%% file: GHH26-leap.tex
\documentclass[a4paper,11pt,reqno]{amsart} 

\usepackage{comment}
\usepackage[utf8]{inputenc}
\usepackage{pifont}
\usepackage[english]{babel}
\usepackage{amsmath}
\usepackage{amsthm}
\usepackage{amsfonts,mathrsfs,relsize,bigints,stmaryrd}
\usepackage{amssymb}
\usepackage{anysize}
\usepackage{hyperref}
\usepackage{esint}
\usepackage{appendix}
\usepackage{mathtools}
\usepackage{graphicx}
\usepackage{tabularx}
\usepackage{xcolor}
\usepackage[normalem]{ulem}
\usepackage{enumitem}

\makeatletter  \@addtoreset{equation}{section} \makeatother
\newtheorem{defi}{Definition}[section]
\newtheorem{lem}{Lemma}[section]
\newtheorem{theo}{Theorem}[section]
\newtheorem{cor}{Corollary}[section]
\newtheorem{pro}{Proposition}[section]

\newtheorem{rem}{Remark}[section]

\DeclareMathOperator{\N}{\mathbb{N}}
\DeclareMathOperator{\R}{\mathbb{R}}
\DeclareMathOperator{\RR}{\mathbb{R}}
\DeclareMathOperator{\C}{\mathbb{C}}
\DeclareMathOperator{\T}{\mathbb{T}}

\DeclareMathOperator{\Z}{\mathbb{Z}}

\newcommand{\cF}{\mathcal{F}}

\newcommand{\cL}{\mathcal{L}}

\allowdisplaybreaks 

\usepackage{subcaption}

\thanks{C.G. has been supported by RYC2022-035967-I (MCIU/AEI/10.13039/501100011033 and
FSE+), and partially by Grants PID2022-140494NA-I00 and PID2022-137228OB-I00 funded
by MCIN/AEI/10.13039/501100011033/FEDER, UE, by Grant C-EXP-265-UGR23 funded by
Consejeria de Universidad, Investigacion e Innovacion \& ERDF/EU Andalusia Program, and
by Modeling Nature Research Unit, project QUAL21-011.
Z.H. has been supported by Grant RYC2023-043706-I funded by MICIU/AEI/10.13039/501100011033 and by ESF+ and by grant PID2022-137228OB-I00 funded
by MCIN/AEI/10.13039/501100011033/FEDER,
T. H. has been supported by Tamkeen under the NYU Abu Dhabi Research Institute grant. Proyecto realizado con la Beca Leonardo a
Investigadores y Creadores Culturales 2024 de la Fundaci\'on BBVA}

\begin{document}

\title[Time-periodic leapfrogging vortex rings in the 3D Euler equations]{Time-periodic leapfrogging vortex rings in \\ the 3D Euler equations}

\author[C. Garc\'ia]{Claudia Garc\'ia}
\address{ Departamento de Matem\'atica Aplicada \& Research Unit ``Modeling Nature'' (MNat), Facultad de Ciencias, Universidad de Granada, 18071 Granada, Spain}
\email{ claudiagarcia@ugr.es}
\author[Z. Hassainia]{Zineb Hassainia}
\address{Departamento de Matem\'atica Aplicada \& Research Unit ``Modeling Nature'' (MNat), Facultad de Ciencias, Universidad de Granada, 18071 Granada, Spain}
\email{zinebhassainia@ugr.es}
\author[T. Hmidi]{Taoufik Hmidi}
\address{New York University Abu Dhabi, Abu Dhabi, United Arab Emirates.}
\email{th2644@nyu.edu}

\begin{abstract}
    We prove the existence of time-periodic leapfrogging vortex rings for the three-dimensional incompressible Euler equations, thereby providing a rigorous realization of a phenomenon first conjectured by Helmholtz (1858). In the leapfrogging motion, two coaxial vortex rings periodically exchange positions, a striking behavior repeatedly observed in experiments and numerical simulations, yet lacking complete mathematical justification.
Our construction relies on a desingularization of two interacting vortex filaments within the contour dynamics formulation, which yields a Hamiltonian description of nearly concentric vortex rings. The main difficulty stems from a singular small-divisor problem arising in the linearized transport dynamics, where the effective time scale degenerates with the ring thickness parameter. To overcome this obstruction, we develop a degenerate KAM-type analysis combined with pseudo-differential operator techniques to control the linearized dynamics around symmetric configurations.
Combining these tools with a Nash–Moser iteration scheme, we construct families of nontrivial time-periodic solutions in an almost uniformly translating frame. This establishes the first rigorous construction of classical leapfrogging motion for axisymmetric Euler flows without swirl, with no restriction on the time interval of existence.
\end{abstract}

\maketitle

\tableofcontents

\section{Introduction}
In this paper, we are concerned with the dynamics of vortex rings, a fundamental  class of coherent structures in fluid dynamics that have been the subject of extensive study for more than a century and a half. A single vortex ring is characterized by the concentration of vorticity around a closed loop, typically a circle in the axisymmetric setting, giving rise to a toroidal region of rotating fluid that propagates steadily through the ambient medium. Vortex rings have been experimentally observed since the nineteenth century, notably in the classical works of Helmholtz \cite{Helmholtz1858,Helmholtz} and Kelvin \cite{Kelvin1867}, and they constitute canonical examples of persistent structures in inviscid flows. Their remarkable stability and self-propelling behavior make them a cornerstone in the study of vortex dynamics.
The reduction of the three-dimensional incompressible Euler equations to the axisymmetric formulation without swirl preserves the Hamiltonian structure of the system, providing a natural and powerful framework to analyze the dynamics of such configurations. Within this setting, the motion of vortex rings can be viewed as the evolution of concentrated vorticity distributions governed by contour dynamics or filament models, depending on the degree of concentration. Among the various interaction patterns between vortex rings, one of the most striking and visually compelling is the leapfrogging phenomenon, in which two coaxial rings of similar strength periodically exchange their positions while maintaining an approximately steady shape. This periodic exchange; where the rear ring contracts, accelerates, and passes through the front one, which simultaneously expands and decelerates, has been extensively documented in experiments and numerical simulations. 
This behavior was first identified by Helmholtz in his pioneering work on vortex dynamics \cite{Helmholtz1858,Helmholtz}, and later modeled by Dyson \cite{Dyson1, Dyson2}, who derived a simplified system describing the interaction of two coaxial vortex rings. In his formulation, each ring is idealized as a Dirac mass in the axisymmetric variables, that is, as a circular vortex filament, leading to the observation of leapfrogging at the level of point-vortex dynamics. The phenomenon has also been reproduced experimentally, for instance, through the generation of air or smoke vortex rings in controlled laboratory conditions \cite{Lim, Maxworthy, Riley}.
The rigorous justification of this phenomenon at the level of classical solutions to the Euler equations remains a challenging problem, and only a few partial results are available in the literature. We refer, for instance, to the works of Dávila–Del Pino–Musso–Wei \cite{delpino}, Buttà–Cavallaro–Marchioro \cite{Butta-Cavallaro-Marchioro}, and Donati–Hientsch–Lacave–Miot \cite{Miot-Lacave}, where leapfrogging dynamics are observed only for local time. Further discussion of these contributions will be provided later. The leapfrogging motion is not specific  to the Euler equations; it also arises in quantum fluid models, most notably in  the Gross–Pitaevskii equation~\cite{smets, smets2}.\\
In contrast, the present work provides the first rigorous analytical construction of long-time leapfrogging motions in the setting of classical solutions to the three-dimensional incompressible Euler equations. More precisely, we construct time-periodic solutions exhibiting leapfrogging behavior within a translating reference frame.
\\
We  consider the dynamics of an incompressible, inviscid fluid in three dimensions, governed by the Euler equations,
\begin{equation}\label{eq:euler}
\partial_t u + u \cdot \nabla u = -\nabla p,\,  
\qquad \nabla \cdot u = 0,\quad x\in\R^3, t\in\R,
\end{equation}
where $u=u(t,x)\in\R^3$ denotes the velocity field and $p=p(t,x)$ the scalar pressure.
A fundamental quantity associated with the Euler equations is the vorticity, defined by
\begin{equation*}
\omega = \nabla \times u.
\end{equation*}
Taking the curl of \eqref{eq:euler} yields the vorticity formulation of the Euler equations,
\begin{equation*}
\partial_t \omega + u \cdot \nabla \omega = \omega \cdot \nabla u,
\end{equation*}
which highlights the transport and stretching mechanisms governing the evolution of vorticity. Notice that the velocity field $u$ can be recovered from the vorticity $\omega = \nabla \times u$ 
through the {Biot--Savart law}
\begin{align*}
u(x) &= \frac{1}{4\pi} \int_{\R^3} \frac{(x-y) \times \omega(y)}{|x-y|^3}\, dy=\int_{\R^3} K(x-y) \times \omega(y)\, dy,
\end{align*}
with $K$ the Biot--Savart kernel
$$
K(x) = \frac{1}{4\pi}\frac{x}{|x|^3}.
$$
The question of global well-posedness for smooth solutions remains open, even though there is numerical evidence of blow up in finite time \cite{Hou2, Hou, Wang, GS-etal}, blow up with boundaries \cite{Chen, Chen2}, and blow up in H\"older spaces \cite{Cordoba, Elgindi}. Nevertheless, an important class of global solutions emerges under the assumption of \emph{axisymmetry without swirl}, that is, when the velocity field is independent of the azimuthal angle $\phi$ in cylindrical coordinates $(r,\phi,z)$ and its angular component vanishes:
$$
u = u^r(r,z,t)\, e_r + u^z(r,z,t)\, e_z.
$$
In this setting, the vorticity possesses only an angular component,
$$
\omega = \omega^\phi(r,z,t)\, e_\phi=(\partial_z u^r-\partial_r u^z)e_\phi,
$$
and the dynamics reduce to the transport equation for the potential vorticity $q = \frac{\omega^\phi}{r}$:
\begin{equation}\label{eq:axisym-vorticity}
\partial_t q 
+ u^r \partial_r q
+ u^z \partial_z q = 0,
\end{equation}
coupled with the velocity reconstruction:
$$
-\Big(\partial_r^2 + \partial_z^2 - \tfrac{1}{r^2}\Big)\psi = r^2 q, 
\qquad 
u^r = -\tfrac{1}{r}\partial_z \psi, 
\qquad 
u^z = \tfrac{1}{r}\partial_r \psi,
$$
where $\psi$ denotes the stream function. 
This property implies the conservation of all the $\|q\|_{L^{p}(\R^{3})}$ norms. 
In~\cite{UkhovskiiYudovich1968}, Ukhoviskii and Yudovich exploited these conservation laws to establish global well-posedness for initial data 
$\omega_{0}, q_{0} \in L^{2}(\R^{3}) \cap L^{\infty}(\R^{3})$. 
More refined results on global existence in various functional settings have subsequently been obtained in~\cite{Abidi, Danchin, Raymond, Shirota-Yanagisawa}.
\\
We now express the three-dimensional axisymmetric Euler equations in Hamiltonian form. This can be achieved by introducing a suitable change of variables defined as follows,
\begin{equation}\label{intro:change:tilde-q}
\varrho=\tfrac{1}{2}r^2, \qquad q(t,r,z)={\bf q}(t,\varrho,z), \qquad \psi(t,r,z)=\Psi(t,\varrho,z).
\end{equation}
In these variables, the derivatives transform as follows
$$
\partial_t q=\partial_t {\bf q},\qquad \partial_r q=r\partial_\varrho {\bf q},\qquad \partial_z q=\partial_z {\bf q}.
$$
Consequently, equation \eqref{eq:axisym-vorticity} takes the form
\begin{equation}\label{intro:eq:tilde-q}
\partial_t {\bf q}+\nabla^\perp \Psi\cdot \nabla {\bf q}=0,
\end{equation}
with $\nabla^\perp=(-\partial_z,\partial_\varrho)$.  
More information about this transformation and the associated Green kernel is given in Section \ref{sec:hamiltonian}.
\\
We recall that in classical studies, see, for instance, the seminal works of Kelvin, Hicks, Lamb, Saffman, and Fraenkel~\cite{Hicks,Kelvin1867,Lamb,Saffman,Fraenkel}, the dynamics of vortex rings with thin cores, whose cross sections are modeled by small circular discs, are often approximated by vortex filaments governed by a coupled system of ODEs.
In contrast, since we work within the Hamiltonian formulation~\eqref{intro:eq:tilde-q}, we shall desingularize these vortex filaments by replacing them with vortex rings possessing, at first approximation, elliptic cores, which is due to  the nonhomogeneous deformation induced by the change of variables.
\\
In the next section, we state our main theorem and present the principal ideas underlying its proof. Later, in Section \ref{intro:previousworks} we give a discussion on previous works related to vortex rings, recalling several classical and fundamental results. Finally, Section \ref{intro:sketch} provides a complete sketch of our main theorem and Section \ref{sec:intro4} gives a comparison of our work with other related works.

\subsection{Main contribution}\label{sec:intromaintheorem}

In this work, we provide a rigorous construction of time--periodic leapfrogging solutions for the 3D Euler equations under the assumption of axisymmetry without swirl. Our analysis builds upon several  techniques from nonlinear analysis and dynamical systems. First, we employ a desingularized formulation through \emph{contour dynamics}, which allows us to replace the singular filament dynamics by an effective system describing the evolution of concentrated vortex patches. This framework enables the derivation of a Hamiltonian system that captures the core interactions of nearly concentric vortex rings.

A central analytical difficulty in our work is to prove the persistence of periodic leapfrogging motions when the system is subject to perturbations. The equations governing the evolution of the vortex boundaries take the form of singular nonlinear transport equations, in which a subtle form of time degeneracy significantly complicates the analysis. A priori, the system involves a coupling between two interacting equations, reflecting the mutual influence of the rings. However, after performing an appropriate symmetry reduction, the problem can be rewritten as a single scalar PDE that includes a time delay term arising from the interaction with the second ring. Within this framework, the vortex rings themselves act as an external periodic forcing mechanism, shaping  the dynamics through the coefficients of the governing PDEs.
Determining whether such a PDE admits time-periodic solutions is precisely the core of our analysis.
Carrying out this program requires overcoming several difficulties related to resonances and small-divisor phenomena, which become even more delicate in a degenerating time regime. These obstructions are intrinsic to the problem and necessitate the careful combination of degenerate KAM techniques, pseudo-differential analysis, and a Nash–Moser scheme.

To state our main result, we begin by introducing an asymptotic model that captures the essential dynamics of the cores of two vortex rings undergoing leapfrogging motion. This reduced formulation isolates the leading‑order interactions responsible for the periodic exchange of positions between the rings and serves as the basis for our analysis. A complete derivation of the model, including the symmetry reductions and scaling arguments, will be provided in Section~ \ref{sec-sym-red} and Section~\ref{section:Symmetry reduction}.
For a small parameter $\varepsilon>0$, the dynamics of the two vortex filaments in the $(\rho,z)$‑plane are represented by the trajectories
$$
t \mapsto \big(P_1(\tau),\, P_2(\tau)\big),\qquad \tau:=|\ln \varepsilon|\, t,\quad P_j=(p_{j,1},p_{j,2}), \quad j=1, 2,
$$
governed by the Hamiltonian system
\begin{equation}\label{intro-eq-P1P2}
|\ln\varepsilon| \dot{P}_j(\tau)=\nabla^\perp_{P_j} H(P_1,P_2),
\end{equation}
where  and the Hamiltonian is given by
\begin{align*}
H(P_1,P_2)&:=\frac{(p_{1,1}p_{2,1})^{\frac14}}{\sqrt{2}}\,
J\!\left(\frac{2(\sqrt{p_{1,1}}-\sqrt{p_{2,1}})^2+(p_{1,2}-p_{2,2})^2}{2\sqrt{p_{1,1} p_{2,1}}}\right)\\ &\qquad\qquad+\sum_{j=1}^2\frac{\sqrt{p_{j,1}}}{\sqrt{8}}\left(|\ln\varepsilon|{-}\frac74{+}\frac54\ln(8){+}\frac34\ln(p_{j,1})\right),
\end{align*}
and the auxiliary function $J$ is defined as
$$
J(s):=\int_0^{\pi}\frac{\cos\theta}{\sqrt{s+2-2\cos\theta}}\, d\theta, 
\qquad s>0.
$$
We assume that the initial configuration  of the vortex filaments  is vertical and
\begin{equation}\label{initial-condition-points}
   P_1(0)=\big(\kappa+{\tfrac12}|\ln\varepsilon|^{-\frac12} (2\kappa)^{\frac12}\lambda,0\big),\quad P_2(0)=\big(\kappa-{\tfrac12}|\ln\varepsilon|^{-\frac12} (2\kappa)^{\frac12}\lambda,0\big),
\end{equation}
for some $\kappa>0$ and $\lambda>0$. As shown in
Section~\ref{sec-sym-red}, the center of mass
$\mathtt{C}=\tfrac12(P_1+P_2)$ evolves along the vertical line
$\rho=\kappa$ and admits the decomposition
$$
\mathtt{C}(\tau)
=
\left(\kappa,0\right)
+
\big(0,U_\varepsilon\,\tau+W_\varepsilon(\tau)\big),
$$
where
$$
U_\varepsilon
=\tfrac{1}{4\sqrt{2\kappa}}
+O\!\left(\tfrac{\ln|\ln\varepsilon|}{\ln\varepsilon}\right),
$$
is a constant drift speed. A better asymptotic expansion of $U_\varepsilon$ is given in Corollary \ref{coro-speed}. In particular, when expressed in the physical
time variable $t$, this velocity exhibits a $|\ln\varepsilon|$-singular
growth. The correction function $W_\varepsilon=O\left(|\ln\varepsilon|^{-1}\right)$ is
$T$-periodic (see
Corollary~\ref{coro-speed}). This period $T=T(\varepsilon,\kappa,\lambda)$ is uniform in $\varepsilon$, and it coincides   with the common period of the quantities
$$
p_{1,1},\qquad p_{2,1},\qquad p_{2,2}-p_{1,2}. 
$$
\begin{figure}
 \centering
\begin{subfigure}{.3\textwidth}
  \centering
  \includegraphics[width=1\linewidth]{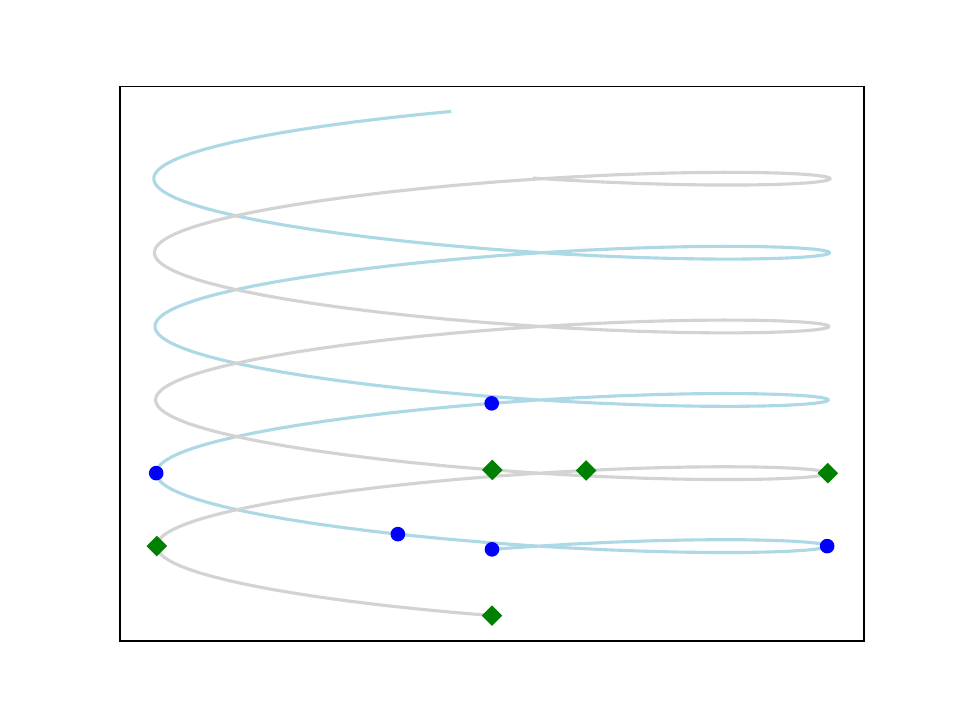}
\end{subfigure}%
\begin{subfigure}{.3\textwidth}
  \centering
  \includegraphics[width=1\linewidth]{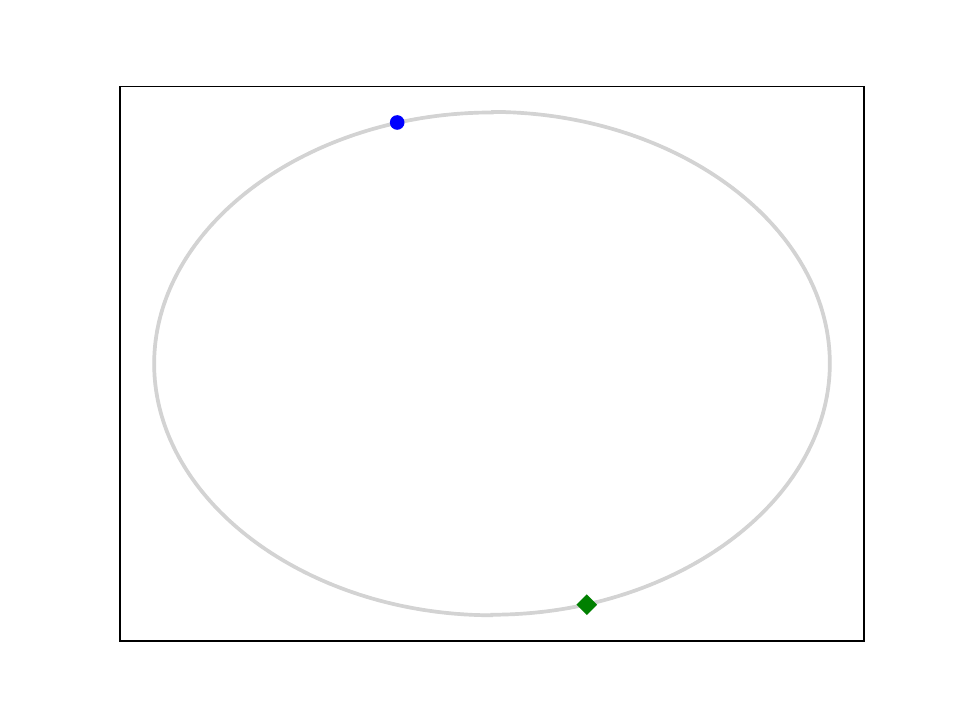}
\end{subfigure}
\caption{Evolution of the points $(P_1,P_2)$ for the parameters $\kappa=0.4$ and $\varepsilon=0.05$.  The green diamonds and blue points represent the vortex ring positions $P_1$, $P_2$, respectively. The right-hand panel displays the evolution of $(P_1, P_2)$
 in the translating frame of reference, where a periodic motion is observed.}
\label{fig-P1P2}
\end{figure}

\begin{figure}
 \centering
\begin{subfigure}{.3\textwidth}
  \centering
  \includegraphics[width=1\linewidth]{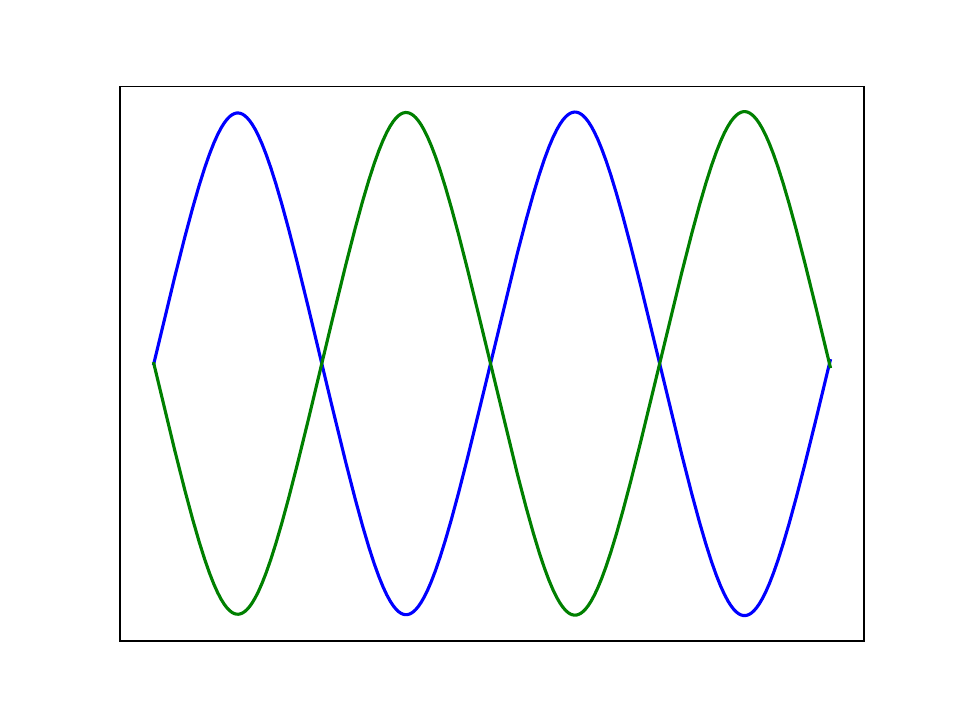}
\end{subfigure}
\begin{subfigure}{.3\textwidth}
  \centering
  \includegraphics[width=1\linewidth]{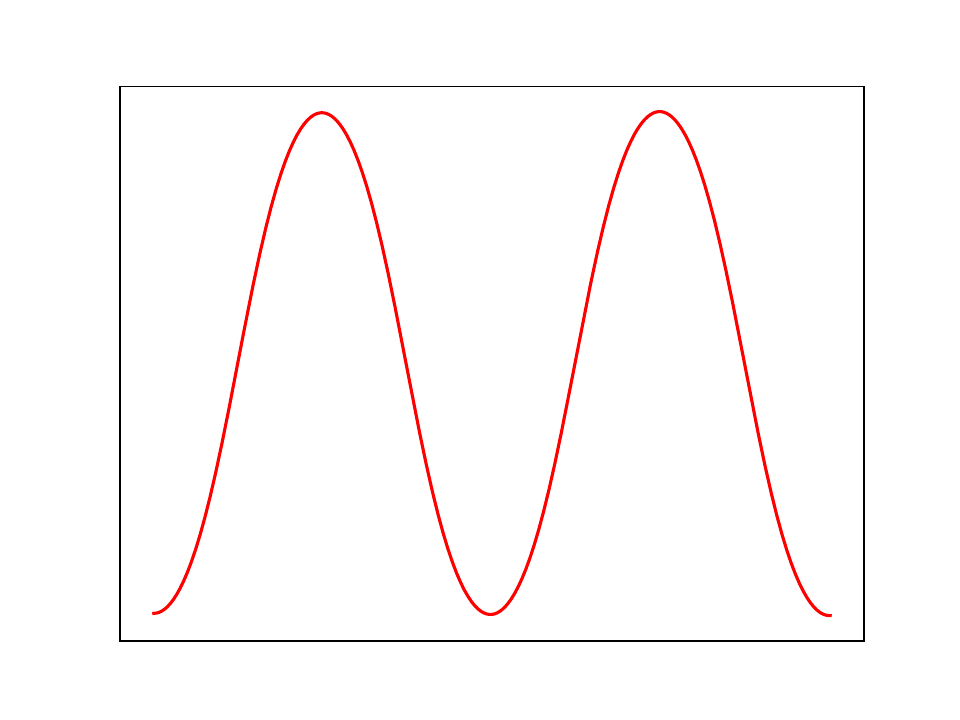}
\end{subfigure}
\caption{Parameters as in Figure~\ref{fig-VR}, namely $\kappa=0.4$ and $\varepsilon=0.05$. he left-hand panel shows the periodic evolution of $\tau\mapsto P_{1,1}(\tau)$ and $\tau\mapsto P_{2,1}(\tau)$while the right-hand panel displays the periodic function $\tau\mapsto (P_{2,2}-P_{1,2})(\tau)$.}
\label{fig-P11P21Z}
\end{figure}
When reverting to the original physical time variable t, this period becomes $\tfrac{T}{|\ln \varepsilon|}.$ This reflects the fact that the two vortex rings are advected with large self-induced velocity and undergo rapid rotational motion.

To illustrate this, we present numerical simulations carried out using the parameters
$\kappa=0.4$ and $\varepsilon=0.05$.
Figure~\ref{fig-P1P2} displays the evolution of the points $(P_1,P_2)$ governed by \eqref{intro-eq-P1P2}. In the left-hand panel, one observes the evolution modulo translation, revealing the characteristic periodic leapfrogging motion. The right-hand panel shows the same trajectories in a translating reference frame.
In Figure~\ref{fig-P11P21Z} we plot the periodic functions $\tau\mapsto p_{1,1}(\tau)$, $\tau\mapsto p_{1,2}(\tau)$ and $\tau\mapsto (p_{2,2}-p_{1,2})(\tau)$.

We now state our main result, which provides a rigorous construction of time-periodic solutions to the axisymmetric Euler equations whose dynamics desingularize and faithfully reproduce the leapfrogging motion of two interacting vortex filaments   $(P_1,P_2)$.

\begin{theo}\label{th-main1}
 Let $0<a<b$ and $\kappa>0$ be fixed. Then there there exists $\varepsilon_0>0$ such that,  every $\varepsilon\in(0,\varepsilon_0)$, there is a Borel set $\mathcal{C}_\varepsilon\subset (a,b)$ of asymptotically full Lebesgue measure, 
 $$
 \lim_{\varepsilon\to 0}|\mathcal{C}_{\varepsilon}|=b-a,
 $$
with the  following property:  for every $\lambda \in \mathcal{C}_\varepsilon$, there exists a global solution of the axisymmetric Euler equation \eqref{intro:eq:tilde-q} of the form
\begin{equation*}
{\bf q}(t,\varrho,z)=\tfrac{1}{\varepsilon^2}{\bf 1}_{\mathcal{D}_1(t)}(\varrho,z)+\tfrac{1}{\varepsilon^2}{\bf 1}_{\mathcal{D}_2(t)}(\varrho,z), 
\end{equation*}
where the two bounded domains are given by
\begin{equation*}
\mathcal{D}_j(t)=\left({P}_j+\varepsilon \mathbf{D}_j+i \varepsilon|\ln\varepsilon|^{-1}\mathbf{V}_j\right)(|\ln \varepsilon|t), \qquad j=1,2.
\end{equation*}
Here, each $\mathbf{D}_j$ is a simply connected domain satisfying the time-periodicity condition
$$
\forall \tau\in\mathbb{R}, \quad \mathbf{D}_j(\tau+T)=\mathbf{D}_j(\tau), 
$$
and whose  boundaries admit the parametrization
\begin{align*}
 \gamma_{1}(\tau,\theta)=&\sqrt{1+2\,\varepsilon\, \mathtt{f}( \tau,\theta) } \Big((2p_{1,1} (\tau))^\frac14\cos(\theta),\, (2p_{1,1} (\tau))^{-\frac14}\sin(\theta)\Big),\\
  \gamma_{2}(\tau,\theta)=&\gamma_{1}\big(\tau+\tfrac{T}{2},\theta\big).
 \end{align*}
The map $\tau \mapsto \mathtt{f}(\tau,\cdot)\in H^s(\mathbb{T})$ is $T$‑periodic for some sufficiently  large Sobolev index $s$,  and the real-valued  functions $\tau \mapsto\mathbf{V}_j:=\mathbf{V}_j(\varepsilon,\mathtt{f})(\tau)$ represent modulation velocities depending on the profile $\mathtt{f}$. Here,  $T$ denotes the period of the vortex‑filament system \eqref{intro-eq-P1P2}.

\end{theo}

\begin{figure}
 \centering
\begin{subfigure}{.19\textwidth}
  \centering
  \includegraphics[width=0.9\linewidth]{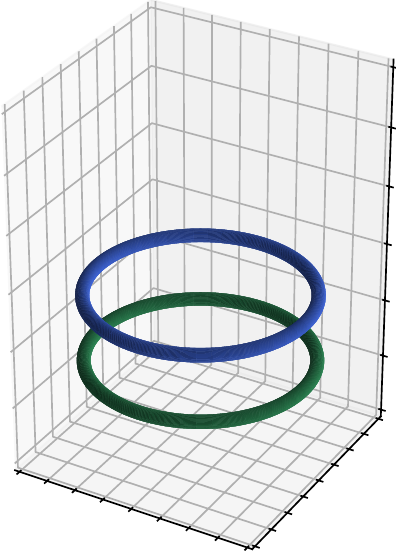}
\end{subfigure}
\begin{subfigure}{.19\textwidth}
  \centering
  \includegraphics[width=.9\linewidth]{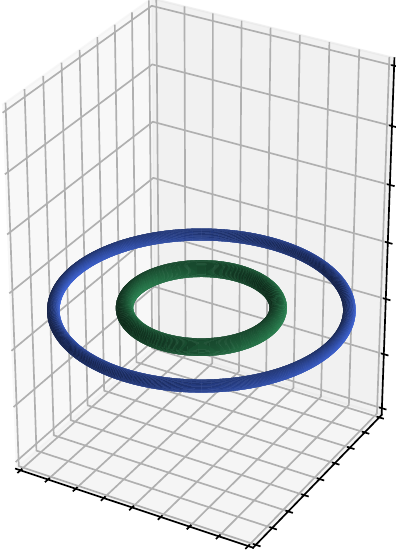}
\end{subfigure}
\begin{subfigure}{.19\textwidth}
  \centering
  \includegraphics[width=.9\linewidth]{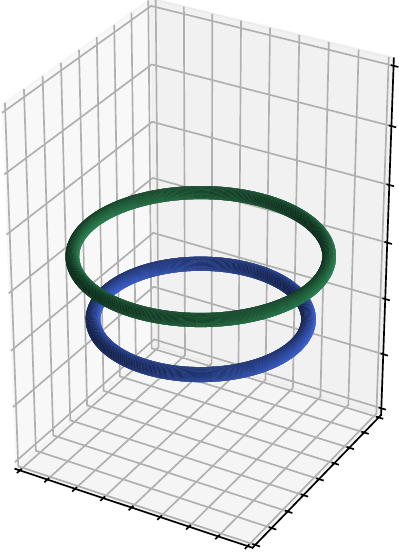}
\end{subfigure}
\begin{subfigure}{.19\textwidth}
  \centering
  \includegraphics[width=.9\linewidth]{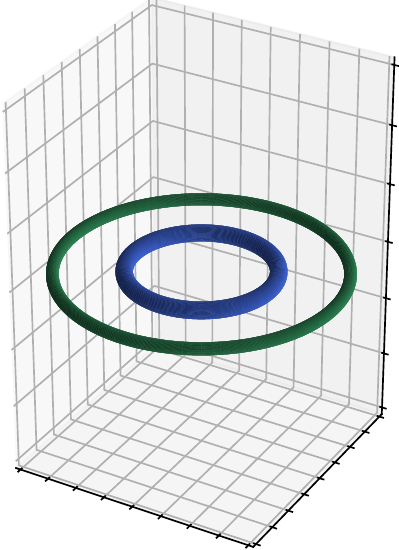}
\end{subfigure}
\begin{subfigure}{.19\textwidth}
  \centering
  \includegraphics[width=.9\linewidth]{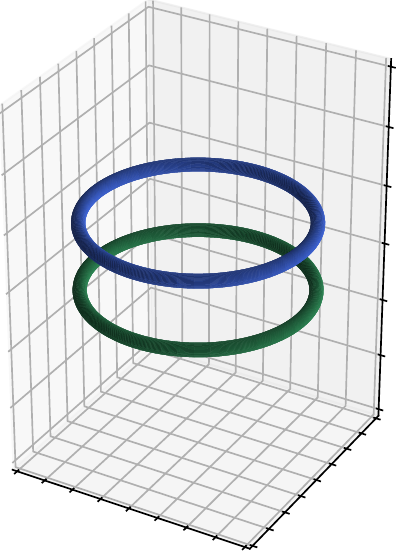}
\end{subfigure}
\caption{Periodic evolution of the two vortex rings described in Theorem~\ref{th-main1}, with parameters $\kappa = 0.4$ and $\varepsilon = 0.05$. The left‑hand panel shows the initial configuration, while the right‑hand panel depicts the state at time $t = T$.}
 \label{fig-VR}
\end{figure}

To illustrate the conclusions of our main theorem, we present numerical simulations carried out for a small value of $\varepsilon$, using the parameters
$\kappa=0.4$ and $\varepsilon=0.05$. In Figure~\ref{fig-P11P21Z} we plot the periodic functions $\tau\mapsto p_{1,1}(\tau)$, $\tau\mapsto p_{1,2}(\tau)$ and $\tau\mapsto (p_{2,2}-p_{1,2})(\tau)$. Figure~\ref{fig-VR} depicts the vortex rings parameterized according to Theorem~\ref{th-main1}, using the choice $\mathtt{f}=\tfrac18 (2p_{1,1})^{-\frac34}\cos(3\theta)$, which corresponds to the leading term in the approximate solution constructed in Section~\ref{sec:approx}.

We now record a few remarks that highlight key features of the result.

\begin{rem}[Non-rigid periodicity]
Although the domain $\mathbf{D}_1$ is time-periodic, it does not rotate rigidly. This is reflected in the  
 structure of the  profile  $\mathtt{f}$, which  admits the asymptotic expansion
$$
\mathtt{f}(\tau,\theta)
=
\tfrac18 (2p_{1,1}(\tau))^{-\frac34}\cos(3\theta)
+
O(\varepsilon|\ln\varepsilon|).
$$
A more refined asymptotic expansion is provided in Theorem $\ref{theo-approx1}.$
\end{rem}

\begin{rem}[Elliptical perturbation]
The domain $\mathbf{D}_1$ is localized around a time-periodic elliptical region whose boundary is parametrized by
$$
(\tau,\theta)\longmapsto 
\Big((2p_{1,1}(\tau))^{\frac14}\cos\theta,\,
(2p_{1,1}(\tau))^{-\frac14}\sin\theta\Big).
$$
This should be contrasted with the classical description of thin-core vortex rings, whose cross-sections are typically modeled by small circular discs (see, e.g., \cite{Hicks,Kelvin1867,Lamb,Saffman,Fraenkel}). 
In our setting, the appearance of an elliptical geometry  is not intrinsic to the physical configuration, but rather results from the nonhomogeneous deformation induced by the change of coordinates~\eqref{intro:change:tilde-q}, which maps at the leading order a disc into an ellipse. This transformation is naturally dictated by the Hamiltonian formulation~\eqref{intro:eq:tilde-q}.
\end{rem}

\begin{rem}[Speed modulation]

Exploiting the translational invariance of the system along the vertical axis, we introduce a modulation of the propagation speed ${\bf V}_j$, depending on $(\varepsilon,\mathtt{f})$, in order to partially remove the degeneracy of the linearized operator at the first Fourier mode. The modulation is chosen so as to eliminate the first sine Fourier mode from the nonlinear equation. 
However, a degeneracy at leading order in $\varepsilon$ persists for the corresponding cosine mode. This difficulty manifests itself in the invertibility properties of the linearized operator. To address this issue, a more refined expansion in $\varepsilon$ is required in order to identify whether the degeneracy can be lifted at smaller scales. In fact, the inversion of the cosine mode reduces to solving a complicated nonlocal integro-differential equation. For that, we do an spectral analysis to conclude that such nonlocal integro-differential equation can be inverted, see Section \ref{Sec-Fundam-lemma}.
\end{rem}

\subsection{Previous results on vortex rings}\label{intro:previousworks}
The first approximation in the study of concentrated vortex rings within the framework of the inviscid three-dimensional Euler equations consists in modeling them as vortex filaments, that is, singular vortex structures supported on one-dimensional curves.
In this setting, analytical developments stemming from the pioneering works of Helmholtz~\cite{Helmholtz1858} and Kelvin~\cite{Kelvin1867}, and pursued by Pocklington~\cite{Pocklington1893}, Love~\cite{Love1894}, and Hicks~\cite{Hicks}, among others, led to approximate models describing the mutual interaction of coaxial vortex filaments. More refined asymptotic and dynamical models for thin vortex filaments were subsequently derived by Moore and Saffman~\cite{Moore-Saffman} and by Klein, Majda, and Damodaran~\cite{KleinMajdaDamodaran1995}. The connection between such reduced filament dynamics and the full three-dimensional Euler system was  established by Marchioro and Negrini~\cite{Marchioro-Negrini}.

A prototypical steady vortex structure in three dimensions is Hill’s spherical vortex~\cite{Hill}.
Subsequent perturbative and variational constructions have produced distinguished families of steadily translating vortex rings, notably those obtained by Fraenkel~\cite{Fraenkel}, Fraenkel–Berger~\cite{Franken1}, and Norbury~\cite{Norbury}.
The stability properties of some of these structures have been investigated in works by Widnall-Sullivan \cite{Widnall}, Choi \cite{Choi}, Cao {\it et al.} \cite{Daomin}, Guo {\it et al} \cite{Guo}.
In perturbative constructions of traveling vortex rings obtained through  the desingularization of circular vortex filament, the propagation speed exhibits a logarithmic singular dependence on the small core radius of the ring.

The dynamics of multiple thin vortex rings and the associated Hamiltonian
structure were developed by Vasilyev~\cite{Vasilyev}, Dyson~\cite{Dyson2},
and Borisov {\it et al.}~\cite{Borisov}. In this framework, the interaction
of idealized rings is described by a highly nonlinear finite-dimensional
Hamiltonian system. Although the resulting dynamics can be complex,
a special regime allows for the emergence of the classical leapfrogging
motion, in which two coaxial rings periodically exchange positions.
Numerical simulations have further supported the persistence of
leapfrogging behavior for sufficiently concentrated vortex rings over
long time intervals~\cite{Acheson1972, ShariffLeonard1992}.

The asymptotic regimes underlying the leapfrogging motion depend in a delicate manner on the relative scaling between the core radius, the filament thickness, and their mutual separation. In the case of two vortex filaments, an appropriate scaling in $\varepsilon$ between the core radius and the mutual separation is required in order to balance the logarithmically singular self-induced speed of each Lamb ring with the interaction induced by the other filament. Such a balance is essential for the emergence of periodic motions. As already noted, such a regime was  recognized in the classical works of Helmholtz and others~\cite{Dyson1,Dyson2,Helmholtz1858,Helmholtz,Hicks}. For more recent developments in this direction, we refer the reader to the works of  Buttà–Cavallaro–Marchioro~\cite{Butta-Cavallaro-Marchioro, Cavallaro-Marchioro}. In particular, in~\cite{Butta-Cavallaro-Marchioro} the authors establish a local-time desingularization of the leapfrogging motion in the thin-core regime, where the filament thickness is of order $\varepsilon$ while the ring radius scales as $|\ln \varepsilon|$. In this asymptotic scaling, the evolution of the two vortex rings is well approximated, at leading order, by a rigid rotational motion in a suitably translating frame.

Another regime for the leapfrogging motion was analyzed by Dávila–Del Pino–Musso–Wei~\cite{delpino} and Donati–Hientzsch–Lacave–Miot~\cite{Miot-Lacave}, where the ring radius is of order one, while the separation between the rings scales as $|\ln \varepsilon|^{-\frac{1}{2}}$. This latter configuration will be considered in the present work as a suitable candidate for the desingularization procedure and the implementation of KAM theory, and it is also consistent with the discussion in the work of  Helmholtz~\cite{Helmholtz1858} and other classical works \cite{Hicks,Kelvin1867, Love1894,Pocklington1893}.

It is important to emphasize that all the aforementioned constructions of leapfrogging are local in time. Extending these results to a global-in-time framework remains a major analytical challenge, primarily due to intrinsic instability mechanisms that are difficult to overcome.

The first long-time result addressing the leapfrogging motion was recently established in the time-periodic setting by the last two authors together with Masmoudi in~\cite{HHM21}. In this work, the authors achieved, for the two-dimensional  Euler equations, the desingularization of two symmetric dipoles corresponding to Love’s pairs, relying on the contour dynamics formulation combined with KAM theory and a Nash–Moser iteration scheme.

In this paper, our main goal is to rigorously establish a long-time result for the leapfrogging motion in the three-dimensional axisymmetric Euler equations. As described in the previous Section \ref{sec:intromaintheorem}, we  desingularize the dynamics of two interacting vortex filaments undergoing leapfrogging motion and  construct time-periodic solutions in an almost uniformly translating frame. This provides the first rigorous derivation of such long-time periodic dynamics in this setting. Our analysis relies on a desingularization procedure achieved through  uniform potential vorticity tubes concentrated around vortex rings, governed by the contour dynamics formulation. To carry out the construction within a space–time periodic configuration, we employ degenerate type of  KAM  theory, combined with a Nash–Moser iterative scheme and pseudo-differential techniques involved in the reducibility process. These analytical tools allow us to overcome the small-divisor and loss of derivative difficulties inherent to this type of problem.

KAM theory has a long and fruitful history in the study of Hamiltonian partial differential equations, where it has been successfully applied to models such as the nonlinear Schr\"odinger equation, the water-wave system, and the two-dimensional Euler equations, leading to the construction of finite-dimensional invariant tori in the vicinity of trivial equilibria. In what follows, we provide a comprehensive  list of references related to applications in fluid dynamics, a relatively recent and rapidly developing area, see for instance
 ~\cite{BCP,BHM,BertiBiascoProcesi2013,BertiMontaltoHaus2017,GIP23,  HHM21, HHR24}. 

To conclude this section, we recall that the dynamics of vortex rings has also been investigated in the context of the axisymmetric Navier–Stokes equations, and we restrict ourselves here to a concise but representative list of contributions. For instance, Bedrossian–Germain–Harrop-Griffiths \cite{Bedrossian-Germain-Harrop-Griffiths} established local well-posedness for initial data corresponding to a smooth vortex filament with small perturbations, while Hidalgo–Gancedo \cite{HidalgoGancedo} proved global well-posedness for helical filaments.
In \cite{Gallay-Sverak-2,Gallay-Sverak}  Gallay–Šverák obtained existence and uniqueness results for solutions issued from a circular filament, together with a detailed analysis of the inviscid limit. A distinct desingularization mechanism was explored by \mbox{Fontelos–Vega \cite{Fontelos-Vega},} who showed that the Navier–Stokes evolution starting from a closed filament remains, for short times, close to the binormal curvature flow; in this setting, the resulting structure corresponds to an axisymmetric solution with swirl rather than a pure vortex ring.
More recently, Gancedo–Hidalgo–Torné– Mengual~\cite{HidalgoGancedo-Mengual} employed convex integration methods to construct infinitely many local-in-time weak Euler solutions with initial data given by a circular vortex filament.

\subsection{Sketch of the proof}\label{intro:sketch}

The proof of Theorem~\ref{th-main1} is technically involved and proceeds through several delicate steps; we therefore limit ourselves to a sketch of the main ideas.

The construction is carried out by means of a Nash--Moser iteration, which is required to compensate for the loss of derivatives generated by the small--divisor mechanism. A central ingredient is the construction of an \emph{approximate} right inverse for the linearized operator within an appropriate scale of function spaces. In the present setting, the linearization yields a degenerate quasilinear transport operator whose time-- and space--periodic coefficients are inherited from the underlying point--vortex dynamics.
In order to control the small divisors and simplify the linearized dynamics, we implement a KAM-type reducibility scheme, realized through four successive conjugations. For the Nash--Moser framework and its adaptation to similar quasilinear problems, we refer to~\cite{BB13,BB10}.

We now outline the main steps in the proof of Theorem~\ref{th-main1}.

\medskip\noindent{\bf $\diamondsuit$ 
{\sc Step 1:} Hamiltonian reformulation of the problem.}
\medskip

We begin by recasting the existence of a leapfrogging motion for the 3D axisymmetric equations in a Hamiltonian form. Instead of the standard axisymmetric variables $(r,z)$, we work with the coordinates introduced above, 
$$(\varrho,z) = \Big(\tfrac{1}{2} r^2, z\Big).$$
In these variables, the evolution of the potential vorticity ${\bf q}$ acquires a Hamiltonian structure: the flow is generated by a Hamiltonian vector field and, in particular, preserves the area form in the $(\varrho,z)$--plane. This reformulation is crucial, since it allows us to impose a natural zero--average condition and thereby remove the Fourier zero mode from the nonlinear functional we will solve.

We seek leapfrogging vortex rings that are close to the leapfrogging vortex filaments. Let $P_1$ and $P_2$ be two vortex filaments solving \eqref{intro-eq-P1P2}. In Section \ref{sec-sym-red} we analyze the evolution of $P_1$ and $P_2$ and obtain a leapfrogging dynamics, provided the parameters satisfy a specific $\varepsilon$--dependence. Accordingly, we assume that the vortex rings are distributed as follows:
\begin{align*}
&p_{j,1}
=\kappa+\tfrac12(-1)^{j+1}\, r_\varepsilon\,(2\kappa)^{\frac14}\,\mathtt{x}_1,
\qquad j=1,2, \\
& p_{1,2}-p_{2,2}
= r_\varepsilon\,(2\kappa)^{-\frac14}\,\mathtt{x}_2,
\qquad 
r_\varepsilon := (2\kappa)^{\frac14}\,|\ln\varepsilon|^{-\frac12},
\end{align*}
with $\kappa>0$. Moreover, the initial configuration for the vortex filaments is given by \eqref{initial-condition-points}. We show that $(\mathtt{x}_1,\mathtt{x}_2)$ are periodic functions of time, and hence $p_{1,1}$, $p_{2,1}$, and $p_{1,2}-p_{2,2}$ are $T$--periodic, for some $T>0$. The corresponding initial configuration for the filaments is prescribed by \eqref{initial-condition-points}.
for $j=1,2$ and $\kappa>0.$ We can show that $(\mathtt{x}_1,\mathtt{x}_2)$ are $T$-periodic functions, which implies that $p_{1,1}$, $p_{2,1}$ and $p_{2,2}-p_{1,2}$ are $T$-periodic functions, for some $T>0$.

\begin{figure}[h!]
\centering
\def\svgwidth{0.6\textwidth} 
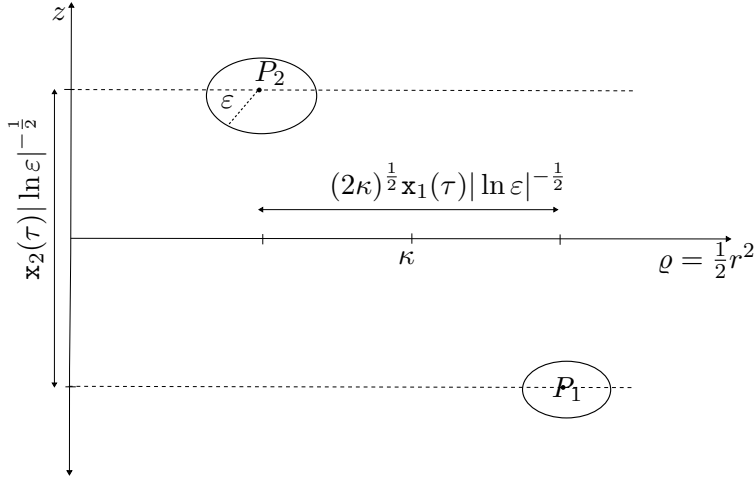
\caption{Illustration of our solution ansatz in the coordinates $(\varrho,z)$. The separation between the centers of the two rings in both the $\varrho$- and $z$‑directions is of order $|\ln\varepsilon|^{-1/2}$, while each vortex‑ring cross‑section has area of order $\varepsilon^{2}$. The corresponding potential vorticity $\mathbf{q}$ attains amplitude of order $\varepsilon^{-2}$ on each ring.}
\label{fig-1}
\end{figure}

We model ${\bf q}$ as the superposition of two thin vortex rings concentrated near the vortex filaments $(P_1,P_2)$. More precisely, we take
$$
{\bf q}(t,\cdot)=\tfrac{1}{\varepsilon^2}{\bf 1}_{\mathcal{D}_1(t)}+\tfrac{1}{\varepsilon^2}{\bf 1}_{\mathcal{D}_2(t)},
$$
where $\mathcal{D}_1(t)$ and $\mathcal{D}_2(t)$ denote the cross-sectional domains of the rings.
Since the velocity induced by a vortex filament is singular, we recall the fast time variable, introduced before, $$\tau=|\ln\varepsilon| t.$$  In terms of $\tau$, the boundaries of the domains are parametrized by $\gamma_j(\tau,\theta)$, $\theta\in[0,2\pi]$, namely
\begin{align*}
 \gamma_{1}(\tau,\theta)=&\,P_1(\tau)+i\,\varepsilon|\ln\varepsilon|^{-1}{\bf V}_1(\tau)
 +\varepsilon\sqrt{1+2\varepsilon \mathtt{f}(\tau,\theta)}
 \Big((2p_{1,1})^\frac14\cos\theta,\,(2p_{1,1})^{-\frac14}\sin\theta\Big),\\
  \gamma_{2}(\tau,\theta)=&\,P_2(\tau)+i\,\varepsilon|\ln\varepsilon|^{-1}{\bf V}_2(\tau)
 +\varepsilon\sqrt{1+2\varepsilon \mathtt{f}\big(\tau+\tfrac{T}{2},\theta\big)}
 \Big((2p_{1,1})^\frac14\cos\theta,\,(2p_{1,1})^{-\frac14}\sin\theta\Big),
\end{align*}
for a shape perturbation $(\tau,\theta)\mapsto \mathtt{f}(\tau,\theta)$, assumed to be $T$-periodic in time and $2\pi$-periodic in space. 
In Section~\ref{SEc-Freq-Norma}, we perform a time reparametrization 
$\tau \mapsto \varphi$ in order to transform $T$-periodic functions 
into $2\pi$-periodic ones in time. Defining
$$
f(\varphi,\theta):=\mathtt{f}(\tau,\theta),
$$
so that $f$ is $2\pi$-periodic with respect to $\varphi$, 
and substituting this expression into the contour dynamics equation 
arising from \eqref{intro:eq:tilde-q}, using the symmetry reduction from Section \ref{Sym-Redu},  we obtain that the pair 
$(\varepsilon,f)$ satisfies the nonlinear functional equation
$$
{\bf F}(\varepsilon,\mathtt{V}_1,\mathtt{V}_2,f)=0,
$$
where 
\begin{align*}
& {\bf F}(\varepsilon,\mathtt{V}_1,\mathtt{V}_2,f)(\varphi,\theta) :=\varepsilon^3 |\ln\varepsilon | \omega \partial_\varphi f(\varphi ,\theta)+\big(2\mathtt{p}_{1,1}(\varphi)\big)^{-\frac12} \partial_\theta \big\{\Psi({\gamma}(\varphi,\theta))\big\}\\ &\quad  -\omega \big( |\ln \varepsilon|\dot{\mathtt{P}_1}(\varphi)+i \,\varepsilon \dot{\mathtt{V}}_1(\varphi)\big)\cdot i\partial_\theta \gamma(\varphi,\theta)+\varepsilon^2 |\ln\varepsilon |\,\omega \tfrac{ \dot{\mathtt{p}}_{1,1}(\varphi)}{8\mathtt{p}_{1,1}(\varphi)}\partial_\theta\big[ (1+2\varepsilon f(\varphi,\theta))\sin(2\theta)],    
\end{align*}
and
\begin{align*}
 &{\pi\sqrt{2}}\, \Psi(\gamma(\varphi,\theta))=\int_{0}^{2\pi}\int_0^{w(\varphi,\eta)}G\big(\mathtt{P}_1(\varphi)+\varepsilon w(\varphi,\theta) \mathcal{Z}(\varphi,\theta),\mathtt{P}_1(\varphi)+\varepsilon \rho   \mathcal{Z}(\varphi,\eta)\big)\rho d\rho d\eta\\
&+\int_{0}^{2\pi}\int_0^{w_\star(\varphi,\eta)}G\big(\mathtt{P}_1(\varphi)+i\,\varepsilon |\ln\varepsilon|^{-1} (\mathtt{V}_1-\mathtt{V}_2)(\varphi)+\varepsilon w(\varphi,\theta) \mathcal{Z}(\varphi,\theta),\mathtt{P}_2(\varphi)+\varepsilon \rho   \mathcal{Z}_\star(\varphi,\eta)\big)\rho  d\rho d\eta.
\end{align*}
Here
\begin{align*}
f_\star(\varphi,\theta)&:=f(\varphi+\pi,\theta),\qquad\,
\mathtt{P}_j(\varphi):=P_j(\tau),\qquad\,
\mathtt{V}_j(\varphi):=\mathbf{V}_j(\tau),
\\
w&:=\sqrt{1+2\varepsilon f},
\qquad
\mathcal{Z}(\varphi,\theta):=\Big((2\mathtt{p}_{1,1})^\frac14\cos\theta,\,(2\mathtt{p}_{1,1})^{-\frac14}\sin\theta\Big),
\end{align*}
and $\omega=\tfrac{2\pi}{T}$ is the frequency of the vortex filament pair $(P_1,P_2).$

The linearization of ${\bf F}$ with respect to $f$ is computed in Proposition~\ref{prop:asymp-lin} and is given by
  \begin{align*}
     &\partial_f {\bf F}(\varepsilon,\mathtt{V}_1,\mathtt{V}_2,f)[h]=\varepsilon^3 |\ln\varepsilon|\omega\partial_\varphi h+ \tfrac{1}{2}\varepsilon \left(\mathcal{H}[h]-\partial_\theta h\right)+\varepsilon^2 \partial_\theta\big\{   V(\varepsilon,f)\, h\big\}+\varepsilon^2\mathcal{H}_{\mathtt{u},0}[h]\\ &  +\varepsilon^2{\partial_\theta}\mathcal{S}[h]+{\varepsilon^3}|\ln\varepsilon|  {\partial_\theta}\mathcal{Q}[h] -\varepsilon^3\partial_\theta\fint_{\mathbb{T}} W({f})(\varphi,\theta,\eta)\ln\big|\sin\big(\tfrac{\theta-\eta}{2}\big)\big|h(\eta) d\eta+\varepsilon^3\partial_\theta\mathscr{R}_{\infty}[h].
  \end{align*}
Here $\mathcal{H}$ denotes the Hilbert transform acting on the spatial variable;  $\mathcal{Q}$ is a projection operator localized on the spatial mode one; $\mathcal{H}_{\mathtt{u},0} \in\mathrm{OPS}^{0}$ and  $\mathcal{S}, \mathscr{R}_\infty \in\mathrm{OPS}^{-\infty}$ are shift operators; and $V(\varepsilon,f)$ and $W({f})$ are smooth functions depending on $f$. The notation $\textnormal{OPS}^{-m}$ means a pseudo differential operator of order $-m$. Precise definitions of all these operators and functions are provided in Section~\ref{sec-Asym-Line-op}.

\medskip\noindent{\bf $\diamondsuit$ {\sc Step 2:} Speed modulation.}
\medskip

A major difficulty in the analysis of the nonlinear functional 
${\bf F}(\varepsilon,\mathtt{V}_1,\mathtt{V}_2,f)$ stems from the degeneracy of the first Fourier mode. 
Indeed, the functions $\sin\theta$ and $\cos\theta$ lie in the kernel of the leading-order operator 
${\mathcal{H}}-\partial_\theta$, which creates a serious obstruction to the invertibility of the linearized operator. The \emph{speed modulation} introduced in Section~\ref{sec:modeone} is tailored precisely to remove this obstruction by exploiting the
translation invariance of the Euler equations in the vertical direction. This is the motivation for
introducing the variables $\mathtt{V}_1$ and $\mathtt{V}_2$, which act as vertical speed
perturbations of the point vortices $\mathtt{P}_1$ and $\mathtt{P}_2$.

More precisely, in Section~\ref{sec:modeone} we choose $\mathtt{V}_1$ and $\mathtt{V}_2$ as functions of
$(\varepsilon,f)$ so as to cancel the first \emph{sine} Fourier mode of the nonlinear equation. Accordingly, we work with the new  functional
$$
{\bf F}_1(\varepsilon,f)
:=\Pi_{1,{\bf s}}^c {\bf F}(\varepsilon, \mathtt{V}_1(f),\mathtt{V}_2(f),f),
$$
where  $\Pi_{1,{\bf s}}^c$ denotes  the projection outside  the first \emph{sine}  mode.
The linearized operator of ${\bf F}_1$ with respect to $f$ has a slightly different structure than that of ${\bf F}$ and, according to \mbox{Corollary \ref{cor:asymp-linF1},} it reads as follows,
 \begin{align*}
      \partial_f {\bf F}_1(\varepsilon,f)[h]
      &={\Pi_{1,{\bf s}}^c}\bigg(\varepsilon^3 |\ln\varepsilon|\omega\partial_\varphi h+ \tfrac{\varepsilon}{2} \big(\mathcal{H}[h]-\partial_\theta h\big)+\varepsilon^2 \partial_\theta\big( h V_1(\varepsilon,f)\big)+\varepsilon^2\mathcal{H}_{\mathtt{u},0}[h]\\ &\quad+\varepsilon^2{\partial_\theta}\mathcal{S}[h]{+\varepsilon^2 \mathcal{M}_1(f)[h]\cos(\theta)+\varepsilon^3|\ln\varepsilon| \mathcal{M}_2(f)[h]\cos(\theta)} +{\varepsilon^3}|\ln\varepsilon|  {\partial_\theta}\mathcal{Q}[h]\\ &\quad-\varepsilon^3\partial_\theta\fint_{\mathbb{T}} W({f})(\varphi,\theta,\eta)\ln\big|\sin\big(\tfrac{\theta-\eta}{2}\big)\big|h(\eta) d\eta +\varepsilon^3|\ln\varepsilon|^{\frac12}\partial_\theta \mathcal{R}_{1,\infty}[h]\bigg),
  \end{align*}
  where the additional new operators $\mathcal{M}_1$ and $\mathcal{M}_2$ arise from the speed modulation constraint.\\
  They act only on the time variable and are nonlocal operators; see Corollary~\ref{cor:asymp-linF1} for their precise definition. The operator $\mathcal{R}_{1,\infty} \in\mathrm{OPS}^{-\infty}$ is smoothing at any order.

\medskip\noindent{\bf $\diamondsuit$ {\sc Step 3:} Construction of a suitable  approximate solution.}
\medskip

The construction of an approximate solution is a central and delicate step in the overall proof. It serves as the starting point for the Nash–Moser iteration, and its accuracy directly determines whether the nonlinear scheme can converge. As we shall explain in Section \ref{sec:approx}, this stage involves several intertwined analytical challenges that stem both from the intrinsic structure of the equation and from the degeneracies specific to the problem.
A major difficulty is that the nonlinear functional exhibits a strong imbalance among its Fourier modes. In particular, the cosine mode one plays a distinguished role: it is the only tangential mode that interacts in a degenerate manner with the geometry of the moving vortex cores, whereas all higher modes correspond to normal oscillations of the boundary. This anisotropy is reflected in the linearized operator. While the higher modes behave regularly, the mode one; even after suppressing its sine component through speed modulation, still satisfies a degenerate first-order ODE whose coefficients depend on the filament dynamics. This equation cannot be solved explicitly, and constructing  a periodic solution requires a nontrivial solvability condition, which is obtained by exploiting the reversibility property of the underlying Hamiltonian system.\\
To be more precise, for any $N\in\N$ fixed, we construct the approximate solution through a finite  expansion in $\varepsilon$, separating at each order the tangential and
normal components:
$$
f_N(\varphi,\theta)
   = \sum_{m=0}^N \varepsilon^m
      \big( g_m(\varphi)\cos\theta + h_m(\varphi,\theta) \big).
$$
The normal components $h_m$ are chosen so as to cancel the projected nonlinear
terms of order $\varepsilon^{m+2}$, making use of the fact that the operator
$\partial_\theta-\mathcal{H}$, with $\mathcal{H}$ being the toroidal Hilbert transform,  is invertible on all modes except the mode one.  The tangential coefficients $g_m$ are then determined by solving an ODE encoding the
residual evolution of the cosine mode one, of the type
\begin{align*}
 \nonumber  g_m^\prime(\varphi)+\check{\beta}_1\big(g_m-g_{m,\star}\big)(\varphi)+&\check{\beta}_2(\varphi)\int_0^\varphi \check{\beta}_3(s)(g_m-g_{m,\star})(s)ds=\, F_m(\varepsilon,\varphi)\\
 +&|\ln\varepsilon|^{-\frac12}G_{1,m}(\varepsilon,\kappa, g_{m},g_{m,\star})(\varphi),
\end{align*}
where $\check{\beta}_j$ are smooth functions inherited from the point-vortex dynamics, and with
$$
\sup_{\varepsilon\in(0,\varepsilon_0)}\|G_{1,m}(\varepsilon,\kappa,g_m,g_{m,\star})\|_s^{{\textnormal{Lip},\nu}}<\infty,
$$
where the norm $||\cdot||_s^{\textnormal{Lip},\nu}$ is the standard $H^s$ norm for $(\varphi,\theta)$, and $(\textnormal{Lip},\nu)$ refers to the regularity with respect to the parameter $\lambda$, see Section \ref{SEc-function-sapces} for more details. The nonlocal term appears as a consequence of the speed modulation. In Section \ref{Sec-Fundam-lemma} we give a fundamental lemma about the invertibility of such integro-differential operators. Later, we choose $\varepsilon$ small such that a fixed point argument can be performed. The existence of a periodic solution
to this ODE is obtained through a
refined analysis using in a crucial way the reversibility of the dynamical system together with the symmetry property of the phase space.
Then, iterating this construction yields an approximate solution $f_N$ satisfying
$$
{\bf F}_1(\varepsilon,f_N) = O(\varepsilon^{N+3}|\ln\varepsilon|).
$$
We point out that this construction does not rely on any Cantor set.
\\
Finally, we rescale the functional ${\bf F}_1$ with the approximate solution and then work with
$$
\mathcal{F}(\varepsilon,g) := \frac{1}{\varepsilon^{1+\mu}}{\bf F}_1\big(\varepsilon, f_N + \varepsilon^\mu g\big),
\qquad 1 \leqslant \mu < 1+N,
$$
where
$$\mathcal{F}(\varepsilon,0) = O\big(\varepsilon^{\,2+N-\mu}|\ln\varepsilon|\big),
$$
and in view of Corollary \ref{prop:asymp-lin-2},
the  linearized operator $\mathscr{L}(\varepsilon,g):=\partial_g\mathcal{F}(\varepsilon,g)$ writes as
\begin{align}\label{Lin-Expand}
    \nonumber  \mathscr{L}(\varepsilon,g)[h]&={\Pi_{1,{\bf{s}}}^c}\bigg(\varepsilon^2 |\ln\varepsilon|\omega(\lambda)\partial_\varphi h+ \tfrac{1}{2} \big(\mathcal{H}[h]-\partial_\theta h\big)+\varepsilon \partial_\theta\big( h \mathcal{V}(\varepsilon,g)\big)+\varepsilon \mathcal{H}_{\mathtt{u},0}[h]+\varepsilon{\partial_\theta}\mathcal{S}[h]\\ 
    \nonumber&\quad +{\varepsilon \,\mathscr{M}_1[h]\cos(\theta)+\varepsilon^2\,|\ln\varepsilon| \mathscr{M}_2[h]\cos(\theta)} +{\varepsilon^2}\,|\ln\varepsilon|  {\partial_\theta}\mathcal{Q}[h]\\ &\quad-\varepsilon^2\partial_\theta\fint_{\mathbb{T}} \mathcal{W}({g})(\varphi,\theta,\eta)\ln\big|\sin\big(\tfrac{\theta-\eta}{2}\big)\big|h(\eta) d\eta+\varepsilon^2\partial_\theta\mathcal{R}_{\infty}(g)[h]\bigg)\,.
  \end{align}

\medskip\noindent{ $\diamondsuit$ \bf {\sc Step 4:} Invertibility of the linearized operator.}
\medskip

A central step in the Nash--Moser scheme is the construction of an approximate right inverse
for the linearized operator associated with the renormalized functional $\cF(\varepsilon,g)$.
This requires a careful analysis of the operator
$
\mathscr{L}(\varepsilon,g) 
$
whose structure is highly nontrivial due to the degeneracy of
the first Fourier mode, the presence of nonlocal pseudo-differential terms, and the various
singular scalings in the small parameter $\varepsilon$. Sections \ref{sec:invertibility} and \ref{sec:rednormalop} develop a systematic procedure to overcome these obstructions and
establish the required invertibility estimates. The main ideas are summarized below.

The operator $\mathscr{L}(\varepsilon,g)$ acts on functions of $(\varphi,\theta)$ and naturally separates
into the mode associated with $\cos\theta$ (the \emph{tangential direction}) and all higher
Fourier modes (the \emph{normal directions}). Writing $\Pi$ for the projector onto the
cosine mode one, and $\Pi^\perp$ for its complement which localized on space modes $|j|\geqslant2$, the operator can be represented as a
$2\times 2$ block matrix operator:
\begin{align*}
\mathbb{M} &\;:=\;\begin{pmatrix}
\Pi\mathscr{L}(\varepsilon,g)\Pi&\Pi\mathscr{L}(\varepsilon,g)\Pi^\perp \\
\Pi^\perp\mathscr{L}(\varepsilon,g)\Pi & \Pi^\perp\mathscr{L}(\varepsilon,g)\Pi^\perp
\end{pmatrix}=:\begin{pmatrix}
\mathcal{L}_{1}&\varepsilon \mathcal{L}_{1,\perp} \\
\varepsilon \mathcal{L}_{\perp,1} & \mathcal{L}_{\perp}
\end{pmatrix},
\end{align*}
where $\cL_1$ is the restriction to the tangential mode and $\cL_\perp$ acts on the higher modes.
The presence of off-diagonal couplings, although of size $O(\varepsilon),$ prevents direct inversion and requires a block reduction.

\medskip\noindent \textbf{$\blacktriangle$ {\sc Step 4.1:} Triangularization up to a smoothing error.}
\medskip

Section \ref{sec:triangularred} develops the first  delicate step toward inverting the linearized operator 
$\mathscr{L}(\varepsilon,g)$, namely a \emph{triangular reduction} that separates the dynamics of the 
degenerate Fourier mode~$1$ from the remaining normal modes. 
The starting point is the block decomposition of the linearized operator, represented by the matrix operator $\mathbb{M}.$
To remove the lower--left block up to higher--order errors, Section \ref{sec:triangularred} introduces a carefully designed 
near–identity linear transformation of the form
$$
\Psi =
\begin{pmatrix}
\mathrm{Id} & 0 \\
\varepsilon \psi & \mathrm{Id}
\end{pmatrix},
\qquad 
\Psi^{-1} =
\begin{pmatrix}
\mathrm{Id} & 0 \\
-\varepsilon \psi & \mathrm{Id}
\end{pmatrix},
$$
where the operator $\psi$ is chosen so that the conjugated operator becomes almost triangular:
$$
\Psi^{-1} \mathbb{M} \Psi
=
\begin{pmatrix}
\mathcal{L}_1 + \varepsilon^{2} \mathcal{L}_{1,\perp}\psi & \varepsilon \mathcal{L}_{1,\perp} \\
0 &
\mathcal{L}_\perp - \varepsilon^{2} \psi \mathcal{L}_{1,\perp}\psi
\end{pmatrix}
+ \varepsilon^{M+2}{\mathbb{P}_M},
$$
where the integer  $M\in\N$ is arbitrary and $\mathbb{P}_M\in \textnormal{OPS}^{-\infty}$. The remarkable feature is that the troublesome coupling term $\mathcal{L}_{\perp,1}$ has been 
eliminated, leaving a block-upper-triangular structure, up to a small and smoothing remainder .
This allows the invertibility problem to be reduced to the two diagonal blocks 
$\mathcal{L}_1$ and $\mathcal{L}_\perp$, which are handled separately in Sections \ref{section-Inver-Mode1} and \ref{sec:inver-normal}.

\medskip

\noindent \textbf{$\blacktriangle$ {\sc Step 4.2:} Invertibility of the tangential operator $\cL_1.$}
\medskip

Section \ref{section-Inver-Mode1} is devoted to the analysis of the operator $\mathcal{L}_1$, which governs the
dynamics of the cosine mode one. This mode is intrinsically delicate: it is the only
tangential direction that survives after the speed modulation, and it carries the residual
degeneracy of the nonlinear problem. After the reductions of Section \ref{section-Inver-Mode1}, the operator takes
the form
$$
\mathcal{L}_1
    \;=\;
    \varepsilon^2|\ln\varepsilon|\,\omega(\lambda)\,\partial_\varphi
    \;+\;
    \mathcal{K}_1(\varepsilon,g),
$$
where $\mathcal{K}_1$ is a zero–order operator with Lipschitz dependence  on the
state~$g$.
At first sight, the small prefactor $\varepsilon^2|\ln\varepsilon|$ appears to obstruct the
invertibility of $\mathcal{L}_1$. The crucial observation, however, is that the perturbation
satisfies the asymptotic relation
$$
\mathcal{K}_1(\varepsilon,g)
    \;\sim\;
    \varepsilon^2|\ln\varepsilon|\,\mathcal{K}_0(g),
$$
so that $\mathcal{K}_1$ vanishes at precisely the same rate as the transport part. Thus
$\mathcal{L}_1$ is not a small operator plus a large perturbation, but rather a
\emph{balanced} sum of two contributions of comparable size. In this sense,
$\mathcal{L}_1$ is a compact perturbation of a first–order transport operator, and this
structure is at the heart of its invertibility.
A more refined understanding relies on the specific form of the approximate solution
$f_N$. When restricted to the cosine mode one, the equation $\mathcal{L}_1[h_1(\varphi)\cos\theta] = k_1(\varphi)\cos(\theta)$ reduces to
a delayed first–order ODE of the form
\begin{align*}
\varepsilon^2 |\ln\varepsilon|\Big(\omega(\lambda) h_1'
        +\check{\beta}_1\, (h_1-h_{1,\star})&+\check{\beta_2}\int_0^\varphi \check{\beta}_3(s)(h_1-h_{1,\star})(s)ds
   \\
   &+ |\ln\varepsilon|^{-\frac12}\big(a_1[g](\varphi) h_1
        + a_2[g](\varphi) h_{1,\star}\big)\Big)
    = k_1,
\end{align*}
where $\check{\beta}_j$  encode lower–order contributions depending on the state~$g$. Indeed, the nonlocal term appears as a consequence of the speed modulation. At the main order, the integrodifferential equation equals to the one studied for the approximate solution. Hence, we refer to Section \ref{Sec-Fundam-lemma} for the general study of the invertibility of such type of operators. Later, solving this equation in the periodic setting relies on a fixed–point argument, exploiting
both the smoothing character of $\mathcal{K}_1$ and the reversibility symmetry of the
underlying Hamiltonian system.  We choose $\varepsilon$ small appropriate to perform the fixed point argument. 
This strategy yields the desired tame estimate
$$
\big\|\mathcal{L}_1^{-1}[h]\big\|_s
    \;\lesssim\;
    \frac{1}{\varepsilon^2|\ln\varepsilon|}
    \Big( \|h\|_{s-1}
        + \|g\|_s \|h\|_{s_0} \Big),
    \qquad 3<s_0\leqslant s,
$$
which provides the quantitative control needed for the Nash--Moser iteration.

\medskip\noindent\textbf{$\blacktriangle$ {\sc Step 4.3:} Invertibility of the normal operator $\mathcal{L}_\perp.$}
\medskip

The construction of an approximate right inverse for the normal operator
$\mathcal{L}_\perp$ is highly delicate due to the intricate structure of the
linearized operator. Indeed, $\mathcal{L}_\perp$ takes the form
\begin{align*}
   \mathcal{L}_\perp (\varepsilon,g) 
   = \Pi^\perp \Big[ \mathcal{T} + \tfrac12 \mathcal{H}
     + \varepsilon \mathcal{H}_{\mathtt{u},0}
     + \varepsilon \partial_\theta \mathcal{S}
     + \varepsilon^2 \partial_\theta \mathcal{R}(g)\Big]\Pi^\perp,
\end{align*}
where
\begin{align*}
   \mathcal{T}(g)[h] 
   := \varepsilon^2 |\ln \varepsilon|\, \omega(\lambda)\,\partial_\varphi h 
   + \partial_\theta\!\left(\Big(-\tfrac12 + \varepsilon \mathcal{V}(\varepsilon,g)\Big)h\right).
\end{align*}
The operators above are defined in Section~\ref{sec-Asym-Line-op}:
$\mathcal{T}$ stands for the transport part, $\mathcal{H}$ is the Hilbert
transform, $\mathcal{H}_{\mathtt{u},0}$ is a zero–order operator (a type of heterogeneous Hilbert transform), $\mathcal{S}$
is a shifting operator. As to the operator $\mathcal{R}$, it  takes the form
\begin{align*}
 \mathcal{R}(g)[h](\varphi,\theta) 
   = \fint_{\T} \mathcal{W}(g)(\varphi,\theta,\eta)\,
      \ln\!\Big|\sin\!\Big(\tfrac{\theta-\eta}{2}\Big)\Big|\,
      h(\varphi,\eta)\, d\eta 
      + \mathcal{R}_{\infty}(g)[h](\varphi,\theta).
\end{align*}
with $\mathcal{W}(g)$ being a function and $\mathcal{R}_{\infty}(g)\in\textnormal{OPS}^{-\infty}.$

Several obstructions to the invertibility of $\mathcal{L}_\perp$ still remain:
a degenerately small transport term in the $\varphi$–direction, a non-diagonal
zero-order contribution, nonlocal pseudo-differential components, and a shifting
structure in the $\theta$–variable produced by the interaction of the two vortex
rings. To overcome these difficulties, we develop a sequence of
conjugations, each designed to remove a specific obstruction and gradually
transform $\mathcal{L}_\perp$ into an operator that is diagonal up to a small
smoothing remainder.
From Proposition~\ref{prop-diagonal-impo}, we obtain the following result.
There exists a reversibility-preserving isomorphism; for the phase space definition see \eqref{compact-notat1},
$$
\Phi:\mathbb{X}_\perp^s \to \mathbb{X}_\perp^s,
$$
satisfying the tame estimates
$$
    \|\Phi^{\pm1}[h]\|_{s}^{\mathrm{Lip},\nu}
    \lesssim \|h\|_{s}^{\mathrm{Lip},\nu}
       + \|g\|_{s+3}^{\mathrm{Lip},\nu}\,\|h\|_{s_0}^{\mathrm{Lip},\nu},
$$
such that for each $n\in\mathbb{N}$ and for all $\lambda$ in the Borel set
$$
{\mathcal{O}_{n}^{1}(g)}
 = \bigcap_{\substack{(\ell,j)\in\mathbb{Z}^{2} \\ 1\leqslant |j|\leqslant N_n}}
   \left\{
   \lambda\in \mathcal{O}:\;
   \big|\varepsilon^2|\ln\varepsilon|\omega(\lambda)\ell
   + j\,\mathtt{c}_1(\lambda,{\varepsilon},g)\big|
   \geqslant \nu |j|^{-\tau}
   \right\},
$$
we have
\begin{align}\label{DIAGONA-LL}
\mathcal{L}_{\perp,1}:=\Phi^{-1}\mathcal{L}_\perp\Phi
&= \Pi^\perp\Big(
     \varepsilon^2|\ln\varepsilon|\,\omega(\lambda)\,\partial_\varphi
     +\mathcal{D}_\varepsilon
   \Big)\Pi^\perp 
 + \varepsilon^2|\ln\varepsilon|^{1/2}\Pi^\perp\partial_\theta\mathcal{R}_{4,-3}\Pi^\perp
 + \mathtt{E}_n^4,
\end{align}
where $\mathcal{D}_\varepsilon$ is a diagonal Fourier multiplier, $\partial_\theta\mathcal{R}_{4,-3}\in\mathrm{OPS}^{-3}$ and $\mathtt{E}_{n}^4$ is a smooth residual operator. The constant $\mathtt{c}_1(\lambda,{\varepsilon},g)$ arises from the successive averaging steps performed in the KAM scheme used to straighten the transport part. 
\\
The construction of $\Phi$ is a key point to invert the normal operator. It will be carried out in Section \ref{sec:rednormalop} and is achieved through four successive conjugations, producing maps
$\Phi_1,\Phi_2,\Phi_3,\Phi_4$ such that
$$
\Phi \;=\; \Phi_1 \circ \Phi_2 \circ \Phi_3 \circ \Phi_4.
$$
Each transformation removes one source of non-diagonal structure or
nonlocality, ultimately yielding a normal operator that is diagonal modulo a
small smoothing remainder.
We now briefly describe the role of these transformations; a more detailed sketch 
of the reduction is given at the beginning of Section~\ref{sec:rednormalop}. 
The first transformation $\Phi_1$ reduces the transport part of the operator to a 
Fourier multiplier up to small zero-order remainders, by means of KAM-type techniques.
The second transformation $\Phi_2$ is introduced to eliminate the non–diagonal 
zero–order terms of the linearized operator. As a result, the zero–order structure 
is significantly simplified, although the remaining nonlocal operators are still 
not diagonal and retain a shifting component. 
The third transformation $\Phi_3$ is then designed to remove this antidiagonal 
contribution of the operator, up to a small remainder.
After performing the preceding transformations, the Hilbert transform appears multiplied by
a coefficient depending only on $\varphi$. In order to obtain an operator that is
diagonal in the variables $(\varphi,\theta)$, up to a small remainder, we introduce
a final transformation $\Phi_4$ which removes the $\varphi$–dependence of this
coefficient. 
The invertibility of $\mathcal{L}_{\perp,1}$ is obtained by first inverting its
principal part and then applying a perturbative argument. Since the leading
operator is a Fourier multiplier, its invertibility reduces to the verification
of suitable Diophantine conditions. In particular, one must avoid the
occurrence of small divisors. This leads to restricting the parameter
$\lambda$ to the following Cantor-type set

\begin{equation*}
\mathcal{O}_{n}^2(g)
=\bigcap_{\substack{(\ell,j)\in\mathbb{Z}^{2}\\ 2\leqslant |j|\leqslant N_{n}}}
\Big\{\lambda\in \mathcal{O}, \,\,\left|\varepsilon^2|\ln\varepsilon| \omega(\lambda)\ell+\mu_{j,2}(\lambda,\varepsilon,g)\right|\geqslant \nu | j|^{-\tau}\Big\},
\end{equation*}
with $\mu_{j,2}$ is the spectrum of diagonal operator $\mathcal{D}_\varepsilon$, given by 
\begin{equation*}
      \mu_{j,2}(\lambda,{\varepsilon}, g)=\mathtt{c}_1(\lambda,{\varepsilon},g)j+\mathtt{c}_2(\lambda,{\varepsilon},g) \, \textnormal{sign} j+\varepsilon^2\mathtt{c}_3(\lambda,{\varepsilon},g)\,\tfrac{\textnormal{sign} j}{j^2-1},
\end{equation*}
where $\mathtt{c}_k$ are constant in time and space,  depending on $g$ and $\varepsilon$. Moreover they are close to universal constants for $\varepsilon$ small enough.

\medskip\noindent{\bf $\diamondsuit$ {\sc Step 5: }  Nash-Moser scheme  and size of Cantor set.}
\medskip

In Section~\ref{sec:nashmoser} we combine all the previous reductions in order to implement a Nash--Moser scheme, following a classical strategy previously used in related settings; see for instance \cite{BCP,BHM,HHM24}. 
This allows us to construct a sequence of approximate solutions $(g_n)$. In this step we make the following choice of  $\nu$, see \eqref{lambda-choice}, 
$$
 \nu:= \varepsilon^2|\ln\varepsilon|^{\delta}.$$
We prove that this sequence converges to an exact solution provided that the parameter $\lambda$ belongs to a final Cantor-like set $\mathtt{C}_\infty$ given by
\begin{equation*}
\mathtt{C}_{\infty}:=\bigcap_{m\in\mathbb{N}}\mathcal{A}_{m}.
\end{equation*}
Here, given $0<\lambda_\star<\lambda^\star$, the sets $\mathcal{A}_n$ are defined recursively by
$$
\mathcal{A}_{0}=(\lambda_\star,\lambda^\star), 
\qquad 
\mathcal{A}_{n+1}=\mathcal{A}_{n}\cap\mathcal{O}_{n}^1(g_n)\cap\mathcal{O}_{n}^2(g_n),
\quad \forall n\in\mathbb{N}.
$$

Finally, in Section~\ref{sec:cantorset} we estimate the measure of the set $\mathtt{C}_\infty$ and prove that
$$
|(\lambda_\star,\lambda^\star)\setminus \mathtt{C}_{\infty}|
\leqslant C|\ln\varepsilon|^{\delta-1},
$$
for some $\delta\in(\tfrac12,1)$. In particular, the set $\mathtt{C}_\infty$ has asymptotically full Lebesgue measure.

\medskip
Finally, we give five appendices that complement this work. In Appendix \ref{apendix:transport} we give the idea of the invertibility of the transport equation, as a toy model of our linear operator, when we have a degeneracy in the time derivative using some appropriate Cantor set. The classical conjugation using KAM techniques of a transport equation with variable coefficients is given in Appendix \ref{appendix:KAM}, which we will use later for our linear operator. Appendix \ref{sec:appendix} refers to the asymptotic analysis of the Green function for the 3D axisymmetric equations. In Appendix \ref{appendix:Fourier} we give different expansions of the linear operator. Some useful integrals needed in this paper are computed in Appendix \ref{appendix:integrals}.

\subsection{Differences and challenges compared to previous works}\label{sec:intro4}
Finally, after outlining the general structure of the proof, we turn to the
discussion of some of the main difficulties that arise in the present work. These
difficulties are specific to the three-dimensional setting and mark a
significant departure from the mechanisms encountered in previous related
results.

\medskip

\noindent $\blacktriangle$ {\bf Small divisor problem and KAM techniques.}

\medskip

Previous works on leapfrogging motion for the $3$d Euler equations,
notably those of Dávila--Del Pino--Musso--Wei~\cite{delpino},
Buttà--Cavallaro--Marchioro~\cite{Butta-Cavallaro-Marchioro},
and Donati--Hientzsch--Lacave--Miot~\cite{Miot-Lacave}
establish the existence of this phenomenon only locally in time.
In particular, the solutions constructed in these works are not periodic in time,
and their validity is restricted to regimes where the interaction of the vortex
filaments can be controlled over finite times.

In contrast, the purpose of the present work is to construct time--periodic
solutions of the three--dimensional Euler equations that exhibit the
leapfrogging motion. Such solutions are necessarily global in time and
therefore require a much finer control of the dynamics. Indeed, when one
attempts to extend the analysis beyond short time scales, new difficulties
arise from the presence of resonances in the linearized dynamics. These
resonances may destroy the periodic structure of the motion and may lead
to complex filamentation phenomena of the vortex tubes.

This issue marks a fundamental transition in the analysis compared with the
aforementioned works. While short--time constructions rely essentially on
perturbative arguments around the filament dynamics, the global--in--time
periodic problem requires controlling small divisors generated by the
resonant structure of the linearized operator. To overcome this difficulty,
we employ KAM--type techniques combined with a Nash--Moser iteration scheme.

\medskip
\noindent $\blacktriangle$ {\bf 2D leapfrogging for point vortices.}
\medskip

In \cite{HHM24}, Hassainia–Hmidi–Masmoudi established the existence of a non-rigid periodic configuration for point vortices in the two-dimensional Euler equations, providing a simplified dynamical model for leapfrogging motion. The general strategy developed in the present work is inspired by that approach. However, the three-dimensional setting and the vortex-ring structure introduce substantial new difficulties and conceptual challenges, which require significantly different analytical tools and refinements of the method. We briefly outline below some of the main challenges.
\subsubsection*{$1)$ Singular self-induced speed.}
In contrast with the two-dimensional case, a single vortex ring in the three-dimensional Euler equations travels with a singular self-induced velocity of order $\ln\varepsilon$. As a consequence, the nonlinear functional cannot be expanded as a regular power series in $\varepsilon$. Instead, one must identify and work within an intermediate asymptotic regime determined by this logarithmic scaling, which substantially complicates the perturbative analysis.

\subsubsection*{$2)$ {\it  Hamiltonian structure and Elliptical shape.}}
The transport equation for the potential vorticity in the axisymmetric three-dimensional Euler equations does not naturally possess a Hamiltonian structure in the original cylindrical variables $(r,z)$. To recover a Hamiltonian formulation, we perform the change of variables $(r,z) \mapsto (\sqrt{2r},z)$. 
However, this transformation has an important geometric consequence. 
An initially small circular cross--section in the $(r,z)$ variables 
is mapped, at leading order, to an ellipse in the new coordinates. 
As a result, the natural geometry of the vortex tube is no longer 
isotropic, and the deformation introduced by the change of variables 
propagates throughout the analysis. 
This geometric distortion is especially reflected in the asymptotic expansions 
carried out in Appendices~\ref{sec:appendix} and~\ref{appendix:Fourier}, 
where the resulting anisotropic structure leads to more involved expansion 
formulas and additional coupling terms.
In particular, the structure of both linearized and  nonlinear operators 
is significantly affected, since their coefficients inherit this anisotropic 
deformation. As a result, the analysis of these operators becomes considerably 
more delicate and requires careful treatment.

  \subsubsection*{$3)$ {\it Leapfrogging filaments.}}
In the two-dimensional Euler setting, the intuition behind leapfrogging motion originates from the existence of non-rigid periodic dynamics of four point vortices, whose evolution is explicitly known and well understood. In contrast, in the three-dimensional case the relevant heuristic picture comes from the leapfrogging motion of two vortex filaments. However, unlike the 2D point-vortex system, the dynamics of interacting vortex filaments is considerably more intricate and far less explicit.

Although classical work asserts the existence of configurations of two vortex filaments exhibiting the leapfrogging motion, our analysis requires a much more precise description of the underlying dynamics, as developed in   Section~\ref{sec:filamentdyn}. In particular, several quantitative properties of the filament system play a crucial role in the subsequent construction. 
One of the key ingredients is the behavior of the period of the leapfrogging motion as a function of the energy. More precisely, we require a monotonicity property of the period function with respect to the energy level of the two–filament system. The explicit expression of this period involves an integral representation. Establishing this monotonic dependence is crucial for our analysis, as it provides the required non–degeneracy of the frequency map and plays a central role in the measure estimates for the final Cantor set arising in the KAM scheme.

Another important aspect concerns the geometry of the leapfrogging trajectories. The structure of the orbits in the reduced point–vortex system exhibits specific symmetry properties which can be exploited in the reformulation of the vortex–patch dynamics. In particular, these geometric features allow for a symmetry reduction that couples the dynamics of the two cross–sections of the vortex rings. As a consequence, the evolution of the two sections can be reduced to a single functional equation, in which the interaction between the rings appears through a delayed term reflecting the phase shift along the leapfrogging orbit.

This refined understanding of the point–filament dynamics therefore constitutes a fundamental preliminary step in our analysis, as it provides the structural information needed to implement the desingularization procedure and to formulate the contour–dynamics equation in a form amenable to the KAM–Nash–Moser scheme.

\subsubsection*{$4)$ {\it Speed modulation and elimination of the first sine mode.}}
In both the 2D and 3D settings, the first Fourier mode constitutes a major difficulty due to a degeneracy in the linearized operator. In the 2D Euler case, this issue was resolved through the analysis of an associated monodromy matrix. In the present framework, however, the increased complexity of the operator makes this approach considerably more delicate. Instead, we exploit the translation invariance of the system to introduce a modulation of the propagation speed. This additional degree of freedom is then chosen so as to completely eliminate the first sine Fourier mode from the nonlinear functional, thereby reducing the problem to solving an ordinary differential equation for the first cosine Fourier mode. Note that the modulation procedure comes at a certain cost. Indeed, it generates additional
nonlocal contributions in the linearized operator, represented by the
operators $\mathscr{M}_j$ in \eqref{Lin-Expand}. These terms appear at
intermediate orders in $\varepsilon$ and therefore contribute in a
non-negligible way to the structure of the linearized operator. As a
result, they cannot be treated as harmless remainders and must be carefully
tracked throughout the analysis.

A key structural property, however, allows us to control these terms.
The most delicate contribution, $\mathscr{M}_1$, possesses an
anti-diagonal structure, while the operator $\mathscr{M}_2$ is activated
only on the first Fourier mode. These features play a crucial role in the
reducibility scheme, as they allow us to isolate the problematic
components and treat them through suitable transformations and
mode-by-mode analysis.

\subsubsection*{$5)$ Degeneracy in $\varepsilon$ of the linearized operator}

A further major difficulty concerns the degeneracy of the linearized operator
with respect to the small parameter $\varepsilon$. This phenomenon already
appears in the desingularization of Love's configuration performed in
\cite{HHM21}, but the structure of the operator in that context is
substantially simpler than in the present three--dimensional setting.
More precisely, in the two--dimensional framework, the linearized operator takes the form
\begin{align*}
\partial_g \mathcal{F}(\varepsilon,g)[h]
&=
\varepsilon^2\omega(\lambda)\,\partial_\varphi h
+\partial_{\theta}\big[\mathcal{V}^\varepsilon(g)h \big]
-\tfrac{1}{2}\mathcal{H}[h]
-\varepsilon^2\partial_\theta\mathcal{Q}_0[h](\varphi,\theta)
+\varepsilon^3 \partial_\theta \mathcal{R}(g)[h],
\end{align*}
where $\mathcal{Q}_0$ is a finite--rank operator localizing on the first
Fourier mode and $\mathcal{R}$ is a smoothing operator of arbitrary order.
In particular, apart from the transport term and the Hilbert transform,
all the additional contributions appear at order $\varepsilon^2$ or higher,
and the remainder is strongly regularizing.
This asymptotic structure is in sharp contrast with the situation arising
in the present three--dimensional problem. Here the linearized operator,
described in \eqref{Lin-Expand}, exhibits a significantly more intricate
behavior. In particular, intermediate regimes of size $\varepsilon$ appear,
involving nonlocal and non--diagonal operators that cannot be treated as
negligible perturbations. Moreover, the zero--order remainder is no longer
smoothing at arbitrary order, which prevents the use of the simpler
reducibility arguments employed in the two--dimensional case.
As a consequence, the diagonalization and reducibility scheme developed in
this work must be considerably refined. 
The presence of these intermediate nonlocal effects requires a more delicate sequence of transformations in order to isolate the principal part of the operator and to construct a suitable approximate right inverse, which constitutes a crucial ingredient in the Nash--Moser scheme.

 \subsubsection*{$6)$ {\it Construction of a refined approximate solution.}}

Another important difficulty concerns the construction of a suitable
approximate solution prior to the implementation of the Nash--Moser
iteration. In the two--dimensional setting considered previously,
the approximate right inverse of the linearized operator exhibits a
loss of regularity together with an inflation of the operator norm of
order $\varepsilon^{-2-\delta}$. As a consequence, the trivial state
$g=0$ cannot serve as a sufficiently accurate initial approximation.
Indeed, starting the Nash--Moser scheme from such a rough state would
generate corrections whose norms grow too rapidly, making the iteration
scheme inapplicable.
To overcome this issue, it is necessary to construct a sufficiently
accurate approximate solution. In the two--dimensional case, this was achieved by building
an approximate solution with an error of order $O(\varepsilon^{3-\mu})$,
which provides enough precision to compensate for the loss of derivatives
and the inflation of the inverse operator.

In the present three--dimensional setting, the situation becomes even
more delicate. The inflation of the approximate right inverse is more
pronounced and behaves as $\varepsilon^{-3}$. This stronger degeneracy
makes the construction of a high--quality approximate solution even
more crucial in order to initiate the Nash--Moser iteration. 
Fortunately, as explained earlier in Step~$3$, we are able to construct an approximate
solution to arbitrary order in $\varepsilon$ through a flexible and robust procedure.

\subsubsection*{$7)$ {\it Cantor type sets and logarithmic regime.}}
As described in \eqref{DIAGONA-LL}, the transformed operator can be written as the sum of a diagonal Fourier multiplier and a sufficiently smoothing remainder of size $\varepsilon^2|\ln\varepsilon|^{\frac12}$. 
The invertibility of this operator is obtained through a perturbative argument, by first inverting its leading diagonal part and then treating the smoothing contribution as a small perturbation.
This strategy is effective provided that the inflation in $\varepsilon$ of the right inverse of the diagonal part can be compensated by the smallness of the remainder. In order to guarantee this balance, we introduce a penalization in the small–divisor conditions defining the Cantor set by choosing the threshold parameter in the form
$$
\nu=\varepsilon^2|\ln\varepsilon|^\delta ,
$$
for some exponent $\delta>\tfrac12$. With this choice, the perturbative argument becomes applicable and the inversion of the full operator can be carried out.
On the other hand, in order to ensure that the resulting Cantor set still has asymptotically full measure, the exponent must also \mbox{satisfy $\delta<1$. } 
This choice of the threshold parameter $\nu$ should be contrasted with the situation arising in the two--dimensional setting. In that lower--dimensional framework, no intermediate logarithmic regime appears and the expansion of the linearized operator is analytic with respect to the small parameter $\varepsilon$. As a consequence, the small--divisor threshold can be selected in the simpler form
$$
\nu=\varepsilon^{2+\delta}.
$$
The presence of the logarithmic scale in the present three--dimensional problem therefore introduces an additional level of complexity, requiring a more delicate balance between the size of the smoothing remainder and the inflation of the right inverse of the diagonal part.

\medskip\noindent $\diamondsuit$ {\bf Quasiperiodic leapfrogging vortex rings.} 

\medskip

Theorem~\ref{th-main1} is stated in the time--periodic setting, which is
already sufficiently intricate from both the analytical and technical
points of view. A natural question is whether the result can be extended
to the quasi--periodic regime by coupling the self--induced frequencies
of the vortex rings with those of the underlying filament system.
Such an extension appears conceivable by adapting the strategy developed
in works such as
\cite{BCP,BHM,BertiBiascoProcesi2013,BertiMontaltoHaus2017,GIP23,HHM21,HHR24}.

However, this generalization would introduce an additional level of
complexity due to the presence of time resonances. In particular, the
small--divisor conditions would no longer involve only the spatial modes
but also the temporal frequencies. As a consequence, the construction of
the Cantor set would require penalization with respect to the time
frequencies rather than the spatial ones, which is the mechanism used in
the present work.

Another significant difference concerns the perturbative simplification
of the linearized operator. In the periodic case studied here, we are able
to partially reduce the operator and treat the remaining terms as a small
and sufficiently smoothing perturbation in the spatial variable. In the
quasi--periodic setting, such a simplification would no longer be
available. Instead, one would need to perform a complete reducibility of
the linearized operator. Nevertheless, it is reasonable to expect that
such a reduction could be achieved by extending the techniques developed
in the present work, combined with the quasi--periodic KAM framework
implemented in the aforementioned references.

\section{Axisymmetric 3D Euler equations}\label{sec:axisym-euler}

A central theme in fluid dynamics is the existence and interaction of coherent structures such as vortex rings. Among these, the \emph{leapfrogging} motion of vortex rings is a paradigmatic phenomenon whose rigorous analysis relies on the symmetries and  geometric structure of the equations. In this work we focus on the axisymmetric setting, which naturally captures vortex-ring dynamics and admits a Hamiltonian reformulation reminiscent of two-dimensional active scalar models.

\subsection{Axisymmetric formulation}

We begin by recalling the structure of axisymmetric vector fields. A vector field $u:\R^3\to \R^3$ is called \emph{axisymmetric} (around the $x_3$-axis) if it is invariant under rotations about the $x_3$-axis. In cylindrical coordinates $(r,\theta,z)$, with
$$
r=\sqrt{x_1^2+x_2^2},\qquad x_1=r\cos\theta,\qquad x_2=r\sin\theta,\qquad x_3=z,
$$
an axisymmetric field can be written as
$$
u(x) \;=\; u^r(r,z)\, e_r(\theta) \;+\; u^\theta(r,z)\, e_\theta(\theta) \;+\; u^z(r,z)\, e_z,
$$
where
$$
e_r(\theta)=(\cos\theta,\sin\theta,0), \qquad e_\theta(\theta)=(-\sin\theta,\cos\theta,0), \qquad e_z=(0,0,1).
$$

When $u^\theta\equiv 0$, the flow is said to be \emph{axisymmetric without swirl}. This reduction is particularly relevant for vortex rings, since the azimuthal component is absent. Under this assumption, the incompressibility condition $\nabla\cdot u=0$ becomes
$$
\partial_r (r u^r) + \partial_z (r u^z)=0.
$$
Consequently, one can introduce a scalar \emph{stream function} $\psi=\psi(r,z)$ such that
$$
u=\nabla\times \phi, \qquad \phi=\frac{\psi}{r}\, e_\theta,
$$
Equivalently,
$$
(u^r,u^z)=\tfrac{1}{r}\,\nabla^\perp_{(r,z)}\psi,
$$
where $\nabla^\perp_{(r,z)}=( -\partial_z,\;\partial_r)$.  

The vorticity field $\omega=\nabla\times u$ is purely azimuthal:
$$
\omega(x)=(\partial_z u^r(r,z) - \partial_r u^z(r,z))\, e_\theta(\theta)=\omega^\theta(r,z) e_\theta(\theta).
$$
It is convenient to introduce the \emph{relative vorticity}
$$
q(r,z):=\frac{\omega^\theta(r,z)}{r},
$$
A direct computation shows that $\psi$ solves the elliptic relation
\begin{align}\label{stream-operat}
-\frac{1}{r^2}\mathcal{L}\psi=q, \qquad \mathcal{L}:=r\partial_r\left(\frac{1}{r}\partial_r\right)+\partial_z^2.
\end{align}
Moreover, $q$ obeys a transport-type active scalar equation: combining the above identities yields
\begin{equation}\label{eq:q}
\partial_t q + u\cdot\nabla q = 0,
\qquad
u=\frac{1}{r}\,\nabla^\perp_{(r,z)}\psi,
\end{equation}
where $u\cdot\nabla q=u^r\partial_r q+u^z\partial_z q$.
Thus $q$ is advected by the velocity field it generates, in close analogy with the
two-dimensional Euler dynamics (with the important geometric feature that $u$ is divergence-free in
$\R^3$, or equivalently divergence-free in $(r,z)$ with respect to the weighted measure $r\,dr\,dz$).

\medskip
\noindent\textbf{Axisymmetric Biot--Savart law.}
To characterize the Biot--Savart law in the axisymmetric setting, we start from
$$
-\Delta \phi=\omega^\theta e_\theta, \qquad 
\phi=\frac{1}{4\pi|x|}\star_{\R^3}(\omega^\theta e_\theta).
$$
Exploiting the rotational invariance of the axisymmetric Biot--Savart law (see, for instance, \cite{BertozziMajda, GallaySverak2014, Gallay-Sverak}), one can derive an explicit kernel representation of the stream function:
\begin{equation*}
\psi(r,z) \;=\; \frac{1}{2\pi}\bigintsss_{\Pi_+} \sqrt{rr'} \; 
J\!\left(\frac{(r-r')^2+(z-z')^2}{rr'}\right)\, q(r',z')\, r'\, dr'\,dz',
\end{equation*}
where $\Pi_+:=\{(r,z)\in\R^2:\ r>0\}$ and the auxiliary function $J:(0,\infty)\to\R$ is defined by
\begin{equation*}
J(s):=\bigintsss_0^{\pi}\frac{\cos\theta}{\sqrt{s+2-2\cos\theta}}\, d\theta, 
\qquad s>0.
\end{equation*}
This function is classically related to the Legendre function of the second kind with index $\tfrac12$. More precisely,
\begin{equation*}
J(s)\;=\; Q_{\frac12}\!\left(\tfrac{s}{2}+1\right),\qquad Q_{\frac12}(x) \;=\; \int_0^{\pi}\frac{\cos\theta}{\sqrt{2x-2\cos\theta}}\, d\theta,
\quad x>1.
\end{equation*}
The function $Q_{\tfrac12}$ solves the Legendre differential equation
\begin{equation*}
(1-x^2)Q_{\frac12}^{\prime\prime}(x) - 2x\, Q_{\frac12}^{\prime}(x) + \tfrac34 Q_{\frac12}(x) = 0,
\end{equation*}
and by straightforward application of the chain rule one deduces that $J$ satisfies 
\begin{equation}\label{eq:J-ODE0}
s(s+4)\,J^{\prime\prime}(s) + 2(s+2)\,J^{\prime}(s) - \tfrac{3}{4}\,J(s) = 0.
\end{equation}
For completeness, we also recall that $Q_{\frac12}$ admits a classical hypergeometric representation:
\begin{equation*}
Q_{\frac12}(x) \;=\; \tfrac{\pi}{4\sqrt{2}}\, x^{-\tfrac32}\, 
{}_2F_1\!\left(\tfrac54,\tfrac34;2;\tfrac{1}{x^2}\right),
\end{equation*}
where ${}_2F_1$ denotes the standard Gauss hypergeometric function. In addition,
$$
 Q_{\frac12}(x)={\sqrt{\tfrac\pi2}}\int_0^\infty e^{-tx}\tfrac{I_1(t)}{\sqrt{t}} dt,\qquad x>1,
$$
where $I_1$ is the modified Bessel function of the first kind.

\subsection{Hamiltonian reformulation}\label{sec:hamiltonian}
For analytical purposes, it is convenient to introduce the change of variables 
$$
(r,z)\longmapsto (\varrho,z),\qquad \varrho:=\tfrac{r^{2}}{2},
$$
and to set
$$
q(t,r,z)=\mathbf{q}(t,\varrho,z),\qquad \psi(t,r,z)=\Psi(t,\varrho,z).
$$
Note that $d\varrho=r\,dr$, so that the weighted measure $r\,dr\,dz$ on $\Pi_+=\{r>0\}$ becomes the Lebesgue measure $d\varrho\,dz$ on $(0,\infty)\times\R$.
By the chain rule,
$$
\partial_t q=\partial_t\mathbf{q},\qquad
\partial_r q=r\,\partial_\varrho\mathbf{q},\qquad
\partial_z q=\partial_z\mathbf{q}.
$$
Recalling that, in the swirl-free case,
$$
u^r=-\frac{1}{r}\partial_z\psi,\qquad u^z=\frac{1}{r}\partial_r\psi,
$$
we obtain in $(\varrho,z)$-variables
$$
u^r=-\frac{1}{r}\partial_z\Psi,\qquad u^z=\partial_\varrho\Psi.
$$
Therefore,
$$
u\cdot\nabla q
=u^r\partial_r q+u^z\partial_z q
=-\partial_z\Psi\,\partial_\varrho\mathbf{q}+\partial_\varrho\Psi\,\partial_z\mathbf{q}
=\nabla^\perp\Psi\cdot\nabla\mathbf{q},
$$
where
$$
\nabla:=(\partial_\varrho,\partial_z),\qquad \nabla^\perp:=(-\partial_z,\partial_\varrho).
$$
Thus \eqref{eq:q} becomes the canonical Hamiltonian transport equation
\begin{equation}\label{eq:tilde-q}
\partial_t \mathbf{q}+\nabla^\perp \Psi\cdot \nabla \mathbf{q}=0,
\qquad\text{equivalently}\qquad
\partial_t\mathbf{q}+\{\Psi,\mathbf{q}\}=0,
\end{equation}
with Poisson bracket $\{a,b\}:=\partial_\varrho a\,\partial_z b-\partial_z a\,\partial_\varrho b$.
In particular, the advecting field $\nabla^\perp\Psi$ is divergence-free in the usual sense.

\medskip
\noindent\textbf{Stream function and Green kernel.}
Transforming the axisymmetric Biot--Savart representation of $\psi$ yields
\begin{equation}\label{eq:stream-Psi}
\Psi(\varrho,z)=\frac{1}{\pi\sqrt{2}}\int_0^\infty\int_{\R} 
G(\varrho,z;\varrho',z')\, \mathbf{q}(\varrho',z')\, d\varrho' dz',
\end{equation}
where the symmetric kernel is given 
\begin{align}\label{eq:G}
G(\varrho,z;\varrho',z') \;:=\; (\varrho\varrho')^{1/4}\,
J\!\left(\tfrac{2(\sqrt{\varrho}-\sqrt{\varrho'})^2+(z-z')^2}{2\sqrt{\varrho \varrho'}}\right).
\end{align}
We will often use the shorthand (abusively identifying $(\varrho,z)$ with the complex number $\varrho+i z$),
$$
G(\varrho,z;\varrho',z')=G(\varrho+i z;\varrho'+ iz').
$$
This formulation shows that the axisymmetric Euler equations without swirl can be viewed as an \emph{active scalar} model on the half-plane $(0,\infty)\times\R$ with a Hamiltonian structure. This perspective is crucial for the study of coherent structures such as leapfrogging vortex rings, and it is well aligned with the symplectic/KAM techniques developed later in the paper.

\medskip
\noindent\textbf{Elliptic operator in $(\varrho,z)$-variables.}
Using \eqref{stream-operat} and the identities
$$
\partial_r = r\,\partial_\varrho,\qquad
\partial_r^2\Psi=\partial_\varrho\Psi+2\varrho\,\partial_\varrho^2\Psi,
$$
one readily checks that \eqref{stream-operat} becomes
\begin{equation*}
\mathbf{L}\Psi(\varrho,z)=\mathbf{q}(\varrho,z),
\qquad
\mathbf{L}:=-\Big(\partial_\varrho^2+\frac{1}{2\varrho}\partial_z^2\Big).
\end{equation*}
In particular, \eqref{eq:stream-Psi} implies that $(\pi\sqrt2)^{-1}G$ is the Green kernel of
$\mathbf{L}$ on $(0,\infty)\times\R$, i.e. in the sense of distributions
$$
\frac{1}{\pi\sqrt2}\,\mathbf{L}_{(\varrho,z)}G(\varrho,z;\varrho',z')
=\delta(\varrho-\varrho')\,\delta(z-z').
$$
Consequently, $G(\cdot,\cdot;\varrho',z')$ is $\mathbf{L}$-harmonic away from its singularity:
\begin{equation}\label{eq:G-harmonic}
\left(2\varrho\,\partial_\varrho^2+\partial_z^2\right)G(\varrho,z;\varrho',z')=0,
\qquad (\varrho,z)\neq(\varrho',z').
\end{equation}
Another way to get this latter identity is to proceed as follows. Let
$$
s=\frac{2(\sqrt{\varrho}-\sqrt{\varrho'})^2+(z-z')^2}{2\sqrt{\varrho\varrho'}}.
$$
A direct computation gives
\begin{equation*}
\left(2\varrho\partial_\varrho^2+\partial_z^2\right)G(\varrho,z;\varrho',z')
=\frac{(\varrho\varrho')^{1/4}}{2\varrho}
\Bigg[{s(s+4)}J''(s)+2({s+2})J'(s)-\tfrac{3}{4}J(s)\Bigg].
\end{equation*}
\medskip
Since $J$ satisfies \eqref{eq:J-ODE0}, the bracket vanishes and \eqref{eq:G-harmonic} follows.

\subsection{Derivation of the contour dynamics equation}

A central theme in vortex dynamics is the study of solutions to the Euler equations that are concentrated on thin regions of vorticity, typically of small cross-section compared to their distance of separation. In the axisymmetric without swirl setting, such configurations correspond to vortex rings (or, in the $(\varrho,z)$-plane, to compact ``patches'' of relative vorticity).  The \emph{contour dynamics method}, introduced in the two-dimensional setting (see \cite{DZ}), provides a powerful tool to reduce the PDE governing vorticity transport to an integro-differential equation for the evolution of the vortex boundaries. This approach will be essential for capturing the fine-scale interactions that give rise to leapfrogging motions.

Throughout the paper, we identify  the half plane $\Pi_:=(0,\infty)\times\R$  with a subset of    $\C$ through the mapping
$(\varrho,z)\leftrightarrow \varrho+i z$. In this identification, the real inner product is
$$
a\cdot b:=\textnormal{Re}(a\,\overline{b}),\qquad a,b\in\C,
$$
and multiplication by $i$ corresponds to a rotation by $\pi/2$.

\medskip
\noindent\textbf{Ansatz for concentrated vorticity.}
Fix a small parameter $\varepsilon\in (0,1)$ and consider solutions to \eqref{eq:tilde-q} that  concentrate in two small regions $\mathcal{D}_1(t)$ and $\mathcal{D}_2(t)$ of characteristic size  $O(\varepsilon)$.  We set
\begin{equation}\label{vorticity}
\mathbf{q}(t,\cdot)=\tfrac{1}{\varepsilon^2}{\bf 1}_{\mathcal{D}_1(t)}+\tfrac{1}{\varepsilon^2}{\bf 1}_{\mathcal{D}_2(t)}.
\end{equation}
The factor $\varepsilon^{-2}$ ensures that each patch carries $O(1)$ total mass in the $(\varrho,z)$ variables. 
We assume that each patch admits the decomposition
\begin{equation*}
\mathcal{D}_j(t)=\mathcal{P}_j(t)+\varepsilon D_j(t)+i \varepsilon|\ln\varepsilon|^{-1}\mathcal{V}_j(t), \qquad j=1,2,
\end{equation*}
where $\mathcal{P}_j=(\mathcal{P}_{j,1},\mathcal{P}_{j,2})\in \Pi_+$ the moving core of the $j$-th vortex, $ \mathcal{V}_j\in\mathbb{R}$  is a vertical displacement (see Section~\ref{sec:modeone}), and $D_j(t)\subset\R^2$ is a simply connected domain of $O(1)$ size, localized near an ellipse.

\medskip
\noindent\textbf{Slow time and boundary parametrization.} The self-induced motion of slender vortex rings involves a logarithmic time scale. We therefore introduce the slow time
$$
\tau=|\ln\varepsilon|\, t,
$$
and write
$$
D_j(t)=\mathbf{D}_j(|\ln\varepsilon|t), \quad \mathcal{V}_j(t)=\mathbf{V}_j(|\ln\varepsilon|t), \quad \mathcal{P}_j(t)=P_j(|\ln\varepsilon|\, t), \quad  
P_j(\tau)=(p_{j,1}(\tau),p_{j,2}(\tau)).
$$
We parametrize the boundary  $\partial\mathbf{D}_j$  by
\begin{equation}\label{parametrization}
\gamma_j(\tau,\theta)=P_j(\tau)+i\varepsilon|\ln\varepsilon|^{-1}\mathbf{V}_j(\tau)+\varepsilon \, w_j(\tau,\theta)\,\mathcal{Z}_j(\tau,\theta),
\end{equation}
where
$$
w_j(\tau,\theta)=\sqrt{1+2\varepsilon f_j(\tau,\theta)},
$$ 
encodes small shape deformations of the core, and
$$
\mathcal{Z}_j(\tau,\theta)=(2p_{j,1}(\tau))^{1/4}\cos(\theta)+ i\,(2p_{j,1}(\tau))^{-1/4}\sin(\theta),
$$
is the anisotropic ``elliptic'' parametrization dictated by the geometry of the kernel $G$ in
\eqref{eq:G}.
 This choice ensures that, in the singular limit regime, the evolution of the vortex rings naturally inherits a Hamiltonian structure. 

\medskip
\noindent\textbf{Contour dynamics equation.}
The evolution of ${\mathbf{q}}$ is governed by the transport equation \eqref{eq:tilde-q}. A standard argument in contour dynamics then yields the following evolution law for each boundary parametrization:
\begin{equation}\label{eq-gamma-tilde0}
|\ln\varepsilon|\, \partial_\tau \gamma_j(\tau,\theta)\cdot i\partial_\theta \gamma_j(\tau,\theta)
= \partial_\theta \big\{ \Psi(\gamma_j(\tau,\theta))\big\},
\end{equation}
where $\Psi$ is the stream function defined in \eqref{eq:stream-Psi}. 
Equation \eqref{eq-gamma-tilde0} expresses that the \emph{normal} velocity of the contour is determined by the induced flow $\Psi$, up to a reparametrization of the boundary. \\
We now compute explicitly the left-hand side of \eqref{eq-gamma-tilde0}. This requires expansions of both $\partial_\tau \gamma_j$ and $i\partial_\theta \gamma_j$. 
From the parametrization \eqref{parametrization}, we compute
\begin{align*}
\partial_\tau \gamma_j(\tau,\theta)
&= \dot{P}_j(\tau)+i\varepsilon|\ln\varepsilon|^{-1} \dot{\mathbf{V}}_j(\tau) 
   + \varepsilon^2 w_j^{-1}(\tau,\theta)\, \partial_\tau f_j(\tau,\theta)\, \mathcal{Z}_j(\tau,\theta) \\
&\quad + \tfrac12 \varepsilon\, \dot{p}_{j,1}(\tau)\, w_j(\tau,\theta) 
      \Big( (2p_{j,1}(\tau))^{-\frac34}\cos(\theta) - i(2p_{j,1}(\tau))^{-\frac54}\sin(\theta)\Big),
\end{align*}
where we use the notation $\dot{P}=\tfrac{d P}{d\tau}.$ For simplicity, we omit the explicit dependence on $(\tau,\theta)$ and rewrite:
$$
\partial_\tau \gamma_j
= \dot{P}_j+i\varepsilon|\ln\varepsilon|^{-1}\dot{\mathbf{V}}_j 
  + \varepsilon^2 w_j^{-1}\, \partial_\tau f_j\, \mathcal{Z}_j
  + \tfrac14 \varepsilon (p_{j,1})^{-1}\dot{p}_{j,1}\, w_j\, \overline{\mathcal{Z}_j}.
$$
Concerning the tangential derivative, differentiating \eqref{parametrization} with respect to $\theta$ gives
\begin{align*}
i\partial_\theta \gamma_j
&= \varepsilon^2 w_j^{-1}\, \partial_\theta f_j\, i\,\mathcal{Z}_j  - \varepsilon w_j \Big( i (2p_{j,1})^{1/4}\sin(\theta) + (2p_{j,1})^{-1/4}\cos(\theta)\Big).
\end{align*}
Equivalently, in a more compact form:
$$
i\partial_\theta \gamma_j
= \varepsilon^2 w_j^{-1}\, \partial_\theta f_j\, i\mathcal{Z}_j
  + \varepsilon w_j\, i\partial_\theta \mathcal{Z}_j.
$$  
Taking the inner complex product, a straightforward computation yields
\begin{align*}
\partial_\tau \gamma_j\cdot i\partial_\theta \gamma_j
&= \left(\dot{P}_j+i\varepsilon|\ln\varepsilon|^{-1}\dot{\mathbf{V}}_j\right)\cdot i\partial_\theta \gamma_j  - \varepsilon^3 \partial_\tau f_j
      - \varepsilon^3 \tfrac{\dot{p}_{j,1}}{4p_{j,1}}\, \partial_\theta f_j\,\sin(2\theta) \\ &\quad - \varepsilon^2 \tfrac{\dot{p}_{j,1}}{4p_{j,1}}\, w_j^2\,\cos(2\theta).
\end{align*}
This can be reorganized in the compact form
\begin{align*}
\partial_\tau \gamma_j\cdot i\partial_\theta \gamma_j
&= \left(\dot{P}_j+i\varepsilon|\ln\varepsilon|^{-1}\dot{\mathbf{V}}_j\right)\cdot i\,\partial_\theta \gamma_j
  - \varepsilon^3 \partial_\tau f_j  - \varepsilon^2 \tfrac{\dot{p}_{j,1}}{8p_{j,1}}\, 
     \partial_\theta\!\left(w_j^2\sin(2\theta)\right).
\end{align*}
Substituting this expansion into \eqref{eq-gamma-tilde0}, we obtain
\begin{align}\label{eq-gamma-tilde}
\nonumber\mathbf{G}_j(\varepsilon,f_1,f_2)(\tau,\theta)
&:= \varepsilon^3|\ln\varepsilon|\, \partial_\tau f_j(\tau,\theta)
   + \partial_\theta\left\{\Psi(\gamma_j(\tau,\theta))\right\}+ \varepsilon^2|\ln\varepsilon|\, \tfrac{\dot{p}_{j,1}(\tau)}{8p_{j,1}(\tau)}\, 
     \partial_\theta\!\left(w_j^2(\tau,\theta)\sin(2\theta)\right)\\ &\quad  - \big(|\ln\varepsilon|\, \dot{P}_j(\tau)+i\varepsilon\, \dot{\mathbf{V}}_j(\tau)\big)\cdot i\partial_\theta \gamma_j(\tau,\theta)
=0.
\end{align}

\medskip
\noindent\textbf{Evaluation of $\Psi(\gamma_j)$.}
Recalling \eqref{vorticity} and \eqref{eq:stream-Psi}, and using the change of variables
$$\rho^\prime+i z^\prime=P_k(\tau)+i\,\varepsilon|\ln\varepsilon|^{-1}\mathbf V_k(\tau)+\varepsilon Y,
$$
we get
\begin{align*}
\Psi(\gamma_j(\tau,\theta))
&=\frac{1}{\pi\sqrt{2}}\sum_{k=1}^2
\int_{\mathbf{D}_k(\tau)}
G\Big(\gamma_j(\tau,\theta),\,P_k(\tau)+i\,\varepsilon|\ln\varepsilon|^{-1}\mathbf{V}_k(\tau)+\varepsilon Y\Big)\,dY.
\end{align*}
For each $k\in\{1,2\},$ we introduce the anisotropic scaling
$$
A_k(\tau):=
\begin{pmatrix}
(2p_{k,1}(\tau))^{1/4} & 0\\
0 & (2p_{k,1}(\tau))^{-1/4}
\end{pmatrix},
\qquad Y=A_k(\tau)\,y,
$$
so that $\det A_k(\tau)=1$ and therefore Lebesgue measure is preserved. Set
$$
\widetilde{\mathbf{D}}_k(\tau):=A_k(\tau)^{-1}\mathbf{D}_k(\tau).
$$
Then
\begin{align*}
\Psi(\gamma_j(\tau,\theta))
&=\frac{1}{\pi\sqrt{2}}\sum_{k=1}^2
\int_{\widetilde{\mathbf{D}}_k(\tau)}
G\Big(\gamma_j(\tau,\theta),\,P_k(\tau)+i\,\varepsilon|\ln\varepsilon|^{-1}\mathbf{V}_k(\tau)
+\varepsilon A_k(\tau)y\Big)\,dy.
\end{align*}
Assuming $\widetilde{\mathbf{D}}_k(\tau)$ is star-shaped with respect to the origin, we parametrize its
boundary as a radial graph
$$
\widetilde{\gamma}_k(\tau,\eta)=w_k(\tau,\eta)\,e^{i\eta},\qquad \eta\in\T,
$$
and pass to polar coordinates $y=\rho e^{i\eta}$. By setting
$$
\mathcal{Z}_k(\tau,\eta):=A_k(\tau)(\cos\eta,\sin\eta)
=(2p_{k,1}(\tau))^{\frac14}\cos\eta+i(2p_{k,1}(\tau))^{-\frac14}\sin\eta,
$$
we obtain
\begin{align}\label{Pol-ah1}
\Psi(\gamma_j(\tau,\theta))
&=\Psi_1(\varepsilon,f_1,f_j)(\tau,\theta)+\Psi_2(\varepsilon,f_2,f_j)(\tau,\theta),
\end{align}
where, for $k=1,2$,
\begin{align*}
\Psi_k(\varepsilon,f_k,f_j)(\tau,\theta)
:=\frac{1}{\pi\sqrt{2}}\int_0^{2\pi}\!\int_0^{w_k(\tau,\eta)}
G\Big(\gamma_j(\tau,\theta),P_k(\tau)+i\varepsilon|\ln\varepsilon|^{-1}\mathbf{V}_k(\tau)
+\varepsilon\rho \mathcal{Z}_k(\tau,\eta)\Big)\rho d\rho d\eta.
\end{align*}
Finally, dilating the radial variable $\rho\mapsto \rho\,w_k(\tau,\eta)$ yields
\begin{align}\label{TH-1}
\nonumber &\Psi_k(\varepsilon,f_k,f_j)(\tau,\theta)
=\\ &\frac{1}{\pi\sqrt{2}}\int_0^{2\pi}\!\int_0^{1}
G\Big(\gamma_j(\tau,\theta), (P_k+i\varepsilon|\ln\varepsilon|^{-1}\mathbf{V}_k)(\tau)
+\varepsilon\rho w_k(\tau,\eta)\mathcal{Z}_k(\tau,\eta)\Big)
w_k^2(\tau,\eta)\rho d\rho d\eta.
\end{align}

\section{Filament dynamics}\label{sec:filamentdyn}

In this section, we show that the dynamics of two interacting vortex filaments can  be formulated as a four–dimensional Hamiltonian system, since each filament is represented by two coordinates. At first sight, this places the problem in a setting that appears too complex to be integrable. However,  the system possesses crucial symmetries. On the one hand, the equations are invariant under vertical translations; on the other, the two filaments are indistinguishable and the dynamics remains unchanged when they are exchanged.
Using these invariances, the four–dimensional dynamics reduces to a planar Hamiltonian system. In this reduced setting, the motion of the filaments is governed by a single two–dimensional Hamiltonian function, and the system becomes completely integrable: every trajectory lies on a level curve of the Hamiltonian, and the long–time dynamics are entirely determined by the geometry of these curves.

\subsection{Symmetry reduction for suitable 4D Hamiltonian systems}\label{sec-sym-red}

The purpose of this section is to present the symmetry reduction in a general Hamiltonian framework, before applying it later to the specific case of interacting vortex filaments.
\\
Let $\mathcal{G}: \R^4\setminus \Delta_4 \to \R$ and $\mathcal{G}_1:\R^3\setminus \Delta_3 \to \R$ be two smooth functions, where
$$
\Delta_4=\{(X,X)\colon X\in \R^2\}, 
\qquad 
\Delta_3=\{(x,x,0)\colon x\in\R\}.
$$
Here $X_1=(x_{1,1},x_{1,2})$ and $X_2=(x_{2,1},x_{2,2})$ denote the positions of the two degrees of freedom.  
We consider solutions $(X_1(\tau),X_2(\tau))$ of the Hamiltonian system
\begin{equation}\label{System-A1}
\dot X_j= \nabla_{X_j}^\perp \mathcal{G}(X_1,X_2),
\qquad j=1,2,
\end{equation}
\medskip
where $\nabla^\perp f:=J\nabla f$ with $J=\begin{psmallmatrix}0&-1\\1&0\end{psmallmatrix}$. 
We impose the following structural assumptions:
\begin{enumerate}[label=\roman*)]
\item \emph{Partial translation invariance:} there exists $\mathcal{G}_0$ such that
$$
\mathcal{G}(X_1,X_2)=\mathcal{G}_0(x_{1,1},x_{2,1},x_{1,2}-x_{2,2}).
$$
\item \emph{Symmetry under exchange and reflection:}
\begin{equation}\label{sym-ident}
\mathcal{G}_0(x_{1,1},x_{2,1},x_2)=\mathcal{G}_0(x_{2,1},x_{1,1},x_2)=\mathcal{G}_0(x_{1,1},x_{2,1},-x_2).
\end{equation}
\end{enumerate}

\medskip
We emphasize that the two properties imply in particular that
$$
\mathcal{G}(X_1,X_2)=\mathcal{G}(X_2,X_1).
$$
Under these hypotheses, the system \eqref{System-A1} can be rewritten as
\begin{equation}\label{eq-sys1}
\begin{cases}
\dot{x}_{1,1}=-(\partial_{3}\mathcal{G}_0)(x_{1,1},x_{2,1},x_{1,2}-x_{2,2}),\\[0.4em]
\dot{x}_{1,2}=(\partial_{1}\mathcal{G}_0)(x_{1,1},x_{2,1},x_{1,2}-x_{2,2}),\\[0.4em]
\dot{x}_{2,1}=(\partial_{3}\mathcal{G}_0)(x_{1,1},x_{2,1},x_{1,2}-x_{2,2}),\\[0.4em]
\dot{x}_{2,2}=(\partial_{2}\mathcal{G}_0)(x_{1,1},x_{2,1},x_{1,2}-x_{2,2}).
\end{cases}
\end{equation}
Introducing the relative vertical coordinate
$$
x_2:=x_{1,2}-x_{2,2},
$$
we obtain the reduced system
\begin{equation}\label{System1}
\begin{cases}
\dot{x}_{1,1}=-(\partial_{3})\mathcal{G}_0)(x_{1,1},x_{2,1},x_2),\\[0.4em]
\dot{x}_{2,1}=(\partial_{3}\mathcal{G}_0)(x_{1,1},x_{2,1},x_2),\\[0.4em]
\dot{x}_2=(\partial_{1}\mathcal{G}_0)(x_{1,1},x_{2,1},x_2)-(\partial_{2}\mathcal{G}_0)(x_{1,1},x_{2,1},x_2).
\end{cases}
\end{equation}
\medskip
From the first two equations, we deduce the conservation law
$$
x_{1,1}(\tau)+x_{2,1}(\tau)=x_{1,1}(0)+x_{2,1}(0)=:2\kappa,
$$
 which expresses a balance constraint along the horizontal direction. In the abstract setting one only has $\kappa\in\mathbb R$; in the filament application below, $\kappa>0$.

Substituting $x_{2,1}=2\kappa-x_{1,1}$ into \eqref{System1}, one finds
$$
\dot{x}_{1,1}=-\partial_{3}\,\mathcal{G}_0(x_{1,1},2\kappa-x_{1,1},x_2), 
\qquad
\dot{x}_2=\partial_{x_{1,1}}\big(\mathcal{G}_0(x_{1,1},2\kappa-x_{1,1},x_2)\big).
$$
By introducing the shifted coordinate $x_{1,1}=\kappa+x_1$, the system takes the canonical Hamiltonian form
\begin{equation}\label{system0}
\dot{x}_1=-\partial_{x_2}H(x_1,x_2), 
\qquad
\dot{x}_2=\partial_{x_1}H(x_1,x_2),
\end{equation}
with the  Hamiltonian

$$
H(x_1,x_2):=\mathcal{G}_0(\kappa+x_1,\kappa-x_1,x_2).
$$
\medskip
The symmetry conditions \eqref{sym-ident} yield
$$
H(-x_1,x_2)=H(x_1,-x_2)=H(x_1,x_2),
$$
so the Hamiltonian is invariant under both reflection in $x_1$ and reflection in $x_2$. These discrete symmetries imply strong reversibility properties for the dynamics.\\ Next, we examine how the symmetry of the Hamiltonian will be reflected on the structure of the orbits. 
\begin{pro}\label{prop-mouh}
Let $(x_1(\tau),x_2(\tau))$ be a periodic solution of \eqref{system0} with minimal period $T$, and assume that its orbit is a simple closed curve surrounding the origin. If
$$
(x_1(0),x_2(0))=(x_{1,0},0),\qquad x_{1,0}>0,
$$
then
$$
x_1(-\tau)=x_1(\tau),\qquad x_2(-\tau)=-x_2(\tau),
$$
$$
x_1(T-\tau)=x_1(\tau),\qquad x_2(T-\tau)=-x_2(\tau),
$$
$$
x_1\!\left(\tau+\tfrac{T}{2}\right)=-x_1(\tau),\qquad
x_2\!\left(\tau+\tfrac{T}{2}\right)=-x_2(\tau),
$$
and
$$
x_1\!\left(\tfrac{T}{2}-\tau\right)=-x_1(\tau),\qquad
x_2\!\left(\tfrac{T}{2}-\tau\right)=x_2(\tau).
$$
Moreover, the functions $\tau\mapsto \dot{x}_{1,2}(\tau)$ and $\tau\mapsto \dot{x}_{2,2}(\tau)$ are $T$-periodic and even.
\end{pro}

\begin{proof}
Since $H$ is even in $x_2$, the transformed curve
$$
\widehat{x}_1(\tau):=x_1(-\tau),\qquad \widehat{x}_2(\tau):=-x_2(-\tau)
$$
solves the same Hamiltonian system as $(x_1,x_2)$ and satisfies
$$
(\widehat{x}_1(0),\widehat{x}_2(0))=(x_{1,0},0).
$$
By uniqueness of the Cauchy problem,
$$
x_1(-\tau)=x_1(\tau),\qquad x_2(-\tau)=-x_2(\tau).
$$
Using $T$-periodicity gives
$$
x_1(T-\tau)=x_1(\tau),\qquad x_2(T-\tau)=-x_2(\tau).
$$
Next, since $H$ is even in both variables, the central reflection
$$
\widetilde{x}_1(\tau):=-x_1(\tau),\qquad \widetilde{x}_2(\tau):=-x_2(\tau)
$$
is again a solution of \eqref{system0}. Its image is the same simple closed orbit, so there exists $s\in(0,T)$ such that
$$
(\widetilde{x}_1(\tau),\widetilde{x}_2(\tau))=(x_1(\tau+s),x_2(\tau+s)).
$$
Applying this identity twice yields
$$
x_j(\tau+2s)=x_j(\tau),\qquad j=1,2.
$$
Since $T$ is the minimal period, we must have $2s=T$. Hence
$$
x_1\!\left(\tau+\tfrac{T}{2}\right)=-x_1(\tau),\qquad
x_2\!\left(\tau+\tfrac{T}{2}\right)=-x_2(\tau).
$$
Combining this with the first symmetry pair gives
$$
x_1\!\left(\tfrac{T}{2}-\tau\right)
=-x_1(-\tau)=-x_1(\tau),
\qquad
x_2\!\left(\tfrac{T}{2}-\tau\right)
=-x_2(-\tau)=x_2(\tau).
$$
Finally,
$$
\dot{x}_{1,2}(\tau)
=\partial_1\mathcal{G}_0\bigl(\kappa+x_1(\tau),\kappa-x_1(\tau),x_2(\tau)\bigr),
$$
and similarly
$$
\dot{x}_{2,2}(\tau)
=\partial_2\mathcal{G}_0\bigl(\kappa+x_1(\tau),\kappa-x_1(\tau),x_2(\tau)\bigr).
$$
Their $T$-periodicity follows from that of $(x_1,x_2)$. Using the evenness of $x_1$, the oddness of $x_2$, and the symmetry
$$
\mathcal{G}_0(x,y,-z)=\mathcal{G}_0(x,y,z),
$$
we obtain
$$
\dot{x}_{1,2}(-\tau)=\dot{x}_{1,2}(\tau),\qquad
\dot{x}_{2,2}(-\tau)=\dot{x}_{2,2}(\tau).
$$
This proves the claim.
\end{proof}
\subsection{Filament dynamics and scaling}\label{sec:ring-dynamics}
The filament motion in the axisymmetric setting under consideration can be described, to leading order, by the dynamics of two Dirac masses in the half-plane  $
\Pi_+ = \{(r,z)\,:\, r>0\},$  
interacting through a Hamiltonian system. This framework has been extensively studied in the literature, and we shall recall its derivation in the following sections. 
More precisely, consider two points $(P_1,P_2)\in\R^4$, with $P_j=(p_{j,1},p_{j,2})$, for $j=1,2$. In order to capture the leapfrogging regime, we shall introduce the following anisotropic scaling
\begin{align}\label{scaling-points}
p_{j,1}=\kappa+\tfrac{ 1}{2}(-1)^{j+1} r_\varepsilon (2\kappa)^\frac14 \mathtt{x}_1, 
\qquad {p_{1,2}-p_{2,2}}=r_\varepsilon(2\kappa)^{-\frac14}\mathtt{x}_2, \qquad r_\varepsilon:=(2\kappa)^\frac14|\ln\varepsilon|^{-\frac12},
\end{align}
for $\kappa>0.$ That scaling is designed to capture the small-amplitude, fast (logarithmically rescaled) dynamics associated with leapfrogging. Motivated by Proposition \ref{prop:G1-trivial}, we impose the following equation for the point vortices:
\begin{equation}\label{eq-points-general}
\begin{aligned}
    |\ln\varepsilon | \dot{P}_j(\tau) =\tfrac{1}{\sqrt{2}}\nabla_{P_j} ^\perp G(P_1,P_{2})&+ \tfrac{i}{4}(2 p_{j,1})^{-\frac12}\Big[|\ln\varepsilon|{+}\tfrac{1}{4}\big(5\ln8+3\ln(p_{j,1})-1\big)\Big],
\end{aligned}
\end{equation}
with $G$ as in \eqref{eq:G}. A direct computation shows that the vector field in \eqref{eq-points-general} has the Hamiltonian structure
$$
|\ln\varepsilon| \dot{P}_j(\tau)=\nabla^\perp_{P_j} H(P_1,P_2),
$$
with Hamiltonian
\begin{align}\label{decompo-Hamil}
    H(P_1,P_2):=&\tfrac{1}{\sqrt{2}}G(P_1, P_2)+\mathtt{G}(P_1,P_2),
\end{align}
where
\begin{align*}
   \mathtt{G}(P_1,P_2):=&\sum_{j=1}^2\frac{(p_{j,1})^\frac12}{2\sqrt{2}}\left(|\ln\varepsilon|{-}\tfrac74{+}\tfrac54\ln(8){+}\tfrac34\ln(p_{j,1})\right).
\end{align*}
We observe that the Hamiltonian satisfies the following symmetry properties.
\begin{enumerate}[label=\roman*)]
\item Partial translation invariance:
$$
H(P_1,P_2)=H_0(p_{1,1},p_{2,1},p_{1,2}-p_{2,2})
$$
\item Symmetry under exchange and reflection: 
\begin{align*}
H_0(x,y,z)=H_0(y,x,z)=H_0(x,y,-z).
\end{align*}

    \end{enumerate}
    These are the same symmetries as in \eqref{System-A1}, and hence the symmetry reduction of Section~\ref{sec-sym-red} applies.
Hence, by using the anisotropic scaling \eqref{scaling-points} we have that, by the chain rule, the system for the point vortices becomes
\begin{equation}\label{system001}
		\left\lbrace\begin{array}{ll}
			\dot{\mathtt{x}_1}=-\partial_{\mathtt{x}_2} \mathcal{H}(\mathtt{x}_1,\mathtt{x}_2),\\
			\dot{\mathtt{x}_2}=\partial_{\mathtt{x}_1} \mathcal{H}(\mathtt{x}_1,\mathtt{x}_2),
		\end{array}\right.   
	\end{equation}
    with
    \begin{align*}
    \mathcal{H}(\mathtt{x}_1,\mathtt{x}_2)
    &=
{\sqrt{\tfrac{2}{\kappa}}}H_0\left(\kappa+{\tfrac12}r_\varepsilon (2\kappa)^{\frac14}\mathtt{x}_1, \kappa-{\tfrac12}r_\varepsilon (2\kappa)^{\frac14}\mathtt{x}_1, r_\varepsilon (2\kappa)^{-\frac14} \mathtt{x}_2\right).
    \end{align*}
    The goal is to expand in $\varepsilon$ the Hamiltonian $\mathcal{H}.$ From \eqref{decompo-Hamil} we infer 
    \begin{align*}
\mathcal{H}(\mathtt{x}_1,\mathtt{x}_2)&=\tfrac{1}{\sqrt{\kappa}}G\big(Z+r_\varepsilon X, Z+r_\varepsilon Y\big)+F_0(\mathtt{x}_1),
    \end{align*}
    with the notation
    $$
Z=(\kappa,0),\quad  X=({\tfrac12}(2\kappa)^{\frac14}\mathtt{x}_1,(2\kappa)^{-\frac14} \mathtt{x}_2),\quad Y=(-\tfrac12(2\kappa)^{\frac14}\mathtt{x}_1,0)
 $$
   and
        \begin{align*}
   F_0(\mathtt{x}_1):=&\sum_{j=1}^2\tfrac{1}{2\sqrt{\kappa}}\left(\kappa+\tfrac{(-1)^{j+1}}{2}r_\varepsilon (2\kappa)^{\frac14}\mathtt{x}_1\right)^\frac12\left(|\ln\varepsilon|{-}\tfrac74{+}\tfrac54\ln(8){+}\tfrac34\ln\left(\kappa+\tfrac{(-1)^{j+1}}{2}r_\varepsilon (2\kappa)^{\frac14}\mathtt{x}_1\right)\right).
\end{align*}
From Taylor expansion
\begin{align}\label{Es-F0}
  \nonumber F_0(\mathtt{x}_1)=&|\ln\varepsilon|{-}\tfrac{1}{4}\left(7-5\ln(8)-3\ln(\kappa)\right)-\tfrac18r_\varepsilon^2 |\ln\varepsilon|(2\kappa)^{-\frac32}\mathtt{x}_1^2+r_\varepsilon^2F_1(\mathtt{x}_1)\\
   =&|\ln\varepsilon|{-}\tfrac{1}{4}\left(7-5\ln(8)-3\ln(\kappa)\right)-\tfrac{1}{16{\kappa}}\mathtt{x}_1^2+r_\varepsilon^2F_1(\mathtt{x}_1),
\end{align}
with
$$
|F_1(\mathtt{x}_1)|\lesssim \mathtt{x}_1^2.$$
Using Lemma \ref{prop-split-1} we get the expansion 
    $$
G\big(Z+r_\varepsilon X, Z+r_\varepsilon Y\big)=|\ln(r_\varepsilon)|\sum_{n=0}^3\mathcal{A}_n(X,Y,Z) r_{\varepsilon}^n+\sum_{n=0}^2\mathcal{B}_n(X,Y,Z) r_{\varepsilon}^n+O(r_{\varepsilon}^3),
$$
where in our case:
\begin{align*}
   \mathcal{A}_0(X,Y,Z)&=\sqrt{\kappa},\quad
   \mathcal{A}_1(X,Y,Z)=\mathcal{A}_3(X,Y,Z)=0,\\ 
   \mathcal{A}_2(X,Y,Z)&=\tfrac{\sqrt{2}}{64\kappa}\left(3\mathtt{x}_2^2-\mathtt{x}_1^2\right),\\
 \mathcal{B}_0(X,Y,Z)&=\sqrt{\kappa}\left(\ln8-2+\tfrac34\ln(2\kappa)\right)-\tfrac{\sqrt{\kappa}}{2}\ln\left({\mathtt{x}_1^2+\mathtt{x}_2^2}\right),\quad \mathcal{B}_1(X,Y,Z)=0,
   \\
   \mathcal{B}_2(X,Y,Z)&=-\tfrac{\mathtt{x}_1^2}{8\sqrt{2} \kappa}\left(\ln8-2+\tfrac34 \ln(2\kappa)-\tfrac12\ln\left(\mathtt{x}_1^2+\mathtt{x}_2^2\right)\right)-\tfrac{\mathtt{x}_1^2}{16\sqrt{2} \kappa }\tfrac{3\mathtt{x}_1^2+2\mathtt{x}_2^2}{\mathtt{x}_1^2+\mathtt{x}_2^2}\\  &\quad+ \left(\mathtt{x}_1^2+\mathtt{x}_2^2\right)\left(\tfrac{9\ln2-1}{32\sqrt{2}\,\kappa}+\tfrac{9}{128\sqrt{2}\,\kappa}\ln(2\kappa)-\tfrac{3}{64\sqrt{2}\kappa}\ln\left(\mathtt{x}_1^2+\mathtt{x}_2^2\right)  \right).
\end{align*}
Consequently,
\begin{align}\label{Es-F1}
\tfrac{1}{\sqrt{\kappa}}G\big(Z+r_\varepsilon X, Z+r_\varepsilon Y\big)&=\mathcal{C}_0(\varepsilon,\kappa)-\tfrac12\,\ln\left({\mathtt{x}_1^2+\mathtt{x}_2^2}\right)+r_{\varepsilon}^2\,F_2(\mathtt{x}_1,\mathtt{x}_2),
\end{align}
with $\mathcal{C}_0$ is a constant and 
$$
|F_2(\mathtt{x}_1,\mathtt{x}_2)|\lesssim (\mathtt{x}_1^2+\mathtt{x}_2^2)|\ln(\mathtt{x}_1^2+\mathtt{x}_2^2)|.
$$
Putting together \eqref{Es-F0} and \eqref{Es-F1} we deduce by removing the constant from the Hamiltonian that
$$
\mathcal{H}(\mathtt{x}_1,\mathtt{x}_2)
=-\tfrac12\ln\bigl(\mathtt{x}_1^2+\mathtt{x}_2^2\bigr)
-\tfrac{1}{16\kappa}\mathtt{x}_1^2
+r_\varepsilon^2 F(\mathtt{x}_1,\mathtt{x}_2),$$
where
\begin{align}\label{F-est1}
|F(\mathtt{x}_1,\mathtt{x}_2)|\lesssim (\mathtt{x}_1^2+\mathtt{x}_2^2)|\ln(\mathtt{x}_1^2+\mathtt{x}_2^2)|\quad\hbox{and}\quad |\nabla F(\mathtt{x}_1,\mathtt{x}_2)|\lesssim (|\mathtt{x}_1|+|\mathtt{x}_2|)|\ln(\mathtt{x}_1^2+\mathtt{x}_2^2)|.
\end{align}
Therefore, the system \eqref{system001} becomes
\begin{align}\label{system002}
		\left\lbrace\begin{array}{ll}
			\dot{\mathtt{x}_1}=\frac{\mathtt{x}_2}{\mathtt{x}_1^2+\mathtt{x}_2^2}-r_{\varepsilon}^2\partial_{\mathtt{x}_2}F(\mathtt{x}_1,\mathtt{x}_2),\\
			\dot{\mathtt{x}_2}=-\frac{\mathtt{x}_1}{\mathtt{x}_1^2+\mathtt{x}_2^2}-\frac{1}{8\kappa}\mathtt{x}_1+r_{\varepsilon}^2\partial_{\mathtt{x}_1} F(\mathtt{x}_1,\mathtt{x}_2).
		\end{array}\right.   
	\end{align}

\subsection{Periodic motion for the unperturbed system}\label{sec:lep-vortex-pairs}

In this section, we study the precise conditions that produce the leapfrogging motion of the limiting system in \eqref{system002} and derive quantitative control of the period together with its analytic dependence on the initial configuration. In our setting, leapfrogging corresponds to a non-rigid time–periodic motion when viewed in a uniformly translating frame. In fact, the limiting system can be identified with the dynamics of two equal vortices in the presence of uniform shear, already investigated in \cite{Kim-Has}.
\\
The limiting system in \eqref{system002} is obtained by setting $r_\varepsilon=0$, namely
\begin{align}\label{system003}
		\left\lbrace\begin{array}{l}
        \dot{\underline{\mathtt{x}}_1}=\frac{\underline{\mathtt{x}}_2}{\underline{\mathtt{x}}_1^2+\underline{\mathtt{x}}_2^2},\\
			\dot{\underline{\mathtt{x}}_2}=-\frac{\underline{\mathtt{x}}_1}{\underline{\mathtt{x}}_1^2+\underline{\mathtt{x}}_2^2}-\tfrac{1}{8\kappa}\underline{\mathtt{x}}_1,\\
             \underline{\mathtt{x}}_1(0)=\lambda>0,\quad \underline{\mathtt{x}}_2(0)=0,
		\end{array}\right.   
	\end{align}
    with the Hamiltonian
    $$
\mathcal{H}_0(\underline{\mathtt{x}}_1,\underline{\mathtt{x}}_2)=-\tfrac12\,\ln\Big({\underline{\mathtt{x}}_1^2+\underline{\mathtt{x}}_2^2}\Big)-\tfrac{1}{16\kappa}\underline{\mathtt{x}}_1^2.$$

\begin{pro}\label{prop:period}
For each $(\lambda,\kappa)\in(0,\infty)^2$, let
$$
\underline{\mathtt{x}}(\tau;\lambda,\kappa)
=
\bigl(\underline{\mathtt{x}}_1(\tau;\lambda,\kappa),\underline{\mathtt{x}}_2(\tau;\lambda,\kappa)\bigr)
$$
be the unique global solution to the system \eqref{system003}. Then the orbit issued from \((\lambda,0)\) is a simple closed curve, and its minimal period
$T_0(\lambda,\kappa)>0$ is given by
    \begin{align}\label{T(z)}
	T_0(\lambda,\kappa)
=
2\lambda^2
\bigintss_0^{1}
\frac{
e^{\alpha(1-s)}
}{
\sqrt{
e^{\alpha(1-s)}-s
}
}
\,\frac{ds}{\sqrt{s}}, \qquad \alpha:=\tfrac{\lambda^2}{8\kappa}\cdot
	\end{align}
 Moreover, the following properties hold.
    \begin{enumerate} 
    \item The trajectories and the period are separately real analytic with respect to
$\lambda$ and $\kappa$. More precisely, for every fixed
$(\tau,\kappa)\in\mathbb R\times(0,\infty)$, the maps
$\lambda\mapsto \underline{\mathtt{x}}_j(\tau;\lambda,\kappa)$ and
$\lambda\mapsto T_0(\lambda,\kappa)$ are real analytic on $(0,\infty)$,
while for every fixed $(\tau,\lambda)\in\mathbb R\times(0,\infty)$, the maps
$\kappa\mapsto \underline{\mathtt{x}}_j(\tau;\lambda,\kappa)$ and
$\kappa\mapsto T_0(\lambda,\kappa)$ are real analytic on $(0,\infty)$.

\item For every fixed $\kappa>0$, the map
$$
\lambda\longmapsto T_0(\lambda,\kappa)
$$
is strictly increasing on $(0,\infty)$. Moreover, for every compact set
$Q\subset(0,\infty)$,
$$
\inf_{\lambda\in Q}\partial_\lambda T_0(\lambda,\kappa)>0.
$$

\item For every $(\lambda,\kappa)\in(0,\infty)^2$,
$$
\frac{2\pi \lambda^2}{\sqrt{1+\alpha}}
\leqslant
T_0(\lambda,\kappa)
\leqslant
2\pi \lambda^2 e^{\alpha/2}.
$$
    \item On the trajectory, we  have the estimates
    $$
 \forall \lambda>0,\,\forall \kappa>0,\,\forall\,\tau\in\mathbb{R},\quad   \lambda^2\leqslant |\underline{\mathtt{x}}^2_1(\tau)+\underline{\mathtt{x}}^2_2(\tau)|\leqslant \lambda^2 e^{\frac{\lambda^2}{8\kappa}}.
    $$
	\end{enumerate}
\end{pro}
\begin{proof}
From the Hamiltonian structure of \eqref{system003} we infer that     
 the orbit is contained in the level set:
	\begin{align}\label{level-set}
	\underline{\mathtt{x}}_1^2(\tau)+\underline{\mathtt{x}}_2^2(\tau)=\lambda^2 e^{-\tfrac{\underline{\mathtt{x}}_1^2(\tau)-\lambda^2}{8\kappa}},\quad \underline{\mathtt{x}}_1(0)=\,\lambda>0.
	\end{align}
	 Consequently,  the orbits are contained in the planar set  
	$$\Big\{(\underline{\mathtt{x}}_1,\underline{\mathtt{x}}_2)\in\R^2,\quad \underline{\mathtt{x}}_1^2+\underline{\mathtt{x}}_2^2=\lambda^2e^{-\tfrac{\underline{\mathtt{x}}_1^2-\lambda^2}{8\kappa} }\Big\}.
	$$ 
	We can check that this set is a simple closed curve and therefore the motion is always periodic. Moreover, the orbit intersects the $\underline{\mathtt{x}}_2$ axis at the point   
	\begin{align}\label{eta-equation0}
{\mathtt{x}}_{2,\mathtt{c}}^2=\lambda^2e^{\tfrac{\lambda^2}{8\kappa}}.
	\end{align}
The orbit is two-fold and it is invariant by reflection with respect to  the two axes. Then, since the system is autonomous,  if $T=T_0(\lambda,\kappa)>0$ denotes the period then at $T/4$ the point on the curve is located at the vertical axis.
    Since
$$
\dot{\underline{\mathtt{x}}}_2(0)
=
-\frac{1}{\lambda}-\frac{\lambda}{8\kappa}<0,
$$
the orbit initially enters the lower half-plane. It follows that
$$
\underline{\mathtt{x}}_2(\tau)<0
\qquad\text{for all }\tau\in(0,T/4),
$$
and therefore
$$
\dot{\underline{\mathtt{x}}}_1(\tau)<0
\qquad\text{for all }\tau\in(0,T/4).
$$
In particular,
$$
\underline{\mathtt{x}}_1(0)=\lambda,
\qquad
\underline{\mathtt{x}}_1(T/4)=0,
\qquad
\underline{\mathtt{x}}_2(T/4)=\mathtt{x}_{2,\mathtt{c}}.
$$
	 Next, we want to find a closed formula for the period, which will be useful to derive some quantitative properties. 
To get it, we first write from \eqref{system003} and \eqref{level-set}
	$$
	\dot{\underline{\mathtt{x}}}_{1}= \frac{\underline{\mathtt{x}}_2}{\underline{\mathtt{x}}_2^2+\underline{\mathtt{x}}_{1}^2}=\frac{ \underline{\mathtt{x}}_2}{\lambda^2}e^{\tfrac{\underline{\mathtt{x}}_{1}^2-\lambda^2}{8\kappa}}\,.
	$$
	Hence, we deduce from \eqref{level-set}
	\begin{align*}
	\dot{\underline{\mathtt{x}}_{1}}(\tau)&=- {\lambda^{-2}}\,{\sqrt{\lambda^2\, e^{-\tfrac{\underline{\mathtt{x}}_{1}^2(\tau)-\lambda^2}{8\kappa}}-\underline{\mathtt{x}}_{1}^2(\tau)}}\,\,e^{\tfrac{\underline{\mathtt{x}}_{1}^2(\tau)-\lambda^2}{8\kappa}}.
	\end{align*}
 Thus, integrating this ODE allows us to get the following formula
	\begin{align*}
	T_0(\lambda,\kappa)&={4\lambda^2}\bigintss_0^{\lambda}\frac{e^{\frac{\lambda^2-\rho^2}{8\kappa}}}{
\sqrt{
\lambda^2
e^{\frac{\lambda^2-\rho^2}{8\kappa}}
-
\rho^2
}
}
\,d\rho.
	\end{align*}
    By making the change of variables $x=\tfrac{\rho^2}{8\kappa},$ we get
    $$
T_0(\lambda,\kappa)
=
2\lambda^2
\int_0^\alpha
\frac{e^{\alpha-x}}{
\sqrt{
\alpha e^{\alpha-x}-x 
}
}
\frac{dx}{\sqrt{x}},\quad \alpha=\tfrac{\lambda^2}{8\kappa}\cdot
$$
Finally, with the change of variables $x=\alpha s$, we find
$$
T_0(\lambda,\kappa)
=
2\lambda^2
\int_0^1
\frac{
e^{\alpha(1-s)}
}{
\sqrt{
e^{\alpha(1-s)}-s
}
}
\,\frac{ds}{\sqrt{s}},
$$
which is exactly \eqref{T(z)}.

\medskip
\noindent{\bf (1)} We rewrite the denominator as
$$
e^{\alpha(1-s)}-s
=
(1-s)\left(
1+\frac{e^{\alpha(1-s)}-1}{1-s}
\right)
=
(1-s)\left(
1+\alpha\int_0^1 e^{\theta\alpha(1-s)}\,d\theta
\right).
$$
Hence
$$
T_0(\lambda,\kappa)=2\lambda^2\bigintss_0^1\frac{e^{\alpha(1-s)}}{\sqrt{1+\alpha\int_0^1 e^{\theta\alpha(1-s)}\,d\theta}}\,\frac{ds}{\sqrt{s(1-s)}}.
$$
For every compact interval $I\subset(0,\infty)$, the function
$$
F(\alpha,s):=
\frac{
e^{\alpha(1-s)}
}{
\sqrt{
1+\alpha\int_0^1 e^{\theta\alpha(1-s)}\,d\theta
}
}
$$
is real analytic in $\alpha$, continuous on $I\times[0,1]$, and all its
$\alpha$-derivatives are bounded on $I\times[0,1]$. Since
$\alpha=\lambda^2/(8\kappa)$, the maps $\lambda\mapsto T_0(\lambda,\kappa)$
$\kappa\mapsto T_0(\lambda,\kappa)$ are real analytic on $(0,\infty)$. 

\medskip
\noindent Next, we shall check the analyticity of the trajectories.
Fix $(\tau_0,\lambda_0,\kappa_0)\in \mathbb R\times(0,\infty)^2$, and write
$$
F(x_1,x_2,\kappa)
:=
\left(
\frac{x_2}{x_1^2+x_2^2},
-\frac{x_1}{x_1^2+x_2^2}-\frac{x_1}{8\kappa}
\right).
$$
Then $F$ is real analytic on
$$
\bigl(\mathbb R^2\setminus\{(0,0)\}\bigr)\times(0,\infty).
$$
From \eqref{level-set}, one has
\begin{equation}\label{inequality1st}
    \underline{\mathtt{x}}_1^2(\tau;\lambda,\kappa)
+
\underline{\mathtt{x}}_2^2(\tau;\lambda,\kappa)
=
\lambda^2
\exp\!\left(
\frac{\lambda^2-\underline{\mathtt{x}}_1^2(\tau;\lambda,\kappa)}{8\kappa}
\right).
\end{equation}
Since the function
$$
u\longmapsto u\, e^{\frac{u}{8\kappa}}
$$
is strictly increasing on $[0,\infty)$, then
$$
\underline{\mathtt{x}}_1^2(\tau;\lambda,\kappa)\leqslant \lambda^2
\qquad\forall \tau\in\mathbb R,
$$
hence
\begin{equation}\label{inequality2nd}
\underline{\mathtt{x}}_1^2(\tau;\lambda,\kappa)
+
\underline{\mathtt{x}}_2^2(\tau;\lambda,\kappa)
\geqslant \lambda^2
\qquad\forall \tau\in\mathbb R.
\end{equation}
The standard theorem on analytic dependence of solutions of ODEs
on parameters and initial data yields that the flow map
$$
(\tau,\lambda,\kappa)\longmapsto
(\underline{\mathtt{x}}_1,\underline{\mathtt{x}}_2)(\tau;\lambda,\kappa)
$$
is real analytic in a neighborhood of every point
$(\tau_0,\lambda_0,\kappa_0)\in \mathbb R\times(0,\infty)^2$.      

 \medskip         \noindent {\bf (2)} Let's check the monotonicity of the period. To this end, we write
$$
T_0(\lambda,\kappa)=16\kappa \int_0^1 K_0(\alpha,s)\,\frac{ds}{\sqrt{s}},
$$
where
$$
K_0(\alpha,s)
:=
\frac{\alpha\,e^{\alpha(1-s)}}{\sqrt{e^{\alpha(1-s)}-s}},
\qquad 0\le s<1.
$$
Hence, it is enough to show that, for every fixed $s\in[0,1)$, the map
$\alpha\mapsto K_0(\alpha,s)$ is strictly increasing on $(0,\infty)$.
A direct differentiation gives
$$
\partial_\alpha K_0(\alpha,s)
=
\frac{e^{\alpha (1-s)}}{(e^{\alpha (1-s)}-s)^{3/2}}
\left[
e^{\alpha (1-s)}-s+\tfrac{\alpha (1-s)}{2}\bigl(e^{\alpha (1-s)}-2s\bigr)
\right].
$$
Therefore it remains to prove that
$$
B:=e^{\alpha (1-s)}-s+\frac{\alpha (1-s)}{2}\bigl(e^{\alpha (1-s)}-2s\bigr)>0.
$$
Writing $t:=\alpha (1-s)\ge 0$, we obtain
$$
B=(e^t-s)+\tfrac{t}{2}(e^t-2s).
$$
Since $e^t\ge 1+t$ for all $t\ge 0$, it follows that
$$
B
\ge
t+1-s+\tfrac{t}{2}(t+1-2s)
=
1-s+\tfrac{t}{2}+\tfrac{t^2}{2}+t(1-s)
>0.
$$
Thus
$$
\partial_\alpha K_0(\alpha,s)>0
\qquad\text{for every }\alpha>0,\; s\in[0,1).
$$
Hence $\lambda\mapsto T_0(\lambda,\kappa)$ is strictly increasing on $(0,\infty)$.

 \medskip        

 \noindent \medskip {\bf (3)}  Now, we shall estimate the period.  Fix $\alpha>0$. Then for every $s\in(0,1)$,
$$
\frac{1}{\sqrt{1+\alpha}}\,
\frac{1}{\sqrt{1-s}}
\;\le\;
\frac{e^{\alpha(1-s)}}{\sqrt{e^{\alpha(1-s)}-s}}\,
\;\le\;
e^{\alpha/2}\,
\frac{1}{\sqrt{1-s}}\cdot
$$
In fact, setting
$
t:=\alpha(1-s)
$
we may write 
$$
\frac{e^{\alpha(1-s)}}{\sqrt{e^{\alpha(1-s)}-s}}\,
=
\frac{e^{t/2}}{\sqrt{1-s e^{-t}}}\cdot
$$
We first estimate the factor $1-s e^{-t}$,
$$
1-s e^{-t}
=
(1-s)+s(1-e^{-t})
\le
(1-s)+(1-e^{-t}).
$$
Using the elementary inequality $1-e^{-x}\leqslant x$ for $x\geqslant 0$, we get
$$
1-s e^{-t}
\leqslant
(1-s)+t
=
(1-s)+\alpha(1-s)
=
(1+\alpha)(1-s).
$$
Therefore,
$$
1-s
\leqslant
1-s e^{-t}
\leqslant
(1+\alpha)(1-s),
$$
On the other hand,  we have
$$
1\leqslant e^{t/2}\leqslant e^{\alpha/2}.
$$
Hence,
$$
\frac{1}{\sqrt{1+\alpha}}\,
\frac{1}{\sqrt{1-s}}
\leqslant
\frac{e^{t/2}}{\sqrt{1-s e^{-t}}}
\leqslant
e^{\alpha/2}\,
\frac{1}{\sqrt{1-s}}.
$$
This proves the claim. Now, from the period formula we get
\begin{align*}
	\frac{2\lambda^2}{\sqrt{1+\alpha}}\bigintsss_0^{1}\frac{d\tau}{\sqrt{1-\tau}\sqrt{\tau}} \leqslant T_0(\lambda,\kappa)&\leqslant {2\lambda^2}e^{\alpha/2}\bigintsss_0^{1}\frac{d\tau}{\sqrt{1-\tau}\sqrt{\tau}}.
	\end{align*} 
    From the classical identity  $$
    \bigintsss_0^{1}\frac{ds}{\sqrt{\tau(1-\tau)}}=\pi,
    $$
    we deduce that
\begin{align*}
	\frac{2\pi \lambda^2}{\sqrt{1+\alpha}}\leqslant T_0(\lambda,\kappa)&\leqslant {2\pi \lambda^2}e^{\alpha/2}.
	\end{align*} 

    \medskip
\noindent{\bf (4)} The last point follows immediately from \eqref{inequality1st} and \eqref{inequality2nd}. 
This completes the proof.
\end{proof}

\subsection{Periodic motion for the perturbed system }\label{sec:lep-vortex-pairs-perturbed}
In this section, we shall explore the existence of periodic solutions to the system \eqref{system002}, which takes the Hamiltonian form 
\begin{align}\label{system004}
		\left\lbrace\begin{array}{ll}
			\dot{\mathtt{x}}_1=-\partial_{\mathtt{x}_{2}} \mathcal{H}(\mathtt{x}_{1},\mathtt{x}_2),\\
            \dot{\mathtt{x}}_{2}=\partial_{\mathtt{x}_1} \mathcal{H}(\mathtt{x}_{1},\mathtt{x}_2),\\
            \mathtt{x}_1(0)={\lambda}>0,\quad\mathtt{x}_{2}(0)=0,
		\end{array}\right.   
	\end{align}
    with the Hamiltonian
    $$
\mathcal{H}(\varepsilon,\mathtt{x}_{1},\mathtt{x}_2)=-\tfrac12\,\ln\left({\mathtt{x}_{1}^2+\mathtt{x}_2^2}\right)-\tfrac{1}{16\kappa}\mathtt{x}_{1}^2+r_{\varepsilon}^2F(\mathtt{x}_{1},\mathtt{x}_2).$$
\begin{pro}
   Let $\kappa>0$ and $\lambda\in [a,b]\subset (0,\infty)$. Then there exists $\varepsilon_0>0$ such that for any $\varepsilon\in (0,\varepsilon_0),$ the solution to \eqref{system004} is a periodic orbit enclosing a set containing the origin and invariant by reflection with respect to the real and imaginary axes. In particular,
   \begin{align*}
\mathtt{x}_1(-t)=\mathtt{x}_1(t), &\quad \mathtt{x}_2(t)=-\mathtt{x}_2(t),\quad
\mathtt{x}_1\big(t+\tfrac{T(\varepsilon,\lambda,\kappa)}{2}\big)=-\mathtt{x}_1(t), \quad \mathtt{x}_2\big(t+\tfrac{T(\varepsilon,\lambda,\kappa)}{2}\big)=-\mathtt{x}_2(t).
\end{align*}
In addition, the  period function $\lambda\in[a,b]\mapsto T(\varepsilon,\lambda,\kappa)$ is $C^1$ and satisfies
   $$
T(\varepsilon,\lambda,\kappa)=T_0(\lambda, \kappa)+O\left(\tfrac{1}{\sqrt{|\ln\varepsilon|}}\right),\quad\hbox{and}\quad
\partial_\lambda T(\varepsilon,\lambda,\kappa)=T_0'(\lambda,\kappa)+O\left(\tfrac{1}{\sqrt{|\ln\varepsilon|}}\right),
   $$
   where $T_0(\lambda,\kappa)$ is defined in \eqref{T(z)}.
\end{pro}
\begin{proof}
The trajectory of the ODE is located on the level set
$$
\mathcal{C}_{\sigma}(\varepsilon):=\{(\mathtt{x}_2,\mathtt{x}_{1}),\, \mathcal{H}(\varepsilon,\mathtt{x}_{1},\mathtt{x}_2)=\mathcal{H}(\varepsilon,0,\lambda)=-\sigma   \}.
$$
Using the polar  coordinates
	$$
	\mathtt{x}_{1}=\sqrt{\mathtt{I}(\theta)}\cos( \theta)\quad\hbox{and}\quad \mathtt{x}_2=\sqrt{\mathtt{I}(\theta)}\sin( \theta), \theta\in \R,
	$$
    we get
    $$
    \ln(I)+\tfrac{1}{8\kappa} I\cos^2(\theta)-2r_\varepsilon^2 F\big(\sqrt{\mathtt{I}(\theta)}\cos( \theta),\sqrt{\mathtt{I}(\theta)}\sin( \theta)\big)=-2\mathcal{H}(\varepsilon,0,\lambda).
    $$
    Observe that, from \eqref{F-est1} we have the estimate
    \begin{align}\label{F-est2}
|F\big(\sqrt{\mathtt{I}(\theta)}\cos( \theta),\sqrt{\mathtt{I}(\theta)}\sin( \theta)\big)|\lesssim I|\ln(I)|.
\end{align}
    Thus the function $\theta\in\mathbb{R}\mapsto I(\theta)$ is the solution to the  nonlinear functional equation
    $$
    \mathcal{F}(\varepsilon, I)(\theta):=I(\theta)e^{\tfrac{1}{8\kappa} I(\theta)\cos^2(\theta)} e^{-2r_\varepsilon^2 F\big(\sqrt{\mathtt{I}(\theta)}\cos( \theta),\sqrt{\mathtt{I}(\theta)}\sin( \theta)\big)}e^{2\mathcal{H}(\varepsilon,0,\lambda)}-1=0.
    $$
    We shall introduce the auxiliary problem
    $$
    \mathcal{F}_0(\varepsilon,\sigma, I)(\theta):=I(\theta)e^{\tfrac{1}{8\kappa} I(\theta)\cos^2(\theta)} e^{-2r_\varepsilon^2 F\big(\sqrt{\mathtt{I}(\theta)}\cos( \theta),\sqrt{\mathtt{I}(\theta)}\sin( \theta)\big)}e^{2\sigma}-1=0.
    $$
    Thus
    $$
    \mathcal{F}(\varepsilon, I)= \mathcal{F}_0(\varepsilon,\mathcal{H}(\varepsilon,0,\lambda), I).$$
    Constructing $2\pi-$periodic solution to \eqref{system004} amounts to  show that the previous equation admits a periodic solution $\theta\in\R\mapsto I(\theta).$  Due to the singular behavior of $r_\varepsilon$ with respect to $\varepsilon$ we shall instead use the variable
    $
    \epsilon=r_\varepsilon$ and introduce the rescaled  functional
    $$
\mathcal{F}_1(\epsilon,\sigma,I):=\mathcal{F}_0(\varepsilon,\sigma,I).
    $$
     For $\varepsilon=0,$ or $\epsilon=0,$ we have constructed a closed simple orbit $\theta\in\R \mapsto I_0(\theta)$ from the previous section and we want to implement the implicit function theorem to construct similar structures for small $\varepsilon,$ or small $\epsilon.$ Introduce the space $\mathcal{X}=\mathcal{C}(\mathbb{T},\R)$ equipped with the uniform norm $\|\cdot\|_{\infty}$ and consider the  open set 
    $$
    \mathcal{U}=\big\{ I\in \mathcal{X}, I>0\}, \quad \mathcal{J}=(-1,1), \quad \sigma_0=\mathcal{H}(0,0,\lambda).
    $$
    Then $\mathcal{F}_1:\mathcal{J}\times (a,b)\times \mathcal{U}\to \mathcal{X}$ is of class $\mathcal{C}^1$ and 
    \begin{align*}
    \partial_I \mathcal{F}_1(0,\sigma_0, I_0)h&=h(\theta)e^{\tfrac{1}{8\kappa} I_0(\theta)\cos^2(\theta)} e^{2\mathcal{H}(0,0,\lambda)}\Big(1+\tfrac{1}{8\kappa} I_0(\theta) \cos^2(\theta)  \Big).
    \end{align*}
As $I_0>0$ then $\partial_I \mathcal{F}(0,\sigma, I_0):\mathcal{X}\to \mathcal{X}$ is an isomorphism. Therefore, using the implicit function theorem we deduce the existence of a family of function $(\epsilon,\sigma)\in (-\epsilon_0,\epsilon_0)\times(\sigma_0-\epsilon_0,\sigma_0+\epsilon_0)\mapsto I(\epsilon,\sigma)\in \mathcal{X}$ that satisfies
$$
\mathcal{F}_1(\epsilon,\sigma, I(\epsilon,\sigma))=0,
$$
with
$$\|I(\epsilon,\sigma)-I(0,\sigma_0)\|_{L^\infty}+\|\partial_\sigma I(\epsilon,\sigma)-\partial_\sigma  I(0,\sigma_0)\|_{L^\infty}\lesssim |\epsilon|+|\sigma-\sigma_0|,
$$
which yields in the original variable a family of solutions
$\varepsilon\in(0,\varepsilon_0)\mapsto I_\varepsilon $ to the original problem
$$
\mathcal{F}_0(\varepsilon,\sigma, I(\varepsilon,\sigma))=0,
$$
 with
\begin{align}\label{varia-act}
\|I(\varepsilon,\sigma)-I(0,\sigma_0)\|_{L^\infty}+\|\partial_\sigma I(\varepsilon,\sigma)-\partial_\sigma  I(0,\sigma_0)\|_{L^\infty}\lesssim |r_\varepsilon|+|\sigma-\sigma_0|.
\end{align}
The area enclosed by the periodic orbit $\mathcal{C}_\sigma$ of energy $-\sigma$ is
\begin{align*}
\mathcal{A}_\sigma(\varepsilon)&=\frac{1}{2}\int_0^{2\pi}I(\varepsilon,\sigma)(\theta) d\theta,\quad r_\varepsilon=\epsilon.
\end{align*}
It follows that from the classical formula on the period
\begin{align*}
\mathtt{T}(\varepsilon,\sigma)&=\partial_\sigma \mathcal{A}_\sigma(\varepsilon)=\frac{1}{2}\int_0^{2\pi}\partial_\sigma I(\varepsilon,\sigma)(\theta) d\theta\\
&=\frac{1}{2}\int_0^{2\pi}\partial_\sigma I(0,\sigma_0)(\theta) d\theta+\frac{1}{2}\int_0^{2\pi}(\partial_\sigma I(\varepsilon,\sigma)-\partial_\sigma I(0,\sigma_0))(\theta) d\theta. 
\end{align*}
Hence
\begin{align*}
\mathtt{T}(\varepsilon,\sigma)&=\mathtt{T}(0,\sigma_0)+\frac{1}{2}\int_0^{2\pi}(\partial_\sigma I(\varepsilon,\sigma)-\partial_\sigma I(0,\sigma_0))(\theta) d\theta.
\end{align*}
Using \eqref{varia-act} we infer
\begin{align*}
|\mathtt{T}(\varepsilon,\sigma)-\mathtt{T}(0,\sigma_0)|\lesssim r_{\varepsilon}+|\sigma-\sigma_0|.
\end{align*}
Notice that the particular values 
\begin{align}\label{link-lambda-x}
\sigma=\frac12\,\ln(\lambda^2)+\frac{\lambda^2}{16\kappa}+r_{\varepsilon}^2 F(0,\lambda),\quad \sigma_0=\frac12\,\ln(\lambda^2)+\frac{\lambda^2}{16\kappa},
\end{align}
 correspond to the orbits  $\mathcal{C}_{\sigma}(\varepsilon)$ and $\mathcal{C}_{\sigma_0}(0)$, successively. Thus, with this choice we get in view of \eqref{F-est2}
 \begin{align*}
|\mathtt{T}(\varepsilon,\sigma)-\mathtt{T}(0,\sigma_0)|\lesssim r_{\varepsilon}+\lambda^2|\ln (\lambda)|r_{\varepsilon}^2.
\end{align*}
By doing the derivative with respect $\sigma$ in the period function one gets
$$
\partial_\sigma \mathtt{T}(\varepsilon,\sigma)=\partial_\sigma \mathtt{T}(0,\sigma_0)+\frac12 \int_0^{2\pi}(\partial^2_\sigma I(\varepsilon,\sigma)-\partial^2_\sigma I(0,\sigma_0))(\theta) d\theta.
$$
Recall that   $\sigma=-\mathcal{H}(\varepsilon,\lambda,0)$. From this, we have
$$\mathtt{T}(\varepsilon,-\mathcal{H}(\varepsilon,\lambda,0))=T(\varepsilon,\lambda, \kappa).
$$
and $\mathtt{T}(0,\sigma)=T_0(\lambda, \kappa)$. By similar analysis, we get the asymptotics of $\partial_\sigma \mathtt{T}(\varepsilon,\sigma).$
    \end{proof}
As a consequence, we get the following result.
\begin{pro}\label{prop-ring}
Let $\kappa>0$, $0<a<b$,  $\lambda\in [a,b]$ and $P_j(0)=(p_{j,1}(0),p_{j,2}(0))\in \R^2$, $j=1,2$, two points such that
$$
p_{1,1}(0)=\kappa+{\tfrac12}r_\varepsilon (2\kappa)^{\frac14}\lambda,\quad  p_{2,1}(0)=\kappa-{\tfrac12}r_\varepsilon (2\kappa)^{\frac14}\lambda \quad p_{1,2}(0)-p_{2,2}(0)=0, \quad r_\varepsilon:=(2\kappa)^\frac14|\ln\varepsilon|^{-\frac12}.
$$
Then, there exists $\varepsilon_0>0$ small enough such that for any $\varepsilon\in (0,\varepsilon_0),$ the system \eqref{eq-points-general} admits a unique  solution
$\tau\in\R\mapsto P_1(\tau),P_2(\tau)$ such that the functions $p_{1,1}, p_{2,1}$ and $p_{1,2}-p_{2,2}$ are periodic with a period $T=T(\varepsilon,\lambda)$  such  that
\begin{enumerate}
    \item Asymptotic expansion
$$
p_{1,1}(\tau)=\kappa+{\tfrac12}r_\varepsilon (2\kappa)^{\frac14}\mathtt{x}_{1}, \quad
p_{2,1}(\tau)=\kappa-{\tfrac12}r_\varepsilon (2\kappa)^{\frac14}\mathtt{x}_{1},\quad p_{1,2}-p_{2,2}= r_\varepsilon (2\kappa)^{-\frac14} \mathtt{x}_2,
$$
where $\mathtt{x}_{1},\mathtt{x}_2$ are periodic solutions to the system \eqref{system004}.
\item The period $T(\varepsilon,\lambda,\kappa)$ takes the form
$$
   T(\varepsilon,\lambda,\kappa)=T_0(\lambda,\kappa)+O\left(\tfrac{1}{\sqrt{|\ln\varepsilon|}}\right),\quad\hbox{and}\quad
\partial_\lambda T(\varepsilon,\lambda,\kappa)=T_0'(\lambda,\kappa)+O\left(\tfrac{1}{\sqrt{|\ln\varepsilon|}}\right).
   $$
   \item Orbital symmetries: we have
\begin{equation*}
p_{1,1}\big(\tau+\tfrac{T}{2}\big)=p_{2,1}(\tau)\quad\hbox{and} \quad (p_{1,2}-p_{2,2})\big(\tau+\tfrac{T}{2}\big)=-(p_{1,2}-p_{2,2})(\tau)=(p_{1,2}-p_{2,2})(-\tau),
\end{equation*}
   and
   \begin{align*}
   p_{1,1}(\tau)&=p_{1,1}(T-\tau)=p_{1,1}(-\tau),\quad 
\dot{P}_1(\tau+\tfrac{T}{2})=\dot{P}_2(\tau),\\
 \dot{p}_{1,2}(\tau)&=\dot{p}_{1,2}(T-\tau)=\dot{p}_{1,2}(-\tau).
\end{align*} 
\end{enumerate}
\end{pro}
We conclude this section with the following corollary, which describes the dynamics of the center of mass. In the case of a single vortex filament, this expression coincides with the classical singular velocity obtained in the standard works on vortex filaments. For two interacting vortex filaments, previous studies only require the first term in the expansion of the center of mass. For completeness, we provide a better asymptotic expansion in the following corollary.

\begin{cor}\label{coro-speed}
Let $\kappa>0$, $0<a<b$,  $\lambda\in [a,b]$ and
$$
\mathtt{C}(\tau)=(\mathtt{C}_1(\tau),\mathtt{C}_2(\tau))
:=\tfrac12\big(P_1(\tau)+P_2(\tau)\big),\quad  P_j(0)=\big(\kappa+{\tfrac12}(-1)^{j+1}|\ln\varepsilon|^{-\frac12} (2\kappa)^{\frac12}\lambda,0\big). 
$$
Then
$$
{\mathtt{C}(\tau)
=
\bigl(\kappa,\,U_\varepsilon\,\tau+W_\varepsilon(\tau)\bigr),}
$$
where the constant drift speed is given by
\begin{equation}\label{defUeps}
U_\varepsilon
:=\frac{|\ln\varepsilon|^{-1}}{2T(\varepsilon,\lambda)}
\int_{0}^{T(\varepsilon,\lambda)}
\Bigl(\tfrac{1}{\sqrt{2}}\partial_{p_{1,1}}G
+\tfrac{1}{\sqrt{2}}\partial_{p_{2,1}}G
+\partial_{p_{1,1}}\mathtt{G}
+\partial_{p_{2,1}}\mathtt{G}\Bigr)(P_1,P_2)(s)\,ds,
\end{equation}
and $W_\varepsilon$ is an odd $T(\varepsilon,\lambda)$-periodic function.
Moreover,
\begin{align*}
U_\varepsilon
&=\tfrac{1}{4\sqrt{2\kappa}}
\biggl[1+\tfrac12|\ln\varepsilon|^{-1}\Bigl|\ln\Bigl(\tfrac{\sqrt{2\kappa}}{|\ln\varepsilon|}\Bigr)\Bigr|
+|\ln\varepsilon|^{-1}\Bigl(\tfrac52\ln8+\tfrac32\ln(\kappa)-\tfrac54\Bigr)\biggr]\\
&\quad+\tfrac{|\ln\varepsilon|^{-1}}{4T(\varepsilon,\lambda)\sqrt{2\kappa}}
\int_{0}^{T(\varepsilon,\lambda)}
\left(\tfrac{\mathtt{x}_1^2}{\mathtt{x}_1^2+\mathtt{x}_2^2}
-\tfrac12\ln\bigl(\mathtt{x}_1^2+\mathtt{x}_2^2\bigr)\right)(s)\,ds
+O\!\left(\tfrac{|\ln(|\ln\varepsilon|)|}{|\ln\varepsilon|^2}\right),
\end{align*}
and
$$
W_\varepsilon=O\!\left(|\ln\varepsilon|^{-1}\right).
$$
\end{cor}

\begin{proof} 
Using \eqref{decompo-Hamil} and the definition of $\mathtt{C}$, we obtain
$$
2|\ln\varepsilon|\,\dot{\mathtt{C}}(\tau)
=\nabla_{P_1}^{\perp}H(P_1,P_2)+\nabla_{P_2}^{\perp}H(P_1,P_2),
$$
where 
$$
H(P_1,P_2)=\tfrac{1}{\sqrt{2}}\,G(P_1,P_2)+\mathtt{G}(P_1,P_2).
$$
Since $\mathtt{G}$ is independent of the second coordinates (see \eqref{decompo-Hamil}),
$$
\partial_{p_{1,2}}\mathtt{G}(P_1,P_2)=\partial_{p_{2,2}}\mathtt{G}(P_1,P_2)=0.
$$
Moreover, by the structural property of the interaction kernel $G$ (cf. \eqref{eq:G}),
it depends on $p_{1,2},p_{2,2}$ only through their difference, hence
$$
\partial_{p_{1,2}}G(P_1,P_2)+\partial_{p_{2,2}}G(P_1,P_2)=0.
$$
Therefore $\dot{\mathtt{C}}_1(\tau)=0$, and
$$
\mathtt{C}_1(\tau)=\tfrac12\bigl(p_{1,1}(0)+p_{2,1}(0)\bigr)=\kappa.
$$
For the second component, we write
\begin{equation}\label{eqC2}
2|\ln\varepsilon | \dot{\mathtt{C}}_2(\tau)
=\tfrac{1}{\sqrt{2}}\big(\partial_{p_{1,1}}G
+ \partial_{p_{2,1}}G\big)(P_1,P_2)
+\big(\partial_{p_{1,1}}\mathtt{G}
+ \partial_{p_{2,1}}\mathtt{G}\big)(P_1,P_2).
\end{equation}
Integrating \eqref{eqC2} and using \eqref{defUeps} yields
$$
 {\mathtt{C}}_2(\tau)
={{\mathtt{C}}_2(0)}+ U_\varepsilon\, \tau + W_\varepsilon(\tau),
$$
where 
$$
W_\varepsilon (\tau)=|\ln\varepsilon|^{-1}\int_{0}^\tau\Big\{\Big(\tfrac{1}{2\sqrt{2}}\partial_{p_{1,1}}G
+ \tfrac{1}{2\sqrt{2}}\partial_{p_{2,1}}G+\tfrac12\partial_{p_{1,1}}\mathtt{G}
+ \tfrac12\partial_{p_{2,1}}\mathtt{G}\Big)(P_1,P_2)(s)-U_\varepsilon \Big\}ds.
$$
By construction, $W_\varepsilon$ is $T(\varepsilon,\lambda)$-periodic, and it is odd thanks to the time-symmetry of the orbit $(P_1,P_2)$ established earlier. 
To obtain the expansion of $U_\varepsilon$, we use \eqref{decompo-Hamil} and
\eqref{scaling-points}. First,
\begin{align*}
\bigl(\partial_{p_{1,1}}\mathtt{G}+\partial_{p_{2,1}}\mathtt{G}\bigr)(P_1,P_2)
&= \tfrac{1}{4}\sum_{m=1}^2(2 p_{m,1})^{-\frac12}
\biggl[|\ln\varepsilon|
+\tfrac{1}{4}\bigl(5\ln8+3\ln(p_{m,1})-1\bigr)\biggr]\\
&=\tfrac{1}{2}(2 \kappa)^{-\frac12}
\biggl[|\ln\varepsilon|
+\tfrac{1}{4}\bigl(5\ln8+3\ln(\kappa)-1\bigr)\biggr]
+O(|\ln\varepsilon|^{-1}).
\end{align*}
Next, applying Lemma~\ref{pro-decomp2-nabla} with
$$
z_1=\kappa,\quad z_2=0,\quad
x_1=\tfrac12 \mathtt{x}_1,\quad
y_1=-\tfrac12\mathtt{x}_1,\quad
x_2-y_2=\mathtt{x}_2,\quad
\epsilon =(2\kappa)^{\frac14}|\ln\varepsilon|^{-\frac12},
$$
we obtain
\begin{align*}
& |\ln\varepsilon|^{-1}\bigl(\partial_{p_{1,1}}G+\partial_{p_{2,1}}G\bigr)(P_1,P_2)
=\tfrac{1}{4 \sqrt{\kappa}}
\Bigl|\ln\Bigl(\tfrac{\sqrt{2\kappa}}{|\ln\varepsilon|}\Bigr)\Bigr|
\,|\ln\varepsilon|^{-1}\\
&\quad+|\ln\varepsilon|^{-1} \tfrac{1}{2\sqrt{\kappa}}
\biggl[
\tfrac54\ln 8 + \tfrac34 \ln(\kappa) - 1
+ \tfrac{\mathtt{x}_1^2}{\mathtt{x}_1^2+\mathtt{x}_2^2}
- \tfrac12 \ln\bigl(\mathtt{x}_1^2+\mathtt{x}_2^2\bigr)
\biggr]+O\!\left(\tfrac{|\ln(|\ln\varepsilon|)|}{|\ln\varepsilon |^2}\right).
\end{align*}
Substituting these expansions into \eqref{defUeps} gives the claimed asymptotic formula for
$U_\varepsilon$, as well as $W_\varepsilon=O(|\ln\varepsilon|^{-1})$.
\end{proof}
\section{Symmetry reduction and time scaling in the patch pairs motion}\label{Sym-Redu}
The analysis of leapfrogging motions in the 3D axisymmetric Euler equations relies crucially on exploiting suitable symmetry invariances. Indeed, the configuration of two interacting vortex patches possesses a natural invariance under the exchange of the two vortices, which should be faithfully reflected in the governing equations. This symmetry allows us to recast the pairwise interaction in terms of a single  equation. In this section, we introduce the appropriate symmetry reduction, showing that the temporal and spatial periodicities imposed on the perturbations of the two patches lead to an equivalence between their governing operators. This reduction serves as the first essential step towards obtaining a well-posed formulation of the leapfrogging motion, and it prepares the ground for the subsequent time rescaling that captures the relevant dynamical regime.

\subsection{Symmetry reduction}\label{section:Symmetry reduction}

In this section, we exploit the symmetry between the point vortices $P_1$ and $P_2$ in order to derive an analogous symmetry relation between ${\bf G}_1$ and ${\bf G}_2$. By exploiting this property we can reduce the motion of two rings to only one PDE containing a delay term.   Recalling that the dynamics of the pair of patches are governed by \eqref{eq-gamma-tilde}, we arrive at the following key result.
\begin{lem}\label{lem-red}
Let $\tau\mapsto \mathbf{V}_k(\tau)\in\mathbb{C}, k=1,2,$ be two functions such that
$$
{\dot{\mathbf{V}}_k(\tau+T)=\dot{\mathbf{V}}_k(\tau),\quad \dot{\mathbf{V}}_2(\tau)=\dot{\mathbf{V}}_1\big(\tau+\tfrac{T}{2}\big)\quad\hbox{and}\quad ({\mathbf{V}}_1-{\mathbf{V}}_2)(\tau+\tfrac{T}{2})=({\mathbf{V}}_2-{\mathbf{V}}_1)(\tau)},
$$
for each $k=1,2.$ Assume that a smooth  function $f_1:\R\times\T\to \R$ is $T$-periodic in time and define  $$
\forall \tau\in\mathbb{R},\, \forall\theta\in\R,\quad f_2(\tau,\theta):=f_1\big(\tau+\tfrac{T}{2},\theta\big).
$$
Then 
 we have the identity
$$
{\bf G}_1(\varepsilon,f_1,f_2)\big(\tau+\tfrac{T}{2},\theta\big)={\bf G}_2(\varepsilon,f_1,f_2)(\tau,\theta),\quad \forall \tau,\theta\in\mathbb{R}.
$$
As a consequence, the two equations in \eqref{eq-gamma-tilde} are equivalent and  the system reduces to 
\begin{align*}
\forall \tau\in\mathbb{R},\, \forall\theta\in\R,\quad{\bf G}_1\Big(\varepsilon,f_1,f_1\big(\cdot+\tfrac{T}{2}\big)\Big)\big(\tau,\theta\big)=0.
\end{align*}   
\end{lem}

\begin{proof}
Invoking \eqref{eq-gamma-tilde} and \eqref{TH-1}, and suppressing the explicit
dependence on the spatial variable~$\theta$ in the notation, we obtain for all $\tau$,
\begin{align*}
&\mathbf{G}_1\big(\varepsilon, f_1, f_2\big)(\tau+\tfrac{T}{2})= \varepsilon^3 |\ln\varepsilon |  \partial_\tau f_1(\tau+\tfrac{T}{2})+\partial_\theta 
\big[ \Psi_1(\varepsilon,f_1,f_1)(\tau+\tfrac{T}{2})+\Psi_2(\varepsilon,f_2,f_1)(\tau+\tfrac{T}{2}) \big]\\& - \big(|\ln\varepsilon | \dot{P}_1+i\,\varepsilon\,\dot{\mathbf{V}}_1)(\tau+\tfrac{T}{2}\big)\cdot i\partial_\theta \gamma_1(\tau+\tfrac{T}{2}) +\varepsilon^2 |\ln\varepsilon | \tfrac{\dot{p}_{1,1}(\tau+\frac{T}{2})}{8p_{1,1}(\tau+\frac{T}{2})} \partial_\theta\left( w_1^2(\tau+\tfrac{T}{2})\sin(2\theta)\right).
\end{align*}
According to Proposition \ref{prop-ring} we have
$$
\dot{P}_1\big(\tau+\tfrac{T}{2}\big)=\dot{P}_2(\tau), \quad\hbox{and}\quad p_{1,1}\big(\tau+\tfrac{T}{2}\big)=p_{2,1}(\tau).
$$ 
Moreover, using the fact that $f_2(\tau,\theta)=f_1\big(\tau+\frac{T}{2},\theta\big)$ together with
$
\partial_\theta \gamma_1\big(\tau+\tfrac{T}{2},\theta\big)=\partial_\theta\gamma_2\big(\tau,\theta\big)
$ we infer
\begin{align}\label{GG-1}
\nonumber &\mathbf{G}_1(\varepsilon, f_1, f_2) \big(\tau+\tfrac{T}{2}\big)= \varepsilon^3 |\ln\varepsilon |  \partial_\tau f_2(\tau)+\partial_\theta\big[ \Psi_1(\varepsilon,f_1,f_1)\big(\tau+\tfrac{T}{2}\big)+\Psi_2(\varepsilon,f_2,f_1)\big(\tau+\tfrac{T}{2}\big) \big]\\&\qquad  -\big( |\ln\varepsilon | \dot{P}_2(\tau)+i\,\varepsilon\,\dot{\mathbf{V}}_2(\tau)\big)\cdot i\partial_\theta \gamma_2(\tau) +\varepsilon^2 |\ln\varepsilon | \tfrac{\dot{p}_{2,1}(\tau)}{8p_{2,1}(\tau)} \partial_\theta\left( w_1^2(\tau)\sin(2\theta)\right),
\end{align}
where we have used the assumption
$
 \dot{\mathbf{V}}_1\big(\tau+\tfrac{T}{2}\big)=\dot{\mathbf{V}}_2(\tau).
$
On the other hand, from \eqref{TH-1} and using the structure of the kernel $G$, given by \eqref{eq:G}, we find
\begin{align*}
\Psi_j(\varepsilon,f_j,f_j)(\tau,\theta)
&= \frac{1}{\pi\sqrt{2}}\int_0^{2\pi}\!\int_0^{1}
   \mathcal{K}_j\big(\tau,\theta,\rho,\eta\big)\rho\, d\rho\, d\eta.
\end{align*}
with
$$
\mathcal{K}_j(\tau,\theta,\rho,\eta):=G\Big(p_{j,1}(\tau)+\varepsilon  w_j(\tau,\theta)\,\mathcal{Z}_j(\tau,\theta), p_{j,1}(\tau)+\varepsilon \rho w_j(\tau,\eta)\,\mathcal{Z}_j(\tau,\eta)\Big)
   w_j^2(\tau,\eta).
$$
It follows that
\begin{align*}
\Psi_1(\varepsilon,f_1,f_1)(\tau+\tfrac{T}{2},\theta)
=\frac{1}{\pi\sqrt{2}}\int_0^{2\pi}\!\int_0^{1}
   \mathcal{K}_1\big(\tau+\tfrac{T}{2},\theta,\rho,\eta\big)\rho\, d\rho\, d\eta.
   \end{align*}
   Applying  Proposition \ref{prop-ring}-(3) gives
   $$
   \mathcal{K}_1\big(\tau+\tfrac{T}{2},\theta,\rho,\eta\big)=\mathcal{K}_2\big(\tau,\theta,\rho,\eta\big)$$
   and therefore
   \begin{equation}\label{psi12t2}
    \Psi_1(\varepsilon,f_1,f_1)\big(\tau+\tfrac{T}{2},\theta\big)=\Psi_2(\varepsilon,f_2,f_2)(\tau,\theta).
\end{equation}
Similarly, from \eqref{TH-1} and  \eqref{eq:G}, we have for $j\neq k$
\begin{align*}
& \Psi_k(\varepsilon,f_k,f_j)(\tau,\theta)=\\
&\frac{1}{\pi\sqrt{2}}\int_0^{2\pi}\!\int_0^{1}
   G\Big(\big(p_{j,1}+i\,\zeta_{jk}\big)(\tau)+\varepsilon( w_j\mathcal{Z}_j)(\tau,\theta), p_{k,1}(\tau)+\varepsilon \rho (w_k\mathcal{Z}_k)(\tau,\eta)\Big)
   w_k^2(\tau,\eta)\rho d\rho d\eta,
\end{align*}
with 
$$
\zeta_{jk}(\tau):=\big(p_{j,2}-p_{k,2}\big)(\tau)+{\varepsilon}{|\ln\varepsilon|^{-1}} \big(\mathbf{V}_j-\mathbf{V}_k\big)(\tau).
$$
By assumptions and using  Proposition \ref{prop-ring}-(3), we have
$$
\zeta_{12}\big(\tau+\tfrac{T}{2},\theta\big)=\zeta_{21}(\tau),\quad\textnormal{and}\quad  p_{2,1}\big(\tau+\tfrac{T}{2}\big)=p_{1,1}(\tau).
$$
It follows that
\begin{align*}
  & G\Big(\big(p_{1,1}+i\zeta_{12}\big)(\tau+\tfrac{T}{2})+\varepsilon( w_1\mathcal{Z}_1)(\tau+\tfrac{T}{2},\theta), p_{2,1}(\tau+\tfrac{T}{2})+\varepsilon \rho (w_2\mathcal{Z}_2)(\tau+\tfrac{T}{2},\eta)\Big)
   \\ 
   &\quad =
   G\Big(\big(p_{2,1}+i\zeta_{21}\big)(\tau)+\varepsilon( w_2\mathcal{Z}_2)(\tau,\theta), p_{1,1}(\tau)+\varepsilon \rho (w_1\mathcal{Z}_1)(\tau,\eta)\Big),
\end{align*}
implying the identity
\begin{align}\label{psi12t22}
& \Psi_2(\varepsilon,f_2,f_1)\big(\tau+\tfrac{T}{2},\theta\big)
= \Psi_1(\varepsilon,f_1,f_2)(\tau,\theta).
\end{align}
Consequently, 
$$
\Psi_1(\varepsilon,f_1,f_1)\big(\tau+\tfrac{T}{2},\theta\big)+\Psi_2(\varepsilon,f_2,f_1)\big(\tau+\tfrac{T}{2},\theta\big)=\Psi_1\big(\varepsilon,f_1,f_2\big)(\tau,\theta)+\Psi_2(\varepsilon,f_2,f_2)(\tau,\theta).
$$
Plugging \eqref{psi12t2} and  \eqref{psi12t22} into \eqref{GG-1} leads to
$$\mathbf{G}_1(\varepsilon, f_1, f_2) (\tau+\tfrac{T}{2},\theta)=\mathbf{G}_2(\varepsilon, f_1, f_2) (\tau,\theta).
$$
This achieves the proof of the desired result.
\end{proof}
As a consequence of Lemma~\ref{lem-red}, the  system \eqref{eq-gamma-tilde} is equivalent to
the single delayed scalar equation
\begin{equation}\label{eq-G1-red}
\mathbf{G}_1\Big(\varepsilon,f_1,f_1(\cdot+\tfrac{T}{2})\Big)(\tau,\theta)=0,
\qquad (\tau,\theta)\in\R\times\T,
\end{equation}
where (recall \eqref{eq-gamma-tilde})
\begin{align*}
\mathbf{G}_1(\varepsilon,f_1,f_2)(\tau,\theta)
&=\varepsilon^3|\ln\varepsilon|\,\partial_\tau f_1(\tau,\theta)
+\partial_\theta\Big(\Psi(\gamma_1(\tau,\theta))\Big)
-\Big(|\ln\varepsilon|\dot P_1(\tau)+i\varepsilon\dot{\mathbf{V}}_1(\tau)\Big)\cdot
i\,\partial_\theta\gamma_1(\tau,\theta)\\
&\quad+\varepsilon^2|\ln\varepsilon|\,
\frac{\dot p_{1,1}(\tau)}{8p_{1,1}(\tau)}\,
\partial_\theta\!\Big(w_1^2(\tau,\theta)\sin(2\theta)\Big),
\end{align*}
and
$$
\gamma_1(\tau,\theta)=P_1(\tau)+i\,\frac{\varepsilon}{|\ln\varepsilon|}\mathbf{V}_1(\tau)
+\varepsilon\,w_1(\tau,\theta)\,\mathcal{Z}_1(\tau,\theta),
\qquad
w_1(\tau,\theta)=\sqrt{1+2\varepsilon f_1(\tau,\theta)},
$$
$$
\mathcal{Z}_1(\tau,\theta)
=(2p_{1,1}(\tau))^{\frac14}\cos\theta+i\,(2p_{1,1}(\tau))^{-\frac14}\sin\theta.
$$

Finally, recalling \eqref{Pol-ah1}--\eqref{TH-1}, we write
\begin{equation}\label{Pol-ah2}
\Psi(\gamma_1(\tau,\theta))=\Psi_1(\varepsilon,f_1)(\tau,\theta)+\Psi_2(\varepsilon,f_1)(\tau,\theta),
\end{equation}
where $\Psi_1(\varepsilon,f_1):=\Psi_1(\varepsilon,f_1,f_1)$ is the self-interaction term
\begin{align*}
\sqrt{2}\pi\,\Psi_1(\varepsilon,f_1)(\tau,\theta)
:=\int_0^{2\pi}\!\int_0^{1}
G\Big(P_1(\tau)+\varepsilon (w_1\mathcal{Z}_1)(\tau,\theta),\,
      P_1(\tau)+\varepsilon\rho (w_1\mathcal{Z}_1)(\tau,\eta)\Big)
w_1^2(\tau,\eta)\rho d\rho d\eta,
\end{align*}
and $\Psi_2(\varepsilon,f_1):=\Psi_2(\varepsilon,f_2,f_1)$ is the interaction term with
$f_2(\tau,\cdot)=f_1(\tau+\tfrac{T}{2},\cdot)$ (hence $w_2(\tau,\cdot)=w_1(\tau+\tfrac{T}{2},\cdot)$),
\begin{align*}
\sqrt{2}\pi\,\Psi_2(\varepsilon,f_1)(\tau,\theta)
:=\int_0^{2\pi}\!\int_0^{1}
G\Big(&\big(P_1+i\tfrac{\varepsilon}{|\ln\varepsilon|}\mathbf{V}_1\big)(\tau)+\varepsilon (w_1\mathcal{Z}_1)(\tau,\theta),\\ &\quad
      \big(P_2+i\tfrac{\varepsilon}{|\ln\varepsilon|}\mathbf{V}_2\big)(\tau)+\varepsilon\rho (w_2\mathcal{Z}_2)(\tau,\eta)\Big)
w_2^2(\tau,\eta)\rho d\rho d\eta.
\end{align*}
\subsection{Reversibility property}
In the search for periodic orbits within Hamiltonian systems, symmetry structures play a decisive role. Among these, reversibility is particularly significant, as it provides a natural mechanism by which time-periodic dynamics may arise. Roughly speaking, reversibility encodes the invariance of the underlying vector field under a transformation that combines a spatial reflection with a time reversal. This structure forces the trajectories of the system to be symmetric with respect to specific { reversal times}, and precisely this symmetry allows periodic motions to be generated.
We now make this precise by examining the action of the symmetry operator $\mathscr{S}$ define by 
\begin{equation}\label{def-S}
(\mathscr{S} u)(\tau,\theta)\;:=\;{u(-\tau,-\theta)},
\end{equation}
on the reduced vector field. The goal is to show that the vector field associated with 
$$
f_1\mapsto \mathbf{G}_1\big(\varepsilon, f_1, f_2\big),
$$
is reversible in the sense of Definition~\ref{Def-Rev}. More precisely, we have the following result.
\begin{pro}\label{prop-reversibility}
Let   ${\mathbf{V}}_1$ and ${\mathbf{V}}_2$ be as in Lemma $\ref{lem-red},$ and assume in addition that both functions are  odd. Then, we have 
$$
\forall \tau,\forall \theta\in\R,\quad \mathbf{G}_1\big(\varepsilon, (\mathscr{S}{f_1}), (\mathscr{S}{f_2}\big)(\tau,\theta)=-\mathbf{G}_1\big(\varepsilon, f_1, f_2\big)(-\tau,-\theta).
$$
\end{pro}

\begin{proof}

By Definition~\ref{Def-Rev}, the vector field $f_1\mapsto \mathbf{G}_1(\varepsilon,f_1,f_1(\cdot+\tfrac{T}{2}))$
is $\mathscr{S}$–reversible if
\begin{equation}\label{eq:RevGoal}
\mathbf{G}_1\!\left(\varepsilon,\,\mathscr{S} f_1,\,(\mathscr{S} f_2)\right)(\tau,\theta)
\;=\;-\,\mathbf{G}_1\!\left(\varepsilon,\,f_1,\,f_2\right)(-\tau,-\theta).
\end{equation}
Differentiation and $\mathscr{S}$ obey to the anti-commutation relation
\begin{equation*}
\partial_\tau(\mathscr{S} u)=-\mathscr{S}(\partial_\tau u),\qquad
\partial_\theta(\mathscr{S} u)=-\mathscr{S}(\partial_\theta u).
\end{equation*}
By Proposition~\ref{prop-ring}, $\tau\mapsto p_{1,1}(\tau)$ is even, hence
\begin{equation}\label{eq:p11-sym}
p_{1,1}(-\tau)=p_{1,1}(\tau),\qquad \dot p_{1,1}(-\tau)=-\dot p_{1,1}(\tau).
\end{equation}
Let $\gamma_{f_1}$ denote the boundary parametrization, instead of $\gamma_1$, which is described through \eqref{parametrization},
$$
\gamma_{f_1}(\tau,\theta)
= P_1(\tau)+ \mathbf{V}_1(\tau)+\varepsilon\,w_{f_1}(\tau,\theta)\,\mathcal Z_1(\tau,\theta).
$$
From \eqref{parametrization}, Proposition~\ref{prop-ring} and by assumptions we have
\begin{equation}\label{eq:geom-sym}
\mathcal Z_1(-\tau,-\theta)=\overline{\mathcal Z_1(\tau,\theta)},\qquad
\dot P_1(-\tau)=-\overline{\dot P_1(\tau)} ,\qquad
{\dot{\mathbf{V}}_1(-\tau)={\dot{ \mathbf{V}}_1(\tau)}}.
\end{equation}
Moreover, the scalar weight transforms as
\begin{equation}\label{eq:w-sym}
w_{\mathscr{S} f_1}(-\tau,-\theta)=w_{f_1}(\tau,\theta),
\qquad
w_{\mathscr{S} f_1}\!\left(-\tau+\tfrac{T}{2},-\theta\right)=w_{f_1}\!\left(\tau+\tfrac{T}{2},\theta\right).
\end{equation}
Using \eqref{eq:geom-sym}–\eqref{eq:w-sym},
\begin{align}
\partial_\theta\gamma_{\mathscr{S} f_1}(\tau,\theta)
&=\partial_\theta\!\Big(\varepsilon\,w_{\mathscr{S} f_1}(\tau,\theta)\,\mathcal Z_1(\tau,\theta)\Big)\nonumber\\
&=\partial_\theta\!\Big(\varepsilon\,\overline{w_{f_1}(-\tau,-\theta)\,\mathcal Z_1(-\tau,-\theta)}\Big)\nonumber\\
&=-\,\partial_\theta\big(\overline{\gamma_{f_1}}\big)(-\tau,-\theta).
\label{eq:dtheta-gammaS}
\end{align}
Consequently, combining \eqref{eq:geom-sym} and \eqref{eq:dtheta-gammaS},
\begin{align*}
&\Big(\tfrac{d}{d\tau}\big(P_1+i\varepsilon|\ln\varepsilon|^{-1} \mathbf{V}_1\big)\cdot i\,\partial_\theta\gamma_{\mathscr{S} f_1}\Big)(-\tau,-\theta)\\
&\quad =-\,\tfrac{d}{d\tau}\big(P_1+i\varepsilon|\ln\varepsilon|^{-1} \mathbf{V}_1\big)(-\tau)\cdot i\,\partial_\theta\big(\overline{\gamma_{f_1}}\big)(\tau,\theta)\nonumber\\
&\quad =\overline{\tfrac{d}{d\tau}\big(P_1+i\varepsilon|\ln\varepsilon|^{-1}\mathbf{V}_1\big)(\tau)}\cdot i\,\partial_\theta\big(\overline{\gamma_{f_1}}\big)(\tau,\theta)\nonumber\\
&\quad =-\,\Big(\big(\dot{P}_1+i\varepsilon|\ln\varepsilon|^{-1} \dot{\mathbf{V}}_1\big)\cdot i\,\partial_\theta\gamma_{f_1}\Big)(\tau,\theta).
\end{align*}
For the explicit drift term, by 
using \eqref{eq:p11-sym}, \eqref{eq:w-sym}, and $\sin(2(-\theta))=-\sin(2\theta)$, we get
$$
\frac{\dot p_{1,1}(\tau)}{8p_{1,1}(\tau)}\,\partial_\theta\!\Big(w_1^2(-\tau,-\theta)\sin(2\theta)\Big)
=-\left[\frac{\dot p_{1,1}}{8p_{1,1}}\,\partial_\theta\!\big(w_1^2\sin(2\theta)\big)\right](-\tau,-\theta),
$$
which is the required oddness for this contribution.

Let us move to  the stream function $\Psi$ defined in~\eqref{Pol-ah2}. The kernel $G$ involved in $\Psi$ is defined by \eqref{eq:G} and it is invariant under complex conjugation:
$$
G(\overline{z_1},\overline{z_2})=G(z_1,z_2),
$$
and depends on the variables $z_2-z_1.$
With \eqref{eq:geom-sym}–\eqref{eq:w-sym} we obtain, for all $\rho\in[0,1]$,
\begin{align*}
& G\big(P_1(-\tau)+\varepsilon\, (w_{f_1}\mathcal Z_1)(-\tau,-\theta),\,P_1(-\tau)+\varepsilon\rho\, (w_{f_1}\mathcal Z_1)(-\tau,-\eta)\big)\\
&\qquad=G\left(P_1(\tau)+\varepsilon w_{f_1}(\tau,\theta)\,\overline{\mathcal Z_1(\tau,\theta)},\;
P_1(\tau)+\varepsilon\rho\, w_{f_1}(\tau,\eta)\,\overline{\mathcal Z_1(\tau,\eta)}\right)\\
&\qquad=G\left(\overline{P_1(\tau)+\varepsilon w_{f_1}(\tau,\theta)\,\mathcal Z_1(\tau,\theta)},\;
\overline{P_1(\tau)+\varepsilon\rho\, w_{f_1}(\tau,\eta)\,\mathcal Z_1(\tau,\eta)}\right)\\
&\qquad=G\left({P_1(\tau)+\varepsilon w_{f_1}(\tau,\theta)\,\mathcal Z_1(\tau,\theta)},\;
{P_1(\tau)+\varepsilon\rho\, w_{f_1}(\tau,\eta)\,\mathcal Z_1(\tau,\eta)}\right).
\end{align*}
The terms in $\Psi$ with the time shift by $T/2$ are handled similarly thanks to the second identity in \eqref{eq:w-sym} and using the fact that ${\mathbf{V}}_1-{\mathbf{V}}_2$ is odd. Therefore
\begin{equation*}
\Psi\big(\gamma_{\mathscr{S} f_1}\big)(-\tau,-\theta)=\Psi\big(\gamma_{f_1}\big)(\tau,\theta).
\end{equation*}
Putting these transformations together yields exactly \eqref{eq:RevGoal}, that is,
$$
\mathbf{G}_1\!\left(\varepsilon,\,\mathscr{S} f_1,\,(\mathscr{S} f_1)(\cdot+\tfrac{T}{2})\right)(\tau,\theta)
=-\,\mathbf{G}_1\!\left(\varepsilon,\,f_1,\,f_1(\cdot+\tfrac{T}{2})\right)(-\tau,-\theta),
$$
and the vector field is $\mathscr{S}$–reversible. This ends the proof.
\end{proof}

\subsection{Functional value at the circular rings}

In this section, we derive an expansion of the functional $\mathbf{G}_1(\varepsilon,0,0)$ associated with the contour dynamics equations \eqref{eq-gamma-tilde}, evaluated on the circular cross-sections of the vortex rings in the original coordinates. This computation allows us to identify the leading-order contribution in the small parameter $\varepsilon$.
We show that the classical point-vortex system \eqref{eq-points-general} naturally arises as the compatibility condition required to cancel the leading-order term in $\varepsilon$ in the expansion of \eqref{eq-G1-red}. In other words, the point-vortex dynamics governs the dominant interaction between the rings when their core size is small. This step is essential for constructing a sufficiently accurate approximate solution that will serve as the starting point of the Nash–Moser iterative scheme developed in the subsequent analysis.
The main input is an asymptotic expansion of the stream function evaluated on the unperturbed elliptic contours.
\\
Recall from 
\eqref{Pol-ah2} that, for $f_1=f_2=0$, one has $w_1=w_2\equiv 1$ and the
stream function along the first contour splits into the self-induced and interaction parts. We
introduce, for $j\in\{1,2\}$,
\begin{align}\label{Psi0j}
\Psi_{0,j}(\tau,\theta)
&:=\frac{1}{\sqrt{2}}\,
\partial_\theta\Psi_j(\varepsilon,0)(\tau,\theta)\\
&=
\partial_\theta\!\int_0^1\fint_{\T}
G\Big(\big(P_1+i\tfrac{\varepsilon}{|\ln\varepsilon|}\mathbf{V}_1\big)(\tau)+\varepsilon\mathcal{Z}_1(\tau,\theta), \big(P_j+i\tfrac{\varepsilon}{|\ln\varepsilon|}\mathbf{V}_j\big)(\tau)+\varepsilon\rho \mathcal{Z}_j(\tau,\eta)\Big) \rho d\rho d\eta, \nonumber
\end{align}
where $\mathcal{Z}_j$ is defined in \eqref{parametrization}.

\subsubsection{Self-induced effect}

We now analyze the asymptotic expansion of the self-induced term $\Psi_{0,1}$ as $\varepsilon \to 0$. Before stating the result, we recall the following notation introduced in \eqref{Pi-proj0}. We denote by $\Pi_{k,\mathbf c}$ the $L^2(\T)$-orthogonal projector onto the mode $\cos(k\theta)$, and by $\Pi_{k,\mathbf s}$ the $L^2(\T)$-orthogonal projector onto $\sin(k\theta)$. We also define
\begin{align}\label{Proj-one-mode}
\Pi_1:=\Pi_{1,\mathbf s}+\Pi_{1,\mathbf s}\quad\hbox{and}\quad \Pi_1^\perp := \mathrm{Id}-\Pi_{1,\mathbf s}-\Pi_{1,\mathbf c},
\end{align}
which denotes the orthogonal projector onto the complement of the first Fourier mode.
\begin{pro}\label{prop-ind}
Let $\tau\in\mathbb{R},$ then  the function $\theta\mapsto \Psi_{0,1}(\tau,\theta)$ is odd. In particular,
\begin{equation}\label{Psi01:symmetry}
\Pi_{k,\mathbf c}\Psi_{0,1}(\tau,\cdot)=0,\qquad \forall\,k\geqslant 0,
\end{equation}
Moreover, we have  the expansions
\begin{align*}
\Pi_{1,\mathbf s}\Psi_{0,1}(\tau,\theta)
&=-\tfrac{\sqrt{2}}{8}\,\varepsilon\,(2p_{1,1})^{-\frac14}
\Big[
|\ln\varepsilon|
+\tfrac54\ln 8+\tfrac34\ln(p_{1,1})-\tfrac14\\ &\qquad\qquad\qquad\qquad\qquad
+\tfrac{3}{16}\varepsilon^2|\ln\varepsilon|(2p_{1,1})^{-\frac32}
+O(\varepsilon^2)
\Big]\sin\theta
\end{align*}
and 
\begin{align*}
\Pi_1^\perp\Psi_{0,1}(\tau,\theta)
&=\tfrac{\sqrt{2}}{16}\,\varepsilon\Big(-(2p_{1,1})^{-\tfrac14}\sin(3\theta)+\tfrac{3}{2}\varepsilon|\ln\varepsilon| (2p_{1,1})^{-1}\sin(2\theta)\Big)+O(\varepsilon^2),
\end{align*}
\end{pro}

\begin{proof}
In the self-induced case $j=1$, the vertical shift
$i\,\varepsilon|\ln\varepsilon|^{-1}\mathbf{V}_1(\tau)$ appears in both arguments of $G$ and cancels
because $G(\varrho,z;\varrho',z')$ depends on $z-z'$ (see \eqref{eq:G}). Hence
\begin{equation}\label{eq:Psi01-noV}
\Psi_{0,1}(\tau,\theta)
=\partial_\theta\int_0^1\fint_{\T}
G\Big(P_1(\tau)+\varepsilon\mathcal{Z}_1(\tau,\theta),\,
      P_1(\tau)+\varepsilon\rho\,\mathcal{Z}_1(\tau,\eta)\Big)\,\rho\,d\rho\,d\eta .
\end{equation}
Moreover, $\mathcal{Z}_1(\tau,-\theta)=\overline{\mathcal{Z}_1(\tau,\theta)}$ and $G$ is invariant
under complex conjugation in both arguments. After the change of variables $\eta\mapsto-\eta$, the
integral in \eqref{eq:Psi01-noV} is even in $\theta$, hence its $\partial_\theta$-derivative is odd.
This proves \eqref{Psi01:symmetry}.
Apply Lemma~\ref{prop:asympt-induced} (kernel expansion) with the substitutions
$$
z_1=p_{1,1}(\tau),\; z_2=p_{1,2}(\tau),\;
X=(x_1,x_2)=(\cos\theta,\sin\theta),\;
Y=(y_1,y_2)=(\rho\cos\eta,\rho\sin\eta),\;
\epsilon=\varepsilon,
$$
we obtain, with a uniform error on $\rho\in[0,1]$ and $\eta\in\T$,
$$
G\big(P_1+\varepsilon\mathcal{Z}_1(\theta),\,P_1+\varepsilon\rho\mathcal{Z}_1(\eta)\big)
=
\sum_{n=0}^{3}\varepsilon^{n}|\ln\varepsilon|\mathscr{A}_n(X,Y)
+\sum_{n=0}^{2}\varepsilon^{n}\mathscr{B}_n(X,Y)
+O(\varepsilon^3).
$$
Substituting this expansion into \eqref{eq:Psi01-noV}, differentiating in $\theta$, and integrating term by term, we find
\begin{align*}
\Psi_{0,1}(\tau,\theta)
&=
\sum_{n=0}^{2}\varepsilon^{n}\,\partial_\theta\left(|\ln\varepsilon|
\fint_{\T}\int_0^{1}\mathscr{A}_n(X,Y)\,\rho\,d\rho\,d\eta
+\,
\fint_{\T}\int_0^{1}\mathscr{B}_n(X,Y)\,\rho\,d\rho\,d\eta\right)
+O(\varepsilon^3|\ln\varepsilon|).
\end{align*}
Since $\mathscr{A}_0=\sqrt{p_{1,1}}$ is $\theta$-independent, its contribution vanishes after
$\partial_\theta$.\\
For $\mathscr{A}_1=\frac{2^{1/4}}{4p_{1,1}^{1/4}}(\cos\theta+\rho\cos\eta)$ we obtain
$$
\partial_\theta\fint_{\T}\int_0^1\mathscr{A}_1(X,Y)\,\rho\,d\rho\,d\eta
=-\frac{\sqrt{2}}{8}(2p_{1,1})^{-1/4}\sin\theta.
$$
    Concerning $\mathscr{A}_2,$ we have
    \begin{align*}
\mathscr{A}_2(X,Y)&=\frac{\sqrt{2}}{64p_{1,1}}\Big(3(\sin\theta-\rho\sin\eta)^2 -2\rho\cos\theta\cos\eta-3\cos^2(\theta)-3\rho^2\cos^2(\eta)\Big),
   \end{align*}
   which implies after integration
$$
\partial_\theta\fint_{\T}\int_0^1\mathscr{A}_2(X,Y)\,\rho\,d\rho\,d\eta
=\frac{3\sqrt{2}}{64p_{1,1}}\sin(2\theta).
$$
For $\mathscr{A}_3$, one has
 $$
\mathscr{A}_3(X,Y)=\frac{2^{3/4}}{256 p_{1,1}^{7/4}}(\cos\theta +\rho\cos\eta )\Big(5 \cos^2\theta -2\rho\cos\theta\cos\eta  + 5\rho^2\cos^2(\eta) -3( \sin\theta -\rho\sin\eta)^2\Big).
  $$
  After averaging in $\eta$ and integrating in $\rho$, we get
$$
\partial_\theta\fint_{\T}\int_0^1\mathscr{A}_3(X,Y)\,\rho\,d\rho\,d\eta
=-\frac{3\sqrt{2}}{128}(2p_{1,1})^{-7/4}\big(\sin\theta+2\sin(3\theta)\big).
$$
The term 
\begin{align*}
\mathscr{B}_0(X,Y)&=\sqrt{p_{1,1}}\left(\frac54\ln8+\frac{3}{4}\ln(p_{1,1})-2\right)-\frac{\sqrt{p_{1,1}}}{2}\ln\big(\mathscr{D}(\theta,\eta,\rho)\big)
\end{align*}
yields no contribution after $\partial_\theta$ because the only
$\theta$-dependence occurs through $\ln\mathscr{D}(\theta,\eta,\rho)$ with
$$\mathscr{D}(\theta,\eta,\rho)=1+\rho^2-2\rho\cos(\eta-\theta),$$ and the change of variables $\eta\mapsto\eta+\theta$ eliminates $\theta$. That is, 
\begin{align*} 
\partial_\theta\fint_{\T}\int_0^{1}\mathscr{B}_0(X,Y)\rho d\rho d\eta &=0.
\end{align*}
As for $\mathscr{B}_1$, we have the expression
\begin{align*}
  \mathscr{B}_1(X,Y)&=\tfrac{\sqrt{2}}{4(2p_{1,1})^{1/4} }\big(\tfrac54\ln8+\tfrac34\ln(p_{1,1})-1\big)\big(\cos(\theta)+\rho \cos(\eta)\big)\\ &\quad+\tfrac{\sqrt{2}}{4 (2p_{1,1})^{1/4}}\big({\mathscr{D}(\theta,\eta,\rho)}\big)^{-1}(\cos(\theta)+\rho \cos(\eta))(\cos(\theta)-\rho \cos(\eta))^2\\ &\quad-\tfrac{\sqrt{2}}{8 (2p_{1,1})^{1/4}}\big(\cos(\theta)+\rho \cos(\eta)\big)\ln\big(\mathscr{D}(\theta,\eta,\rho)\big).
\end{align*}
Using the integral identities \eqref{id1}--\eqref{id14}, one obtains
$$
\partial_\theta\fint_{\T}\int_0^{1}\mathscr{B}_1(X,Y)\,\rho\,d\rho\,d\eta
=
-\frac{\sqrt{2}}{16}(2p_{1,1})^{-1/4}
\Big[\sin(3\theta)+\tfrac12\big(5\ln 8+3\ln(p_{1,1})-1\big)\sin\theta\Big].
$$
Finally, the $\mathscr{B}_2$ term is given by
\begin{align*}
    \tfrac{16 p_{1,1}}{\sqrt{2}}\mathscr{B}_2(X,Y)&=\left(\tfrac54 \ln(8)-2+\tfrac34 \ln(p_{1,1})\right)\mathscr{B}_{2,1}(\theta,\eta,\rho)+\tfrac{1}{4}\left(\tfrac{15}{4}\ln 8-1+\tfrac94 \ln(p_{1,1})\right)\mathscr{D}(\theta,\eta,\rho)\\ &\quad -\tfrac{1}{2}\left(\mathscr{B}_{2,1}(\theta,\eta,\rho)+\tfrac34\mathscr{D}(\theta,\eta,\rho)\right)\ln\big(\mathscr{D}(\theta,\eta,\rho)\big)+\tfrac{1}{2}\mathscr{B}_{2,2}(\theta,\eta,\rho),
\end{align*}
with
\begin{align*}
   \mathscr{B}_{2,1}(\theta,\eta,\rho)&:=\rho\cos(\theta)\cos(\eta)-\tfrac32 \cos^2(\theta)-\tfrac32 \rho^2\cos^2(\eta),\\
    \mathscr{B}_{2,2}(\theta,\eta,\rho)&:=\tfrac{1}{4 \mathscr{D}(\theta,\eta,\rho)^2}\left(\cos(\theta) - \rho \cos(\eta)\right)^2 \Big(-11 - 44 \rho^2 - 11 \rho^4 - 2 (1 + 3 \rho^2) \cos(2\theta) \\
    &\quad+ \cos(4\theta)+ 
   6 \rho \cos(3\theta-\eta)  + 40 \rho \cos(\theta-\eta) + 40 \rho^3 \cos(\theta-\eta) - 
   14 \rho^2 \cos(2\theta-2\eta) \\ &\quad - 6 \rho^2 \cos(2\theta) \cos(2\eta) - 
   2 \rho^4 \cos(2 \eta) + \rho^4 \cos(4 \eta) + 
   2 \rho \cos(\theta+\eta) + 2 \rho^3 \cos(\theta+\eta) \\ &\quad- 
   2 \rho^2 \cos(2\theta+2\eta) + 6 \rho^3 \cos(\theta-3\eta)\Big).
\end{align*}
One can readily verify that $\mathscr{D}$, $\mathscr{B}_{2,1}$ and  $\mathscr{B}_{2,2}$ are all invariant under the shift $(\theta,\eta)\mapsto (\theta+\pi,\eta+\pi)$. Consequently, $\mathscr{B}_{2}$ inherits this symmetry. It follows that the integral
$$
\theta\mapsto \partial_\theta\fint_{\T}\int_0^{1} \mathscr{B}_{2}(X,Y)\rho d\rho d\eta,
$$ 
is $\pi$-periodic. Thus,  the coefficient of the $\sin\theta$  term in its Fourier expansion vanishes.
Putting together the previous identities and grouping the $\sin\theta$-mode and the higher modes gives exactly the stated expansions for $\Pi_{1,\mathbf s}\Psi_{0,1}$ and $\Pi_1^\perp\Psi_{0,1}$.
\end{proof}
\subsubsection{Interaction effect} We now derive the asymptotic expansion of the interaction contribution
$\Psi_{0,2}(\tau,\theta)$ defined in \eqref{Psi0j}. A key simplification comes from the fact that, away from the diagonal $(\varrho,z)=(\varrho',z')$, the kernel $G$ is harmonic with respect to the degenerate elliptic operator $\mathbf{L}$  as shown  in  \eqref{eq:G-harmonic}. In particular, in view of this identity, we obtain 
\begin{equation}\label{eq:G-harmonic-2}
\Big(2p_{1,1}\,\partial_{p_{1,1}}^{2}+\partial_{p_{1,2}}^{2}\Big)G(P_1,P_2)=0,
\qquad
\Big(2p_{2,1}\,\partial_{p_{2,1}}^{2}+\partial_{p_{2,2}}^{2}\Big)G(P_1,P_2)=0,
\qquad P_1\neq P_2.
\end{equation}
\begin{pro}\label{prop-int}
With the notation \eqref{Proj-one-mode}, the interaction term $\Psi_{0,2}(\tau,\theta)$ admits the expansion
\begin{align*}
\nonumber\Pi_1 \Psi_{0,2}(\tau, \theta) &=  \tfrac{1}{2} \varepsilon \Big(\nabla_1 G(P_1, P_2)+\varepsilon |\ln\varepsilon|^{-1}(\mathbf{V}_1-\mathbf{V}_2)\, \partial_{p_{1,2}}\nabla_1 G(P_1,P_2)
\Big)\cdot \partial_\theta \mathcal{Z}_1(\tau, \theta)\\ &\quad \quad+\tfrac{1}{\sqrt{2}}\varepsilon^3|\ln\varepsilon|\,(2p_{1,1})^{\frac14}\mathcal{C}_0 \sin(\theta) 
 + O\left( \varepsilon^3 |\ln \varepsilon|^{-\frac12} \right)
 \end{align*}
 and
 \begin{align*}
\Pi_1^\perp \Psi_{0,2}(\tau,\theta)&= \tfrac{1}{2}\varepsilon^{2}\partial_{p_{1,1}}\partial_{p_{1,2}}G(P_1,P_2)\cos(2\theta)\nonumber\\
    &\quad-\tfrac{1}{4}\varepsilon^{2} (2p_{1,1})^{-\frac12}\Big(2p_{1,1}\partial_{p_{1,1}}^2G(P_1,P_2) -\partial_{p_{1,2}}^2G(P_1,P_2)\Big)\sin(2\theta)+ O\left( \varepsilon^3 |\ln \varepsilon|^{\frac{3}{2}} \right), \nonumber
\end{align*}
where
\begin{align}\label{def:C0}
\mathcal{C}_0
&:=- \tfrac{1}{16\sqrt{p_{1,1}}}|\ln\varepsilon|^{-1}\Big[(\partial_{p_{1,1}}^2G)- (2p_{1,1})^{-1}(\partial_{p_{1,2}}^2G)  \Big].
\end{align}
\end{pro}
\begin{proof}
Recall from \eqref{Psi0j} that 
\begin{align*}
    \Psi_{0,2}(\tau,\theta)&=\partial_\theta \int_0^1 \fint_{\T}G(P_1+i\varepsilon |\ln\varepsilon|^{-1}\mathbf{V}_1+\varepsilon \mathcal{Z}_1(\tau,\theta), P_2+i\varepsilon |\ln\varepsilon|^{-1}\mathbf{V}_2+\varepsilon {\rho}\mathcal{Z}_2(\tau,\eta))\ \rho d\rho d\eta.
\end{align*}
Since $G$ depends on the vertical variables only through $z-z'$, we may shift both arguments by
$-i\varepsilon|\ln\varepsilon|^{-1}\mathbf{V}_2(\tau)$ and obtain
\begin{align*}
\Psi_{0,2}(\tau,\theta)
&= \partial_\theta \int_0^1 \fint_{\T}G(P_1+i\varepsilon |\ln\varepsilon|^{-1}(\mathbf{V}_1-\mathbf{V}_2)+\varepsilon \mathcal{Z}_1(\tau,\theta), P_2+\varepsilon {\rho}\mathcal{Z}_2(\tau,\eta))\ \rho d\rho d\eta,
\end{align*}
To simplify the notation, we define
\begin{align*}
 \mathtt{Z}_1&:= i|\ln\varepsilon|^{-1}(\mathbf{V}_1-\mathbf{V}_2)(\tau)+\mathcal{Z}_1(\tau,\theta),\qquad \mathtt{Z}_2:= \rho \mathcal{Z}_2(\tau,\eta),
\end{align*}
Expanding $G(\cdot,\cdot)$ around $(P_1,P_2)$ using Taylor's formula gives
\begin{align}\label{G-expand}
&G(P_1+\varepsilon\,\mathtt{Z}_1,P_2+\varepsilon\,\mathtt{Z}_2)
=G(P_1,P_2)+\varepsilon\,\nabla_1G(P_1,P_2)\cdot\mathtt{Z}_1
+\varepsilon\,\nabla_2G(P_1,P_2)\cdot\mathtt{Z}_2 \nonumber\\
&\quad
+\frac{\varepsilon^2}{2}\sum_{j,k\in\{1,2\}}\nabla_{j,k}^3G(P_1,P_2)[\mathtt{Z}_j,\mathtt{Z}_k]
+\frac{\varepsilon^3}{6}\sum_{n,j,k\in\{1,2\}}\nabla_{n,j,k}^3G(P_1,P_2)[\mathtt{Z}_n,\mathtt{Z}_j,\mathtt{Z}_k]+\varepsilon^4\,G_{r,\varepsilon},
\end{align}
where $G_{r,\varepsilon}$ involves the fourth derivatives of $G$. In the regime considered here
($|P_1-P_2|\sim |\ln\varepsilon|^{-1/2}$), one has the rough bound
\begin{equation*}
\varepsilon^4G_{r,\varepsilon}=O(\varepsilon^4|\ln\varepsilon|^{2}).
\end{equation*}
Using $\fint_{\T}\mathcal{Z}_2(\tau,\eta)\,d\eta=0$, we find
$$
\fint_{\T}\int_0^1 \nabla_2G(P_1,P_2)\cdot\mathtt{Z}_2\,\rho\,d\rho\,d\eta=0,
\qquad
\fint_{\T}\int_0^1 \nabla_{1,2}^2G(P_1,P_2)[\mathtt{Z}_1,\mathtt{Z}_2]\,\rho\,d\rho\,d\eta=0.
$$
Moreover, after averaging with respect to $\eta$, the quantity 
$\nabla_2^2 G(P_1,P_2)[\mathtt{Z}_2,\mathtt{Z}_2]$ becomes independent of $\theta$, and therefore it vanishes upon differentiation with respect to $\theta$.
Therefore, after inserting \eqref{G-expand} into the definition of $\Psi_{0,2}$ and performing the
$\rho$-integration, we obtain
\begin{align}\label{psi02-1}
\Psi_{0,2}(\tau,\theta)
&=\frac{\varepsilon}{2}\,\nabla_1G(P_1,P_2)\cdot \partial_\theta\mathcal{Z}_1
+\frac{\varepsilon^2}{4}\,\partial_\theta\nabla_1^2G(P_1,P_2)[\mathcal{Z}_1,\mathcal{Z}_1]\nonumber\\
&\quad+\frac{\varepsilon^2}{2}|\ln\varepsilon|^{-1}(\mathbf{V}_1-\mathbf{V}_2)\,
\partial_\theta\nabla_1^2G(P_1,P_2)[i,\mathcal{Z}_1]\nonumber\\
&\quad+\frac{\varepsilon^3}{6}\,\partial_\theta\int_0^1\fint_{\T}
\sum_{n,j,k\in\{1,2\}}\nabla_{n,j,k}^3G(P_1,P_2)[\mathtt{Z}_n,\mathtt{Z}_j,\mathtt{Z}_k]\,
\rho\,d\rho\,d\eta
+O(\varepsilon^4|\ln\varepsilon|^{2}).
\end{align}
Notice that one readily obtains the following, 
\begin{align*}
\nabla_1^2G(P_1,P_2)[i,\mathcal{Z}_1]
&=\partial_{p_{1,2}}\nabla_1G(P_1,P_2)\cdot\mathcal{Z}_1.
\end{align*}
Then, by writing $\mathcal{Z}_1=A_{1,1}\cos\theta+i A_{1,2}\sin\theta$ with
$A_{1,1}=(2p_{1,1})^{1/4}$ and $A_{1,2}=(2p_{1,1})^{-1/4}$, it yields
\begin{align*}
\partial_\theta\nabla_1^2G(P_1,P_2)[\mathcal{Z}_1,\mathcal{Z}_1]
&=\Big((\partial_{p_{1,2}}^2G)A_{1,2}^2-(\partial_{p_{1,1}}^2G)A_{1,1}^2\Big)\sin(2\theta)
+2A_{1,1}A_{1,2}(\partial_{p_{1,1}}\partial_{p_{1,2}}G)\cos(2\theta),
\end{align*}
which gives the  desired $\sin(2\theta)$ and $\cos(2\theta)$ contributions in $\Pi_1^\perp \Psi_{0,2}$. 
In the cubic term in \eqref{psi02-1}, all contributions involving an odd number of 
$\mathtt{Z}_2$ vanish after averaging with respect to $\eta$. Consequently, the only potentially nonzero terms are
$\nabla_1^3G[\mathtt{Z}_1,\mathtt{Z}_1,\mathtt{Z}_1]$ and $\nabla_1\nabla_2^2G[\mathtt{Z}_1,\mathtt{Z}_2,\mathtt{Z}_2]$.
A direct expansion gives
\begin{align}\label{third-der}
\partial_\theta\!\int_0^1\fint_{\T}\!\sum_{n,j,k}\nabla_{n,j,k}^3G[\mathtt{Z}_n,\mathtt{Z}_j,\mathtt{Z}_k]\rho\,d\rho\,d\eta
&=\frac12\,\partial_\theta\nabla_1^3G(P_1,P_2)[\mathcal{Z}_1,\mathcal{Z}_1,\mathcal{Z}_1]\nonumber\\
&\quad+3\,\partial_\theta\!\int_0^1\fint_{\T}\nabla_1\nabla_2^2G(P_1,P_2)[\mathcal{Z}_1,\mathcal{Z}_2,\mathcal{Z}_2]\rho^3\,d\rho\,d\eta\nonumber\\
&\quad+\frac{3}{2}|\ln\varepsilon|^{-1}(\mathbf{V}_1-\mathbf{V}_2)\,
\partial_\theta\nabla_1^3G(P_1,P_2)[i,\mathcal{Z}_1,\mathcal{Z}_1]\nonumber\\
&\quad+\frac{3}{2}|\ln\varepsilon|^{-2}(\mathbf{V}_1-\mathbf{V}_2)^2\,
\partial_\theta\nabla_1^3G(P_1,P_2)[i,i,\mathcal{Z}_1].
\end{align}
The mixed term with $\nabla_1\nabla_2^2G$ vanishes identically by \eqref{eq:G-harmonic-2} (in the
second argument), 
\begin{equation}\label{third-deriv2}
\partial_\theta\!\int_0^1\fint_{\T}\nabla_1\nabla_2^2G(P_1,P_2)[\mathcal{Z}_1,\mathcal{Z}_2,\mathcal{Z}_2]\rho^3\,d\rho\,d\eta=0.
\end{equation}
In fact, one has
$$
\begin{aligned}
\nabla_1\nabla_2^2G(P_1,P_2)[\mathcal{Z}_1,\mathcal{Z}_2,\mathcal{Z}_2] = &
A_{1,1}\Big[  A_{2,1}^2  \, \partial_{p_{1,1}} \partial_{p_{2,1}}^2 G \cos^2(\eta)+  A_{2,2}^2  \partial_{p_{1,1}}  \partial_{p_{2,2}}^2 G \sin^2(\eta) \\
& + 2 A_{2,1} A_{2,2}   \partial_{p_{1,1}}\partial_{p_{2,1}}\partial_{p_{2,2}} G\cos(\eta) \sin(\eta)\Big] \cos(\theta) \\
& + A_{1,2}\Big[ A_{2,1}^2  \, \partial_{p_{1,2}} \partial_{p_{2,1}}^2G \,\cos^2(\eta) +  A_{2,2}^2\, \partial_{p_{1,2}}\partial_{p_{2,2}}^2 G \sin^2(\eta) \\
& +  2A_{2,1} A_{2,2} \, \partial_{p_{1,2}}\partial_{p_{2,1}}\partial_{p_{2,2}} G  \, \cos(\eta) \sin(\eta) \Big]\sin(\theta), 
\end{aligned}
$$
with
$$A_{j,1}=(2p_{j,1})^{1/4} \qquad A_{j,2}=(2p_{j,1})^{-1/4}.$$
Integrating in $\eta$ gives
$$
\begin{aligned}
\frac{1}{2\pi}\int_{0}^{2\pi}\nabla_1\nabla_2^2G(P_1,P_2)[\mathcal{Z}_1,\mathcal{Z}_2,\mathcal{Z}_2]d\eta = &
\tfrac12 A_{1,1}\partial_{p_{1,1}} \Big[  A_{2,1}^2  \,  \partial_{p_{2,1}}^2 G +  A_{2,2}^2   \partial_{p_{2,2}}^2 G \Big] \cos(\theta) \\
& +\tfrac12 A_{1,2}\partial_{p_{1,2}}\Big[ A_{2,1}^2  \, \partial_{p_{2,1}}^2G +  A_{2,2}^2\, \partial_{p_{2,2}}^2 G  \Big]\sin(\theta). 
\end{aligned}
$$
which is identically zero according to \eqref{eq:G-harmonic-2}.
Now decomposing $\nabla_1^3G(P_1,P_2)[\mathcal{Z}_1,\mathcal{Z}_1,\mathcal{Z}_1]$ into Fourier modes gives
$$
\nabla_1^3G(P_1,P_2)[\mathcal{Z}_1,\mathcal{Z}_1,\mathcal{Z}_1]
=\mathcal{C}_1\cos\theta+\mathcal{C}_2\cos(3\theta)+\mathcal{C}_3\sin\theta+\mathcal{C}_4\sin(3\theta),
$$
with
\begin{align*}
\mathcal{C}_1\cos(\theta)&+\mathcal{C}_2\cos(3\theta)=A_{1,1}^3(\partial_{p_{1,1}}^3G) \cos^3(\theta)+3A_{1,1} A_{1,2}^2(\partial_{p_{1,1}}\partial_{p_{1,2}}^2G) \cos(\theta)\sin^2(\theta)\\
=&A_{1,1}^3(\partial_{p_{1,1}}^3G) \Big(\frac34 \cos(\theta)+\frac14\cos(3\theta)\Big)+\frac34A_{1,1} A_{1,2}^2(\partial_{p_{1,1}}\partial_{p_{1,2}}^2G )\Big(\cos(\theta)-\cos(3\theta)\Big),
\end{align*}
and
\begin{align*}
\mathcal{C}_3\sin(\theta)&+\mathcal{C}_4\sin(3\theta)=(A_{1,2}^3\partial_{p_{1,2}}^3G) \sin^3(\theta)+3A_{1,2} A_{1,1}^2(\partial_{p_{1,1}}^2\partial_{p_{1,2}}G) \sin(\theta)\cos^2(\theta)\\
=&A_{1,2}^3(\partial_{p_{1,2}}^3G) \Big(\frac34 \sin(\theta)-\frac14\sin(3\theta)\Big)+\frac34A_{1,2} A_{1,1}^2(\partial_{p_{1,1}}^2\partial_{p_{1,2}}G) \Big(\sin(\theta)+\sin(3\theta)\Big).
\end{align*}
From these identities together with the expressions of $A_{j,1},A_{j,2}$ we conclude  that
\begin{align*}
\mathcal{C}_1
&=\tfrac34(2p_{1,1})^{\frac14}\Big[(2p_{1,1})^{\frac12}(\partial_{p_{1,1}}^3G) +(2p_{1,1})^{-\frac12}(\partial_{p_{1,1}}\partial_{p_{1,2}}^2G) \Big]\\
&=\tfrac34(2p_{1,1})^{\frac14}\partial_{p_{1,1}}\Big[ (2p_{1,1})^{\frac12}(\partial_{p_{1,1}}^2G)+ (2p_{1,1})^{-\frac12}(\partial_{p_{1,2}}^2G)  \Big]\\ 
& \quad -\tfrac34\Big[ (2p_{1,1})^{-\frac14}(\partial_{p_{1,1}}^2G)- (2p_{1,1})^{-\frac54}(\partial_{p_{1,2}}^2G)  \Big]\\
&=-\tfrac34\Big[ (2p_{1,1})^{-\frac14}(\partial_{p_{1,1}}^2G)- (2p_{1,1})^{-\frac54}(\partial_{p_{1,2}}^2G)  \Big], 
\end{align*}
where in the last expression, we used \eqref{eq:G-harmonic-2}. Similar computations allow to get
\begin{align*}
\mathcal{C}_3
&=\tfrac34\Big[ A_{1,2} A_{1,1}^2(\partial_{p_{1,1}}^2\partial_{p_{1,2}}G)+ A_{1,2}^3(\partial_{p_{1,2}}^3G)  \Big]\\
&=\tfrac34(2p_{1,1})^{-\frac14}\Big[ (2p_{1,1})^{\frac12}(\partial_{p_{1,1}}^2\partial_{p_{1,2}}G)+ (2p_{1,1})^{-\frac12}(\partial_{p_{1,2}}^3G)  \Big]\\
&=\tfrac34(2p_{1,1})^{-\frac14}\partial_{p_{1,2}}\Big[(2p_{1,1})^{\frac12}(\partial_{p_{1,1}}^2G)+  (2p_{1,1})^{-\frac12}(\partial_{p_{1,2}}^2G)  \Big]\\
&=0,
\end{align*}
Concerning the terms $\mathcal{C}_2$ and $\mathcal{C}_4$, their explicit formulas are not required; however, their scaling with respect to  $\varepsilon$ is essential.
Consequently,
\begin{equation}\label{third-deriv1}
\frac12\,\partial_\theta\nabla_1^3G(P_1,P_2)[\mathcal{Z}_1,\mathcal{Z}_1,\mathcal{Z}_1]
=\frac{1}{2}\mathcal{C}_1\sin\theta-\frac{3}{2}\mathcal{C}_2\sin(3\theta)+\frac{3}{2}\mathcal{C}_4\cos(3\theta).
\end{equation}
Similarly we get
\begin{align*}
   \partial_\theta\nabla_1^3G(P_1,P_2)[i,i,\mathcal{Z}_1]&= (2p_{1,1})^{-\frac{1}{4}}(\partial_{p_{1,2}}^3G) \cos(\theta)-(2p_{1,1})^{\frac{1}{4}}(\partial_{p_{1,1}}\partial_{p_{1,2}}^2G)\sin(\theta) \\
   &= \partial_{p_{1,2}}^2\nabla_1 G \cdot\partial_\theta \mathcal{Z}_1,
\end{align*}
and
\begin{align*}
   \partial_\theta\nabla_1^3G(P_1,P_2)[i,\mathcal{Z}_1,\mathcal{Z}_1]&= \Big((2p_{1,1})^{-\frac{1}{2}} (\partial_{p_{1,2}}^3G) -(2p_{1,1})^{\frac{1}{2}}(\partial_{p_{1,1}}^2\partial_{p_{1,2}}G)\Big)\sin(2\theta)\\ &\quad+ 2 (\partial_{p_{1,1}}\partial_{p_{1,2}}^2G)\cos(2\theta).
\end{align*}
Inserting \eqref{third-deriv2}, \eqref{third-deriv1} and the  later two  identities into \eqref{third-der} gives 
\begin{align*}
&\frac{\varepsilon^3}{6}\,\partial_\theta\!\int_0^1\fint_{\T}\!\sum_{n,j,k}\nabla_{n,j,k}^3G[\mathtt{Z}_n,\mathtt{Z}_j,\mathtt{Z}_k]\rho\,d\rho\,d\eta
=\frac{1}{\sqrt{2}}\varepsilon^3|\ln\varepsilon|(2p_{1,1})^{\frac14}\mathcal{C}_0\sin\theta-\frac{\varepsilon^3}{4}\,\mathcal{C}_2\sin(3\theta)\nonumber\\
&\quad+\frac{\varepsilon^3}{4}\,\mathcal{C}_4\cos(3\theta)+\frac{\varepsilon^3}{4}\,|\ln\varepsilon|^{-1}(\mathbf{V}_1-\mathbf{V}_2)\,
\Big((2p_{1,1})^{-\frac{1}{2}} (\partial_{p_{1,2}}^3G) -(2p_{1,1})^{\frac{1}{2}}(\partial_{p_{1,1}}^2\partial_{p_{1,2}}G)\Big)\sin(2\theta)\nonumber\\
&\quad+\frac{\varepsilon^3}{4}\,|\ln\varepsilon|^{-1}(\mathbf{V}_1-\mathbf{V}_2)\,
2 (\partial_{p_{1,1}}\partial_{p_{1,2}}^2G)\cos(2\theta)+\frac{\varepsilon^3}{4}\,|\ln\varepsilon|^{-2}(\mathbf{V}_1-\mathbf{V}_2)^2\,
\partial_{p_{1,2}}^2\nabla_1 G \cdot\partial_\theta \mathcal{Z}_1,
\end{align*}
where $\mathcal{C}_0$ is given by \eqref{def:C0}.
Since the  third derivatives  scale like $|P_1-P_2|^{-3}\sim |\ln\varepsilon|^{3/2}$, then
the projection by $\Pi_{1}^\perp$ yields the remainder $O(\varepsilon^3|\ln\varepsilon|^\frac32)$.
\end{proof}

\begin{pro}\label{prop:G1-trivial}
Assume that the core trajectories $P_1,P_2$ satisfy \eqref{eq-points-general}. 
Then, we have the asymptotics,
 \begin{align*}
\Pi_{1,{\bf{c}}} \mathbf{G}_1(\varepsilon, 0, 0)(\tau,\theta)& =\bigg( \frac{1}{\sqrt{2}} \frac{\varepsilon^2}{ |\ln\varepsilon|}  (2p_{1,1})^{-\frac{1}{4}} (\mathbf{V}_1-\mathbf{V}_2)  \partial_{p_{1,2}}^2 G(P_1,P_2)
+ O\left( \varepsilon^3 |\ln \varepsilon|^{-\frac12} \right)\bigg)\cos(\theta) ,
\end{align*}
and
\begin{align*}
\Pi_{1,{\bf{s}}} \mathbf{G}_1(\varepsilon, 0, 0)(\tau,\theta)& = \varepsilon^2 (2p_{1,1})^\frac14\bigg(  \, \dot{\mathbf{V}}_1(\tau) -\frac{1}{\sqrt{2}|\ln\varepsilon|}  (\mathbf{V}_1-\mathbf{V}_2) \partial_{p_{1,2}}\partial_{p_{1,1}} G(P_1,P_2)\\ &\qquad\qquad\qquad  
+\varepsilon|\ln\varepsilon|\Big(\mathcal{C}_0-\frac{3}{256 p_{1,1}^2} \Big)+ O\left( \varepsilon^3 |\ln \varepsilon|^{-\frac12} \right)\bigg)\sin( \theta),
\end{align*}
where $\mathcal{C}_0$ is given by \eqref{def:C0}.
Moreover,  
\begin{align*}
\Pi_1^\perp\mathbf{G}_1(\varepsilon,0,0)(\tau,\theta)
&=-\frac{1}{8}\varepsilon(2p_{1,1})^{-\frac14}\sin(3\theta)+\frac{1}{4}\varepsilon^2|\ln\varepsilon|
\bigg(
\frac{3}{8p_{1,1}}
-2\sqrt{p_{1,1}}\,|\ln\varepsilon|^{-1}\,\partial_{p_{1,1}}^{2}G(P_1,P_2) \\
&
\quad +\frac{1}{\sqrt{p_{1,1}}}\,|\ln\varepsilon|^{-1}\,\partial_{p_{1,2}}^{2}G(P_1,P_2)
\bigg)\sin(2\theta) \\
&\quad +\varepsilon^2|\ln\varepsilon|
\bigg(
\frac{\dot p_{1,1}}{4p_{1,1}}
+\frac{1}{\sqrt{2}}|\ln\varepsilon|^{-1}\,\partial_{p_{1,1}}\partial_{p_{1,2}}G(P_1,P_2)
\bigg)\cos(2\theta)
+O(\varepsilon^{2}).
\end{align*}
\end{pro}

\begin{proof}
We start from the reduced equation \eqref{eq-G1-red}. Setting $f_1=f_2=0$ (so that $w_1\equiv 1$) and
using \eqref{Pol-ah2}--\eqref{Psi0j}, we obtain
\begin{align}\label{eq:G1-00-decomp}
\mathbf{G}_1(\varepsilon,0,0)(\tau,\theta)
&=-\varepsilon\Big(|\ln\varepsilon|\,\dot P_1(\tau)+i\varepsilon\,\dot{\mathbf{V}}_1(\tau)\Big)\cdot
i\,\partial_\theta\mathcal{Z}_1(\tau,\theta)
+\varepsilon^2|\ln\varepsilon|\,\frac{\dot p_{1,1}(\tau)}{8p_{1,1}(\tau)}\,\partial_\theta\!\big(\sin(2\theta)\big) \nonumber\\
&\quad+\sqrt{2}\,\Psi_{0,1}(\tau,\theta)+\sqrt{2}\,\Psi_{0,2}(\tau,\theta).
\end{align}
Projecting \eqref{eq:G1-00-decomp} onto $\mathrm{span}\{\cos\theta,\sin\theta\}$ and using
Propositions~\ref{prop-ind} and \ref{prop-int} (self-induced and interaction expansions), we find
\begin{align*}
\Pi_1\mathbf{G}_1(\varepsilon,0,0)
&=-\varepsilon\Big(|\ln\varepsilon|\,\dot P_1+i\varepsilon\,\dot{\mathbf{V}}_1\Big)\cdot
i\,\partial_\theta\mathcal{Z}_1
\\
&\quad+\frac{\varepsilon}{\sqrt{2}}\Big(\nabla_1G(P_1,P_2) +\frac{\varepsilon}{|\ln\varepsilon|}(\mathbf{V}_1-\mathbf{V}_2)
\partial_{p_{1,2}}\nabla_1G(P_1,P_2)
\Big)\cdot\partial_\theta\mathcal{Z}_1\\
&\quad -\,\frac{\varepsilon}{4}(2p_{1,1})^{-\frac14}
\Big(|\ln\varepsilon|+\tfrac54\ln 8+\tfrac34\ln(p_{1,1})-\tfrac14\Big)\sin(\theta)
+\varepsilon^3|\ln\varepsilon|(2p_{1,1})^{\frac14}\mathcal{C}_0\sin(\theta) \\
&\quad -\frac{3}{64}\varepsilon^3 |\ln \varepsilon| (2p_{1,1})^{-\frac74}+O\left( \varepsilon^3 |\ln \varepsilon|^{-\frac12} \right)\cos(\theta)
+O\left( \varepsilon^3 |\ln \varepsilon|^{-\frac12} \right)\sin(\theta).
\end{align*}
 By imposing \eqref{eq-points-general},   the leading-order terms in $\varepsilon$ cancels, leaving
\begin{align*}
&\Pi_1 \mathbf{G}_1(\varepsilon, 0, 0) = {\varepsilon^2 (2p_{1,1})^\frac14 \, \dot{\mathbf{V}}_1\sin(\theta)}+\frac{1}{\sqrt{2}}\frac{\varepsilon^2}{ |\ln\varepsilon|} (\mathbf{V}_1-\mathbf{V}_2)\,  \partial_{p_{1,2}}\nabla_1 G(P_1,P_2)
\cdot \partial_\theta \mathcal{Z}_1\\ &\qquad+\varepsilon^3|\ln\varepsilon|(2p_{1,1})^\frac14\Big(\mathcal{C}_0-\frac{3}{256 p_{1,1}^2} \Big)\sin(\theta)+ O\left( \varepsilon^3 |\ln \varepsilon|^{-\frac12} \right)\sin(\theta) 
+ O\left( \varepsilon^3 |\ln \varepsilon|^{-\frac12} \right)\cos(\theta).
\end{align*}
Extracting separately the $\cos\theta$ and $\sin\theta$ components gives the stated formulas for
$\Pi_{1,\mathbf c}\mathbf{G}_1$ and $\Pi_{1,\mathbf s}\mathbf{G}_1$.
Finally, projecting \eqref{eq:G1-00-decomp} onto $\Pi_1^\perp$, 
\begin{align*}
\Pi_1^\perp\mathbf{G}_1(\varepsilon,0,0)
&=\varepsilon^2|\ln\varepsilon|\,\frac{\dot p_{1,1}}{4p_{1,1}}\cos(2\theta)
+\sqrt{2}\,\Pi_1^\perp\Psi_{0,1}+\sqrt{2}\,\Pi_1^\perp\Psi_{0,2}
+O(\varepsilon^2),
\end{align*}
 and substituting the explicit expansions for $\Pi_1^\perp\Psi_{0,1}$ and $\Pi_1^\perp\Psi_{0,2}$ in  Propositions~\ref{prop-ind} and \ref{prop-int}, respectively, yields the claimed expansion.
\end{proof}

\subsection{Frequency renormalization}\label{SEc-Freq-Norma}

The time reparameterization introduced in this section will play a central role in Section~\ref{section-linearization}. Its purpose is to remove the leading factor
$\sqrt{2p_{1,1}}$ that appears in the transport part of the linearized operator. This new parameterization ensures that the dynamics is naturally $2\pi$--periodic in the new time variable, thereby allowing us to place the forthcoming analysis, and in particular the linearization in Section~\ref{section-linearization}, into a more canonical form.

Let $p_{1,1}:\R\to(0,\infty)$ be the function provided by Proposition~\ref{prop-ring}.
In particular, $p_{1,1}$ is $T(\varepsilon,\lambda,\kappa)$--periodic, where
$T(\varepsilon,\lambda,\kappa)$ is the period given in Proposition~\ref{prop-ring}.  
Introduce the variable
\begin{equation}\label{def:phi-change}
\varphi=\varphi(\tau):=\omega\int_0^\tau \sqrt{2p_{1,1}(s)}\,ds. 
\end{equation}
This reparametrization is $2\pi$--periodic (i.e. $\varphi(\tau+T)=\varphi(\tau)+2\pi$) if and only if
\begin{equation}\label{lambda--1}
\omega=\omega(\varepsilon,\lambda)
=\frac{2\pi}{\displaystyle\int_0^{T(\varepsilon,\lambda)}\sqrt{2p_{1,1}(s)}\,ds}\,\cdot
\end{equation}
Given a function $F_1:(\tau,\theta)\in\R\times\T\to\R$, we define $F:(\varphi,\theta)\in\R\times\T\to\R$ by
$$
F_1(\tau,\theta) \;=\; F\Big( \omega(\varepsilon,\lambda) \int_0^\tau \sqrt{2p_{1,1}(s)}\,ds,\,\theta\Big).
$$
Then $F_1$ is $T$--periodic if and only if $F$ is $2\pi$--periodic. Under this reparametrization, we introduce the renormalized unknowns that will be used throughout the remainder of the paper.
\begin{align*}
f_1(\tau,\theta)&=f(\varphi,\theta),\qquad w_1(\tau,\theta)=w(\varphi,\theta),\qquad  \gamma_1(\tau,\theta)=\gamma(\varphi,\theta).
\end{align*}
Dividing \eqref{eq-G1-red} (see also \eqref{eq-gamma-tilde}) by $(2p_{1,1})^{1/2}$ and applying the above change of variables, we obtain the equivalent formulation 
\begin{equation}\label{F-def}
\begin{aligned}
\mathbf{F}(\varepsilon,\mathtt V_1,\mathtt V_2,f)(\varphi,\theta)
&:=\varepsilon^3|\ln\varepsilon|\,\omega\,\partial_\varphi f(\varphi,\theta)
+(2\mathtt p_{1,1}(\varphi))^{-\frac12}\,\partial_\theta\!\big\{\Psi(\gamma(\varphi,\theta))\big\}\\
&\quad-\omega\Big(|\ln\varepsilon|\,\dot{\mathtt P}_1(\varphi)+i\varepsilon\,\dot{\mathtt V}_1(\varphi)\Big)\cdot
i\,\partial_\theta\gamma(\varphi,\theta)\\
&\quad+\varepsilon^2|\ln\varepsilon|\,\omega\,\frac{\dot{\mathtt p}_{1,1}(\varphi)}{8\mathtt p_{1,1}(\varphi)}\,
\partial_\theta\!\Big((1+2\varepsilon f(\varphi,\theta))\sin(2\theta)\Big),
\end{aligned}
\end{equation}
and the renormalized equation to consider is $\mathbf{F}(\varepsilon,\mathtt V_1,\mathtt V_2,f)=0$.
Here, we have used the notation
$$
\mathtt{P}_j(\varphi):=P_j(\tau),\quad \mathtt{P}_j(\varphi):=\big(\mathtt{p}_{j,1}(\varphi),\mathtt{p}_{j,2}(\varphi)\big),
\quad \mathtt{V}_j(\varphi):=\mathbf{V}_j(\tau),\quad \mathcal{Z}_1(\tau,\theta)=\mathcal{Z}(\varphi,\theta).
$$
The constraints arising from  Lemma~\ref{lem-red} and Proposition~\ref{prop-reversibility} can be written as
\begin{equation}\label{assumo-speed}
\begin{aligned}
\mathtt{V}_1(-\varphi)&=-\mathtt{V}_1(\varphi),\qquad \mathtt{V}_2(-\varphi)=-\mathtt{V}_2(\varphi),\\  \dot{\mathtt{V}}_{1}(\varphi+\pi)&=\dot{\mathtt{V}}_2(\varphi), \qquad
\mathtt{U}_{-}(\varphi+\pi)=-\mathtt{U}_{-}(\varphi)\quad\hbox{with}\quad \mathtt{U}_{-}(\varphi):=\big(\mathtt{V}_{1}-\mathtt{V}_{2}\big)(\varphi)\,.
\end{aligned}
\end{equation}
In view of \eqref{Pol-ah1}, the stream function splits as
\begin{align}\label{Psi-form} &\Psi(\gamma(\varphi,\theta))=\sqrt{2}\fint_{\T}\int_0^{w(\varphi,\eta)}G\big(\mathtt{P}_1(\varphi)+\varepsilon w(\varphi,\theta) \mathcal{Z}(\varphi,\theta),\mathtt{P}_1(\varphi)+\varepsilon \rho   \mathcal{Z}(\varphi,\eta)\big)\rho d\rho d\eta\\
\nonumber&\quad+\sqrt{2}\fint_{\T}\int_0^{w_\star(\varphi,\eta)}G\big(\mathtt{P}_1(\varphi)+i\,\varepsilon |\ln\varepsilon|^{-1} \mathtt{U}_{-}(\varphi)+\varepsilon w(\varphi,\theta) \mathcal{Z}(\varphi,\theta),\mathtt{P}_2(\varphi)+\varepsilon \rho   \mathcal{Z}_\star(\varphi,\eta)\big)\rho  d\rho d\eta,
\end{align}
where, for any $2\pi$--periodic function $u(\varphi,\theta)$, we use the time-shift notation,  defined in the notation paragraph (see Section \ref{sec-notat}) by
$$
u_\star(\varphi,\theta) := u(\varphi+\pi,\theta)\,.
$$
Next, we derive an asymptotic expansion of the frequency $\omega(\varepsilon,\lambda)$ with respect to the small parameter $\varepsilon$, which will be used at several stages of the analysis below.
\begin{lem}\label{lem-lam}
The frequency $\omega$ introduced in \eqref{lambda--1} satisfies
$$
\omega(\varepsilon,\lambda)=\frac{2\pi}{\sqrt{2\kappa}\,T_0(\lambda,\kappa)}+O\left({{|\ln\varepsilon|^{-\frac12}}}\right)\quad\hbox{and}\quad \partial_\lambda\omega(\varepsilon,\lambda)=-\frac{2\pi T_0'(\lambda)}{\sqrt{2\kappa}\,T^2_0(\lambda)}+O\left(\tfrac{1}{\sqrt{|\ln\varepsilon|}}\right).$$

\end{lem}

\begin{proof}
The asymptotic expansion of $\omega(\varepsilon,\lambda)$ as $\varepsilon \to 0$, for any fixed $\lambda$, follows directly from Proposition~\ref{prop-ring}. Differentiating the relation \eqref{lambda--1} with respect to $\lambda$ and using once again Proposition~\ref{prop-ring} yields the stated expansion for $\partial_\lambda \omega(\varepsilon,\lambda)$.
\end{proof}

From \eqref{eq-points-general} and the change of variables \eqref{def:phi-change}, the system for $(\mathtt P_1,\mathtt P_2)$ becomes
\begin{equation}\label{eq-points-generalL}
\begin{aligned}
 \omega |\ln\varepsilon |\,  \dot{\mathtt{P}}_j(\varphi) =&\frac{1}{2\sqrt{\mathtt{p}_{j,1}}}(\nabla_{\mathtt{P}_j} ^\perp G)(\mathtt{P}_1,\mathtt{P}_2) +\frac{1}{\sqrt{2\mathtt{p}_{j,1}}}(\nabla_{\mathtt{P}_j} ^\perp \mathtt{G})(\mathtt{P}_1,\mathtt{P}_2),
\end{aligned}
\end{equation}
with
\begin{equation}\label{nablaGtt}
(\nabla_{\mathtt{P}_j} ^\perp \mathtt{G})(\mathtt{P}_1,\mathtt{P}_2)=\frac{i}{4\sqrt{2\mathtt{p}_{j,1}}}\Big(|\ln\varepsilon|{+}\tfrac{1}{4}\big(5\ln8+3\ln(\mathtt{p}_{j,1})-1\big)\Big).
\end{equation}
In terms of the scaled variables (cf. \eqref{scaling-points}),
\begin{align}\label{p-decompos}
\mathtt{p}_{j,1}=\kappa+\tfrac{ 1}{2}(-1)^{j+1} r_\varepsilon (2\kappa)^\frac14 \mathtt{y}_1,\quad j=1,2, 
\quad \mathtt{p}_{1,2}-\mathtt{p}_{2,2}=r_\varepsilon(2\kappa)^{-\frac14}\mathtt{y}_2, \quad r_\varepsilon=\frac{(2\kappa)^\frac14}{{|\ln\varepsilon|^{\frac12}}},
\end{align}
the pair $(\mathtt y_1,\mathtt y_2)(\varphi):= {({\mathtt{x}}_1,{\mathtt{x}}_2)}(\tau)$ has some symmetry obtained a consequence of Proposition \ref{prop-mouh}, which is stated in the following corollary.
\begin{cor}\label{cor:symmetryy}
The pair $(\mathtt{y}_1,\mathtt{y}_2)(\varphi)$ is $2\pi$-periodic and satisfies
\begin{align*}
\mathtt{y}_1(-\varphi)=\mathtt{y}_1(\varphi), &\qquad \mathtt{y}_2(\varphi)=-\mathtt{y}_2(\varphi),\quad
\mathtt{y}_1(\varphi+\pi)=-\mathtt{y}_1(\varphi), & \mathtt{y}_2(\varphi+\pi)=-\mathtt{y}_2(\varphi).
\end{align*}
Moreover,    the functions $\mathtt{y}_j:(0,\infty)\times \mathbb{T}\mapsto \mathtt{y}_j(\lambda,\varphi)$, $j=1,2$, are real analytic. 
\end{cor}
\noindent We introduce the following $2\pi$--periodic functions, which provide a normalized version of those appearing in Lemma~\ref{prop-asym-derivatives-G}:
\begin{align} \label{list-functions}
  \nonumber\mathtt{f}_{2}(\varphi)&:=(2\mathtt{p}_{1,1})^{-\frac12} \tfrac{1}{\sqrt{2}}|\ln\varepsilon|^{-1}(\partial_{\mathtt{p}_{1,1}}\partial_{\mathtt{p}_{1,2}}G)(\mathtt{P}_1,\mathtt{P}_2),\quad
  \mathtt{h}_2(\varphi):=-\tfrac{1}{\sqrt{\mathtt{p}_{1,1}}} |\ln\varepsilon|^{-1}(\partial_{\mathtt{p}_{1,2}}^2G)(\mathtt{P}_1,\mathtt{P}_2),\\
 \nonumber
 \mathtt{g}_{2}(\varphi)&:=\tfrac{1}{4}(2\mathtt{p}_{1,1})^{-\frac12}\Big(\tfrac{3}{8p_{1}}-2\sqrt{\mathtt{p}_{1,1}}|\ln\varepsilon|^{-1}(\partial_{\mathtt{p}_{1,1}}^2G)(\mathtt{P}_1,\mathtt{P}_2)+\tfrac{1}{\sqrt{\mathtt{p}_{1,1}}}|\ln\varepsilon|^{-1}(\partial_{\mathtt{p}_{1,2}}^2 G)(\mathtt{P}_1,\mathtt{P}_2)\Big),
\\  
 \mathtt{g}_{3}(\varphi)&:={\tfrac18 (2\mathtt{p}_{1,1})^{-\frac34}},\end{align}
 together with
\begin{align*}\mathtt{q}_1(\varphi)&:=-\tfrac{1}{\sqrt{2}}(2\mathtt{p}_{1,1})^{-\frac14}(2\mathtt{p}_{2,1})^{\frac14}|\ln\varepsilon|^{-1}(\partial_{\mathtt{p}_{1,1}}\partial_{\mathtt{p}_{2,1}}G)(\mathtt{P}_1,\mathtt{P}_2),
      \\ \nonumber \mathtt{q}_2(\varphi)& :=-\tfrac{1}{\sqrt{2}}(2\mathtt{p}_{1,1})^{-\frac14}(2\mathtt{p}_{2,1})^{-\frac14}|\ln\varepsilon|^{-1}(\partial_{\mathtt{p}_{1,1}}\partial_{\mathtt{p}_{2,2}}G)(\mathtt{P}_1,\mathtt{P}_2),
       \\ \nonumber \mathtt{q}_3(\varphi) &:=\tfrac{1}{\sqrt{2}}(2\mathtt{p}_{1,1})^{-\frac34}(2\mathtt{p}_{2,1})^{\frac14}|\ln\varepsilon|^{-1}(\partial_{\mathtt{p}_{1,2}}\partial_{\mathtt{p}_{2,1}}G)(\mathtt{P}_1,\mathtt{P}_2),
      \\ \nonumber \mathtt{q}_4(\varphi)&:=\tfrac{1}{\sqrt{2}}(2\mathtt{p}_{1,1})^{-\frac34}(2\mathtt{p}_{2,1})^{-\frac14}|\ln\varepsilon|^{-1}(\partial_{\mathtt{p}_{1,2}}\partial_{\mathtt{p}_{2,2}}G)(\mathtt{P}_1,\mathtt{P}_2).
  \end{align*}
  The next result follows immediately from Lemma~\ref{prop-asym-derivatives-G}.
\begin{lem}\label{lem-functions}
The following expansions hold:
\begin{align} \label{asymt-list1}
 \nonumber \mathtt{f}_{2}(\varphi) 
 &=\tfrac{1}{\sqrt{2\kappa}}
 \tfrac{\mathtt{y}_1 \mathtt{y}_2}{(\mathtt{y}_1^2+\mathtt{y}_2^2)^2}  +O(|\ln\varepsilon|^{-\frac12})=:\frac{1}{\sqrt{2\kappa}}\check{\alpha}(\varphi)+O(|\ln\varepsilon|^{-\frac12}),
 \\
 \nonumber \mathtt{h}_{2}(\varphi) 
 &=
 \tfrac{\mathtt{y}_1^2- \mathtt{y}_2^2}{(\mathtt{y}_1^2+\mathtt{y}_2^2)^2}  +O(|\ln\varepsilon|^{-\frac12})=:\check{\mathtt{h}}_2(\varphi)+O(|\ln\varepsilon|^{-\frac12}),
 \\
 \mathtt{g}_{2}(\varphi) &=\tfrac{3}{16}(2\kappa)^{-\frac32} -{\tfrac12}(2\kappa)^{-\frac12}\tfrac{\mathtt{y}_1^2-\mathtt{y}_2^2}{(\mathtt{y}_1^2+\mathtt{y}_2^2)^2}  +O(|\ln\varepsilon|^{-\frac12}),\\
   \nonumber\mathtt{g}_{3}(\varphi) &=\tfrac{1}{8 (2\kappa)^{3/4}}+O(|\ln\varepsilon|^{-\frac12}),
     \end{align} 
and    \begin{align*}\mathtt{q}_1(\varphi)&=\tfrac{1}{2\sqrt{2\kappa}}\tfrac{\mathtt{y}_1^2-\mathtt{y}_2^2}{(\mathtt{y}_1^2+\mathtt{y}_2^2)^2}+O(|\ln\varepsilon|^{-\frac12});\quad\mathtt{q}_2(\varphi) =\tfrac{1}{\sqrt{2\kappa}}\tfrac{\mathtt{y}_1\mathtt{y}_2}{(\mathtt{y}_1^2+\mathtt{y}_2^2)^2}+O(|\ln\varepsilon|^{-\frac12}),
       \\ \nonumber \mathtt{q}_3(\varphi) &=-\tfrac{1}{\sqrt{2\kappa}}\tfrac{\mathtt{y}_1\mathtt{y}_2}{(\mathtt{y}_1^2+\mathtt{y}_2^2)^2}+O(|\ln\varepsilon|^{-\frac12})=:-\frac{1}{\sqrt{2\kappa}}\check{\alpha}(\varphi)+O(|\ln\varepsilon|^{-\frac12}),
      \\ \nonumber\mathtt{q}_4(\varphi)&=\tfrac{1}{2\sqrt{2\kappa}}\tfrac{\mathtt{y}_1^2-\mathtt{y}_2^2}{(\mathtt{y}_1^2+\mathtt{y}_2^2)^2}+O(|\ln\varepsilon|^{-\frac12}).
  \end{align*}
  \end{lem}
  As a corollary of Proposition~\ref{prop:G1-trivial}, we obtain the expansion of the renormalized nonlinear functional at $f=0$.
\begin{cor}\label{cor-F0}
The asymptotics of $\mathbf{F}$, defined in \eqref{F-def}, at $f=0$ takes the form
\begin{align*}
  \Pi_{1,{\bf{s}}} {\bf F}(\varepsilon,\mathtt{V}_1,\mathtt{V}_2,0)(\varphi,\theta)&=\varepsilon^2 (2\mathtt{p}_{1,1}(\varphi))^\frac14 \Big({  \, \omega\,\dot{\mathtt{V}}_1(\varphi)-  \mathtt{U}_{-}(\varphi)\,  \mathtt{f}_2(\varphi)}- \varepsilon|\ln\varepsilon|\mathtt{b}(\varphi)\Big)\sin(\theta) \\
&\quad+O\big({ \varepsilon^3}|\ln\varepsilon|^{\frac12}\big)\sin(\theta),
   \\
    \Pi_{1,{\bf{c}}} {\bf F}(\varepsilon,\mathtt{V}_1,\mathtt{V}_2,0)(\varphi,\theta)&=-{\varepsilon^2 \tfrac{1}{2} (2\mathtt{p}_{1,1}(\varphi))^{-\frac14}\mathtt{U}_{-}(\varphi)\,\mathtt{h}_2(\varphi)}\cos(\theta)+O\big({ \varepsilon^3}|\ln\varepsilon|^{\frac12}\big)\cos(\theta),
      \\
      \Pi_1^\perp {\bf F}(\varepsilon,\mathtt{V}_1,\mathtt{V}_2,0)(\varphi,\theta)&=-\varepsilon \mathtt{g}_{3}(\varphi)\sin(3\theta)+\varepsilon^2|\ln\varepsilon|\big(\mathtt{f}_{2}(\varphi)+{\tfrac{\omega \dot{\mathtt{p}}_{1,1}(\varphi)}{4\mathtt{p}_{1,1}(\varphi)}}\big)\cos(2\theta)\\ &\quad + \varepsilon^2|\ln\varepsilon|\mathtt{g}_{2}(\varphi)\sin(2\theta)+O(\varepsilon^2),
\end{align*}
where
\begin{equation}\label{definitionb}
    \mathtt{b}(\varphi):=\frac{1}{16\kappa}\Big(\frac{3}{16 \kappa}+\check{\mathtt{h}}_2(\varphi)\Big).
\end{equation}
Moreover, $\mathtt f_2$ and $\check{\alpha}$ are smooth and odd in $\varphi$, while $\mathtt h_2$, $\check{\mathtt{h}}_2(\varphi)$, $\mathtt g_2$ and $\mathtt g_3$ are smooth and even. In addition, $\check{\alpha}$ and $\check{\mathtt{h}}_2(\varphi)$  have the symmetry property
$$
\check{\alpha}_\star=\check{\alpha}, \qquad \check{\mathtt{h}}_{2,\star}(\varphi)=\check{\mathtt{h}}_2(\varphi).
$$
\end{cor}
\begin{proof}
This follows by applying the change of variables \eqref{def:phi-change} to
Proposition~\ref{prop:G1-trivial} and then dividing by $\sqrt{2\mathtt p_{1,1}}$ as in
\eqref{F-def}. The parity properties follow from the symmetry of the core orbit
(Proposition~\ref{prop-ring}), together with the structure of $G$ in \eqref{eq:G}
and the definition \eqref{list-functions}.
\end{proof}

\section{Toroidal pseudo-differential operators}

The purpose of this section is to collect several tools that will be used repeatedly
in the sequel. We begin by introducing the appropriate function spaces, with special
attention to even/odd symmetries and orthogonality constraints that reflect the
reversible structure of the problem. We then recall some classical properties of
toroidal pseudo-differential operators, including their symbol calculus, operator
topologies, and continuity estimates. In addition, we establish commutator bounds
and smoothing properties, which will later play a crucial role in the Nash--Moser
iteration scheme and in controlling the nonlinear terms of the equation.
This preparatory material is standard in spirit, but we adapt it to the present
setting to ensure that the subsequent analysis is self-contained. We start by fixing
some notation and introducing the function spaces that will serve as the backbone of
our estimates.

 \subsection{ Notation}\label{sec-notat}
Here we list some notation used throughout this paper.

\begin{itemize}[label=$\diamond$]
    \item We denote
$$ 
\qquad\,\mathbb{N}= \{0,1,2,\ldots\},\;\;
\mathbb{N}^*= \mathbb{N}\setminus\{0\},\;\;
\mathbb{Z}=\{\ldots,-1,0,1,\ldots\},\;\;
\mathbb{Z}^*= \mathbb{Z}\setminus\{0\},\;\;
\mathbb{T}= \mathbb{R}/2\pi\mathbb{Z}.
$$
\item We introduce a sequence that will later be used to define cutoff projectors,
which play an essential role in both the reduction procedure and the Nash--Moser scheme:
\begin{equation}\label{definition of Nm}
N_{-1}:= 1,\qquad
\forall n\in\mathbb{N},\qquad
N_{n}:= N_{0}^{\left(\frac{3}{2}\right)^{n}},\qquad
N_0\geqslant 2.
\end{equation}
\item We consider a list of real numbers satisfying
\begin{equation}\label{cond1}
\nu \in(0,1),\qquad 1<\tau\leqslant \tfrac32, \qquad 3<s_0\leqslant s\leqslant s_{\textnormal{up}},
\end{equation}
where $s_{\textnormal{up}}$ is a large exponent
\item For any function $h:\T\to\R$ we define its average by
$$
\langle h\rangle_\theta :=\fint_{\mathbb{T}}h(\theta)\,d\theta
=\frac{1}{2\pi }\int_0^{2\pi}h(\theta)\,d\theta.
$$
\item For any functional $g\mapsto F(g)$ (with $g$ belonging to a set of functions),
we define the difference
$$
\Delta_{12}F:=F(g_1)-F(g_2).
$$
\item For a $2\pi$-periodic function $f:\R\to\R$ we define
$$
\forall \theta\in\R,\quad f_\star(\theta)=f(\theta+\pi).
$$
\end{itemize}

\subsection{Function spaces}\label{SEc-function-sapces}
We introduce various function spaces that will be used frequently throughout this paper.
The first one is the classical Sobolev space $H^{s}(\mathbb{T}\times\T;\C)$, namely the set of
all complex-valued periodic functions $h:\mathbb{T}\times\T\to\C$ such that
$$
h(\varphi,\theta)=\sum_{(l,j)\in \Z^{2}}h_{l,j}\,\mathbf{e}_{l,j}(\varphi,\theta),
\qquad h_{l,j}\in\C,\qquad
\|h\|_{H^s(\T^2)}^2=\sum_{(l,j)\in \Z^{2}}\langle l,j\rangle^{2s}|h_{l,j}|^2,
$$
where
$$
\mathbf{e}_{l,j}(\varphi,\theta):=e^{ i (l\varphi+j\theta)}
\quad\hbox{and}\quad
\langle l,j\rangle:=\sqrt{1+|l|^2+|j|^2}.
$$
For reasons connected with the reversibility of the system we distinguish the subspaces
\begin{align}\label{Hs-even}
H^s_{\textnormal{even}}=\Big\{h\in H^{s}(\mathbb{T}^{2};\RR)\;:\;
h(-\varphi,-\theta)=h(\varphi,\theta),\ \forall (\varphi,\theta)\in\T^{2}\Big\},
\end{align}
and
$$
H^s_{\textnormal{odd}}=\Big\{h\in H^{s}(\mathbb{T}^{2};\RR)\;:\;
h(-\varphi,-\theta)=-h(\varphi,\theta),\ \forall (\varphi,\theta)\in\T^{2}\Big\}.
$$
We denote $H^\infty:= \cap_{s\in\RR} H^s$ and similarly define the subspaces
$H^\infty_{\textnormal{even}}$ and $H^\infty_{\textnormal{odd}}$.

\smallskip
Three additional subspaces will be used repeatedly throughout this paper, namely
$H^{s}_\star(\mathbb{T}^{2},\mathbb{R})$, $H^{s}_{\circ}(\mathbb{T}^{2},\mathbb{R})$ and
$H^{s}_\perp(\mathbb{T}^{2},\mathbb{R})$, defined as follows:
$$
 H^{s}_\star:=H^{s}_\star(\mathbb{T}^{2},\mathbb{R})=\left\{ h\in H^{s}(\mathbb{T}^{2},\mathbb{R})\ :\ \int_{\T^2}h(\varphi,\theta)\cos(\theta)\, d\theta\,d\varphi=0,\ \forall \varphi\in\T,\, \int_{\T}h(\varphi,\theta) d\theta=0\right\},
$$
$$
H^{s}_\circ:=H^{s}_{\circ}(\mathbb{T}^{2},\mathbb{R})
=\left\{ h\in H^{s}_\star(\mathbb{T}^{2},\mathbb{R})\ :\ 
\forall \varphi\in\T,\ \int_{\T}h(\varphi,\theta)\sin(\theta)\, d\theta=0\right\},
$$
and
$$
H^{s}_\perp:=H^{s}_\perp(\mathbb{T}^{2},\mathbb{R})
=\left\{ h\in H^{s}_{\circ}(\mathbb{T}^{2},\mathbb{R})\ :\ 
\forall \varphi\in\T,\ \int_{\T}h(\varphi,\theta)\cos(\theta)\, d\theta=0\right\}.
$$
We also define the subspaces
$$
H^{s}_{\circ,\textnormal{even}}=H^{s}_{\textnormal{even}}\cap H^s_{\circ},
\qquad
H^{s}_{\circ,\textnormal{odd}}=H^{s}_{\textnormal{odd}}\cap H^s_{\circ},
$$
and, similarly,
$$
H^{s}_{\perp,\textnormal{even}}=H^{s}_{\textnormal{even}}\cap H^s_{\perp},
\qquad
H^{s}_{\perp,\textnormal{odd}}=H^{s}_{\textnormal{odd}}\cap H^s_{\perp}.
$$
Note that any element $h\in H^{s}_{\circ}$ decomposes as
\begin{align*}
h(\varphi,\theta)&=h_1(\varphi) \cos(\theta)+\sum_{|j|\geqslant 2}h_j(\varphi) e^{ i j\theta}.
\end{align*}
This leads to the direct sum
$$
H^s_{\circ}=H^{s}_{\circ,\bf{c}}\oplus H^s_\perp.
$$
We introduce the projectors: for $h\in H^s_{\circ}$ and $j\geqslant 1$, 
\begin{align}\label{Pi-proj0}
 \nonumber \Pi_{j,{\bf{c}}}h(\varphi,\theta)&={\frac{1}{\pi}}\int_{\T}h(\varphi,\eta) \cos(j\eta)d\eta\,\, \cos(j\theta),\\
 \Pi_{j,\bf{s}}h(\varphi,\theta)&=\frac{1}{\pi}\int_{\T}h(\varphi,\eta) \sin(j\eta)d\eta\,\, \sin(j\theta).
 \end{align}
Sometimes we use the notation: for $h\in H^s_{\circ}$, 
\begin{align*}
 \Pi h:=\Pi_{1,{\bf{c}}}h\quad\hbox{and} \quad \Pi^\perp h=h-\Pi h.
\end{align*}

To implement the Nash--Moser scheme for the nonlinear equation we need to measure the dependence of solutions with respect to the external parameter $\lambda\in\mathbb{R}$; for this purpose we introduce weighted Sobolev (Lipschitz) spaces.
\begin{def}\label{Def-WS}
Let $\nu\in(0,1), s\in\RR$ and let $\mathcal{O}$ be an open bounded subset of $\mathbb{R}$.
We define
$$
\textnormal{Lip}_\nu(\mathcal{O},H^{s}(\T^{2}))=
\Big\{ h:\mathcal{O}\rightarrow H^{s}(\T^2)\;:\;\|h\|_{s}^{\textnormal{Lip},\nu}<\infty\Big\},
$$
with
\begin{equation*}
\|h\|_{s}^{\textnormal{Lip},\nu}:= \sup_{\lambda\in{\mathcal{O}}}\|h(\lambda,\cdot)\|_{H^{s}(\T^{2})}+\nu\sup_{\lambda_1\neq\lambda_2\in{\mathcal{O}}}\tfrac{\|h(\lambda_1,\cdot)-h(\lambda_2,\cdot)\|_{H^{s-1}(\T^{2})}}{|\lambda_1-\lambda_2|}\cdot
\end{equation*}
\end{def}
Later in the Nash--Moser scheme we will also use the notation
$$\|h\|_{s,\mathcal{O}}^{\textnormal{Lip},\nu}= \|h\|_{s}^{\textnormal{Lip},\nu}.$$
The following notation will be frequently used
\begin{align}\label{compact-notat1}
    \mathbb{X}^s:=\textnormal{Lip}_\nu(\mathcal{O},H^{s}_{\circ,\textnormal{even}}(\T^{2}))\quad\hbox{and}\quad \mathbb{Y}^s:=\textnormal{Lip}_\nu(\mathcal{O},H^{s}_{\circ,\textnormal{odd}}(\T^{2})),
\end{align}
and
\begin{align}  \label{compact-notat2}\mathbb{X}_\perp^s:=\textnormal{Lip}_\nu(\mathcal{O},H^{s}_{\perp,\textnormal{even}}(\T^{2}))\quad\hbox{and}\quad \mathbb{Y}_\perp^s:=\textnormal{Lip}_\nu(\mathcal{O},H^{s}_{\perp,\textnormal{odd}}(\T^{2})).
\end{align}
Next we recall some classical results related to law products and interpolation laws.
\begin{lem}\label{Law-prodX1}
 Let $s\geqslant 1,   s_{0}>2$ and  $f,g\in \textnormal{Lip}_\nu(\mathcal{O},H^{s}(\T^{2}))$. Then:
\begin{enumerate}
    \item 
 The product $fg$ belongs to $ \textnormal{Lip}_\nu(\mathcal{O},H^{s}(\T^{2}))$ and
$$\|fg\|_{s}^{\textnormal{Lip},\nu}\lesssim\|f\|_{s_{0}}^{\textnormal{Lip},\nu}\|g\|_{s}^{\textnormal{Lip},\nu}+\|f\|_{s}^{\textnormal{Lip},\nu}\|g\|_{s_{0}}^{\textnormal{Lip},\nu}.$$
\item Let $n\in\N$. Then
\begin{align*}
\|\partial_\theta^n(fg)\|_{s}^{\textnormal{Lip},\nu}\lesssim\;&
\|f\|_{s_{0}+n}^{\textnormal{Lip},\nu}\|g\|_{s}^{\textnormal{Lip},\nu}
+\|f\|_{s+n}^{\textnormal{Lip},\nu}\|g\|_{s_{0}}^{\textnormal{Lip},\nu}\\
&+\|f\|_{s_{0}}^{\textnormal{Lip},\nu}\|\partial_\theta^ng\|_{s}^{\textnormal{Lip},\nu}
+\|f\|_{s}^{\textnormal{Lip},\nu}\|\partial_\theta^ng\|_{s_{0}}^{\textnormal{Lip},\nu}.
\end{align*}
\item  Let $1<s_{1}\leqslant s_{3}\leqslant s_{2}$ and $\overline{\theta}\in[0,1],$ with  $s_{3}=\overline{\theta} s_{1}+(1-\overline{\theta})s_{2}.$			If $f\in \textnormal{Lip}_\nu(\mathcal{O},H^{s_2}(\T^{2})))$, then  $\rho\in \textnormal{Lip}_\nu(\mathcal{O},H^{s_3}(\T^{2})))$ and
				$$\|f\|_{s_{3}}^{\textnormal{Lip},\nu}\lesssim\left(\|f\|_{s_{1}}^{\textnormal{Lip},\nu}\right)^{\overline{\theta}}\left(\|f\|_{s_{2}}^{\textnormal{Lip},\nu}\right)^{1-\overline{\theta}}.$$

\end{enumerate}
\end{lem}
\begin{proof}  
The first and second results are standard and whose proof can be found in \cite{BertiMontalto} and the references therein.
The interpolation result uses Kato-Ponce inequality, see \cite{DongLi}.
\end{proof}
Next, we introduce a cutoff frequency operator which has the advantage of providing smoothing effects that will be used in KAM and Nash-Moser schemes.
For $N\in\mathbb{N}^{*}$ we define the orthogonal projections on $H^s(\T^{2})$ by
$$
\Pi_{N} h=\sum_{\substack{(l,j)\in\mathbb{Z}\times\mathbb{Z}\\ \langle l,j\rangle\leqslant N}}
h_{l,j}\,\mathbf{e}_{l,j},
\qquad
\Pi^{\perp}_{N}=\textnormal{Id}-\Pi_{N}.
$$
The following result is elementary and can be  checked by a straightforward computation.
\begin{lem}\label{orthog-Lem1}
Let $N\in\mathbb{N}^{*}, s\in\RR$ and $\mu\in\mathbb{R}_{+}.$ Then 
\begin{equation*}
\|\Pi_{N}h\|_{H^{s+\mu}}\leqslant N^{\mu}\|h\|_{H^{s}},\quad \|\Pi_{N}^{\perp}h\|_{H^{s}}\leqslant N^{-\mu}\|h\|_{H^{s+\mu}},
\end{equation*}
and
\begin{equation*}
\|\Pi_{N}h\|_{s+\mu}^{\textnormal{Lip},\nu}\leqslant N^{\mu}\|h\|_{s}^{\textnormal{Lip},\nu},\quad \|\Pi_{N}^{\perp}h\|_{s}^{\textnormal{Lip},\nu}\leqslant N^{-\mu}\|h\|_{s+\mu}^{\textnormal{Lip},\nu}.
\end{equation*}
\end{lem}
The proof of the following lemma can be found in \cite[Lemma A.2]{HR22}:
\begin{lem}\label{lemma HR22}
Let $(\nu,s_{0},s)$ satisfy \eqref{cond1} and let
$f\in \textnormal{Lip}_\nu(\mathcal{O},H^{s}(\mathbb{T}^{2}))$.
Define $g:\mathcal{O}\times\mathbb{T}^3\rightarrow\mathbb{C}$ by
	$$g(\lambda,\varphi,\theta,\eta)=\left\lbrace\begin{array}{ll}
		\tfrac{f(\lambda,\varphi,\eta)-f(\lambda,\varphi,\theta)}{{\tan}(\frac{\eta-\theta}{2})} & \textnormal{if }\theta\neq \eta,\\
		2\partial_{\theta}f(\lambda,\varphi,\theta)& \textnormal{if }\theta=\eta.
	\end{array}\right.$$
	Then 
$$\|g\|_{s}^{\textnormal{Lip},\nu}\lesssim\|f\|_{s+1}^{\textnormal{Lip},\nu}.$$ 
\end{lem}
We point out that in Lemma \ref{lemma HR22}, the function $g$ depends on three variables
$\varphi$, $\theta$, and $\eta$, and is $2\pi$-periodic in each of them. The corresponding
Sobolev spaces are therefore those defined on the torus $\mathbb{T}^3$.
\subsection{Operators and Reversibility}
In this section, we collect some notations for $\lambda-$dependent families of operators and recall the notions of reversibility that will be used throughout the paper.

\medskip

Let $\mathcal{O}\subset\R$ be an open bounded set and let
$$
\mathcal{A}:\mathcal{O}\to \mathcal{L}\big(H^s(\T^{2});H^s(\T^{2})\big),\qquad
\lambda\mapsto \mathcal{A}(\lambda),
$$
be a (say, smooth) family of bounded linear operators on $H^s(\T^2)$.

\smallskip
\noindent\textbf{Matrix representation.}
With respect to the Fourier basis $\mathbf{e}_{l,j}(\varphi,\theta):=e^{ i(l\varphi+j\theta)}$,
each $\mathcal{A}(\lambda)$ can be represented by an infinite matrix
$\big(\mathcal{A}^{\,l,j}_{l_0,j_0}(\lambda)\big)_{(l,j),(l_0,j_0)\in\Z^2}$ defined by
$$
\mathcal{A}(\lambda)\mathbf{e}_{l_0,j_0}
=\sum_{(l,j)\in\Z^2}\mathcal{A}^{\,l,j}_{l_0,j_0}(\lambda)\,\mathbf{e}_{l,j}.
$$
We use the normalized $L^2(\T^2)$ inner product
$$
\langle u,v\rangle_{L^2(\T^2)}:=\frac{1}{(2\pi)^2}\int_{\T^2}u(\varphi,\theta)\,\overline{v(\varphi,\theta)}\,d\varphi\,d\theta,
$$
so that, since $\overline{\mathbf{e}_{l,j}}=\mathbf{e}_{-l,-j}$,
\begin{align*}
\mathcal{A}^{\,l,j}_{l_0,j_0}(\lambda)
&=\big\langle \mathcal{A}(\lambda)\mathbf{e}_{l_0,j_0},\mathbf{e}_{l,j}\big\rangle_{L^{2}(\T^{2})} \\
&=\frac{1}{(2\pi)^{2}}
\int_{\T^{2}}
\big(\mathcal{A}(\lambda)\mathbf{e}_{l_0,j_0}\big)(\varphi,\theta)\,
\mathbf{e}_{-l,-j}(\varphi,\theta)\,d\varphi\, d\theta.
\end{align*}
\smallskip
Throughout this paper, operators and functions may depend on the same external parameter $\lambda$.
Accordingly, if $h\in \textnormal{Lip}_{\nu}(\mathcal{O},H^{s}(\T^{2}))$ and
$\mathcal{A}=(\mathcal{A}(\lambda))_{\lambda\in\mathcal{O}}$, we define the pointwise action by
$$
(\mathcal{A}h)(\lambda,\varphi,\theta):=\big(\mathcal{A}(\lambda)\,h(\lambda,\cdot)\big)(\varphi,\theta),
$$
i.e. both objects are evaluated at the same $\lambda$. This is justified by the fact that  with respect to the  equations that we consider here, $\lambda$ is an external parameter in the model.

\smallskip
We now collect some definitions and properties concerning reversible operators; see for instance
\cite{Baldi-berti,BFM21,BFM,BertiMontalto}. Recall the involution $\mathscr{S}$ introduced in
\eqref{def-S}:
$$
(\mathscr{S}h)(\varphi,\theta):=h(-\varphi,-\theta),\qquad (\varphi,\theta)\in\T^2.
$$
\begin{defi}\label{Def-Rev}
Let $\mathcal{A}=(\mathcal{A}(\lambda))_{\lambda\in\mathcal{O}}$ be a $\lambda$--dependent family of operators.
We say that $\mathcal{A}$ is
\begin{enumerate}[label=\roman*)]
\item \emph{real} if, for every $\lambda\in\mathcal{O}$ and every $h\in L^2(\T^2,\C)$,
$$
\overline{h}=h \quad\Longrightarrow\quad \overline{\mathcal{A}(\lambda)h}=\mathcal{A}(\lambda)h.
$$
\item \emph{reversible} if, for every $\lambda\in\mathcal{O}$,
$$
\mathcal{A}(\lambda)\circ\mathscr{S}=-\,\mathscr{S}\circ\mathcal{A}(\lambda).
$$
\item \emph{reversibility preserving} if, for every $\lambda\in\mathcal{O}$,
$$
\mathcal{A}(\lambda)\circ\mathscr{S}=\mathscr{S}\circ\mathcal{A}(\lambda).
$$
\end{enumerate}
\end{defi}
\subsection{Symbol class topology}\label{Se-Toroidal pseudo-differential operators}

In this subsection we recall some classical facts about pseudo-differential operators in the periodic setting. Let $\mathcal{A}:C^\infty(\T)\to C^\infty(\T)$ be a linear operator. Its (toroidal) symbol
$\sigma_{\mathcal{A}}$ is defined by
$$
\forall\,(\theta,\xi)\in\T\times\Z,\qquad
\sigma_{\mathcal{A}}(\theta,\xi):=\mathbf{e}_{-\xi}(\theta)\,\big(\mathcal{A}\mathbf{e}_{\xi}\big)(\theta),
\qquad \mathbf{e}_{\xi}(\theta):=e^{ i \xi\theta}.
$$
Thus, if $h(\theta)=\sum_{\xi\in\Z}h_\xi e^{ i \xi\theta}$, then
$$
\mathcal{A}h(\theta)=\sum_{\xi\in\Z}\sigma_{\mathcal{A}}(\theta,\xi)\,h_\xi\,e^{ i \xi\theta}.
$$
The operator also admits the kernel representation
$$
\mathcal{A}h(\theta)=\int_{\T}K(\theta,\eta)\,h(\eta)\,d\eta,
\qquad
K(\theta,\eta):=\frac{1}{2\pi}\sum_{\xi\in\Z}\sigma_{\mathcal{A}}(\theta,\xi)\,e^{ i(\theta-\eta)\xi}.
$$
Moreover, the symbol can be recovered from the kernel by Fourier inversion:
$$
\sigma_{\mathcal{A}}(\theta,\xi)=\int_{\T}K(\theta,\theta+\eta)\,e^{ i \eta\xi}\,d\eta.
$$
\smallskip
\noindent\textbf{Difference operators.}
Given a discrete function $h:\Z\to\C$, we define the forward difference operator $\Delta_\xi$ by
$$
\Delta_\xi h(\xi):=a(\xi+1)-a(\xi),\qquad
\Delta_\xi^\ell:=\underbrace{\Delta_\xi\circ\cdots\circ\Delta_\xi}_{\ell\ \text{times}}.
$$
The backward difference $\overline{\Delta}_\xi$ is defined by
$$
\overline{\Delta}_\xi h(\xi):=h(\xi)-h(\xi-1),\qquad
\overline{\Delta}_\xi^\ell:=\underbrace{\overline{\Delta}_\xi\circ\cdots\circ\overline{\Delta}_\xi}_{\ell\ \text{times}}.
$$
We recall some useful identities (see e.g. \cite[Chapter~3]{Ruzhansky}). First, the discrete Leibniz rule:
$$
\forall\, \gamma\in\N,\quad\Delta_\xi^\gamma(fg)(\xi)=\sum_{0\leqslant \beta\leqslant\gamma}\left(_\beta^\gamma\right)\Delta_\xi^\beta f(\xi)\,\Delta_\xi^{\gamma-\beta}g(\xi+\beta).
$$
Second, summation by parts:
$$
\sum_{\xi\in\Z} f(\xi)\,\Delta_\xi^\gamma g(\xi)
=(-1)^{\gamma}\sum_{\xi\in\Z}\big(\overline{\Delta}_\xi^\gamma f\big)(\xi)\,g(\xi).
$$
Finally, by induction we can check that
$$
\Delta_\xi^\gamma e^{ i \xi \eta}=\big(e^{ i \eta}-1\big)^\gamma e^{ i \xi \eta}.
$$
\smallskip
\noindent\textbf{Symbols with parameters.}
We will frequently consider families of operators $\varphi\mapsto \mathcal{A}(\varphi)$ acting on the $\theta$--variable. If
$$
h(\varphi,\theta)=\sum_{\xi\in\Z}h_\xi(\varphi)\,e^{ i \xi\theta},
$$
we define the (one--parameter) symbol $\sigma_{\mathcal{A}}(\varphi,\theta,\xi)$ by requiring
\begin{equation*}
\mathcal{A}(\varphi)h(\varphi,\theta)
=\sum_{\xi\in\Z}\sigma_{\mathcal{A}}(\varphi,\theta,\xi)\,h_\xi(\varphi)\,e^{ i \xi\theta}.
\end{equation*}
The symbol can be recovered from the operator as follows
$$
\sigma_{\mathcal{A}}(\varphi,\theta,\xi)
:=\mathbf{e}_{-\xi}(\theta)\,\big(\mathcal{A}(\varphi)\mathbf{e}_{\xi}\big)(\theta),
\qquad \mathbf{e}_\xi(\theta)=e^{ i \xi\theta},
$$
and $\mathcal{A}(\varphi)$ admits the kernel representation
$$
\mathcal{A}(\varphi)h(\varphi,\theta)=\int_{\T} K(\varphi,\theta,\eta)\,h(\varphi,\eta)\,d\eta,
\qquad
K(\varphi,\theta,\eta):=\frac{1}{2\pi}\sum_{\xi\in\Z}\sigma_{\mathcal{A}}(\varphi,\theta,\xi)\,e^{ i(\theta-\eta)\xi}.
$$
By the Fourier inversion formula we get
$$
\sigma_{\mathcal{A}}(\varphi,\theta,\xi)=\int_{\T}K(\varphi,\theta,\theta+\eta)\,e^{ i \eta\xi}\,d\eta.
$$
Let $s,m\in\R$ and $\gamma\in\N$. Writing $\langle\xi\rangle:=\max\{1,|\xi|\}$, we define the norm over the symbol class
\begin{equation}\label{Def-Norm-M1}
\interleave \mathcal{A}\interleave_{m,s,\gamma}
:=\sup_{\xi\in\Z}\ \sup_{0\leqslant \ell\leqslant\gamma}
\langle \xi\rangle^{-m+\ell}\,
\big\|\Delta_\xi^\ell\sigma_{\mathcal{A}}(\cdot,\cdot,\xi)\big\|_{H^{s}(\T^{2})}.
\end{equation}
\smallskip
We also consider families $(\lambda,\varphi)\in \mathcal{O}\times\T\mapsto \mathcal{A}(\lambda,\varphi)$.
If
$$
h(\lambda,\varphi,\theta)=\sum_{\xi\in\Z}h_\xi(\lambda,\varphi)\,e^{ i \xi\theta},
$$
we define the multi-parameter symbol of $\mathcal{A}$ by
\begin{equation}\label{A-phi-lambda}
\mathcal{A}(\lambda,\varphi)h(\lambda,\varphi,\theta)
=\sum_{\xi\in\Z}\sigma_{\mathcal{A}}(\lambda,\varphi,\theta,\xi)\,h_\xi(\lambda,\varphi)\,e^{ i \xi\theta},
\end{equation}
and the kernel representation becomes
\begin{equation}\label{kernel-phi-lambda}
\mathcal{A}(\lambda,\varphi)h(\lambda,\varphi,\theta)=\int_{\T}K(\lambda,\varphi,\theta,\eta)\,h(\lambda,\varphi,\eta)\,d\eta,
\end{equation}
with
$$
K(\lambda,\varphi,\theta,\eta)=\frac{1}{2\pi}\sum_{\xi\in\Z}e^{ i(\theta-\eta)\xi}\,\sigma_{\mathcal{A}}(\lambda,\varphi,\theta,\xi).
$$
Given $m,s\in\R$, $\gamma\in\N$ and $\nu\in(0,1)$, we define the weighted norm
$$
\interleave \mathcal{A}\interleave_{m,s,\gamma}^{\textnormal{Lip},\nu}
:= \sup_{\lambda\in\mathcal{O}}\interleave \mathcal{A}(\lambda,\cdot)\interleave_{m,s,\gamma}
+\nu\sup_{\lambda_1\neq\lambda_2\in\mathcal{O}}
\frac{\interleave \mathcal{A}(\lambda_1,\cdot)-\mathcal{A}(\lambda_2,\cdot)\interleave_{m,s-1,\gamma}}{|\lambda_1-\lambda_2|}\cdot
$$
We say that $\mathcal{A}\in \textnormal{OPS}^m$ if for every $s\in\R$ and $\gamma\in\N$,
$$
\interleave \mathcal{A}\interleave_{m,s,\gamma}^{\textnormal{Lip},\nu}<\infty.
$$
In the particular case $m=\gamma=0$ we write
$$
\interleave \mathcal{A}\interleave_{s}^{\textnormal{Lip},\nu}
:= \interleave \mathcal{A}\interleave_{0,s,0}^{\textnormal{Lip},\nu}.
$$

\noindent {\bf Continuity and law products.}
We recall the following classical result, whose proof can be found, for instance, in \cite{BertiMontalto,Ruzhansky}.
	
\begin{lem}\label{Lem-Rgv1LM}
 Let $\mathcal{A}$ be a pseudo-differential operator as in \eqref{A-phi-lambda} and $s\geqslant 1, s_0>\frac{5}{2}$. Then 
$$
\|\mathcal{A}h\|_{s}^{\textnormal{Lip},\nu}\lesssim \interleave \mathcal{A}\interleave_{s_0}^{\textnormal{Lip},\nu} \|h\|_{s}^{\textnormal{Lip},\nu}+\interleave\mathcal{A}\interleave_{s}^{\textnormal{Lip},\nu} \|h\|_{s_0}^{\textnormal{Lip},\nu}.
$$
\end{lem}
The following lemma shows that an integral operator with a smooth kernel is smoothing of arbitrary order.
\begin{lem}[Lemma~2.32 in~\cite{BertiMontalto}]\label{lem:Int}
Assume that $K(\lambda,\cdot)\in C^\infty(\T^3)$ for every $\lambda\in\mathcal{O}$, and let $\mathcal{A}$ be the integral operator defined as in \eqref{kernel-phi-lambda}, namely
$$
\big(\mathcal{A}(\lambda,\varphi)h\big)(\lambda,\varphi,\theta)
=\int_{\T}K(\lambda,\varphi,\theta,\eta)\,h(\lambda,\varphi,\eta)\,d\eta .
$$
Then $\mathcal{A}\in \mathrm{OPS}^{-\infty}$. Moreover, for every $m\in\R$, $\gamma\in\N$, $s\geqslant 1$, and $s_0>\tfrac{5}{2}$, there exists a constant $C=C(m,s,\gamma)>0$ such that
$$
\interleave \mathcal{A} \interleave_{-m,s,\gamma}^{\textnormal{Lip},\nu}
\le C(m,s,\gamma)\,\|K\|_{s+m+\gamma}^{\textnormal{Lip},\nu}.
$$
\end{lem}

We also recall the following composition and commutator estimates (see, e.g., \cite{BertiMontalto}).

\begin{lem}\label{comm-pseudo1}
Let $\mathcal{A}$ and $\mathcal{B}$ be as in \eqref{A-phi-lambda}, and assume that the parameters satisfy \eqref{cond1}. Let $m_1,m_2\in\R$ and $\gamma\in\N$. Then
\begin{align*}
\interleave \mathcal{A}\mathcal{B}\interleave_{m_1+m_2,s,\gamma}^{\textnormal{Lip},\nu}
\lesssim\;&
\interleave \mathcal{A}\interleave_{m_1,s,\gamma}^{\textnormal{Lip},\nu}\,
\interleave \mathcal{B}\interleave_{m_2,s_0+\gamma+|m_1|,\gamma}^{\textnormal{Lip},\nu}\\
&+
\interleave \mathcal{A}\interleave_{m_1,s_0,\gamma}^{\textnormal{Lip},\nu}\,
\interleave \mathcal{B}\interleave_{m_2,s+\gamma+|m_1|,\gamma}^{\textnormal{Lip},\nu}.
\end{align*}
In particular (for $m_1=m_2=\gamma=0$) one has
$$
\interleave \mathcal{A}\mathcal{B}\interleave_{s}^{\textnormal{Lip},\nu}
\lesssim
\interleave \mathcal{A}\interleave_{s_0}^{\textnormal{Lip},\nu}\,
\interleave \mathcal{B}\interleave_{s}^{\textnormal{Lip},\nu}
+\interleave \mathcal{A}\interleave_{s}^{\textnormal{Lip},\nu}\,
\interleave \mathcal{B}\interleave_{s_0}^{\textnormal{Lip},\nu}.
$$
Moreover, the commutator $[\mathcal{A},\mathcal{B}]:=\mathcal{A}\mathcal{B}-\mathcal{B}\mathcal{A}$ belongs to
$\mathrm{OPS}^{m_1+m_2-1}$ and satisfies
\begin{align*}
\interleave [\mathcal{A},\mathcal{B}]\interleave_{m_1+m_2-1,s,\gamma}^{\textnormal{Lip},\nu}
\lesssim\;&
\interleave \mathcal{A}\interleave_{m_1,s+2+m_2+\gamma,\gamma+1}^{\textnormal{Lip},\nu}\,
\interleave \mathcal{B}\interleave_{m_2,s_0+\gamma+|m_1|,\gamma+1}^{\textnormal{Lip},\nu}\\
&+
\interleave \mathcal{A}\interleave_{m_1,s_0+2+m_2+\gamma,\gamma+1}^{\textnormal{Lip},\nu}\,
\interleave \mathcal{B}\interleave_{m_2,s+\gamma+|m_1|,\gamma+1}^{\textnormal{Lip},\nu}.
\end{align*}
\end{lem} 
\subsection{Inhomogeneous Fourier multipliers}
We shall consider integral operators of the form
\begin{equation*}
(\mathcal{A}h)(\lambda,\varphi,\theta)
:=\int_{\T}K(\lambda,\varphi,\theta,\eta)\,h(\lambda,\varphi,\eta)\,d\eta,
\qquad
K(\lambda,\varphi,\theta,\eta)=f(\lambda,\varphi,\theta)\,\mathtt{K}(\theta-\eta),
\end{equation*}
where $f$ is smooth and $\mathtt{K}\in L^1(\T)$.
Equivalently, 
$$
(\mathcal{A}h)(\lambda,\varphi,\theta)
=f(\lambda,\varphi,\theta)\,(\mathtt{K}*h)(\lambda,\varphi,\theta),
\qquad
(\mathtt{K}*h)(\lambda,\varphi,\theta):=\int_{\T}\mathtt{K}(\theta-\eta)\,h(\lambda,\varphi,\eta)\,d\eta.
$$
In this case, the symbol factorizes as
\begin{align*}
\sigma_{\mathcal{A}}(\lambda,\varphi,\theta,\xi)
&=\int_{\T}K(\lambda,\varphi,\theta,\theta+\eta)\,e^{ i \eta\xi}\,d\eta
=f(\lambda,\varphi,\theta)\,\widehat{\mathtt{K}}(\xi),
\end{align*}
where $\widehat{\mathtt{K}}(\xi)$ denotes the $\xi$-th Fourier coefficient of $\mathtt{K}$.
Therefore,
\begin{equation*}
\interleave \mathcal{A}\interleave_{m,s,\gamma}
=\|f\|_{H^s(\T^2)}\,
\sup_{\xi\in\Z}\ \sup_{0\leqslant\ell\leqslant\gamma}
\langle \xi\rangle^{-m+\ell}\,|\Delta_\xi^\ell\widehat{\mathtt{K}}(\xi)|.
\end{equation*}
We now briefly recall the Hilbert transform on $\T$.
Let $h:\T\to\R$ be continuous with zero average. We define
\begin{equation}\label{Hilbert1}
(\mathcal{H}h)(\theta):=\fint_{\T}h(\eta)\,\cot\!\Big(\tfrac{\eta-\theta}{2}\Big)\,d\eta.
\end{equation}
It is well known that $\mathcal{H}$ is a Fourier multiplier and
$$
\forall j\in\Z^*,\qquad
\mathcal{H}\mathbf{e}_j(\theta)= i\,\operatorname{sign}(j)\,\mathbf{e}_j(\theta).
$$
Next, we  introduce auxiliary operators naturally arising in the linearization and conjugation procedures.
For $m\in\N$ we define the (homogeneous) Fourier multiplier
\begin{equation}\label{FM-n-1}
(\Lambda_m h)(\lambda,\varphi,\theta)
:=-\fint_{\T}\Big|\sin\!\Big(\tfrac{\theta-\eta}{2}\Big)\Big|^m
\ln\!\Big|\sin\!\Big(\tfrac{\theta-\eta}{2}\Big)\Big|\,
h(\lambda,\varphi,\eta)\,d\eta,
\end{equation}
together with its inhomogeneous variant
\begin{equation}\label{eq:Lambda-m-f}
(\Lambda_{m,f} h)(\lambda,\varphi,\theta)
:=-\fint_{\T}\Big|\sin\!\Big(\tfrac{\theta-\eta}{2}\Big)\Big|^m
\ln\!\Big|\sin\!\Big(\tfrac{\theta-\eta}{2}\Big)\Big|\,
f(\lambda,\varphi,\theta,\eta)\,h(\lambda,\varphi,\eta)\,d\eta,
\end{equation}
where $f$ is a smooth function.

\begin{lem}\label{Lem-homo-inhomo}
Let $m,\gamma\in\N$ and $s\in\R$. Then, the following assertions hold.
\begin{enumerate}
\item If $m$ is even, then $\Lambda_m\in \textnormal{OPS}^{-m-1}$ and
$$
\interleave\Lambda_m\interleave_{-m-1,s,\gamma}<\infty.
$$
If $m$ is odd, then for every $\varepsilon>0$ one has $\Lambda_m\in \textnormal{OPS}^{-m-1+\varepsilon}$ and
$$
\interleave\Lambda_m\interleave_{-m-1+\varepsilon,s,\gamma}<\infty.
$$

\item Assume in addition that $f\in \textnormal{Lip}_\nu(\mathcal{O},H^{s+m+\gamma+3}(\T^3))$.
If $m$ is even, then $\Lambda_{m,f}\in \textnormal{OPS}^{-m-1}$ and
$$
\interleave \Lambda_{m,f}\interleave_{-m-1,s,\gamma}^{\textnormal{Lip},\nu}
\lesssim \|f\|_{s+m+\gamma+3}^{\textnormal{Lip},\nu}.
$$
If $m$ is odd, then for every $\varepsilon>0$ one has $\Lambda_{m,f}\in \textnormal{OPS}^{-m-1+\varepsilon}$ and
$$
\interleave \Lambda_{m,f}\interleave_{-m-1+\varepsilon,s,\gamma}^{\textnormal{Lip},\nu}
\lesssim \|f\|_{s+m+\gamma+3}^{\textnormal{Lip},\nu}.
$$
\end{enumerate}
\end{lem}

\begin{proof}
We first study the Fourier coefficients of the kernel. Define, for $m\in\N$ and $n\in\Z$,
\begin{equation*}
I_{m,n}:=-\fint_{\T}\Big|\sin\!\Big(\tfrac{\eta}{2}\Big)\Big|^{m}
\ln\!\Big|\sin\!\Big(\tfrac{\eta}{2}\Big)\Big|\,e^{ i n\eta}\,d\eta.
\end{equation*}
Since the integrand is even, we have
$$
I_{m,n}=-\frac{1}{\pi}\int_{0}^{\pi}\sin^{m}\!\Big(\tfrac{\eta}{2}\Big)
\ln\!\Big(\sin\!\big(\tfrac{\eta}{2}\big)\Big)\cos(n\eta)\,d\eta.
$$
For $\Re(\alpha)>-1$, a classical formula (see \cite[(3.631)]{Gradshteyn}) yields
$$
F(\alpha,n):=\frac{1}{\pi}\int_{0}^\pi\Big|\sin\!\Big(\tfrac{\eta}{2}\Big)\Big|^{\alpha}e^{ i n\eta}\,d\eta
=\frac{(-1)^n}{2^{\alpha}}
   \frac{\Gamma(\alpha+1)}
        {\Gamma\!\big(\frac{\alpha}{2}-n+1\big)\Gamma\!\big(\frac{\alpha}{2}+n+1\big)}.
$$
Since
$$
\frac{\partial}{\partial\alpha}\Big|\sin\!\Big(\tfrac{\eta}{2}\Big)\Big|^{\alpha}
=
\Big|\sin\!\Big(\tfrac{\eta}{2}\Big)\Big|^{\alpha}\ln\Big|\sin\!\Big(\tfrac{\eta}{2}\Big)\Big|,
$$
we obtain
$$
I_{m,n}=-\partial_\alpha F(m,n).
$$
From the explicit form of $F(\alpha,n)$ we infer
$$
\partial_\alpha\ln F(\alpha,n)
=-\ln2+\psi(\alpha+1)
-\tfrac{1}{2}\psi\!\Big(\tfrac{\alpha}{2}-n+1\Big)
-\tfrac{1}{2}\psi\!\Big(\tfrac{\alpha}{2}+n+1\Big),
$$
where $\psi(x)=\Gamma'(x)/\Gamma(x)$ is the digamma function. Hence,
\begin{equation}\label{form:Inm}
\begin{aligned}
I_{m,n}
&=-\frac{(-1)^n}{2^{m}}
  \frac{\Gamma(m+1)}
       {\Gamma(\tfrac{m}{2}-n+1)\Gamma(\tfrac{m}{2}+n+1)}\\
&\quad\times
\Big[
-\ln2+\psi(m+1)
-\tfrac12\psi\!\Big(\tfrac{m}{2}-n+1\Big)
-\tfrac12\psi\!\Big(\tfrac{m}{2}+n+1\Big)
\Big].
\end{aligned}
\end{equation}

Using the reflection formula
$$
\frac{1}{\Gamma(\tfrac{\alpha}{2}-n+1)}
=\frac{\sin\!\big(\pi(n-\tfrac{\alpha}{2})\big)}{\pi}\,
  \Gamma\!\Big(n-\tfrac{\alpha}{2}\Big)
=\frac{(-1)^n\sin(\pi\alpha/2)}{\pi}\,
  \Gamma\!\Big(n-\tfrac{\alpha}{2}\Big),
$$
we obtain
$$
F(\alpha,n)
=\frac{2^{-\alpha}\Gamma(\alpha+1)\sin(\pi\alpha/2)}{\pi}\,
  \frac{\Gamma(n-\alpha/2)}{\Gamma(n+\alpha/2+1)}\cdot
$$
For large $n\to\infty$, we have the asymptotics
$$
\frac{\Gamma(n-\alpha/2)}{\Gamma(n+\alpha/2+1)}
\sim n^{-\alpha-1},
$$
and therefore
$$
F(\alpha,n)\sim C(\alpha)\,n^{-\alpha-1},
\qquad
C(\alpha):=\frac{2^{-\alpha}\Gamma(\alpha+1)\sin(\pi\alpha/2)}{\pi}\cdot
$$
Differentiating the asymptotic form gives
$$
I_{m,n}\sim \big(C(m)\ln n-C'(m)\big)\,n^{-m-1}.
$$
Moreover,
$$
\frac{C'(\alpha)}{C(\alpha)}
=-\ln 2+\psi(\alpha+1)+\frac{\pi}{2}\cot\!\Big(\frac{\pi\alpha}{2}\Big).
$$
For odd $m=2k+1$ one obtains, as $n\to\infty$,
\begin{align}\label{odd-asy}
 I_{m,n}
&\sim
\tfrac{\Gamma(2k+2)(-1)^k}{\pi\, 2^{2k+1}}\,
n^{-2k-2}
\left[
\ln n+\ln 2-\psi(2k+2)
\right]
\sim
\frac{\Gamma(2k+2)(-1)^{k}}{\pi\, 2^{2k+1}}\,
\frac{\ln n}{n^{2k+2}}\cdot
\end{align}
If $m=2k$ is even, then $\sin(\pi m/2)=0$ and the leading term vanishes; using $C'(2k)$ yields
\begin{align}\label{even-asy}
I_{m,n}
\sim
2^{-(2k+1)}\Gamma(2k+1)(-1)^k\,
n^{-(2k+1)},
\qquad n\to\infty.
\end{align}
We now estimate the symbol of $\Lambda_{m,f}$. Writing $\eta\mapsto \theta+\eta$ in \eqref{eq:Lambda-m-f},
\begin{align*}
\sigma_{\Lambda_{m,f}}(\lambda,\varphi,\theta,n)
&=-\fint_{\T}\Big|\sin\!\Big(\tfrac{\eta}{2}\Big)\Big|^m
\ln\!\Big|\sin\!\Big(\tfrac{\eta}{2}\Big)\Big|\,
f(\lambda,\varphi,\theta,\theta+\eta)\,e^{ i n\eta}\,d\eta\\
&=f(\lambda,\varphi,\theta,\theta)\,I_{m,n}
\\ &\quad-\fint_{\T}\Big|\sin\!\Big(\tfrac{\eta}{2}\Big)\Big|^m
\ln\!\Big|\sin\!\Big(\tfrac{\eta}{2}\Big)\Big|\,
\big(f(\lambda,\varphi,\theta,\theta+\eta)-f(\lambda,\varphi,\theta,\theta)\big)\,e^{ i n\eta}\,d\eta\\
&=: \sigma_{1,m}(\lambda,\varphi,\theta,n)+\sigma_{2,m}(\lambda,\varphi,\theta,n).
\end{align*}
Introduce
$$
g(\lambda,\varphi,\theta,\eta)
:=\frac{f(\lambda,\varphi,\theta,\theta+\eta)-f(\lambda,\varphi,\theta,\theta)}{\sin(\eta/2)}\cdot
$$
Then $g$ is smooth and
\begin{equation*}
\sigma_{2,m}(\lambda,\varphi,\theta,n)
=-\fint_{\T}\Big|\sin\!\Big(\tfrac{\eta}{2}\Big)\Big|^m
\sin\!\Big(\tfrac{\eta}{2}\Big)\,
\ln\!\Big|\sin\!\Big(\tfrac{\eta}{2}\Big)\Big|\,
g(\lambda,\varphi,\theta,\eta)\,e^{ i n\eta}\,d\eta.
\end{equation*}
Using \eqref{odd-asy} and \eqref{even-asy} and standard bounds on discrete differences of sequences with polynomial/logarithmic decay, we infer by straightforward estimates (combined with Sobolev embeddings) that:
\begin{enumerate}[label=\roman*)]
\item If $m$ is even,
$$
\sup_{n\in\Z}\ \sup_{0\le\ell\le\gamma}
\langle n\rangle^{m+1+\ell}\,
\big\|\Delta_n^\ell\sigma_{1,m}(\cdot,\cdot,n)\big\|_{H^{s}(\T^{2})}
\lesssim \|f\|_{H^{s+1}(\T^{3})}.
$$
\item If $m$ is odd, then for any $\varepsilon>0$,
$$
\sup_{n\in\Z}\ \sup_{0\le\ell\le\gamma}
\langle n\rangle^{m+1-\varepsilon+\ell}\,
\big\|\Delta_n^\ell\sigma_{1,m}(\cdot,\cdot,n)\big\|_{H^{s}(\T^{2})}
\lesssim \|f\|_{H^{s+1}(\T^{3})}.
$$
\end{enumerate}
On the other hand, using $\Delta_n^\ell e^{ i n\eta}=(e^{ i \eta}-1)^\ell e^{ i n\eta}$ and
$\partial_\eta^{m+1+\ell}e^{ i n\eta}=( i n)^{m+1+\ell}e^{ i n\eta}$, we write
\begin{align*}
&(i n)^{m+1+\ell}\Delta_n^\ell\,\sigma_{2,m}(\lambda,\varphi,\theta,n)=\\
&-\fint_{\T}\Big|\sin\!\Big(\tfrac{\eta}{2}\Big)\Big|^m
\sin\!\Big(\tfrac{\eta}{2}\Big)\,
\ln\!\Big|\sin\!\Big(\tfrac{\eta}{2}\Big)\Big|\,
g(\lambda,\varphi,\theta,\eta)\,(e^{ i \eta}-1)^\ell\,
\partial_\eta^{m+1+\ell}\!\big[e^{ i n\eta}\big]\,d\eta.
\end{align*}
Integrating by parts $(m+1+\ell)$ times yields
$$
\langle n\rangle^{m+1+\ell}\,
\big\|\Delta_n^\ell\sigma_{2,m}(\cdot,\cdot,n)\big\|_{H^s(\T^2)}
\lesssim \|g\|_{H^{s+m+\ell+1}(\T^3)}
\lesssim \|f\|_{H^{s+m+\ell+3}(\T^3)},
$$
where in the last estimate we used that $g$ is a difference quotient and Sobolev embeddings.
Collecting the bounds for $\sigma_{1,m}$ and $\sigma_{2,m}$ gives the claimed estimates in $\textnormal{OPS}^{-m-1}$ (if $m$ is even) and $\textnormal{OPS}^{-m-1+\varepsilon}$ (if $m$ is odd).
\end{proof}
Given a smooth function $f:\T^2\to\R$ and an integer $m\geqslant 0$, we define the operator
\begin{equation}\label{Hf-n}
\mathcal{H}_{f,-m}:=\partial_\theta\big(f\Lambda_m+\Lambda_m f\big),
\end{equation}
where products denote composition, right or left, of $\Lambda_m$ with the multiplication operator by $f$.

\begin{lem}\label{Commutators-hilbert}
Let $\nu\in(0,1)$ and $s_0>2$ be as in \eqref{cond1}, and let $f_1,f_2$ be two smooth (possibly $\lambda$--dependent) functions on $\T^2$. Then, the following assertions hold.
\begin{enumerate}
\item The operators $f_1\mathcal{H}$ and $\mathcal{H}_{f_1,0}$ belong to $\textnormal{OPS}^0$ and, for all
$s\in\R$ and $\gamma\in\N$,
$$
\interleave f_1\mathcal{H}\interleave_{0,s,\gamma}^{\textnormal{Lip},\nu}
+\interleave \mathcal{H}_{f_1,0}\interleave_{0,s,\gamma}^{\textnormal{Lip},\nu}
\lesssim \|f_1\|_{s}^{\textnormal{Lip},\nu}.
$$

\item The commutator $[f_1\mathcal{H},f_2\mathcal{H}]$ belongs to $\textnormal{OPS}^{-\infty}$ and, for all
$m,\gamma\in\N$ and $s>1$,
$$
\interleave[f_1\mathcal{H},f_2\mathcal{H}]\interleave_{-m,s,\gamma}^{\textnormal{Lip},\nu}
\lesssim
\|f_1\|_{s_0+\gamma+m+1}^{\textnormal{Lip},\nu}\,\|f_2\|_{s+\gamma+m+1}^{\textnormal{Lip},\nu}
+\|f_1\|_{s+\gamma+m+1}^{\textnormal{Lip},\nu}\,\|f_2\|_{s_0+\gamma+m+1}^{\textnormal{Lip},\nu}.
$$

\item The commutators $[\mathcal{H},\mathcal{H}_{f_1,0}]\in\textnormal{OPS}^{-\infty}$ and
$[\mathcal{H}_{f_1,0},\mathcal{H}_{f_2,0}]\in\textnormal{OPS}^{-3}$, with
$$
\interleave[\mathcal{H},\mathcal{H}_{f_1,0}]\interleave_{-m,s,\gamma}^{\textnormal{Lip},\nu}
\lesssim \|f_1\|_{s+\gamma+m+1}^{\textnormal{Lip},\nu},
$$
and
$$
\interleave[\mathcal{H}_{f_1,0},\mathcal{H}_{f_2,0}]\interleave_{-3,s,\gamma}^{\textnormal{Lip},\nu}
\lesssim
\|f_1\|_{s+\gamma+7}^{\textnormal{Lip},\nu}\,\|f_2\|_{s_0+\gamma+7}^{\textnormal{Lip},\nu}
+\|f_1\|_{s_0+\gamma+7}^{\textnormal{Lip},\nu}\,\|f_2\|_{s+\gamma+7}^{\textnormal{Lip},\nu}.
$$

\item Let $\Pi:H^s(\T^2)\to H^s(\T^2)$ be a  finite-rank operator, in the variable $\theta$, and $\mathcal{A}\in\textnormal{OPS}^m$ for some
$m\in\R$. Then $\Pi\mathcal{A}$ and $\mathcal{A}\Pi$ belong to $\textnormal{OPS}^{-\infty}$.
\end{enumerate}
\end{lem}

\begin{proof}
In the proofs below, to simplify the notation, we will omit the variable $\varphi$, since the operators act only on the spatial variable $\theta$.\\
\medskip\noindent\textbf{(1)} Since $\mathcal{H}$ is a Fourier multiplier with bounded symbol, it belongs to
$\textnormal{OPS}^0$. Multiplication by $f_1$ also belongs to $\textnormal{OPS}^0$, hence
$f_1\mathcal{H}\in\textnormal{OPS}^0$ and the stated bound follows from \eqref{Def-Norm-M1}.
Moreover, $\Lambda_0\in\textnormal{OPS}^{-1}$ by Lemma~\ref{Lem-homo-inhomo} (case $m=0$), hence
$\mathcal{H}_{f_1,0}=\partial_\theta(f_1\Lambda_0+\Lambda_0 f_1)\in \textnormal{OPS}^0$.

\medskip\noindent\textbf{(2)} We first show that $[f_1,\mathcal{H}]\in \textnormal{OPS}^{-\infty}$.
Indeed, for $h=h(\theta)$,
$$
[f_1,\mathcal{H}]h(\theta)
=\mathrm{p.v.}\,\fint_{\T}\frac{f_1(\theta)-f_1(\eta)}{\tan\!\big(\frac{\eta-\theta}{2}\big)}\,h(\eta)\,d\eta
=\fint_{\T}K(\theta,\eta)\,h(\eta)\,d\eta,
$$
where the quotient defining $K(\theta,\eta)$ extends smoothly across $\eta=\theta$.
By Lemma~\ref{lemma HR22} we have
$$
\|K\|_{s}^{\textnormal{Lip},\nu}\lesssim\|f_1\|_{s+1}^{\textnormal{Lip},\nu},
$$
and Lemma~\ref{lem:Int} implies $[f_1,\mathcal{H}]\in\textnormal{OPS}^{-\infty}$ together with
\begin{equation}\label{commu-Hilb-f}
\forall m,\gamma\in\N,\ \forall s>1,\qquad
\interleave [f_1,\mathcal{H}]\interleave_{-m,s,\gamma}^{\textnormal{Lip},\nu}
\lesssim \|f_1\|_{s+m+\gamma+1}^{\textnormal{Lip},\nu}.
\end{equation}
Now, using
$$
[f_1\mathcal{H},f_2\mathcal{H}]
=f_1[\mathcal{H},f_2]\mathcal{H}-f_2[\mathcal{H},f_1]\mathcal{H},
$$
and applying \eqref{commu-Hilb-f} together with Lemma~\ref{comm-pseudo1}, we conclude that
$[f_1\mathcal{H},f_2\mathcal{H}]\in\textnormal{OPS}^{-\infty}$ and obtain the stated estimate.

\medskip\noindent\textbf{(3)} By definition,
\begin{equation}\label{id-Hf0}
\mathcal{H}_{f_1,0}
=\tfrac12 f_1\mathcal{H}+\tfrac12 \mathcal{H}f_1+(\partial_\theta f_1)\Lambda_{0}.
\end{equation}
Hence,
$$
[\mathcal{H},\mathcal{H}_{f_1,0}]
=\tfrac12[\mathcal{H}, f_1\mathcal{H}]+\tfrac12[\mathcal{H}, \mathcal{H}f_1]
+[\mathcal{H},(\partial_\theta f_1)]\Lambda_{0},
$$
where we used that $\mathcal{H}$ and $\Lambda_0$ commute (both are Fourier multipliers).
Since $[\mathcal{H},f_1]\in\textnormal{OPS}^{-\infty}$ by \eqref{commu-Hilb-f}, we infer
$[\mathcal{H},\mathcal{H}_{f_1,0}]\in\textnormal{OPS}^{-\infty}$ with the claimed bound.
\\
Next, from \eqref{id-Hf0} we also have
$$
\mathcal{H}_{f_1,0}
=f_1\mathcal{H}-\tfrac12[f_1,\mathcal{H}]+(\partial_\theta f_1)\Lambda_{0}.
$$
Using again that $[f_1,\mathcal{H}]\in\textnormal{OPS}^{-\infty}$, we obtain
$$
[\mathcal{H}_{f_1,0},\mathcal{H}_{f_2,0}]
=[f_1\mathcal{H},(\partial_\theta f_2)\Lambda_0]-[(\partial_\theta f_1)\Lambda_0,f_2\mathcal{H}]
+[(\partial_\theta f_1)\Lambda_0,(\partial_\theta f_2)\Lambda_0]+\textnormal{OPS}^{-\infty}.
$$
Moreover,
$$
[f_1\mathcal{H},(\partial_\theta f_2)\Lambda_0]
=(\partial_\theta f_2)[f_1,\Lambda_0]\mathcal{H}+\textnormal{OPS}^{-\infty},
\qquad
[(\partial_\theta f_1)\Lambda_0,f_2\mathcal{H}]
=(\partial_\theta f_1)[f_2,\Lambda_0]\mathcal{H}+\textnormal{OPS}^{-\infty},
$$
so that
$$
[\mathcal{H}_{f_1,0},\mathcal{H}_{f_2,0}]
=\Big((\partial_\theta f_2)[f_1,\Lambda_0]-(\partial_\theta f_1)[f_2,\Lambda_0]\Big)\mathcal{H}
+[(\partial_\theta f_1)\Lambda_0,(\partial_\theta f_2)\Lambda_0]+\textnormal{OPS}^{-\infty}.
$$
A direct computation shows that
$$
\Big((\partial_\theta f_2)[f_1,\Lambda_0]-(\partial_\theta f_1)[f_2,\Lambda_0]\Big)[h]
=\Lambda_{2,K}[h],
$$
where $\Lambda_{2,K}$ is defined as in \eqref{eq:Lambda-m-f} with $m=2$ and coefficient $K=K(\theta,\eta)$ given by
$$
K(\theta,\eta)=
\begin{cases}
\dfrac{(\partial_\theta f_2)(\theta)\big(f_1(\theta)-f_1(\eta)\big)-(\partial_\theta f_1)(\theta)\big(f_2(\theta)-f_2(\eta)\big)}
{\sin^2\!\big(\frac{\theta-\eta}{2}\big)},
& \theta\neq \eta,\\
2\Big((\partial_\theta f_1)(\partial_\theta^2 f_2)-(\partial_\theta f_2)(\partial_\theta^2 f_1)\Big),
& \theta=\eta.
\end{cases}
$$
The function $K(\theta,\eta)$ is smooth and satisfies
$$
\|K\|_{s}^{\textnormal{Lip},\nu}
\lesssim
\|f_1\|_{s_0+2}^{\textnormal{Lip},\nu}\|f_2\|_{s+2}^{\textnormal{Lip},\nu}
+\|f_1\|_{s+2}^{\textnormal{Lip},\nu}\|f_2\|_{s_0+2}^{\textnormal{Lip},\nu}.
$$
Therefore, Lemma~\ref{Lem-homo-inhomo} (with $m=2$) yields $\Lambda_{2,K}\in\textnormal{OPS}^{-3}$, hence
$\Lambda_{2,K}\mathcal{H}\in\textnormal{OPS}^{-3}$. Finally, since $(\partial_\theta f_i)\Lambda_0\in\textnormal{OPS}^{-1}$,
Lemma~\ref{comm-pseudo1} implies
$$
[(\partial_\theta f_1)\Lambda_0,(\partial_\theta f_2)\Lambda_0]\in\textnormal{OPS}^{-3}.
$$
This proves $[\mathcal{H}_{f_1,0},\mathcal{H}_{f_2,0}]\in\textnormal{OPS}^{-3}$ and the stated bound.

\medskip\noindent\textbf{(4)} If $\Pi$ has finite rank, then $\Pi$ is smoothing of arbitrary order; hence
$\Pi\mathcal{A}$ and $\mathcal{A}\Pi$ have smooth kernels and belong to $\textnormal{OPS}^{-\infty}$.
\end{proof}

\subsection{A fundamental Lemma}\label{Sec-Fundam-lemma}
In this section we establish an auxiliary result that will play a key role in the subsequent analysis. The lemma addresses the solvability and compactness properties of a class of integral operators that naturally arise in the treatment of the degeneracy of the first Fourier mode. This result will be used at two different stages of the argument: first, in the construction of an approximate solution, and second, in the study of the invertibility of the linearized operator restricted to the first mode.\\
Throughout this section, we consider three non-zero, smooth, periodic, even, real-valued functions $\varrho_j:\T\to\R$, $j=1,2,3$, satisfying the symmetry condition
$$
\varrho_{j,\star}=\varrho_j,\quad\textnormal{where}\quad \varrho_{j,\star}(\varphi):=\varrho_{j}(\varphi+\pi). 
$$
These functions will act as coefficients in a family of nested integral operators that encode the structure of the reduced equations.
Given a parameter $b\in\R$, we introduce the linear operators $\mathscr{P}$ and $\mathscr{T}$ defined below:
$$
\mathscr{P}[g](\varphi):=\varrho_1(\varphi)\left[b+\int_0^\varphi\varrho_2(\tau)\left(\int_0^\tau\varrho_3(s) g(s)ds\right) d\tau  \right],
$$
and
$$
 \mathscr{T}[g](\varphi):=\varrho_1(\varphi)\int_0^\varphi\varrho_2(\tau)\left(\int_0^\tau\varrho_3(s) g(s)ds\right) d\tau.
 $$
 The operator $\mathscr{P}$ incorporates a constant contribution together with a double integral term, while $\mathscr{T}$ corresponds to the purely integral component. 
To capture the symmetry properties of the problem, we work in the subspace
$$
H^s_{\textnormal{even},*}(\T)
:=\big\{ g\in H^s_{\textnormal{even}}(\T,\R),\ g_\star=-g \big\},
$$
which consists of even periodic functions in the Sobolev space $H^s(\T,\R)$ satisfying the additional compatibility condition imposed by the $\star$-symmetry. Equivalently, these are functions whose Fourier expansion contains only odd cosine modes, namely
$$
\forall \varphi\in\R,\qquad 
g(\varphi)=\sum_{n\in\N} a_n \cos\big((2n+1)\varphi\big).
$$
The following lemma gathers several structural properties of the operators $\mathscr{P}$ and $\mathscr{T}$. In particular, we identify the compatibility condition on the parameter $b$ ensuring that $\mathscr{P}$ preserves the space $H^s_{\textnormal{even},*}(\T)$, establish compactness properties of the operators involved, and derive an invertibility criterion for the operator $\mathrm{Id}-\mathscr{P}$. These results will be very useful in Section \ref{sec:approx} and Section \ref{section-Inver-Mode1}. 
\begin{lem}\label{lem-fundamental}
    The following assertions hold true. Let  $s\geqslant 0.$ Then
    \begin{enumerate}
        \item The operator $\mathscr{P}:H^s_{\textnormal{even},*}(\T)\to H^s_{\textnormal{even},*}(\T)$ is well-defined if and only if
        \begin{align}\label{b-choice}
       b= -\frac12\int_0^{\pi}\varrho_2(\tau)\left(\int_0^\tau\varrho_3(s) g(s)ds\right)d\tau.\end{align}
       \item The operator $\mathscr{T}:L^\infty([0,2\pi];\R)\to L^\infty([0,2\pi];\R)$ is compact and  $\textnormal{Id}-\mathscr{T}$ has a bounded inverse.
        \item With the assumption \eqref{b-choice}, the operator $\mathscr{P}:H^s_{\textnormal{even},*}(\T)\to H^s_{\textnormal{even},*}(\T)$ is a compact operator. In addition,  assume that
        $$
 2+\int_0^{\pi}\varrho_2(\tau)\left(\int_0^\tau\varrho_3(s) (\textnormal{Id}-\mathscr{T})^{-1}[\varrho_1](s)ds\right)d\tau\neq0,
 $$
 then, the operator $\textnormal{Id}-\mathscr{P}$ has a bounded inverse $(\textnormal{Id}-\mathscr{P})^{-1}:H^s_{\textnormal{even},*}(\T)\to H^s_{\textnormal{even},*}(\T)$.
    \end{enumerate}
\end{lem}
\begin{proof}
    {\bf{(1)}} Let us first show that $\mathscr{P}$ is well-defined. Consider  $g\in H^s_{\textnormal{even},*}(\T,\R)$. We intend to show that the map 
    $$F:\tau\mapsto \int_0^\tau\varrho_3(s) g(s)ds, 
    $$ is $2\pi-$periodic. As the function  $\varrho_3 g$ is $2\pi-$periodic, then $F$ is $2\pi-$periodic if and only if
    $$
    F(2\pi)=0.
    $$
    From the assumptions, we have
    $$
    \varrho_{3,\star}(s)g_\star(s)=-\varrho_{3}(s)g(s).
    $$
    Then the Fourier decomposition of $\varrho_{3}g$  takes the form
    \begin{align}\label{Phi3-dec}
    \varrho_{3}(s)g(s)=\sum_{n\in\mathbb{N}}a_n \cos\big((2n+1) s\big).
    \end{align}
    Integrating this sum yields  
$$
F(2\pi)=\int_{0}^{2\pi}\sum_{n\in\mathbb{N}}a_n \cos\big((2n+1) s\big)ds=0.
$$
On the other hand, since $\varrho_{3}g$ is even, we deduce that $F$ is an odd periodic function. Consequently, as $\varrho_2$ is an even periodic function, the product $\varrho_2 F$ is also odd and periodic. Integrating then yields
$$
\varphi\mapsto \int_0^\varphi \varrho_2(\tau)\,F(\tau)d\tau
$$
is an even $2\pi-$periodic function. Consequently, we obtain that for each $b\in\R$ and $g\in H^s_{\textnormal{even},*}(\T)$, the function $\mathscr{P}[g]$ is also an even $2\pi-$periodic function. Now, let us check
the condition 
$$(\mathscr{P}[g])_\star=-\mathscr{P}[g].$$
Since $\varrho_{j,\star}=\varrho_j$, we can write
\begin{align}\label{Ss-star}
   (\mathscr{P}[g])_\star(\varphi)
   =\varrho_1(\varphi)\left[b+\int_0^{\varphi+\pi}\varrho_2(\tau)F(\tau)\, d\tau  \right].
\end{align}
Performing the change of variables $\tau\rightsquigarrow\tau+\pi$ yields
\begin{align}\label{FF-st}
\int_0^{\varphi+\pi}\varrho_2(\tau)F(\tau)\, d\tau
  = \int_{-\pi}^{\varphi}\varrho_2(\tau)F_\star(\tau)\, d\tau .
\end{align}
Moreover, a direct computation using the change of variables $s\rightsquigarrow s+\pi$ gives
\begin{align*}
F_\star(\tau) &=\int_0^{\tau+\pi}\varrho_3(s) g(s)ds =\int_{-\pi}^{\tau}\varrho_{3,\star}(s)g_\star(s) ds\\
&=-\int_{-\pi}^{\tau}\varrho_{3}(s)g(s) ds=-\int_{0}^{\tau}\varrho_{3}(s)g(s) ds-\int_{-\pi}^{0}\varrho_{3}(s)g(s) ds\\
&=-F(\tau)-\int_{-\pi}^{0}\varrho_{3}(s)g(s) ds.
\end{align*}
Applying \eqref{Phi3-dec} yields
$$
\int_{-\pi}^{0}\varrho_{3}(s)g(s) ds=0.
$$
Hence
$$
F_\star=-F.
$$
Plugging this identity into \eqref{FF-st} leads to 

\begin{align*}
\int_0^{\varphi+\pi}\varrho_2(\tau)F(\tau) d\tau & = -\int_{-\pi}^{\varphi}\varrho_2(\tau)F(\tau) d\tau\\
&=-\int_{-\pi}^{0}\varrho_2(\tau)F(\tau) d\tau-\int_{0}^{\varphi}\varrho_2(\tau)F(\tau) d\tau.
\end{align*}
Since $F$ is odd and $\varrho_2$ is even, the product $\varrho_2 F$ is odd. In particular, it has zero average over a period, and therefore
$$
\int_{-\pi}^{0}\varrho_2(\tau)F(\tau)\,d\tau
=-\int_{0}^{\pi}\varrho_2(\tau)F(\tau)\,d\tau.
$$
Consequently,
\begin{align*}
\int_0^{\varphi+\pi}\varrho_2(\tau)F(\tau) d\tau 
&=\int_0^{\pi}\varrho_2(\tau)F(\tau) d\tau-\int_{0}^{\varphi}\varrho_2(\tau)F(\tau) d\tau.
\end{align*}
Inserting this relation into \eqref{Ss-star} implies
\begin{align*}
   (\mathscr{P}[g])_\star(\varphi)=\varrho_1(\varphi)\left[b+\int_0^{\pi}\varrho_2(\tau)F(\tau) d\tau-\int_0^{\varphi}\varrho_2(\tau)F(\tau) d\tau  \right]. 
\end{align*}
Then 
\begin{align*}
   (\mathscr{P}[g])_\star(\varphi)=-(\mathscr{P}[g])(\varphi)\Longleftrightarrow \varrho_1(\varphi)\left[b+\int_0^{\pi}\varrho_2(\tau)F(\tau) d\tau  \right] =-b\, \varrho_1(\varphi).
\end{align*}
As $\varrho_1$ is not identically zero, then this is  equivalent to
\begin{align*}
   (\mathscr{P}[g])_\star=-\mathscr{P}[g]\Longleftrightarrow b=-\frac12\int_0^{\pi}\varrho_2(\tau)F(\tau) d\tau.
\end{align*}
As for the regularity property, the fact that $g\in H^s(\T)$ implies $\mathscr{P}[g]\in H^s(\T)$ follows readily from the product rules in Sobolev spaces and the smoothing effect of the integral operators. This is a classical fact, and we omit the details.
\\
{\bf{(2)}} We can easily show that
\begin{align*}
    \|\mathscr{T}[g]\|_{L^\infty}\lesssim \, \|g\|_{L^\infty}\prod_{j=1}^3\|\varrho_j\|_{L^\infty}\,.
\end{align*}
and by differentiation
\begin{align*}
    \|(\mathscr{T}[g])^\prime\|_{L^\infty}\lesssim \, \|g\|_{L^\infty}\|\varrho_1\|_{\textnormal{Lip}}\prod_{j=2}^3\|\varrho_j\|_{L^\infty}\,.
\end{align*}
Therefore,
\begin{align*}
    \|\mathscr{T}[g]\|_{\textnormal{Lip}}\lesssim \, \|g\|_{L^\infty}\|\varrho_1\|_{\textnormal{Lip}}\prod_{j=2}^3\|\varrho_j\|_{L^\infty}\,.
\end{align*}
Here the Lipschitz nor is given by
$$
\|g\|_{\textnormal{Lip}}=\|g\|_{L^\infty}+\sup_{x\neq y}\tfrac{|g(x)-g(y)}{|x-y|}\cdot
$$
It follows that $\mathscr{T}: L^\infty([0,2\pi],\R)\to \textnormal{Lip}([0,2\pi],\R)$ is continuous. Since the embedding  $$\textnormal{Lip}([0,2\pi],\R)\hookrightarrow L^\infty([0,2\pi],\R)
$$ is compact, we deduce that $\mathscr{T}: L^\infty([0,2\pi],\R)\to L^\infty([0,2\pi],\R)$ is compact.\\ It remains to show that $\textnormal{Id}- \mathscr{T}:L^\infty([0,2\pi],\R)\to L^\infty([0,2\pi],\R)$ has a bounded inverse. Since $\mathscr{T}$ is compact, then according to Fredholm alternative the inverse $(\hbox{Id}-\mathscr{T})^{-1}:L^\infty([0,2\pi];\R)\to L^\infty([0,2\pi];\R)$ exists and is  bounded iff the kernel of $\hbox{Id}-\mathscr{T}$ is trivial. To this end, assume that a function $g\in L^\infty([0,2\pi];\R)$ satisfies
 \begin{align}\label{kernel-gL}
g(\varphi)=\varrho_1(\varphi)\int_0^\varphi\varrho_2(\tau)\left(\int_0^\tau\varrho_3(s) g(s)ds\right) d\tau.
 \end{align}
 We claim that, for any $n\in\N$ 
 $$
 \forall \varphi\in[0,2\pi],\quad |g(\varphi)|\leqslant \|g\|_{L^\infty} \frac{C^n}{(2n)!} \varphi^{2n}\quad\hbox{with}\quad C:=\prod_{j=1}^3\|\varrho_j\|_{L^\infty}.
 $$
 We prove this estimate  by induction. The case $n=0$ is immediate. Assume that the bound holds for some $n\in\N$ and let us show that it also holds for $n+1$. Using \eqref{kernel-gL} we find
 \begin{align*}
|g(\varphi)|&\leqslant C \int_0^\varphi\left(\int_0^\tau |g(s)|ds\right) d\tau\\
&\leqslant C \|g\|_{L^\infty} \frac{C^n}{(2n)!} \int_0^\varphi\left(\int_0^\tau s^{2n}ds\right) d\tau\\
&\leqslant  \|g\|_{L^\infty} \frac{C^{n+1}}{(2n+2)!}\varphi^{2n+2}\,.
 \end{align*}
 This completes the induction step and hence the proof of the claim. Letting $n\to\infty$ in the above estimate, we obtain $g\equiv 0$. 
\\
 {\bf{(3)}} 
 Note that
$
\mathscr{P}=\varrho_1\,\mathscr{T}.$
Then, using the first and the second point, together with the fact that multiplication by the smooth function $\varrho_1$ is bounded on $H^s(\T)$, we infer that
$$
\mathscr{P}:H^s_{\textnormal{even},*}(\T)\to H^s_{\textnormal{even},*}(\T)
$$
is a compact operator.
 Thus $\textnormal{Id}- \mathscr{P}:H^s_{\textnormal{even},*}(\T)\to H^s_{\textnormal{even},*}(\T)$ is Fedholm and it is invertible if and only if its kernel is trivial. We will study the kernel equation and assume that we have a function $g\in H^s_{\textnormal{even},*}(\T)$ such that
 \begin{align}\label{Eq-gg}
 \forall \varphi\in[0,2\pi],\quad g(\varphi)&=\varrho_1(\varphi)\left[b+\int_0^\varphi\varrho_2(\tau)\left(\int_0^\tau\varrho_3(s) g(s)ds\right) d\tau  \right]\\
 \nonumber b&=-\frac12\int_0^{\pi}\varrho_2(\tau)F(\tau) d\tau.
 \end{align}
 Now, coming back to the equation \eqref{Eq-gg}
which can be written in the form
\begin{align*}
 \forall \varphi\in[0,2\pi],\quad (\hbox{Id}-\mathscr{T})[g](\varphi)&=b\,\varrho_1(\varphi),  \quad b=-\frac12\int_0^{\pi}\varrho_2(\tau)F(\tau) d\tau.
 \end{align*}
 Thus
 \begin{align}\label{Eq-bno}
 \forall \varphi\in[0,2\pi],\quad g(\varphi)&=b\,(\hbox{Id}-\mathscr{T})^{-1}[\varrho_1](\varphi)  \\
 \nonumber b&=-\frac12\int_0^{\pi}\varrho_2(\tau)\left(\int_0^\tau\varrho_3(s) g(s)ds\right)d\tau.
 \end{align}
 Therefore the compatibility condition on $b$ becomes
 $$
 b\left[1+\frac12\int_0^{\pi}\varrho_2(\tau)\left(\int_0^\tau\varrho_3(s) (\hbox{Id}-\mathscr{T})^{-1}[\varrho_1](s)ds\right)d\tau\right]=0.
 $$
 It follows that under the assumption
 $$
 2+\int_0^{\pi}\varrho_2(\tau)\left(\int_0^\tau\varrho_3(s) (\hbox{Id}-\mathscr{T})^{-1}[\varrho_1](s)ds\right)d\tau\neq0
 $$
 we get $b=0$. This yields in view of \eqref{Eq-bno} that $g\equiv 0$. Hence the kernel equation \eqref{Eq-gg} is trivial.
 THis completes the proof of the fundamental lemma.
 \end{proof}
\section{Linearization}\label{section-linearization}
As stated in Lemma \ref{lem-red}, in the symmetric setting the nonlinear dynamics of two
concentrated vortex rings reduce to the Hamiltonian equation \eqref{F-def}.
To understand the behavior of this reduced system under small perturbations,
we study its linearization, which provides valuable insight into the dynamics. This is a crucial step: the
linearized operator captures the leading-order action of the nonlinearity,
reveals the mechanisms driving stability and bifurcation, and organizes the
hierarchy of $\varepsilon$--dependent contributions. In particular, it
separates the dominant singular part (generated by the kernel $G$) from
lower-order terms, where logarithmic corrections interact with algebraic
scales.
The main difficulty is that $\mathbf{F}$ admits a delicate asymptotic
expansion near the small state. The singular structure of $G$, combined with
the geometry encoded by $\mathcal{Z}$ and the $\varepsilon$--dependent weights,
produces a genuinely nonlocal operator that must be decomposed carefully.

We begin by studying the regularity and deriving the linearized operator and 
for convenience, let us recall the basic unknowns:
$$
\mathtt{P}_1=\big(\mathtt{p}_{1,1},\mathtt{p}_{1,2}\big),\quad w=\sqrt{1+2\varepsilon f}, \quad \mathcal{Z}(\theta)= (2\mathtt{p}_{1,1})^{\frac14}\cos(\theta)+{ i}(2\mathtt{p}_{1,1})^{-\frac14}\sin(\theta).
$$
For $r>0, s\in\R,$ define the closed ball
\begin{align}\label{ball-space}
B_{r,\textnormal{even}}(s):=\big\{f\in \textnormal{Lip}_\nu(\mathcal{O},H^{s}_{\star,\textnormal{even}}(\T^2)),\,\|f\|_s^{\textnormal{Lip},\nu}\leqslant r \big\}.
\end{align}

\begin{lem}\label{prop:asymp-lin0}
Assume that \eqref{cond1} holds, and use the notation introduced in
\eqref{ball-space}--\eqref{F-def}. Then there exists $\varepsilon_0>0$ such that,
for every $\varepsilon\in(0,\varepsilon_0)$, the following statements hold.
\begin{enumerate}
\item The map $${\bf F}(\varepsilon,\mathtt{V}_1,\mathtt{V}_2,\cdot):B_{1,\textnormal{even}}(s_0)\cap \textnormal{Lip}_\nu(\mathcal{O},H^{s}_{\star,\textnormal{even}}(\T^2))\to \textnormal{Lip}_\nu(\mathcal{O},H^{s-1}_{\star,\textnormal{odd}}(\T^2))$$ is well-defined.
\item The map $f\mapsto {\bf F}(\varepsilon,\mathtt{V}_1,\mathtt{V}_2,f)$ is  of class $C^1$. More precisely, for every
$$
f\in B_{1,\mathrm{even}}(s_0)
\quad\text{and}\quad
h\in \mathrm{Lip}_\nu(\mathcal{O},H^s_{\star,\mathrm{even}}(\T^2)),
$$
one has  
 \begin{align*}
\nonumber  \partial_f {\bf F}(\varepsilon,\mathtt{V}_1,\mathtt{V}_2,f)[h] (\varphi,\theta)&=\varepsilon^3 |\ln\varepsilon | \omega \partial_\varphi h(\varphi ,\theta) +\varepsilon\,    \,\partial_\theta\big\{ h(\varphi,\theta)U(\varepsilon,f)(\varphi,\theta) \big\}\\ &\quad +\tfrac{\varepsilon}{\sqrt{\mathtt{p}_{1,1}(\varphi)}}  \,\big(\mathcal{I}_1(\varepsilon,f)[h]+\mathcal{I}_2(\varepsilon,f)[h]\big)(\varphi,\theta).
\end{align*}
Here, the operators $\mathcal{I}_1(f)$ and $\mathcal{I}_2(f)$ are defined by
 \begin{align*}
\qquad \mathcal{I}_1(\varepsilon,f)[h](\varphi,\theta)&:=\partial_\theta \fint_{\mathbb{T}} G\left(\Gamma_1(f)(\varphi,\theta,1), \Gamma_1(f)(\varphi,\eta,1)\right) h(\varphi,\eta)d\eta,\\
\qquad \mathcal{I}_2(\varepsilon,f)[h](\varphi,\theta)&:=\partial_\theta\fint_{\T}G\left(\Gamma_1(f)(\varphi,\theta,1)+\tfrac{i\varepsilon \mathtt{U}_{-}(\varphi)}{ |\ln\varepsilon|} ,\Gamma_2(f)(\varphi,\eta,1)\right)h_\star(\varphi,\eta)  d\eta, 
\end{align*}
while the function $U(\varepsilon,f)(\varphi,\theta)$ is given by
  \begin{equation*}
\begin{aligned}
    U(\varepsilon,f)(\varphi,\theta)&:=\tfrac{1}{\sqrt{\mathtt{p}_{1,1}(\varphi)}}\big( U_1(\varepsilon,f)+ U_2(\varepsilon,f)\big)(\varphi,\theta)-\varepsilon |\ln \varepsilon| \omega  \dot{\mathtt{P}_1}(\varphi)\cdot i\frac{\mathcal{Z}(\varphi,\theta)}{w(\varphi,\theta)}\\ & \qquad  -i\varepsilon^2 \omega\,  \dot{\mathtt{V}}_1(\varphi)\cdot i\frac{\mathcal{Z}(\varphi,\theta)}{w(\varphi,\theta)} +\varepsilon^2 |\ln\varepsilon |\,  \frac{ \omega \,\dot{\mathtt{p}}_{1,1}(\varphi)}{4\mathtt{p}_{1,1}(\varphi)}\sin(2\theta),
    \end{aligned}
\end{equation*}
with
\begin{align*}
\quad\qquad    U_1(\varepsilon,f)(\varphi,\theta)
&:= \varepsilon\!
\fint_{\T}\!\int_{0}^{1}\!
\nabla_1G\left(\Gamma_1(f)(\varphi,\theta,1),
\Gamma_1(f)(\varphi,\eta,\rho)\right)\cdot \mathcal{Z}(\varphi,\theta)\tfrac{w^2(\varphi, \eta)}{w(\varphi,\theta)}
\rho d\rho d\eta,\\
\quad\qquad U_2(\varepsilon,f)(\varphi,\theta)
&:= \varepsilon\!\fint_{\T}\!\int_0^{1}\!\nabla_1G\left(\Gamma_1(f)(\varphi,\theta,1)+\tfrac{i\varepsilon\, \mathtt{U}_{-}}{ |\ln\varepsilon|} ,\Gamma_2(f)(\varphi,\eta,\rho)\right)\cdot \mathcal{Z}(\varphi,\theta)\tfrac{w^2(\varphi, \eta)}{w(\varphi,\theta)} d\rho d\eta,
\end{align*}
and 
\begin{equation}\label{def:Gamma12}
\begin{aligned}
\Gamma_1(f)(\varphi,\eta,\rho)&:=\mathtt{P}_1(\varphi)+\varepsilon\rho\, w(\varphi,\eta)\mathcal{Z}(\varphi,\eta),\\
\Gamma_2(f)(\varphi,\eta,\rho) &:=\mathtt{P}_2(\varphi)+\varepsilon\rho\, w_\star(\varphi,\eta)\mathcal{Z}_\star(\varphi,\eta).
\end{aligned}
\end{equation}
\end{enumerate}

\end{lem}
\begin{proof}
${\bf{(1)}}$ Let's first check the symmetry. For $f\in B_{1,\textnormal{even}}(s),$ it is easy to check from \eqref{F-def} that an integration by parts allows to get
$$\int_{\mathbb{T}}{\bf F}(\varepsilon,\mathtt{V}_1,\mathtt{V}_2,f)(\varphi,\theta)d\theta=\omega\,\varepsilon^2|\ln\varepsilon|\partial_\varphi\int_{\mathbb{T}} f(\varphi,\theta)d\theta=0.$$
On the other hand, the equation is reversible in the sense that for $f\in H^s_{\textnormal{even}}$, we have that ${\bf F}(\varepsilon,\mathtt{V}_1,\mathtt{V}_2,f)$ is odd, that is,
$$
{\bf F}(\varepsilon,\mathtt{V}_1,\mathtt{V}_2,f)(-\varphi,-\theta)=-{\bf F}(\varepsilon,\mathtt{V}_1,\mathtt{V}_2,f)(\varphi,\theta).
$$
Therefore
$$
\int_{\mathbb{T}}{\bf F}(\varepsilon,\mathtt{V}_1,\mathtt{V}_2,f)(\varphi,\theta)\sin(\theta)d\varphi d\theta=0.
$$
The regularity property ${\bf F}(\varepsilon,\mathtt{V}_1,\mathtt{V}_2,f)\in H^{s-1}$ and ${\bf F}(\varepsilon,\mathtt{V}_1,\mathtt{V}_2,f)\in\textnormal{Lip}_\nu(\mathcal{O},H^{s-1}_{\star,\textnormal{odd}}(\T^2))$ can be done using classical tools of products law and composition.

\medskip

\noindent
${\bf{(2)}}$ Differentiating \eqref{F-def} with respect to $f$ in the direction $h$ gives
\begin{align}\label{diff-F-def}
 & \partial_f {\bf F}(\varepsilon,\mathtt{V}_1,\mathtt{V}_2,f)[h] =\varepsilon^3 |\ln\varepsilon | \omega \partial_\varphi h(\varphi ,\theta)+\big(2\mathtt{p}_{1,1}(\varphi)\big)^{-\frac12} \partial_\theta \big\{ \partial_f \Psi({\gamma}(\varphi,\theta))[h]\big\}\\ &\quad  -\varepsilon^2 \omega \big( |\ln \varepsilon|\dot{\mathtt{P}_1}(\varphi)+i \,\varepsilon \dot{\mathtt{V}}_1(\varphi)\big)\cdot i\partial_\theta\Big\{\tfrac{h(\varphi,\theta)}{w(\varphi,\theta)}\mathcal{Z}(\varphi,\theta)\Big\}+\varepsilon^3 |\ln\varepsilon |\,\omega\,  \tfrac{ \dot{\mathtt{p}}_{1,1}(\varphi)}{4\mathtt{p}_{1,1}(\varphi)}\partial_\theta\big\{ h(\varphi,\theta)\sin(2\theta)\big\}.  \nonumber
\end{align}
Next, from  \eqref{Psi-form}, we decompose the stream function along the boundary as
$$
\Psi(\gamma(\varphi,\theta))=\Psi_1(\varepsilon,f)+\Psi_2(\varepsilon,f),
$$
where $\Psi_1$ and $\Psi_2$ denote, respectively, the self-interaction and interaction contributions, 
\begin{align*} 
\Psi_1(\varepsilon,f)&=\sqrt{2}\fint_{\T}\int_0^{w(\varphi,\eta)}G\left(\Gamma_1(f)(\varphi,\theta,1),
\Gamma_1(0)(\varphi,\eta,\rho)\right)\rho d\rho d\eta\\
\Psi_2(\varepsilon,f)&=\sqrt{2}\int_{\T}\int_0^{w_\star(\varphi,\eta)}G\left(\Gamma_1(\varphi,\theta,1)+i\tfrac{\varepsilon }{ |\ln\varepsilon|} \mathtt{U}_{-}(\varphi),\Gamma_2(0)(\varphi,\eta,\rho)\right)\rho  d\rho d\eta,
\end{align*}
From \eqref{def:Gamma12} one has
\begin{equation*}
\begin{aligned}
\partial_f w[h](\varphi,\eta)&=\varepsilon \frac{h(\varphi,\eta)}{w(\varphi,\eta)},\\
\partial_f\Gamma_1(f)[h](\varphi,\eta,\rho)&=\varepsilon\rho\, \frac{h(\varphi,\eta)}{w(\varphi,\eta)}\mathcal{Z}(\varphi,\eta),\\
\partial_f\Gamma_2(f)[h](\varphi,\eta,\rho) &=\varepsilon\rho\, \frac{h_\star(\varphi,\eta)}{w_\star(\varphi,\eta)}\mathcal{Z}_\star(\varphi,\eta).
\end{aligned}
\end{equation*}
Thus,  linearizing $\Psi_1(\varepsilon,f)$ gives  
\begin{align*}
 & \partial_\theta \big\{ \partial_f\Psi_1(\varepsilon,f)[h]\big\}(\varphi,\theta)=\sqrt{2}{\varepsilon}\partial_\theta\fint_{\T}G\big(\mathtt{P}_1(\varphi)+\varepsilon (w \mathcal{Z})(\varphi,\theta),\mathtt{P}_1(\varphi)+\varepsilon (w   \mathcal{Z})(\varphi,\eta)\big)h(\varphi,\eta)  d\eta\\
&+\sqrt{2} \varepsilon^2\partial_\theta\bigg(\tfrac{h(\varphi,\theta)}{w(\varphi,\theta)}\fint_{\T}\int_0^{w(\varphi,\eta)}\nabla_1G\big(\mathtt{P}_1(\varphi)+\varepsilon (w \mathcal{Z})(\varphi,\theta),\mathtt{P}_1(\varphi)+\varepsilon \rho   \mathcal{Z}(\varphi,\eta)\big)\cdot\mathcal{Z}(\varphi,\theta)\rho d\rho d\eta\bigg).
\end{align*}
Using the change of variable $\rho\mapsto \rho w(\varphi,\eta)$ in the last integral  we find 
\begin{align*}
 \partial_\theta \big\{ \partial_f\Psi_1(\varepsilon,f)&[h]\big\}(\varphi,\theta)=\sqrt{2} \varepsilon \,\mathcal{I}_1(\varepsilon,f)[h](\varphi,\theta)+\sqrt{2} \varepsilon\, \partial_\theta\big\{h(\varphi,\theta) U_1(\varepsilon,f)(\varphi,\theta)\big\}.
\end{align*}
In a similar way, we get
\begin{align*}
 & \partial_\theta \big\{ \partial_f\Psi_2(\varepsilon,f)[h]\big\}(\varphi,\theta)=\sqrt{2} \varepsilon \,\mathcal{I}_2(\varepsilon,f)[h](\varphi,\theta)+\sqrt{2} \varepsilon\, \partial_\theta\big\{h(\varphi,\theta)U_2(\varepsilon,f)(\varphi,\theta)\big\}.
\end{align*}
Substituting these two identities into \eqref{diff-F-def} yields the claimed formula.
Finally, the $C^1$-regularity of $f\mapsto {\bf F}(\varepsilon,\mathtt{V}_1,\mathtt{V}_2,f)$
follows from the product and composition estimates.
\end{proof}

To unravel the fine structure of the linearized operator given by Lemma  \ref{prop:asymp-lin0}  we begin by expanding the contribution arising from the \emph{induced effect}, followed by that of the interaction term. This allows us to identify the leading singular contributions and isolate the terms that will be crucial in the subsequent reducibility and invertibility analysis.

\subsection{Induced effect expansion}

For now, our goal is to construct an explicit decomposition of the operator $\mathcal{I}_1(f)$, defined in Lemma  \ref{prop:asymp-lin0}. 
We  proceed with the definition of the main operators that will be used in the expansion.
The first one is the \emph{toroidal Hilbert transform}, previously introduced in \eqref{Hilbert1}, which admits the equivalent representation
\begin{align}\label{Hilbert1alt}
\mathcal{H}[h](\varphi,\theta)&=-\frac{1}{\pi} \partial_\theta\int_0^{2\pi}\ln(|e^{i\theta}-  e^{i\eta}|)h(\varphi,\eta) d\eta=\fint_{\mathbb{T}} h(\varphi,\eta) \cot\left(\tfrac{\eta-\theta}{2}\right) d\eta.
\end{align}
The second operator  localizes on the spatial mode~$1$,  
\begin{align}\label{def:Q1}
      \mathcal{Q}_1[h](\varphi,\theta):=&-\tfrac{1}{8 }(2\mathtt{p}_{1,1})^{-\frac32}\Bigg( \cos\theta\fint_{\mathbb{T}} h(\eta) \cos\eta d\eta+3\sin\theta\fint_{\mathbb{T}} h(\eta) \sin\eta d\eta\Bigg).
  \end{align}
We also recall the \emph{shift operator} introduced in \eqref{Hf-n},
$$
\mathcal{H}_{\mathtt{u},0}=\partial_\theta\left(\mathtt{u}\Lambda_{0}+\Lambda_{0}\mathtt{u}\right),
$$
where the operator $\Lambda_{0}$ is a Fourier multiplier defined by \eqref{FM-n-1} with $m=0$, and the function $\mathtt{u}$ is given by
\begin{equation}\label{func-mathttu}
    \mathtt{u}(\varphi,\theta):=\tfrac{1}{2}   (2\mathtt{p}_{1,1})^{-\frac34}\cos(\theta).
\end{equation}
Thus,
\begin{align}\label{def-Hu0}
   \mathcal{H}_{\mathtt{u},0}[h](\varphi,\theta)=-\tfrac12(2\mathtt{p}_{1,1})^{-\frac34}\partial_\theta \fint_{\T} \ln\!\Big|\sin\!\Big(\tfrac{\theta-\eta}{2}\Big)\Big|\,\big(\cos(\theta)+\cos(\eta)\big) h(\varphi,\eta)\, d\eta.
\end{align}
Another \emph{shift operator} localizing on the spatial Fourier modes $1$ and $2$, is defined by
\begin{align}\label{shift-operator1}
\mathcal{S}[h](\varphi,\theta)&:= 
-\tfrac{1}{8}   (2\mathtt{p}_{1,1})^{-\frac34}\fint_{\mathbb{T}}\Big({\cos(\theta+2\eta)+\cos(2\theta+\eta)}\Big) h(\varphi,\eta)d\eta.
\end{align}
The main result reads as follows.
\begin{pro}\label{prop-induc25}
Under assumption \eqref{cond1}, there exists $\varepsilon_0>0$ such that for
every $\varepsilon\in(0,\varepsilon_0)$, every smooth
$f\in B_{1,\mathrm{even}}(s_0)$, and every
$h\in\mathrm{Lip}_\nu(\mathcal{O},H^{s}_{\star,\mathrm{even}}(\T^2))$, one has
\begin{align*}
\frac{1}{\sqrt{\mathtt{p}_{1,1}(\varphi)}}&\mathcal{I}_1(\varepsilon,f)[h](\varphi,\theta)
=
\Big(\tfrac12\mathcal{H}[h]
+\varepsilon\mathcal{H}_{\mathtt{u},0}[h]
+\varepsilon \partial_\theta\mathcal{S}[h]
+\varepsilon^2|\ln\varepsilon|\,\partial_\theta\mathcal{Q}_1[h]\Big)(\varphi,\theta)\\
&\quad+\varepsilon^2\partial_\theta\fint_{\T}
W(f)(\varphi,\theta,\eta)\,
\ln\Big|\sin\Big(\tfrac{\theta-\eta}{2}\Big)\Big|\,
h(\varphi,\eta)\,d\eta
+\varepsilon^2\partial_\theta\mathcal{R}_1(\varepsilon,f)[h](\varphi,\theta),
\end{align*}
where $W(f):\mathcal{O}\times\T^3\to\R$ is smooth and satisfies
$$
\|W(f)\|^{\mathrm{Lip},\nu}_{s}\lesssim 1+\|f\|^{\mathrm{Lip},\nu}_{s},
\qquad
W(f)(-\varphi,-\theta,-\eta)=W(f)(\varphi,\theta,\eta)
=W(f)(\varphi,\eta,\theta).
$$
Moreover, the operator $\mathcal{R}_1(\varepsilon,f)$ is smoothing: for every $N\in\N$ and every
$s\ge s_0$,
$$
\|\partial_\theta^N\mathcal{R}_1(\varepsilon,f)[h]\|^{\mathrm{Lip},\nu}_s
\lesssim (1+\|f\|_{s_0+N}^{\mathrm{Lip},\nu})\|h\|_s^{\mathrm{Lip},\nu}
+\|f\|_{s+N+1}^{\mathrm{Lip},\nu}\|h\|_{s_0}^{\mathrm{Lip},\nu}.
$$
\end{pro}
\begin{proof} We note that, through this proof, the $\varphi$–dependence is omitted in order to simplify the notation.     Using Lemma \ref{lem:kernel-expansion} with
$$
X=(x_1,x_2)=w(\theta)\mathcal{Z}(\theta),\quad
Z=\mathtt{P}_1=(\mathtt{p}_{1,1},\mathtt{p}_{1,2}),\quad Y=(y_1,y_2)=w(\eta) \mathcal{Z}(\eta),
$$we infer through straightforward computations that
    $$
G\big(\mathtt{P}_1+\varepsilon w(\theta) \mathcal{Z}(\theta),\mathtt{P}_1+\varepsilon w(\eta)   \mathcal{Z}(\eta)\big)
=|\ln(\varepsilon)|\sum_{n=0}^2\mathcal{A}_n\varepsilon^n+\sum_{n=0}^1\mathcal{B}_n \varepsilon^n+\varepsilon^2 \mathcal{C},
$$
where
\begin{align*}
   &\mathcal{A}_0=\sqrt{\mathtt{p}_{1,1}},\qquad
   \mathcal{A}_1=\frac{\sqrt{2}}{4}(2\mathtt{p}_{1,1})^{-\frac14}\Big(w(\theta)\cos(\theta)+w(\eta) \cos(\eta)\Big),\\ 
  & \mathcal{A}_2=-\frac{\sqrt{2}}{32}(2\mathtt{p}_{1,1})^{-1}\Big(3w^2(\theta)\cos(2\theta)+3w^2(\eta) \cos(2\eta)+2w(\eta) w(\theta)\cos(\theta) \cos(\eta)\\ &\qquad \qquad\qquad\qquad \qquad\ +6w(\eta) w(\theta)\sin(\theta) \sin(\eta)\Big),
\end{align*}
and
\begin{align*}
 \mathcal{B}_0&=C_0(\mathtt{p}_{1,1})-\frac{\sqrt{\mathtt{p}_{1,1}}}{2}\ln\Big(\big(w(\theta)-w(\eta)\big)^2
+4 w(\theta)w(\eta)\sin^2\Big(\tfrac{\theta-\eta}{2}\Big)\Big).
   \end{align*}
   The coefficient $\mathcal{B}_1$ will be evaluated later.
   As $\langle h\rangle_\theta=0,$ then
   \begin{align*}
     \mathcal{I}_1(\varepsilon,f)[h](\varphi,\theta)&=\varepsilon^2|\ln\varepsilon|\partial_\theta\fint_{\mathbb{T}} \mathcal{A}_2(\varphi,\theta,\eta)h(\varphi,\eta) d\eta+\partial_\theta\fint_{\mathbb{T}} \mathcal{B}_0(\varphi,\theta,\eta)h(\varphi,\eta) d\eta\\
    &\quad +\varepsilon\partial_\theta\fint_{\mathbb{T}} \mathcal{B}_1(\varphi,\theta,\eta)h(\varphi,\eta) d\eta+\varepsilon^2\partial_\theta\fint_{\mathbb{T}} \mathcal{C}(\varphi,\theta,\eta)h(\varphi,\eta) d\eta.
   \end{align*}
   Since $w=\sqrt{1+2\varepsilon f}$ then
    \begin{align*}
    &\partial_\theta\fint_{\mathbb{T}} \mathcal{A}_2(\varphi,\theta,\eta)h(\varphi,\eta) d\eta\\ &=-\frac{\sqrt{2}}{16}(2\mathtt{p}_{1,1})^{-1}\partial_\theta\fint_{\mathbb{T}}\Big(w(\eta) w(\theta)\cos(\theta) \cos(\eta)+3w(\eta) w(\theta)\sin(\theta) \sin(\eta)\Big)h(\varphi,\eta)d\eta\\
    &=-\frac{\sqrt{2}}{16}(2\mathtt{p}_{1,1})^{-1}\partial_\theta\fint_{\mathbb{T}}\Big(\cos(\theta) \cos(\eta)+3\sin(\theta) \sin(\eta)\Big)h(\varphi,\eta)d\eta+\varepsilon\partial_\theta\mathcal{R}_{\mathcal{A}_2}[h],
   \end{align*}
   where 
   \begin{align*}
    \mathcal{R}_{\mathcal{A}_2}[h])(\varphi,\theta)&=-\frac{\sqrt{2}}{16}(2\mathtt{p}_{1,1})^{-1}\fint_{\mathbb{T}}\big(w(\eta) w(\theta)-1\big)\Big(\cos(\theta) \cos(\eta)+3\sin(\theta) \sin(\eta)\Big)h(\varphi,\eta)d\eta\\
    &=:\fint_{\mathbb{T}}\mathcal{K}_{\mathcal{A}_2}(\varphi,\theta,\eta)h(\varphi,\eta)d\eta.
   \end{align*}
   We note that
   $\mathcal{R}_{\mathcal{A}_2}$ is smoothing at any order and belongs to the class $ \textnormal{OPS}^{-\infty}$. It  satisfies the estimates
$$
\forall N\in\mathbb{N},\;\forall s\geqslant s_0, \quad \|\partial_\theta^N\mathcal{R}_{\mathcal{A}_2}[h]\|^{\textnormal{Lip},\nu}_s\lesssim \|f\|^{\textnormal{Lip},\nu}_{s_0+N}\|h\|_s +\|f\|^{\textnormal{Lip},\nu}_{s+N}\|h\|^{\textnormal{Lip},\nu}_{s_0}.$$
On the other hand, we observe that
   \begin{align*}
 \mathcal{B}_0
 &=C_0(\mathtt{p}_{1,1})-\frac{\sqrt{\mathtt{p}_{1,1}}}{2}\ln(4 w(\theta)w(\eta))-{\sqrt{\mathtt{p}_{1,1}}}\ln\big|\sin\big(\tfrac{\theta-\eta}{2}\big)\big|+\varepsilon^2 \mathcal{C}_1,
   \end{align*}
   with $\mathcal{C}_1$ a smooth function that
   $$
   \|\mathcal{C}_1\|^{\textnormal{Lip},\nu}_{s}\lesssim \|f\|^{\textnormal{Lip},\nu}_{s+1}.
   $$
 Therefore,
   \begin{align*}
    \partial_\theta\fint_{\mathbb{T}} \mathcal{B}_0(\varphi,\theta,\eta)h(\varphi,\eta) d\eta&=-\sqrt{\mathtt{p}_{1,1}}\,\partial_\theta\fint_{\mathbb{T}}\ln\big|\sin\big(\tfrac{\theta-\eta}{2}\big)\big|h(\varphi,\eta)d\eta+\varepsilon^2\partial_\theta\mathcal{R}_{\mathcal{B}_0}[h],
   \end{align*}
   with $\mathcal{R}_{\mathcal{B}_0}\in \textnormal{OPS}^{-\infty}$ and satisfies
   $$
\forall N\in\mathbb{N},\;\forall s\geqslant s_0, \quad \|\partial_\theta^N\mathcal{R}_{\mathcal{B}_0}[h]\|^{\textnormal{Lip},\nu}_s\lesssim \|f\|^{\textnormal{Lip},\nu}_{s_0+N}\|h\|_s +\|f\|^{\textnormal{Lip},\nu}_{s+1+N}\|h\|^{\textnormal{Lip},\nu}_{s_0}.$$ 
 Return to $\mathcal{B}_1$, which takes the form
   \begin{align*}
   \mathcal{B}_1(\varepsilon)&=\frac{1}{4 \sqrt{z_1}}(x_1+y_1)\bigg(\ln8-2-\tfrac12\ln\Big(\frac{(x_1-y_1)^2+2z_1(x_2-y_2)^2}{(2z_1)^2}\Big)\bigg)\\  &\quad +\frac{x_1+y_1}{2\sqrt{z_1}}\Big(\frac{(x_1-y_1)^2+z_1(x_2-y_2)^2}{(x_1-y_1)^2+2z_1(x_2-y_2)^2}\Big).
\end{align*}
It follows that
\begin{align*}
  &\partial_\theta\fint_{\mathbb{T}}\mathcal{B}_1(0) \,h(\eta)
  =-\tfrac{\sqrt{2}}{4}(2\mathtt{p}_{1,1})^{-\frac14}\partial_\theta\int_0^{2\pi}\Big(\cos\theta+\cos\eta\Big)\ln\big|\sin\big(\tfrac{\theta-\eta}{2}\big)\big|h(\eta) d\eta
  \\ &\qquad +\tfrac{1}{\sqrt{2} }(2\mathtt{p}_{1,1})^{-\frac14}\partial_\theta\fint_{\mathbb{T}}\Big({\cos\theta+\cos\eta}\Big)\frac{(\cos\theta-\cos\eta)^2+\frac12(\sin\theta-\sin\eta)^2}{(\cos\theta-\cos\eta)^2+(\sin\theta-\sin\eta)^2}h(\eta)d\eta.
   \end{align*}
Using the identity
\begin{align*}
\frac{(\cos\theta-\cos\eta)^2+\frac12(\sin\theta-\sin\eta)^2}{(\cos\theta-\cos\eta)^2+(\sin\theta-\sin\eta)^2}&=1-\frac18\frac{(\sin\theta-\sin\eta)^2}{\sin^2(\frac{\theta-\eta}{2})}\\
&=1-\frac12\cos^2\Big(\tfrac{\theta+\eta}{2}\Big)\\
&=\frac34-\frac14\cos(\theta+\eta),
\end{align*}
and applying classical trigonometric identities, we obtain
\begin{align*}
  \frac{1}{\sqrt{2} }(2\mathtt{p}_{1,1})^{-\frac14}&\partial_\theta\fint_{\mathbb{T}}\Big({\cos\theta+\cos\eta}\Big)\Big(\frac{(\cos\theta-\cos\eta)^2+\tfrac12(\sin\theta-\sin\eta)^2}{(\cos\theta-\cos\eta)^2+(\sin\theta-\sin\eta)^2}\Big)h(\eta)d\eta\\
  &=-\frac{1}{8\sqrt{2} }(2\mathtt{p}_{1,1})^{-\frac14}\partial_\theta\fint_{\mathbb{T}}\Big({\cos(\theta+2\eta)+\cos(2\theta+\eta)}\Big) h(\eta)d\eta.
   \end{align*}
   Consequently,
   \begin{align*}
  \partial_\theta\int_{0}^{2\pi}\mathcal{B}_1(0) \,h(\eta)
  &=-\frac{1}{2\sqrt{2} }(2\mathtt{p}_{1,1})^{-\frac14}\partial_\theta\fint_{\mathbb{T}}\Big(\cos\theta+\cos\eta\Big)\ln|\sin\big(\tfrac{\theta-\eta}{2}\big)|h(\eta) d\eta
  \\  &\quad-\frac{1}{8\sqrt{2} }(2\mathtt{p}_{1,1})^{-\frac14}\partial_\theta\fint_{\mathbb{T}}\Big({\cos(\theta+2\eta)+\cos(2\theta+\eta)}\Big) h(\eta)d\eta.
   \end{align*}
   By expanding $\mathcal{B}_1(\varepsilon)$ we get
   \begin{align*}
  \partial_\theta\int_{0}^{2\pi}\mathcal{B}_1(\varepsilon) \,h(\eta)
  =&-\frac{1}{2\sqrt{2} }(2\mathtt{p}_{1,1})^{-\frac14}\partial_\theta\fint_{\mathbb{T}}\Big(\cos\theta+\cos\eta\Big)\ln\big|\sin\big(\tfrac{\theta-\eta}{2}\big)\big|h(\eta) d\eta
  \\ \quad &-\frac{1}{8\sqrt{2} }(2\mathtt{p}_{1,1})^{-\frac14}\partial_\theta\fint_{\mathbb{T}}\Big({\cos(\theta+2\eta)+\cos(2\theta+\eta)}\Big) h(\eta)d\eta\\
  &+\varepsilon \partial_\theta\fint_{\mathbb{T}}W_{0,1}(\varphi,\theta,\eta)\ln\big|\sin\big(\tfrac{\theta-\eta}{2}\big)\big|h(\eta) d\eta+\varepsilon\partial_\theta\mathcal{R}_{\mathcal{B}_1}[h],
   \end{align*}
   where $W_{0,1}$ is a smooth function which satisfies the estimates
   $$
   \|W_{0,1}\|^{\textnormal{Lip},\nu}_{s}\lesssim \|f\|^{\textnormal{Lip},\nu}_{s},
   $$
   and $\mathcal{R}_{\mathcal{B}_1}\in\textnormal{OPS}^{-\infty}$  is a smoothing operator with
   $$
\|\partial_\theta^N\mathcal{R}_{\mathcal{B}_1}[h]\|^{\textnormal{Lip},\nu}_s\lesssim \|f\|^{\textnormal{Lip},\nu}_{s_0+N}\|h\|_s +\|f\|^{\textnormal{Lip},\nu}_{s+1+N}\|h\|^{\textnormal{Lip},\nu}_{s_0}.$$
   Notice that the structure of  $\mathcal{C}$ is quite similar to $\mathcal{B}_0 $ and $\mathcal{B}_1$ in the sense that it decomposes as
   $$
   \mathcal{C}(\varphi,\theta,\eta)=W_{0,2} (\varphi,\theta,\eta)\ln\big|\sin\big(\tfrac{\theta-\eta}{2}\big)\big|+W_{0,3}(\varphi,\theta,\eta),
   $$
   where $W_{0,2}, W_{0,3}$ are smooth functions with
   $$
   \|W_{0,2}\|^{\textnormal{Lip},\nu}_{s}\lesssim 1+ \|f\|^{\textnormal{Lip},\nu}_{s},\quad \|W_{0,3}\|^{\textnormal{Lip},\nu}_{s}\lesssim 1+\|f\|^{\textnormal{Lip},\nu}_{s+1}.
   $$
   Putting altogether the preceding identities ends the proof of the desired expansion.
\end{proof}

We now turn to the expansion of the function $U_1(\varepsilon,f)(\varphi,\theta)$ defined in Lemma \ref{prop:asymp-lin0},
corresponding to the self-induced velocity within the transport component of the linearized operator.
\begin{lem}\label{prop-V1}
  Consider a  smooth function
$f\in B_{1,\mathrm{even}}(s_0)$. Then, the function $U_1(\varepsilon,f)$ admits the expansion
\begin{align*}
\frac{1}{\sqrt{\mathtt{p}_{1,1}(\varphi)}}U_1(\varepsilon,f)(\varphi,\theta)
&=-\frac12+ \frac{\varepsilon}{\sqrt{2\mathtt{p}_{1,1}}}(\nabla_{\mathtt{P}_1} ^\perp \mathtt{G})(\mathtt{P}_1,\mathtt{P}_2)\cdot i \,\frac{\mathcal{Z}(\varphi,\theta)}{w(\varphi,\theta)}+\varepsilon\mathtt{g}_3(\varphi)\left(\cos(3\theta)-2  \cos(\theta)\right)
\\ &\quad+\frac\varepsilon2 f(\varphi,\theta) -\frac{3}{16}\varepsilon^2|\ln\varepsilon|\,(2\mathtt{p}_{1,1})^{-3/2}\cos(2\theta)+\varepsilon^2 \widetilde{U}_1
(\varepsilon,f)(\varphi,\theta).
\end{align*}
where $\mathtt{g}_3$ and $(\nabla_{\mathtt{P}_1} ^\perp \mathtt{G})(\mathtt{P}_1,\mathtt{P}_2)$ are given by \eqref{list-functions} and \eqref{nablaGtt}, respectively.  
Moreover the function   $\widetilde{U}_1(\varepsilon,f)$ is a smooth and satisfies 
$$
   \|\widetilde{U}_1(\varepsilon,f)\|^{\textnormal{Lip},\nu}_{s}\lesssim 1+ \|f\|^{\textnormal{Lip},\nu}_{s}.
   $$

\end{lem}
\begin{proof}
Applying Lemma~\ref{pro-decomp2-nabla} with
$$
X=(x_1,x_2)=w(\theta)(\cos\theta,\sin\theta), \quad Y=(y_1,y_2)=\rho w(\eta)(\cos\eta,\sin\eta), \quad z_1=\mathtt{p}_{1,1},
$$
we obtain
\begin{align*}
&\varepsilon\,\nabla_1G\big(P+\varepsilon w(\theta)\mathcal{Z}(\theta),\,P+\varepsilon\rho w(\eta)\mathcal{Z}(\eta)\big)\cdot\mathcal{Z}(\theta)\\
&\qquad=
|\ln\varepsilon|\Big(\varepsilon\,\mathscr C_1+\varepsilon^2\,\mathscr C_2\Big)(X,Y,Z)\cdot\mathcal{Z}(\theta)
+\Big(\mathscr D_0+\varepsilon\,\mathscr D_1\Big)(X,Y,Z)\cdot\mathcal{Z}(\theta)
+O(\varepsilon^2).
\end{align*}
Consequently,
\begin{align}\label{eq:V1-split}
\frac{1}{\sqrt{\mathtt{p}_{1,1}}}U_1(\varepsilon,f)(\theta)
&=\frac{1}{\sqrt{\mathtt{p}_{1,1}}}
\fint_{\T}\int_0^{1}
|\ln\varepsilon|\Big(\varepsilon\,\mathscr C_1+\varepsilon^2\,\mathscr C_2\Big)\cdot\mathcal{Z}(\theta)\,
\frac{w^2(\eta)}{w(\theta)}\,\rho\,d\rho\,d\eta \nonumber\\
&\quad+\frac{1}{\sqrt{\mathtt{p}_{1,1}}}
\fint_{\T}\int_0^{1}
\Big(\mathscr D_0+\varepsilon\,\mathscr D_1\Big)\cdot\mathcal{Z}(\theta)\,
\frac{w^2(\eta)}{w(\theta)}\,\rho\,d\rho\,d\eta
+O(\varepsilon^2).
\end{align}
Using $\mathscr C_1(X,Y,Z)=\frac{1}{4\sqrt{\mathtt{p}_{1,1}}}(1,0)$ and the definition of $\mathcal Z$, we get
$$
\mathscr C_1(X,Y,Z)\cdot\mathcal Z(\theta)=\frac{1}{2\sqrt2}(2\mathtt{p}_{1,1})^{-1/4}\cos(\theta).
$$
Since $\fint_{\T}f(\eta)\,d\eta=0$, we have
$$\fint_{\T}w^2(\eta)\,d\eta=1.
$$
Therefore,
\begin{align*}
\frac{1}{\sqrt{\mathtt{p}_{1,1}}}\fint_{\T}\int_0^{1}
\mathscr C_1\cdot\mathcal Z(\theta)\,\frac{w^2(\eta)}{w(\theta)}\,\rho\,d\rho\,d\eta
&=(2\mathtt{p}_{1,1})^{-3/4}\frac{1}{4 w(\theta)}\cos(\theta).
\end{align*}
From Lemma~\ref{pro-decomp2-nabla} one computes
$$
\mathscr C_2(X,Y,Z)\cdot\mathcal Z(\theta)
=\frac{\sqrt2}{64\mathtt{p}_{1,1}}\Big(-6w(\theta)\cos(2\theta)-2\rho w(\eta)\cos\theta\cos\eta-6\rho w(\eta)\sin\theta\sin\eta\Big).
$$
At the order $\varepsilon^2|\ln\varepsilon|$ we may set $w\equiv1$ in this contribution,
so that the $\eta$--oscillatory parts average to $0$ and only the $-6\cos(2\theta)$ term remains. Therefore
\begin{align*}
\frac{1}{\sqrt{\mathtt{p}_{1,1}}}\fint_{\T}\int_0^{1}
\mathscr C_2\cdot\mathcal Z(\theta)\,\frac{w^2(\eta)}{w(\theta)}\,\rho\,d\rho\,d\eta
=-\frac{3}{16}(2\mathtt{p}_{1,1})^{-3/2}\cos(2\theta)+O(\varepsilon).
\end{align*}
Using the explicit formula for $\mathscr D_0$,
$$
\mathscr D_0(X,Y,Z)\cdot\mathcal Z(\theta)
=-\sqrt{\mathtt{p}_{1,1}}\,
\frac{w(\theta)-\rho w(\eta)\cos(\theta-\eta)}
{w^2(\theta)+\rho^2w^2(\eta)-2\rho w(\theta)w(\eta)\cos(\theta-\eta)},
$$
gives
\begin{align*}
&\frac{1}{\sqrt{\mathtt{p}_{1,1}}}\fint_{\T}\int_0^{1}
\mathscr D_0\cdot\mathcal Z(\theta)\,\frac{w^2(\eta)}{w(\theta)}\,\rho\,d\rho\,d\eta\\
&\qquad
=-\frac{1}{ w(\theta)}\fint_{\T}\int_0^{1}
\frac{w(\theta)-\rho w(\eta)\cos(\theta-\eta)}
{w^2(\theta)+\rho^2w^2(\eta)-2\rho w(\theta)w(\eta)\cos(\theta-\eta)}
\,w^2(\eta)\,\rho\,d\rho\,d\eta.
\end{align*}
Introduce
$$
\mathcal{J}(\alpha):=\fint_{\T}\int_0^{1}
\ln\Big(\big|\alpha w(\theta)e^{i\theta}-\rho w(\eta)e^{i\eta}\big|^2\Big)\,w^2(\eta)\,\rho\,d\rho\,d\eta.
$$
A direct differentiation gives
$$
\frac{1}{\sqrt{\mathtt{p}_{1,1}}}\fint_{\T}\int_0^{1}\mathscr D_0\cdot\mathcal Z(\theta)\,\frac{w^2(\eta)}{w(\theta)}\,\rho\,d\rho\,d\eta
=-\frac{1}{2 w^2(\theta)}\mathcal{J}'(\alpha)\big|_{\alpha=1}.
$$
Expanding $w(\theta)=1+\varepsilon f(\theta)+O(\varepsilon^2)$ and using the change of variable $\rho\to \alpha\rho $ , one obtains
\begin{align*} 
\frac{\mathcal{J}(\alpha)}{2 w^2(\theta)}  &= {\alpha^2}\fint_{\T} \int_0^{1/\alpha} \left[ \ln(\alpha) + \ln\left( |e^{i\theta} - \rho e^{i\eta}| \right) \right] \rho \left( 1 + 2\varepsilon(f(\eta) - f(\theta)) \right)d\rho d\eta \\ &\quad +  \varepsilon \alpha^2 \fint_{\T} \int_0^{1/\alpha} \text{Re} \left( \frac{e^{i\theta}f(\theta) - \rho e^{i\eta}f(\eta)}{e^{i\theta} - \rho e^{i\eta}} \right) \rho  d\rho d\eta + O(\varepsilon^2). \end{align*}
Taking the derivative with respect to $\alpha$ at $\alpha = 1$, we get  
\begin{align*}
     \frac{\mathcal{J}'(1)}{2 w^2(\theta)}
     &=\frac12(1-2\varepsilon f(\theta))+  4\varepsilon \fint_{\T}\int_0^{1}  f(\eta)\ln\left(| e^{i\theta}-\rho e^{i\eta}|\right)\rho d\rho d\eta
     -2\varepsilon  \fint_{\T}  \ln\left(| e^{i\theta}- e^{i\eta}|\right) f(\eta) d\eta 
     \\
     &\quad + 2\varepsilon  \fint_{\T}\int_0^{1} \text{Re} \left( \frac{e^{i\theta}f(\theta) - \rho e^{i\eta}f(\eta)}{e^{i\theta} - \rho e^{i\eta}} \right)\rho  d\rho d\eta
     - \varepsilon\fint_{\T}\text{Re} \left( \frac{e^{i\theta}f(\theta) -  e^{i\eta}f(\eta)}{e^{i\theta} -  e^{i\eta}} \right)  d\eta.
\end{align*}
Using the identity
\begin{align}\label{identity re}
\text{Re} \left( \frac{e^{i\theta}f(\theta) - \rho e^{i\eta}f(\eta)}{e^{i\theta} - \rho e^{i\eta}} \right)&=\frac{1-\rho\cos(\theta-\eta)}{1+\rho^2-2\rho \cos(\theta-\eta)}f(\theta)+\frac{\rho-\cos(\theta-\eta)}{1+\rho^2-2\rho \cos(\theta-\eta)}\rho f(\eta),
\end{align}
and the integral formula \eqref{id:integralV11} we get
\begin{align*}
   2 \fint_{\T}\int_0^{1}  \text{Re} \left( \frac{e^{i\theta}f(\theta) - \rho e^{i\eta}f(\eta)}{e^{i\theta} - \rho e^{i\eta}} \right)\rho  d\rho d\eta&=  f(\theta)+  \fint_{\T}\int_0^{1}  \frac{2\rho-2\cos(\theta-\eta)}{1+\rho^2-2\rho \cos(\theta-\eta)} f(\eta)\rho^2  d\rho d\eta.
\end{align*}
Then, integrating by parts  yields
\begin{align*}
    2  \fint_{\T}\int_0^{1}  \text{Re} \left( \frac{e^{i\theta}f(\theta) - \rho e^{i\eta}f(\eta)}{e^{i\theta} - \rho e^{i\eta}} \right)\rho  d\rho d\eta    &=f(\theta)+  \fint_{\T}\int_0^{1}\partial_\rho \ln\left(| e^{i\theta}-\rho e^{i\eta}|^2\right)   f(\eta)\rho^2  d\rho d\eta\\
    &=f(\theta)+ 2 \fint_{\T} \ln\left(| e^{i\theta}- e^{i\eta}|\right)   f(\eta)  d\rho d\eta\\ 
    &- 4\fint_{\T}\int_0^{1} \ln\left(| e^{i\theta}-\rho e^{i\eta}|\right)   f(\eta)\rho  d\rho d\eta.
\end{align*}
Using again the identity \eqref{identity re} we find   
$$
\fint_{\T} \text{Re} \left( \frac{e^{i\theta}f(\theta) -  e^{i\eta}f(\eta)}{e^{i\theta} -  e^{i\eta}} \right) d\eta=\frac12 f(\theta).
$$
Bringing the computations above together, we obtain
\begin{align*}
     \frac{\mathcal{J}'(1)}{ 2 w^2(\theta)}
     &=  \frac12-\frac\varepsilon2 f(\theta)+O(\varepsilon^2).
\end{align*}
Hence
\begin{align*}
     \frac{1}{\sqrt{p_1} }\fint_{\T}\int_0^{1}\mathscr{D}_0\cdot\mathcal{Z}(\theta) \frac{w^2(\eta)}{w(\theta)}\rho d\rho d\eta
     &=- \frac12+\frac\varepsilon2 f(\theta)+O(\varepsilon^2).
\end{align*}
Finally, from the explicit formula of $\mathscr{D}_1$ we have
\begin{align*}
   \mathscr{D}_1&(X,Y,Z)\cdot\mathcal{Z}(\theta)
   =\frac{(2\mathtt{p}_{1,1})^{\frac14}}{2\sqrt{2} }\Big(\tfrac54\ln8-1+\tfrac34\ln(z_1)\Big)\cos(\theta)-\frac{(2\mathtt{p}_{1,1})^{\frac14}}{4\sqrt{2} }\cos(\theta) \ln\Big(|e^{i\theta}-\rho e^{i\eta}|^2\Big)\\ \quad &+\frac{(2\mathtt{p}_{1,1})^{\frac14}}{2\sqrt{2} }\bigg( \cos(\theta)\frac{(\cos(\theta)-\rho\cos(\eta))^2}{|e^{i\theta}-\rho e^{i\eta}|^2}- \cos(\theta)\frac{(\cos(\theta)+\rho \cos(\eta))(\cos(\theta)-\rho\cos(\eta))^3}{|e^{i\theta}-\rho e^{i\eta}|^4}\\
   \ &+ \cos(\theta)\frac{(\cos(\theta)+\rho \cos(\eta))(\cos(\theta)-\rho\cos(\eta))(\sin(\theta)-\rho\sin(\eta))^2}{|e^{i\theta}-\rho e^{i\eta}|^4}\\ & -\sin(\theta)\frac{(\cos(\theta)+\rho \cos(\eta))(\sin(\theta)-\rho\sin(\eta))}{|e^{i\theta}-\rho e^{i\eta}|^2}\\ &-2 \sin(\theta)\frac{(\cos(\theta)+\rho \cos(\eta))(\cos(\theta)-\rho\cos(\eta))^2(\sin(\theta)-\rho\sin(\eta))}{|e^{i\theta}-\rho e^{i\eta}|^4}\bigg)+O(\varepsilon).
\end{align*}
Then, from the integral identities in Lemma \ref{lem:integrals}, we obtain
\begin{align*}
   &\frac{1}{ \sqrt{\mathtt{p}_{1,1}}} \fint_{\T}\int_0^{1}\mathscr{D}_1\cdot\mathcal{Z}(\theta) \frac{w^2(\eta)}{w(\theta)}\rho d\rho d\eta =\frac{(2\mathtt{p}_{1,1})^{-\frac34}}{4}\Big(\tfrac54\ln8-1+\tfrac34\ln(z_1)\Big)\cos(\theta)\\ &\quad+\frac{(2\mathtt{p}_{1,1})^{-\frac34}}{4}\bigg(\big(\tfrac{1}{2}+\tfrac{1}{4}\cos(2\theta)\big)\cos(\theta)
   -\big(\tfrac{3}{8}+\tfrac12\cos(2\theta)\big)\cos(\theta)\\ &\quad+\big(\tfrac{1}{8}+\tfrac14\cos(2\theta)\big)\cos(\theta)-\tfrac{3}{4}\sin(2\theta)\sin(\theta)-\tfrac{1}{4}\sin(2\theta)\sin(\theta)\bigg),
\end{align*}
which simplifies using trigonometric identities  to
\begin{align*}
   \frac{1}{ \sqrt{\mathtt{p}_{1,1}}} \fint_{\T}\int_0^{1}\mathscr{D}_1\cdot\mathcal{Z}(\theta) \frac{w^2(\eta)}{w(\theta)}\rho d\rho d\eta 
   &=\frac{(2\mathtt{p}_{1,1})^{-\frac34}}{16}\Big(5\ln8-5+3\ln(\mathtt{p}_{1,1})\Big)\cos(\theta)\\
  &\quad +\frac{(2\mathtt{p}_{1,1})^{-\frac34}}{8}\cos(3\theta).
\end{align*}
Substituting the preceding identities into \eqref{eq:V1-split} yields
\begin{align*}
\frac{1}{\sqrt{\mathtt{p}_{1,1}}}U_1(\varepsilon,f)(\theta)
&=-\frac12+\frac\varepsilon2 f(\theta)+\varepsilon|\ln\varepsilon|(2\mathtt{p}_{1,1})^{-3/4}\frac{1}{4 w(\theta)}\cos(\theta)
\\ &\quad +\frac{\varepsilon}{16}(2\mathtt{p}_{1,1})^{-3/4}\left(5\ln 8+3\ln (\mathtt{p}_{1,1})-5\right)\cos(\theta)\\
&\quad +\frac{\varepsilon}{8}(2\mathtt{p}_{1,1})^{-3/4}\cos(3\theta)
-\frac{3}{16}\varepsilon^2|\ln\varepsilon|\,(2\mathtt{p}_{1,1})^{-3/2}\cos(2\theta)
+O(\varepsilon^2).
\end{align*}
Note that from the asymptotics
$$  w(\theta)^{-1}=1-\varepsilon f(\theta)+O(\varepsilon^2),$$
one may write
\begin{align*}
\frac{1}{\sqrt{\mathtt{p}_{1,1}}}U_1(\varepsilon,f)(\theta)
&=-\frac12+\frac\varepsilon2 f(\theta)+\varepsilon(2\mathtt{p}_{1,1})^{-3/4}\left(|\ln\varepsilon|+\frac{1}{4}\left(5\ln 8+3\ln (\mathtt{p}_{1,1})-1\right)\right)\frac{1}{4w(\theta)}\cos(\theta)
\\ &\quad -\frac{\varepsilon}{4}(2\mathtt{p}_{1,1})^{-3/4}\cos(\theta)+\frac{\varepsilon}{8}(2\mathtt{p}_{1,1})^{-3/4}\cos(3\theta)  
-\frac{3}{16}\varepsilon^2|\ln\varepsilon|\,(2\mathtt{p}_{1,1})^{-3/2}\cos(2\theta)
\\ &\quad+O(\varepsilon^2).
\end{align*}
From the expressions of $(\nabla_{\mathtt{P}_1}^\perp \mathtt{G})(\mathtt{P}_1,\mathtt{P}_2)$ and $\mathcal{Z}$ one has
$$
\frac{1}{4}(2\mathtt{p}_{1,1})^{-3/4}\left(|\ln\varepsilon|+\frac{1}{4}\left(5\ln 8+3\ln (\mathtt{p}_{1,1})-1\right)\right)\cos(\theta)= \frac{1}{\sqrt{2\mathtt{p}_{1,1}}}(\nabla_{\mathtt{P}_1}^\perp \mathtt{G})(\mathtt{P}_1,\mathtt{P}_2) \cdot i \mathcal{Z}(\theta)
$$
It follows that
\begin{align*}
\frac{1}{\sqrt{\mathtt{p}_{1,1}}}U_1(\varepsilon,f)(\theta)
&=-\frac12+ \frac{\varepsilon}{\sqrt{2\mathtt{p}_{1,1}}}(\nabla_{\mathtt{P}_1}^\perp \mathtt{G})(\mathtt{P}_1,\mathtt{P}_2)\cdot i \frac{\mathcal{Z}(\theta)}{w(\theta)}+\frac{\varepsilon}{8}(2\mathtt{p}_{1,1})^{-3/4}\left(\cos(3\theta)-2  \cos(\theta)\right) 
\\ &\quad+\frac\varepsilon2 f(\theta)
-\frac{3}{16}\varepsilon^2|\ln\varepsilon|\,(2\mathtt{p}_{1,1})^{-3/2}\cos(2\theta)+O(\varepsilon^2),
\end{align*}
which is exactly the claimed expansion once  the remaining $O(\varepsilon^2)$ contribution is absorbed into
$\varepsilon^2\widetilde{U}_1(\varepsilon,f)$.
\end{proof}
\subsection{Interaction term  expansion}
We now derive the asymptotic expansion, as $\varepsilon\to0$, of the interaction
contribution encoded in $\mathcal I_2(\varepsilon,f)$ and $U_2(\varepsilon,f)$.
Recall that
\begin{align*}
\mathcal{I}_2(\varepsilon,f)[h](\varphi,\theta)
&=\partial_\theta\fint_{\T}
G\!\left(\Gamma_1(f)(\varphi,\theta,1)+i\tfrac{\varepsilon \mathtt{U}_{-}(\varphi)}{ |\ln\varepsilon|} ,\Gamma_2(f)(\varphi,\eta,1)\right)
h_\star(\varphi,\eta)\,d\eta,
\\
U_2(\varepsilon,f)(\varphi,\theta)
&= \varepsilon\fint_{\T}\!\int_0^{1}\!
\nabla_1G\!\left(\Gamma_1(f)(\varphi,\theta,1)+i\tfrac{\varepsilon\, \mathtt{U}_{-}(\varphi)}{ |\ln\varepsilon|} ,\Gamma_2(f)(\varphi,\eta,\rho)\right)\cdot
\tfrac{\mathcal{Z}(\varphi,\theta)}{w(\varphi,\theta)}\,w_\star^2(\varphi,\eta)\,\rho\, d\rho d\eta,
\end{align*}
where
\begin{equation*}
\begin{aligned}
\Gamma_1(f)(\varphi,\eta,\rho)&:=\mathtt{P}_1(\varphi)+\varepsilon\rho\, w(\varphi,\eta)\mathcal{Z}(\varphi,\eta),\\
\Gamma_2(f)(\varphi,\eta,\rho)&:=\mathtt{P}_2(\varphi)+\varepsilon\rho\, w_\star(\varphi,\eta)\mathcal{Z}_\star(\varphi,\eta).
\end{aligned}
\end{equation*}
\begin{lem}\label{prop-int25}
There exists $\varepsilon_0>0$ such that, for every $\varepsilon\in(0,\varepsilon_0)$
and every smooth $f\in B_{1,\mathrm{even}}(s_0)$ and $h\in \mathrm{Lip}_\nu(\mathcal{O},H^s_{\star,\mathrm{even}}(\T^2))$,
\begin{align*}
\frac{1}{\sqrt{\mathtt p_{1,1}(\varphi)}}\mathcal I_2(\varepsilon,f)[h](\varphi,\theta)
&=
\varepsilon^2|\ln\varepsilon|\,\partial_\theta \mathcal Q_2[h](\varphi,\theta)
+\varepsilon^3|\ln\varepsilon|\,\partial_\theta \mathcal R_2(\varepsilon,f,\mathtt{U}_-)[h](\varphi,\theta),
\end{align*}
where the operator $\mathcal{Q}_2$ is given by
\begin{align}\label{def:Q2}
      \mathcal{Q}_2[h](\varphi,\theta)&:= \tfrac{1}{\sqrt{\mathtt{p}_{1,1}(\varphi)}}|\ln\varepsilon|^{-1}\fint_{\mathbb{T}} h_\star(\varphi,\eta) \nabla_{1,2}^2G(\mathtt{P}_1,\mathtt{P}_2)[\mathcal{Z}(\varphi,\theta),\mathcal{Z}_\star(\varphi,\eta)]d\eta,
  \end{align}
  and the operator $\mathcal R_2$ is smoothing: for every $N\in\N$ and every $s\ge s_0$,
  $$
  \|\partial_\theta^N\mathcal{R}_2[h]\|^{\mathrm{Lip},\nu}_s
\lesssim (1+\|f\|_{s_0+N}^{\mathrm{Lip},\nu})\|h\|_s^{\mathrm{Lip},\nu}
+\|f\|_{s+N+1}^{\mathrm{Lip},\nu}\|h\|_{s_0}^{\mathrm{Lip},\nu}.
$$
  Moreover,
\begin{align*}
\tfrac{1}{\sqrt{\mathtt{p}_{1,1}(\varphi)}}U_2(\varepsilon,f)(\varphi,\theta)
&=\tfrac{\varepsilon}{2\sqrt{\mathtt{p}_{1,1}(\varphi)}} \,\, \nabla_1^\perp G(\mathtt{P}_1(\varphi),\mathtt{P}_2(\varphi))\cdot i\, \tfrac{\mathcal{Z}(\varphi,\theta)}{w(\varphi,\theta)}+\varepsilon^2|\ln\varepsilon|\mathtt{f}_2(\varphi) \sin(2\theta)
\\ &\quad+\varepsilon^2|\ln\varepsilon|\Big(\tfrac{3}{16}(2\mathtt{p}_{1,1}(\varphi))^{-\frac32}-\mathtt{g}_2(\varphi)\Big)\cos(2\theta)  
 +\varepsilon^2 \widetilde{U}_2(\varepsilon,f,\mathtt{U}_-)(\varphi,\theta),
\end{align*}
where the functions $\mathtt{f}_2$, $\mathtt{h}_2$, $\mathtt{g}_2$ are given by \eqref{list-functions} and the function $\widetilde{U}_2(\varepsilon,f,\mathtt{U}_-)$ is smooth and satisfies 
$$
   \|\widetilde{U}_2(\varepsilon,f,\mathtt{U}_-)\|^{\textnormal{Lip},\nu}_{s}\lesssim 1+ \|f\|^{\textnormal{Lip},\nu}_{s}.
   $$

\end{lem}
\begin{proof}
We first expand the kernel of $\mathcal I_2$. A second-order Taylor formula at $(\mathtt{P}_1,\mathtt{P}_2)$ yields
\begin{align*}
 \nonumber G\big(\mathtt{P}_1&+\varepsilon w \mathcal{Z}+i\tfrac{\varepsilon}{|\ln\varepsilon|}\mathtt{U}_{-},\mathtt{P}_2+\varepsilon  w_\star   \mathcal{Z}_\star\big)= G(\mathtt{P}_1,\mathtt{P}_2)+\varepsilon w\nabla_1G(\mathtt{P}_1,\mathtt{P}_2)[\mathcal{Z}]+\varepsilon  w_\star \nabla_2G(\mathtt{P}_1,\mathtt{P}_2)[\mathcal{Z}_\star]
\\ &+\varepsilon |\ln\varepsilon|^{-1}\mathtt{U}_{-} \, \nabla_1G(\mathtt{P}_1,\mathtt{P}_2)[i]+\tfrac12\varepsilon^2w^2\nabla_1^2G(\mathtt{P}_1,\mathtt{P}_2)[\mathcal{Z},\mathcal{Z}]+\tfrac12\varepsilon^2 w_\star^2\nabla_2^2G(\mathtt{P}_1,\mathtt{P}_2)[\mathcal{Z}_\star,\mathcal{Z}_\star]
 \\ \nonumber&+\varepsilon^2 |\ln\varepsilon|^{-1}\mathtt{U}_{-} w \, \nabla_1^2G(\mathtt{P}_1,\mathtt{P}_2)[i,\mathcal{Z}] +\tfrac12 \varepsilon^2 |\ln\varepsilon|^{-2}\mathtt{U}_{-}^2 \, \nabla_1^2G(\mathtt{P}_1,\mathtt{P}_2)[i,i] \\
&+\varepsilon^2 w w_\star\,\nabla_{1,2}^2G(\mathtt{P}_1,\mathtt{P}_2)[\mathcal{Z},\mathcal{Z}_\star]+\varepsilon^2 |\ln\varepsilon|^{-1}\mathtt{U}_{-}\,  w_\star \nabla_{1,2}^2G(\mathtt{P}_1,\mathtt{P}_2)[i,\mathcal{Z}_\star]+{\varepsilon^3}|\ln\varepsilon|^{\frac32}\mathcal{R}(\varepsilon,f,\mathtt{U}_-),
\end{align*}
with 
$$
\mathcal{R}(\varepsilon,f,\mathtt{U}_-):=\frac{1}{2}|\ln\varepsilon|^{-\frac32}\!\sum_{i,j,k=1}^2\int_0^1\!(1-t)^2\nabla^3_{i,j,k}G(P+t\varepsilon w\mathcal{Z}+i \tfrac{\varepsilon t}{|\ln\varepsilon|}\mathtt{U}_{-},P_\star+t\varepsilon  w_\star\mathcal{Z}_\star)[Z_i,Z_j,Z_k]dt,$$
and $Z_i,Z_j, Z_k\in\{w\mathcal{Z}+i|\ln\varepsilon|^{-1}\mathtt{U}_{-}, w_\star\mathcal{Z}_\star\}$. 
As the average in space of $h$ is zero,  then
\begin{align*}
\mathcal{I}_2(\varepsilon,f)[h](\varphi,\theta)
&=\partial_\theta\fint_{\T}G\big(\mathtt{P}_1+\varepsilon w(\theta) \mathcal{Z}(\theta)+i\varepsilon|\ln\varepsilon|^{-1}\mathtt{U}_{-},\mathtt{P}_2+\varepsilon w_\star(\eta)   \mathcal{Z}_\star(\eta)\big)h_\star(\eta) d\eta
\\
& =\varepsilon^2\partial_\theta\bigg(w(\theta)\fint_{\T} w_\star(\eta )\,\nabla_{1,2}^2G(\mathtt{P}_1,\mathtt{P}_2)[\mathcal{Z}(\theta),\mathcal{Z}_\star(\eta)]h_\star(\eta) d\eta\bigg)\\ &\quad +{\varepsilon^3}|\ln\varepsilon|^{\frac32}\partial_\theta\fint_{\T}\mathcal{R}(f)h_\star(\eta) d\eta.
\end{align*}
Moreover from the asymptotics
$$  w(\theta)=1+\varepsilon f(\theta)+O(\varepsilon^2),$$
we conclude that
\begin{align*}
\mathcal{I}_2(\varepsilon,f)[h](\varphi,\theta)
& =\varepsilon^2\partial_\theta\fint_{\T} \nabla_{1,2}^2G(\mathtt{P}_1,\mathtt{P}_2)[\mathcal{Z}(\theta),\mathcal{Z}_\star(\eta)]h_\star(\eta) d\eta \ +{\varepsilon^3}|\ln\varepsilon|\partial_\theta\mathcal{R}_2(\varepsilon,f,\mathtt{U}_-)[h](\theta),
\end{align*}
where
\begin{align*}
{\varepsilon^3}|\ln\varepsilon|\partial_\theta\mathcal{R}_2&(\varepsilon,f,\mathtt{U}_-)[h](\theta)=\varepsilon^2\partial_\theta\bigg(w(\theta)\fint_{\T} w_\star(\eta )\,\nabla_{1,2}^2G(\mathtt{P}_1,\mathtt{P}_2)[\mathcal{Z}(\theta),\mathcal{Z}_\star(\eta)]h_\star(\eta) d\eta\\ &-\fint_{\T} \nabla_{1,2}^2G(\mathtt{P}_1,\mathtt{P}_2)[\mathcal{Z}(\theta),\mathcal{Z}_\star(\eta)]h_\star(\eta) d\eta\bigg) +{\varepsilon^3}|\ln\varepsilon|^{\frac32}\partial_\theta\fint_{\T}\mathcal{R}(\varepsilon,f,\mathtt{U}_-)h_\star(\eta) d\eta.
\end{align*}
Notice that  the operator $\mathcal{R}_2(\varepsilon,f,\mathtt{U}_-)\in \textnormal{OPS}^{-\infty}$ is smoothing and satisfies  the claimed estimate.
\\
We now turn to $V_2$. Expanding $\nabla_1G$ at $(P_1,P_2)$ gives
\begin{align*}
 \nonumber &\varepsilon\nabla_1G\big(\mathtt{P}_1+\varepsilon w \mathcal{Z}+i\varepsilon|\ln\varepsilon|^{-1}\mathtt{U}_{-},\mathtt{P}_2+\varepsilon \rho w_\star   \mathcal{Z}_\star\big)
 =\varepsilon \nabla_1G(\mathtt{P}_1,\mathtt{P}_2)+\varepsilon^2 w\nabla_1^2 G(\mathtt{P}_1,\mathtt{P}_2)[\mathcal{Z},\cdot]\\ &\quad +\varepsilon^2 |\ln\varepsilon|^{-1}\mathtt{U}_{-} \, \nabla_1^2 G(\mathtt{P}_1,\mathtt{P}_2)[i,\cdot]+\varepsilon^2 \rho w_\star \nabla_{1,2}^2G(\mathtt{P}_1,\mathtt{P}_2)[\cdot,\mathcal{Z}_\star]+ O\left(\varepsilon^3|\ln\varepsilon|^\frac32\right).
\end{align*}
Substituting into the expression of $U_2(\varepsilon,f)(\varphi,\theta)$, and using
$$
\fint_{\T}w_\star(\varphi,\eta)^2\,d\eta=1,
\qquad
\fint_{\T}\mathcal Z_\star(\varphi,\eta)\,d\eta=0,
$$
we obtain
\begin{align*}
U_2(\varepsilon,f)(\varphi,\theta)
&=\frac12 \varepsilon\nabla_1G(\mathtt{P}_1,\mathtt{P}_2)\cdot \frac{\mathcal{Z}(\varphi,\theta)}{w(\varphi,\theta)}+\frac12\varepsilon^2\nabla_1^2 G(\mathtt{P}_1,\mathtt{P}_2)[\mathcal{Z}(\theta),\mathcal{Z}(\theta)] +\varepsilon^2 \widetilde{U}_2(\varepsilon,f,\mathtt{U}_-)(\varphi,\theta).
\end{align*}
Here, the mixed term involving $\nabla_{1,2}^2G$ is absorbed into the remainder because its
leading part vanishes by the zero average of $\mathcal Z_\star$.
Direct computations lead to 
\begin{align*}
  \nonumber \nabla_{1}^2G(\mathtt{P}_1,\mathtt{P}_2)[\mathcal{Z}(\theta),\mathcal{Z}(\theta)]&=(2\mathtt{p}_{1,1})^{\frac12}(\partial_{\mathtt{p}_{1,1}}^2G)\cos^2(\theta)+(2\mathtt{p}_{1,1})^{-\frac12}(\partial_{\mathtt{p}_{1,2}}^2G)\sin^2(\theta) \\ &\quad \nonumber +2(\partial_{\mathtt{p}_{1,1}}\partial_{\mathtt{p}_{1,2}}G)\,\cos(\theta)\sin(\theta)\\
  \nonumber &=\tfrac12\Big[(2\mathtt{p}_{1,1})^{\frac12}(\partial_{\mathtt{p}_{1,1}}^2G)-(2\mathtt{p}_{1,1})^{-\frac12}(\partial_{\mathtt{p}_{1,2}}^2G)\Big] \cos(2\theta)\\ &\quad +(\partial_{\mathtt{p}_{1,1}}\partial_{\mathtt{p}_{1,2}}G)\,\sin(2\theta)\\
   &=2\sqrt{\mathtt{p}_{1,1}}|\ln\varepsilon|\left\{\tfrac{3}{16(2\mathtt{p}_{1,1})^\frac32}-\mathtt{g}_2\right\}\cos(2\theta)+2\sqrt{\mathtt{p}_{1,1}}\mathtt{f}_2|\ln\varepsilon| \sin(2\theta),
\end{align*}
where we used \eqref{list-functions} in the last identity. Therefore,
\begin{align*}
\tfrac{1}{\sqrt{\mathtt{p}_{1,1}(\varphi)}}U_2(\varepsilon,f)(\varphi,\theta)
&=\tfrac{\varepsilon}{2\sqrt{\mathtt{p}_{1,1}}} \,\, \nabla G(\mathtt{P}_1,\mathtt{P}_2)\cdot \, \tfrac{\mathcal{Z}(\theta)}{w(\theta)}+\varepsilon^3|\ln\varepsilon|\,\Big(\tfrac{3}{16}(2\mathtt{p}_{1,1})^{-\frac32}-\mathtt{g}_2\Big)\cos(2\theta)
\\ &\quad +\varepsilon^3|\ln\varepsilon|\mathtt{f}_2 \sin(2\theta) 
 +\varepsilon^2 \widetilde{U}_2(\varepsilon,f,\mathtt{U}_-)(\varphi,\theta).
\end{align*}
Since $\nabla_1G\cdot \mathcal Z=\nabla_1^\perp G\cdot i\,\mathcal Z$, this is exactly the
desired expansion.
\end{proof}
\subsection{Asymptotic expansion of the linearized operator}\label{sec-Asym-Line-op}

The goal of this section is  twofold. First, we identify the 
principal building blocks of the linearization by introducing a set of 
auxiliary operators that naturally arise in the decomposition of the nonlinear 
effects. These include the toroidal Hilbert transform, which captures the 
nonlocal singular structure of the kernel, a projection operator localizing 
on the first Fourier mode, and a shift operator that encodes certain coupling 
mechanisms between harmonics. Second, we combine these ingredients to derive 
a precise asymptotic description of the linearized operator, organized in 
powers of $\varepsilon$ and $|\ln \varepsilon|$.\\  
Define a \emph{projection localized on the spatial mode~$1$} as  
\begin{align}\label{operator-mode1}
      \mathcal{Q}[h](\varphi,\theta):=&-\tfrac{1}{8 }(2\mathtt{p}_{1,1})^{-\frac32}\Bigg( \cos\theta\fint_{\mathbb{T}} h(\eta) \cos\eta d\eta+3\sin\theta\fint_{\mathbb{T}} h(\eta) \sin\eta d\eta\Bigg)\nonumber\\ 
      &+\tfrac{1}{\sqrt{\mathtt{p}_{1,1}}}|\ln\varepsilon|^{-1}\fint_{\mathbb{T}} h_\star(\eta) \nabla_{1,2}^2G(\mathtt{P}_1,\mathtt{P}_2)[\mathcal{Z}(\theta),\mathcal{Z}_\star(\eta)]d\eta.
  \end{align}
We also introduce the functions
\begin{align}\label{V-001}
  V_0(\varepsilon,f)(\varphi,\theta)&:=\mathtt{g}_{3}(\varphi)\big(\cos(3\theta)-2\cos(\theta)\big)-\varepsilon|\ln\varepsilon| \mathtt{g}_{2}(\varphi)\cos(2\theta)\\
  \nonumber &\quad +\varepsilon|\ln\varepsilon| \left(\mathtt{f}_{2}(\varphi)+ \frac{\omega \,\dot{\mathtt{p}}_{1,1}}{4\mathtt{p}_{1,1}}\right)\sin(2\theta)+\tfrac12 f(\varphi,\theta),
\end{align}
where the auxiliary functions $\mathtt{f}_{2}, \mathtt{g}_{2}$ and $ \mathtt{g}_{3}$  are defined in \eqref{list-functions}. Also note
\begin{equation*}
V_0(\varepsilon,f)=V_0(\varepsilon,0)+\tfrac12 f.
\end{equation*}
The main result of this section deals with the asymptotic structure of the linearized operator around a small state.

\begin{pro}\label{prop:asymp-lin}
Assume \eqref{cond1} and adopt the notation of 
\eqref{ball-space}--\eqref{F-def}. Let $N \in \mathbb{N}$. 
Then there exists $\varepsilon_0 > 0$ such that, for all 
$\varepsilon \in (0,\varepsilon_0)$ and every $ f$ smooth with zero average in space and satisfying 
		\begin{equation}\label{smallnessf}
		\|f\|_{s_0}^{{\textnormal{Lip},\nu}}\leqslant 1,
		\end{equation}
the following hold.
\begin{enumerate}
\item The linearized operator satisfy the expansion
  \begin{align*}
    \qquad  \partial_f {\bf F}(\varepsilon,\mathtt{V}_1,\mathtt{V}_2,f)[h]&=\varepsilon^3 |\ln\varepsilon|\omega\partial_\varphi h+ \partial_\theta\big\{\big( \tfrac{1}{2}\varepsilon+  \varepsilon^2 V(\varepsilon,f,\mathtt{U}_-)\big)\, h\big\}\\ & \quad  +\tfrac{1}{2}\varepsilon \mathcal{H}[h]+\varepsilon^2\mathcal{H}_{\mathtt{u},0}[h]+\varepsilon^2{\partial_\theta}\mathcal{S}[h]+{\varepsilon^3}|\ln\varepsilon|  {\partial_\theta}\mathcal{Q}[h] \\ & \quad +\varepsilon^3\partial_\theta\fint_{\mathbb{T}} W({f})(\varphi,\theta,\eta)\ln\big|\sin\big(\tfrac{\theta-\eta}{2}\big)\big|h(\eta) d\eta+\varepsilon^3\partial_\theta{R}_{\infty}(\varepsilon,f,\mathtt{U}_-)[h],
  \end{align*}
where the operators $\mathcal{H}$, $\mathcal{H}_{\mathtt{u},0}$, $\mathcal{S}$, $\mathcal{Q}[h]$ are defined by \eqref{Hilbert1alt}, \eqref{def-Hu0}, \eqref{shift-operator1}, \eqref{operator-mode1}, respectively, the function $W(f):\mathcal{O}\times\mathbb{T}^3\to\mathbb{R}$ is given by Proposition \ref{prop-induc25} and    $V:\mathcal{O}\times\mathbb{T}^2\to\mathbb{R}$ is smooth and even, with 
\begin{align*}
V(\varepsilon,f,\mathtt{U}_-)&:=V_0(\varepsilon,f)+\varepsilon {V}_{-1}(\varepsilon,f,\mathtt{U}_-), \quad  \|V_{-1}(\varepsilon,f,\mathtt{U}_-)\|^{\textnormal{Lip},\nu}_s\lesssim 1+\|f\|^{\textnormal{Lip},\nu}_s.
\end{align*}
In addition, the  operator $\partial_\theta{R}_{\infty}(\varepsilon,f,\mathtt{U}_-)\in \textnormal{OPS}^{-\infty}$ is reversible and satisfies  the estimates
$$
\quad \|\partial_\theta^N{R}_{\infty}(\varepsilon,f,\mathtt{U}_-)[h]\|^{\textnormal{Lip},\nu}_{s}\lesssim (1+\|f\|^{\textnormal{Lip},\nu}_{s_0+N})\|h\|_s +\|f\|^{\textnormal{Lip},\nu}_{s+1+N}\|h\|^{\textnormal{Lip},\nu}_{s_0}.
$$
\item The higher  derivatives take the form
 \begin{align*}
      d_{f}^2 {\bf F}(\varepsilon,\mathtt{V}_1,\mathtt{V}_2,f)[h_1,h_2]=&\tfrac12 \varepsilon^2 \partial_\theta\big\{ h_1h_2\big\}+\varepsilon^3 \mathcal{E}_2(\varepsilon,f)[h_1,h_2],
  \end{align*}
and 
\begin{align*}
      &d_{f}^3 {\bf F}(\varepsilon,\mathtt{V}_1,\mathtt{V}_2,f)[h_1, h_2, h_3]=\varepsilon^{3} \mathcal{E}_3(\varepsilon,f)[h_1,h_2,h_3],
  \end{align*}
  where
\begin{align*}
			\| \mathcal{E}_2(\varepsilon,f)[h_1,h_2]\|_{s}^{{\textnormal{Lip},\nu}}&\lesssim\|h_{1}\|_{s_{0}+2}^{{\textnormal{Lip},\nu}}\|h_{2}\|_{s+2}^{{\textnormal{Lip},\nu}}+\Big(\|h_{1}\|_{s+2}^{{\textnormal{Lip},\nu}}+\| f\|_{s+2}^{{\textnormal{Lip},\nu}}\|h_{1}\|_{s_{0}+2}^{{\textnormal{Lip},\nu}}\Big)\|h_{2}\|_{s_{0}+2}^{{\textnormal{Lip},\nu}},\\
			\| \mathcal{E}_3(\varepsilon,f)[h,h,h]\|_{s}^{{\textnormal{Lip},\nu}}&\lesssim \big(\|h\|_{s_{0}+2}^{{\textnormal{Lip},\nu}}\big)^2\Big(\|h\|_{s+2}^{{\textnormal{Lip},\nu}}+\|f\|_{s+2}^{{\textnormal{Lip},\nu}}\|h\|_{s_0+2}^{{\textnormal{Lip},\nu}}\Big).
			\end{align*}	
\end{enumerate}

\end{pro}
\begin{proof}
${\bf{(1)}}$  According to Lemma  \ref{prop:asymp-lin0} the linearized operator of the functional ${\bf F}$ with respect to  $f$ in the direction $h$ is given by
 \begin{align*}
\nonumber  \partial_f {\bf F}(\varepsilon,\mathtt{V}_1,\mathtt{V}_2,f)[h] (\varphi,\theta)&=\varepsilon^3 |\ln\varepsilon | \omega \partial_\varphi h(\varphi ,\theta) +\varepsilon\,    \,\partial_\theta\big\{ h(\varphi,\theta)U(\varepsilon,f)(\varphi,\theta)\big\}\\ &\quad +\tfrac{\varepsilon}{\sqrt{\mathtt{p}_{1,1}(\varphi)}}  \,\big(\mathcal{I}_1(f)[h]+\mathcal{I}_2(f)[h]\big)(\varphi,\theta), 
\end{align*}
with
\begin{align*}
U(\varepsilon,f)(\varphi,\theta)&=\tfrac{1}{\sqrt{\mathtt{p}_{1,1}(\varphi)}}\big( U_1(\varepsilon,f)+ U_2(\varepsilon,f)\big)(\varphi,\theta)-\varepsilon \omega  |\ln \varepsilon|\dot{\mathtt{P}_1}(\varphi)\cdot i\, \tfrac{\mathcal{Z}(\varphi,\theta)}{w(\varphi,\theta)}\\ & \quad  -i\varepsilon^2 \omega\,  \dot{\mathtt{V}}_1(\varphi)\cdot i\tfrac{\mathcal{Z}(\varphi,\theta)}{w(\varphi,\theta)} +\varepsilon^2 |\ln\varepsilon |\,  \tfrac{ \omega \,\dot{\mathtt{p}}_{1,1}(\varphi)}{4\mathtt{p}_{1,1}(\varphi)}\sin(2\theta). 
\end{align*}
In view of  Lemma  \ref{prop-int25} and Lemma \ref{prop-V1}, one has
\begin{align*}
U(\varepsilon,f)(\varphi,\theta)&=
 \varepsilon\bigg(\tfrac{1}{2\sqrt{\mathtt{p}_{1,1}}} \,\, \nabla_1^\perp G(\mathtt{P}_1,\mathtt{P}_2)+ \tfrac{1}{\sqrt{2\mathtt{p}_{1,1}}}(\nabla_{\mathtt{P}_1}^\perp \mathtt{G})(\mathtt{P}_1,\mathtt{P}_2)- \omega  |\ln \varepsilon|\dot{\mathtt{P}_1}\bigg)\cdot i\, \tfrac{\mathcal{Z}(\varphi,\theta)}{w(\varphi,\theta)} 
\\ &\quad-\frac12+\frac\varepsilon2 f(\varphi,\theta)+\varepsilon\mathtt{g}_3(\varphi)\left(\cos(3\theta)-2  \cos(\theta)\right)
-\varepsilon^2|\ln\varepsilon|\mathtt{g}_2(\varphi)\cos(2\theta)
\\ &\quad+\varepsilon^2|\ln\varepsilon|\Big(\mathtt{f}_2(\varphi)+ \tfrac{ \omega \,\dot{\mathtt{p}}_{1,1}(\varphi)}{4\mathtt{p}_{1,1}(\varphi)}\Big) \sin(2\theta) -i\varepsilon^2 \omega\,  \dot{\mathtt{V}}_1(\varphi)\cdot i\tfrac{\mathcal{Z}(\varphi,\theta)}{w(\varphi,\theta)}\\ &\quad+\varepsilon^2 \widetilde{U}_1
(\varepsilon,f)(\varphi,\theta) +\varepsilon^2 \widetilde{U}_2(\varepsilon,f,\mathtt{U}_-)(\varphi,\theta). 
\end{align*}
Recall from \eqref{eq-points-generalL} that 
\begin{align*}
\omega |\ln\varepsilon |\,  \dot{\mathtt{P}}_1(\varphi) -\tfrac{1}{2\sqrt{\mathtt{p}_{1,1}}}(\nabla_{1} ^\perp G)(\mathtt{P}_1,\mathtt{P}_2)-\tfrac{1}{\sqrt{2\mathtt{p}_{1,1}}}(\nabla_{\mathtt{P}_1}^\perp \mathtt{G})(\mathtt{P}_1,\mathtt{P}_2)=0.
\end{align*}
It follows that
\begin{align*}
U(\varepsilon,f)(\varphi,\theta)&=-\frac12+\varepsilon V_{0}(\varepsilon,f)(\varphi,\theta)+\varepsilon^2 V_{-1}(\varepsilon,f)(\varphi,\theta), 
\end{align*}
where $V_{0}(\varepsilon,f)$ is given by \eqref{V-001}, and
\begin{align*}
V_{-1}(\varepsilon,f,\mathtt{U}_-)(\varphi,\theta)&=\widetilde{U}_1
(\varepsilon,f)(\varphi,\theta) + \widetilde{U}_2(\varepsilon,f,\mathtt{U}_-)(\varphi,\theta)-i \omega\,  \dot{\mathtt{V}}_1(\varphi)\cdot i\tfrac{\mathcal{Z}(\varphi,\theta)}{w(\varphi,\theta)}.
\end{align*}
Using the decompositions in Proposition \ref{prop-induc25} and Lemma  \ref{prop-int25} we get 
\begin{align*}
&\tfrac{1}{\sqrt{\mathtt{p}_{1,1}}}  \,\big(\mathcal{I}_1(\varepsilon,f)[h]+\mathcal{I}_2(\varepsilon,f)[h]\big)
=\tfrac12\mathcal{H}[h]+\varepsilon\mathcal{H}_{\mathtt{u},0}[h]+\varepsilon\partial_\theta\mathcal{S}[h]+\varepsilon^2|\ln\varepsilon|\partial_\theta \mathcal{Q}[h]
\\
&+\varepsilon^2\partial_\theta\fint_{\T}W(f)(\varphi,\theta,\eta)\ln\big|\sin\big(\tfrac{\theta-\eta}{2}\big)\big|h(\eta) d\eta
+{\varepsilon^2}\partial_\theta{R}_{\infty}(\varepsilon,f,\mathtt{U}_-)[h](\theta),
\end{align*}
with
\begin{align*}
 \partial_\theta\mathcal{Q}[h]& = \partial_\theta\mathcal{Q}_1[h]+\partial_\theta \mathcal{Q}_2[h],\\
 \partial_\theta{R}_{\infty}(\varepsilon,f,\mathtt{U}_-)[h]&=\partial_\theta\mathcal{R}_1(\varepsilon,f)[h]+{\varepsilon}|\ln\varepsilon|^{\frac32}\partial_\theta\mathcal{R}_2(\varepsilon,f,\mathtt{U}_-)[h],
 \end{align*}
 where the operators $\partial_\theta\mathcal{Q}_1[h]$ and $\partial_\theta \mathcal{Q}_2[h]$ are given by \eqref{def:Q1} and \eqref{def:Q2}, respectively. 
Combining the previous computations, we obtain the asymptotic expression of the linearized operator appearing in the statement of this proposition. Finally, higher-order derivatives of ${\bf{F}}$ with respect to $f$ are obtained by differentiating the asymptotics of $\partial_f {\bf{F}}$, and the corresponding estimates follow directly. This completes the proof of \mbox{Proposition \ref{prop:asymp-lin}.}
\end{proof}

\section{Speed modulation and approximate solution}
 As revealed by the detailed linearization performed in Proposition~\ref{prop:asymp-lin}, there is a degeneracy associated with the first Fourier mode in the angular variable. This degeneracy reflects the symmetry invariance of Euler equations under vertical translations, and manifests itself through the vanishing of the leading-order operator on the sine component of the first mode. In order to restore solvability at the linear level, it is therefore necessary to exploit this symmetry by introducing a carefully designed modulation of the reference frame speed.

The purpose of this section is twofold: first, to eliminate from the nonlinear problem the problematic first sine mode through an appropriate modulation of the vertical drift; and second, to construct an approximate solution that captures the essential nonlinear interactions up to sufficiently high order in the small parameter $\varepsilon.$ This preparation is crucial for the Nash–Moser scheme developed in the subsequent sections, as it ensures that the linearized operator admits a tame right inverse once projected onto the appropriate functional subspaces.

\subsection{Suppressing  the first sine mode  by speed modulation}\label{sec:modeone}
In this section, we implement a speed modulation for $\mathtt{V}_1$ and $\mathtt{V}_2$ as functions of the shape parameter $f$, with the aim of removing the $\sin(\theta)$ Fourier mode from the full nonlinear equation ${\bf F}$ defined in \eqref{F-def}. The modulation is designed so that the resulting dynamics satisfy the compatibility conditions imposed on $\mathtt{V}_1$ and $\mathtt{V}_2$, as stated in Lemma~\ref{lem-red} and Proposition~\ref{prop-reversibility}. 
{To state our main result, we need to introduce the following functionals,
\begin{align*}
    &\mathscr{N}_1(f)(\varphi)
    =-\tfrac{1}{\omega}e^{\frac{1}{\omega}\int_0^\varphi\mathtt({f}_2+\mathtt{f}_{2,\star})(\tau)d\tau}\int_0^\varphi e^{-\frac{1}{\omega}\int_0^s(\mathtt{f}_2+\mathtt{f}_{2,\star})(\tau)d\tau}\big(\mathbf{A}_1-\mathbf{A}_{1,\star}\big)(s)ds,
\end{align*}
with
\begin{align}\label{A1-f}
   \nonumber \mathbf{A}_1(f)&:=2 (2\mathtt{p}_{1,1})^{-\frac14}\fint_{\T}\bigg( \partial_\theta\big\{\big({\mathtt{g}_{3}}\big[\cos(3\theta)-2\cos(\theta)\big]+\tfrac{1}{4}f\big) f\big\}+\mathcal{H}_{\mathtt{u},0}[f]+{\partial_\theta}\mathcal{S}[f]\bigg)\sin(\theta)d\theta\\
   &\quad +2\varepsilon|\ln\varepsilon|(2\mathtt{p}_{1,1})^{-\frac14}\fint_{\T}\bigg(\partial_\theta\big\{\big(\mathtt{f}_2\sin(2\theta)-\mathtt{g}_2\cos(2\theta)\big)f\big\}+\partial_\theta\mathcal{Q}[f]\bigg)\sin(\theta)d\theta.
\end{align}
Fore the different functions and operators used above, see \eqref{list-functions},\eqref{func-mathttu}, \eqref{shift-operator1} and \eqref{operator-mode1}. On the other hand, according to Section~\ref{sec-notat}, the notation $f_\star$ is defined for a $2\pi$-periodic function $f$ by
$$
f_\star(\varphi)=f(\varphi+\pi).
$$

\begin{pro}\label{prop:U}
Assume \eqref{cond1} and let  $f\in \textnormal{Lip}_\nu(\mathcal{O},H^s_{\circ,\textnormal{even}}(\T^2))$ satisfy \eqref{smallnessf}. Then there exist 
$\mathtt{V}_1(f), \mathtt{V}_2(f)\in \textnormal{Lip}_\nu(\mathcal{O},H^s_{\textnormal{odd}}(\R))$  such that  
$$
\Pi_{1,{\bf s}}\Big( {\bf F}\big(\varepsilon,\mathtt{V}_1(f), \mathtt{V}_2(f),f\big)\Big)=0.
$$
Moreover, the function $\mathtt{U}_{-}(f):=\mathtt{V}_1(f)- \mathtt{V}_2(f)$ satisfies 
\begin{align*}
    \|\mathtt{U}_{-}(f)- \mathscr{N}_1(f)\|_{s}^{\textnormal{Lip},\nu}\lesssim  \varepsilon|\ln\varepsilon|^{\frac12}.
\end{align*}
\end{pro}
\begin{proof}
We start by expanding 
${\bf F}$ in a Taylor series to obtain
$$
{\bf F}(\varepsilon,\mathtt{V}_1,\mathtt{V}_2,f)={\bf F}(\varepsilon,\mathtt{V}_1,\mathtt{V}_2, 0)+d_f{\bf F}(\varepsilon,\mathtt{V}_1,\mathtt{V}_2,0)[f]+\tfrac{1}{2}d^2_f {\bf F}(\varepsilon,\mathtt{V}_1,\mathtt{V}_2,0)[f,f]+\mathbf{Q}(f),
$$
where $\mathbf{Q}(f)$ is at least cubic in $f$.
From Corollary \ref{cor-F0}, by naming the error term of the first sine mode projection by $\mathscr{C}$, which could depend on $\mathtt{U}_-$, we get
\begin{align*}
   \Pi_{1,{\bf s}}{\bf F}(\varepsilon,\mathtt{V}_1,\mathtt{V}_2,0)=&\varepsilon^2 (2\mathtt{p}_{1,1})^\frac14\,\Big(    \omega\,\dot{\mathtt{V}}_1(\varphi)-\mathtt{U}_{-}(\varphi)  \mathtt{f}_2(\varphi)-\varepsilon|\ln\varepsilon|\mathtt{b}+\varepsilon{|\ln\varepsilon|^{\frac12}}\, {\mathscr{C}(\mathtt{U}_{-})}(\varphi)\Big)\sin(\theta), 
\end{align*}
where $\mathtt{b}$ is given by \eqref{definitionb}.
By Proposition \ref{prop:asymp-lin} we obtain 
\begin{align*}
   &\Pi_{1,{\bf s}}\Big(d_f{\bf F}(\varepsilon,\mathtt{V}_1,\mathtt{V}_2,0)[f]+{\tfrac12}d_f^2{\bf F}(\varepsilon,\mathtt{V}_1,\mathtt{V}_2,0)[f,f]+\mathbf{Q}(f)\Big)\\
   &=\varepsilon^2 \,\Pi_{1,{\bf s}}\Big(\mathtt{g}_{3}\partial_\theta\big\{\big(\cos(3\theta)-2\cos(\theta)\big) f\big\}+\mathcal{H}_{\mathtt{u},0}[f]+{\partial_\theta}\mathcal{S}[f]+\tfrac14\partial_\theta\big\{ f^2 \big\}\\ &\quad +\varepsilon|\ln\varepsilon|\partial_\theta\big\{\big(\mathtt{f}_2\sin(2\theta)-\mathtt{g}_2\cos(2\theta)\big)f+\mathcal{Q}[f]\big)\big\}+O\big(\varepsilon|\ln\varepsilon|^{\frac12}\big)\Big)\\
   &=:\varepsilon^2(2\mathtt{p}_{1,1})^\frac14\,\Big(\mathbf{A}_1(f)+\varepsilon|\ln\varepsilon|^{\frac12} \mathbf{A}_2(\mathtt{U}_{-}, f)\Big)\sin(\theta), 
\end{align*}
where we recall from \eqref{assumo-speed}, 
$$ \dot{\mathtt{V}}_2=\dot{\mathtt{V}}_{1,\star},\quad\mathtt{U}_{-}=\mathtt{V}_1-\mathtt{V}_2,\quad\hbox{and} \quad\mathtt{U}_{-,\star}=-\mathtt{U}_{-} .
$$
Consequently, the equation
\begin{align*}
   \Pi_{1,{\bf s}}{\bf F}(\varepsilon,\mathtt{V}_1,\mathtt{V}_2,f)&=0,
\end{align*}
is equivalent to 
\begin{align}\label{eq-VV1}
      \omega\,\dot{\mathtt{V}}_1(\varphi)= \mathtt{U}_{-}(\varphi)   \mathtt{f}_2(\varphi)+\varepsilon|\ln\varepsilon|\mathtt{b}{-\varepsilon|\ln\varepsilon|^\frac12\mathscr{C}(\mathtt{U}_{-})(\varphi)}-\mathbf{A}_1(f)(\varphi)-\varepsilon|\ln\varepsilon|^{\frac12} \mathbf{A}_2(\mathtt{U}_{-},f)(\varphi).
\end{align}
Let us now explain how to solve equation \eqref{eq-VV1} in the periodic setting. From the definition of $\mathtt{b}$ in \eqref{definitionb} and the symmetry of the points in Section \ref{sec:filamentdyn} we get that $\mathtt{b}_{\star}=\mathtt{b}$,  and  then define
\begin{equation}\label{Form-Ident5}
    \begin{aligned}
 \mathbf{B}_1(f)&:=   {-}\mathbf{A}_1(f){-}\mathbf{A}_{1,\star}(f)+2\varepsilon|\ln\varepsilon| \mathtt{b}, \quad\mathbf{B}_2(\mathtt{U}_{-},f):= {-}\mathbf{A}_2(\mathtt{U}_{-}, f){-}\mathbf{A}_{2,\star}(\mathtt{U}_{-}, f)  , \\  \mathbf{C}_1(f)&:= {-}  \mathbf{A}_1(f){+}\mathbf{A}_{1,\star}(f),\quad
 \mathbf{C}_2(\mathtt{U}_{-},f):={-} \mathbf{A}_2(\mathtt{U}_{-}, f){+}\mathbf{A}_{2,\star}(\mathtt{U}_{-},f){-\mathscr{C}(\mathtt{U}_{-})+\mathscr{C}_\star(\mathtt{U}_{-})}.
\end{aligned}
\end{equation}
Notice that the assumptions in \eqref{assumo-speed}  imply that
\begin{align}\label{assumo-speed-U2}
\dot{\mathtt{V}}_{2,\star}(\varphi)=\dot{\mathtt{V}}_1(\varphi),\quad \mathtt{U}_{-}(-\varphi)=-\mathtt{U}_{-}(\varphi),  \quad\hbox{and}\quad  \mathtt{U}_{-,\star}=-\mathtt{U}_{-}.
\end{align}
From \eqref{eq-VV1} and \eqref{assumo-speed-U2} we infer that
\begin{align}\label{Eq-U2}
\omega\,   \dot{\mathtt{U}}_{-}-(\mathtt{f}_2+\mathtt{f}_{2,\star})\, \mathtt{U}_{-}=\mathbf{C}_1(f)+\varepsilon|\ln\varepsilon|^{\frac12}\mathbf{C}_2(\mathtt{U}_{-},f).
\end{align}
Defining \begin{align}\label{f3-math}
\mathtt{f}_3(\varphi):=\frac{1}{\omega}\int_0^\varphi(\mathtt{f}_2+\mathtt{f}_{2,\star})(\tau)d\tau,
\end{align}
we obtain after  integration
\begin{align}\label{Fixed-point1} {\mathtt{U}_{-}}(\varphi)&=\tfrac{1}{\omega}e^{\mathtt{f}_3(\varphi)}\bigg(  \int_0^\varphi e^{-\mathtt{f}_3(s)}\mathbf{C}_1(f)(s)ds+{\varepsilon|\ln\varepsilon|^{\frac12}}\int_0^\varphi e^{-\mathtt{f}_3(s)}\mathbf{C}_2(\mathtt{U}_{-},f)(s)ds\bigg)\\ \nonumber &=:\mathscr{N}(\mathtt{U}_{-},f).
\end{align}
 The construction of a solution to this equation follows a fixed point argument that will be done in several steps. First, we need to show that the right hand side functional in \eqref{Fixed-point1} sends $2\pi-$periodic function to a $2\pi-$periodic function. For this purpose, we need to show that 
if  $\mathbf{A}:\mathbb{T}\to \RR,$ is periodic,  then the following function
 \begin{align}\label{G-phi-per}
 G(\varphi)=\int_0^\varphi e^{-\mathtt{f}_3(s)}\big[\mathbf{A}(s)-\mathbf{A}_{\star}(s)\big]ds,
 \end{align}
 is also periodic. Indeed, from \eqref{f3-math} and  Corollary \ref{cor-F0} we infer  that  $\mathtt{f}_2+\mathtt{f}_{2,\star}$ is periodic and odd. Therefore, the map $$s\in\R\mapsto e^{-\mathtt{f}_3(s)},
 $$ is $2\pi-$ periodic and even. Consequently, the function $G$ is periodic if and only if
 $$
 G(2\pi)=\int_0^{2\pi} e^{-\mathtt{f}_3(s)}\big[\mathbf{A}(s)-\mathbf{A}(s+\pi)\big]ds=0.
 $$
Using a change of variables, we infer
 \begin{align*}
     G(2\pi)&=\int_0^{2\pi} e^{-\mathtt{f}_3(s)}\mathbf{A}(s)ds-\int_{\pi}^{3\pi} e^{-\frac{1}{\omega}\int_0^{s-\pi}(\mathtt{f}_2+\mathtt{f}_{2,\star})(\tau)d\tau}\mathbf{A}(s)ds\\
   &=\int_0^{2\pi} e^{-\mathtt{f}_3(s)}\mathbf{A}(s)ds-\int_{0}^{2\pi} e^{-\frac{1}{\omega}\int_0^{s-\pi}(\mathtt{f}_2+\mathtt{f}_{2,\star})(\tau)d\tau}\mathbf{A}(s)ds. 
 \end{align*}
 Now,  we write
 \begin{align*}
    \int_0^{s-\pi}(\mathtt{f}_2+\mathtt{f}_{2,\star})(\tau)d\tau= \int_0^{s}(\mathtt{f}_2+\mathtt{f}_{2,\star})(\tau)d\tau+ \int_s^{s-\pi}(\mathtt{f}_2+\mathtt{f}_{2,\star})(\tau)d\tau.
 \end{align*}
As $(\mathtt{f}_2+\mathtt{f}_{2,\star})_\star=(\mathtt{f}_2+\mathtt{f}_{2,\star}),$ we deduce that
 $$
 \frac{d}{ds}\int_s^{s-\pi}(\mathtt{f}_2+\mathtt{f}_{2,\star})(\tau)d\tau =(\mathtt{f}_2+\mathtt{f}_{2,\star})(s-\pi)-(\mathtt{f}_2+\mathtt{f}_{2,\star})(s)=0.
 $$
 Then,  by the fact that  $(\mathtt{f}_2+\mathtt{f}_{2,\star})$ is odd, we get
 $$
 \int_s^{s-\pi}(\mathtt{f}_2+\mathtt{f}_{2,\star})(\tau)d\tau=\int_{\frac\pi2}^{-\frac\pi2}(\mathtt{f}_2+\mathtt{f}_{2,\star})(\tau)d\tau=0.
 $$
 It follows that
 \begin{align*}
    \int_0^{s-\pi}(\mathtt{f}_2+\mathtt{f}_{2,\star})(\tau)d\tau= \int_0^{s}(\mathtt{f}_2+\mathtt{f}_{2,\star})(\tau)d\tau.
 \end{align*}
 Hence  $G(2\pi)=0$ and thus the map $G$ is $2\pi-$ periodic.  
This allows to show that the functional ${ \mathscr{N}(\mathtt{U}_{-},f)}$  defined in \eqref{Fixed-point1}  stabilizes periodic functions. According to \eqref{Fixed-point1}, one writes  
\begin{align*}
{ \mathscr{N}(\mathtt{U}_{-},f)(\varphi)}=&\tfrac{1}{\omega}e^{\mathtt{f}_3(\varphi)} \int_0^\varphi e^{-\mathtt{f}_3(s)}\mathbf{C}_1(f)(s)ds+{\varepsilon|\ln\varepsilon|^{\frac12}}\tfrac{1}{\omega}e^{\mathtt{f}_3(\varphi)}\int_0^\varphi e^{-\mathtt{f}_3(s)}\mathbf{C}_2(\mathtt{U}_{-},f)ds\\
  =:&\mathscr{N}_1(f)(\varphi)+\varepsilon|\ln\varepsilon|^{\frac12}\mathscr{N}_2(f,\mathtt{U}_{-})(\varphi).
\end{align*}
Define the closed  balls
 $$ B_1^s:=\big\{f\in \textnormal{Lip}_\nu(\mathcal{O},H^s_{\circ,\textnormal{even}}(\T^2)), \quad \|f\|_{s}^{\textnormal{Lip},\nu}\leqslant 1\big\},
$$
and
$$
B_{2}^s:=\big\{\mathtt{U}_{-}\in \textnormal{Lip}_\nu(\mathcal{O},H_{\textnormal{odd}}^s(\T)),\quad\|\mathtt{U}_{-}-\mathscr{N}_1(f)\|_{s}^{\textnormal{Lip},\nu}\leqslant1  \big\}.
$$
Notice that for $f\in B_1^s$ we can show the existence of $\varepsilon_0>0$ small enough but  independent of $f$ such that for any $\varepsilon\in(0,\varepsilon_0)$ we have
$$\mathtt{U}_{-}\in B_2^s\mapsto\varepsilon|\ln\varepsilon|^{\frac12}\mathscr{N}(\mathtt{U}_{-},f)\in B_2^s,
$$
is a well-defined contracting map. Indeed, the oddness of the function $\mathscr{N}(\mathtt{U}_{-})$ follows from the reversibility of the functional defining the contour dynamics equation. The Lipschitz dependence of the map $\mathtt{U}_{-}\in B_2^s\mapsto\mathscr{N}(\mathtt{U}_{-},f)$  follows  from  straightforward analysis based on the fact that  the functional $\mathtt{U}_{-}\mapsto {\bf F}(\varepsilon,\mathtt{V}_1,\mathtt{V}_2,f)$ is  smooth enough. Therefore the fixed point problem \eqref{Fixed-point1}   admits a unique solution $\mathtt{U}_{-}=\mathtt{U}_{-}(f)\in B_2^s$ satisfying 
$$
\mathtt{U}_{-}=\mathscr{N}_1(f)+O(\varepsilon|\ln\varepsilon|^{\frac12}).
$$
In addition, the map $f\in B_2^s\mapsto \mathtt{U}_{-}(f)$ is smooth. Let's check the assumption
$$
\mathtt{U}_{-,\star}=-\mathtt{U}_{-},
$$
as claimed in \eqref{assumo-speed-U2}. The oddness of $\mathtt{U}_{-}$ is a consequence of the construction as it belongs to $B_2^s.$ For the second property, we define
$$
\mathtt{W}_{-}:=-\mathtt{U}_{-,\star},
$$
and we intend to prove  that
$$
 \mathtt{W}_{-}=\mathtt{U}_{-}
$$
By taking the $\star$ operation in  \eqref{Eq-U2}, we find that
\begin{align*}
\omega\,   \dot{\mathtt{W}}_{-}(\varphi)-(\mathtt{f}_2+\mathtt{f}_{2,\star})\,(\varphi) \mathtt{W}_{-}(\varphi)&=-\mathbf{C}_{1,\star}(f)(\varphi)-\varepsilon|\ln\varepsilon|^{\frac12}\mathbf{C}_{2,\star}(\mathtt{U}_{-},f)(\varphi).
\end{align*}
Applying \eqref{Form-Ident5}, using in particular that for a periodic function $f$ we have $(f_\star)_\star=f$, yields
$$
   \mathbf{C}_{j,\star}=   -\mathbf{C}_{j}.
 $$
 Therefore,
 \begin{align*}
\omega\,   \dot{\mathtt{W}}_{-}(\varphi)-(\mathtt{f}_2+\mathtt{f}_{2,\star})\,(\varphi) \mathtt{W}_{-}(\varphi)&=\mathbf{C}_{1}(f)(\varphi)+\varepsilon|\ln\varepsilon|^{\frac12}\mathbf{C}_{2}(\mathtt{U}_{-},f)(\varphi).
\end{align*}
From this we deduce that the periodic function ${\eta}_2:=\mathtt{W}_{-}-\mathtt{U}_{-}$ satisfies the linear ODE
$$
\omega\,   \dot{\eta}_2(\varphi)-(\mathtt{f}_2+\mathtt{f}_{2,\star})\,(\varphi) \eta_2(\varphi)=0,
$$
which can be solved explicitly as
$$
\eta_2(\varphi)=\mathtt{c}\, e^{\frac{1}{\omega}\int_0^\varphi(\mathtt{f}_2+\mathtt{f}_{2,\star})(\tau)d\tau}.
$$
As by construction $\langle \mathtt{U}_{-}\rangle=0$ we get $\langle \mathtt{U}_{-,\star}\rangle=0$ and therefore $\langle \eta_2\rangle=0$. It follows that $\mathtt{c}=0$ and $\eta_2=0$.  Hence 
$
\mathtt{W}_{-}=\mathtt{U}_{-}.$
Define now
\begin{align*}
\mathtt{U}_{+}:=\mathtt{{V}}_1+\mathtt{{V}}_2.
\end{align*}
Coming back to \eqref{eq-VV1} and \eqref{Form-Ident5}, then we can write
 \begin{align}\label{Eq-U1}
&\omega\,   \dot{\mathtt{U}}_{+}(\varphi)=\mathbf{B}_1(f) (\varphi)+\varepsilon|\ln\varepsilon|^{\frac12}\mathbf{B}_2(\mathtt{U}_{-},f)(\varphi).
\end{align}
By setting 
 \begin{align*}
\nonumber \,\mathtt{U}_0&=\tfrac{1}{2\omega}\langle \mathbf{B}_1\rangle+\tfrac{1}{2\omega}\varepsilon|\ln\varepsilon|^{\frac12} \langle \mathbf{B}_2(\mathtt{U}_{-},f)\rangle\\ &= -\tfrac1\omega\langle  \mathbf{A}_1\rangle+\tfrac{1}{\omega}\varepsilon|\ln\varepsilon|\langle\mathtt{b}_\varepsilon\rangle+\tfrac1\omega\varepsilon|\ln\varepsilon|^{\frac12} \langle  \mathbf{A}_2(\mathtt{U}_{-},f)\rangle,
\end{align*}
and integrating \eqref{Eq-U1}, we get
$$
\mathtt{U}_{+}(\varphi)=2\mathtt{U}_0\varphi+\mathtt{W}(\varphi),
$$
with $\mathtt{W}$ an odd $2\pi-$periodic function.
Hence, one has 
\begin{align*}
    {\mathtt{V}_1}&=\tfrac12({\mathtt{U}}_{+}+{\mathtt{U}}_{-})= \mathtt{U}_0\varphi+\mathtt{W}_{+}(\varphi),
\end{align*}
and
\begin{align*}
    {\mathtt{V}_2}&=\tfrac12({\mathtt{U}}_{+}-{\mathtt{U}}_{-})= \mathtt{U}_0\varphi+\mathtt{W}_{-}(\varphi),
\end{align*}
with $\mathtt{W}_{+}$ and $\mathtt{W}_{-}$ two  $2\pi-$periodic functions. This achieves the proof of the proposition.
\end{proof}
Given $\mathtt{V}_1$ and $\mathtt{V}_2$ as in Proposition \ref{prop:U} yields
$$
{\bf F}(\varepsilon, \mathtt{V}_1(f),\mathtt{V}_2(f),f)=\Pi_{1,{\bf s}}^c {\bf F}(\varepsilon, \mathtt{V}_1(f),\mathtt{V}_2(f),f),
 $$}
where the whole first sine mode is removed. Thus, equation \eqref{F-def} is equivalent to 
 \begin{align}\label{Eq-modify1}
{\bf F}_1(\varepsilon,f)
:=\Pi_{1,{\bf s}}^c {\bf F}(\varepsilon, \mathtt{V}_1(f),\mathtt{V}_2(f),f)=0, 
\end{align}
and we will work now with ${\bf F}_1$ instead of ${\bf F}$. For later  use, we need to introduce the following operators. For $j=1,2$, we define the operators
\begin{equation}\label{defMj}
 \begin{aligned}
     \mathcal{M}_j(f)[h](\varphi)&:=\tfrac{\mathtt{h}_2(\varphi) e^{\mathtt{f}_3(\varphi)}} {2\omega (2\mathtt{p}_{1,1})^{\frac14}}\, \int_0^\varphi e^{-\mathtt{f}_3(s)} \big(\mathtt{L}_j[h]-(\mathtt{L}_{j}[h])_{\star}\big)(s)ds +(2-j)\mathcal{M}_{3}[f,h],\\
     \mathcal{M}_{3}[f,h](\varphi)(\varphi)&:=\tfrac{\mathtt{h}_2(\varphi) e^{\mathtt{f}_3(\varphi)}} {2\omega (2\mathtt{p}_{1,1})^{\frac14}}\, \int_0^\varphi e^{-\mathtt{f}_3(s)} \big(\mathtt{Q}_{1}[f,h]-(\mathtt{Q}_{1}[f,h])_{\star}\big)(s)ds,
\end{aligned}   
\end{equation}
with
\begin{equation}\label{defL1Q1L2Q2}
\begin{aligned}
  &\mathtt{L}_1[h]:=-\frac{1}{8\mathtt{p}_{1,1}}\fint_{\T}\big(\cos(3\theta)-2\cos(\theta)\big) h(\theta)\cos(\theta)d\theta\\ &\qquad \qquad+2(2\mathtt{p}_{1,1})^{-\frac14}\fint_{\T}\big(\mathcal{H}_{\mathtt{u},0}[h](\theta)+{\partial_\theta}\mathcal{S}[h](\theta)\big)\sin(\theta)d\theta,\\ &\mathtt{L}_2[h]:=2(2\mathtt{p}_{1,1})^{-\frac14}\bigg(\fint_{\T}\partial_\theta\mathcal{Q}[h](\theta)\sin(\theta)d\theta -\fint_{\T}\big(\mathtt{f}_2\sin(2\theta)-\mathtt{g}_2\cos(2\theta)\big)h(\theta)\cos(\theta)d\theta\bigg),
\\ &\mathtt{Q}_1[f,h]:= -(2\mathtt{p}_{1,1})^{-\frac14}\fint_{\T} f(\theta)h(\theta) \cos(\theta)d\theta.  
\end{aligned}
\end{equation}
We, next examine the regularity properties and the linearization of the functional ${\bf F}_1$.
\begin{cor}\label{cor:asymp-linF1}
Assume \eqref{cond1} and let $ N\in\mathbb{N}$. There exists $\varepsilon_0>0$ such that for all $\varepsilon\in(0,\varepsilon_0)$, the functional $${\bf F}_1(\varepsilon,\cdot):B_{1,\textnormal{even}}(s_0)\cap \textnormal{Lip}_\nu(\mathcal{O},H^{s}_{\circ,\textnormal{even}}(\T^2))\to \textnormal{Lip}_\nu(\mathcal{O},H^{s-1}_{\circ,\textnormal{odd}}(\T^2)),
$$  is well-defined and of class $C^1$. In addition, we have 
  \begin{align*}
      \partial_f {\bf F}_1(\varepsilon,f)[h]
      &={\Pi_{1,{\bf s}}^c}\bigg(\varepsilon^3 |\ln\varepsilon|\omega\partial_\varphi h+ \tfrac{\varepsilon}{2} \big(\mathcal{H}[h]-\partial_\theta h\big)+\varepsilon^2 \partial_\theta\big( h V_1(\varepsilon,f)\big)+\varepsilon^2\mathcal{H}_{\mathtt{u},0}[h]\\ &\quad+\varepsilon^2{\partial_\theta}\mathcal{S}[h]{+\varepsilon^2 \mathcal{M}_1(f)[h]\cos(\theta)+\varepsilon^3|\ln\varepsilon| \mathcal{M}_2(f)[h]\cos(\theta)} +{\varepsilon^3}|\ln\varepsilon|  {\partial_\theta}\mathcal{Q}[h]\\ &\quad+\varepsilon^3\partial_\theta\fint_{\mathbb{T}} W({f})(\varphi,\theta,\eta)\ln\big|\sin\big(\tfrac{\theta-\eta}{2}\big)\big|h(\eta) d\eta +\varepsilon^3|\ln\varepsilon|^{\frac12}\partial_\theta \mathcal{R}_{1,\infty}[h]\bigg),
  \end{align*}
where 
\begin{align*}
   V_1(\varepsilon,f)&=V_0(\varepsilon,f)+O(\varepsilon),
\end{align*}
with
$$
\|V_1(f)\|^{\textnormal{Lip},\nu}_s\lesssim 1+\|f\|^{\textnormal{Lip},\nu}_s\quad\hbox{and}\quad \|W(f)\|^{\textnormal{Lip},\nu}_s\lesssim 1+\|f\|^{\textnormal{Lip},\nu}_{s}.
$$
Moreover,  $\mathcal{R}_{1,\infty}\in \textnormal{OPS}^{-\infty}$ and  satisfies  in particular the estimates
$$
\quad \|\partial_\theta^N\mathcal{R}_{1,\infty}h\|^{\textnormal{Lip},\nu}_{s}\lesssim \|f\|^{\textnormal{Lip},\nu}_{s_0+N}\|h\|_s +\|f\|^{\textnormal{Lip},\nu}_{s+1+N}\|h\|^{\textnormal{Lip},\nu}_{s_0}.
$$
 Furthermore, 
\begin{equation*}
\begin{aligned}
    \partial^2_f {\bf F}_1(\varepsilon,f)[h_1,h_2]
&=\tfrac12\varepsilon^2 {\Pi_{1,{\bf s}}^c}\partial_\theta \big(h_1h_2 \big) +\varepsilon^2 \mathcal{M}_3[h_1,h_2]\cos(\theta)\\ &\quad
+\varepsilon^3|\ln\varepsilon|^{\frac12} {\bf E}_2(\varepsilon,f)[h_1,h_2],
\end{aligned}
\end{equation*}
and
\begin{equation*}
\partial^3_f {\bf F}_1(\varepsilon,f)[h_1,h_2,h_3]
={\varepsilon^3 |\ln\varepsilon|^{\frac12}} {\bf E}_3(\varepsilon,f)[h_1,h_2,h_3],
\end{equation*}
  where
\begin{align*}
			&\| {\bf E}_2(\varepsilon,f)[h_1,h_2]\|_{s}^{{\textnormal{Lip},\nu}}\lesssim\|h_{1}\|_{s_{0}+2}^{{\textnormal{Lip},\nu}}\|h_{2}\|_{s+2}^{{\textnormal{Lip},\nu}}+\Big(\|h_{1}\|_{s+2}^{{\textnormal{Lip},\nu}}+\| f\|_{s+2}^{{\textnormal{Lip},\nu}}\|h_{1}\|_{s_{0}+2}^{{\textnormal{Lip},\nu}}\Big)\|h_{2}\|_{s_{0}+2}^{{\textnormal{Lip},\nu}},\\
            &\| {\bf E}_3(\varepsilon,f)[h_1,h_2,h_3]\|_{s}^{{\textnormal{Lip},\nu}}\lesssim\|h_{1}\|_{s_{0}+2}^{{\textnormal{Lip},\nu}}\|h_{2}\|_{s_0+2}^{{\textnormal{Lip},\nu}}\|h_{3}\|_{s+2}^{{\textnormal{Lip},\nu}}\\
            &+\Big(\|h_{1}\|_{s+2}^{{\textnormal{Lip},\nu}}\|h_{2}\|_{s_{0}+2}^{{\textnormal{Lip},\nu}}+\Big(\|h_{2}\|_{s+2}^{{\textnormal{Lip},\nu}}+\| f\|_{s+2}^{{\textnormal{Lip},\nu}}\|h_{2}\|_{s_{0}+2}^{{\textnormal{Lip},\nu}}\Big)\|h_{1}\|_{s_{0}+2}^{{\textnormal{Lip},\nu}}\Big)\|h_{3}\|_{s_{0}+2}^{{\textnormal{Lip},\nu}}.
            \end{align*}
Recall that the function $W$ was introduced in Proposition $\ref{prop:asymp-lin}$,  $V_0$ in \eqref{V-001}.
\end{cor}
\begin{proof}
From the choice of $\mathtt{V}_1(f)$ and $\mathtt{V}_2(f)$ the sine-mode 1 is absent from the nonlinear functional, which implies that the range of ${\bf F}_1$ is in $\textnormal{Lip}_\nu(\mathcal{O},H^{s-1}_\circ(\T^2))$. To get the $C^1$ regularity of the mapping we simply  use Proposition \ref{prop:asymp-lin} and Proposition \ref{prop:U}.\\ 
As ${\bf F}(\varepsilon, \mathtt{V}_1,\mathtt{V}_2,f)$ depends on $(\mathtt{U}_{-},f)$, then differentiating \eqref{Eq-modify1} with respect to $f$ in the direction $h$ gives
\begin{align*}
    \partial_f {\bf F}_1(\varepsilon,f)[h]=\Pi_{1,{\bf{s}}}^c (\partial_f {\bf F})(\varepsilon, \mathtt{V}_1,\mathtt{V}_2,f)[h]+\Pi_{1,{\bf{s}}}^c \partial_{\mathtt{U}_{-}} {\bf F}(\varepsilon, \mathtt{V}_1,\mathtt{V}_2,f){\partial_f \mathtt{U}_{-}}(\varepsilon,f)[h].
\end{align*}
The asymptotics of the first term $\partial_f {\bf F}$ are given by Proposition \ref{prop:asymp-lin}. For the second term, we shall first compute  $ \partial_f \mathtt{U}_{-}.$
Using Proposition \ref{prop:U}, we obtain
\begin{align*}
(\partial_f \mathtt{U}_{-})(f)[h]&=(\partial_f\mathscr{N}_1)(f)[h]+O\big(\varepsilon|\ln\varepsilon|^{\frac12}\big)\\
&=-\tfrac{1}{\omega}e^{\mathtt{f}_3(\varphi)}\int_0^\varphi e^{-\mathtt{f}_3(s)} \partial_f\big(\mathbf{A}_1(f)-\mathbf{A}_{1,\star}(f)\big)[h](s)ds+O\big(\varepsilon|\ln\varepsilon|^{\frac12}\big).
\end{align*}
Here $O\big(\varepsilon|\ln\varepsilon|^{\frac12}\big)$ denotes a smoothing operator of order $\varepsilon|\ln\varepsilon|^{\frac12}.$
Linearizing \eqref{A1-f} and integrating by parts yields 

\begin{align*}
\partial_f& \mathbf{A}_1(f)[h]
=\;-2 (2\mathtt{p}_{1,1})^{-\frac14}\fint_{\mathbb{T}}\Big(\mathtt{g}_{3}\big[\cos(3\theta)-2\cos(\theta)\big]+\tfrac12 f\Big)h\cos(\theta)\,d\theta \\
&+\;2 (2\mathtt{p}_{1,1})^{-\frac14}\fint_{\mathbb{T}}\big(
\mathcal{H}_{\mathtt{u},0}[h]
+\partial_\theta\mathcal{S}[h]\big)
\sin(\theta)\,d\theta +2\varepsilon|\ln\varepsilon|(2\mathtt{p}_{1,1})^{-\frac14}\fint_{\mathbb{T}}\mathcal{Q}[h]
\sin(\theta)\,d\theta\\
&-2\varepsilon|\ln\varepsilon|(2\mathtt{p}_{1,1})^{-\frac14}\fint_{\mathbb{T}}
\Big(\mathtt{f}_2\sin(2\theta)-\mathtt{g}_2\cos(2\theta)\Big)h
\cos(\theta)\,d\theta.
\end{align*}
By applying \eqref{defL1Q1L2Q2}, we find
\begin{align*}
 \partial_f \mathbf{A}_1(f)(\varphi)[h]&={\mathtt{L}_1[h]}+\varepsilon|\ln\varepsilon|\mathtt{L}_2[h]+\mathtt{Q}_1[f,h].
\end{align*}
 Thus, from \eqref{defMj}, we conclude that
\begin{align}\label{partialfU}
(\partial_f \mathtt{U}_{-})(f)[h]&=-{2}\tfrac{(2\mathtt{p}_{1,1})^{\frac14}}{\mathtt{h}_2}\Big(  \mathcal{M}_1[h]+\varepsilon|\ln\varepsilon| \mathcal{M}_2[h]\Big)+O\big(\varepsilon|\ln\varepsilon|^{\frac12}\big).
\end{align}
Note that the operator $h\mapsto \partial_f \mathtt{U}_{-}[h]$ is  smoothing at any order.
By Taylor expansion, we infer
\begin{align*}
    {\bf F}(\varepsilon, \mathtt{V}_1,\mathtt{V}_2,f)={\bf F}(\varepsilon, \mathtt{V}_1,\mathtt{V}_2,0)+\int_0^1( \partial_f {\bf F})(\varepsilon,\mathtt{V}_1,\mathtt{V}_2,sf)[f]ds.
\end{align*}
From Proposition \ref{prop:asymp-lin} and \eqref{partialfU} we obtain  
$$
 \Pi_{1,{\bf{s}}}^c \bigg(\partial_{\mathtt{U}_{-}}\int_0^1( \partial_f {\bf F})(\varepsilon,\mathtt{V}_1,\mathtt{V}_2,sf)[f]ds\bigg)\partial_f \mathtt{U}_{-}(\varepsilon,f)[h]={\varepsilon^3} \partial_\theta\mathscr{R}_{1,\infty}[h],
$$
where
$$
\quad \|\partial_\theta^N\mathscr{R}_{1,\infty}h\|^{\textnormal{Lip},\nu}_{s}\lesssim \|f\|^{\textnormal{Lip},\nu}_{s_0+N}\|h\|_s +\|f\|^{\textnormal{Lip},\nu}_{s+1+N}\|h\|^{\textnormal{Lip},\nu}_{s_0}.
$$
By Corollary \ref{cor-F0}, it follows that
\begin{align*}
    \Pi_{1,{\bf{s}}}^c \partial_{\mathtt{U}_{-}} {\bf F}(\varepsilon, \mathtt{V}_1,\mathtt{V}_2,0)\partial_f \mathtt{U}_{-}[h]&=-\varepsilon^2\tfrac12 (2\mathtt{p}_{1,1})^{-\frac14}\mathtt{h}_2\, \partial_f \mathtt{U}_{-}(\varepsilon,f)[h] \cos(\theta)+O\big(\varepsilon^3|\ln\varepsilon|^\frac12\big).
\end{align*}
Applying once again \eqref{partialfU}, we find
\begin{align*}
   \Pi_{1,{\bf{s}}}^c \partial_{\mathtt{U}_{-}} {\bf F}(\varepsilon, \mathtt{V}_1,\mathtt{V}_2,0)\partial_f \mathtt{U}_{-}[h]&=\varepsilon^2 \mathcal{M}_1[h]\cos(\theta)+\varepsilon^3|\ln\varepsilon| \mathcal{M}_2[h]\cos(\theta)\\ &\quad+\varepsilon^3|\ln\varepsilon|^\frac12\partial_\theta\mathcal{R}_{\mathcal{M},\infty}[h],
\end{align*}
where
$$
\quad \|\partial_\theta^N\mathcal{R}_{\mathcal{M},\infty} h\|^{\textnormal{Lip},\nu}_{s}\lesssim \|f\|^{\textnormal{Lip},\nu}_{s_0+N}\|h\|_s +\|f\|^{\textnormal{Lip},\nu}_{s+1+N}\|h\|^{\textnormal{Lip},\nu}_{s_0}.
$$
Putting everything together, we obtain
  \begin{align*}
      \partial_f {\bf F}_1(\varepsilon,f)[h]
      &={\Pi_{1,{\bf s}}^c}\bigg(\varepsilon^3 |\ln\varepsilon|\omega\partial_\varphi h+ \tfrac{\varepsilon}{2} \big(\mathcal{H}[h]-\partial_\theta h\big)+\varepsilon^2 \partial_\theta\big( h V_1(\varepsilon,f)\big)+\varepsilon^2\mathcal{H}_{\mathtt{u},0}[h]\\ &\quad+\varepsilon^2{\partial_\theta}\mathcal{S}[h]{+\varepsilon^2 \mathcal{M}_1[h]\cos(\theta)+\varepsilon^3|\ln\varepsilon| \mathcal{M}_2[h]\cos(\theta)} +{\varepsilon^3}|\ln\varepsilon|  {\partial_\theta}\mathcal{Q}[h]\\ &\quad+\varepsilon^3\partial_\theta\fint_{\mathbb{T}} W({f})(\varphi,\theta,\eta)\ln\big|\sin\big(\tfrac{\theta-\eta}{2}\big)\big|h(\eta) d\eta +\varepsilon^3|\ln\varepsilon|^{\frac12}\partial_\theta \mathcal{ R}_{1,\infty}[h]\bigg),
  \end{align*}
  with 
$\partial_\theta\mathcal{R}_{1,\infty}:=\partial_\theta{R}_{\infty}+\partial_\theta\mathcal{R}_{\mathcal{M},\infty}\in \textnormal{OPS}^{-\infty}$ and
\begin{align*}
{V}_1(\varepsilon,f)
&=V_0(\varepsilon,f)+O(\varepsilon),
\end{align*}
where we have used Proposition \ref{prop:U}. Finally, higher-order derivatives of ${\bf F}_1$ with respect to $f$ are obtained by differentiating the asymptotics of $\partial_f {\bf F}_1$, and the corresponding estimates follow directly. This completes the proof.
\end{proof}

\subsection{Approximate solution}\label{sec:approx}
In this section, we construct an approximate solution to the nonlinear equation \eqref{Eq-modify1} in the periodic setting, valid up to any prescribed polynomial order in the small parameter $\varepsilon$.  
The procedure is based on an inductive iterative scheme in which each step improves the accuracy of the approximation by canceling the residual at a higher order \mbox{in $\varepsilon$.}  
Starting from a suitably chosen leading-order profile, we generate successive corrections by solving linearized equations around the previously computed approximation.  
At each iteration, the error term is pushed to a strictly higher power of $\varepsilon$, so that after $n$ steps the residual is of size $\varepsilon^{n+3}|\ln\varepsilon|$.  
More precisely, we construct a sequence of profiles  
$$
f_n(\varphi,\theta) = \sum_{m=0}^n \varepsilon^m\big(g_m(\varphi)\cos\theta + h_m(\varphi,\theta)\big),
$$
satisfying both the prescribed symmetry conditions and a precise hierarchy of estimates on ${\bf F}_1(\varepsilon,f_n)$.  
The method relies on decomposing each perturbation into a normal component $h_m$ and a tangential component $g_m$: the former is obtained explicitly from the complementary projection $\Pi_1^c {\bf F}_1$, while the latter is determined by eliminating the mode-one term in the projection $\Pi_1 {\bf F}_1$, remind that the sine part is absent from the nonlinear functional ${\bf F}_1$.  
At each step, this decomposition leads to an ordinary differential equation in $\varphi$ for $g_m$, whose resolution exploits reversibility  properties and a contraction argument for small $\varepsilon$.  
This procedure, initiated for $n=0$ and then extended inductively, produces a profile $f_n$ that preserves the symmetry structure while improving the approximation order of the solution. To formulate our result we need to introduce some definitions. We set
\begin{align}\label{Def-Phij}
\varrho_1(\varphi):=e^{-\frac{T_0}{\pi}\int_0^\varphi\check{\alpha}(\tau)d\tau},\quad  \varrho_2(\varphi):=\tfrac{ T_0^2\check{\mathtt{h}}_2(\varphi) e^{\check{\mathtt{f}}_3(\varphi)}} {64\pi^2\kappa}e^{\frac{T_0}{\pi}\int_0^\varphi\check{\alpha}(\tau)d\tau},\quad\hbox{and}\quad \varrho_3(\varphi):=e^{-\check{\mathtt{f}}_3(s)},
\end{align}
together with the operator $\mathscr{T}:L^\infty([0,2\pi];\R)\to L^\infty([0,2\pi];\R)$ is defined by
\begin{align}\label{Def-TT}
 \mathscr{T}[g](\varphi):=\varrho_1(\varphi)\int_0^\varphi\varrho_2(\tau)\left(\int_0^\tau\varrho_3(s) g(s)ds\right) d\tau.
 \end{align}
 For the auxiliary functions involved in \eqref{Def-Phij}, we refer to Lemma \ref{lem-functions}.
From Lemma \ref{lem-fundamental}, the operator $\textnormal{Id}-\mathscr{T}$ has a bounded inverse. Now we define the open set
\begin{align}\label{Condi-Fund0}
\mathcal{O}_0&:=\Big\{(\lambda,\kappa)\in  (0,\infty)^2\quad\textnormal{s.t.}\quad \mathcal{P}(\lambda,\kappa)\neq 0\Big\},\\
\nonumber\mathcal{P}(\lambda,\kappa)&:=1+\int_0^{\pi}\varrho_2(\tau)\left(\int_0^\tau\varrho_3(s) (\textnormal{Id}-\mathscr{T})^{-1}[\varrho_1](s)ds\right)d\tau.
\end{align}
Now, we state our main result.
\begin{theo}\label{theo-approx1}
Assume \eqref{cond1} and let $(\lambda,\kappa) \in \mathcal{O}_0$ defined by \eqref{Condi-Fund0}.  For any $n\geqslant 1,$ there exists  {$\varepsilon_0=\varepsilon_{0}(n,\lambda,\kappa)>0$  such that for any  $\varepsilon\in(0,\varepsilon_{0})$ there exists } $f_n\in \textnormal{Lip}_\nu(\mathcal{O},H^{s}_{\circ,\textnormal{even}}(\T^2))$, for any $s\in\R,$ satisfying
$$
{\bf F}_1(\varepsilon,f_n)=O(\varepsilon^{3+n}|\ln\varepsilon|),\quad \Pi_{1,{\bf{c}}}{\bf F}_1(\varepsilon,f_n)=0=\Pi_{1,{\bf{s}}}{\bf F}_1(\varepsilon,f_n).
$$
Moreover,
$$
f_n(\varphi,\theta)=\sum_{m= 0}^n \varepsilon^m\big(g_m(\varphi)\cos(\theta)+h_m(\varphi,\theta)\big), 
$$
where all the $h_m$ are even with
\begin{align*}
  h_0(\varphi,\theta)=& \mathtt{g}_{3}(\varphi)\cos(3\theta)-2 \varepsilon|\ln\varepsilon|\mathtt{g}_{2}(\varphi)\cos(2\theta)+2\varepsilon|\ln\varepsilon|\mathtt{f}_{2}(\varphi)\sin(2\theta)+O(\varepsilon|\ln\varepsilon|^{\frac12}),
\end{align*}
and  all the $g_m$ are even with zero average.
The projections $\Pi_{m,{\bf{c}}}$ are  defined in \eqref{Pi-proj0}, see also \eqref{asymt-list1} for the involved functions. 
\end{theo}
\begin{proof}
We aim to prove by induction the existence of a sequence of smooth functions $$f_n\in \textnormal{Lip}_\nu\big(\mathcal{O},H^{s}_{\circ,\textnormal{even}}(\T^2)\big),$$
such that
 $$
 \Pi_{1,{\bf{c}}}{\bf F}_1(\varepsilon,f_{n})=0,\qquad \Pi_{1}^c{\bf F}(\varepsilon,f_n)=O(\varepsilon^{n+3}|\ln\varepsilon|).
 $$
{$\diamond$ \bf  Initialization ($n=0$).}
We first construct
 $f_0\in \textnormal{Lip}_\nu\big(\mathcal{O},H^{s}_{\circ,\textnormal{even}}(\T^2)\big)$
such that
\begin{align}\label{Initia-0}
\Pi_{1,{\bf{c}}}{\bf F}_1(\varepsilon,f_0)=0,\quad {\bf F}_1(\varepsilon,f_0)=O(\varepsilon^{3}|\ln\varepsilon|).
\end{align}
By Corollary~\ref{cor-F0} we have $\mathbf{F}(\varepsilon,0)=O(\varepsilon)$, and hence
\begin{align}\label{Fep0}
{\bf F}_1(\varepsilon,0)=O(\varepsilon).
\end{align}
We seek $f_0$ in the form
$$f_0(\varphi,\theta)=g_0(\varphi)\cos(\theta)+h_0(\varphi,\theta), \qquad \textnormal{with}\qquad \Pi_1^c h_0=h_0,$$ and impose the symmetry and normalization conditions
 $$
\int_{\mathbb{T}} g_0(\varphi)d\varphi=0,\quad  g_0(-\varphi)=g_0(\varphi),\quad\hbox{and}\quad h_0(-\varphi,-\theta)=h_0(\varphi,\theta).
 $$
 In the first step $g_0$ is fixed smooth function and satisfies
 $$
 \|g_0\|_{s_0}^{{\textnormal{Lip},\nu}}\leqslant 1.
 $$
According to Corollary \ref{cor:asymp-linF1}, the linearized operator satisfies
$$
\partial_f {\bf F}_1(\varepsilon,f)=\varepsilon^3 |\ln\varepsilon|\omega\partial_\varphi h+ \tfrac{\varepsilon}{2} \big(\mathcal{H}[h]-\partial_\theta h\big)+\varepsilon^2 {\Pi_{1,{\bf{s}}}^c} \mathfrak{L}(f)[h],
$$
where we have used the notation
\begin{align*}
   \nonumber   &\mathfrak{L}(f)[h]:= \partial_\theta\big( h V_1(\varepsilon,f)\big)+\mathcal{H}_{\mathtt{u},0}[h]+{\partial_\theta}\mathcal{S}[h]+{\varepsilon}|\ln\varepsilon|  {\partial_\theta}\mathcal{Q}[h]+\mathcal{M}_1(f)[h]\cos(\theta) \\ &\quad+\varepsilon|\ln\varepsilon|\mathcal{M}_2(f)[h] \cos(\theta)-\varepsilon\partial_\theta\fint_{\mathbb{T}} W({f})(\varphi,\theta,\eta)\ln\big|\sin\big(\tfrac{\theta-\eta}{2}\big)\big|h(\eta) d\eta+\varepsilon\partial_\theta\mathcal{R}_{1,\infty}[h].
  \end{align*}
Using the  Taylor expansion of $\mathbf{F}_1(\varepsilon,f)$ around $f=0$ up to the second order, we obtain
\begin{align}\label{FM0}
 \nonumber{\bf F}_1(\varepsilon,f_0)&=\Pi_{1,{\bf c}}{\bf F}_1(\varepsilon,0)+\Pi_{1,{\bf c}}\partial_f {\bf F}_1(\varepsilon,0)[f_0]+{\tfrac{1}{2}}\Pi_{1,{\bf{c}}} \partial^2_f {\bf F}_1(\varepsilon,0)[f_0,f_0]+\varepsilon^3{|\ln\varepsilon|^{\frac12}}\Pi_{1,{\bf c}} \mathcal{N}_0(f_0)\\
 \nonumber&\quad +\Pi_1^c{\bf F}_1(\varepsilon,0)+\varepsilon^3|\ln\varepsilon| \omega\partial_\varphi h_0+\tfrac{\varepsilon}{2}\{\mathcal{H}-\partial_\theta\}[h_0]+\varepsilon^2\Pi_1^c \mathfrak{L}(0)[f_0]\\
 &\quad+\tfrac{1}{2}\Pi_1^c \partial^2_f {\bf F}_1(\varepsilon,0)[f_0,f_0]+\varepsilon^3 {|\ln\varepsilon|^{\frac12}}\Pi_1^c\mathcal{N}_0(f_0),
\end{align}
where $\mathcal{N}_0(f_0)$ collects the terms that are at least cubic in $f_0$. 
The prefactor $\varepsilon^3 |\ln\varepsilon|^{\frac12}$ follows from the control of the third derivative given in Corollary~\ref{cor:asymp-linF1}. 
We also recall that $\partial_f^2 {\bf F}_1$ is of order $\varepsilon^2$, again by Corollary~\ref{cor:asymp-linF1}.
Our strategy is to first determine the normal component $h_0$, and then the tangential component $g_0$. We define $h_0$ as the solution of
\begin{align}\label{h0-L}
\Pi_1^c{\bf F}_1(\varepsilon,0)+\tfrac{\varepsilon}{2}\{\mathcal{H}-\partial_\theta\}[h_0]=- \Pi_{1}^c \big(\varepsilon^2 \mathfrak{L}(0)(f_0)+\tfrac12 \partial^2_f {\bf F}_1(\varepsilon,0)[f_0,f_0]\big).
\end{align}
Equivalently,
\begin{align}\label{Eq-NonL}
 \nonumber h_0&=2\varepsilon^{-1}\{\partial_\theta-\mathcal{H}\}^{-1}\Pi_1^c{\bf F}_1(\varepsilon,0)+\varepsilon   \{\partial_\theta-\mathcal{H}\}^{-1}\Pi_1^cQ_0(\varepsilon, g_0, h_0)\\
 &:=h_{0,\textnormal{in}}+\varepsilon Q_1(\varepsilon, g_0, h_0),
\end{align}
where $Q_0$ is a quadratic operator defined by
\begin{align*}
Q_0(\varepsilon, g_0, h_0)&:=\varepsilon^{-2}\partial^2_f {\bf F}_1(\varepsilon,0)\big[f_0,f_0\big]+\Pi_{1}^c  \mathfrak{L}(0)(f_0).
\end{align*}
Using Corollary~\ref{cor:asymp-linF1}, we find
$$
\sup_{0<\varepsilon\leqslant \varepsilon_0}\|Q_0(\varepsilon, g_0, h_0)\|_s^{{\textnormal{Lip},\nu}}\leqslant C(\|g_0\|_{s+1}^{{\textnormal{Lip},\nu}}+\|h_0\|_{s+1}^{{\textnormal{Lip},\nu}})^2.
$$
On the other hand
$$
\|\{\mathcal{H}-\partial_\theta\}^{-1}\Pi_1^c h\|_{s+1}^{{\textnormal{Lip},\nu}}\leqslant C\|h\|_s^{{\textnormal{Lip},\nu}}.
$$
Therefore
$$
\sup_{0<\varepsilon\leqslant \varepsilon_0}\|Q_1(\varepsilon, g_0, h_0)\|_s^{{\textnormal{Lip},\nu}}\leqslant C(\|g_0\|_{s}^{{\textnormal{Lip},\nu}}+\|h_0\|_{s}^{{\textnormal{Lip},\nu}})^2.
$$
By virtue of  \eqref{Fep0} we get
$$
\sup_{0<\varepsilon\leqslant \varepsilon_0}\|h_{0,\textnormal{in}}\|_s^{{\textnormal{Lip},\nu}}=\sup_{0<\varepsilon\leqslant \varepsilon_0}2\varepsilon^{-1}\|\{\partial_\theta-\mathcal{H}\}^{-1}\Pi_1^c{\bf F}_1(\varepsilon,0)\|_s^{{\textnormal{Lip},\nu}}<\infty.
$$
Consider the closed  ball $B_{1,\textnormal{in}}$ 
$$ B^s_{1,\textnormal{in}}:=\big\{h_0\in \textnormal{Lip}_\nu(\mathcal{O},H^{s}_{\perp,\textnormal{even}}(\T^2)), \quad \|h_0-h_{0,\textnormal{in}}\|_{s}^{\textnormal{Lip},\nu}\leqslant 1\big\}.
$$
The Lipschitz dependence of the map $h_0\in B^s_{1,\textnormal{in}}\mapsto Q_1\in \textnormal{Lip}_\nu(\mathcal{O},H^{s}_{\perp,\textnormal{even}}(\T^2))$ can be established in a straightforward way.
Therefore, for sufficiently small $\varepsilon$, we may apply the Banach fixed-point theorem with $s=s_0$ to deduce that the nonlinear equation \eqref{Eq-NonL} admits a unique even solution $h_0=h_0(\varepsilon,g_0)$. 
Moreover, by a persistence-of-regularity argument, this solution is in fact smoother and belongs to 
$
\textnormal{Lip}_\nu\big(\mathcal{O},H^{s}_{\perp,\textnormal{even}}(\T^2)\big),
$
for any $s$. In this way, we ensure that the smallness condition on $\varepsilon$ can be chosen independently of $s$.
 The expansion of $h_0$ at the leading order  follows from \eqref{Eq-NonL} and reads as
\begin{align}\label{h-0-four1}
   h_0(\varphi,\theta)&= 2\varepsilon^{-1}\{\partial_\theta-\mathcal{H}\}^{-1}\Pi_1^c{\bf F}_1(\varepsilon,0)+O(\varepsilon)\\
  \nonumber &=\mathtt{g}_{3}(\varphi)\cos(3\theta)-2 \varepsilon|\ln\varepsilon|\big(\mathtt{g}_{2}(\varphi)\cos(2\theta)-\mathtt{f}_{2}(\varphi)\sin(2\theta)\big)+\varepsilon|\ln\varepsilon|^\frac12 H_0(\varepsilon,g_0),
\end{align}
with, 
$$
 \|g_0\|_{s}^{{\textnormal{Lip},\nu}}\leqslant 1\Longrightarrow \sup_{0<\varepsilon\leqslant \varepsilon_0}\|H_0(\varepsilon,g_0)\|_{s}^{{\textnormal{Lip},\nu}}<\infty,
 $$
 and we have a smooth dependence of $g_0\in B_1^s\mapsto H_0(\varepsilon,g_0)\in \textnormal{Lip}_\nu(\mathcal{O},H^{s}_{\circ,\textnormal{even}}(\T)),$ where
  the closed  ball $B_1^s$ is defined by 
 \begin{align*}
 B_1^s:=\big\{g_0\in \textnormal{Lip}_\nu(\mathcal{O},H^{s}_{\circ,\textnormal{even}}(\T)), \quad \|g_0\|_{s}^{\textnormal{Lip},\nu}\leqslant 1\big\}.
\end{align*}
We point out that the space $H^{s}_{\circ,\textnormal{even}}(\T)$ used here is closely related to the spaces introduced in Section~\ref{SEc-function-sapces} defined on the torus $\T^2$. It consists of even functions on the one-dimensional torus that belong to $H^s(\T)$ and have zero average.\\
Another point used to obtain the error of size $\varepsilon |\ln \varepsilon|^{\frac12}$ in \eqref{h-0-four1} is the decomposition \eqref{p-decompos}, which allows us to obtain
$$
\left|{\tfrac{\omega \dot{\mathtt{p}}_{1,1}(\varphi)}{4\mathtt{p}_{1,1}(\varphi)}}\right|=O(|\ln\varepsilon|^{-\frac12}).
$$
 With the choice of $h_0=h_0(\varepsilon,g_0)$ dictated by \eqref{h0-L}, the identity \eqref{FM0} becomes
\begin{align}\label{first-app1}
 \nonumber{\bf F}_1(\varepsilon,f_0)=&\Pi_{1,{\bf c}}{\bf F}_1(\varepsilon,0)+\Pi_{1,{\bf c}}\partial_f {\bf F}_1(\varepsilon,0)[f_0]+{\tfrac{1}{2}}\Pi_{1,{\bf c}} \partial^2_f {\bf F}_1(\varepsilon,0)[f_0,f_0]+\varepsilon^3 \Pi_{1,{\bf c}} {|\ln\varepsilon|^{\frac12}}\mathcal{N}_0(f_0)\\
 &+\varepsilon^3|\ln\varepsilon| \omega\partial_\varphi h_0+\varepsilon^3 {|\ln\varepsilon|^{\frac12}}\Pi_1^c\mathcal{N}_0(f_0).
\end{align}
We now choose $g_0$ so as to eliminate the entire first  mode. More precisely, we impose
\begin{align}\label{First-mode1}
 \Pi_{1,{\bf c}}{\bf F}_1(\varepsilon,0)+\Pi_{1,{\bf c}}\partial_f {\bf F}_1(\varepsilon,0)[f_0]+{\tfrac{1}{2}}\Pi_{1,{\bf c}} \partial^2_f {\bf F}_1(\varepsilon,0)[f_0,f_0]+\varepsilon^3|\ln\varepsilon|^\frac12 \Pi_{1,{\bf c}} \mathcal{N}_0(f_0)=0.
\end{align}
By Corollary \ref{cor-F0} and Proposition \ref{prop:U}, we have
$$\Pi_{1,{\bf{c}}}{\bf F}_1(\varepsilon,0)=O\big({ \varepsilon^3}|\ln\varepsilon|^{\frac12}\big).
$$ 
Now, the Fourier expansion of $h_0$ in $\theta$ takes the form 
$$
h_0(\varphi,\theta)=\sum_{m\geqslant 2}\big(h_{0,m,c}(\varphi)\cos(\theta)+h_{0,m,s}(\varphi)\sin(\theta)\big).
$$
Then, from \eqref{h-0-four1}, we infer that the coefficients in the low modes satisfy
\begin{align*}
h_{0,2,s}&=2\varepsilon|\ln\varepsilon|\mathtt{f}_2+O(\varepsilon|\ln\varepsilon|^\frac12),\quad h_{0,2,c}=-2\varepsilon|\ln\varepsilon|\mathtt{g}_2+O(\varepsilon|\ln\varepsilon|^\frac12),\\
h_{0,3,s}&=O(\varepsilon|\ln\varepsilon|^\frac12), 
 \qquad h_{0,3,c}=\mathtt{g}_3+O(\varepsilon|\ln\varepsilon|^\frac12),\qquad h_{0,4,s}=O(\varepsilon|\ln\varepsilon|^\frac12),
\end{align*}
where the  $O(\varepsilon|\ln\varepsilon|^\frac12)$ remainders depend smoothly on $g_0$ (by \eqref{h-0-four1}).
Using Lemma~\ref{prop:dF-computations3}--(3) together with the asymptotics in Lemma \ref{lem-functions}, we obtain
\begin{align*}
 \Pi_{1,{\bf{c}}}(\partial_f {\bf F}_1)(\varepsilon,0)[f_0]&+\tfrac12 \Pi_{1,{\bf{c}}}(\partial_f^2 {\bf F}_1)(\varepsilon,0)[f_0,f_0]= \varepsilon^3 |\ln\varepsilon| \omega\Big( g_0^\prime+\tfrac{T_0 \check{\alpha}}{2\pi} (g_0- g_{0,\star}) \\ &+\,G_1[g_{0}-g_{0,\star}](\varphi)+|\ln\varepsilon|^{-\frac12}G_2(\varepsilon,\kappa, g_0,g_{0,\star})\Big)\cos(\theta),
\end{align*}
with 
\begin{align*}
  G_1[g](\varphi):= \tfrac{ T_0^2\check{\mathtt{h}}_2(\varphi) e^{\check{\mathtt{f}}_3(\varphi)}} {64\kappa\pi^2}\, \int_0^\varphi e^{-\check{\mathtt{f}}_3(s)} g(s)ds,
\end{align*}
where $\check{\mathtt{f}}_3$ is defined in \eqref{check-f3} and $\check{\mathtt{h}}_2$ in Lemma  \ref{lem-functions}. The functional  $G_2$ satisfies for each $\kappa>0$ and $\varepsilon_0>0$ small enough,
$$
 \|g_0\|_{s}^{{\textnormal{Lip},\nu}}\leqslant 1\Longrightarrow \sup_{0<\varepsilon\leqslant \varepsilon_0}\|G_2(\varepsilon,\kappa, g_0,g_{0,\star})\|_{s}^{{\textnormal{Lip},\nu}}<\infty,
 $$
 and we have a smooth dependence of $g_0\in B_1^s\mapsto G_2(\varepsilon,\kappa, g_0,g_{0,\star})\in \textnormal{Lip}_\nu(\mathcal{O},H^{s}_{\textnormal{odd}}(\T)).$ 
Then $g_0$ must solve the following ODE:
\begin{align}\label{Eq-Lun1}
    g_0^\prime+\tfrac{T_0}{2\pi}\check{\alpha}\, (g_0- g_{0,\star}) +\,G_1[g_{0}-g_{0,\star}]=|\ln\varepsilon|^{-\frac12}G_2(\varepsilon,\kappa, g_0,g_{0,\star}).
\end{align}
As in the study of \eqref{G-phi-per}, we can show that $G_1[g_{0}-g_{0,\star}]$ is $2\pi$ periodic. Moreover, by the reversibility of ${\bf{F}}_1$ and its regularity, we get 
\begin{align}\label{sym-rever}
g_0 \in H_{\textnormal{even}}^s(\T),\quad \langle g_0\rangle_\varphi=0\quad\Longrightarrow\quad  G_1[g_{0}-g_{0,\star}]\quad\hbox{and}\quad \,G_2(\varepsilon,\kappa,g_0,g_{0,\star})\in H_{\textnormal{odd}}^s(\T)\,.
\end{align}
Let us solve the equation in the periodic setting. To this end, define 
$$
g_{\pm}:=g_0\pm\,g_{0,\star}.
$$
From the last point in Corollary \ref{cor-F0}, we have $$\check{\alpha}_{\star}=\check{\alpha},\quad\,\check{\mathtt{h}}_{2,\star}=\check{\mathtt{h}}_{2},\quad\,\check{\mathtt{f}}_{3,\star}=\check{\mathtt{f}}_{3}.
$$It follows from \eqref{Eq-Lun1} that the unknown $g_{+}$ satisfies the ODE
\begin{align*}
 g_{+}^\prime&= G_1[g_{-}]+G_{1,\star}[g_{-}]+|\ln\varepsilon|^{-\frac12}\big(G_2(\varepsilon,\kappa, g_0,g_{0,\star})+G_{2,\star}(\varepsilon,\kappa,g_0,g_{0,\star})\big)\\
 &=:{G}_{1,+}[g_{-}] +|\ln\varepsilon|^{-\frac12}{G}_{2,+}(\varepsilon,\kappa,g_{+},g_{-}),\end{align*}
and  $g_{-}$ satisfies 
\begin{align*}
g_{-}^\prime+\tfrac{{\check{\alpha}\,T_0}}{\pi}  g_{-}&= \big(G_1[g_{-}]-G_{1,\star}[g_{-}]\big)+|\ln\varepsilon|^{-\frac12}\big(G_2(\varepsilon,\kappa,g_0,g_{0,\star})-G_{2,\star}(\varepsilon,\kappa,g_0,g_{0,\star})\big)\\
&=:{G}_{1,-}[g_{-}] +|\ln\varepsilon|^{-\frac12}{G}_{2,-}(\varepsilon,\kappa,g_{+},g_{-}).
\end{align*}
Integrating both equations, we  obtain  the equivalent integral formulation
\begin{align}\label{Eq--01}
  \nonumber g_{+}(\varphi)&= a+ \int_0^\varphi {G}_{1,+}[g_{-}](s) ds+|\ln\varepsilon|^{-\frac12}\int_0^\varphi {G}_{2,+}(\varepsilon,\kappa,g_{+},g_{-})(s) ds\\
&=:\mathcal{S}_+(g_+,g_-)(\varphi),
\end{align}
and
\begin{align}\label{g--Eq}
\nonumber g_{-}(\varphi)&= \mathscr{P}[g_{-}](\varphi)+|\ln\varepsilon|^{-\frac12}e^{-\frac{T_0}{\pi}\int_0^\varphi\check{\alpha}(\tau)d\tau}\left(b_2+\int_0^\varphi e^{\frac{T_0}{\pi}\int_0^s\check{\alpha}(\tau)d\tau}{G}_{2,-}(\varepsilon,\kappa,g_{+},g_{-})(s)ds\right)\\
&=:\mathscr{P}[g_{-}](\varphi)+|\ln\varepsilon|^{-\frac12}\mathcal{S}_-(g_+,g_-)(\varphi),
\end{align}
where
\begin{align}\label{mathSS}
 \mathscr{P}[g]:=e^{-\frac{T_0}{\pi}\int_0^\varphi\check{\alpha}(\tau)d\tau} \left(b_1+\int_0^\varphi e^{\frac{T_0}{\pi}\int_0^s\check{\alpha}(\tau)d\tau}{G}_{1,-}[g](s)ds\right),
\end{align}
and the constants $b_1$ and $b_2$ to be determined later.
We intend to solve this system with 
$$g_+\in \textnormal{Lip}_\nu(\mathcal{O},H^{s,+}_{\textnormal{even},\circ}(\T))\quad\hbox{and} \quad g_-\in \textnormal{Lip}_\nu(\mathcal{O},H^{s,-}_{\textnormal{even}}(\T)).
$$
Here, 
$$g_+\in \textnormal{Lip}_\nu(\mathcal{O},H^{s,+}_{\textnormal{even},\circ}(\T))\Longleftrightarrow g_+\in \textnormal{Lip}_\nu(\mathcal{O},H^{s}_{\textnormal{even},\circ}(\T))\quad\textnormal{and}  \quad g_{+,\star}=g_+
$$
and 
$$g_-\in \textnormal{Lip}_\nu(\mathcal{O},H^{s,-}_{\textnormal{even}}(\T))\Longleftrightarrow g_-\in \textnormal{Lip}_\nu(\mathcal{O},H^{s}_{\textnormal{even}}(\T))\quad\textnormal{and}  \quad g_{-,\star}=-g_{-}.
$$
From Corollary \ref{cor-F0}, we know that  $\check{\alpha}$ is periodic and odd, which implies that  $\displaystyle{\varphi\mapsto \int_0^\varphi\check{\alpha}(\tau)d\tau}$ is periodic and even.  Moreover, ${G}_{1,+}[g_{-}]$ and ${G}_{2,+}(\varepsilon,\kappa,g_{+},g_{-})$ are odd functions and invariant under the $\star-$transformation. Consequently,  $\mathcal{S}_{+}(g_+,g_-)$ is periodic, and  has zero average if and only if $a=0$. 
Hence \eqref{Eq--01} becomes
\begin{align}\label{Eq--02}
  \nonumber g_{+}(\varphi)&= \int_0^\varphi {G}_{1,+}[g_{-}](s) ds+|\ln\varepsilon|^{-\frac12}\int_0^\varphi {G}_{2,+}(\varepsilon,\kappa,g_{+},g_{-})(s) ds\\
&=\mathcal{S}_+(g_+,g_-)(\varphi),
\end{align}
Using the definition of $\mathcal{S}_{+}$, together with the fact that $({G}_{j,+})_\star={G}_{j,+}$, we readily verify that
 $$
\mathcal{S}_{+,\star}(g_+,g_-)=\mathcal{S}_{+}(g_+,g_-)\Longleftrightarrow \int_0^\pi {G}_{1,+}[g_{-}](s) ds=0=\int_0^\pi {G}_{2,+}(\varepsilon,\kappa,g_{+},g_{-})(s) ds.
$$
These identities follow from the fact that the functions ${G}_{j,+}$ are odd and satisfy $({G}_{j,+})_\star={G}_{j,+}$. Thus, their Fourier expansion contains only even sine modes, namely,
$$
{G}_{j,+}(\varphi)=\sum_{n\in\N} a_n \sin(2n\varphi).
$$
In particular, their average over a half period vanishes.
Consequently, using the previous symmetry properties   with  the stability of the regularity, which can be checked by straightforward estimates,  we obtain the existence of  
$\varepsilon_0>0,$ such that for any $\varepsilon\in (0,\varepsilon_0)$, the mapping 
$$
(g_+,g_-)\in B_{1,+}^s\times B_{1,-}^s\mapsto \mathcal{S}_+(g_+,g_-)\in \textnormal{Lip}_\nu(\mathcal{O},H^{s,+}_{\textnormal{even},\circ}(\T))
$$
is well-defined and smooth, where 
\begin{align*}
 B_{1,+}^s&=\big\{g_+\in \textnormal{Lip}_\nu(\mathcal{O},H^{s,+}_{\textnormal{even},\circ}(\T)),\quad  \|g_+\|_{s}^{\textnormal{Lip},\nu}\leqslant1  \big\},\\
 B_{1,-}^s&=\big\{g_-\in \textnormal{Lip}_\nu(\mathcal{O},H^{s,-}_{\textnormal{even}}(\T)),\quad  \|g_-\|_{s}^{\textnormal{Lip},\nu}\leqslant1\big\}.
\end{align*}
Next, we turn to equation~\eqref{g--Eq}. 
We begin by simplifying the expression of ${G}_{1,-}[g]$ for an arbitrary function $g\in H^{s,-}_{\textnormal{even}}(\T)$.
Using \eqref{Def-Phij}, Corollary \ref{cor:symmetryy} and Lemma \ref{lem-functions}, we can show that all the $\varrho_j$ are even and 
$$
\varrho_{j,\star}=\varrho_j, \quad \forall\, j=1,2,3.$$
Then we can write
\begin{align*}
    {G}_{1,-}[g](\varphi)&=G_1[g](\varphi)-G_{1,\star}[g](\varphi)=\varrho_2(\varphi)\, \int_{\varphi+\pi}^\varphi \varrho_3(s) g(s)ds\\
    &:=\varrho_2(\varphi)\,F(\varphi).
\end{align*}
Differentiating $F$ and using $g_\star=-g$, we obtain
\begin{align*}
F^\prime(\varphi)&=\varrho_3(\varphi)g(\varphi)-\varrho_{3,\star}(\varphi)g_\star(\varphi)=2\varrho_3(\varphi)g(\varphi).
\end{align*}
Applying Taylor's formula, we infer
\begin{align*}
F(\varphi)&=-\int_0^\pi\varrho_3(\varphi)g(\varphi)d\varphi+2\int_0^\varphi\varrho_3(\tau)g(\tau)d\tau.
\end{align*}
Since $\varrho_3g$ is even and $(\varrho_3g)_\star=-\varrho_3g$, then its Fourier expansion takes the form
$$
\varrho_3(\varphi)g(\varphi)=\sum_{n\in\N}a_n\cos\big((2n+1)\varphi\big).
$$
Then, integrating yields
$$
\int_0^\pi\varrho_3(\varphi)g(\varphi)d\varphi=0.
$$
It follows that
$$
F(\varphi)=2\int_0^\varphi\varrho_3(\tau)g(\tau)d\tau.
$$
Consequently,
\begin{align*}
    {G}_{1,-}[g](\varphi)
    &:=2\varrho_2(\varphi)\,\int_0^\varphi\varrho_3(\tau)g(\tau)d\tau.
\end{align*}
Plugging this identity into \eqref{mathSS} yields the following.
\begin{align}\label{mathSSLL}
 \mathscr{P}[g](\varphi)=\varrho_1(\varphi) \left(b_1+2\int_0^s\varrho_2(s)\,\int_0^\varphi\varrho_3(\tau)g(\tau)d\tau ds\right).
\end{align}
Applying Lemma \ref{lem-fundamental} (we may change there $\varrho_2$ by $2\varrho_2$) we infer that for the choice 
$$
b_1= -\int_0^{\pi}\varrho_2(\tau)\left(\int_0^\tau\varrho_3(s) g(s)ds\right)d\tau,
$$
we get that the operator $\mathscr{P}:H^s_{\textnormal{even},*}(\T)\to H^s_{\textnormal{even},*}(\T)$ is a compact, and under the assumption
\begin{align}\label{Condi-Fund}
1+\int_0^{\pi}\varrho_2(\tau)\left(\int_0^\tau\varrho_3(s) (\textnormal{Id}-\mathscr{T})^{-1}[\varrho_1](s)ds\right)d\tau\neq0,
\end{align}
the operator $\textnormal{Id}-\mathscr{P}$ has a bounded inverse $(\textnormal{Id}-\mathscr{P})^{-1}:H^s_{\textnormal{even},*}(\T)\to H^s_{\textnormal{even},*}(\T)$. Here, the operator $\mathscr{T}:L^\infty([0,2\pi];\R)\to L^\infty([0,2\pi];\R)$ is defined by
$$
 \mathscr{T}[g](\varphi):=2\varrho_1(\varphi)\int_0^\varphi\varrho_2(\tau)\left(\int_0^\tau\varrho_3(s) g(s)ds\right) d\tau.
 $$
 For $\mathcal{S}_-(g_+,g_-)$, one can show in a similar way that it maps $H^s_{\textnormal{even},*}(\T)$ into itself if and only if
 $$
 b_2=-\frac12\int_0^{\pi} e^{\frac{T_0}{\pi}\int_0^s\check{\alpha}(\tau)d\tau}{G}_{2}(\varepsilon,\kappa,g_{+},g_{-})(s)ds.
 $$
 Hence, under this assumption and the condition \ref{Condi-Fund}, the equation \eqref{g--Eq} writes
 \begin{align}\label{g--Eq1}
 g_{-}(\varphi)&= |\ln\varepsilon|^{-\frac12}(\textnormal{Id}-\mathscr{P})^{-1}\mathcal{S}_-(g_+,g_-)(\varphi).
\end{align}
Plugging this equation into \eqref{Eq--02} yields
\begin{align}\label{g--Eq02}
  \nonumber g_{+}(\varphi)&= |\ln\varepsilon|^{-\frac12}\int_0^\varphi {G}_{1,+}[(\textnormal{Id}-\mathscr{P})^{-1}\mathcal{S}_-(g_+,g_-)](s) ds+|\ln\varepsilon|^{-\frac12}\int_0^\varphi {G}_{2,+}(\varepsilon,\kappa,g_{+},g_{-})(s) ds\\
&=:\mathcal{S}_+(g_+,g_-)(\varphi),
\end{align}
Here, to simplify the notation, we use the same symbol $\mathcal{S}_+$ as introduced in \eqref{Eq--02}.
Performing a straightforward analysis, we can show that 
under the assumption \eqref{Condi-Fund}, there exists 
$\varepsilon_0>0,$ such that for any $\varepsilon\in (0,\varepsilon_0)$, the mapping 
$$
(g_+,g_-)\in B_{1,+}^s\times B_{1,-}^s\mapsto (\mathcal{S}_+,\mathcal{S}_-)\in B_{1,+}^s\times B_{1,-}^s
$$
is a contraction. Thus, the Banach fixed-point theorem ensures the existence of a unique solution $(g_+,g_-)\in B_{1,+}^s\times B_{1,-}^s$ to the system \eqref{g--Eq02}–\eqref{g--Eq1}.
 To recover $g_0,$ we write
$$
g_0=\tfrac12(g_{+}+g_{-}),\quad  g_{0,\star}=\tfrac12(g_{+}-g_{-}).
$$
The compatibility condition $(g_{0})_{\star}=g_{0,\star}$ follows from the structure as $g_{+,\star}=g_+$
and $g_{-,\star}=-g_{-}.$
 This completes the construction of $g_0$. Therefore, with
$f_0=g_0\cos(\theta)+h_0$, we have
$$
\Pi_{1}   {\bf F}_1(\varepsilon,f_0)=0.
$$
Finally, combining together \eqref{first-app1} and  \eqref{First-mode1}, we conclude
\begin{align*}
 {\bf F}_1&(\varepsilon,f_0)=\varepsilon^3|\ln\varepsilon| \omega\partial_\varphi h_0+\varepsilon^3|\ln\varepsilon|^{\frac12}\Pi_1^c\mathcal{N}_0(f_0)=O(\varepsilon^3|\ln\varepsilon|),
\end{align*}
which achieves \eqref{Initia-0}.\\
{$\diamond$ \bf  Induction step.}  
Assume that, for some $n\geqslant 1$, we have constructed 
$$f_{n-1}(\varphi,\theta):=\sum_{m= 0}^{n-1} \varepsilon^{m}(g_{m}(\varphi)\cos(\theta)+h_{m}(\varphi,\theta)),\quad\hbox{with}\quad \Pi_1^c f_{n-1}=\sum_{m= 0}^{n-1}{\varepsilon^m}\,h_m,
$$ 
where the coefficients satisfy the symmetry conditions
\begin{align*}
\int_{\mathbb{T}} g_m(\varphi)d\varphi=0,\qquad  g_m(-\varphi)=g_m(\varphi),\quad\hbox{and}\quad h_m(-\varphi,-\theta)=h_m(\varphi,\theta),
 \end{align*}
and the inductive bounds
\begin{equation}\label{ind-asump}
\Pi_1  {\bf F}_1(\varepsilon,f_{n-1})=0,\qquad \Pi_1^c  {\bf F}_1(\varepsilon,f_{n-1})= O(\varepsilon^{2+n}|\ln\varepsilon|).
\end{equation}
This holds for $n=1$ by the construction of $f_0$.
We seek a correction $$w_n:= g_n\cos(\theta)+h_n,$$  and set
$$
f_n:=f_{n-1} +\varepsilon^n w_n.
$$
Using Taylor expansion of ${\bf F}_1$ at $f_{n-1}$ up to the first order, as in the initialization step, yields
\begin{align*}
 {\bf F}_1(\varepsilon,f_{n})&=\varepsilon^n \Pi_{1,{\bf{c}}} \partial_f{\bf F}_1(\varepsilon,f_{n-1})[w_{n}]+\varepsilon^{2+2n}\Pi_{1,{\bf{c}}} \mathcal{N}_n(w_n)+\Pi_1^c{\bf F}_1(\varepsilon,f_{n-1})\\
 &\quad+\varepsilon^{3+n}|\ln\varepsilon| \omega\partial_\varphi h_n+\tfrac12{\varepsilon^{1+n}}\{\mathcal{H}-\partial_\theta\}[h_n]+\varepsilon^{2+n}\Pi_1^c\mathfrak{L}(f_{n-1})[w_n]+{\varepsilon^{2+2n}}\Pi_1^c\mathcal{N}_n(w_n).
\end{align*}
where $\mathcal{N}_n$ is at least quadratic in $w_n$.
We define $h_n$ by solving
$$
\Pi_1^c{\bf F}_1(\varepsilon,f_{n-1})+\tfrac12{\varepsilon^{1+n}}\{\mathcal{H}-\partial_\theta\}[h_n]=-\varepsilon^{2+n} \Pi_1^c\mathfrak{L}(f_{n-1})[w_n].
$$
As in the base step, this equation is uniquely solvable by a contraction argument as for $h_0$. By reversibility of ${\bf F}_1(\varepsilon,f_{n-1})$ the solution $h_n$ is even.
Using the inductive assumption \eqref{ind-asump}, there exists $\varepsilon_{0,n}>0$ such that, for $0<\varepsilon\leqslant \varepsilon_{0,n}$, the solution $h_n$ can be decomposed as
\begin{equation}\label{Hn}
h_n(\varphi)=\varepsilon|\ln\varepsilon| H_{n,1}(\varepsilon)(\varphi)+\varepsilon H_{n,2}(\varepsilon,g_0)(\varphi),
\end{equation}
with
$$
\sup_{\varepsilon\in(0,\varepsilon_{0,n})}\|H_{n,1}(\varepsilon)\|_s^{{\textnormal{Lip},\nu}}<\infty,\quad \sup_{\varepsilon\in(0,\varepsilon_{0,n})}\|H_{n,2}(\varepsilon,g_0)\|_s^{{\textnormal{Lip},\nu}}<\infty.
$$
Substituting $h_n$  into the expansion of ${\bf F}_1(\varepsilon,f_{n})$ yields
\begin{align}\label{Fn-stepn}
\nonumber {\bf F}_1(\varepsilon,f_{n})=&\varepsilon^n \Pi_{1,{\bf{c}}} \partial_f{\bf F}_1(\varepsilon,f_{n-1})[w_{n}]+\varepsilon^{2+2n}\Pi_{1,{\bf{c}}} \mathcal{N}_n(w_n)+\varepsilon^{3+n}|\ln\varepsilon| \omega\partial_\varphi h_n\\
&+{\varepsilon^{2+2n}}\Pi_1^c\mathcal{N}_n(w_n).
\end{align}
We now choose $g_n$ so as to eliminate the entire mode $1$ (only the cosine component is present), namely
\begin{align*}
  \Pi_{1,{\bf{c}}} \partial_f{\bf F}_1(\varepsilon,f_{n-1})[g_{n}\cos(\theta)+h_{n}]
&+\varepsilon^{2+n}\Pi_{1,{\bf{c}}} \mathcal{N}_n(w_n)=0.
\end{align*}
Using Lemma~\ref{prop:dF-computations3}--(1) and the structure $f_{n-1}=f_0+O(\varepsilon)$ 
we get
\begin{align*}
 \Pi_{1,{\bf{c}}}(\partial_f {\bf F}_1)&(\varepsilon,f_{n-1})[g_{n}\cos(\theta)+h_{n}]= \Big\{\varepsilon^3 |\ln\varepsilon|\Big(\omega g_n'+\mathtt{q}_3 g_{n,\star}+ \mathtt{f}_2 g_{n}-\tfrac12 \mathtt{g}_3 g_0 h_{n,3,s}- \mathtt{f}_2 h_{n,3,c}\\
 &-\mathtt{g}_2 h_{n,3,s}\Big)+\tfrac14\varepsilon^2  \mathtt{g}_3 ( h_{n,2,s}+3h_{n,4,s})+\tfrac14\varepsilon^2 g_0 h_{n,2,s} +\varepsilon^2 \mathcal{M}_1(f_0)[g_{n}\cos(\theta)+h_{n}]
\\
 & +\varepsilon^3|\ln\varepsilon| \mathcal{M}_2(f_0)[g_{n}\cos(\theta)+h_{n}]+O(\varepsilon^3|\ln\varepsilon|^{\frac12}h)\Big\}\cos(\theta).
\end{align*}
where the terms involved in $\mathcal{M}_1(f_0)[f_n]+\varepsilon^3|\ln\varepsilon| \mathcal{M}_2(f_0)[f_n]$ have the following Fourier expansion:
\begin{align*}
   &\mathtt{L}_1[f_n]+   \mathtt{Q}_1[f_0,f_n]+\varepsilon|\ln\varepsilon| \mathtt{L}_2[f_n]= -\tfrac14 (2\mathtt{p}_{1,1})^{-\frac14}g_0 h_{n,2,c}+\tfrac{1}{64\mathtt{p}_{1,1}}(h_{n,4,c}-7 h_{n,2,c})\\
   &+(2\mathtt{p}_{1,1})^{-\frac14}\varepsilon|\ln\varepsilon|\Big(\mathtt{g}_2 h_{n,3,c}+\mathtt{g}_2g_n-\mathtt{f}_2 h_{n,3,s}+\tfrac{1}{16}(2\mathtt{p}_{1,1})^{-\frac32}g_n-\mathtt{q}_1g_{n,\star} \Big)+O(\varepsilon|\ln\varepsilon|^{\frac12} ).
\end{align*}
Using \eqref{Hn}, Lemma \ref{lem-lam}, Lemma \ref{asymt-list1} and similarly to the first step,  we find that 
\begin{align*}
  & g_n^\prime(\varphi)+\tfrac{\check{\alpha}T_0}{2\pi}\big(g_n-g_{n,\star}\big)+\tfrac{\check{\mathtt{h}}_2 e^{\check{\mathtt{f}}_3} T_0^2(\lambda)}{64\pi^2\kappa}\int_0^\varphi e^{-\check{\mathtt{f}}_3(s)}(g_n-g_{n,\star})(s)ds \\
  &\quad=F_n(\varepsilon,\varphi)+|\ln\varepsilon|^{-\frac12}G_{2,n}(\varepsilon,\kappa, g_{n},g_{n,\star})(\varphi).
\end{align*}
This equation can be written in the form
\begin{align*}
 \nonumber  g_n^\prime(\varphi)+\tfrac{\check{\alpha}T_0}{2\pi}\big(g_n-g_{n,\star}\big)(\varphi)+G_{1}[g_{n}-g_{n,\star}](\varphi)=\, F_n(\varepsilon,\varphi)&+|\ln\varepsilon|^{-\frac12}G_{2,n}(\varepsilon,\kappa, g_{n},g_{n,\star})(\varphi)
\end{align*}
where $G_1$ is the same functional as \eqref{check-f3} and  
the functional $G_{2,n}$ satisfies for each $\kappa>0, n\in\N$ and $\varepsilon_0>0$ small enough,
$$
\|g_n\|_s^{{\textnormal{Lip},\nu}}\leqslant1\Longrightarrow \sup_{\varepsilon\in(0,\varepsilon_0)}\|G_{2,n}\|_s^{{\textnormal{Lip},\nu}}<\infty.
$$
Proceeding as in the first step in the construction of $g_0$, and under the assumption \eqref{Condi-Fund}, we construct a solution $g_n$ by means of a fixed-point argument so that it satisfies the required properties. Consequently, the new approximation
$$
f_n(\varphi,\theta)=f_{n-1}(\varphi,\theta)+\varepsilon^{n}w_n(\varphi,\theta)=f_{n-1}(\varphi,\theta)+\varepsilon^{n}\big(g_{n}(\varphi)\cos(\theta)+h_n(\varphi,\theta)\big).
$$
satisfies by construction,
$$
 \Pi_{1,{\bf{c}}}{\bf F}_1(\varepsilon,f_{n})=0,
 $$
 and from \eqref{Fn-stepn} and \eqref{Hn}, we deduce
\begin{align*}
{\bf F}_1&(\varepsilon,f_{n})=
\varepsilon^{3+n}|\ln\varepsilon| \omega\partial_\varphi h_n+{\varepsilon^{2+2n}}\Pi_1^c\mathcal{N}_n(w_n)=O(\varepsilon^{n+3}|\ln\varepsilon|).
\end{align*}
This ends the induction step and achieves the proof of the result.
\end{proof}
\subsection{Analysis of the non-resonant condition \texorpdfstring{\eqref{Condi-Fund0}}{condition (Condi-Fund0)}}
The main objective of this section is to analyze the non-resonance condition \eqref{Condi-Fund0}, which plays a crucial role in the construction of the approximate solution in Theorem~\ref{theo-approx1}, as well as in the invertibility of the mode-one operator established in Proposition~\ref{pro-L1-iso}.
 In particular, this condition ensures the solvability of certain auxiliary ODEs governing the mode-one component and prevents the occurrence of degeneracies in the reduced system.
For convenience, we recall here the precise formulation of this assumption. We define
\begin{align*}
F_0
&:=\Big\{(\lambda,\kappa)\in  (0,\infty)^2\ \textnormal{s.t.}\ 
\mathcal{P}(\lambda,\kappa)= 0\Big\},\\
\nonumber
\mathcal{P}(\lambda,\kappa)
&=1+\int_0^{\pi}\varrho_2(\tau)\left(\int_0^\tau\varrho_3(s)\,(\textnormal{Id}-\mathscr{T})^{-1}[\varrho_1](s)\,ds\right)d\tau,
\end{align*}
where the involved functions are defined through \eqref{Def-Phij} and from \eqref{Def-TT}, the operator
$\mathscr{T}$ is defined by 
\begin{align*}
 \mathscr{T}[g](\varphi)=\varrho_1(\varphi)\int_0^\varphi\varrho_2(\tau)\left(\int_0^\tau\varrho_3(s) g(s)ds\right) d\tau.
 \end{align*}
In what follows, we investigate the topological structure of the open set $\mathcal{O}_0$.

\begin{pro}
 The following assertions hold.
   \begin{enumerate}
   \item Given $\kappa=\kappa_0>0$, then the set
$$
F_0\cap \big\{(\lambda,\kappa_0), \lambda>0 \big\},
$$
is discrete and it does not concentrate around $\lambda=0.$
\item Given $\lambda=\lambda_0>0$, then the set
$$
F_0\cap \big\{(\lambda_0,\kappa), \kappa>0 \big\},
$$
is discrete.

   \end{enumerate}
\end{pro}
\begin{proof}
    Recall that
    \begin{align*}
\varrho_1(\varphi)=e^{-\frac{T_0}{\pi}\int_0^\varphi\check{\alpha}(\tau)d\tau},\quad  \varrho_2(\varphi)=\tfrac{ T_0^2\check{\mathtt{h}}_2(\varphi) e^{\check{\mathtt{f}}_3(\varphi)}} {64\pi^2\kappa\, \varrho_1(\varphi)}\quad\hbox{and}\quad \varrho_3(\varphi)=e^{-\check{\mathtt{f}}_3(s)}\,.
\end{align*}
From the real analyticity of the involved functions we get that $(\lambda,\kappa)\mapsto \mathcal{P}(\lambda,\kappa)$ is real analytic.  In addition, according to Lemma \ref{lem-functions} and Proposition \ref{prop:period}, we infer for any $\varphi\in[0,\pi]$
\begin{align*}
    \tfrac{T_0}{\pi}\left|\int_0^\varphi\check{\alpha}(\tau)d\tau\right|&\leqslant {2 \lambda^2}e^{\frac{\lambda^2}{16\kappa}}\int_0^\pi\frac{d\tau}{y_1^2(\tau)+y_2^2(\tau)}\leqslant 2 \pi e^{\frac{\lambda^2}{16\kappa}}, 
\end{align*}
which implies that
$$
\forall \varphi\in[0,\pi],\quad|\varrho_1^{\pm1}(\varphi)|\leqslant e^{2 \pi e^{\frac{\lambda^2}{16\kappa}}}.
$$
As to the functions $\check{\mathtt{h}}_2$ and $\check{\mathtt{f}}_3$, they are  defined by Lemma  \ref{lem-functions} and \eqref{check-f3}. Their estimates follow the same lines as $\varrho_1$ and one gets
$$
\forall \varphi\in[0,\pi],\quad|\varrho_3^{\pm1}(\varphi)|\leqslant e^{2 \pi e^{\frac{\lambda^2}{16\kappa}}}\quad\hbox{and}\quad |\check{\mathtt{h}}_2(\varphi)|\leqslant \lambda^{-2},
$$
and similarly
\begin{align*}
\forall \varphi\in[0,\pi],\quad|\varrho_2(\varphi)|&\leqslant e^{4 \pi e^{\frac{\lambda^2}{16\kappa}}}\lambda^{-2}\tfrac{T_0^2}{64\pi^2\kappa}\leqslant e^{4 \pi e^{\frac{\lambda^2}{16\kappa}}}e^{\frac{\lambda^2}{16\kappa}}\tfrac{\lambda^2}{16\kappa}\leqslant e^{5 \pi e^{\frac{\lambda^2}{16\kappa}}}\tfrac{\lambda^2}{16\kappa}\cdot
\end{align*}
It follows that 
\begin{align*}
\|\mathscr{T}[g]\|_{L^\infty}&\leqslant \pi^2\|g\|_{L^\infty}\prod_{j=1}^3\|\varrho_j\|_{L^\infty}\leqslant \pi^2e^{9 \pi e^{\frac{\lambda^2}{16\kappa}}}\tfrac{\lambda^2}{16\kappa} \|g\|_{L^\infty}.
\end{align*}
Therefore, by taking 
$$
\pi^2e^{9 \pi e^{\frac{\lambda^2}{16\kappa}}}\tfrac{\lambda^2}{16\kappa}\leqslant \tfrac12,
$$
we get
$$
\|(\textnormal{Id}-\mathscr{T})^{-1}[g]\|_{L^\infty}\leqslant 2\|g\|_{L^\infty}.
$$
Hence, it follows that
\begin{align*}
|\mathcal{P}(\lambda,\kappa)-1|&\leqslant 2\pi^2\prod_{j=1}^3\|\varrho_j\|_{L^\infty}\leqslant \pi^2e^{9 \pi e^{\frac{\lambda^2}{16\kappa}}}\tfrac{\lambda^2}{8\kappa}. 
\end{align*}
This shows that, for any fixed $\kappa>0$, the function $\lambda\in(0,\infty)\mapsto \mathcal{P}(\lambda,\kappa)$ is not identically zero; indeed, it suffices to consider the limit $\lambda\to 0$. Consequently, any zeros of $\mathcal{P}(\cdot,\kappa)$, if they exist, must lie in an interval of the form $[\lambda_0,\infty)$ for some $\lambda_0>0$. 
Moreover, by analyticity with respect to $\lambda$, the set of zeros is discrete and therefore finite on any compact interval $[a,b]\subset [0,\infty)$.
\\
Fixing now $\lambda>0$ and viewing $\mathcal{P}$ as a function of $\kappa$, the same arguments apply. In this case, the relevant asymptotic regime corresponding to $\lambda\to 0$ is given by $\kappa\to\infty$.
\end{proof}
\subsection{Rescaling and linearization}

Fix an arbitrary integer $N$ and consider an approximate solution $f_N$ as provided by Theorem~\ref{theo-approx1}. In order to analyze the nonlinear equation in a neighborhood of $f_N$, we introduce a rescaled unknown and  rescale the functional ${\bf F}_1$. More precisely, we define
\begin{equation}\label{def:mathcalF}
\mathcal{F}(\varepsilon,g) := \frac{1}{\varepsilon^{1+\mu}}{\bf F}_1\big(\varepsilon, f_N + \varepsilon^\mu g\big),
\qquad 1 \leqslant \mu < 1+N.
\end{equation}
This rescaling isolates the leading-order behavior of the functional around the approximate solution and shows that, at the level of the PDE, only $f_N$ contributes at the time-degeneracy scale, while the perturbation $g$ acts only at smaller scales. This observation is important, as the degeneracy of the first mode will be handled  by the leading term, whose structure is well understood.
By Theorem~\ref{theo-approx1}, the residual associated with the approximate solution $f_N$ is sufficiently small. In terms of the rescaled functional, this translates into the estimate
\begin{equation*}
  \mathcal{F}(\varepsilon,0) = O\big(\varepsilon^{\,2+N-\mu}|\ln\varepsilon|\big),
\end{equation*}
which shows that the defect of $f_N$ in solving the nonlinear equation is of higher order with respect to the rescaled variables. 
\\
The next step consists in studying the linearization of $\mathcal{F}$ with respect to $g$, which takes the form
\begin{align}\label{Linea-G}
\partial_g\mathcal{F}(\varepsilon,g)[h] 
=\frac{1}{\varepsilon}\,\partial_f{\bf F}_1\big(\varepsilon, f_N + \varepsilon^\mu g\big)[h].
\end{align}
Using Corollary~\ref{cor:asymp-linF1}, we derive an explicit asymptotic expansion of the associated linearized operator together with precise estimates for each of its components. 
The resulting structure is described below and summarized in the following result.
\begin{cor}\label{prop:asymp-lin-2}
Assume \eqref{cond1} and let  $ n\in\mathbb{N}^\star$. There exists $\varepsilon_0>0$ such that for all $\varepsilon\in(0,\varepsilon_0)$, the functional
$$\mathcal{F}(\varepsilon,\cdot):B_{1,\textnormal{even}}(s_0)\cap \textnormal{Lip}_\nu(\mathcal{O},H^{s}_{\circ,\textnormal{even}}(\T^2))\to \textnormal{Lip}_\nu(\mathcal{O},H^{s-1}_{\circ,\textnormal{odd}}(\T^2)),
$$  is well-defined and of class $C^1$ with  
  \begin{align*}
      \partial_g\mathcal{F}(\varepsilon,g)[h]&={\Pi_{1,{\bf{s}}}^c}\bigg(\varepsilon^2 |\ln\varepsilon|\omega(\lambda)\partial_\varphi h+ \tfrac{1}{2} \big(\mathcal{H}[h]-\partial_\theta h\big)+\varepsilon \partial_\theta\big(  \mathcal{V}(\varepsilon,g)h\big)+\varepsilon \mathcal{H}_{\mathtt{u},0}[h]+\varepsilon{\partial_\theta}\mathcal{S}[h]\\ &\quad {+\varepsilon \mathscr{M}_1[h]\cos(\theta)+\varepsilon^2|\ln\varepsilon| \mathscr{M}_2[h]\cos(\theta)} +{\varepsilon^2}|\ln\varepsilon|  {\partial_\theta}\mathcal{Q}[h]\\ &\quad+\varepsilon^2\partial_\theta\fint_{\mathbb{T}} \mathcal{W}({g})(\varphi,\theta,\eta)\ln\big|\sin\big(\tfrac{\theta-\eta}{2}\big)\big|h(\eta) d\eta+\varepsilon^2\partial_\theta\mathcal{R}_{\infty}(g)[h]\bigg),
  \end{align*}
where
\begin{align*}
\mathscr{M}_j[h]:=    \mathcal{M}_j(f_0)[h],
\end{align*}
\begin{align}\label{mathcalV}
   \mathcal{V}(\varepsilon,g):=&V_1(\varepsilon,f_N+\varepsilon^\mu g)\\
   \nonumber=&\tfrac32 \mathtt{g}_{3}(\varphi)\cos(3\theta)+\big(\tfrac12 g_0(\varphi)-2\mathtt{g}_{3}(\varphi)\big)\cos(\theta)\\
   \nonumber&+2\varepsilon|\ln\varepsilon|\big(\mathtt{f}_{2}(\varphi)\sin(2\theta)- \mathtt{g}_{2}\cos(2\theta)\big)+O(\varepsilon|\ln\varepsilon|^\frac12),
\end{align}
with
$$
\|\mathcal{V}(\varepsilon,g)\|^{\textnormal{Lip},\nu}_s\lesssim 1+\|g\|^{\textnormal{Lip},\nu}_s\quad\hbox{and}\quad \|\mathcal{W}(g)\|^{\textnormal{Lip},\nu}_s\lesssim 1+\|g\|^{\textnormal{Lip},\nu}_{s}.
$$
In addition, the operator $\mathcal{R}_{\infty}\in \textnormal{OPS}^{-\infty}$ and satisfies  in particular the estimates
$$
\quad \|\partial_\theta^n\mathcal{R}_{\infty}h\|^{\textnormal{Lip},\nu}_{s}\lesssim (1+\|g\|^{\textnormal{Lip},\nu}_{s_0+n})\|h\|_s +\|g\|^{\textnormal{Lip},\nu}_{s+1+n}\|h\|^{\textnormal{Lip},\nu}_{s_0}.
$$
 Furthermore
\begin{equation*}
\partial^2_g {\mathcal F}(\varepsilon,g)[h_1,h_2]
 \varepsilon^{1+\mu} {\rm E}_2(\varepsilon,g)[h_1,h_2],
\end{equation*}
  where
\begin{align*}
			\| {\rm E}_2(\varepsilon,g)[h_1,h_2]\|_{s}^{{\textnormal{Lip},\nu}}&\lesssim\|h_{1}\|_{s_{0}+2}^{{\textnormal{Lip},\nu}}\|h_{2}\|_{s+2}^{{\textnormal{Lip},\nu}}+\Big(\|h_{1}\|_{s+2}^{{\textnormal{Lip},\nu}}+\|g\|_{s+2}^{{\textnormal{Lip},\nu}}\|h_{1}\|_{s_{0}+2}^{{\textnormal{Lip},\nu}}\Big)\|h_{2}\|_{s_{0}+2}^{{\textnormal{Lip},\nu}}.
            \end{align*}
The  functions $\mathtt{f}_2$, $\mathtt{g}_2$, $\mathtt{g}_3$   are defined in  \eqref{list-functions}.
The operators $\mathcal{H}$, $\mathcal{Q}$, $\mathcal{H}_{\mathtt{u},0}$ and $\mathcal{S}$ are introduced in \eqref{Hilbert1alt},   \eqref{operator-mode1}, \eqref{def-Hu0} and \eqref{shift-operator1}, respectively. 
\end{cor}
\section{Invertibility of the linearized operator}\label{sec:invertibility}
To analyze the invertibility of the  linearized operator described in Corollary \ref{prop:asymp-lin-2}, we begin by decomposing the phase space into two complementary subspaces.  
The first is the \emph{tangential part}, associated with the fundamental cosine mode $1$, while the second is the \emph{normal part}, consisting of all Fourier modes with index $|n|\geqslant 2$.  
This decomposition allows us to separate the delicate low-frequency behavior from the higher modes, and to treat each contribution within its natural framework.  
With this splitting, the problem of invertibility reduces to the study of a matrix-valued operator that couples the tangential and normal components.  
The key idea is to simplify this coupling by means of a triangularization procedure, implemented through a carefully chosen zero order transformation.  
This step reorganizes the structure of the operator in such a way that the tangential and normal directions can be treated almost independently.  

After this reduction, the question of invertibility is further transformed into the invertibility of the linearized operator restricted separately to the tangential part and to the normal part.  
In other words, the global problem decomposes into two more manageable subproblems: the analysis of a finite-dimensional matrix acting on the tangential modes, and the study of the infinite-dimensional operator governing the normal directions.  
This strategy lies at the heart of the reduction scheme and paves the way for establishing the right-invertibility of the full linearized operator.

\subsection{Triangular reduction up to a smoothing small error}\label{sec:triangularred}

The linearized operator described in Corollary~\ref{prop:asymp-lin-2}, exhibits a rich structure combining a transport part, a singular Hilbert transform contribution, and several non-diagonal perturbative terms at different scales on $\varepsilon$ and {the state $g$. } 
Although the operator acts naturally on the functional spaces $\mathbb{X}^s$ and $\mathbb{Y}^{s-1}$, introduced in \eqref{compact-notat1}, its action mixes tangential and normal modes in a nontrivial way.  
In order to make progress on the invertibility problem, it is convenient to reformulate the action of the linearized operator in terms of a block matrix operator $\mathbb{M}$, reflecting the decomposition of the phase space into tangential and normal components.  
The key point is to simplify this matrix structure by means of a suitable conjugation.  
By introducing an explicit lower–triangular transformation, we aim to eliminate the dominant off–diagonal terms and reduce $\mathbb{M}$ to a triangular form.  
In this new representation, the coupling between tangential and normal parts becomes negligible, up to a smoothing remainder of arbitrarily high order. \\
We begin by recalling the structure of the linearized operator given in Corollary~\ref{prop:asymp-lin-2}. We set
$$
\mathscr{L}(\varepsilon,g)[h]:=\partial_g\mathcal{F}(\varepsilon,g)[h],
$$
then
\begin{align}\label{LL-MM}
     \mathscr{L}(\varepsilon,g)[h] &= {\Pi_{1,{\bf{s}}}^c}\Big[\varepsilon^2 |\ln\varepsilon|\omega(\lambda)\partial_\varphi h+ \tfrac{1}{2} \left\{\mathcal{H}[h]-\partial_\theta h\right\}+\varepsilon \partial_\theta\Big\{ h(\theta) \mathcal{V}(\varepsilon,g)\Big\}+\varepsilon \mathcal{H}_{\mathtt{u},0}[h]\\ 
     \nonumber &+\varepsilon {\partial_\theta}\mathcal{S}[h]{+\varepsilon \mathscr{M}_1[h]\cos(\theta)+\varepsilon^2|\ln\varepsilon| \mathscr{M}_2[h]\cos(\theta)} +{\varepsilon^2}|\ln\varepsilon|  {\partial_\theta}\mathcal{Q}[h]+\varepsilon^2{\partial_\theta}\mathcal{R}(\varepsilon,g)\Big],
\end{align}
with
\begin{align}\label{RId-10}
\mathcal{R}(\varepsilon,g)[h]&=\fint_{\mathbb{T}} \mathcal{W}({g})(\varphi,\theta,\eta)\ln\big|\sin\big(\tfrac{\theta-\eta}{2}\big)\big|h(\eta)\, d\eta+\mathcal{R}_{\infty}(g)[h].
\end{align}
This operator satisfies the mapping property
\begin{align*}
\mathscr{L}(\varepsilon,g):\mathbb{X}^s\to \mathbb{Y}^{s-1},
\end{align*}
where the spaces are defined in \eqref{compact-notat1}. 
In particular, thanks to the speed modulation performed in Section~\ref{sec:modeone}, we have
$$
\forall h\in H_\circ^{s}(\T^{2}),\qquad 
\Pi_{1,{\bf s}}\mathscr{L}(\varepsilon,g)[h]=0,
$$
where $\Pi_{1,{\bf s}}$ denotes the projector defined in \eqref{Pi-proj0}.
   Moreover, the operator ${\partial_\theta}\mathcal{R}(\varepsilon,g)$ is of order zero. 
We now decompose the phase space as, using the notation \eqref{compact-notat2}, 
$$
\mathbb{X}^s=\textnormal{Lip}_\nu(\mathcal{O}, H_{\circ,\textnormal{even},{\bf c}}^{s}(\T^{2}))\oplus \mathbb{X}_\perp^s \quad\hbox{with}\quad \mathbb{X}_\perp^s:=\textnormal{Lip}_\nu(\mathcal{O}, H_{\textnormal{even},\perp}^{s}(\T^{2})),
$$
where the tangential space $\textnormal{Lip}_\nu(\mathcal{O}, H_{\circ,\textnormal{even},{\bf c}}^{s}(\T^{2}))$ coincides with the set of functions of the form $f(\varphi)\cos(\theta)$, subject to the conditions
$$
\langle f\rangle_\varphi=0,\quad f(-\varphi)=f(\varphi),\quad \|f\|_{s}^{\textnormal{Lip},\nu}<\infty.
$$
Similarly, we decompose $\mathbb{Y}^s$ as follows
$$
\mathbb{Y}^s=\textnormal{Lip}_\nu(\mathcal{O}, H_{\circ,\textnormal{odd},{\bf c}}^{s}(\T^{2}))\oplus \mathbb{Y}_\perp^s \quad\hbox{with}\quad \mathbb{Y}_\perp^s:=\textnormal{Lip}_\nu(\mathcal{O}, H_{\textnormal{odd},\perp}^{s}(\T^{2})).
$$
With this decomposition, the action of the operator $\mathscr{L}(\varepsilon,g)$ can be equivalently represented by the matrix operator
$$
\mathbb{M}: \mathbb{E}^{s}\to \mathbb{F}^{s-1},
$$
where
\begin{align}\label{OP-MM}
\mathbb{M}&:=\begin{pmatrix}
\Pi\mathscr{L}(\varepsilon,g)\Pi&\Pi\mathscr{L}(\varepsilon,g)\Pi^\perp \\
\Pi^\perp\mathscr{L}(\varepsilon,g)\Pi & \Pi^\perp\mathscr{L}(\varepsilon,g)\Pi^\perp
\end{pmatrix},
\end{align}
and
\begin{align}\label{product-spaces}
\mathbb{E}^s&:= \textnormal{Lip}_\nu(\mathcal{O}, H_{\circ,\textnormal{even},{\bf c}}^{s}(\T^{2}))\times \mathbb{X}_\perp^s,\quad
 \mathbb{F}^s:= \textnormal{Lip}_\nu(\mathcal{O}, H_{\circ,\textnormal{odd},{\bf c}}^{s}(\T^{2}))\times \mathbb{Y}_\perp^s.
\end{align}
Here, $\Pi = \Pi_{1,{\bf c}}$ denotes the projection onto the cosine mode $\mathtt{c}_1$, as defined in~\eqref{Pi-proj0}.  
Since the sine mode does not belong to the phase space, this projection coincides with the projection onto the entire Fourier mode $n=1$.  
Consequently, $\Pi^\perp$ corresponds to the projection onto the subspace spanned by all Fourier modes with index $|n|\geqslant 2$. \\ 
We compute the off–diagonal components:  
\begin{align*}
    \Pi\mathscr{L}(\varepsilon,g)\Pi^\perp&=\varepsilon\Pi \partial_\theta\Big\{  \mathcal{V}(\varepsilon,g)\Pi^\perp h\Big\}+\varepsilon \Pi{\partial_\theta}\mathcal{S}[\Pi^\perp h]+\varepsilon  \Pi\mathcal{H}_{\mathtt{u},0}[\Pi^\perp h]\\  &\quad {+\varepsilon \Pi\mathscr{M}_1[\Pi^\perp h]\cos(\theta)+\varepsilon^2|\ln\varepsilon| \Pi\mathscr{M}_2[\Pi^\perp h]\cos(\theta)}  
    +\varepsilon^2 \Pi{\partial_\theta}\mathcal{R}(\varepsilon,g)[\Pi^\perp h]\\
      &=:\varepsilon \mathcal{L}_{1,\perp}[h],
\end{align*}
and
\begin{eqnarray}\label{Lperp-1}
    \Pi^\perp\mathscr{L}(\varepsilon,g)\Pi \nonumber h&=&\varepsilon \Pi^\perp\partial_\theta\Big\{  \mathcal{V}(\varepsilon,g)\Pi h\Big\}+\varepsilon \Pi^\perp{\partial_\theta}\mathcal{S}[\Pi h]+\varepsilon  \Pi^\perp\mathcal{H}_{\mathtt{u},0}[\Pi h]
    +\varepsilon^2 \Pi^\perp{\partial_\theta}\mathcal{R}(\varepsilon,g)[\Pi h]\\
      &=:&\varepsilon \mathcal{L}_{\perp,1}[h].
\end{eqnarray}
We note  that the operator ${\partial_\theta}\mathcal{Q}$ does not contribute to the off-diagonal operators as it is localized on the Fourier mode one according to the last point of Lemma \ref{prop:SandQ}. Moreover, 
 both $\mathcal{L}_{1,\perp}$ and $\mathcal{L}_{\perp,1}$ are finite-rank operators and therefore belong to $\textnormal{OPS}^{-\infty}$.  
Accordingly, the matrix operator can be written in the  form
\begin{align*}
\mathbb{M}
&=\begin{pmatrix}
\mathcal{L}_{1}&\varepsilon \mathcal{L}_{1,\perp} \\
\varepsilon \mathcal{L}_{\perp,1} & \mathcal{L}_{\perp}
\end{pmatrix},
\end{align*}
where $\mathcal{L}_1:=\Pi\mathscr{L}(\varepsilon,g)\Pi$ and $\mathcal{L}_\perp:=\Pi^\perp \mathscr{L}(\varepsilon,g)\Pi^\perp$.

The central idea  is to simplify the structure of the matrix operator $\mathbb{M}$ by conjugating it with a carefully designed transformation. The goal is to obtain a triangular form in which the tangential and normal components are weakly coupled, leaving only a smoothing remainder of arbitrarily high order.
This triangularization is a key step and prepares the ground for the subsequent analysis of right-invertibility.
Our main result reads as follows.
\begin{pro}\label{pro-matrixope}
    Assume \eqref{cond1} and let $M\in\N$. Consider $g\in \mathbb{X}^s\cap B_{1,\textnormal{even}}(s_0), $ there exists a  matrix-valued operator
\begin{align*}
\Psi=\Psi[g]
&:=\begin{pmatrix}
\textnormal{Id}&0 \\
\varepsilon \psi & \textnormal{Id}
\end{pmatrix},
\end{align*}
with  $ \psi:\textnormal{Lip}_\nu(\mathcal{O},H^s_{\circ,\textnormal{even},{\bf c}})\to \mathbb{X}^s_\perp$ a bounded smoothing operator  such that
\begin{enumerate}
    \item The map $\psi$ is reversibility preserving, smoothing  and satisfies the tame estimates: 
    $$
  \forall n\in\mathbb{N},\quad \|\partial_\theta^n\psi[h]\|^{\textnormal{Lip},\nu}_{s}\lesssim (1+\|g\|_{s_0+n+{M}}^{\textnormal{Lip},\nu})\| h\|_{s}^{\textnormal{Lip},\nu}+\|g\|_{s+n+1+{M}}^{\textnormal{Lip},\nu}\| h\|_{s_0}^{\textnormal{Lip},\nu}.
  $$
  \item $\Psi$ is an isomorphism of $\mathbb{E}^s$ with
  \begin{align*}
\nonumber \Psi^{-1} \mathbb{M}\Psi
&=\begin{pmatrix}
\mathcal{L}_{1}+\varepsilon^2 \mathcal{L}_{1,\perp}\psi& \varepsilon \mathcal{L}_{1,\perp} \\
0 & \mathcal{L}_{\perp}-\varepsilon^2{\psi\mathcal{L}_{1,\perp}}
\end{pmatrix}+\varepsilon^{M+2}\mathbb{P}_M,
\end{align*}
where  $\mathbb{P}_M\in \textnormal{OPS}^{-\infty}$ and satisfies
\begin{align*}
 \forall n\in\mathbb{N},\quad \|\partial_\theta^n\mathbb{P}_M[h]\|^{\textnormal{Lip},\nu}_{s}\lesssim (1+\|g\|_{s_0+n+{M}}^{\textnormal{Lip},\nu})\| h\|_{s}^{\textnormal{Lip},\nu}+\|g\|_{s+n+1+{M}}^{\textnormal{Lip},\nu}\| h\|_{s_0}^{\textnormal{Lip},\nu}.
  \end{align*}
\end{enumerate}
\end{pro}

\begin{proof}
The scheme of the proof proceeds as follows. We begin by introducing a lower--triangular operator $\Psi$ depending on a smoothing map $\psi$ yet to be determined. The inverse of $\Psi$ is explicit, which allows us to compute directly the conjugated operator $\Psi^{-1}\mathbb{M}\Psi$. At this stage, the lower--left block still contains a complicated term $\mathcal{N}\psi$, which mixes tangential and normal directions.
The key step is to choose $\psi$ in such a way that $\mathcal{N}\psi$ becomes negligible, up to a smoothing remainder. To achieve this, we expand $\psi$ as a finite series in powers of $\varepsilon$
and solve recursively a sequence of homological--type equations for the \mbox{coefficients $\psi_m$.}
To start, let us introduce 
\begin{align*}
\Psi
&:=\begin{pmatrix}
\hbox{Id}&0 \\
\varepsilon \psi & \hbox{Id}
\end{pmatrix},
\end{align*}
with  $ \psi:\textnormal{Lip}_\nu(\mathcal{O},H^s_{\circ,\textnormal{even},{\bf c}})\to\mathbb{X}^s_\perp$ a linear continuous map to be determined later. Then  the matrix operator $\Psi:\mathbb{X}^s\to \mathbb{X}^s$ is invertible with
\begin{align*}
\Psi^{-1}
&=\begin{pmatrix}
\hbox{Id}&0 \\
-\varepsilon \psi & \hbox{Id}
\end{pmatrix}.
\end{align*}
Straightforward computations give
\begin{align}\label{conjug-MM}
\Psi^{-1} \mathbb{M}\Psi
&=\begin{pmatrix}
\mathcal{L}_{1}& \varepsilon \mathcal{L}_{1,\perp} \\
\varepsilon \mathcal{N}\psi & \mathcal{L}_{\perp}
\end{pmatrix}+\varepsilon^2\begin{pmatrix}
 \mathcal{L}_{1,\perp}\psi& 0 \\
0 & -\psi\mathcal{L}_{1,\perp}
\end{pmatrix},
\end{align}
with
$$
\mathcal{N}\psi:=\mathcal{L}_{\perp,1}+\mathcal{L}_{\perp}\psi-\psi\mathcal{L}_{1}-\varepsilon^2\psi\mathcal{L}_{1,\perp}\psi.
$$
Applying Lemma \ref{prop:dF-computations3}--(2) yields 
\begin{align*}
     \mathcal{L}_{1}&=\varepsilon^2 |\ln\varepsilon|\omega(\lambda)\partial_\varphi +\varepsilon^2|\ln\varepsilon| \mathtt{R}_{1},
  \end{align*}
where $\mathtt{R}_{1}:\textnormal{Lip}_\nu(\mathcal{O}, H_{\circ,\textnormal{even},{\bf c}}^{s}(\T^{2}))\to\textnormal{Lip}_\nu(\mathcal{O}, H_{\circ,\textnormal{even},{\bf c}}^{s}(\T^{2}))$ is a bounded operator. More details will be given later in Section \ref{section-Inver-Mode1}. In addition, in  view of Corollary \ref{prop:asymp-lin-2} 
\begin{eqnarray}\label{Eq-Lin-OP}
   \mathcal{L}_\perp (\varepsilon,g)[h] 
 \nonumber  &=& \Pi^\perp \Big[ \varepsilon^2 |\ln \varepsilon| \omega(\lambda) \,\partial_\varphi h 
   + \partial_\theta\Big(\big(-\tfrac12 + \varepsilon \mathcal{V}(\varepsilon,g)\big)h\Big)\Big]\Pi^\perp \\
   &&+ \Pi^\perp \Big[\tfrac12 \mathcal{H}+  \varepsilon \partial_\theta \mathcal{S}_{-1} 
   + \varepsilon^2 \partial_\theta \mathcal{R}(g)\Big]\Pi^\perp,
\end{eqnarray}
where we have using the notation
$$
\partial_\theta \mathcal{S}_{-1}:=\partial_\theta \mathcal{S}+\mathcal{H}_{\mathtt{u},0}.
$$
Thus
\begin{align*}
    \mathcal{N}\psi=& \mathcal{L}_{\perp,1}+\tfrac12\Big(\mathcal{H}-\partial_\theta\Big)\psi+\varepsilon^2 |\ln\varepsilon|\omega(\lambda)\partial_\varphi \psi+\varepsilon \Pi^\perp\partial_\theta\Big\{  \mathcal{V}(\varepsilon,g)\psi\Big\}+\varepsilon \Pi^\perp{\partial_\theta}\mathcal{S}_0\psi\\
      &+{\varepsilon^2}|\ln\varepsilon|  \Pi^\perp{\partial_\theta}\mathcal{Q}\psi-\varepsilon^2|\ln\varepsilon| \psi\mathtt{R}_{1}+\varepsilon^2\Pi^\perp\partial_\theta\mathcal{R}\psi-\varepsilon^2\psi\mathcal{L}_{1,\perp}\psi\\
      :=& \mathcal{L}_{\perp,1}+\tfrac12\Big(\mathcal{H}-\partial_\theta\Big)\psi+\varepsilon \mathtt{T}(\psi)-\varepsilon^2\psi\mathcal{L}_{1,\perp}\psi,
  \end{align*}
  where $\psi\mapsto \mathtt{T}(\psi)$ is linear. Note that the operator $\partial_\varphi\psi$ is the commutator $[\partial_\varphi,\psi].$
  Now, we solve approximately this operator equation as follows. We write $$\psi=\sum_{m=0}^M\varepsilon^m \psi_m,
  $$ 
  and we impose the recursive equations
  \begin{align*}
     \tfrac{1}{2}\Big(\partial_\theta- \mathcal{H}\Big)\psi_0=\mathcal{L}_{\perp,1},
  \end{align*}
  and
   \begin{align*}
    m=1,..,M,\quad \tfrac{1}{2}\Big(\partial_\theta- \mathcal{H}\Big)\psi_m=\mathtt{T}(\psi_{m-1})-{\varepsilon}\sum_{j+\ell=m-1}\psi_j\mathcal{L}_{1,\perp}\psi_\ell.
  \end{align*}
  Therefore, we find
  \begin{equation}\label{Npsi}
\mathcal{N}\psi=\varepsilon^{M+1}\mathtt{T}(\psi_{M})-\sum_{j+\ell\geqslant M}\varepsilon^{2+\ell+j}\psi_j\mathcal{L}_{1,\perp}\psi_\ell.
  \end{equation}
  It follows that
  \begin{align}\label{psi0-defi}
      \psi_0=2\left(\partial_\theta- \mathcal{H}\right)^{-1}\Pi^\perp\mathcal{L}_{\perp,1}.
  \end{align}
  As already mentioned, the operator $\mathcal{L}_{\perp,1}$ has finite rank. 
Therefore, by Lemma~\ref{Commutators-hilbert}-(4), it belongs to $\textnormal{OPS}^{-\infty}$. 
Moreover, using the triangle inequality with \eqref{Lperp-1}, we get
  \begin{align*}
    \|\partial_\theta^n\mathcal{L}_{\perp,1}[h]\|_{s}^{\textnormal{Lip},\nu}\leqslant&\| \partial_\theta^{n+1}(\mathcal{V}(\varepsilon,g)\Pi h)\|_{s}^{\textnormal{Lip},\nu}+\| \partial_\theta^{n+1}\mathcal{S}_{-1}[\Pi h]\|_{s}^{\textnormal{Lip},\nu}+\| \partial_\theta^{n+1}\mathcal{R}(\varepsilon,g)[\Pi h]\|_{s}^{\textnormal{Lip},\nu}.
\end{align*}
To estimate these terms we mainly use Lemma \ref{Law-prodX1} and Lemma \ref{orthog-Lem1}. For instance, for the first transport term, we get
\begin{align*}
     \|\partial_\theta^{n+1}(\mathcal{V}(\varepsilon,g)\Pi h)\|_{s}^{\textnormal{Lip},\nu}\lesssim& \|\mathcal{V}(\varepsilon,g)\|_{s+n+1}^{\textnormal{Lip},\nu}\|\Pi h\|_{s_0}^{\textnormal{Lip},\nu}+\|\mathcal{V}(\varepsilon,g)\|_{s}^{\textnormal{Lip},\nu}\|\partial_\theta^{n+1}\Pi h\|_{s_0}^{\textnormal{Lip},\nu}\\
     &+\|\partial_\theta^{n+1}\Pi h\|_{s}^{\textnormal{Lip},\nu}\|\mathcal{V}(\varepsilon,g)\|_{s_0}^{\textnormal{Lip},\nu}+\|\Pi h\|_{s}^{\textnormal{Lip},\nu}\|\mathcal{V}(\varepsilon,g)\|_{s_0+n+1}^{\textnormal{Lip},\nu}.
\end{align*}
Combined with the estimates of $\mathcal{V}$ in  Corollary \ref{prop:asymp-lin-2} and  Lemma \ref{orthog-Lem1} we deduce that
\begin{align*}
     \|\partial_\theta^{n+1}(\mathcal{V}(\varepsilon,g)\Pi h)\|_{s}^{\textnormal{Lip},\nu}&\lesssim (1+\|g\|_{s_0+n+1}^{\textnormal{Lip},\nu})\| h\|_{s}^{\textnormal{Lip},\nu}+(1+\|g\|_{s+n+1}^{\textnormal{Lip},\nu})\| h\|_{s_0}^{\textnormal{Lip},\nu}\\
&\lesssim
     (1+\|g\|_{s_0+n+1}^{\textnormal{Lip},\nu})\| h\|_{s}^{\textnormal{Lip},\nu}+\|g\|_{s+n+1}^{\textnormal{Lip},\nu}\| h\|_{s_0}^{\textnormal{Lip},\nu}.
\end{align*}
Next, we turn to the term $\mathcal{S}_0$ and focus first on  $\mathcal{H}_{\mathtt{u},0}$ which is  defined by \eqref{def-Hu0}. 
Using \eqref{eq:Lambda-m-f}, it can be written in the form
\begin{align*}
\nonumber \mathcal{H}_{\mathtt{u},0}[h]&=-\tfrac{1}{2}   (2\mathtt{p}_{1,1})^{-\frac34} \partial_\theta\big((\cos(\theta)\Lambda_{0}+\Lambda_{0}(\cos(\cdot))[h].
\end{align*}
Notice that $\Lambda_0$ is an operator of order $-1$ and Lemma \ref{Lem-homo-inhomo} applies. Therefore, we get 
\begin{align*}
     \|\partial_\theta^{n}(\mathcal{H}_{\mathtt{u},0}[\Pi h])\|_{s}^{\textnormal{Lip},\nu}&\lesssim  \|\partial_\theta^{n+1}(\cos(\cdot)\Lambda_{0}(\Pi h))\|_{s}^{\textnormal{Lip},\nu}+\|\partial_\theta^{n+1}(\Lambda_0( \cos(\cdot)\Pi h))\|_{s}^{\textnormal{Lip},\nu}\\
     &\lesssim  \|\partial_\theta^{n}(\Pi h))\|_{s}^{\textnormal{Lip},\nu}\lesssim \| h\|_{s}^{\textnormal{Lip},\nu}. 
\end{align*}
The estimates of $\mathcal{S}$ defined by \eqref{shift-operator1} are quite similar. For the term associated with the  remainder $\partial_\theta\mathcal{R}$ defined through \eqref{RId-10}, it has a structure similar to $\mathcal{H}_{\mathtt{u},0}$. Hence, using the estimates of Corollary \ref{prop:asymp-lin-2} we find
\begin{align*}
     \|\partial_\theta^{n+1}(\mathcal{R}[\Pi h])\|_{s}^{\textnormal{Lip},\nu}&\lesssim (1+\|g\|_{s_0+n+1}^{\textnormal{Lip},\nu})\| h\|_{s}^{\textnormal{Lip},\nu}+\|g\|_{s+n+2}^{\textnormal{Lip},\nu}\| h\|_{s_0}^{\textnormal{Lip},\nu}.
\end{align*}
Putting together these estimates yields
\begin{align}\label{Lperp-est}
  \forall n\in\mathbb{N},\quad \|\partial_\theta^n\mathcal{L}_{\perp,1}[h]\|^{\textnormal{Lip},\nu}_{s}\lesssim (1+\|g\|_{s_0+n+1}^{\textnormal{Lip},\nu})\| h\|_{s}^{\textnormal{Lip},\nu}+\|g\|_{s+n+2}^{\textnormal{Lip},\nu}\| h\|_{s_0}^{\textnormal{Lip},\nu}.
\end{align}
Returning to $\psi_0$ introduced in \eqref{psi0-defi} and using the foregoing  estimates, we find
$$
  \forall n\in\mathbb{N},\quad \|\partial_\theta^n\psi_0[h]\|^{\textnormal{Lip},\nu}_{s}\lesssim (1+\|g\|_{s_0+n}^{\textnormal{Lip},\nu})\| h\|_{s}^{\textnormal{Lip},\nu}+\|g\|_{s+n+1}^{\textnormal{Lip},\nu}\| h\|_{s_0}^{\textnormal{Lip},\nu}.
  $$
  We remark that one may use Lemma \ref{comm-pseudo1} and Lemma \ref{Commutators-hilbert} and gets that $\psi_0\in\textnormal{OPS}^{-\infty}$.
  As $g\in \mathbb{X}^s$ and $\mathcal{L}_{\perp,1}$ is reversible, then $\psi_0$ is reversibility preserving.
  Concerning the next term of the asymptotics $\psi_1$, it satisfies 
  \begin{align*}
   \tfrac{1}{2}\Big(\partial_\theta- \mathcal{H}\Big)\psi_1&=\mathtt{T}(\psi_{0})-{\varepsilon}\,\psi_0\mathcal{L}_{1,\perp}\psi_0\\
   &=\Pi^\perp\big(\mathtt{T}(\psi_{0})-{\varepsilon}\,\psi_0\mathcal{L}_{1,\perp}\psi_0\big)\in \textnormal{OPS}^{-\infty}.
  \end{align*}
  Here we have used the fact that the operators $\mathtt{T}$ and $\mathcal{L}_{1,\perp}$ act continuously on $\textnormal{OPS}^{-\infty},$ and  this latter class is stable by composition.
  and therefore
  \begin{align*}
   \psi_1=2\left(\partial_\theta- \mathcal{H}\right)^{-1}\Pi^\perp\big(\mathtt{T}(\psi_{0})-\varepsilon^{2}\psi_0\mathcal{L}_{1,\perp}\psi_0\big)\in \textnormal{OPS}^{-\infty}.
  \end{align*}
  Similar arguments allow us to construct 
   $\psi_m$  recursively and for each $m=0,..,M$ with
   $
   \psi_m\in \textnormal{OPS}^{-\infty}.
   $
   The map $\psi_m$ is reversibility preserving and proceeding as for $\psi_0$ we find the tame estimates,
  $$
  \forall n\in\mathbb{N},\quad \|\partial_\theta^n\psi_m[h]\|^{\textnormal{Lip},\nu}_{s}\lesssim (1+\|g\|_{s_0+n+{m}}^{\textnormal{Lip},\nu})\| h\|_{s}^{\textnormal{Lip},\nu}+\|g\|_{s+n+1+{m}}^{\textnormal{Lip},\nu}\| h\|_{s_0}^{\textnormal{Lip},\nu}.
  $$
  Summing up yields
  \begin{align}\label{psi-final0}
  \forall n\in\mathbb{N},\quad \|\partial_\theta^n\psi[h]\|^{\textnormal{Lip},\nu}_{s}\lesssim (1+\|g\|_{s_0+n+{M}}^{\textnormal{Lip},\nu})\| h\|_{s}^{\textnormal{Lip},\nu}+\|g\|_{s+n+1+{M}}^{\textnormal{Lip},\nu}\| h\|_{s_0}^{\textnormal{Lip},\nu}.
  \end{align}
  Together with \eqref{Npsi} we infer that $\mathcal{N}\psi\in\textnormal{OPS}^{-\infty} $ and 
  $$
  \forall n\in\N,\quad \|\partial_\theta^n\mathcal{N}\psi[h]\|^{\textnormal{Lip},\nu}_{s}\lesssim \varepsilon^{M+1}\Big((1+\|g\|_{s_0+n+{M}}^{\textnormal{Lip},\nu})\| h\|_{s}^{\textnormal{Lip},\nu}+\|g\|_{s+n+1+{M}}^{\textnormal{Lip},\nu}\| h\|_{s_0}^{\textnormal{Lip},\nu}\Big).
  $$
  Finally, it follows from \eqref{conjug-MM} that
  \begin{align}\label{Inver-matrix89}
\nonumber \Psi^{-1} \mathbb{M}\Psi
&=\begin{pmatrix}
\mathcal{L}_{1}+\varepsilon^2 \mathcal{L}_{1,\perp}\psi& \varepsilon \mathcal{L}_{1,\perp} \\
0 & \mathcal{L}_{\perp}-\varepsilon^2\psi\mathcal{L}_{1,\perp}
\end{pmatrix}+\begin{pmatrix}
0& 0 \\
\varepsilon \mathcal{N}\psi & 0
\end{pmatrix}\\
&:=\mathbb{M}_1+\varepsilon^{M+2}\mathbb{P}_M,
\end{align}
with $\mathbb{P}_M\in \textnormal{OPS}^{-\infty}$ and 
\begin{align*}
 \forall n\in\N,\quad \|\partial_\theta^n\mathbb{P}_{N}[h]\|^{\textnormal{Lip},\nu}_{s}\lesssim (1+\|g\|_{s_0+n+{M}}^{\textnormal{Lip},\nu})\| h\|_{s}^{\textnormal{Lip},\nu}+\|g\|_{s+n+1+{M}}^{\textnormal{Lip},\nu}\| h\|_{s_0}^{\textnormal{Lip},\nu}.
  \end{align*}
  We emphasize once again that the required symmetry properties are preserved throughout the reduction scheme. 
As a result of the proof, the matrix-valued operator is transformed into a triangular form, modulo a smoothing remainder of arbitrarily small size. 
This triangularization achieves a weak decoupling between the tangential and normal components of the operator.
  \end{proof}

\subsection{Invertibility of the first-mode operator \texorpdfstring{$\mathcal{L}_1$}{L1}}\label{section-Inver-Mode1}
The aim of this subsection is to isolate and analyze the tangential dynamics carried by the first Fourier mode. We focus on the one–dimensional block extracted by the projector $\Pi=\Pi_{1,{\bf c}}$, which singles out the cosine component of frequency~$1$. On this block, the linearized operator reflects the parity switch induced by the $\partial_\varphi,\partial_\theta$–structure, mapping even functions into odd ones with a loss of one derivative. After performing the normal–tangential splitting and the triangularization introduced earlier, the operator under consideration reduces to
$$
\mathcal{L}_1 = \Pi\,\mathscr{L}(\varepsilon,g)\,\Pi.
$$
This operator admits the decomposition
$$
\mathcal{L}_1 \;=\; \varepsilon^2|\ln\varepsilon| \omega(\lambda)\,\partial_\varphi \;+\; \mathcal{K}_1(\varepsilon,g),
$$
where the first part encodes the principal transport strength in the mode~$1$, and $\mathcal{K}_1$ is a zero–order perturbation that smoothly depends on the state $g$. The smallness of the first term in $\varepsilon$ is, at first glance, an obstruction to invertibility. However, as we shall see below,  this degeneracy is precisely balanced by the structure of $\mathcal{K}_1$, which satisfies 
$$
\mathcal{K}_1(\varepsilon,g) \;=\; {\varepsilon^2|\ln\varepsilon|}\,\mathtt{R}_1,
$$
where $\mathtt{R}_{1}$ is a bounded operator. 
In particular, $\mathcal{K}_1(\varepsilon,g)$ is a compact perturbation of the transport part, vanishing at the same critical rate as $\varepsilon^2|\ln\varepsilon|$. This compensation mechanism ensures that $\mathcal{L}_1$ remains perturbatively invertible. The proof relies on a contraction-mapping argument based on Lemma \ref{lem-fundamental}.
The main  result is the isomorphism estimate stated below, which forms a cornerstone for reconstructing the inverse of the full operator.

\begin{pro}\label{pro-L1-iso}
    Assume \eqref{cond1} and consider $g\in \mathbb{X}^s\cap B_{1,\textnormal{even}}(s_0). $ Then, the operator 
$$\mathcal{L}_1:\textnormal{Lip}_\nu(\mathcal{O}, H_{\circ,\textnormal{even},{\bf c}}^{s}(\T^{2}))\to \textnormal{Lip}_\nu(\mathcal{O}, H_{\circ,\textnormal{odd},{\bf c}}^{s-1}(\T^{2})),
    $$ is an isomorphism with
      $$
      \|\mathcal{L}_1^{-1}[h]\|^{\textnormal{Lip},\nu}_{s}\lesssim \tfrac{1}{\varepsilon^2|\ln \varepsilon|} \Big(\|h\|^{\textnormal{Lip},\nu}_{s-1}+\|g\|_s\|h\|^{\textnormal{Lip},\nu}_{s_0}\Big).
      $$
\end{pro}

\begin{proof}
We begin by computing the action of the operator on the first Fourier mode.  
For a real-valued function $h_1(\varphi)$, we infer from \eqref{Linea-G}
\begin{align*}
\mathcal{L}_1\big[h_1(\varphi)\,\cos(\theta)\big]&=\Pi\mathscr{L}(\varepsilon,g)\big[h_1(\varphi)\,\cos(\theta)\big]\\ &=\tfrac{1}{\varepsilon}\Pi(\partial_f {\bf F}_1)(\varepsilon,f_N+\varepsilon^\mu g)\big[h_1(\varphi)\cos(\theta)\big].
\end{align*}
Then, applying Lemma~\ref{prop:dF-computations3}--(2) one finds
\begin{align*}
\mathcal{L}_1\big[h_1(\varphi)\,\cos(\theta)\big]
   = {\omega(\lambda)} \varepsilon^2 |\ln\varepsilon|  \bigg[& h_1^\prime   +  \tfrac{T_0(\lambda,\kappa)}{2\pi} \check{\alpha} ( h_1  -h_{1,\star})+  G_1[h_1-h_{1,\star}]\\
&\qquad - |\ln\varepsilon|^{-\frac12}  G_2(\varepsilon,\kappa,g)[h_1]\bigg]\cos(\theta),
\end{align*}
where 
$$
G_1[h]=  \tfrac{T_0^2(\lambda)\check{\mathtt{h}}_2 e^{\check{\mathtt{f}}_3(\varphi)}} {64\kappa\pi^2}\, \int_0^\varphi e^{-\check{\mathtt{f}}_3(s)}(h_1-h_{1,\star}) (s)ds 
$$
 and the linear maps $h\in  H_{\circ,\textnormal{even}}^{s}(\T)\mapsto G_2(\varepsilon,\kappa,g)[h]\in \textnormal{Lip}_\nu(\mathcal{O}, H_{\textnormal{odd}}^{s-1}(\T))$ is  continuous and satisfies the tame estimates
 $$
 \|g\|_{s}^{{\textnormal{Lip},\nu}}\leqslant 1\Longrightarrow \sup_{\varepsilon\in(0,\varepsilon_0)}\|G_2(\varepsilon,\kappa,g)[h]\|_{s}^{{\textnormal{Lip},\nu}}\lesssim 
 \|h\|_{s}^{{\textnormal{Lip},\nu}}+\|h\|_{s_0}^{{\textnormal{Lip},\nu}}\|g\|_{s}^{{\textnormal{Lip},\nu}}.
 $$
 Now, let us explore the invertibility of this operator $\mathcal{L}_1$.
Given $g_1=q_1(\varphi)\,\cos(\theta)$ with $q_1\in \textnormal{Lip}_\nu(\mathcal{O}, H_{\textnormal{odd}}^{s-1}(\T))$, we seek to solve the equation
$$
\mathcal{L}_1\big[h_1(\varphi)\,\cos(\theta)\big]=q_1(\varphi)\,\cos(\theta),
$$
with $h_1\in H_{\circ,\textnormal{even}}^{s}(\T)).$
It is equivalent to
\begin{align*}
 h_1^\prime   +  \tfrac{T_0(\lambda,\kappa)}{2\pi} \check{\alpha}\, ( h_1  -h_{1,\star})+ G_1[h_1-h_{1,\star}]=|\ln\varepsilon|^{-\frac12}  G_2(\varepsilon,\kappa,g)[h_1]
 +\varepsilon^{-2} |\ln\varepsilon|^{-1}\tfrac{q_1}{\omega(\lambda)}.
\end{align*}
Notice that this ODE is closely related to \eqref{Eq-Lun1} and can therefore be analyzed in the same way; we omit the details. Under the assumption \eqref{Condi-Fund}, we deduce that, for $\varepsilon$ sufficiently small, the ODE admits a unique solution $h_1$ satisfying for all $s\geqslant s_0$,
$$
\|h_1\|_s^{\textnormal{Lip},\nu}
   \;\lesssim\; \tfrac{1}{\varepsilon^2|\ln\varepsilon|}\Big(
      \|q_1\|_{s-1}^{\textnormal{Lip},\nu}
      + \|q_1\|_{s_0}^{\textnormal{Lip},\nu}\,\|g\|_{s}^{\textnormal{Lip},\nu}\Big).
$$
This concludes the proof of the desired result.
\end{proof}

\subsection{Right-invertibility of the normal operator \texorpdfstring{$\mathcal{L}_\perp$}{L⊥}}\label{sec:inver-normal}

The main goal of this section is to state a key result concerning the right-invertibility of the normal operator defined in \eqref{Eq-Lin-OP},
$$
\mathcal{L}_\perp:\mathbb{X}_\perp^s \;\to\; \mathbb{Y}_\perp^{\,s-1}.
$$
A detailed analysis of this operator will be carried out in Section~\ref{sec:rednormalop}. Here, we focus on the main results and their implications for the invertibility of the full linearized operator.\\
We shall  construct a suitable isomorphism $\Phi:\mathbb{X}_\perp^s \to \mathbb{X}_\perp^s$ that conjugates the normal operator to a Fourier multiplier, up to a smoothing operator of order $-3$, which is the key step to be able to find an approximate right inverse to the normal operator. The main result is stated below.

\begin{pro}\label{prop-diagonal-impo}
Let $(\nu, \tau,s_0,s_{\textnormal{up}})$ as in \eqref{cond1}. Fix $\mu_2\in\R$ such that
\begin{equation*}
\mu_{2}\geqslant 4\tau+3,
\end{equation*}
and let $g\in \mathbb{X}^{s_\textnormal{up}}$.
There exists ${\epsilon}_0>0$ such that if 
\begin{align}\label{small-C2-0}
\|g\|_{\frac{3}{2}\mu_{2}+2s_{0}+2\tau+{4}}^{\textnormal{Lip},\nu}\leqslant1, \quad \varepsilon^3|\ln\varepsilon|^{-\frac12}\nu^{-1}\leqslant 1\quad\textnormal{and}\quad\nu \varepsilon^{{-}2}|\ln\varepsilon|^{-1}+ N_{0}^{\mu_{2}}\varepsilon^{{5}} {\nu^{-1}}\leqslant{\epsilon}_0,
\end{align}
then we can construct a reversibility-preserving isomorphism $\Phi:\mathbb{X}_\perp^s\to\mathbb{X}_\perp^s$ such that the following hold. \begin{enumerate}
\item[\rm (i)]
$($\emph{Tame bounds.}$)$ The map $\Phi$ and its inverse $\Phi^{-1}$ satisfy, for every $s_0\leqslant s\leqslant s_{up}$
 $$
    \|\Phi^{\pm1}[h]\|_{s}^{\textnormal{Lip},\nu}
    \lesssim \|h\|_{s}^{\textnormal{Lip},\nu}
       + \|g\|_{s+3}^{\textnormal{Lip},\nu}\,\|h\|_{s_0}^{\textnormal{Lip},\nu}.
$$

\item[\rm (ii)]
$($\emph{Diagonalization up to a smoothing remainder.}$)$
Let $n\in\N$. For every $\lambda\in\mathcal{O}_{n}^{1}$ (as in Proposition~$\ref{Second-change}$), one has
\begin{align*}
\Phi^{-1}\mathcal{L}_\perp\Phi&=
\varepsilon^2|\ln\varepsilon|\,\omega(\lambda)\,\partial_\varphi
     + \mathtt{c}_1\partial_\theta+\mathtt{c}_2
     \,\mathcal{H}+\varepsilon^2\mathtt{c}_3 \partial_\theta \Lambda_2
     \\ &\quad + \varepsilon^2|\ln\varepsilon|^{\frac12}\Pi^\perp\partial_\theta\mathcal{R}_{4,-3}\Pi^\perp
     + \mathtt{E}_n^4,
\end{align*}
where the following properties hold. 
\begin{enumerate}
\item[\rm (a)]
The operator $\partial_\theta\Lambda_2$ is a Fourier multiplier. More precisely, for all $|j|\geqslant 2$,
$$
\partial_\theta\Lambda_2[e^{ij\theta}]= i\, \mathtt{d}_j e^{ij\theta},\qquad \mathtt{d}_j:=
        -\tfrac{j}{4(j^2-1)|j|}\cdot     
$$
 \item[\rm (b)]
The operator $\partial_\theta\mathcal{R}_{4,-3}\in\mathrm{OPS}^{-3}$ satisfies the tame estimates
$$
  \|\partial_\theta^4\mathcal{R}_{4,-3}[h]\|_{s}^{\textnormal{Lip},\nu}
   \lesssim \|h\|_{s}^{\textnormal{Lip},\nu}
      + \|g\|_{s+2\tau+6}^{\textnormal{Lip},\nu}\|h\|_{s_0}^{\textnormal{Lip},\nu}.
$$
\item[\rm (c)]
The remainder $\mathtt{E}_{n}^4$ satisfies
\begin{equation*}
\|\mathtt{E}_{n}^4[h]\|_{s_0}^{\textnormal{Lip},\nu}\lesssim \varepsilon^{{5}} N_{0}^{\mu_{2}}N_{n+1}^{-\mu_{2}}\|h\|_{s_{0}+2}^{\textnormal{Lip},\nu}+\varepsilon N_n^{s_0-s}\|  g\|_{s+2\tau+6}^{\textnormal{Lip},\nu}\|h\|_{s_0}^{\textnormal{Lip},\nu}.
\end{equation*}
The constants $\mathtt{c}_j=\mathtt{c}_j(\varepsilon,g)$, $j=1,2,3$  satisfy the Lipschitz dependence
\begin{align*}
\|\Delta_{12}\mathtt{c}_1\|^{\textnormal{Lip},\nu}&\lesssim \varepsilon^2|\ln\varepsilon|^{\frac12}
                \| \Delta_{12}g\|_{2{s}_{0}+2\tau+3}^{\textnormal{Lip},\nu},\\
				\|\Delta_{12}\mathtt{c}_2\|^{\textnormal{Lip},\nu}&\lesssim \varepsilon^2
                \| \Delta_{12}g\|_{2{s}_{0}+2\tau+3}^{\textnormal{Lip},\nu},\\ \quad \|\Delta_{12}\mathtt{c}_3\|^{\textnormal{Lip},\nu}&\lesssim \| \Delta_{12}g\|_{2{s}_{0}+2\tau+3}^{\textnormal{Lip},\nu}.
			\end{align*}

\end{enumerate}
    \end{enumerate}
\end{pro}

\begin{proof}
 The proof of Proposition \ref{prop-diagonal-impo} follows easily by setting
$$
\Phi=\Phi_1\circ\Phi_2\circ\Phi_3\circ\Phi_4,$$
where $\Phi_1, \Phi_2, \Phi_3$ and $\Phi_4$ are defined in Sections \ref{sec:Phi1}, \ref{sec:Phi2}, \ref{sec:Phi3} and \ref{sec:Phi4}, respectively. The identity
$$
\partial_\theta\Lambda_2[e^{ij\theta}]= i\, \mathtt{d}_j e^{ij\theta},     
$$
follows from \eqref{FM-n-1} and  \eqref{form:Inm} and the limit
$$
\lim_{z\to 2} \frac{\psi(z-j)}{\Gamma(z-j)}=
        (-1)^{j-1}(j-2)! \quad \forall j\geqslant 2,
	$$
leading to
$$
\Lambda_2[e^{i\,j\theta}]=-I_{2,j} e^{i\,j\theta},
$$
where
$$
I_{2,j}= \tfrac{1}{4(j^2-1)j}, \quad \forall j\geqslant 2. 
$$
The case $j\leqslant -2$ can be obtained using the fact that $I_{2,j}$ is even in $j$. This ends the proof of the proposition.
\end{proof}
From the previous proposition we have
\begin{align}\label{Lperp-def}
\nonumber\Phi^{-1}\mathcal{L}_\perp\Phi&=
\varepsilon^2|\ln\varepsilon|\,\omega(\lambda)\,\partial_\varphi
     + \mathtt{c}_1\partial_\theta+\mathtt{c}_2
     \,\mathcal{H}+\varepsilon^2\mathtt{c}_3 \partial_\theta \Lambda_2
     + \varepsilon^2|\ln\varepsilon|^{\frac12}\Pi^\perp\partial_\theta\mathcal{R}_{4,-3}\Pi^\perp
     + \mathtt{E}_n^4\\
     &=:\mathbb{L}_\perp+ \varepsilon^2|\ln\varepsilon|^{\frac12}\Pi^\perp\partial_\theta\mathcal{R}_{4,-3}\Pi^\perp
     + \mathtt{E}_n^4.
\end{align}
After this reduction scheme, the operator is brought into a nearly diagonal form with an explicitly computable principal part: a transport operator in $\varphi$, a constant-coefficient differential operator in $\theta$, and a constant multiple of the Hilbert transform.  
The remaining terms are  smoothing and sufficiently small so that they can be treated perturbatively.  Then, now the  difficulty is reduced to the invertibility of the leading Fourier multiplier $\mathbb{L}_\perp$  which will be the goal of the next result.
We denote by $\{\mu_{j,2}, \; |j|\geqslant2\}$ the spectrum set of  the space diagonal Fourier multiplier of $\mathbb{L}_\perp,$ given by
\begin{equation}\label{def:muj2}
      \mu_{j,2}(\lambda,g)=\mathtt{c}_1(\lambda,g)j+\mathtt{c}_2(\lambda,g)  \textnormal{sign} j+\varepsilon^2\mathtt{c}_3(\lambda,g)\mathtt{d}_j. 
\end{equation}
The main goal is to invert the operator $\mathbb{L}_\perp.$
The proof of the following result can be carried out along the same lines as \cite[Proposition~5.3]{HHM24}. However, due to minor variations and for the sake of completeness, we provide the proof here.
\begin{pro}\label{prop-perp-orth-1}
Assume \eqref{cond1}. There exists $\varepsilon_0>0$ small enough  such that for any  $\varepsilon\in(0,\varepsilon_0),$
  there exists a family of  linear operators $\big(\mathtt{T}_n\big)_{n\in\mathbb{N}}$  with  the estimates
$$
\forall \,N\geqslant 0,  s\in[s_0,s_\textnormal{up}],\quad \sup_{n\in\mathbb{N}}\|\partial_\theta^N\mathtt{T}_nh\|_{s}^{\textnormal{Lip},\nu}\leqslant C\nu^{-1}\||\partial_\theta|^{2\tau+N}h\|_{s}^{\textnormal{Lip},\nu},
$$
and  by restricting the parameter $\lambda$ on  the Cantor set
\begin{equation*}
	\mathcal{O}_{n}^2(g)=\bigcap_{\substack{(\ell,j)\in\mathbb{Z}^{2} \ 2\leqslant |j|\leqslant N_{n}}}\Big\{\lambda\in \mathcal{O}, \,\,\left|\varepsilon^2|\ln\varepsilon| \omega(\lambda)\ell+\mu_{j,2}(\lambda,g)\right|\geqslant {\nu}{| j|^{-\tau}}\Big\},	
\end{equation*}
we get
$$
\mathbb{L}_{\perp}\mathtt{T}_n=\textnormal{Id}+{\mathbb{E}_{n}},
$$
with
\begin{align*}
\|\mathbb{E}_{n}[h]\|_{s_0}^{\textnormal{Lip},\nu}
 &\leqslant  C\nu^{-1}N_n^{s_0-s} \||\partial_\theta|^{2\tau+1}h\|_{{s}}^{\textnormal{Lip},\nu}.\end{align*}
\end{pro}

\begin{proof}
We shall use the splitting
\begin{align}\label{dekomp1}
\nonumber\mathbb{L}_{\perp}&=\underbrace{\varepsilon^2|\ln\varepsilon|\,\omega(\lambda)\,\partial_\varphi
     +\Pi_{N_n}\big( \mathtt{c}_1\partial_\theta+\mathtt{c}_2
     \,\mathcal{H}+\varepsilon^2\mathtt{c}_3 \partial_\theta \Lambda_2\big)}_{:=\,\mathtt{L}_n}\\ &\quad +\Pi_{N_n}^\perp\big( \mathtt{c}_1\partial_\theta+\mathtt{c}_2
     \,\mathcal{H}+\varepsilon^2\mathtt{c}_3 \partial_\theta \Lambda_2\big),
\end{align}
where the projector $\Pi_{N}$ is defined by
$$
h=\sum_{\substack{\ell\in\mathbb{Z}\\ 2\leqslant |j|}}h_{\ell,j}{\bf e}_{\ell,j},\quad \Pi_{N}h =\sum_{\substack{\ell\in\mathbb{Z}\\ 2\leqslant |j|\leqslant N }}h_{\ell,j}{\bf e}_{\ell,j},\quad {\bf e}_{\ell,j}:=e^{i(\ell\varphi+j\theta)}.
$$
By definition
$$
\mathtt{L}_n{\bf e}_{\ell,j}={\bf e}_{\ell,j}\left\lbrace \begin{array}{rcl}  i\big(\varepsilon^2|\ln\varepsilon|\,\omega(\lambda)\,\ell+\mu_{j,2}(\lambda,g)\big);& \hbox{if}& \quad  2\leqslant |j|\leqslant N_n,\quad \ell\in\mathbb{Z},\\
 i\,\varepsilon^2|\ln\varepsilon|\,\omega(\lambda)\,\ell;& \hbox{if} & \quad |j|> N_n,\quad \ell\in\mathbb{Z}.
\end{array}\right.$$
Define the diagonal  operator  $\mathtt{T}_n$ by 
\begin{eqnarray*}
\mathtt{T}_{n}h(\varphi,\theta)&:=&
- i\sum_{\substack{ 2\leqslant |j|\leqslant N_n \\
\ell\in\mathbb{Z}}}\tfrac{\chi\left((\varepsilon^2|\ln\varepsilon|\,\omega(\lambda)\,\ell+\mu_{j,2}(\lambda))\nu^{-1}| j|^{\tau}\right)}{\varepsilon^2|\ln\varepsilon|\,\omega(\lambda)\,\ell+\mu_{j,2}(\lambda)}h_{\ell,j}\,{\bf e}_{\ell,j}(\varphi,\theta),\\
\end{eqnarray*}
where $\chi\in\mathscr{C}^\infty(\mathbb{R},[0,1])$ is an even positive cut-off function  such that 
\begin{equation*}
	\chi(\xi)=\left\{ \begin{array}{ll}
		0\quad \hbox{if}\quad |\xi|\leqslant\frac13,&\\
		1\quad \hbox{if}\quad |\xi|\geqslant\frac12.
	\end{array}\right.
\end{equation*}
Thus, in  the Cantor set $\mathcal{O}_{n,2}$ one has
\begin{align}\label{LT11}
\mathtt{L}_n\mathtt{T}_{n}=\hbox{Id}-\Pi_{N_n}^\perp.
\end{align}
One can easily  check  from Fourier side  that for any $N\geqslant 0$
$$
\sup_{n\in\mathbb{N}}\||\partial_\theta|^N\mathtt{T}_nh\|_{s}\leqslant C\nu^{-1}\||\partial_\theta|^{\tau+N}h\|_{s}.
$$
On the other hand, denote by
\begin{align*}
g_{\ell,j}(\lambda)&=\tfrac{\chi\left((\varepsilon^2|\ln\varepsilon|\,\omega(\lambda)\,\ell+\mu_{j,2}(\lambda))\nu^{-1}| j|^{\tau}\right)}{\varepsilon^2|\ln\varepsilon|\,\omega(\lambda)\,\ell+\mu_{j,2}(\lambda)}=a_{j}\,\widehat{\chi}\big(a_{j}A_{\ell,j}(\lambda)\big),
\end{align*}
with the definition:
\begin{align*}\widehat{\chi}(x):=\,\tfrac{\chi(x)}{x},\quad A_{\ell,j}(\lambda)&:=\,\varepsilon^2|\ln\varepsilon|\,\omega(\lambda)\,\ell+\mu_{j,2}(\lambda),\quad a_{j}:=\,\lambda^{-1}| j|^{\tau}.\nonumber
				\end{align*}
				Notice that $\widehat{\chi}$ is $C^{\infty}$ with bounded derivatives and $\widehat{\chi}(0)=0.$  Then
				$$
				\|g_{\ell,j}\|_{L^\infty}\lesssim  \nu^{-1}| j|^{\tau}.
				$$
				Taking the Lipschitz norm in $\lambda$,  we get
				$$
				\|A_{\ell,j}\|_{\textnormal{Lip}}\lesssim \langle \ell,j\rangle,
				$$
				and therefore
				$$
				\|g_{\ell,j}\|_{\textnormal{Lip}}\lesssim a_{j}^2\|A_{\ell,j}\|_{\textnormal{Lip}}\lesssim \nu^{-2}| j|^{2\tau}  \langle \ell,j\rangle. 
				$$
				It follows that
				\begin{align*}
				\sup_{\lambda_1\neq\lambda_2\in(a,b)}\frac{\|(\mathtt{T}_nh)(\lambda_1)-(\mathtt{T}_nh)(\lambda_2)\|_{H^{s-1}}}{|\lambda_1-\lambda_2|}&\lesssim \nu^{-1} \sup_{\lambda_1\neq\lambda_2\in (a,b)}\frac{\||\partial_\theta|^\tau h(\lambda_1)-h(\lambda_2)\|_{H^{s-1}}}{|\lambda_1-\lambda_2|}\\
				&\quad+\nu^{-2} \sup_{\lambda\in (a,b)}\||\partial_\theta|^{2\tau}h(\lambda)\|_{H^{s}}.
				\end{align*}
				Consequently
				\begin{align*}
				\sup_{n\in\mathbb{N}}\|\mathtt{T}_nh\|_{s}^{\textnormal{Lip},\nu}&\leqslant C\nu^{-1}\||\partial_\theta|^{2\tau}h\|_{s}^{\textnormal{Lip}, \nu}.
				\end{align*}
				Similarly, we find for any $N\geqslant 0$
				\begin{align}\label{Tn-01}
				\sup_{n\in\mathbb{N}}\||\partial_\theta|^N\mathtt{T}_nh\|_{s}^{\textnormal{Lip}, \nu}&\leqslant C\nu^{-1}\||\partial_\theta|^{2\tau+N}h\|_{s}^{\textnormal{Lip},\nu}.
				\end{align}
				Putting together \eqref{dekomp1} with \eqref{LT11} yields on   the Cantor set $\mathcal{O}_{n,2}$ to the identity
\begin{align*}
\mathbb{L}_{1,\perp}\mathtt{T}_{n}=\hbox{Id}\underbrace{-\Pi_{N_n}^\perp-\Pi_{N_n}^\perp\big( \mathtt{c}_1\partial_\theta+\mathtt{c}_2
     \,\mathcal{H}+\varepsilon^2\mathtt{c}_3 \partial_\theta \Lambda_2\big)\mathtt{T}_{n}}_{:= \mathtt{E}_{n}^{2}}.
\end{align*}
From straightforward estimates, using in particular \eqref{Tn-01}, we find
\begin{align*}
\|\Pi_{N_n}^\perp\big( \mathtt{c}_1\partial_\theta+\mathtt{c}_2
     \,\mathcal{H}+\varepsilon^2\mathtt{c}_3 \partial_\theta \Lambda_2\big)\mathtt{T}_nh\|_{s_0}^{\textnormal{Lip},\nu}&\leqslant CN_n^{s_0-s}\|\partial_\theta\mathtt{T}_nh\|_{s}^{\textnormal{Lip},\nu}\\
&\leqslant CN_n^{s_0-s}\nu^{-1}\||\partial_\theta|^{2\tau+1}h\|_{s}^{\textnormal{Lip},\nu}.
\end{align*}
This ends the proof of the desired result.
\end{proof}

As a consequence, we get our final estimate on the existence of a right inverse operator \mbox{for $\mathcal{L}_{\perp}.$} 
\begin{pro}\label{prop-inverseKL} 
There exists $\epsilon_0>0$ such that under the assumptions \eqref{cond1} and
\begin{align}\label{hao11MM}
&\|g\|_{\frac{3}{2}\mu_{2}+2s_{0}+2\tau+{4}}^{\textnormal{Lip},\nu}\leqslant1, \quad\varepsilon^{2}|\ln\varepsilon|^{\frac12} \nu^{-1}+\nu \varepsilon^{{-}2}|\ln\varepsilon|^{-1}+ N_{0}^{\mu_{2}}\varepsilon^{{5}} {\nu^{-1}}\leqslant{\epsilon}_0,
\end{align}
 there exists a family of linear  operators $\big({\mathcal{T}}_{n,\perp}\big)_{n\in\mathbb{N}}$ satisfying
\begin{equation*}
				\forall \, s\in\,[ s_0, S],\quad\sup_{n\in\mathbb{N}}\|{\mathcal{T}}_{n,\perp}h\|_{s}^{\textnormal{Lip},\nu}\leqslant C\nu^{-1}\Big(\||\partial_\theta|^{2\tau}h\|_{s}^{\textnormal{Lip},\nu}+\|g\|_{s+{2\tau+6}}^{\textnormal{Lip},\nu}\|\partial_\theta|^{2\tau}h\|_{s_0}^{\textnormal{Lip},\nu}\Big),
			\end{equation*}
			and such that in the Cantor set $\mathcal{O}_n^{1}(g)\cap \mathcal{O}_n^{2}(g)$
			we have
			$$
			{\mathcal{L}_\perp}\,{\mathcal{T}}_{n,\perp}=\textnormal{Id}_{\mathbb{X}_\perp^s}+{\mathcal{E}}_{n,\perp},
			$$
			with  the following estimate
			\begin{align*}
				 \|\mathcal{E}_{n,\perp}h\|_{s_0}^{\textnormal{Lip},\nu}
				\nonumber&\leqslant C\nu^{-1}N_n^{s_0-s}\big( \|h\|_{s+2\tau+1}^{\textnormal{Lip},\nu}+\| \rho\|_{s+4\tau{+6}}^{\textnormal{Lip},\nu}\|h\|_{s_{0}+2\tau}^{\textnormal{Lip},\nu}\big)\\ &\quad +C\varepsilon^5\nu^{-1}N_{0}^{\mu_{2}}N_{n+1}^{-\mu_{2}}\|h\|_{s_{0}+2\tau+2}^{\textnormal{Lip},\nu}.
			\end{align*}

\end{pro}
\begin{proof}
In view of Propositions \ref{prop-perp-orth-1} and \ref{prop:Phi4}, and as $\tau\leqslant \frac32$, we get for all   $s\in[s_0,s_\textnormal{up}]$ 
\begin{align}\label{Tn-induc}
\nonumber
\varepsilon^2|\ln\varepsilon|\|\mathtt{T}_n \Pi^\perp\partial_\theta\mathcal{R}_{4,-3}\Pi^\perp h\|_{s_0}^{\textnormal{Lip},\nu} &\leqslant  
C\nu^{-1}\varepsilon^2|\ln\varepsilon|\||\partial_\theta|^{2\tau}\Pi^\perp\partial_\theta\mathcal{R}_{4,-3}\Pi^\perp h \|_{s_0}^{\textnormal{Lip},\nu}
\\ 
&\leqslant C\nu^{-1}\varepsilon^2|\ln\varepsilon|\Big( \|h\|_{s_0}^{\textnormal{Lip},\nu}
      + \|g\|_{s_0+2\tau+6}^{\textnormal{Lip},\nu}\|h\|_{s}^{\textnormal{Lip},\nu}\Big).
\end{align}
It follows, under the smallness condition \eqref{hao11MM}, that
\begin{align*}
\varepsilon^2|\ln\varepsilon|\|\mathtt{T}_n \Pi^\perp\partial_\theta\mathcal{R}_{4,-3}\Pi^\perp h\|_{s_0}^{\textnormal{Lip},\nu} &\leqslant  
\tfrac12\|h\|_{s_0}^{\textnormal{Lip},\nu}.
\end{align*}
Thus, the operator $\varepsilon^2|\ln\varepsilon|\mathtt{T}_n \Pi^\perp\partial_\theta\mathcal{R}_{4,-3}\Pi^\perp$ is a contraction. Hence, using  Neumann series, we conclude that the operator 
$$\mbox{Id}_{\mathbb{X}_\perp^{s_0}}+\varepsilon^2|\ln\varepsilon|\mathtt{T}_n \Pi^\perp\partial_\theta\mathcal{R}_{4,-3}\Pi^\perp,
$$
 is invertible  with
$$
\big\|\left(\mbox{Id}+\varepsilon^2|\ln\varepsilon|\mathtt{T}_n \Pi^\perp\partial_\theta\mathcal{R}_{4,-3}\Pi^\perp\right)^{-1}h\big\|_{s_0}^{\textnormal{Lip},\nu}\leqslant 2 \|h\|_{s_0}^{\textnormal{Lip},\nu}.
$$
As to the invertibility for $s\in[s_0,s_\textnormal{up}]$, one can check by induction from \eqref{Tn-induc}, $\forall m\geqslant 1,$
\begin{align*}
\|(\varepsilon^2|\ln\varepsilon|&\mathtt{T}_n \Pi^\perp\partial_\theta\mathcal{R}_{4,-3}\Pi^\perp)^mh\|_{s}^{\textnormal{Lip},\nu} \leqslant \Big(C\varepsilon^2|\ln\varepsilon|\nu^{-1}\Big)^m \|h\|_{{s}}^{\textnormal{Lip},\nu}\\
&+C\varepsilon^2|\ln\varepsilon|\nu^{-1}(m+1) 2^{m-2}\Big(C\varepsilon^2|\ln\varepsilon|\nu^{-1}\Big)^{m-1}\|g\|_{s+2\tau+6}^{\textnormal{Lip},\nu} \|h\|_{{s_0}}^{\textnormal{Lip},\nu}.
\end{align*}
Consequently, we get under the smallness condition \eqref{hao11MM} 
\begin{align*}
\sum_{m\geqslant 0}\|(\varepsilon^2|\ln\varepsilon|\mathtt{T}_n \Pi^\perp\partial_\theta\mathcal{R}_{4,-3}\Pi^\perp)^mh \|_{s}^{\textnormal{Lip},\nu} &\leqslant C \|h\|_{{s}}^{\textnormal{Lip},\nu}+C\|g\|_{s+2\tau+6}^{\textnormal{Lip},\nu}  \|h\|_{{s_0}}^{\textnormal{Lip},\nu}.
\end{align*}
It follows that
\begin{align}\label{jeud-2}
\|(\hbox{Id}_{\mathbb{X}_\perp^s}+\varepsilon^2|\ln\varepsilon|\mathtt{T}_n \Pi^\perp\partial_\theta\mathcal{R}_{4,-3}\Pi^\perp)^{-1}h\|_{s}^{\textnormal{Lip},\nu}&\leqslant C \|h\|_{{s}}^{\textnormal{Lip},\nu}+C\|g\|_{s+2\tau+6}^{\textnormal{Lip},\nu} \|h\|_{{s_0}}^{\textnormal{Lip},\nu}.
\end{align}
Set
\begin{equation}\label{Tn-def}
   {\mathcal{T}}_{n,\perp}:=\left(\mbox{Id}_{\mathbb{X}_\perp^s}+\varepsilon^2|\ln\varepsilon|\mathtt{T}_n \Pi^\perp\partial_\theta\mathcal{R}_{4,-3}\Pi^\perp\right)^{-1} \mathtt{T}_n.
\end{equation}
From  \eqref{jeud-2}, Proposition \ref{prop-perp-orth-1} and \eqref{hao11MM}, we deduce that
\begin{align}\label{jeud-3}
\nonumber \|{\mathcal{T}}_{n,\perp}h\|_{s}^{\textnormal{Lip},\nu} &\leqslant C \|\mathtt{T}_nh\|_{{s}}^{\textnormal{Lip},\nu}+C\|g\|_{s+2\tau+6}^{\textnormal{Lip},\nu}  \|\mathtt{T}_nh\|_{{s_0}}^{\textnormal{Lip},\nu}\\
&\leqslant C\nu^{-1}\Big(\|h\|_{s,2\tau}^\lambda+\|g\|_{s+2\tau+6}^{\textnormal{Lip},\nu}\|h\|_{s_0,2\tau}^{\textnormal{Lip},\nu}\Big).
\end{align}
By virtue of \eqref{Lperp-def} and \eqref{Tn-def}, we find that on the Cantor set $ \mathcal{O}_n^2$ 
\begin{align*}
\nonumber{\mathcal{L}_\perp}\,{\mathcal{T}}_{n,\perp}&=\left(\mathbb{L}_\perp+ \varepsilon^2|\ln\varepsilon|^{\frac12}\Pi^\perp\partial_\theta\mathcal{R}_{4,-2}\Pi^\perp
     + \mathtt{E}_n^4\right)\left(\mbox{Id}_{\mathbb{X}_\perp^s}+\varepsilon^2|\ln\varepsilon|\mathtt{T}_n \Pi^\perp\partial_\theta\mathcal{R}_{4,-2}\Pi^\perp\right)^{-1} \mathtt{T}_n   \\
\nonumber
&=\mbox{Id}_{\mathbb{X}_\perp^s}+\left(\mathtt{E}_n^4-\varepsilon^2|\ln\varepsilon|^{\frac12}\mathbb{E}_n\Pi^\perp\partial_\theta\mathcal{R}_{4,-3}\Pi^\perp
\right){\mathcal{T}}_{n,\perp}
\\
&=:\mbox{Id}_{\mathbb{X}_\perp^s}+{\mathcal{E}}_{n,\perp}.
\end{align*}
Combining  \eqref{jeud-3}, Proposition \ref{prop-perp-orth-1} and Proposition \ref{prop:Phi4}-(2) gives the estimate on ${\mathcal{E}}_{n,\perp}$. This ends the proof of Proposition \ref{prop-inverseKL}.
\end{proof}

\subsection{Right-invertibility of the full operator}\label{sec:rightinvertibility}

Having reduced the operator to triangular form, up to a smoothing remainder of arbitrarily small size, 
the next step is to investigate the right-invertibility of the full operator. 
This can be achieved by perturbative arguments, which show that the triangularization preserves the essential invertibility properties. 
In fact, the triangular reduction yields a weak decoupling between the tangential and normal components, 
so that the invertibility of the full matrix operator can be deduced once suitable approximate inverses for the diagonal blocks are available. 
We emphasize, however, that proving the invertibility of the diagonal block constitutes the most delicate part of the analysis.\\
Our first result in this section states as follows.
  \begin{lem}\label{lem-pertub-inv}
Under  \eqref{cond1}, assume that \eqref{hao11MM} occurs. Let $\mathcal{L}_1^{-1}$ be the inverse of $\mathcal{L}_1$ constructed in Proposition $\ref{pro-L1-iso}$ and $\mathcal{T}_{n,\perp}$ the approximate right inverse of $\mathcal{L}_\perp$ constructed in Proposition $\ref{prop-inverseKL}$. Then  for $\varepsilon$ small enough,
  the following holds.
      \begin{enumerate}
          \item The operator 
     $\mathcal{L}_1+\varepsilon^2 \mathcal{L}_{1,\perp}\psi: \,\textnormal{Lip}_\nu(\mathcal{O}, H_{\circ,\textnormal{even},{\bf c}}^{s}(\T^{2}))\to \textnormal{Lip}_\nu(\mathcal{O}, H_{\circ,\textnormal{odd},{\bf c}}^{s-1}(\T^{2}))$ is an isomorphism with 
      $$
      \|(\mathcal{L}_1+\varepsilon^2 \mathcal{L}_{1,\perp}\psi)^{-1}[h]\|^{\textnormal{Lip},\nu}_{s}\lesssim \tfrac{1}{\varepsilon^2|\ln \varepsilon|} \Big(\|h\|^{\textnormal{Lip},\nu}_{s-1}+\|g\|^{\textnormal{Lip},\nu}_s\|h\|^{\textnormal{Lip},\nu}_{s_0}\Big).$$
     \item The operator $\mathcal{L}_\perp-\varepsilon^2\psi\mathcal{L}_{1,\perp}\psi \,:\, \mathbb{X}^s_\perp\to\mathbb{Y}^{s-1}_\perp$ admits an approximate right inverse denoted by $\mathcal{T}_{n,\perp}^{\varepsilon}$. More precisely, there exists two operators $\mathcal{T}_{n,\perp}^{\varepsilon}$ and $\mathcal{E}_{n,\perp}^{\varepsilon}$ satisfying the estimates
$$
\sup_{n\in\N}\|\mathcal{T}_{n,\perp}^{\varepsilon}[h]\|^{\textnormal{Lip},\nu}_{s}\leqslant C\nu^{-1}\Big(\||\partial_\theta|^{2\tau}h\|_{s}^{\textnormal{Lip},\nu}+\|g\|_{s+{2\tau+6}}^{\textnormal{Lip},\nu}\|h\|_{s_0+2\tau}^{\textnormal{Lip},\nu}\Big),$$ 
      and
    \begin{align*}
				 \|\mathcal{E}_{n,\perp}^\varepsilon [h]\|_{s_0}^{\textnormal{Lip},\nu}
				\nonumber&\leqslant C\nu^{-1}N_n^{s_0-s}\big( \|h\|_{s+2\tau+1}^{\textnormal{Lip},\nu}+\| g\|_{s+4\tau{+6}}^{\textnormal{Lip},\nu}\|h\|_{s_{0}+2\tau}^{\textnormal{Lip},\nu}\big)\\ &\quad +C\varepsilon^5\nu^{-1}N_{0}^{\mu_{2}}N_{n+1}^{-\mu_{2}}\|h\|_{s_{0}+2\tau+2}^{\textnormal{Lip},\nu},
			\end{align*}
            such that for any $\lambda\in\mathcal{O}_n^1\cap \mathcal{O}_n^1 $ we have
            $$(\mathcal{L}_\perp-\varepsilon^2\psi\mathcal{L}_{1,\perp})\mathcal{T}_{n,\perp}^{\varepsilon}=\textnormal{Id}+\mathcal{E}_{n,\perp}^{\varepsilon}.
      $$
       \end{enumerate}
  \end{lem}

\begin{proof}

\medskip\noindent {\bf (1)}
The proof of this first assertion is based on a perturbative argument. To establish the invertibility of 
$$
\mathcal{L}^\varepsilon_1:=\mathcal{L}_1+\varepsilon^2 \mathcal{L}_{1,\perp}\psi,
$$
we factorize the operator as
$$
\mathcal{L}^\varepsilon_1
   = \mathcal{L}_1\Big(\mathrm{Id}+\varepsilon^2 \mathcal{L}_1^{-1}\mathcal{L}_{1,\perp}\psi\Big).
$$
The operator 
$$
\mathcal{L}_{1,\perp}\psi:\textnormal{Lip}_\nu(\mathcal{O}, H_{\circ,\textnormal{even},{\bf c}}^{s}(\T^{2}))
   \;\to\; \textnormal{Lip}_\nu(\mathcal{O}, H_{\circ,\textnormal{odd},{\bf c}}^{s-1}(\T^{2})),
$$
is well-defined and bounded. Indeed, we already observed that $\mathcal{L}_{1,\perp}\in\textnormal{OPS}^{-\infty}$, being of finite rank. Moreover, combining \eqref{Lperp-est} and \eqref{psi-final0}, and using the smallness condition \eqref{hao11MM} (implying $\|g\|_{s_0+1}^{\textnormal{Lip},\nu}\lesssim 1$), we obtain
\begin{align}\label{l1-perp-1}
\| \mathcal{L}_{1,\perp}\psi[h]\|_{s}^{\textnormal{Lip},\nu}
&\lesssim  \| h\|_{s}^{\textnormal{Lip},\nu}+\|g\|_{s+1}^{\textnormal{Lip},\nu}\| h\|_{s_0}^{\textnormal{Lip},\nu}.
\end{align}
On the other hand, Proposition~\ref{pro-L1-iso} yields the tame estimate
$$
\| \mathcal{L}_{1}^{-1}[h]\|_{s}^{\textnormal{Lip},\nu}
\;\lesssim\; \frac{1}{\varepsilon^2|\ln\varepsilon|}
\Big(\|h\|_{s}^{\textnormal{Lip},\nu}+\|g\|_s^{\textnormal{Lip},\nu}\|h\|_{s_0}^{\textnormal{Lip},\nu}\Big).
$$
By composing these two bounds and invoking the smallness condition, we infer
$$
\varepsilon^2\| \mathcal{L}_{1}^{-1}\mathcal{L}_{1,\perp}\psi[h]\|_{s}^{\textnormal{Lip},\nu}
\;\lesssim\; \frac{1}{|\ln\varepsilon|}
\Big(\|h\|_{s}^{\textnormal{Lip},\nu}
   + \|g\|_{s+1}^{\textnormal{Lip},\nu}\|h\|_{s_0}^{\textnormal{Lip},\nu}\Big).
$$
Hence, for $\varepsilon$ sufficiently small, the Neumann series shows that  
$$
\mathrm{Id}+\varepsilon^2\mathcal{L}_1^{-1}\mathcal{L}_{1,\perp}\psi,
$$
is invertible. Moreover, for $s\in[s_0,s_{\textnormal{up}}]$,  
$$
\|(\mathrm{Id}+\varepsilon^2\mathcal{L}_{1}^{-1}\mathcal{L}_{1,\perp}\psi)^{-1}[h]\|_{s}^{\textnormal{Lip},\nu}
\;\lesssim\; \|h\|_{s}^{\textnormal{Lip},\nu}
   + \|g\|_{s+1}^{\textnormal{Lip},\nu}\|h\|_{s_0}^{\textnormal{Lip},\nu}.
$$
By composition, this implies that $\mathcal{L}_1+\varepsilon^2\mathcal{L}_{1,\perp}\psi$ is invertible, with
$$
\|(\mathcal{L}^\varepsilon_1)^{-1}[h]\|_{s}^{\textnormal{Lip},\nu}
\;\lesssim\; \|h\|_{s}^{\textnormal{Lip},\nu}
   + \|g\|_{s+1}^{\textnormal{Lip},\nu}\|h\|_{s_0}^{\textnormal{Lip},\nu}.
$$
More details of this argument are given in \cite[Proposition~5.4]{HHM24}.

\medskip\noindent {\bf (2)} For the second assertion, we show that an approximate right inverse of the perturbed operator $$\mathcal{L}_\perp^\varepsilon:=\mathcal{L}_\perp-\varepsilon^2\psi\mathcal{L}_{1,\perp},
$$ is given by
\begin{align}\label{TN-InV}
\mathcal{T}_{n,\perp}^\varepsilon:=\big(\textnormal{Id}-\varepsilon^2\mathcal{T}_{n,\perp}\psi\mathcal{L}_{1,\perp} \big)^{-1}\mathcal{T}_{n,\perp}.
\end{align}
The inverse of $\textnormal{Id}-\varepsilon^2\mathcal{T}_{n,\perp}\psi\mathcal{L}_{1,\perp}$ is well-defined. Indeed, by Proposition $\ref{prop-inverseKL}$, we have
\begin{align*}
    \varepsilon^2\|\mathcal{T}_{n,\perp}\psi\mathcal{L}_{1,\perp}[h]\|_{s_0}^{\textnormal{Lip},\nu}
    &\lesssim \varepsilon^2\nu^{-1}
    \||\partial_\theta|^{2\tau}\psi\mathcal{L}_{1,\perp}[h]\|_{s_0}^{\textnormal{Lip},\nu}.
\end{align*}
Since $\psi\mathcal{L}_{1,\perp}\in\textnormal{OPS}^{-\infty}$ and satisfies the same estimate \eqref{psi-final0}, we deduce
\begin{align*}
    \varepsilon^2\|\mathcal{T}_{n,\perp}\psi\mathcal{L}_{1,\perp}[h]\|_{s_0}^{\textnormal{Lip},\nu}
    &\lesssim \varepsilon^2\nu^{-1}\|h\|_{s_0}^{\textnormal{Lip},\nu}.
\end{align*}
From \eqref{hao11MM}, we know that $\varepsilon^2\nu^{-1}$ is sufficiently small, so that $\varepsilon^2 \mathcal{T}_{n,\perp}\psi\mathcal{L}_{1,\perp}$ is a contraction. Therefore $\textnormal{Id}-\varepsilon^2\mathcal{T}_{n,\perp}\psi\mathcal{L}_{1,\perp}$ is invertible.  \\
It remains to check that $\mathcal{T}_{n,\perp}^\varepsilon$ is an approximate right inverse of $\mathcal{L}_\perp^\varepsilon.$ Using
$$
\textnormal{Id}={\mathcal{L}_\perp}\,{\mathcal{T}}_{n,\perp}-{\mathcal{E}}_{n,\perp},
$$
we obtain
\begin{align*}
   \mathcal{L}_\perp-\varepsilon^2\psi\mathcal{L}_{1,\perp}
   &=\mathcal{L}_\perp-\varepsilon^2{\mathcal{L}_\perp}\,{\mathcal{T}}_{n,\perp}\psi\mathcal{L}_{1,\perp}+\varepsilon^2{\mathcal{E}}_{n,\perp}\psi\mathcal{L}_{1,\perp}\\
   &=\mathcal{L}_\perp\Big(\textnormal{Id}-\varepsilon^2\,{\mathcal{T}}_{n,\perp}\psi\mathcal{L}_{1,\perp}\Big)+\varepsilon^2{\mathcal{E}}_{n,\perp}\psi\mathcal{L}_{1,\perp}.
\end{align*}
Combining together the latter identity with \eqref{TN-InV} we infer
\begin{align*}
   \mathcal{L}_\perp^\varepsilon\mathcal{T}_{n,\perp}^\varepsilon
   &=\mathcal{L}_\perp\mathcal{T}_{n,\perp}+\varepsilon^2{\mathcal{E}}_{n,\perp}\psi\mathcal{L}_{1,\perp}\big(\textnormal{Id}-\varepsilon^2\mathcal{T}_{n,\perp}\psi\mathcal{L}_{1,\perp}\big)^{-1}\mathcal{T}_{n,\perp}\\
   &=\textnormal{Id}+{\mathcal{E}}_{n,\perp}+\varepsilon^2{\mathcal{E}}_{n,\perp}\psi\mathcal{L}_{1,\perp}\big(\textnormal{Id}-\varepsilon^2\mathcal{T}_{n,\perp}\psi\mathcal{L}_{1,\perp}\big)^{-1}\mathcal{T}_{n,\perp}\\
   &=:\textnormal{Id}+{\mathcal{E}}_{n,\perp}^\varepsilon.
\end{align*}
We remark that the operator $\varepsilon^2\psi\mathcal{L}_{1,\perp}\big(\textnormal{Id}-\varepsilon^2\mathcal{T}_{n,\perp}\psi\mathcal{L}_{1,\perp} \big)^{-1}\mathcal{T}_{n,\perp}$ is bounded. Therefore, ${\mathcal{E}}_{n,\perp}^\varepsilon$ satisfies the same estimates as ${\mathcal{E}}_{n,\perp}$ stated in Proposition $\ref{prop-inverseKL}$, namely
\begin{align}\label{En-epsilon}
				 \nonumber\|\mathcal{E}_{n,\perp}^\varepsilon [h]\|_{s_0}^{\textnormal{Lip},\nu}
				&\leqslant C\nu^{-1}N_n^{s_0-s}\big( \|h\|_{s+2\tau+1}^{\textnormal{Lip},\nu}+\| \rho\|_{s+4\tau{+6}}^{\textnormal{Lip},\nu}\|h\|_{s_{0}+2\tau}^{\textnormal{Lip},\nu}\big)\\ 
                &\quad +C\varepsilon^5\nu^{-1}N_{0}^{\mu_{2}}N_{n+1}^{-\mu_{2}}\|h\|_{s_{0}+2\tau+2}^{\textnormal{Lip},\nu}.
\end{align}
This completes the proof. 
\end{proof}
The next result establishes the existence of an approximate  right inverse for both the matrix operator $\mathbb{M}$ and the scalar operator $\mathscr{L}(\varepsilon,g)$, defined by \eqref{OP-MM} and \eqref{LL-MM}, respectively.
\begin{lem}\label{lem-matrix-opF}
 Assume \eqref{cond1} and  \eqref{hao11MM} and recall \eqref{Inver-matrix89}. Then the following holds.
      \begin{enumerate}
          
      \item There exists a sequence of two matrix operators $(\mathbb{S}_n)_{n\in\N}$ and $(\mathbb{E}_n)_{n\in\N}$ satisfying 
      $$
\sup_{n\in\mathbb{N}}\|\mathbb{S}_n[h]\|_s^{\textnormal{Lip},\nu}\lesssim \frac{1}{\varepsilon|\ln\varepsilon|\nu}\Big(\||\partial_\theta|^{2\tau}h\|_{s}^{\textnormal{Lip},\nu}+\|g\|_{s+{2\tau+6}}^{\textnormal{Lip},\nu}\|h\|_{s_0+2\tau}^{\textnormal{Lip},\nu}\Big),
$$
and
\begin{align*}
				 \|\mathbb{E}_{n} [h]\|_{s_0}^{\textnormal{Lip},\nu}
				&\leqslant C\nu^{-1}N_n^{s_0-s}\big( \|h\|_{s+2\tau+1}^{\textnormal{Lip},\nu}+\| g\|_{s+4\tau{+6}}^{\textnormal{Lip},\nu}\|h\|_{s_{0}+2\tau}^{\textnormal{Lip},\nu}\big)\\ 
                &\quad +C\varepsilon^5\nu^{-1}N_{0}^{\mu_{2}}N_{n+1}^{-\mu_{2}}\|h\|_{s_{0}+2\tau+2}^{\textnormal{Lip},\nu},
\end{align*}
 such that for any $\lambda\in\mathcal{O}_n^1\cap\mathcal{O}_n^2$ we have 
     $$
      \mathbb{M}\mathbb{S}_n=\mathbb{I}+\mathbb{E}_n.
      $$
        \item The operator $\mathscr{L}(\varepsilon,g):\mathbb{X}^s\to \mathbb{Y}^{s-1}$ admits an approximate  right inverse denoted by $\mathcal{T}_n$. More precisely, there exist two linear operators $\mathcal{T}_n$ and $\mathcal{E}_n$ satisfying the estimates
       $$
\sup_{n\in\mathbb{N}}\|\mathcal{T}_n[h]\|_s^{\textnormal{Lip},\nu}\lesssim \frac{1}{\varepsilon|\ln\varepsilon|\nu}\Big(\||\partial_\theta|^{2\tau}h\|_{s}^{\textnormal{Lip},\nu}+\|g\|_{s+{2\tau+6}}^{\textnormal{Lip},\nu}\|h\|_{s_0+2\tau}^{\textnormal{Lip},\nu}\Big),
$$
       and 
       \begin{align*}
				 \|\mathcal{E}_{n} [h]\|_{s_0}^{\textnormal{Lip},\nu}
				&\leqslant C\nu^{-1}N_n^{s_0-s}\big( \|h\|_{s+2\tau+1}^{\textnormal{Lip},\nu}+\| g\|_{s+4\tau{+6}}^{\textnormal{Lip},\nu}\|h\|_{s_{0}+2\tau}^{\textnormal{Lip},\nu}\big)\\ 
                &\quad +C\varepsilon^5\nu^{-1}N_{0}^{\mu_{2}}N_{n+1}^{-\mu_{2}}\|h\|_{s_{0}+2\tau+2}^{\textnormal{Lip},\nu},
\end{align*}
such that for any $\lambda\in\mathcal{O}_n^1\cap\mathcal{O}_n^2$ we have 
       $$
\mathscr{L}\mathcal{T}_n=\textnormal{Id}+\mathcal{E}_n.
       $$
        \end{enumerate}
\end{lem}
 
\begin{proof}
{\bf (1)} The proof of the first point will be carried out in two steps. First, we establish the result for the triangular matrix $\mathbb{M}_1$, and then extend it to the full matrix $\mathbb{M}$ by means of perturbative arguments. 
We recall from \eqref{Inver-matrix89} that
\begin{align*}
 \mathbb{M}_1
&=\begin{pmatrix}
\mathcal{L}^\varepsilon_{1}& \varepsilon \mathcal{L}_{1,\perp} \\
0 & \mathcal{L}^\varepsilon_{\perp}
\end{pmatrix},
\end{align*}
and we plan to check that the matrix 
\begin{align*}
 \mathbb{N}_n
&=\begin{pmatrix}
(\mathcal{L}^\varepsilon_{1})^{-1}& -\varepsilon (\mathcal{L}^\varepsilon_{1})^{-1}\mathcal{L}_{1,\perp}\mathcal{T}^\varepsilon_{\perp} \\
0 & \mathcal{T}^\varepsilon_{\perp}
\end{pmatrix},
\end{align*}
is an approximate right inverse of $\mathbb{M}_1$. Indeed, a direct computation yields
\begin{align}\label{Inver-matrix890}
 \mathbb{M}_1\mathbb{N}_n
&=\mathbb{I}+\begin{pmatrix}
0& 0 \\
0 & {\mathcal{E}}_{n,\perp}^\varepsilon
\end{pmatrix}:=\mathbb{I}+\mathbb{E}_{n,1}.
\end{align}
The estimate of $\mathbb{E}_n$ follows immediately from \eqref{En-epsilon}. However, the estimate of $\mathbb{N}_n$ is a consequence  of Lemma \ref{lem-pertub-inv} and \eqref{l1-perp-1}, which allow to get
\begin{align}\label{est-mathbbN}
\|\mathbb{N}_n[h]\|_s^{\textnormal{Lip},\nu}\lesssim \frac{1}{\varepsilon|\ln\varepsilon|\nu}\Big(\||\partial_\theta|^{2\tau}h\|_{s}^{\textnormal{Lip},\nu}+\|g\|_{s+{4\tau+5}}^{\textnormal{Lip},\nu}\|h\|_{s_0+2\tau}^{\textnormal{Lip},\nu}\Big).
\end{align}
For the second step, we first recall that
\begin{align*}
\nonumber \Psi^{-1} \mathbb{M}\Psi
=\mathbb{M}_1+\varepsilon^{N+2}\mathbb{P}_N,
\end{align*}
with 
\begin{align*}
 \forall n\in\N,\quad \|\partial_\theta^n\mathbb{P}_{N}[h]\|^{\textnormal{Lip},\nu}_{s}\lesssim (1+\|g\|_{s_0+n}^{\textnormal{Lip},\nu})\| h\|_{s}^{\textnormal{Lip},\nu}+\|g\|_{s+n+1}^{\textnormal{Lip},\nu}\| h\|_{s_0}^{\textnormal{Lip},\nu}.
  \end{align*}
We denote 
  \begin{align*}
\nonumber  \mathbb{S}_{n,1}
=\big(\mathbb{I}+\varepsilon^{N+2}\mathbb{N}_n\mathbb{P}_N\big)^{-1}\mathbb{N}_n. 
\end{align*}
This operator is well-defined since $\mathbb{I}+\varepsilon^{N+2}\mathbb{N}_n\mathbb{P}_N$ is invertible. Indeed, for $N\geqslant1$ we have
\begin{align*}
    \varepsilon^{N+2}\|\mathbb{N}_n\mathbb{P}_N[h]\|_{s_0+2\tau}^{\textnormal{Lip},\nu}&\lesssim\tfrac{\varepsilon^{N+1}}{|\ln\varepsilon|\nu}\Big(\||\partial_\theta|^{2\tau}\mathbb{P}_N[h]\|_{s_0+2\tau}^{\textnormal{Lip},\nu}+\|g\|_{s_0+{6\tau+5}}^{\textnormal{Lip},\nu}\|\mathbb{P}_N[h]\|_{s_0+2\tau}^{\textnormal{Lip},\nu}\Big)\\
    &\lesssim\tfrac{\varepsilon^{2}}{|\ln\varepsilon|\nu}\big(1+\|g\|_{s_0+{6\tau+5}}^{\textnormal{Lip},\nu}\big)\|h\|_{s_0+2\tau}^{\textnormal{Lip},\nu}.
\end{align*}
By the smallness condition \eqref{hao11MM}, we obtain
\begin{align*}
    \varepsilon^{N+2}\|\mathbb{N}_n\mathbb{P}_N[h]\|_{s_0+2\tau}^{\textnormal{Lip},\nu}
    &\leqslant \tfrac12\|h\|_{s_0+2\tau}^{\textnormal{Lip},\nu}.
\end{align*}
This shows that $\varepsilon^{N+2}\mathbb{N}_n\mathbb{P}_N$ is a contraction. Therefore, using Neumann series, we deduce that $\mathbb{I}+\varepsilon^{N+2}\mathbb{N}_n\mathbb{P}_N$ is invertible with
$$
\big\|\big(\mathbb{I}+\varepsilon^{N+2}\mathbb{N}_n\mathbb{P}_N\big)^{-1}[h]\big\|_{s_0+2\tau}^{\textnormal{Lip},\nu}\leqslant 2 \|h\|_{s_0+2\tau}^{\textnormal{Lip},\nu}.
$$
Using once again Neumann series and the product law detailed in Lemma \ref{Law-prodX1}, we also find
$$
\big\|\big(\mathbb{I}+\varepsilon^{N+2}\mathbb{N}_n\mathbb{P}_N\big)^{-1}[h]\big\|_{s}^{\textnormal{Lip},\nu}\lesssim \|h\|_{s}^{\textnormal{Lip},\nu}+\|g\|_{s+{4\tau+5}}^{\textnormal{Lip},\nu}\|h\|_{s_0+2\tau}^{\textnormal{Lip},\nu}.
$$
Combining this estimate with \eqref{est-mathbbN}, and fixing $N=1$, we obtain
$$
\big\|\mathbb{S}_{n,1}[h]\big\|_{s}^{\textnormal{Lip},\nu}\lesssim \frac{1}{\varepsilon|\ln\varepsilon|\nu}\Big(\||\partial_\theta|^{2\tau}h\|_{s}^{\textnormal{Lip},\nu}+\|g\|_{s+{4\tau+5}}^{\textnormal{Lip},\nu}\|h\|_{s_0+4\tau}^{\textnormal{Lip},\nu}\Big).
$$
In view of \eqref{Inver-matrix890} we may write
\begin{align*}
\mathbb{M}_1+\varepsilon^{N+2}\mathbb{P}_N&=\mathbb{M}_1+\varepsilon^{N+2} \mathbb{M}_1\mathbb{N}_n\mathbb{P}_N-\varepsilon^{N+2} \mathbb{E}_{n,1}\mathbb{P}_N=\mathbb{M}_1\Big(\mathbb{I}+\varepsilon^{N+2} \mathbb{N}_n\mathbb{P}_N\Big)-\varepsilon^{N+2} \mathbb{E}_{n,1}\mathbb{P}_N.
\end{align*}
Therefore,
\begin{align*}
\Psi^{-1}\mathbb{M}\Psi\mathbb{S}_{n,1}&=\mathbb{M}_1\mathbb{N}_n-\varepsilon^{N+2} \mathbb{E}_{n,1}\mathbb{P}_N\mathbb{S}_{n,1}\\
&=\mathbb{I}+\mathbb{E}_{n,1}-\varepsilon^{N+2} \mathbb{E}_{n,1}\mathbb{P}_N\mathbb{S}_{n,1}\\
&=:\mathbb{I}+\mathbb{E}_{n,2}.
\end{align*}
The estimate of $\mathbb{E}_{n,2}$ can be performed in a straightforward way, leading to an estimate analogous to that of $\mathbb{E}_{n,1}$. By setting
$$
\mathbb{S}_n:=\Psi\mathbb{S}_{n,1}\Psi^{-1},
$$
we deduce
\begin{align*}
\mathbb{M}\mathbb{S}_{n}&=\mathbb{I}+\Psi\mathbb{E}_{n,2}\Psi^{-1}:=\mathbb{I}+\mathbb{E}_{n}.
\end{align*}
Finally, combining the estimates of $\Psi^{\pm1}$ given in Proposition \ref{pro-matrixope} with the bounds on $\mathbb{E}_{n,2}$, we obtain the desired control of $\mathbb{E}_n$ as stated in Lemma \ref{lem-matrix-opF}. The estimate for $\mathbb{S}_n$ follows similarly from that of $\mathbb{S}_{n,1}$. 

\medskip\noindent {\bf (2)} To solve the equation
$$
\mathscr{L}(\varepsilon,g)h = g,
$$
it is equivalent to consider the matrix formulation
$$
\mathbb{M}\begin{pmatrix}
h_1 \\
h_\perp
\end{pmatrix}
=
\begin{pmatrix}
g_1 \\
g_\perp
\end{pmatrix},
$$
where $\begin{pmatrix} h_1 \\ h_\perp \end{pmatrix} \in \mathbb{E}^s$ and 
$\begin{pmatrix} g_1 \\ g_\perp \end{pmatrix} \in \mathbb{F}^s$, the spaces being defined in \eqref{product-spaces}. Since $\mathbb{S}_n$ is an approximate right inverse of $\mathbb{M}$, it follows that an approximate right inverse of the scalar operator $\mathscr{L}(\varepsilon,g)$ can be obtained directly by taking suitable linear combinations of the components of $\mathbb{S}_n$. From this observation we deduce all the required results.
This achieves the proof.
\end{proof}

\section{Reducibility of the normal operator}\label{sec:rednormalop}

The proof of Proposition \ref{prop-diagonal-impo} is based on a sequence of successive conjugation steps for the operator $\mathcal{L}_\perp$ that will be developed in the next sections. Before proceeding, we provide a sketch of the proof, which relies on the construction of a sequence of auxiliary transformations $\Phi_1$, $\Phi_2$, $\Phi_3$, and $\Phi_4$, from which the final transformation will be defined as
$$
\Phi:=\Phi_1\circ\Phi_2\circ\Phi_3\circ\Phi_4.$$
The starting point is the operator $\mathcal{L}_\perp$ which, in view of Corollary \ref{prop:asymp-lin-2} and \eqref{Eq-Lin-OP}, reads as

\begin{equation}\label{Lperp}
   \mathcal{L}_\perp (\varepsilon,g) 
   = \Pi^\perp \Big[ \mathcal{T} + \tfrac12 \mathcal{H}+\varepsilon \mathcal{H}_{\mathtt{u},0} + \varepsilon \partial_\theta \mathcal{S} 
   + \varepsilon^2 \partial_\theta \mathcal{R}(g)\Big]\Pi^\perp,
\end{equation}
where $\mathcal{T}$ stands for the transport part of the operator:
\begin{align}\label{def:T}
   \mathcal{T}(g)[h] 
   := \varepsilon^2 |\ln \varepsilon| \omega(\lambda) \,\partial_\varphi h 
   + \partial_\theta\!\left(\Big(-\tfrac12 + \varepsilon \mathcal{V}(\varepsilon,g)\Big)h\right).
\end{align}

\medskip\noindent $\blacktriangle$ {\bf {\sc Step 1: }}{\bf Construction of $\Phi_1.$}
\medskip

The operator $\mathcal{L}_\perp$ consists of a dominant diagonal component together with lower-order perturbative terms involving variable coefficients. Our first step is to diagonalize its transport part. Specifically, we straighten the transport component through a suitable change of coordinates $\mathscr{B}$, constructed in Proposition~\ref{Second-change},  so that it becomes a Fourier multiplier up to small remainders. This transformation is carefully designed to eliminate resonances within a Cantor set of parameters and is complicated by the presence of time degeneracy. The procedure relies on KAM techniques in the spirit of \cite{FGMP19, HHM21, HHM24, HR21}.  
Subsequently, we conjugate the nonlocal terms $\mathcal{H}$, $\partial_\theta \mathcal{S}$, and $\partial_\theta \mathcal{R}$ with the transformation $\mathscr{B}$. We then introduce  
$$
\Phi_1 := \Pi^\perp \mathscr{B} \Pi^\perp,
$$  
to fully conjugate the operator $\mathcal{L}_\perp$, obtaining a new operator whose transport part has constant coefficient $\mathtt{c}_1$:  
$$
\begin{aligned}
\mathcal{L}_{\perp,1} = \Phi_1^{-1}\mathcal{L}_\perp\Phi_1 = & \;\varepsilon^2|\ln\varepsilon|\omega(\lambda)\,\partial_\varphi + \mathtt{c}_1(\lambda,g)\,\partial_\theta + \tfrac12\mathcal{H} + \varepsilon^2|\ln\varepsilon|\,\partial_\theta \mathcal{S}_1 \Pi^\perp \\
& + \varepsilon \Pi^\perp \mathcal{H}_{\mathtt{u}_1,0}\Pi^\perp + \varepsilon^2 \Pi^\perp \mathcal{H}_{\mathtt{u}_2,-2}\Pi^\perp + \varepsilon^2 |\ln\varepsilon|^{\frac12}\,\partial_\theta \mathcal{R}_{1,-4} + \mathtt{E}_n^1,
\end{aligned}
$$
where $\partial_\theta \mathcal{R}_{1,-4}$ is an operator of order $-4$. All operators appearing above are defined in Section~\ref{sec:Phi1}.  
Since this step uses KAM theory to reduce the transport part, the conjugation can only be performed on the Cantor set $\mathcal{O}_n^1$ defined in Proposition~\ref{Second-change} by,  
$$
\mathcal{O}_n^1(g) = \bigcap_{\substack{(\ell,j)\in\mathbb{Z}^2 \\ 1\leqslant |j|\leqslant N_n}} \Big\{\lambda\in\mathcal{O} : \big|\varepsilon^2|\ln\varepsilon|\omega(\lambda)\ell + j\mathtt{c}_1(\lambda,g)\big| \geqslant \nu |j|^{-\tau}\Big\}, \qquad n\in\mathbb{N}.
$$
The remainder $\mathtt{E}_n^1$, arising from the KAM procedure, depends on $n$ but it is smoothing, small in $\varepsilon$ and decaying in $n$:
$$
\|\mathtt{E}_n^1 h\|_{s_0}^{\mathrm{Lip},\nu} \lesssim \varepsilon^5 N_0^{\mu_2} N_{n+1}^{-\mu_2} \|h\|_{s_0+2}^{\mathrm{Lip},\nu}.
$$

\medskip\noindent $\blacktriangle$ {\bf {\sc Step 2: }}{\bf Construction of $\Phi_2.$}
\medskip

In Section \ref{sec:Phi2} we define $\Phi_2$ to reduce the zero order operator and the $-2$ order operator of $\mathcal{L}_{\perp,1}$:
$$
\varepsilon\Pi^\perp \mathcal{H}_{\mathtt{u}_1,0}\Pi^\perp+\varepsilon^2\Pi^\perp \mathcal{H}_{\mathtt{u}_2,-2}\Pi^\perp.
$$
Then, applying $\Phi_2$ we find in Section \ref{sec:Phi2}:
\begin{align*}
\nonumber \mathcal{L}_{\perp,2} 
   & = \Phi_2^{-1}\mathcal{L}_{\perp,1}\Phi_2\\
   &= \varepsilon^2|\ln\varepsilon|\lambda \partial_\varphi
      + \mathtt{c}_1\partial_\theta+\tfrac12\mathcal{H}
      + \varepsilon^2\mathtt{w}_1 \,\mathcal{H}+\varepsilon^2 \mathtt{w}_2\, \partial_\theta\Lambda_{2}+ \varepsilon \Pi^\perp\partial_\theta\mathcal{S}_{2}\Pi^\perp   \nonumber
       \\
   &\quad+ \varepsilon^2|\ln\varepsilon|\,\Pi^\perp \partial_\theta\mathcal{S}_1 \Pi^\perp
      + \varepsilon^2|\ln\varepsilon|^{\frac12}\Pi^\perp\partial_\theta\mathcal{R}_{2,-3}\Pi^\perp
      + \mathtt{E}_n^2,
\end{align*}
for some $\mathtt{w}_1$ and $\mathtt{w}_1$ even functions in $\varphi$ and constant in $\theta$, and where the operator $\partial_\theta\mathcal{R}_{2,-3}$ is of order $-3$ and $\mathtt{E}_n^2$ is smooth and small in $\varepsilon$. The main difference from $\mathcal{L}_{\perp,1}$ is that now the operators in the main part are smoothing, i.e., $\mathcal{S}_2$ belongs to $\textnormal{OPS}^{-\infty}$.

\medskip\noindent $\blacktriangle$ {\bf {\sc Step 3: }}{\bf Construction of $\Phi_3.$}
\medskip

After Step 2, we still have antidiagonal terms in the principal part of the operator. Hence, in Section \ref{sec:Phi3} these terms will be removed by solving a suitable homological equation in Fourier variables, thereby ensuring that the coupling between different Fourier modes disappeared up to a smoothing remainder. Then, we obtain:
\begin{align*}
 \mathcal{L}_{\perp,3}
  \;=\; \Phi_3^{-1}\,\mathcal{L}_{\perp,2}\,\Phi_3
 \nonumber =& \varepsilon^2|\ln\varepsilon|\,\omega(\lambda)\,\partial_\varphi
     \,+\, \mathtt{c}_1\partial_\theta+\tfrac12\mathcal{H}
     \;+\; \varepsilon^2 \mathtt{w}_1\,\mathcal{H}
     \;+\; \varepsilon^2 \mathtt{w}_2\,\partial_\theta\,\Lambda_2\\
    \nonumber &+ \varepsilon^2|\ln\varepsilon|^{\frac12}\Pi^\perp\partial_\theta\mathcal{R}_{3,-3}\Pi^\perp
     \;+\; \mathtt{E}_n^3,
\end{align*}
where
$\partial_\theta\mathcal{R}_{3,-3}$ is of order $-3$ and   $\mathtt{E}_n^3$ is smoothing and small in $\varepsilon$. 
We emphasize that the coefficients $\mathtt{W}_j$ depend on the time variable, and the associated operators must therefore be reduced. This is precisely the purpose of the next construction.

\medskip\noindent $\blacktriangle$ {\bf {\sc Step 4: }}{\bf Construction of $\Phi_4.$}
\medskip

At this stage, the operator is transformed into a nearly diagonal form whose principal part can be explicitly identified: a transport operator in $\varphi$, a constant-coefficient operator in $\theta$, a multiple of the Hilbert transform and a multiple of the $\partial_\theta \Lambda_2$ operator. The final step in Section \ref{sec:Phi4} is about the residual $\varphi$–dependence of the coefficient in front of the Hilbert transform and the $\partial \Lambda_2$ operator: by conjugating with a propagator of the form $e^{\rho(\varphi)\mathcal{H}}$, the variable coefficients can be replaced by its $\varphi$–average.  Namely, we get
\begin{align}\label{L4-1}
 \nonumber \mathcal{L}_{\perp,4}
  \;=\; \Phi_4^{-1}\,\mathcal{L}_{\perp,3}\,\Phi_4
 \nonumber =& \varepsilon^2|\ln\varepsilon|\,\omega(\lambda)\,\partial_\varphi
     \,+\, \mathtt{c}_1\partial_\theta+\big(\tfrac12
     + \varepsilon^2 \langle \mathtt{w}_1\rangle_{\varphi} \big)\,\mathcal{H}
     \;+\; \varepsilon^2 \langle \mathtt{w}_2\rangle_{\varphi}\,\partial_\theta\,\Lambda_2\\
    \nonumber &+ \varepsilon^2|\ln\varepsilon|^{\frac12}\Pi^\perp\partial_\theta\mathcal{R}_{4,-3}\Pi^\perp
     \;+\; \mathtt{E}_n^4,
\end{align}
where $\partial_\theta\mathcal{R}_{4,-3}$ is of order $-3$ and $\mathtt{E}_n^4$ is smooth and small in $\varepsilon$. Note that the constants $\mathtt{c}_2$ and $\mathtt{c}_3$ appearing in Proposition \ref{prop-diagonal-impo} coincide with,
$$
\mathtt{c}_2=\tfrac12+\varepsilon^2 \langle \mathtt{w}_1\rangle_{\varphi},\quad \mathtt{c}_3=\langle \mathtt{w}_2\rangle_{\varphi}.
$$
The next sections are devoted to the detailed construction of the transformations $(\Phi_j)_{j=1,\dots,4}$ introduced above.
\subsection{Straightening of the transport part: \texorpdfstring{$\Phi_1$}{Phi1}}\label{sec:Phi1}

The central goal of this section is to introduce a change of coordinates $\Phi_1$ that diagonalizes the transport part of the operator $\mathcal{L}_\perp$. The construction is organized into three steps, detailed in Sections~\ref{sec:transportB}--\ref{sec:phi1-2}. 

First, we conjugate the transport operator $\mathcal{T}$, introduced in \eqref{def:T}, into a Fourier multiplier, up to a controlled error, by means of the change of coordinates $\mathscr{B}$. This reduction is feasible only when the external parameter lies in a suitably chosen Cantor set, which guarantees the elimination of resonances through KAM techniques. This step is carried out in Section~\ref{sec:transportB}. 

Next, in Section~\ref{section-CFLO}, we conjugate the non-local operators appearing in the linearized operator using the same transformation $\mathscr{B}$. Finally, in Section~\ref{sec:phi1-2}, we address localization effects in the normal directions and define 
$$
\Phi_1 := \Pi^\perp \mathscr{B} \Pi^\perp,
$$
which simplifies the operator $\mathcal{L}_\perp$. The construction of $\Phi_1$ is delicate but essential, as it prepares the ground for subsequent diagonalization of the full operator.

\subsubsection{Reduction of the transport part}\label{sec:transportB}

It is important to emphasize that implementing the full KAM scheme at this stage, as in \cite{FGMP19,HHM21,HR21}, is  not feasible due to multiple degeneracies affecting both the linearized operator and the Cantor sets. A key obstruction lies in the fact that the perturbation in the transport operator $\mathcal{T}$ is of size $\varepsilon$, which is larger than the time degeneracy scale $\varepsilon^2|\ln \varepsilon|$. Such a mismatch fundamentally breaks down the KAM scheme.
To overcome this difficulty, we follow the strategy developed in \cite[Sec.~4]{HHM24}. We begin with an auxiliary step based on a suitable change of variables acting solely on the spatial component. This transformation reduces the advection vector field to a constant one, up to an error term sufficiently small in $\varepsilon$.
 Interestingly, the time degeneracy plays a favorable role in this step. With this new structure in place, a KAM-type reduction can then be applied.
We intend to establish the following result.
 \begin{pro}\label{First-change0}
  Assume \eqref{cond1}, there exists ${\varepsilon}_0>0$ 
such that if $g\in \mathbb{X}^{s_\textnormal{up}}$ with
		$$\textnormal{and}\quad\|g\|_{s_0+2}^{\textnormal{Lip},\nu}\leqslant 1,$$ 
then for any $\varepsilon\in (0,\varepsilon_0),$  there exists  $\beta_1:\T^2\to\R$ a smooth odd function in the form
$$
\beta_1(\varphi,\theta):={b}_0(\varphi)+\varepsilon {b}_1(\varphi,\theta)+\varepsilon^2|\ln\varepsilon|^{\frac12}{b}_{1,e}(\varphi,\theta),$$
 and 
\begin{equation*}
\begin{aligned}
\mathscr{B}_1 h(\varphi,\theta) & := \big(1+\partial_{\theta}\beta_1(\varphi,\theta)\big) h\big(\varphi,\theta+\beta_1(\varphi,\theta)\big), 
\end{aligned}
\end{equation*}
such that 
\begin{align*}
\mathcal{T}_1(g):=\mathscr{B}_1^{-1}\mathcal{T}(g) \mathscr{B}_1=\varepsilon^2|\ln\varepsilon|\omega(\lambda)\partial_\varphi+\partial_\theta\Big(\big(\mathtt{c}_0(\lambda,g)+\varepsilon^5 {\mathcal{V}}_1(\varepsilon,g)\big)\cdot\Big),
\end{align*}
with the following properties:
\begin{enumerate}
\item The functions $b_k$ take  the form
\begin{align*}
 b_0(\varphi)&= O(|\ln\varepsilon|^{-\frac12}),
\\
    b_1(\varphi,\theta)&=\sin(\theta)\big(g_0-4\mathtt{g}_3\big)+ \mathtt{g}_3\sin(3\theta)-2\varepsilon|\ln\varepsilon|\big(\mathtt{g}_2 \sin(2\theta)+\mathtt{f}_2\cos(2\theta) \big)+O(\varepsilon|\ln\varepsilon|^{-\frac12}).
\end{align*}
For the definitions of ${\mathtt{f}}_2$ ${\mathtt{g}}_2$ and ${\mathtt{g}}_3$, see \eqref{asymt-list1}.
\item The constant $\mathtt{c}_0$ is given by
\begin{equation*}
\mathtt{c}_0:= -\tfrac12+\varepsilon^2|\ln\varepsilon|^\frac12\mathtt{r}_1, \quad \| \mathtt{r}_1 \|^{\textnormal{Lip},\nu}\lesssim 1.\end{equation*}
\item  The functions ${b}_0$, ${b}_1$ and  ${b}_{1,e}$ are odd and  ${\mathcal{V}}_1$ is even and they satisfy: for all $s\in[s_{0},s_\textnormal{up}],$ 
\begin{equation*}
\begin{aligned}
\|{b}_0\|_{s}^{\textnormal{Lip},\nu} +\|{b}_1\|_{s}^{\textnormal{Lip},\nu} +\|{b}_{1,e}\|_{s}^{\textnormal{Lip},\nu} +\|{\mathcal{V}_1}\|_{s}^{\textnormal{Lip},\nu} \lesssim 1+\|g\|_{s+{3}}^{\textnormal{Lip},\nu}.
\end{aligned}
\end{equation*}
\item The operators $\mathscr{B}_1^{\pm 1}$ are reversibility preserving on $\mathbb{X}^s$ and   satisfy: for all $s\in[s_{0},s_\textnormal{up}],$
\begin{equation*}
\begin{aligned}
\|\mathscr{B}^{\pm 1}_1 h\|_{s}^{\textnormal{Lip},\nu}
 \lesssim \|h\|_{s}^{\textnormal{Lip},\nu}+\varepsilon\| g\|_{s+{3}}^{\textnormal{Lip},\nu}\|h\|_{s_{0}}^{\textnormal{Lip},\nu}.
\end{aligned}
\end{equation*} 
	\item Given two small states $g_{1}$ and $g_{2}$, then 
				\begin{align*}
				\|\Delta_{12}\mathtt{c}_0\|^{\textnormal{Lip},\nu}&\lesssim\varepsilon^2|\ln\varepsilon|^\frac12\| \Delta_{12}g\|_{s_0+2}^{\textnormal{Lip},\nu},
				\\ 
				\|\Delta_{12}\mathcal{V}_1\|_{s}^{\textnormal{Lip},\nu}
&\lesssim\Big(\|\Delta_{12} g\|_{s+3}^{\textnormal{Lip},\nu}+{ \|\Delta_{12} g\|_{s_0+3}^{\textnormal{Lip},\nu}
\max_{j=1,2} \|g_j\|_{s+3}^{\textnormal{Lip},\nu}}
\Big).
			\end{align*}
\end{enumerate}
\end{pro}

\begin{proof}
   The change of coordinates $\mathscr{B}_1$  will be used here to reduce the size of the function $\mathcal{V}$. As stated in the proposition, we consider a symplectic periodic change of coordinates taking the form
\begin{align*}
    \mathscr{B}_1 h(\varphi,\theta):=(1+\partial_\theta \beta_1(\varphi,\theta))h(\varphi,\theta+\beta_1(\varphi,\theta)),
\end{align*}
for some $\beta_1$ that we will defined appropriately. Define also
$$
\mathtt{B}_1h(\varphi,\theta):=h(\varphi,\theta+\beta_1(\varphi,\theta)), 
$$
whose inverse takes the form
$$
\mathtt{B}^{-1}_1h(\varphi,\theta):=h(\varphi,\theta+\widehat{\beta}_1(\varphi,\theta)).
$$
Then, using Lemma \ref{algeb1} we infer
\begin{align*}
\mathscr{B}_1^{-1}\mathcal{T}\mathscr{B}_1=\varepsilon^2|\ln\varepsilon|\omega\partial_\varphi+\partial_\theta\Big\{(\mathtt{B}_1^{-1}\widetilde{V}_1)\cdot\Big\},
\end{align*}
where 
\begin{align*}
    \widetilde{V}_1(\varepsilon,g):=&\varepsilon^2|\ln\varepsilon|\omega\partial_\varphi \beta_1-\tfrac12(1+\partial_\theta \beta_1)+\varepsilon \mathcal{V}(\varepsilon,g)(1+\partial_\theta \beta_1)\\
=&\varepsilon^2|\ln\varepsilon|\omega\partial_\varphi \beta_1+\left( -\tfrac12 +\varepsilon\mathcal{V}(\varepsilon,g)\right)+\partial_\theta \beta_1 \left( -\tfrac12 +\varepsilon\mathcal{V}(\varepsilon,g)\right).
\end{align*}
Next, we want to reduce the size of ${V}_1(\varepsilon,f)$ around its average by making suitable ansatz. For this aim, we   define 
$$
\beta_1(\varphi,\theta):= b_0(\varphi)+\sum_{k=1}^4 \varepsilon^k b_k(\varphi,\theta),
$$
and therefore
\begin{align*}
   \widetilde{V}_1(\varepsilon,g)    &=\varepsilon^2|\ln\varepsilon|\omega\partial_\varphi b_0+\sum_{k=1}^4 \varepsilon^{2+k}|\ln\varepsilon|\omega\partial_\varphi b_k+\left( -\tfrac12 +\varepsilon\mathcal{V}(\varepsilon,g)\right)+\sum_{k=1}^4 \varepsilon^k  \partial_\theta b_k \left( -\tfrac12 +\varepsilon\mathcal{V}(\varepsilon,g)\right)\\  &=\varepsilon^2|\ln\varepsilon|\omega\partial_\varphi b_0 -\tfrac12 +\varepsilon\big(\mathcal{V}(\varepsilon,g)-\tfrac12  \partial_\theta b_1 \big) \\
  &\quad +\sum_{k=2}^4 \varepsilon^{k}\Big(\varepsilon\ln\varepsilon|\omega\partial_\varphi b_{k-1}-\tfrac12\partial_\theta b_k+\partial_\theta b_{k-1} \mathcal{V}(\varepsilon,g)\Big)+\varepsilon^{6}|\ln\varepsilon|\omega\partial_\varphi b_{4}+\varepsilon^5  \partial_\theta b_{4} \mathcal{V}(\varepsilon,g).
\end{align*}
Now, we shall impose the conditions
$$
\tfrac12 \partial_\theta b_1=\mathcal{V}(\varepsilon,g)-\langle\mathcal{V}(\varepsilon,g)\rangle_\theta,
$$
and  for $k=2,3,4,$
$$
\tfrac12 \partial_\theta b_k=\varepsilon |\ln\varepsilon| \omega\partial_\varphi b_{k-1}+\partial_\theta b_{k-1} \mathcal{V}(\varepsilon,g)-\langle\varepsilon |\ln\varepsilon| \omega\partial_\varphi b_{k-1}+\partial_\theta b_{k-1} \mathcal{V}(\varepsilon,g)\rangle_\theta.
$$
 These equations can be solved uniquely  with the constraints $\langle b_k\rangle_\theta=0.$ As $\langle \partial_\varphi b_k\rangle_\theta=0,$ then we deduce that 
\begin{align}\label{V1-jeu}
   \widetilde{V}_1(\varepsilon,g)    =&\varepsilon^2|\ln\varepsilon| \omega(\lambda)\partial_\varphi b_0-\tfrac12+\varepsilon\langle\mathcal{V}\rangle_\theta+\sum_{k=1}^3 \varepsilon^{k+1}  \langle\partial_\theta b_k \mathcal{V}(\varepsilon,g)\rangle_\theta\\
   \nonumber &+\varepsilon^5 \partial_\theta b_4 \mathcal{V}(\varepsilon,g)+\varepsilon^6 |\ln\varepsilon|\omega\partial_\varphi b_4.
\end{align}
Now, we define $b_0$ as the unique solution with zero average to the equation
\begin{align*}
\varepsilon^2|\ln\varepsilon| \omega\partial_\varphi b_0=&  \varepsilon\Big(\langle\mathcal{V}\rangle_{\varphi,\theta}-\langle\mathcal{V}\rangle_\theta\Big)+\sum_{k=1}^3 \varepsilon^{k+1}  \Big(\langle\partial_\theta b_k \mathcal{V}\rangle_{\varphi,\theta}-\langle\partial_\theta b_k \mathcal{V}\rangle_\theta\Big),
\end{align*}
so that
\begin{align*}
    \widetilde{V}_1(\varepsilon,g)    &=-\frac12+\varepsilon\langle\mathcal{V}\rangle_{\varphi,\theta}+\sum_{k=1}^3 \varepsilon^{k+1}  \langle\partial_\theta b_k \mathcal{V}(\varepsilon,g)\rangle_{\varphi,\theta}+\varepsilon^5 \partial_\theta b_4 \mathcal{V}(\varepsilon,g)+\varepsilon^6 |\ln\varepsilon|\omega\partial_\varphi b_4.
\end{align*}
From \eqref{mathcalV} and \eqref{asymt-list1} we deduce that
\begin{align*}
\langle  \mathcal{V}\rangle_\theta
   =&O(\varepsilon|\ln\varepsilon|^\frac12),\quad \textnormal{and}\quad \langle  \mathcal{V}\rangle_{\varphi,\theta}
   =O(\varepsilon|\ln\varepsilon|^\frac12).
\end{align*}
Thus,
\begin{align*}
\partial_\varphi b_0=& O(|\ln\varepsilon|^{-\frac12}),
\end{align*}
which implies that the unique periodic solution with zero average is given by 
\begin{align*}
 b_0(\varphi)=&  O(|\ln\varepsilon|^{-\frac12}) .
\end{align*}
Set
\begin{align*}
\mathtt{c}_0&:=-\frac12+\varepsilon\langle\mathcal{V}\rangle_{\varphi,\theta}+\sum_{k=1}^3 \varepsilon^{k+1}  \langle\partial_\theta b_k \mathcal{V}\rangle_{\varphi,\theta}\quad\hbox{and} \quad \widehat{\mathcal{V}}_1:=\partial_\theta b_4\, \mathcal{V}+\varepsilon |\ln\varepsilon|\omega\partial_\varphi b_4,
    \end{align*}
then $\mathtt{c}_0$  admits the asymptotics
\begin{align*}
\mathtt{c}_0&=-\frac12+O(\varepsilon^2|\ln\varepsilon|^\frac12),
    \end{align*}
    and
    \begin{align*}
\mathtt{B}_1^{-1}\widetilde{V}_1=\mathtt{c}_0+\varepsilon^5 \mathcal{V}_1.
\end{align*}
It follows that
\begin{align*}
\mathcal{T}_1:=\mathscr{B}_1^{-1}\mathcal{T}\mathscr{B}_1=\varepsilon^2|\ln\varepsilon|\omega\partial_\varphi+\partial_\theta\Big\{\big(\mathtt{c}_0+\varepsilon^5 \mathcal{V}_1\big)\cdot\Big\}.
\end{align*}
It remains to 
get an explicit expression of $b_1$ up to the order $O(\varepsilon|\ln\varepsilon|)$. Recall that
$$
\partial_\theta b_1=2\mathcal{V}(\varepsilon,g)-2\langle\mathcal{V}(\varepsilon,g)\rangle_\theta.
$$
Then, using the expression of $\mathcal{V}$ in \eqref{mathcalV} we infer
\begin{align*}
 \partial_\theta b_1=&3 \mathtt{g}_{3}(\varphi)\cos(3\theta)+\big( g_0-4\mathtt{g}_{3}\big)\cos(\theta)-4\varepsilon|\ln\varepsilon| \mathtt{g}_{2}\cos(2\theta)\\
 &+4\varepsilon|\ln\varepsilon| \mathtt{f}_{2}(\varphi)\sin(2\theta)+O(\varepsilon|\ln\varepsilon|^\frac12).
\end{align*}
Thus by integration we get
\begin{align*}
    b_1(\varphi,\theta)=&\sin(\theta)\big(g_0-4\mathtt{g}_3\big)+ \mathtt{g}_3\sin(3\theta)-2\varepsilon|\ln\varepsilon|\big(\mathtt{g}_2 \sin(2\theta)+\mathtt{f}_2\cos(2\theta)\big) +O(\varepsilon|\ln\varepsilon|^\frac12),
\end{align*}
which concludes the proof.
\end{proof}
In Proposition~\ref{First-change0}, we carried out a preliminary transformation that reduced the transport 
part to an operator whose coefficients are variable but of higher order, namely of size $\varepsilon^5$. The new operator takes the form
\begin{align*}
\mathcal{T}_1(g)
= \varepsilon^2 |\ln \varepsilon| \omega(\lambda) \, \partial_\varphi
+ \partial_\theta \Big( \big( \mathtt{c}_0 + \varepsilon^5 \mathcal{V}_1(\varepsilon, g) \big) \cdot \Big).
\end{align*}  
The natural continuation of this procedure is a complete reduction, where one seeks to fully diagonalize 
the resulting operator.
This step is achieved through a KAM-type scheme, inspired by the pioneering strategies of~\cite{FGMP19, HHM21, HR21}. 
However, unlike the classical setting, our problem exhibits a degeneracy in the time direction, 
which requires specific adjustments in the iterative procedure. These refinements have been implemented in detail 
in~\cite{HHM24}. For the purposes of this work, we shall not reproduce the technical construction, 
but simply state the main result that follows from this program.  \\
The following result is  an immediate consequence of \cite[Proposition 4.2]{HHM24}.

\begin{pro}\label{Second-change}
Let $(\nu,\tau,s_0,s_{\textnormal{up}})$ as in \eqref{cond1}. Fix $\mu_2\in\R$ such that
\begin{equation*}
\mu_{2}\geqslant 4\tau+3
\end{equation*}
and let $g\in \mathbb{X}^{s_\textnormal{up}}$.
There exists ${\epsilon}_0>0$ such that if \eqref{small-C2-0} holds, then, we can construct $\mathtt{c}_1:= \mathtt{c}_1(\lambda,g)\in \textnormal{Lip}_\nu(\mathcal{O},\mathbb{R})$ and $\beta_2\in \mathbb{Y}^{s_\textnormal{up}}$
such that with $\mathscr{B}_2$ as in \eqref{definition symplectic change of variables} one gets the following results.
\begin{enumerate}
\item  The function $\lambda\mapsto\mathtt{c}_1(\lambda,g)$ $($which is constant with respect to the time-space variables$)$  satisfies the following estimate, 
\begin{equation*}
\| \mathtt{c}_1-\mathtt{c}_0 \|^{\textnormal{Lip},\nu}\lesssim {\varepsilon}^{{5}},
\end{equation*}
where $\mathtt{c}_0$ is defined in Proposition $\ref{First-change0}.$
\item The transformations $\mathscr{B}_2^{\pm 1}$ are reversibility preserving on $\mathbb{X}^s$ and the functions  $ {\beta_2}, \widehat{\beta}_2$ are odd and they  satisfy, for all $s\in[s_{0},s_{\textnormal{up}}],$ 
\begin{equation*}
\|\mathscr{B}_2^{\pm 1}h\|_{s}^{\textnormal{Lip},\nu}
 \lesssim\|h\|_{s}^{\textnormal{Lip},\nu}+{\varepsilon^{{5}}\nu ^{-1}}\| g\|_{s+2\tau+{5}}^{\textnormal{Lip},\nu}\|h\|_{s_{0}}^{\textnormal{Lip},\nu},
\end{equation*}
and 
\begin{equation*}
\|\widehat{\beta}_2\|_{s}^{\textnormal{Lip},\nu}\lesssim\|\beta_2\|_{s}^{\textnormal{Lip},\nu}\lesssim \varepsilon^{{5}}\nu^{-1}\left(1+\| g\|_{s+2\tau+{4}}^{\textnormal{Lip},\nu}\right).
\end{equation*}
\item Let $n\in\N$, then on  the Cantor set
\begin{equation*}
{\mathcal{O}_{n}^{1}(g)}=\bigcap_{\substack{(\ell,j)\in\mathbb{Z}^{2}\\ 1\leqslant |j|\leqslant N_{n}}}\Big\lbrace \lambda\in \mathcal{O};\;\, \big|\varepsilon^2\,|\ln\varepsilon|\omega(\lambda)  \ell+j \mathtt{c}_1(\lambda,g)\big|\geqslant{\nu}{| j|^{-\tau}}\Big\rbrace,
\end{equation*}
we have 
\begin{align*}
\mathcal{T}_2:=\mathscr{B}_2^{-1}\mathcal{T}_1(g)\mathscr{B}_2=\varepsilon^2|\ln\varepsilon|\omega(\lambda)\partial_\varphi+ \mathtt{c}_1(\lambda,g)\partial_{\theta}+\mathtt{E}_{n},
\end{align*}
with $\mathtt{E}_{n}$ a linear operator satisfying
\begin{equation*}
\|\mathtt{E}_{n}[h]\|_{s_0}^{\textnormal{Lip},\nu}\lesssim \varepsilon^{{5}} N_{0}^{\mu_{2}}N_{n+1}^{-\mu_{2}}\|h\|_{s_{0}+2}^{\textnormal{Lip},\nu}.
\end{equation*}
	\item Given two functions $g_{1}$ and $g_{2}$ both satisfying \eqref{small-C2-0}, we have 
			\begin{align*}
				\|\Delta_{12}\mathtt{c}_1\|^{\textnormal{Lip},\nu}&\lesssim \varepsilon^2|\ln\varepsilon|^\frac12\| \Delta_{12}g\|_{2{s}_{0}+2\tau+3}^{\textnormal{Lip},\nu}.
			\end{align*}
            \item Let $\mathscr{B}:= \mathscr{B}_1\mathscr{B}_2$, then $\mathscr{B}$ is reversibility preserving on $\mathbb{X}^s$, invertible and satisfies the estimate, for any $s\in[s_0,s_{\textnormal{up}}]$
\begin{align*}
\|\mathscr{B}^{\pm 1} h\|_{s}^{\textnormal{Lip},\nu}& \lesssim   \|h\|_{s}^{\textnormal{Lip},\nu}+ \| g\|_{s+2\tau+{5}}^{\textnormal{Lip},\nu}\|h\|_{s_0}^{\textnormal{Lip},\nu}.
\end{align*} 
In addition, for $\lambda \in {\mathcal{O}_{n}^{1}(g)}$
\begin{align*}
  \mathscr{B}^{-1}\mathcal{T}(g)\mathscr{B}= \varepsilon^2|\ln\varepsilon|\omega\partial_\varphi+ \mathtt{c}_1(\lambda,g)\partial_{\theta}+\mathtt{E}_{n}.
\end{align*}
Moreover,
\begin{equation*}
\mathscr{B}h(\lambda,\varphi,\theta)=\big(1+\partial_{\theta}\beta(\lambda,\varphi,\theta)\big)  h\big(\lambda,\varphi,\theta+\beta(\lambda,\varphi,\theta)\big),
\end{equation*}
where $\beta(\lambda,\varphi,\theta)$ is given by
\begin{align*}
  \beta(\lambda,\varphi,\theta)&=\beta_1(\lambda,\varphi,\theta)+\beta_2\left(\lambda,\varphi,\theta+\beta_1(\lambda,\varphi,\theta)\right) \\ 
  &={b}_0(\varphi)+\varepsilon {b}_1(\varphi,\theta)+\varepsilon^2|\ln\varepsilon|^{\frac12}\mathtt{b}(\varphi,\theta)
\end{align*}
and satisfies the estimate
\begin{equation*}
\|\beta\|_{s}^{\textnormal{Lip},\nu}+\|\mathtt{b}\|_{s}^{\textnormal{Lip},\nu}\lesssim  
1+\| g\|_{s+2\tau+{4}}^{\textnormal{Lip},\nu},
\end{equation*}
where the function $b_0$ and $b_1$ are defined in Proposition $\ref{First-change0}.$
\end{enumerate}
\end{pro}

\subsubsection{Conjugation of the nonlocal operators}\label{section-CFLO}

With the transport part addressed, we proceed to compute the impact of the transformation $\mathscr{B}$  on the nonlocal terms. Recall from \eqref{Lperp} that
\begin{align*}
   \mathcal{L}_\perp (\varepsilon,g) 
   = { \Pi^\perp }\Big[\mathcal{T} + \tfrac12 \mathcal{H} +\varepsilon \mathcal{H}_{\mathtt{u},0}
   + \varepsilon \partial_\theta \mathcal{S} 
   + \varepsilon^2 \partial_\theta \mathcal{R}(g)\Big]{\Pi^\perp }.
\end{align*}
Here, $\mathcal{H}$, $\mathcal{H}_{\mathtt{u},0}$ and $\mathcal{S}$ denote lower-order components whose precise structure 
has already been described in \eqref{Hilbert1alt}, \eqref{def-Hu0} and \eqref{shift-operator1}, respectively, while $\mathcal{R}(g)$ represents a remainder operator in~\eqref{RId-10},
\begin{align*}
 \mathcal{R}(g)[h](\varphi,\theta) 
   = \fint_{\T} \mathcal{W}(g)(\varphi,\theta,\eta)\,
      \ln\!\Big|\sin\!\Big(\tfrac{\theta-\eta}{2}\Big)\Big|\,
      h(\varphi,\eta)\, d\eta 
      + \mathcal{R}_{\infty}(g)[h](\varphi,\theta).
\end{align*}
This operator is characterized by a singular integral kernel involving the logarithm of a sine function, 
supplemented by a smoother contribution $\mathcal{R}_\infty(g)$.  
In what follows, we analyze the asymptotic structure of the nonlocal operators and their conjugation with the transformation $
   \mathscr{B} = \mathscr{B}_2 \mathscr{B}_1,$
constructed in Proposition \ref{Second-change}.  

Before stating the main result, we introduce an operator that arises in the analysis of the conjugation of the various components of $\mathcal{L}_\perp$.
  \begin{align}\label{S1-defi1}
&\mathcal{S}_0[h](\varphi,\theta):=  - 2\,\mathtt{g}_3\fint_{\mathbb{T}}\Big(\cos\big(2\theta+\eta-3{b_0}\big)+\cos\big(\theta+2\eta-3b_0\big)\Big)h\big(\varphi,\eta\big)\,d\eta\\
&\nonumber-2\varepsilon|\ln\varepsilon| \left(\mathtt{g}_2\fint_{\mathbb{T}}\cos\big(\theta+\eta-2b_0\big)h\big(\varphi,\eta\big)\,d\eta-\mathtt{f}_2\fint_{\mathbb{T}}\sin\big(\theta+\eta-2b_0\big)h\big(\varphi,\eta\big)\,d\eta\right),
\end{align}
where $\mathtt{f}_2$, $\mathtt{g}_2$, $\mathtt{g}_3$  are defined in \eqref{list-functions} and $b_0$ is introduced in Proposition \ref{First-change0}.
This operator is finite-dimensional and acts by localizing onto Fourier modes of order less than~$2$.
Our main result is the following.
  \begin{pro}\label{prop-asym-brut}
Assume the notation and hypotheses of Proposition~$\ref{Second-change}$ and suppose that
      $$
\varepsilon^3|\ln\varepsilon|^{-\frac12}\nu^{-1}\leqslant 1.
$$
Then, for every $\lambda$ in the Cantor set $\mathcal{O}_{n}^{1}(g)$, the conjugated operator
$$
\mathcal{L}_{\mathscr{B}}
:=\mathscr{B}^{-1}\Big(\mathcal{T} + \tfrac12 \mathcal{H} +\varepsilon \mathcal{H}_{\mathtt{u},0}
   + \varepsilon \partial_\theta \mathcal{S} 
   + \varepsilon^2 \partial_\theta \mathcal{R}(g)\Big)\mathscr{B}
$$  
admits the decomposition
      \begin{align*}
\mathcal{L}_{\mathscr{B}}
&=\varepsilon^2|\ln\varepsilon|\omega(\lambda)\partial_\varphi
+\mathtt{c}_1(\lambda,g)\partial_{\theta}
+\tfrac12\mathcal{H}
+\varepsilon \partial_\theta\mathcal{S}_0
+\varepsilon\mathcal{H}_{\mathtt{u}_1,0}
+\varepsilon^2\mathcal{H}_{\mathtt{u}_2,-2}
+\varepsilon^2\partial_\theta\Lambda_{4,\mathcal{W}_3}\\
&\quad+\varepsilon^2|\ln\varepsilon|^{\frac12}\partial_\theta\mathcal{R}_{0,\infty}
+\mathtt{E}_{n}.
    \end{align*}
Moreover, the following properties hold:
    \begin{enumerate}
    \item The functions $\mathtt{u}_1$ and $\mathtt{u}_2$ are even and satisfy
        $$
    \|\mathtt{u}_1\|_{s}^{\textnormal{Lip},\nu}\lesssim 1+\|g\|_{s+1}^{\textnormal{Lip},\nu}\qquad \|\mathtt{u}_2\|_{s}^{\textnormal{Lip},\nu}\lesssim 1+\|g\|_{s+3}^{\textnormal{Lip},\nu},
    $$
  {{and}{
 $$
    \|\Delta_{12}\mathtt{u}_1\|_{s}^{\textnormal{Lip},\nu}\lesssim \|\Delta_{12}g\|_{s+1}^{\textnormal{Lip},\nu}\qquad \|\Delta_{12}\mathtt{u}_2\|_{s}^{\textnormal{Lip},\nu}\lesssim \|\Delta_{12}g\|_{s+3}^{\textnormal{Lip},\nu} .
    $$ 
  }}
    \item The operator $\partial_\theta\Lambda_{4,\mathcal{W}_3}\in \mathrm{OPS}^{-4}$ is a reversible, inhomogeneous Fourier multiplier of the form \eqref{eq:Lambda-m-f}, associated with a function $\mathcal{W}_3$, and it satisfies
\begin{align*}
\|\partial_\theta\Lambda_{4,\mathcal{W}_3}[h]\|_{s}^{\mathrm{Lip},\nu}
&\lesssim \|h\|_{s}^{\mathrm{Lip},\nu}
+\|g\|_{s+6}^{\mathrm{Lip},\nu}\,\|h\|_{s_0}^{\mathrm{Lip},\nu},\\
\|\mathcal{W}_3\|_{s}^{\mathrm{Lip},\nu}
&\lesssim 1+\|g\|_{s+5}^{\mathrm{Lip},\nu}.
\end{align*}
\item The operator $\partial_\theta\mathcal{R}_{0,\infty}$ is reversible and belongs to $\mathrm{OPS}^{-\infty}$. More precisely, 
    $$
     \forall n\geqslant1,\quad\|\partial_\theta^n\mathcal{R}_{0,\infty}[h]\|_{s}^{\textnormal{Lip},\nu}\lesssim \|h\|_{s}^{\textnormal{Lip},\nu}+\|g\|_{s+n}^{\textnormal{Lip},\nu}\|h\|_{s_0}^{\textnormal{Lip},\nu}.
    $$
    \end{enumerate}

      \end{pro}
      \begin{proof}
      
  Let us conjugate with $\mathscr{B}_1$ each term.
Note that the inverse diffeomorphism $\mathscr{B}_1^{-1}$ admits the form
\begin{equation*}
\mathscr{B}_1^{-1}h(\varphi,y)=\big(1+\partial_y\widehat{\beta_1}(\varphi,y)\big)h\big(\varphi,y+\widehat{\beta_1}(\varphi,y)\big).
\end{equation*}
According to Proposition \ref{First-change0} one gets
\begin{align*}
\beta_1(\varphi,\theta)&={b}_0(\varphi)+\varepsilon {b}_1(\varphi,\theta)+\varepsilon^2|\ln\varepsilon|^{\frac12}{b}_{1,e}(\varphi,\theta).
\end{align*}
Applying Lemma \ref{beta-inv-asym} leads to
\begin{align}\label{Expand-M100}
\nonumber\widehat{\beta_1}(\varphi,\theta)=&-{b}_0(\varphi)-\varepsilon b_1(\varphi,\theta-b_0(\varphi))+\varepsilon^{2}|\ln\varepsilon|^{\frac12}\mathfrak{r}_1\big(\varphi,\theta\big)=:-{b}_0(\varphi)-\varepsilon\rho_1,
\\
&\qquad\textnormal{with}\qquad  \|\mathfrak{r}_1\|_{s}^{{\textnormal{Lip}, \nu}}\lesssim 1+\|g\|_{{s+3}}^{{\textnormal{Lip}, \nu}}.
\end{align}
Using Proposition \ref{Second-change} we find
\begin{align}\label{comp-id-com}
\mathscr{B}^{-1} \Big[\mathcal{T} + \tfrac12 \mathcal{H}&+\varepsilon \mathcal{H}_{\mathtt{u},0} 
   + \varepsilon \partial_\theta \mathcal{S} 
   + \varepsilon^2 \partial_\theta \mathcal{R}(g)\Big] \mathscr{B}=\varepsilon^2|\ln\varepsilon|\omega(\lambda) \partial_\varphi+\mathtt{c}_1(\lambda,g)\partial_{\theta}+\mathtt{E}_{n}\\ &+\tfrac12\mathscr{B}^{-1}\mathcal{H} \mathscr{B}+\varepsilon\mathscr{B}^{-1} \mathcal{H}_{\mathtt{u},0}\mathscr{B}
\nonumber+\varepsilon \mathscr{B}^{-1}{\partial_\theta}\mathcal{S}\mathscr{B}+\varepsilon^2 \mathscr{B}^{-1}{\partial_\theta}\mathcal{R}(\varepsilon,g)\mathscr{B}.
    \end{align}
 Let us start with the conjugation of the Hilbert transform term.  We employ the following decomposition:
\begin{align}\label{comm H}
\nonumber\mathscr{B}^{-1}\mathcal{H}\,\mathscr{B}-\mathcal{H}
&= \mathscr{B}_2^{-1}\mathscr{B}_1^{-1}\mathcal{H}\,\mathscr{B}_2\,\mathscr{B}_1-\mathcal{H} \\
\nonumber &= \mathscr{B}_2^{-1}\big[\mathscr{B}_1^{-1}\mathcal{H}\,\mathscr{B}_1-\mathcal{H}\big]\mathscr{B}_2 + \mathscr{B}_2^{-1}\mathcal{H}\,\mathscr{B}_2-\mathcal{H} \\
&= {\mathscr{B}_1^{-1}\mathcal{H}\,\mathscr{B}_1-\mathcal{H}}
+ \left(\mathscr{B}_2^{-1}\big[{\mathscr{B}_1^{-1}\mathcal{H}\,\mathscr{B}_1-\mathcal{H}}\big]\,\mathscr{B}_2-\big[{\mathscr{B}_1^{-1}\mathcal{H}\,\mathscr{B}_1-\mathcal{H}}\big]\right) 
\\ &\quad \nonumber+ \mathscr{B}_2^{-1}\mathcal{H}\,\mathscr{B}_2-\mathcal{H}.
\end{align}
Using the general formula,
\begin{align}\label{general-formula}
\mathscr{B}_1^{-1}\left(\partial_\theta\fint_{\T} K(\varphi,\theta,\eta)h(\varphi,\eta)d\eta \right)\mathscr{B}_1&=\partial_\theta\fint_{\T} K\big(\varphi,\theta+\widehat{\beta}_1(\varphi,\theta),\eta+\widehat{\beta}_1(\varphi,\eta)\big)h(\varphi,\eta)d\eta,
\end{align}
together with the expansion of $\widehat{\beta_1}$, we can get 
\begin{align*}
 \mathcal{H}_1[h](\varphi,\theta)&:=\tfrac12\left(\mathscr{B}_1^{-1}\mathcal{H} \mathscr{B}_1-\mathcal{H}\right)[h](\varphi,\theta)=\partial_\theta\fint_{\mathbb{T}}\mathcal{K}_1(\varphi,\theta,\eta)h(\varphi,\eta)\,d\eta,
\end{align*}
where
\begin{align*}
 \mathcal{K}_1(\varphi,\theta,\eta)&:=\ln\left|\frac{e^{i(\theta+\widehat{\beta}_1(\lambda,\varphi,\theta))}-e^{i (\eta+\widehat{\beta}_1(\lambda,\varphi,\eta))}}{e^{i\theta}-e^{i\eta}} \right|=\varepsilon \mathcal{K}_{1,1}(\varphi,\theta,\eta)+\varepsilon^{2}{\mathcal  K}_{1,2},
\end{align*}
with
\begin{align*}
\mathcal{K}_{1,1}(\varphi,\theta,\eta)&:=\hbox{Im}\bigg\{\frac{e^{i\theta}{b_1\big(\varphi,\theta-b_0(\varphi)\big)}-e^{i\eta}b_1(\varphi,\eta-{b_0(\varphi)})}{e^{i\theta} -e^{i\eta}} \bigg\}\\
&=\tfrac12\left(b_1\left(\varphi,\theta-b_0(\varphi)\right)-b_1\left(\varphi,\eta-b_0(\varphi)\right)\right)\cot\big(\tfrac{\eta-\theta}{2}\big)
\end{align*}
and
\begin{align*}
 \| {\mathcal  K}_{1,2}\|_{s}^{{\textnormal{Lip},\nu}}\lesssim_s 1+\| g\|_{s+3}^{\textnormal{Lip},\nu}.
\end{align*}
Therefore, making a change of variables we get
\begin{align*}
\mathcal{H}_1[h](\varphi,\theta+b_0(\varphi))=&\varepsilon\partial_\theta\fint_{\mathbb{T}}\mathcal{K}_{1,1}\left(\varphi,\theta+b_0(\varphi),\eta+b_0(\varphi)\right)h\big(\varphi,\eta+b_0(\varphi)\big)\,d\eta\\ &+\varepsilon^2\fint_{\mathbb{T}}  {\mathcal  K}_{1,2}( \varphi, \theta,\eta) h( \varphi, \eta)\,d\eta,
\end{align*}
and
\begin{align*}
\mathcal{K}_{1,1}\left(\varphi,\theta+b_0(\varphi),\eta+b_0(\varphi)\right)
&=\tfrac12\left(b_1\big(\varphi,\theta\big)-b_1\big(\varphi,\eta\big)\right)\cot\big(\tfrac{\eta-\theta}{2}\big).
\end{align*}
Recall that
\begin{align*}
    b_1(\varphi,\theta)=&\sin(\theta)\big(g_0-4\mathtt{g}_3\big)+ \mathtt{g}_3\sin(3\theta)-2\varepsilon|\ln\varepsilon|\big(\mathtt{g}_2 \sin(2\theta)+\mathtt{f}_2\cos(2\theta)\big)+O(\varepsilon|\ln\varepsilon|^\frac12).
\end{align*}
Using the trigonometric identities
\begin{align*}
    (\sin(\theta)-\sin(\eta))\cot\big(\tfrac{\eta-\theta}{2}\big)&=-{(\cos \theta+\cos\eta}),\\
   (\sin(2\theta)-\sin(2\eta))\cot\big(\tfrac{\eta-\theta}{2}\big)&=2\cos(\theta+\eta)+\cos(2\theta)+\cos(2\eta),
\end{align*}
together with
\begin{align*}
   (\cos(2\theta)-\cos(2\eta))\cot\big(\tfrac{\eta-\theta}{2}\big)&=-2\sin(\theta+\eta)-\sin(2\theta)-\sin(2\eta),\\
   (\sin(3\theta)-\sin(3\eta))\cot\big(\tfrac{\theta-\eta}{2}\big)&=\cos(3\theta)+\cos(3\eta)+2\cos(2\theta+\eta)+2\cos(\theta+2\eta),
\end{align*}
and
\begin{align*}
    (\sin(4\theta) - \sin(4\eta))  \cot\big( \tfrac{\theta - \eta}{2} \big)
&= 2\cos(2\theta + 2\eta) + \cos(4\theta) + \cos(4\eta) + 2\cos(3\theta + \eta) \\ &\quad + 2\cos(\theta + 3\eta),
\end{align*}
and the fact that $\langle h\rangle_\theta=0$, we get 
\begin{align*}
\mathcal{H}_1[h](\varphi,\theta+b_0(\varphi))&=  -\varepsilon  \mathtt{g}_3\partial_\theta\fint_{\mathbb{T}}(\cos(2\theta+\eta)+\cos(\theta+2\eta))h\big(\varphi,\eta+b_0(\varphi)\big)\,d\eta\\
&\quad-2\varepsilon^2|\ln\varepsilon| \mathtt{g}_2(\varphi)\partial_\theta\fint_{\mathbb{T}}\cos(\theta+\eta)h\big(\varphi,\eta+b_0(\varphi)\big)\,d\eta\\
&\quad+2\varepsilon^2|\ln\varepsilon| \mathtt{f}_2(\varphi)\partial_\theta\fint_{\mathbb{T}}\sin(\theta+\eta)h\big(\varphi,\eta+b_0(\varphi)\big)\,d\eta\\
&\quad+\varepsilon^2|\ln\varepsilon|^{\frac12}\fint_{\mathbb{T}}  {\mathcal  K}_{1,3}( \varphi, \theta,\eta) h( \varphi, \eta)\,d\eta,
\end{align*}
where 
\begin{align*}
 \| {\mathcal  K}_{1,3}\|_{s}^{{\textnormal{Lip},\nu}}\lesssim_s 1+\| g\|_{s+3}^{\textnormal{Lip},\nu}.
\end{align*}
Using a change of variables, it follows that
\begin{align*}
\mathcal{H}_1[h](\varphi,\theta)&=  - \varepsilon\mathtt{g}_3\partial_\theta\fint_{\mathbb{T}}\Big(\cos\big(2\theta+\eta-3{b_0}\big)+\cos\big(\theta+2\eta-3b_0\big)\Big)h\big(\varphi,\eta\big)\,d\eta\\
&\quad-2\varepsilon^2|\ln\varepsilon| \mathtt{g}_2\partial_\theta\fint_{\mathbb{T}}\cos\big(\theta+\eta-2b_0\big)h\big(\varphi,\eta\big)\,d\eta\\
&\quad+2\varepsilon^2|\ln\varepsilon| \mathtt{f}_2\partial_\theta\fint_{\mathbb{T}}\sin\big(\theta+\eta-2b_0\big)h\big(\varphi,\eta\big)\,d\eta\\
&\quad+\varepsilon^2|\ln\varepsilon|^{\frac12}\fint_{\mathbb{T}}  {\mathcal  K}_{1,4}( \varphi, \theta,\eta) h( \varphi, \eta)\,d\eta,
\end{align*}
where
\begin{align*}
 \| {\mathcal  K}_{1,4}\|_{s}^{{\textnormal{Lip},\nu}}\lesssim_s 1+\| g\|_{s+3}^{\textnormal{Lip},\nu}.
\end{align*}
Inserting this identity into \eqref{comm H} gives
\begin{align*}
\tfrac12\left(\mathscr{B}^{-1}\mathcal{H} \mathscr{B}-\mathcal{H}\right)= \mathcal{H}_1+\varepsilon^2|\ln\varepsilon|^{\frac12}\mathcal{H}_3,
\end{align*} 
where
\begin{align*}
\nonumber\mathcal{H}_3[h](\varphi,\theta)
=&\fint_{\mathbb{T}}  {\mathcal  K}_{1,4}( \varphi, \theta,\eta) h( \varphi, \eta)\,d\eta+\varepsilon^{-2}|\ln\varepsilon|^{-\frac12}\big(\mathscr{B}_2^{-1}\mathcal{H}_1\,\mathscr{B}_2-\mathcal{H}_1\big)[h](\varphi,\theta) 
\\ &+ \varepsilon^{-2}|\ln\varepsilon|^{-\frac12}\big(\mathscr{B}_2^{-1}\mathcal{H}\,\mathscr{B}_2-\mathcal{H}\big)[h](\varphi,\theta).
\end{align*}
Notice, in view of  Proposition \ref{Second-change}, that
\begin{equation}\label{est-beta-r0}
\|\widehat{\beta}_2\|_{s}^{\textnormal{Lip},\nu}\lesssim\|\beta_2\|_{s}^{\textnormal{Lip},\nu}\lesssim \varepsilon^5\nu^{-1}\left(1+\| g\|_{s+2\tau+{4}}^{\textnormal{Lip},\nu}\right).
\end{equation}
In view of of this estimate,  Lemma \cite[Lemma 2.36]{BertiMontalto},  \eqref{est-beta-r0} and \eqref{small-C2-0}, the operator
\begin{align*}
		\big(\mathscr{B}_2^{-1}\mathcal{H}\,\mathscr{B}_2-\mathcal{H}\big) h(\varphi, \theta) 
		&= \fint_{\mathbb{T}}  {\mathcal  K}_{1,5}( \varphi, \theta,\eta) h( \varphi, \eta)\,d\eta
\end{align*}
		defines  an integral operator whose kernel satisfies the   estimates : for all $s\in[ s_0,s_{\textnormal{up}}],$  
\begin{align*}
 \| {\mathcal  K}_{1,5}\|_{s}^{{\textnormal{Lip},\nu}}\lesssim_s \varepsilon^{{5}}\nu^{-1}\left(1+\| g\|_{s+2\tau+{5}}^{\textnormal{Lip},\nu}\right).
\end{align*}
Moreover, according to  \cite[Lemma 2.3]{BertiMontalto}, \eqref{est-beta-r0} and \eqref{small-C2-0},  we have
\begin{align*}
\nonumber				\big(\mathscr{B}_2^{-1}\mathcal{H}_1\mathscr{B}_2-\mathcal{H}_1\big)h(\varphi,\theta)&=\fint_{\mathbb{T}}h(\varphi,{\eta}){\mathcal  K}_{1,6}(\varphi,\theta,{\eta})d{\eta},
		\end{align*}
		with
        \begin{align*}
 \| {\mathcal  K}_{1,6}\|_{s}^{{\textnormal{Lip},\nu}}\lesssim_s \varepsilon^{{5}}\nu^{-1}\left(1+\| g\|_{s+2\tau+{5}}^{\textnormal{Lip},\nu}\right).
\end{align*}
Combining the previous kernel estimates with Lemma \ref{lemma HR22}, we deduce that the integral operator  $\mathcal{H}_3$ is an integral operator  in  $\textnormal{OPS}^{-\infty}$.
\\
Let us now turn to the conjugation of $\mathcal{H}_{\mathtt{u},0}$, defined in \eqref{def-Hu0}.
Combining Lemma \ref{beta-inv-asym} with Proposition \ref{Second-change}-(5) yields 
\begin{align}\label{Expand-M1002}
\nonumber\widehat{\beta}(\varphi,\theta)=&-{b}_0(\varphi)-\varepsilon b_1(\varphi,\theta-b_0(\varphi))+\varepsilon^{2}|\ln\varepsilon|^{\frac12}\mathfrak{r}_2\big(\varphi,\theta\big)=:-{b}_0(\varphi)-\varepsilon\rho_2(\varphi,\theta),
\\
&\qquad\textnormal{with}\qquad  \|\mathfrak{r}_2\|_{s}^{{\textnormal{Lip}, \nu}}\lesssim 1+\|g\|_{{s+3}}^{{\textnormal{Lip}, \nu}}.
\end{align}
A direct computation, using \eqref{general-formula} and \eqref{Expand-M1002} yields
\begin{align*}
\mathscr{B}^{-1}\mathcal{H}_{\mathtt{u},0}\,\mathscr{B}[h](\varphi,\theta)
&= -\tfrac{1}{2}(2\mathtt{p}_{1,1})^{-\frac{3}{4}} \fint_{\mathbb{T}} K_{\mathrm{I}}(\varphi,\theta,\eta)\,h(\varphi,\eta)\,d\eta,
\end{align*}
with
\begin{align*}
{K}_{\textnormal{I}}(\varphi,\theta,\eta)&:=\Big(\cos\big(\theta-b_0-\varepsilon\rho_2(\theta)\big)+\cos\big(\eta-b_0-\varepsilon\rho_2(\eta)\big)\Big)\ln(|e^{i\theta-i\varepsilon\rho_2(\theta)}-  e^{i\eta-i\varepsilon\rho_2(\eta)}|)\\
&=\Big(\cos\big(\theta-b_0-\varepsilon\rho_2(\theta)\big)+\cos\big(\eta-b_0-\varepsilon\rho_2(\eta)\big)\Big)\ln(|e^{i\theta}-  e^{i\eta}|)\\
&\quad +\Big(\cos\big(\theta-b_0-\varepsilon\rho_2(\theta)\big)+\cos\big(\eta-b_0-\varepsilon\rho_2(\eta)\big)\Big)\ln\Bigg(\frac{|e^{i\theta-i\varepsilon\rho_2(\theta)}-  e^{i\eta-i\varepsilon\rho_2(\eta)}|}{|e^{i\theta}-e^{i\eta}|}\Bigg).
\\&=:{K}_{\textnormal{I},0}(\varphi,\theta,\eta)\ln(|e^{i\theta}-  e^{i\eta}|)+\varepsilon{K}_{\textnormal{I},\infty}(\varphi,\theta,\eta).
\end{align*}
Notice that ${K}_{\textnormal{I},\infty}$ is smooth and satisfies the estimate
        \begin{align*}
 \| {K}_{\textnormal{I},\infty}\|_{s}^{{\textnormal{Lip},\nu}}\lesssim_s 1+\| g\|_{s+1}^{\textnormal{Lip},\nu}.
\end{align*}
Moreover, by using Taylor expansion in $\varepsilon$ gives
  \begin{align*}
     K_{\textnormal{I},0}(\varphi,\theta,\eta)&=\cos(\theta-b_0(\varphi))+\cos(\eta-b_0(\varphi))+\varepsilon\rho_2(\varphi,\theta)\int_0^1 \sin\big(\theta-b_0(\varphi)-\varepsilon \tau \rho_2(\varphi,\theta)\big)d\tau\\
  &\quad+\varepsilon \rho_2(\varphi,\eta)\int_0^1 \sin\big(\eta-b_0(\varphi)-\varepsilon \tau \rho_2(\varphi,\eta)\big)d\tau. 
  \end{align*}
 It follows that
  \begin{align*}
\nonumber\mathscr{B}^{-1}\mathcal{H}_{\mathtt{u},0}\mathscr{B}[h]&=
{\partial_\theta}\big(\mathtt{u}_0\Lambda_{0}+\Lambda_{0}\mathtt{u}_0\big)[h]+\varepsilon\textnormal{OPS}^{-\infty},
\end{align*}
with
\begin{equation}\label{def-u0}
    \mathtt{u}_0(\varphi,\theta):=\tfrac{1}{2}   (2\mathtt{p}_{1,1})^{-\frac34} \left(\cos(\theta-b_0(\varphi))+ \varepsilon \rho_2(\varphi,\theta)\int_0^1 \sin\big(\theta-b_0(\varphi)-\varepsilon \tau \rho_2(\varphi,\theta)\big)d\tau\right).
\end{equation}
Next, we address the conjugation of $\partial_\theta \mathcal{S}$, defined in \eqref{shift-operator1}. We employ the general formula for smoothing operators (see, for instance, \cite[Lemma~2.3]{BertiMontalto}), together with the expansion of $\widehat{\beta_1}$ in \eqref{Expand-M100}:
\begin{align*}
\mathscr{B}^{-1}\partial_\theta \fint_{\mathbb{T}} K(\varphi,\theta,\eta)\,h(\varphi,\eta)\,d\eta\,\mathscr{B}
&= \partial_\theta \fint_{\mathbb{T}} K\big(\varphi,\theta+\widehat{\beta}_1(\varphi,\theta),\eta+\widehat{\beta}_1(\varphi,\eta)\big)\,h(\varphi,\eta)\,d\eta \\
&= \partial_\theta \fint_{\mathbb{T}} K\big(\varphi,\theta - b_0(\varphi),\eta - b_0(\varphi)\big)\,h(\varphi,\eta)\,d\eta + \varepsilon\,\mathrm{OPS}^{-\infty}.
\end{align*}
This can be applied to $\mathcal{S}$, which yields
\begin{align*}
\mathscr{B}^{-1}\partial_\theta \mathcal{S}\,\mathscr{B}[h](\varphi,\theta)
&= \tfrac{1}{8}(2\mathtt{p}_{1,1})^{-\frac{3}{4}} \fint_{\mathbb{T}} \Big(\cos(\theta+2\eta-3b_0)+\cos(2\theta+\eta-3b_0)\Big)\,h(\varphi,\eta)\,d\eta \\ &\quad + \varepsilon\, \partial_\theta \mathcal{S}_\infty, 
\end{align*}
with $\partial_\theta \mathcal{S}_\infty \mathrm{OPS}^{-\infty}$. 
Finally, in view of \eqref{RId-10}, one has
\begin{align*}
 \mathscr{B}^{-1}{\partial_\theta}\mathcal{R}\mathscr{B}[h] &=\partial_\theta\fint_{\mathbb{T}} \mathcal{W}_0(g)(\varphi,\theta,\eta)\ln\big|\sin\big(\tfrac{\theta-\eta}{2}\big)\big|h(\eta) d\eta+\partial_\theta\mathcal{R}_{0,\infty}, 
\end{align*}
with
$$
\mathcal{W}_0(g)(\varphi,\theta,\eta):=\mathcal{W}(g)\big(\varphi,\theta+\widehat{\beta}_1(\varphi,\theta),\eta+\widehat{\beta}_1(\varphi,\eta)\big),
$$
 and  $\partial_\theta\mathcal{R}_{0,\infty}\in\textnormal{OPS}^{-\infty}.$ 
 Notice that
$$
\mathcal{W}_0(g)(\varphi,\theta,\eta)=\mathcal{W}_0(g)(\varphi,\eta,\theta).
$$
Then, using \cite[Lemma 4.6]{HHM21} together with composition and law products
we get
$$
\mathcal{W}_0(\varphi,\theta,\eta)={\mathcal{W}}_1(\varphi,\theta)+{\mathcal{W}}_1(\varphi,\eta)+\sin^2\big(\tfrac{\theta-\eta}{2}\big){\mathcal{W}}_2(\varphi,\theta,\eta),
$$
where ${\mathcal{W}}_1$ is given by
\begin{equation}\label{defcalw1}
    {\mathcal{W}}_1(g)(\varphi,\theta):=\tfrac12\mathcal{W}_0(g)(\varphi,\theta,\theta).
\end{equation}
The function ${\mathcal{W}}_2$ is symmetric, 
$$
{\mathcal{W}}_2(g)(\varphi,\theta,\eta)={\mathcal{W}}_2(g)(\varphi,\eta,\theta),
$$
and we have the estimates 
\begin{equation}\label{est:w12}
\|{\mathcal{W}}_1(g)\|_{s}^{\textnormal{Lip},\nu}\lesssim 1+\|g\|_{s}^{\textnormal{Lip},\nu}\quad\hbox{and}\quad \|{\mathcal{W}}_2(g)\|_{s}^{\textnormal{Lip},\nu}\lesssim 1+\|g\|_{s+2+\epsilon}^{\textnormal{Lip},\nu}.
\end{equation}
Denoting 
\begin{equation}\label{def:ttu2}
   \mathtt{u}_2(g)(\varphi,\theta):=\tfrac12{\mathcal{W}}_2(g)(\varphi,\theta,\theta), 
\end{equation}
 and applying \cite[Lemma 4.6]{HHM21} once again yields
\begin{align*}
\mathcal{W}_0(\varphi,\theta,\eta)={\mathcal{W}}_1(\varphi,\theta)+{\mathcal{W}}_1(\varphi,\eta)+\sin^2\big(\tfrac{\theta-\eta}{2}\big)\big(\mathtt{u}_2(\varphi,\theta)+\mathtt{u}_2(\varphi,\eta)\big)
+\sin^4\big(\tfrac{\theta-\eta}{2}\big)\mathcal{W}_3(\varphi,\theta,\eta),
  \end{align*}
  with
$$
 \|{\mathcal{W}}_3(g)\|_{s}^{\textnormal{Lip},\nu}\lesssim 1+\|g\|_{s+{5}}^{\textnormal{Lip},\nu}.
$$
This gives the decomposition:
  \begin{align}\label{Iden-hmid-2}
 \nonumber \partial_\theta\fint_{\mathbb{T}} \mathcal{W}_0(g)(\varphi,\theta,\eta)\ln\big|\sin\big(\tfrac{\theta-\eta}{2}\big)\big|h(\eta) d\eta&={\partial_\theta}\big({\mathcal{W}}_1\Lambda_{0}+\Lambda_{0}{\mathcal{W}}_1\big)[h](\varphi,\theta)\\ \nonumber&\quad+{\partial_\theta}\big(\mathtt{u}_2\Lambda_{2}+\Lambda_{2}\mathtt{u}_2\big)[h](\varphi,\theta) +\partial_\theta\Lambda_{4,\mathcal{W}_3}[h](\varphi,\theta).
  \end{align}
In view of Lemma \ref{Lem-homo-inhomo}-(2), the operator $\partial_\theta\Lambda_{4,\mathcal{W}_3}$ belongs to  $ \textnormal{OPS}^{-4}.$ 
Putting everything together in \eqref{comp-id-com} we get
\begin{align*}
\mathscr{B}^{-1} \Big[\mathcal{T} + \tfrac12 \mathcal{H}  +\varepsilon \mathcal{H}_{\mathtt{u},0}
   + \varepsilon \partial_\theta \mathcal{S} 
   &+ \varepsilon^2 \partial_\theta \mathcal{R}(g)\Big]\mathscr{B}=\varepsilon^2|\ln\varepsilon|\omega(\lambda) \partial_\varphi+\mathtt{c}_1(\lambda,g)\partial_{\theta}+\tfrac12\mathcal{H} +\varepsilon {\partial_\theta}\mathcal{S}_0\nonumber\\ &+\varepsilon\mathcal{H}_{\mathtt{u}_1,0}+\varepsilon^2\mathcal{H}_{\mathtt{u}_2,-2}+\varepsilon^2\partial_\theta\Lambda_{4,\mathcal{W}_3}+\varepsilon^2|\ln\varepsilon|^{\frac12}\partial_\theta\mathcal{R}_{0,\infty}+\mathtt{E}_{n},
    \end{align*}
with
 \begin{equation}\label{def:ttu1}
\mathtt{u}_1(\varphi,\theta):= \mathtt{u}_0(\varphi,\theta)+\varepsilon {\mathcal{W}}_1,
  \end{equation}
  and $\mathtt{u}_0$ is given by \eqref{def-u0} and ${\mathcal{W}}_1$ is defined in \eqref{defcalw1}.
 Notice that all the involved kernels in the previous integral operators are symmetric and reversible, that is,
$$
g\in H^s_{\circ,\textnormal{even}}\Longrightarrow K(-\varphi,-\theta,-\eta)=K(\varphi,\theta,\eta)=K(\varphi,\eta,\theta).
$$
Finally, from \eqref{def:ttu1}, \eqref{def:ttu2}, \eqref{def-u0}
we conclude that
$$
\|\Delta_{12}\mathtt{u}_1\|_{s}^{\textnormal{Lip},\nu}\lesssim \|\Delta_{12}g\|_{s+1}^{\textnormal{Lip},\nu}\qquad \|\Delta_{12}\mathtt{u}_2\|_{s}^{\textnormal{Lip},\nu}\lesssim \|\Delta_{12}g\|_{s+3}^{\textnormal{Lip},\nu} .
    $$ 
This achieves the proof of Proposition \ref{prop-asym-brut}.
\end{proof}

\subsubsection{Localized change of coordinates and conjugation of the full linearized operator}\label{sec:phi1-2}
As established in the previous section, the transformation $\mathscr{B}$ acts continuously on $\textnormal{Lip}_\nu(\mathcal{O}, H^{s}(\T^{2}))$ and preserves reversibility.  
 However, it leaves invariant neither the subspace $H_\circ^{s}(\T^{2}))$ nor  $H_\perp^{s}(\T^{2}))$ defined in~\eqref{SEc-function-sapces}. 
Recall from Corollary~\ref{prop:asymp-lin-2} that the linearized operator satisfies
$$
\partial_g {\mathcal F}(\varepsilon, g) : \textnormal{Lip}_\nu(\mathcal{O}, H_\circ^{s}(\T^{2})) \;\longrightarrow\; \textnormal{Lip}_\nu(\mathcal{O}, H_\circ^{s-1}(\T^{2})).
$$
To construct a transformation compatible with this structure, we localize $\mathscr{B}$ by projecting away  in the orthogonal complement of the mode one, and define
$$
\Phi_1 := \Pi^\perp \,\mathscr{B}\, \Pi^\perp.
$$
Here, $\Pi$ denotes the projection onto the Fourier modes $\pm 1$, namely
\begin{align*}
h(\varphi,\theta)=\sum_{n\in\mathbb{Z}^\star}h_n(\varphi)e^{ i n\theta}
\quad \Longrightarrow \quad
\Pi h(\varphi,\theta)=\sum_{{|n|=1}}h_n(\varphi)e^{ i n\theta}, 
\qquad \Pi^\perp = \mathrm{Id}-\Pi.
\end{align*}
It is important to observe that $\mathscr{B}$ is symplectic, and this property ensures the preservation of the spatial zero average. More precisely,
$$
\int_{\mathbb{T}} h(\varphi, \theta)\, d\theta = 0
\quad \Longrightarrow \quad
\int_{\mathbb{T}} \mathscr{B}h(\varphi, \theta)\, d\theta = 0.
$$
Thus, if the Fourier mode $n=0$ is absent in $h$, it remains absent after the action of $\mathscr{B}$. In other words, $\mathscr{B}$ does not generate a zero mode when none is initially present.  \\
The following result, proved in~\cite{HR21}, deals with the main properties of $\Phi_1.$
\begin{pro}\label{prop-ortho}
The operator $\Phi_1:\textnormal{Lip}_\nu(\mathcal{O}, H_\perp^{s}(\T^{2}))\to \textnormal{Lip}_\nu(\mathcal{O}, H_\perp^{s}(\T^{2}))$ is an isomorphism, reversibility preserving  and satisfies the estimates, for any $s\in[s_0,s_{\textnormal{up}}]$
\begin{align*}
\|\Phi_1^{\pm 1} h\|_{s}^{\textnormal{Lip},\nu}& \lesssim   \|h\|_{s}^{\textnormal{Lip},\nu}+ \| g\|_{s+2\tau+{5}}^{\textnormal{Lip},\nu}\|h\|_{s_0}^{\textnormal{Lip},\nu}.
\end{align*} 
\end{pro}
In the following, we aim to conjugate $\mathcal{L}_{\perp}$ with $\Phi_1$. Our main result is the following.

\begin{pro}\label{prop-conju-orth}
  Let $\Phi_1$ as in Proposition $\ref{prop-ortho}$ and $\mathcal{L}_{\perp,1}$ as in Proposition $\ref{prop-asym-brut}.$    Let $n\in\N$, then on  the Cantor set $\mathcal{O}_{n}^{1}$ defined in Proposition $\ref{Second-change}$ we get

  \begin{align*}
  \mathcal{L}_{\perp,1}:= \Phi_1^{-1}\mathcal{L}_\perp\Phi_1&=\varepsilon^2|\ln\varepsilon|\omega(\lambda) \partial_\varphi+ \mathtt{c}_1(\lambda,g)\partial_{\theta}+\tfrac12\mathcal{H}+\varepsilon\Pi^\perp\mathcal{H}_{\mathtt{u}_1,0}\Pi^\perp+\varepsilon^2\Pi^\perp\mathcal{H}_{\mathtt{u}_2,-2}\Pi^\perp\\
  &\quad+\varepsilon^2 |\ln\varepsilon|^{\frac12} \partial_\theta\mathcal{R}_{1,-4}+\mathtt{E}_{n}^1,
  \end{align*}
  where the operator $\partial_\theta\mathcal{R}_{1,-4}$ is of order $-4$,
    $$
    \|\partial_\theta^5\mathcal{R}_{1,-4}[h]\|_{s}^{\textnormal{Lip},\nu}\lesssim \|h\|_{s}^{\textnormal{Lip},\nu}+\|g\|_{s+4}^{\textnormal{Lip},\nu}\|h\|_{s_0}^{\textnormal{Lip},\nu},
    $$
    and $\mathtt{E}_{n}^1$ a linear operator satisfying
\begin{equation*}
\|\mathtt{E}_{n}^1[h]\|_{s_0}^{\textnormal{Lip},\nu}\lesssim \varepsilon^{{5}} N_{0}^{\mu_{2}}N_{n+1}^{-\mu_{2}}\|h\|_{s_{0}+2}^{\textnormal{Lip},\nu}.
\end{equation*}
\end{pro}
\begin{proof}

First, note that by definition of $\mathscr{B}$, in particular due to Proposition \ref{First-change0}, we have
\begin{align}\label{beta-ort-p}
\Pi \mathscr{B}\Pi^\perp=O(\varepsilon)\quad\hbox{and}\quad \Pi^\perp \mathscr{B} \Pi=O(\varepsilon).
\end{align}
Now, take $h\in H^s_\perp (\T^2)$, in part $\Pi^\perp h=h$, and note that
\begin{align*}
 \Big[\mathcal{T} + \tfrac12 \mathcal{H}  +\varepsilon \mathcal{H}_{\mathtt{u},0}
   + \varepsilon \partial_\theta \mathcal{S} 
   + \varepsilon^2& \partial_\theta \mathcal{R}(g)\Big]\Phi_1 h=\Big[\mathcal{T} + \tfrac12 \mathcal{H}  +\varepsilon \mathcal{H}_{\mathtt{u},0}
   + \varepsilon \partial_\theta \mathcal{S} 
   + \varepsilon^2 \partial_\theta \mathcal{R}(g)\Big]\Pi^\perp \mathscr{B}h\\
   =&\Big[\mathcal{T} + \tfrac12 \mathcal{H}  +\varepsilon \mathcal{H}_{\mathtt{u},0}
   + \varepsilon \partial_\theta \mathcal{S} 
   + \varepsilon^2 \partial_\theta \mathcal{R}(g)\Big]\mathscr{B}h\\
   &-\Big[\mathcal{T} + \tfrac12 \mathcal{H}  +\varepsilon \mathcal{H}_{\mathtt{u},0}
   + \varepsilon \partial_\theta \mathcal{S} 
   + \varepsilon^2 \partial_\theta \mathcal{R}(g)\Big]\Pi \mathscr{B}h\\
   =&\mathscr{B}\mathcal{L}_{\mathscr{B}}h-\Big[\mathcal{T} + \tfrac12 \mathcal{H}  +\varepsilon \mathcal{H}_{\mathtt{u},0}
   + \varepsilon \partial_\theta \mathcal{S} 
   + \varepsilon^2 \partial_\theta \mathcal{R}(g)\Big]\Pi \mathscr{B}\Pi ^\perp h,
\end{align*}
where the last term is at least of order $\varepsilon$. Now, we apply $\Pi^\perp$:
\begin{align*}
\Pi^\perp \Big[\mathcal{T} + \tfrac12 \mathcal{H}  +\varepsilon \mathcal{H}_{\mathtt{u},0}
   + \varepsilon \partial_\theta \mathcal{S} 
   &+ \varepsilon^2 \partial_\theta \mathcal{R}(g)\Big]\Phi_1 h
   =\Pi^\perp \mathscr{B}\mathcal{L}_{\mathscr{B}}h\\ &-\Pi^\perp \Big[\mathcal{T} + \tfrac12 \mathcal{H}  +\varepsilon \mathcal{H}_{\mathtt{u},0}
   + \varepsilon \partial_\theta \mathcal{S} 
   + \varepsilon^2 \partial_\theta \mathcal{R}(g)\Big]\Pi \mathscr{B}\Pi ^\perp h,
\end{align*}
with
$$
\Pi^\perp \Big[\mathcal{T} + \tfrac12 \mathcal{H}  +\varepsilon \mathcal{H}_{\mathtt{u},0}
   + \varepsilon \partial_\theta \mathcal{S} 
   + \varepsilon^2 \partial_\theta \mathcal{R}(g)\Big]\Pi \mathscr{B}\Pi ^\perp h=O(\varepsilon^2),
$$
due to the structure of the linear operator, whose leading term is a Fourier multiplier, and \eqref{beta-ort-p}. Moreover, this latter operator  is of finite rank and therefore it belongs to $\textnormal{OPS}^{-\infty}$. Then
\begin{align*}
\Pi^\perp \Big[\mathcal{T} + \tfrac12 \mathcal{H}  +\varepsilon \mathcal{H}_{\mathtt{u},0}
   + \varepsilon \partial_\theta \mathcal{S} 
   + \varepsilon^2 \partial_\theta \mathcal{R}(g)\Big]\Phi_1 h
   =&\Pi^\perp \mathscr{B}\Pi^\perp \mathcal{L}_{\mathscr{B}}h+\Pi^\perp \mathscr{B}\Pi\mathcal{L}_{\mathscr{B}}h+O(\varepsilon^2)\\
   =&\Pi^\perp \mathscr{B}\Pi^\perp \Pi^\perp \mathcal{L}_{\mathscr{B}} \Pi^\perp h+\Pi^\perp \mathscr{B}\Pi\mathcal{L}_{\mathscr{B}}\Pi^\perp h+O(\varepsilon^2)\\
   =&\Phi_1 \Pi^\perp \mathcal{L}_{\mathscr{B}} \Pi^\perp h+O(\varepsilon^2),
\end{align*}
where,  we are using the structure of $\mathcal{L}_{\mathscr{B}}$ for the asymptotics. Now, applying $\Phi_1^{-1}$ we find
\begin{align*}
\Phi_1^{-1}\mathcal{L}_\perp \Phi_1 h&=\Phi_1^{-1}\Pi^\perp \Big[\mathcal{T} + \tfrac12 \mathcal{H}  +\varepsilon \mathcal{H}_{\mathtt{u},0}
   + \varepsilon \partial_\theta \mathcal{S} 
   + \varepsilon^2 \partial_\theta \mathcal{R}(g)\Big]\Phi_1 h\\
   &= \Pi^\perp \mathcal{L}_{\mathscr{B}} \Pi^\perp h+\varepsilon^2 \partial_\theta \mathcal{R}_{\mathscr{B},\infty},
\end{align*}
where the previous $\partial_\theta \mathcal{R}_{\mathscr{B},\infty}$ is in $\textnormal{OPS}^{-\infty}$ since it is a finite rank operator.\\
Next, let us compute $\Pi^\perp \mathcal{L}_{\mathscr{B}}\Pi^\perp.$ From the expression of $\mathcal{S}_0$  stated in \eqref{S1-defi1}, we infer 
\begin{align*}
\tfrac12\Pi^\perp\partial_\theta \mathcal{S}_0\Pi^\perp h&=0.
\end{align*}
Then, we obtain
    \begin{align*}
  \mathcal{L}_{\perp,1}:= \Phi_1^{-1}\mathcal{L}_\perp\Phi_1&=\varepsilon^2|\ln\varepsilon|\omega(\lambda) \partial_\varphi+ \mathtt{c}_1(\lambda,g)\partial_{\theta}+\tfrac12\mathcal{H}+\varepsilon\Pi^\perp\mathcal{H}_{\mathtt{u}_1,0}\Pi^\perp
+\varepsilon^2\Pi^\perp\mathcal{H}_{\mathtt{u}_2,-2}\Pi^\perp\\
  &\quad+\varepsilon^2  \Pi^\perp \partial_\theta\Lambda_{4,\mathcal{W}_3}\Pi^\perp+\varepsilon^2 |\ln\varepsilon|^{\frac12} \Pi^\perp\partial_\theta \mathcal{R}_{0,\infty}\Pi^\perp+\varepsilon^2\partial_\theta \mathcal{R}_{\mathscr{B},\infty}+\Pi^\perp \mathtt{E}_{n}\Pi^\perp\\
  &=: \varepsilon^2|\ln\varepsilon|\omega(\lambda) \partial_\varphi+ \mathtt{c}_1(\lambda,g)\partial_{\theta}+\tfrac12\mathcal{H}+\varepsilon\Pi^\perp\mathcal{H}_{\mathtt{u}_1,0}\Pi^\perp+\varepsilon^2\Pi^\perp\mathcal{H}_{\mathtt{u}_2,-2}\Pi^\perp\\
  &\quad+\varepsilon^2 |\ln\varepsilon|^\frac12 \partial_\theta\mathcal{R}_{1,-4}+\mathtt{E}_n^1,
  \end{align*}
    where the operator $\partial_\theta\mathcal{R}_{1,-4}$ is of order $-4$ and we are including there also the $\textnormal{OPS}^{-\infty}$ operators.\\
    This ends the proof of the result.
\end{proof}

\subsection{Reduction of the zero--order term: \texorpdfstring{$\Phi_2$}{Phi2}}\label{sec:Phi2}

In this step of the analysis, our goal is to simplify the structure of the operator acting on the normal direction by eliminating its zero–order contribution.  
Indeed,  the  operator $\mathcal{L}_{\perp,1}$ still contains terms of the form 
$$
\varepsilon\Pi^\perp\mathcal{H}_{\mathtt{u}_1,0}\Pi^\perp+\varepsilon^2\Pi^\perp\mathcal{H}_{\mathtt{u}_2,-2}\Pi^\perp,
$$
which carry a non-trivial zero-order component generated by singular integral operators.  
These terms prevent the operator from being in a canonical transport–diagonal form and therefore must be conjugated  into a Fourier multiplier.  
The strategy is classical: one introduces a carefully chosen transformation 
\begin{align}\label{phi-2-aid}
\Phi_2 = e^{\varepsilon \Psi_2},
\end{align}where the generator $\Psi_2$ is built from the singular integral structure involving  operators of type \eqref{Hf-n}, and takes the form 
\begin{align}\label{psi-2-aid}
\Psi_2={\Pi^\perp}\partial_\theta(\varrho_0\Lambda_{0}+\Lambda_{0} \varrho_0){\Pi^\perp}+\varepsilon {\Pi^\perp}\partial_\theta(\varrho_1\Lambda_{2}+\Lambda_{2} \varrho_1){\Pi^\perp}.
\end{align}
The role of this transformation is to transfer the zero–order part into a more regular perturbative class, while preserving the leading transport structure of the operator.  
Concretely, the functions $\varrho_0,\varrho_1 $ are determined by solving an auxiliary transport equations that balance the commutator terms arising from the conjugation.  
This reduction procedure ultimately yields a new operator $\mathcal{L}_{\perp,2}$ in which the zero–order component is replaced by an averaged term involving $\mathtt{w}_1 \mathcal{H}$, and the remaining contributions belong either to smoothing remainders or higher–order perturbations.  
\\
Our main result reads as follows.
\begin{pro}\label{prop:phi2}
   Assume \eqref{cond1} and let $g \in \mathbb{X}^{s_{\mathrm{up}}}$ be such that \eqref{small-C2-0} holds. There exists a reversibility-preserving isomorphism $\Phi_2:\mathbb{X}_\perp^s\to\mathbb{X}_\perp^s$ in the form \eqref{phi-2-aid} and \eqref{psi-2-aid}, such that 
    \begin{enumerate}
        \item The map $\Phi_2$ and its inverse $\Phi_2^{-1}$ satisfy the tame estimates
        $$
    \|\Phi_2^{\pm1}[h]\|_{s}^{\textnormal{Lip},\nu}
    \lesssim \|h\|_{s}^{\textnormal{Lip},\nu}
       + \|g\|_{s+2\tau+6}^{\textnormal{Lip},\nu}\,\|h\|_{s_0}^{\textnormal{Lip},\nu}.
$$

\item Let $n\in\N$, then for $\lambda\in\mathcal{O}_{n}^{1}$ defined in Proposition $\ref{Second-change}$, we get 
\begin{align*}
\nonumber\qquad  \mathcal{L}_{\perp,2} 
    := \Phi_2^{-1}\mathcal{L}_{\perp,1}\Phi_2&= \varepsilon^2|\ln\varepsilon|\lambda \partial_\varphi
      + \mathtt{c}_1(\lambda,g)\partial_\theta+\tfrac12\mathcal{H}
      + \varepsilon^2\mathtt{w}_1 \,\mathcal{H}+\varepsilon^2 \mathtt{w}_2\, \partial_\theta\Lambda_{2}+ \varepsilon \Pi^\perp\partial_\theta\mathcal{S}_{1}\Pi^\perp   \nonumber
       \\
   &\quad
      + \varepsilon^2|\ln\varepsilon|^{\frac12}\Pi^\perp\partial_\theta\mathcal{R}_{2,-3}\Pi^\perp
      + \mathtt{E}_n^2,
\end{align*}
with the following properties:
\begin{enumerate}
\item  The functions $(\lambda,\varphi)\mapsto\mathtt{w}_j(\lambda,\varphi), j=1,2$  $($which are constant with respect to the space variable$)$ are even in $\varphi$ and satisfy the following estimates, 
\begin{align*}
\| \mathtt{w}_1\|^{\textnormal{Lip},\nu}_s\lesssim 1+\| g\|^{\textnormal{Lip},\nu}_s, \quad 
\| \mathtt{w}_2 \|^{\textnormal{Lip},\nu}_s\lesssim 1+\| g\|^{\textnormal{Lip},\nu}_{s+3},
\end{align*}
{{and}
\begin{align*}
\| \Delta_{12}\mathtt{w}_1\|^{\textnormal{Lip},\nu}_s\lesssim \|\Delta_{12} g\|^{\textnormal{Lip},\nu}_s, \quad 
\| \Delta_{12}\mathtt{w}_2 \|^{\textnormal{Lip},\nu}_s\lesssim \| \Delta_{12}g\|^{\textnormal{Lip},\nu}_{s+3}.
\end{align*}
}
\item The functions $\varrho_0, \varrho_1$ satisfy $\langle \varrho_j\rangle_\theta=0$ and 
\begin{align*}
\|\varrho_0\|_{s}^{\textnormal{Lip},\nu}+\|\varrho_1\|_{s}^{\textnormal{Lip},\nu}
&\lesssim 1+\|g\|_{s+2\tau+6}^{\textnormal{Lip},\nu}.
\end{align*}
    \item  
 The operator $\partial_\theta\mathcal{R}_{2,-3}\in\textnormal{OPS}^{-3}$ satisfies the tame estimates
$$
  \|\partial_\theta^4\mathcal{R}_{2,-3}[h]\|_{s}^{\textnormal{Lip},\nu}
   \lesssim \|h\|_{s}^{\textnormal{Lip},\nu}
      + \|g\|_{s+2\tau+6}^{\textnormal{Lip},\nu}\|h\|_{s_0}^{\textnormal{Lip},\nu}.
$$
\item The operator  $\mathtt{E}_{n}^2$ satisfies
\begin{equation*}
\|\mathtt{E}_{n}^2[h]\|_{s_0}^{\textnormal{Lip},\nu}\lesssim \varepsilon^{{5}} N_{0}^{\mu_{2}}N_{n+1}^{-\mu_{2}}\|h\|_{s_{0}+2}^{\textnormal{Lip},\nu}+\varepsilon N_n^{s_0-s}\|  g\|_{s+2\tau+6}^{\textnormal{Lip},\nu}\|h\|_{s_0}^{\textnormal{Lip},\nu}.
\end{equation*}
\item The operator $\partial_\theta\mathcal{S}_{1}$ belongs to $\textnormal{OPS}^{-\infty}$ and is given by
$$
\partial_\theta\mathcal{S}_{1}:=\tfrac12[\mathcal{H},{\Pi^\perp}\mathcal{H}_{\varrho_0,0}{\Pi^\perp}].
$$
\end{enumerate}
    \end{enumerate}
\end{pro}
\begin{proof}

We plan to eliminate the zero-order part  of  the operator of $\mathcal{L}_{\perp,1}$.  Let $\Pi_{N_n}$ the projection over the spatial frequencies $|j|\leqslant N_n,$ then
\begin{align*}
\Pi^\perp\mathcal{H}_{\mathtt{u}_1,0}\Pi^\perp[h]
   =& \Pi^\perp\big((\Pi_{N_n}\mathtt{u}_1)\Lambda_{0}+\Lambda_{0} (\Pi_{N_n} \mathtt{u}_1)\big)\Pi^\perp[h]\\
   &+\Pi^\perp\big((\Pi^\perp_{N_n}\mathtt{u}_1)\Lambda_{0}+\Lambda_{0} (\Pi^\perp_{N_n} \mathtt{u}_1)\big)\Pi^\perp[h].
\end{align*}
We introduce the linear propagator
$$
\Phi_2 := e^{\varepsilon \Psi_2}, 
   \qquad 
   \Psi_2=\Psi_{2,0}+\varepsilon\Psi_{2,1},\qquad \Psi_{2,j}:=\partial_\theta\Pi^\perp(\varrho_j\Lambda_{2j}+\Lambda_{2j} \varrho_j)\Pi^\perp, \quad j=0,1,
$$
with $\varrho_j$ to be determined. 
A direct computation yields
$$
\partial_\varphi\Psi_2=[\partial_\varphi,\Psi_2]
   = \partial_\theta\Pi^\perp\left((\partial_\varphi \varrho_0)\Lambda_{0}+\Lambda_{0}(\partial_\varphi \varrho_0)\right)\Pi^\perp
   +\varepsilon \partial_\theta\Pi^\perp\left((\partial_\varphi \varrho_1)\Lambda_{2}+\Lambda_{2}(\partial_\varphi \varrho_1)\right)\Pi^\perp,
$$
and 
\begin{align*}
[\mathtt{c}_1\partial_\theta,\Psi_2]
   &= \mathtt{c}_1\partial_\theta\Pi^\perp\left((\partial_\theta \varrho_0)\Lambda_{0}+\Lambda_{0}(\partial_\theta \varrho_0)\right)\Pi^\perp+\varepsilon \mathtt{c}_1\partial_\theta\Pi^\perp\left((\partial_\theta \varrho_1)\Lambda_{2}+\Lambda_{2}(\partial_\theta \varrho_1)\right)\Pi^\perp.
\end{align*}
It follows that
\begin{align*}
&\varepsilon^2|\ln\varepsilon|\omega(\lambda)\partial_\varphi\Psi_2
   + [\mathtt{c}_1\partial_\theta,\Psi_2]
   + \Pi^\perp\mathcal{H}_{\mathtt{u}_1,0}\Pi^\perp+ \varepsilon \Pi^\perp\mathcal{H}_{\mathtt{u}_2,-2}\Pi^\perp
   = \partial_\theta(k_{n,0}\Lambda_{0}+\Lambda_{0} k_{n,0})\\ &\; +\varepsilon \partial_\theta(k_{n,1}\Lambda_{2}+\Lambda_{2} k_{n,1})
+\partial_\theta\big((\Pi_{N_n}^\perp \mathtt{u}_1)\Lambda_{0}+\Lambda_{0} (\Pi_{N_n}^\perp \mathtt{u}_1)
+\varepsilon(\Pi_{N_n}^\perp \mathtt{u}_2)\Lambda_{2}+\varepsilon\Lambda_{2} (\Pi_{N_n}^\perp \mathtt{u}_2)\big),
\end{align*}
where
$$
k_{n,j} := \varepsilon^2|\ln\varepsilon|\omega(\lambda)\,\partial_\varphi \varrho_j
   + \mathtt{c}_1 \partial_\theta \varrho_j + \Pi_{N_n}\mathtt{u}_{j+1}.
$$
Using the Baker–Campbell–Hausdorff expansion, we find
\begin{align*}
\Phi_2^{-1}\mathcal{L}_{\perp,1,\mathtt{p}}\Phi_2
   &= \mathcal{L}_{\perp,1,\mathtt{p}}
      + \varepsilon[\mathcal{L}_{\perp,1,\mathtt{p}},\Psi_2]
      + \tfrac12 \varepsilon^2\left[[\mathcal{L}_{\perp,1,\mathtt{p}},\Psi_2],\Psi_2\right]
      + \varepsilon^3\partial_\theta\mathcal{R}\\
      \mathcal{L}_{\perp,1,\mathtt{p}}&:=\varepsilon^2|\ln\varepsilon|\omega(\lambda) \partial_\varphi+ \mathtt{c}_1(\lambda,g)\partial_{\theta}+\tfrac12\mathcal{H}.
\end{align*}
Notice that the remainder $\partial_\theta\mathcal{R}$ can be described by Taylor expansion as an integral form with a commutator of order three. This allows to show by virtue of Lemma \ref{comm-pseudo1} that $\partial_\theta\mathcal{R}\in \textnormal{OPS}^{-3}.$
Hence, putting together the preceding identities and Proposition \ref{prop-conju-orth}, we infer
\begin{align*}
\Phi_2^{-1}\mathcal{L}_{\perp,1}\Phi_2
   &= \varepsilon^2|\ln\varepsilon|\omega(\lambda) \partial_\varphi+ \mathtt{c}_1(\lambda,g)\partial_{\theta}+\tfrac12\mathcal{H}+ \varepsilon[\tfrac12\mathcal{H},\Psi_2] +\varepsilon \partial_\theta(k_{n,0}\Lambda_{0}+\Lambda_{0} k_{n,0}) \\
  &\quad+\varepsilon^2 \partial_\theta(k_{n,1}\Lambda_{2}+\Lambda_{2} k_{n,1})+\varepsilon^2 |\ln\varepsilon|^{\frac12} \Pi^\perp\partial_\theta\tilde{\mathcal{R}}_{2,-3}\Pi^\perp+\mathtt{E}_n^2,
\end{align*}
where $\partial_\theta\tilde{\mathcal{R}}_{2,-3}$ is given by
\begin{align*}
\partial_\theta\tilde{\mathcal{R}}_{2,-3}&:=\Phi_2^{-1}\partial_\theta\mathcal{R}_{1,-4}\Phi_2 +\varepsilon^{-1}|\ln\varepsilon|^{-\frac12} \big(\Phi_2^{-1}\Pi^\perp \mathcal{H}_{\mathtt{u}_1,0} \Pi^\perp\Phi_2-\Pi^\perp \mathcal{H}_{\mathtt{u}_1,0} \Pi^\perp)\\ &\quad+|\ln\varepsilon|^{-\frac12} \big(\Phi_2^{-1}\Pi^\perp \mathcal{H}_{\mathtt{u}_2,-2} \Pi^\perp\Phi_2-\Pi^\perp \mathcal{H}_{\mathtt{u}_2,-2} \Pi^\perp) \\ &\quad 
      + \tfrac12 |\ln\varepsilon|^{-\frac12}
\left[[\mathtt{c}_1\partial_\theta+\tfrac12\mathcal{H},\Psi_2],\Psi_2\right]+ \varepsilon |\ln\varepsilon|^{-\frac12} \mathcal{R},
\end{align*}
and
$$
\mathtt{E}_n^2:=\Phi_2^{-1}\mathtt{E}_n^1\Phi_2+\varepsilon\partial_\theta\big((\Pi_{N_n}^\perp \mathtt{u}_1)\Lambda_{0}+\Lambda_{0} (\Pi_{N_n}^\perp \mathtt{u}_1)
+\varepsilon(\Pi_{N_n}^\perp \mathtt{u}_2)\Lambda_{2}+\varepsilon\Lambda_{2} (\Pi_{N_n}^\perp \mathtt{u}_2)\big).
$$
We impose that $\varrho_0$ and $\varrho_1$ solve
\begin{equation}\label{eq-mahma1}
\begin{aligned}
\varepsilon^2|\ln\varepsilon|\omega(\lambda)\,\partial_\varphi \varrho_0
   + \mathtt{c}_1 \partial_\theta \varrho_0 
   = \langle \mathtt{u}_1\rangle_{\theta}-\Pi_{N_n}\mathtt{u}_1,\\
   \varepsilon^2|\ln\varepsilon|\omega(\lambda)\,\partial_\varphi \varrho_1
   + \mathtt{c}_1 \partial_\theta \varrho_1
   = \langle \mathtt{u}_2\rangle_{\theta}-  \Pi_{N_n}\mathtt{u}_2.
   \end{aligned}
\end{equation}
In view of \eqref{def:ttu1},  \eqref{def-u0} and  \eqref{defcalw1}, one has
$$
\langle \mathtt{u}_1\rangle_{\theta}
   =\langle \mathtt{u}_0\rangle_{\theta}+ \varepsilon\langle {\mathcal{W}_1}\rangle_{\theta}=:\varepsilon\mathtt{w}_1.
$$
Moreover, from \eqref{def:ttu2}, one has
$$
\langle \mathtt{u}_2\rangle_{\theta}
   = \tfrac12\langle \mathcal{W}_2\rangle_{\theta}=:\mathtt{w}_2.
$$
By Corollary~\ref{Lemma-transp-inver}, applied with $M=2$ and  under the assumption \eqref{small-C2-0}, there exist  smooth functions $\varrho_0$ and $\varrho_1$, independent of $n$, which solve \eqref{eq-mahma1} when $\lambda\in \mathcal{O}_n^1$. This Cantor set is defined in Proposition \ref{Second-change}. Moreover, $\varrho_0$ and $\varrho_1$ satisfy the estimates
\begin{align*}
\|\varrho_j\|_{s}^{\textnormal{Lip},\nu}&\lesssim 1+\|\mathcal{W}_j\|_{s+2\tau+3}^{\textnormal{Lip},\nu}.
\end{align*}
Applying the estimates of ${\mathcal{W}}_j$  in  \eqref{est:w12}, we obtain
\begin{align}\label{rho-to-g}
\|\varrho_j\|_{s}^{\textnormal{Lip},\nu}
&\lesssim 1+\|g\|_{s+2\tau+6}^{\textnormal{Lip},\nu}.
\end{align}
Substituting back, we deduce that for any $\lambda\in \mathcal{O}_n^1$
\begin{align*}
\Phi_2^{-1}\mathcal{L}_{\perp,1}\Phi_2
   &= \varepsilon^2|\ln\varepsilon|\omega(\lambda) \partial_\varphi+ \mathtt{c}_1(\lambda,g)\partial_{\theta}+\tfrac12\mathcal{H}+\varepsilon^2 \mathtt{w}_1\mathcal{H}+\varepsilon^2 \mathtt{w}_2\partial_\theta\Lambda_{2}+ \varepsilon[\tfrac12\mathcal{H},\Psi_2] \\
  &\quad+\varepsilon^2 |\ln\varepsilon|^{\frac12} \Pi^\perp\partial_\theta\mathcal{R}_{2,-3}\Pi^\perp+\mathtt{E}_n^2.
\end{align*}
As $\Psi_2$ is an operator of zero order, then using the law products from Lemma \ref{Law-prodX1} and the definition of $\Psi_2$, we get
\begin{align*}
   \|\Psi_2[h]\|_{s}^{\textnormal{Lip},\nu}&\leqslant C \|(\varrho_0\Lambda_{0}+\Lambda_{0} \varrho_0)[h] \|_{s+1}^{\textnormal{Lip},\nu}+C\varepsilon \|(\varrho_1\Lambda_{2}+\Lambda_{2} \varrho_1)[h] \|_{s+1}^{\textnormal{Lip},\nu}\\
   &\leqslant C(\|\varrho_0\|^{\textnormal{Lip},\nu}_{s}+\varepsilon\|\varrho_1\|^{\textnormal{Lip},\nu}_{s})\|h\|^{\textnormal{Lip},\nu}_{s_0}+C(\|\varrho_0\|^{\textnormal{Lip},\nu}_{s_0}+\varepsilon\|\varrho_1\|^{\textnormal{Lip},\nu}_{s_0})\|h\|^{\textnormal{Lip},\nu}_{s}.
\end{align*}
Iterating this estimate and using $\|\rho_0\|_{s_0}^{\textnormal{Lip},\nu}+\|\rho_1\|_{s_0}^{\textnormal{Lip},\nu}\leqslant C$, we infer
\begin{align*}
   \|\Psi_2^n[h]\|_{s}^{\textnormal{Lip},\nu}&\leqslant C^{2n}\Big((\|\varrho_0\|^{\textnormal{Lip},\nu}_{s}+\varepsilon\|\varrho_1\|^{\textnormal{Lip},\nu}_{s})\|h\|^{\textnormal{Lip},\nu}_{s_0}+\|h\|^{\textnormal{Lip},\nu}_{s}\Big).
\end{align*}
Hence, by applying Taylor series and \eqref{rho-to-g}, we find
\begin{align*}
   \|\Phi_2^{\pm1}[h]\|_{s}^{\textnormal{Lip},\nu}&\leqslant \sum_{n\in\N} \tfrac{\varepsilon^n}{n!}\|\Psi_2^n[h]\|_{s}^{\textnormal{Lip},\nu}\\
   &\leqslant \Big((\|\varrho_0\|^{\textnormal{Lip},\nu}_{s}+\varepsilon\|\varrho_1\|^{\textnormal{Lip},\nu}_{s})\|h\|^{\textnormal{Lip},\nu}_{s_0}+\|h\|^{\textnormal{Lip},\nu}_{s}\Big)\sum_{n\in\N} \tfrac{C^{2n}\varepsilon^n}{n!}\\
   &\lesssim \|g\|^{\textnormal{Lip},\nu}_{s+2\tau+6}\|h\|^{\textnormal{Lip},\nu}_{s_0}+\|h\|^{\textnormal{Lip},\nu}_{s}.
\end{align*}
This shows that the two linear mappings $\Phi_2^{\pm1}$ satisfy the required  tame estimates.\\
Applying these estimates together with   Proposition \ref{prop-conju-orth}, \eqref{orthog-Lem1} and Lemma \ref{Law-prodX1}, we obtain
\begin{align*}
\|\mathtt{E}_n^2[h]\|_{s_0}^{\textnormal{Lip},\nu}&\lesssim \varepsilon^{{5}} N_{0}^{\mu_{2}}N_{n+1}^{-\mu_{2}}\|h\|_{s_{0}+2}^{\textnormal{Lip},\nu}+\|\Pi_N^\perp \mathtt{u}_1\|_{s_0+1}^{\textnormal{Lip},\nu}\|h\|_{s_0}^{\textnormal{Lip},\nu}\\
&\lesssim \varepsilon^{{5}} N_{0}^{\mu_{2}}N_{n+1}^{-\mu_{2}}\|h\|_{s_{0}+2}^{\textnormal{Lip},\nu}+N_n^{s_0-s}\| \Pi_N^\perp \mathtt{u}_1\|_{s+1}^{\textnormal{Lip},\nu}\|h\|_{s_0}^{\textnormal{Lip},\nu}.
    \end{align*}
From the expression of $\mathcal{W}_1$ we infer that $\Pi_N^\perp \mathtt{u}_1=\varepsilon\Pi_N^\perp \mathcal{W}_1$ implying in view of Proposition \ref{prop-asym-brut} that
\begin{align*}
\|\mathtt{E}_n^2[h]\|_{s_0}^{\textnormal{Lip},\nu}
&\lesssim \varepsilon^{{5}} N_{0}^{\mu_{2}}N_{n+1}^{-\mu_{2}}\|h\|_{s_{0}+2}^{\textnormal{Lip},\nu}+\varepsilon N_n^{s_0-s}\|  g\|_{s+2\tau+6}^{\textnormal{Lip},\nu}\|h\|_{s_0}^{\textnormal{Lip},\nu}.
    \end{align*}
    Recall that
    \begin{align*}
\Psi_2&={\Pi^\perp}\partial_\theta(\varrho_0\Lambda_{0}+\Lambda_{0} \varrho_0){\Pi^\perp}+\varepsilon {\Pi^\perp}\partial_\theta(\varrho_1\Lambda_{2}+\Lambda_{2} \varrho_1){\Pi^\perp}\\ & ={\Pi^\perp}\mathcal{H}_{\varrho_0,0}{\Pi^\perp}+\varepsilon {\Pi^\perp}\mathcal{H}_{\varrho_1,-2}{\Pi^\perp}.
    \end{align*}
According to  Lemma~\ref{Commutators-hilbert} and Lemma \ref{comm-pseudo1}, one obtains
\begin{align*}
  \tfrac12[\mathcal{H},\Psi_2]  & =\tfrac12[\mathcal{H},{\Pi^\perp}\mathcal{H}_{\varrho_0,0}{\Pi^\perp}]+ \tfrac\varepsilon2[\mathcal{H}, {\Pi^\perp}\mathcal{H}_{\varrho_1,-2}{\Pi^\perp}]\\
  &=: \Pi^\perp \partial_\theta\mathcal{S}_{1}\Pi^\perp+\tfrac\varepsilon2[\mathcal{H}, {\Pi^\perp}\mathcal{H}_{\varrho_1,-2}{\Pi^\perp}],
\end{align*}
with 
$$\mathcal{S}_{1}\in\textnormal{OPS}^{-\infty}\quad\hbox{and}\quad [\mathcal{H}, {\Pi^\perp}\mathcal{H}_{\varrho_1,-2}{\Pi^\perp}]\in\textnormal{OPS}^{-3}.
$$  
Similarly, one has
$$
\left[[\mathtt{c}_1\partial_\theta+\tfrac12\mathcal{H},\Psi_2],\Psi_2\right]\in\textnormal{OPS}^{-3}.
$$
Arguing as before,  we finally get
$\partial_\theta{\mathcal{R}}_{2,-3}\in\textnormal{OPS}^{-3}$ and satisfying tame estimates.
By setting
$$
\partial_\theta{\mathcal{R}}_{2,-3}:=\partial_\theta\tilde{\mathcal{R}}_{2,-3}+\tfrac12|\ln\varepsilon|^{-\frac12}[\mathcal{H}, {\Pi^\perp}\mathcal{H}_{\varrho_1,-2}{\Pi^\perp}].
$$
It follows from Lemma \ref{comm-pseudo1} that 
$$
\partial_\theta{\mathcal{R}}_{2,-3}\in \textnormal{OPS}^{-3},
$$
and we have the tame estimate, that can be checked in a straightforward way,
$$
\|\partial_\theta^4\mathcal{R}_{2,-3}[h]\|_{s}^{\textnormal{Lip},\nu}
   \;\lesssim\; \|h\|_{s}^{\textnormal{Lip},\nu}
      + \|g\|_{s+2\tau+6}^{\textnormal{Lip},\nu}\|h\|_{s_0}^{\textnormal{Lip},\nu}.
$$
The estimates of $\Delta_{12}\mathtt{w}_j$ follow from $\Delta_{12}\mathtt{u}_j$ stated in Proposition \ref{prop-asym-brut}.
This ends the proof.
\end{proof}

\subsection{Elimination of the anti-diagonal terms: \texorpdfstring{$\Phi_3$}{Phi3}}\label{sec:Phi3}

In the previous reduction step, we managed to simplify the operator by isolating its principal diagonal contributions while controlling the higher-order remainders.  
However, the transformed operator $\mathcal{L}_{\perp,2}$ still contains undesirable off–diagonal components that couple different Fourier modes in the normal directions.  
These so-called \emph{anti-diagonal terms} prevent a clean diagonalization and would obstruct the subsequent analysis of invertibility.  
To overcome this difficulty, we perform a further conjugation by means of a carefully chosen transformation $\Phi_3 = e^{\varepsilon \Psi_3}$.  
The map $\Psi_3$ is designed so that the commutator with the diagonal part of the operator cancels precisely the anti-diagonal contribution.  
In Fourier variables, this amounts to solving a homological equation that prescribes the Fourier coefficients of $\Psi_3$.  
Crucially, the diagonal coefficients vanish, which ensures the solvability of the equation without small divisor issues.

\begin{pro}\label{prop-L3-perp}
    Assume \eqref{cond1} and let $g \in \mathbb{X}^{s_{\mathrm{up}}}$ be such that \eqref{small-C2-0} holds. There exists a reversibility-preserving isomorphism $\Phi_3:\mathbb{X}_\perp^s\to\mathbb{X}_\perp^s$ such that 
    \begin{enumerate}
        \item The map $\Phi_3$ and its inverse $\Phi_3^{-1}$ satisfy the tame estimates
        $$
    \|\Phi_3^{\pm1}[h]\|_{s}^{\textnormal{Lip},\nu}
    \lesssim \|h\|_{s}^{\textnormal{Lip},\nu}
       + \|g\|_{s+2\tau+6}^{\textnormal{Lip},\nu}\,\|h\|_{s_0}^{\textnormal{Lip},\nu}.
$$

\item Let $n\in\N$, then on  the Cantor set $\mathcal{O}_{n}^{1}$ defined in Proposition \ref{Second-change} we get 
\begin{align}\label{L3-1}
 \mathcal{L}_{\perp,3}
  \;:=\; \Phi_3^{-1}\,\mathcal{L}_{\perp,2}\,\Phi_3
 \nonumber =& \varepsilon^2|\ln\varepsilon|\,\omega(\lambda)\,\partial_\varphi
     \,+\, \mathtt{c}_1\partial_\theta+\tfrac12\mathcal{H}
     \;+\; \varepsilon^2  \mathtt{w}_1\,\mathcal{H}
     \;+\; \varepsilon^2 \mathtt{w}_2\,\partial_\theta\,\Lambda_2\\
    \nonumber &+ \varepsilon^2|\ln\varepsilon|^{\frac12}\Pi^\perp\partial_\theta\mathcal{R}_{3,-3}\Pi^\perp
     \;+\; \mathtt{E}_n^3,
\end{align}
 with the following properties:
\begin{enumerate}
    \item  
 The operator $\partial_\theta\mathcal{R}_{3,-3}\in\textnormal{OPS}^{-3}$ satisfies the tame estimates
$$
  \|\partial_\theta^4\mathcal{R}_{3,-3}[h]\|_{s}^{\textnormal{Lip},\nu}
   \lesssim \|h\|_{s}^{\textnormal{Lip},\nu}
      + \|g\|_{s+2\tau+6}^{\textnormal{Lip},\nu}\|h\|_{s_0}^{\textnormal{Lip},\nu}.
$$
\item The operator  $\mathtt{E}_{n}^3$ satisfies
\begin{equation*}
\|\mathtt{E}_{n}^3[h]\|_{s_0}^{\textnormal{Lip},\nu}\lesssim \varepsilon^{{5}} N_{0}^{\mu_{2}}N_{n+1}^{-\mu_{2}}\|h\|_{s_{0}+2}^{\textnormal{Lip},\nu}+\varepsilon N_n^{s_0-s}\|  g\|_{s+2\tau+6}^{\textnormal{Lip},\nu}\|h\|_{s_0}^{\textnormal{Lip},\nu}.
\end{equation*}

\end{enumerate}
    \end{enumerate}
\end{pro}

\begin{proof}
Introduce the principal part of the operator $\mathcal{L}_{\perp,2}$:
\begin{align*}
  \mathcal{L}_{\perp,2,\mathtt{p}}
  &:= \varepsilon^2|\ln\varepsilon|\,\omega(\lambda)\,\partial_\varphi \,+ \mathtt{c}_1\,\partial_\theta\, +\tfrac12\mathcal{H}
     \,+ \varepsilon^2 \mathtt{w}_1\,\mathcal{H}\,+ \varepsilon^2 \mathtt{w}_2\, \partial_\theta\Lambda_{2}
     \,+ \varepsilon\,\Pi^\perp\partial_\theta \mathcal{S}_{1}\Pi^\perp,
\end{align*}
obtaining
\begin{equation}\label{Lperp21}
  \mathcal{L}_{\perp,2}
  \;=\;
  \mathcal{L}_{\perp,2,\mathtt{p}}
  \;+\; \varepsilon^2|\ln\varepsilon|^{\frac12}\,\Pi^\perp\partial_\theta\mathcal{R}_{2,-3}\Pi^\perp
  \;+\; \mathtt{E}_n^2.
\end{equation}
Let $\Psi_3:H_\perp^s\to H_\perp^s$ be a continuous linear map to be fixed later and set
$$
  \Phi_3:=e^{\varepsilon\Psi_3}.
$$
Using the Baker--Campbell--Hausdorff expansion we get
\begin{align*}
  \Phi_3^{-1}\,&\mathcal{L}_{\perp,2,\mathtt{p}}\,\Phi_3
  = \mathcal{L}_{\perp,2,\mathtt{p}}
     \;+\; \varepsilon[\mathcal{L}_{\perp,2,\mathtt{p}},\Psi_3]
     \;+\;  \varepsilon^2\int_0^1(1-s)e^{-\varepsilon  s\Psi_3}\big[[\mathcal{L}_{\perp,2,\mathtt{p}},\Psi_3],\Psi_3\big]e^{\varepsilon s\Psi_3}ds\\
  =& \,\varepsilon^2|\ln\varepsilon|\,\omega(\lambda)\,\partial_\varphi \,+ \mathtt{c}_1\,\partial_\theta\, +\tfrac12\mathcal{H}
     \,+ \varepsilon^2 \mathtt{w}_1\,\mathcal{H}\,+ \varepsilon^2 \mathtt{w}_2\, \partial_\theta\Lambda_{2}
     \,+ \varepsilon\,\Pi^\perp\partial_\theta \mathcal{S}_{1}\Pi^\perp\\
  +& \varepsilon^3|\ln\varepsilon|\,\omega(\lambda)\,\partial_\varphi \Psi_3
     \;+\; \varepsilon[\mathtt{c}_1\partial_\theta+\tfrac12\mathcal{H},\Psi_3]
     \;+\; \varepsilon^2\,[\Pi^\perp\partial_\theta\mathcal{S}_{1}\Pi^\perp,\Psi_3]\\
  +& \varepsilon^3[\mathtt{w}_1\,\mathcal{H},\Psi_3]\;+\; \varepsilon^3[\mathtt{w}_2\,\partial_\theta \Lambda_{2},\Psi_3]
    +\;  \varepsilon^2\int_0^1(1-s)e^{-\varepsilon  s\Psi_3}\big[[\mathcal{L}_{\perp,2,\mathtt{p}},\Psi_3],\Psi_3\big]e^{\varepsilon s\Psi_3}ds.
\end{align*}
We choose $\Psi_3$ to solve the homological equation that removes $\partial_\theta\mathcal{S}_{1}$, namely,
\begin{equation}\label{Comuta-1}
  [\mathtt{c}_1\partial_\theta+\tfrac12\mathcal{H},\Psi_3] \;=\; -\,\Pi^\perp\partial_\theta\mathcal{S}_{1}\Pi^\perp.
\end{equation}
We recall that
$$
\partial_\theta\mathcal{S}_{1}=\tfrac12[\mathcal{H},{\Pi^\perp}\mathcal{H}_{\varrho_0,0}{\Pi^\perp}].
$$
Let us use the Fourier basis $\{\mathbf{e}_n\}_{|n|\geqslant 2}$ of the normal subspace $H_\perp^s$ and write
$$
  \Psi_3\mathbf{e}_n=\sum_{|j|\geqslant 2}\psi_j^n(\varphi)\,\mathbf{e}_j,
  \qquad
  \Pi^\perp\mathcal{S}_{1}\Pi^\perp\,\mathbf{e}_n
  =\sum_{|j|\geqslant 2} S_j^n(\varphi)\,\mathbf{e}_j.
$$
Since $(\mathtt{c}_1\partial_\theta+\tfrac12\mathcal{H})\mathbf{e}_n=\lambda_n \mathbf{e}_n$ with
$$
  \lambda_n =i\big( \mathtt{c}_1\, n + \tfrac12\,\mathrm{sign}(n)\big),
$$
then we get
$$
  [\mathtt{c}_1\partial_\theta+\tfrac12\mathcal{H},\Psi_3]\mathbf{e}_n
  = \sum_{|j|\geqslant 2} (\lambda_j-\lambda_n)\,\psi_j^n(\varphi)\,\mathbf{e}_j.
$$
It follows that \eqref{Comuta-1} is equivalent to
$$
  (\lambda_j-\lambda_n)\,\psi_j^n(\varphi) \;=\; -\,S_j^n(\varphi),\qquad |j|,|n|\geqslant 2.
$$
To see that the diagonal coefficients cause no obstruction, note first that
\begin{align*}
  S_n^n
  &= \big\langle \partial_\theta\mathcal{S}_{1}\,\mathbf{e}_n,\mathbf{e}_n\big\rangle
  = i n\big\langle \mathcal{S}_{1}\,\mathbf{e}_n,\mathbf{e}_n\big\rangle.
\end{align*}
Recall that
\begin{align*}
  2\,\partial_\theta\mathcal{S}_{1}\,[\mathbf{e}_n]
  &= \mathcal{H}\mathcal{H}_{\varrho_0,0}[\mathbf{e}_n] - \mathcal{H}_{\varrho_0,0}\mathcal{H}[\mathbf{e}_n].
\end{align*}
Taking the inner product with $\mathbf{e}_n$ yields
\begin{align*}
  2\langle \partial_\theta\mathcal{S}_{1}\mathbf{e}_n,\mathbf{e}_n\rangle
  &= \langle \mathcal{H}\mathcal{H}_{\varrho_0,0}\mathbf{e}_n,\mathbf{e}_n\rangle
   - \langle \mathcal{H}_{\varrho_0,0}\mathcal{H}\mathbf{e}_n,\mathbf{e}_n\rangle.
\end{align*}
Since $\mathcal{H}^\ast=-\mathcal{H}$ and $\mathcal{H}\mathbf{e}_n=\pm i\,\mathbf{e}_n$, it follows from the definition of $\Psi_2$ that
\begin{align*}
  \langle \mathcal{H}\mathcal{H}_{\varrho_0,0}\mathbf{e}_n,\mathbf{e}_n\rangle
  &= \pm \langle \mathcal{H}_{\varrho_0,0}\mathbf{e}_n,\mathbf{e}_n\rangle\\
  &= \pm i\,n\,\big(\langle (\varrho_0\Lambda_{0}+\Lambda_{0}\varrho_0)\mathbf{e}_n,\mathbf{e}_n\rangle\\
  &= \mu_n \fint_{\mathbb{T}}\varrho_0(\varphi,\theta)\,d\theta,
\end{align*}
for some $\mu_n\in\mathbb{C}$. According to Proposition \ref{prop:phi2}, the average $\langle\varrho_0\rangle_\theta=0$ and therefore  
$$
  \langle \mathcal{H}\Psi_2\mathbf{e}_n,\mathbf{e}_n\rangle = 0.
$$
A similar computation gives
$$
  \langle \mathcal{H}_{\varrho_0,0}\mathcal{H}\mathbf{e}_n,\mathbf{e}_n\rangle = 0.
$$
It follows that,
$$
  \langle \partial_\theta\mathcal{S}_{1}\mathbf{e}_n,\mathbf{e}_n\rangle = 0,
$$
and hence
$$
  S_n^n=0,\qquad \forall\,|n|\geqslant 2.
$$
Therefore we can define, for $|j|,|n|\geqslant 2$,
\begin{equation*}
  \psi_j^n(\varphi)
  :=
  \begin{cases}
    \dfrac{S_j^n(\varphi)}{\lambda_n-\lambda_j}, & j\ne n,\\[2mm]
    0, & j=n.
  \end{cases}
\end{equation*}
Since $\partial_\theta \mathcal{S}_{1}\in \mathrm{OPS}^{-\infty}$, we obtain $\Psi_3\in \mathrm{OPS}^{-\infty}$ as well.
With this choice,
\begin{align*}
  \Phi_3^{-1}\,\mathcal{L}_{\perp,2,\mathtt{p}}\,\Phi_3
  &= \varepsilon^2|\ln\varepsilon|\,\omega(\lambda)\,\partial_\varphi
     \;+\; \mathtt{c}_1\partial_\theta+\tfrac12\mathcal{H}
     \;+\; \varepsilon^2 \mathtt{w}_1\,\mathcal{H}\;+\; \varepsilon^2 \mathtt{w}_2\,\partial_\theta\,\Lambda_2
     \;+\; \varepsilon^2\,\mathrm{OPS}^{-3},
\end{align*}
where we used $[\mathtt{c}_1\partial_\theta+\tfrac12\mathcal{H},\Psi_3]=-\Pi^\perp\partial_\theta\mathcal{S}_{1}\Pi^\perp$ to cancel the anti-diagonal part and the fact that all remaining commutators are smoothing.
Finally, conjugating \eqref{Lperp21} yields
\begin{align*}
 \nonumber \mathcal{L}_{\perp,3}
  \;:=\; \Phi_3^{-1}\,\mathcal{L}_{\perp,2}\,\Phi_3
 \nonumber =& \varepsilon^2|\ln\varepsilon|\,\omega(\lambda)\,\partial_\varphi
     \;+\; \mathtt{c}_1\partial_\theta+\tfrac12\mathcal{H}
     \;+\; \varepsilon^2 \mathtt{w}_1\,\mathcal{H}
     \;+\; \varepsilon^2 \mathtt{w}_2\,\partial_\theta\,\Lambda_2\\
    \nonumber &+\; \varepsilon^2|\ln\varepsilon|^{\frac12}\Pi^\perp\partial_\theta\mathcal{R}_{3,-3}\Pi^\perp
     \;+\; \mathtt{E}_n^3.
\end{align*}
Notice that the symmetry is not broken along this construction and one gets in particular that the functions $\mathtt{w}_j$ are even.
The operator
$\partial_\theta\mathcal{R}_{3,-3}\in \textnormal{OPS}^{-3}$  satisfies the tame estimates,
$$
\|\partial_\theta^4\mathcal{R}_{3,-3}[h]\|_{s}^{\textnormal{Lip},\nu}
   \;\lesssim\; \|h\|_{s}^{\textnormal{Lip},\nu}
      + \|g\|_{s+2\tau+6}^{\textnormal{Lip},\nu}\|h\|_{s_0}^{\textnormal{Lip},\nu}.
$$
That concludes the proof.
\end{proof}

\subsection{Elimination of the zero-order time-dependence: \texorpdfstring{$\Phi_4$}{Phi4}}\label{sec:Phi4}

In the previous step, we reduced the operator $\mathcal{L}_{\perp,2}$ to a form that is essentially diagonal, up to smoothing remainders, but still retained a zero-order coefficient $\mathtt{w}_1$ multiplying the Hilbert transform and a coefficient $\mathtt{w}_2$ multiplying the $\Lambda_2$ operator.  
The presence of these $\varphi$–dependent coefficients  prevents us from achieving the desired autonomous structure of the operator.  
To overcome this difficulty, we apply a further conjugation with a transformation of the form $\Phi_4 = e^{\Psi_4}$, where $\Psi_4= \rho(\varphi)\,\mathcal{H}$. 
As we shall see, the conjugation removes the time dependence of the coefficient in front of $\mathcal{H}$, while leaving only a constant average contribution.  
As a result, the new operator has a simplified structure in which the Hilbert transform appears with a constant prefactor.

\begin{pro}\label{prop:Phi4}
     Assume \eqref{cond1} and let $g \in \mathbb{X}^{s_{\mathrm{up}}}$ be such that \eqref{small-C2-0} holds. There exists a reversibility-preserving isomorphism $\Phi_4:\mathbb{X}_\perp^s\to\mathbb{X}_\perp^s$ such that 
    \begin{enumerate}
        \item The map $\Phi_4$ and its inverse $\Phi_4^{-1}$ satisfy the tame estimates
        $$
    \|\Phi_4^{\pm1}[h]\|_{s}^{\textnormal{Lip},\nu}
    \lesssim \|h\|_{s}^{\textnormal{Lip},\nu}
       + \|g\|_{s+2\tau+6}^{\textnormal{Lip},\nu}\,\|h\|_{s_0}^{\textnormal{Lip},\nu}.
$$

\item Let $n\in\N$, then on  the Cantor set $\mathcal{O}_{n}^{1}$ defined in Proposition \ref{Second-change} we get 
\begin{align}\label{L4-1}
 \nonumber \mathcal{L}_{\perp,4}
  \;:=\; \Phi_4^{-1}\,\mathcal{L}_{\perp,3}\,\Phi_4
 \nonumber =& \varepsilon^2|\ln\varepsilon|\,\omega(\lambda)\,\partial_\varphi
     \,+\, \mathtt{c}_1(\lambda,g)\partial_\theta+\mathtt{c}_2(\lambda,g)\,\mathcal{H}
     \\
    \nonumber &+\, \varepsilon^2 \mathtt{c}_3(\lambda,g)\,\partial_\theta\Lambda_2+ \varepsilon^2|\ln\varepsilon|^{\frac12}\Pi^\perp\partial_\theta\mathcal{R}_{4,-3}\Pi^\perp
     \;+\; \mathtt{E}_n^4,
\end{align}
with the following properties:
\begin{enumerate}
\item  The functions $\lambda\mapsto\mathtt{c}_2(\lambda,g)$ and $\lambda\mapsto\mathtt{c}_3(\lambda,g)$ $($which are constant with respect to the time-space variables$)$  satisfies the following estimates, 
\begin{align*}
\| {\mathtt{c}_2-\tfrac12 \|^{\textnormal{Lip},\nu}}\lesssim \varepsilon^2,\qquad 
\| {\mathtt{c}_3 \|^{\textnormal{Lip},\nu}}\lesssim 1.
\end{align*}
Moreover, given two functions $g_{1}$ and $g_{2}$ both satisfying \eqref{small-C2-0},  we have 
			\begin{align*}
				\|\Delta_{12}\mathtt{c}_2\|^{\textnormal{Lip},\nu}&\lesssim \varepsilon^2
                \| \Delta_{12}g\|_{{s}_{0}+3}^{\textnormal{Lip},\nu},\\ \quad \|\Delta_{12}\mathtt{c}_3\|^{\textnormal{Lip},\nu}&\lesssim \| \Delta_{12}g\|_{{s}_{0}+3}^{\textnormal{Lip},\nu}.
			\end{align*}
    \item  
 The operator $\partial_\theta\mathcal{R}_{4,-3}\in\textnormal{OPS}^{-3}$ satisfies the tame estimates
$$
\|\partial_\theta^4\mathcal{R}_{4,-3}[h]\|_{s}^{\textnormal{Lip},\nu}
   \lesssim \|h\|_{s}^{\textnormal{Lip},\nu}
      + \|g\|_{s+2\tau+6}^{\textnormal{Lip},\nu}\|h\|_{s_0}^{\textnormal{Lip},\nu}.
$$
\item The operator  $\mathtt{E}_{n}^4$ satisfies
\begin{equation*}
\|\mathtt{E}_{n}^4[h]\|_{s_0}^{\textnormal{Lip},\nu}\lesssim \varepsilon^{{5}} N_{0}^{\mu_{2}}N_{n+1}^{-\mu_{2}}\|h\|_{s_{0}+2}^{\textnormal{Lip},\nu}+\varepsilon N_n^{s_0-s}\|  g\|_{s+2\tau+6}^{\textnormal{Lip},\nu}\|h\|_{s_0}^{\textnormal{Lip},\nu}.
\end{equation*}

\end{enumerate}
    \end{enumerate}
\end{pro}

\begin{proof}
Recall from Proposition \ref{prop-L3-perp} the structure of $\mathcal{L}_{\perp, 3}$, defined by
\begin{align}\label{Lperp4p}
 \nonumber \mathcal{L}_{\perp,3}
  &= \varepsilon^2|\ln\varepsilon|\,\omega(\lambda)\,\partial_\varphi
     \;+ \mathtt{c}_1\partial_\theta+\tfrac12\mathcal{H}
     \;+ \varepsilon^2\mathtt{w}_1\,\mathcal{H}+\varepsilon^2 \mathtt{w}_2\,\partial_\theta\Lambda_2
     + \varepsilon^2|\ln\varepsilon|^{\frac12}\Pi^\perp\partial_\theta\mathcal{R}_{3,-3}\Pi^\perp
     \;+ \mathtt{E}_n^3\\
     &=:  \mathcal{L}_{\perp,3,\mathtt{p}}+\; \varepsilon^2|\ln\varepsilon|^{\frac12}\Pi^\perp\partial_\theta\mathcal{R}_{3,-3}\Pi^\perp
     \;+\; \mathtt{E}_n^3.
\end{align}
Let us consider the propagator 
$$
  \Phi_4 := e^{\Psi_4}, 
  \qquad \Psi_4 := \rho_1(\varphi)\,\mathcal{H}+ \rho_2(\varphi)\,\partial_\theta\Lambda_2,
$$
where $\rho_1$ and $\rho_2$  are scalar functions to be determined.  
Using Taylor expansion we obtain
\begin{align*}
  \Phi_4^{-1}\,\mathcal{L}_{\perp,3,\mathtt{p}}\,\Phi_4
  &= \mathcal{L}_{\perp,3,\mathtt{p}}
     + [\mathcal{L}_{\perp,3,\mathtt{p}},\Psi_4]
     + \int_0^1(1-s) e^{-s\Psi_4}\big[[\mathcal{L}_{\perp,3,\mathtt{p}},\Psi_4],\Psi_4\big]e^{s\Psi_4}\,ds.
\end{align*}
Since 
$$
  [\mathtt{w}_1\mathcal{H}+\mathtt{w}_2\,\partial_\theta\Lambda_2,\Psi_4]=0=[\mathtt{c}_1\partial_\theta+\tfrac12\mathcal{H},\Psi_4],
$$
we deduce
$$
  [\mathcal{L}_{\perp,3,\mathtt{p}},\Psi_4]
  = \varepsilon^2|\ln\varepsilon|\,\omega(\lambda)\,\rho_1'(\varphi)\,\mathcal{H}+\varepsilon^2|\ln\varepsilon|\,\omega(\lambda)\,\rho_2'(\varphi)\,\partial_\theta\Lambda_2,
  \qquad
  \big[[\mathcal{L}_{\perp,3,\mathtt{p}},\Psi_4],\Psi_4\big]=0.
$$
Therefore
\begin{align*}
  \Phi_4^{-1}\,\mathcal{L}_{\perp,3,\mathtt{p}}\,\Phi_4
  &= \varepsilon^2|\ln\varepsilon|\,\omega(\lambda)\,\partial_\varphi
     +\mathtt{c}_1\partial_\theta+\tfrac12\mathcal{H}
     + \varepsilon^2\Big(\mathtt{w}_1+|\ln\varepsilon|\,\omega(\lambda)\,\rho_1'(\varphi)\Big)\,\mathcal{H}\\ &\quad +\varepsilon^2\Big(\mathtt{w}_2+|\ln\varepsilon|\,\omega(\lambda)\,\rho_2'(\varphi)\Big)\,\partial_\theta\Lambda_2.
\end{align*}
We now choose $\rho_1$ and $\rho_2$ to eliminate the $\varphi$–dependence by imposing
\begin{align*}
  |\ln\varepsilon|\,\omega(\lambda)\,\rho_1'(\varphi) &= \langle \mathtt{w}_1\rangle_{\varphi} - \mathtt{w}_1,\quad |\ln\varepsilon|\,\omega(\lambda)\,\rho_2'(\varphi) = \langle \mathtt{w}_2\rangle_{\varphi} - \mathtt{w}_2,
\end{align*}
whose unique zero-average solution is
\begin{align*}
  \rho_1(\varphi) &= \frac{1}{|\ln\varepsilon|}\,\partial_\varphi^{-1}\Big(\langle \mathtt{w}_1\rangle_{\varphi}-\mathtt{w}_1\Big),\quad \rho_2(\varphi)= \frac{1}{|\ln\varepsilon|}\,\partial_\varphi^{-1}\Big(\langle \mathtt{w}_2\rangle_{\varphi}-\mathtt{w}_2\Big).
\end{align*}
Theses functions are  smooth, odd, and small. Using Proposition \ref{prop:phi2}, we obtain the estimate
\begin{align}\label{estim-rho-9}
 \|\rho_1\|_{s}^{\textnormal{Lip},\nu}
&\lesssim \frac{1}{|\ln\varepsilon|}\,\|\mathtt{w}_1\|_{s-1}^{\textnormal{Lip},\nu}\lesssim \frac{1}{|\ln\varepsilon|}\Big(1+\|g\|_{s}^{\textnormal{Lip},\nu}\Big).
\end{align}
Similarly, we get
\begin{align}\label{estim-rho-92}
\nonumber \|\rho_2\|_{s}^{\textnormal{Lip},\nu}
&\lesssim \frac{1}{|\ln\varepsilon|}\Big(1+\|g\|_{s+3}^{\textnormal{Lip},\nu}\Big).
\end{align}
With this choice, we find
\begin{align*}
  \Phi_4^{-1}\,\mathcal{L}_{\perp,3,\mathtt{p}}\,\Phi_4
  &= \varepsilon^2|\ln\varepsilon|\,\omega(\lambda)\,\partial_\varphi
     + \mathtt{c}_1\partial_\theta+\big(\tfrac12
     + \varepsilon^2 \langle \mathtt{w}_1\rangle_{\varphi}\big)\,\mathcal{H}+ \varepsilon^2 \langle \mathtt{w}_2\rangle_{\varphi}\,\partial_\theta\Lambda_2.
\end{align*}
Coming back to \eqref{Lperp4p}, this yields
\begin{align*}
  \mathcal{L}_{\perp,4}&:=\Phi_4^{-1}\,\mathcal{L}_{\perp,3}\,\Phi_4\\
  &= \varepsilon^2|\ln\varepsilon|\,\omega(\lambda)\,\partial_\varphi
     + \mathtt{c}_1\partial_\theta+ \mathtt{c}_2\,\mathcal{H}+ \varepsilon^2 \mathtt{c}_3\,\partial_\theta\Lambda_2+ \varepsilon^2|\ln\varepsilon|^{\frac12}\Pi^\perp\partial_\theta\mathcal{R}_{4,-3}\Pi^\perp
     + \mathtt{E}_n^4,
\end{align*}
where 
$$
\mathtt{c}_2:=\tfrac12
     + \varepsilon^2\langle \mathtt{w}_1\rangle_{\varphi}\quad \mathtt{c}_3:=\langle \mathtt{w}_2\rangle_{\varphi},
$$
and $\partial_\theta\mathcal{R}_{4,-3}\in\textnormal{OPS}^{-3}$ satisfies the tame estimate
$$
  \|\partial_\theta^4\mathcal{R}_{4,-3}[h]\|_{s}^{\textnormal{Lip},\nu}
   \lesssim \|h\|_{s}^{\textnormal{Lip},\nu}
      + \|g\|_{s+2\tau+6}^{\textnormal{Lip},\nu}\|h\|_{s_0}^{\textnormal{Lip},\nu}.
$$
The estimates of $\Delta_{12}\mathtt{c}_j$, $j=2,3$ follow from $\Delta_{12}\mathtt{w}_j$ stated in Proposition \ref{prop:phi2}-(a).

Finally, since $\rho_1$ and $\rho_2$ are odd, then the map $\Phi_4$ preserves reversibility. Moreover, by writing
$$
  {\Psi_4}=\Psi_{4,0}+\Psi_{4,1}, 
  \qquad \Psi_{4,0} := \rho_1(\varphi)\,\mathcal{H} \qquad \Psi_{4,1} :=  \rho_2(\varphi)\,\partial_\theta\Lambda_2,
$$ 
the Neumann series expansions
$$
  e^{\Psi_{4,0}} = \sum_{n\in\mathbb{N}} \frac{1}{n!}\,\rho_1^n\,\mathcal{H}^n,\qquad e^{\Psi_{4,1}} = \sum_{n\in\mathbb{N}} \frac{1}{n!}\,\rho_2^n\,(\partial_\theta\Lambda_2)^n,
$$
and the unitary structure  of the operators $\mathcal{H}$ and $\partial_\theta\Lambda_2$ on $\mathbb{X}^s$, together with Lemma \ref{Law-prodX1}, \eqref{estim-rho-9}, and the bound $\|g\|_{s_0}^{\textnormal{Lip},\nu}\leqslant 1$, we deduce
\begin{align*}
    \|e^{\Psi_4}[h]\|_{s}^{\textnormal{Lip},\nu}
    &\lesssim \|h\|_{s}^{\textnormal{Lip},\nu}
       + \|g\|_{s+3}^{\textnormal{Lip},\nu}\,\|h\|_{s_0}^{\textnormal{Lip},\nu}.
\end{align*}
Thus both $\Phi_4$ and $\Phi_4^{-1}$ satisfy the tame estimates
$$
    \|\Phi_4^{\pm1}[h]\|_{s}^{\textnormal{Lip},\nu}
    \lesssim \|h\|_{s}^{\textnormal{Lip},\nu}
       + \|g\|_{s+3}^{\textnormal{Lip},\nu}\,\|h\|_{s_0}^{\textnormal{Lip},\nu}.
$$
This concludes the proof.
\end{proof}

\section{Main result: construction of the solutions}\label{N-M-S1}

In this section, our main goal is to prove the result stated in \mbox{Theorem \ref{th-main1}}. The proof is carried out in several steps. We begin by constructing approximate solutions through a modified Nash–Moser iteration scheme, following the framework developed in \cite{Baldi-berti,BB13}. A key ingredient in this construction is Lemma  \ref{lem-matrix-opF}, which provides an approximate right inverse with tame estimates. These estimates play a crucial role in the inductive step of the scheme and ensure the convergence of the iterative procedure. Next, we analyze the convergence of the scheme and establish the existence of solutions for values of the external parameter $\lambda,$ which is related to the energy of limiting vortex filament in $(\rho,z)-$plane, belonging to a Cantor-like set. Finally, we estimate the measure of this Cantor set and show that its Lebesgue measure is asymptotically full.

\subsection{Nash-Moser scheme}\label{sec:nashmoser}
The main aim of this section is to construct solutions to the nonlinear equation 
\begin{equation}\label{main-eq1}
\mathcal{F}(g):=\mathcal{F}(\varepsilon,g) = \tfrac{1}{\varepsilon^{1+\mu}}{\bf F}_1\big(\varepsilon, f_N + \varepsilon^\mu g\big)=0,
\end{equation}
where  $f_N$ is given by Theorem~\ref{theo-approx1}, the functional ${\bf F}_1$ is introduced in \eqref{Eq-modify1} and $\mathcal{F}$ is defined in \eqref{def:mathcalF}. The parameter $\mu$ and the integer $N$ satisfy the first and the second constraint in \eqref{Assump-DRP1}, that is,
$$  N\geqslant 5, \qquad 1 \leqslant \mu < N-2.
$$
The main idea is to implement a modified Nash–Moser iteration scheme, in the spirit of \cite{Baldi-berti,BB13,BB10,BertiMontalto}. Specifically, we will construct a recursive and explicit scheme in which, at each step, we generate an approximate solution belonging to the finite-dimensional space
$$
E_{n}:=\Big\{h: \mathcal{O}\times\T^2\to\R;\,\quad \Pi_nh=h\Big\},$$
where $\mathcal{O}=(\lambda_*,\lambda^*)$  and  $\Pi_{n}$ is the projector defined   
$$
h(\varphi,\theta)=\sum_{\substack{\ell\in\Z\\ j\in\Z^\star}}h_{\ell,j}e^{ i(\ell \varphi+j\theta)},\quad\Pi_{n}h(\varphi,\theta)=\sum_{|\ell|+|j|\leqslant N_{n}}h_{\ell,j}e^{ i(\ell \varphi+j\theta)},
$$
{and  consider the sequence of numbers  $(N_n)_{n}$  defined  in \eqref{definition of Nm}.  Here we will use the parameters introduced in \eqref{cond1} and  the following additional quantity
\begin{align}\label{choice-f1}
 \quad \mathtt{b}_0=N-\mu+1.
\end{align}
The number $N_0$ appearing in \eqref{definition of Nm}, as well as the parameter 
$\nu$ involved in the Cantor sets $\mathcal{O}_{n}^1(g)$ and $\mathcal{O}_{n}^2(g)$, 
will be fixed in terms of $\varepsilon$ as specified below,
\begin{equation}\label{lambda-choice}
 \nu:= \varepsilon^2|\ln\varepsilon|^{\delta}, \quad N_{0}:= \varepsilon^{-\upsilon}, \quad \tfrac12<\delta<1, \quad 1<\tau\leqslant \tfrac32.
\end{equation}
 Moreover, we shall  impose the following constraints required along the  Nash-Moser scheme,
\begin{equation}\label{Assump-DRP1}\left\lbrace\begin{array}{rcl}
						\frac52&<&a_2,\\ {\tfrac{3}{2}} s_0+30+\tfrac32 a_2&<& a_1,\\
												\frac{2}{3}a_{1}+5&< & \mu_{2},\\
	
    0&<&\upsilon\,<\,\min\left(\tfrac{3}{\mu_2},\tfrac{N-4}{a_1}\right),\\
						15+\tfrac{9}{\upsilon}&<&\mu_1,
						
					\end{array}\right.
				\end{equation}
                and $b_1$ satisfies
                $$
            {\max\left(s_{0}+6+\tfrac{2}{3}\mu_{1}+a_{1} +{\tfrac5\upsilon},2s_0+7+\tfrac32\tfrac{\mu_2(a_2+\mu_1)}{a_2-1}\right)< b_{1}}.$$
               
                Note that these parameters can be chosen as follows.  
Given the values of $s_0, \tau,\mu$ specified in \eqref{cond1}, we determine the remaining
parameters successively in the order  
$N, a_2, a_1, \mu_2, \upsilon, \mu_1,$ and $b_1$,  
ensuring that each of them satisfies the corresponding conditions in 
\eqref{Assump-DRP1} in the prescribed order.
         This flexibility permits a wide range of admissible parameter choices, and any such selection is sufficient to obtain the results established in this section.\\
         We now turn to our central result, which concerns the implementation of a Nash–Moser scheme to construct approximate solutions to \eqref{main-eq1}. Later, we shall prove that it converges to an exact solution.
\begin{pro}[Nash-Moser scheme]\label{Nash-Moser}
Assume the conditions \eqref{cond1}, \eqref{choice-f1},\eqref{lambda-choice} and   \eqref{Assump-DRP1}. There exist $C_{\ast}>0$ and ${\varepsilon}_0>0$ such that for any $\varepsilon\in[0,\varepsilon_0]$   we get  for all $n\in\mathbb{N}$ the following properties,
\begin{enumerate}[label=\roman*)]
\item  $(\mathcal{P}1)_{n}$ There exists a Lipschitz function 
$$g_{n}:\begin{array}[t]{rcl}
\mathcal{O} & \rightarrow &  E_{n-1},\\
\lambda & \mapsto & g_n,
\end{array}$$
satisfying 
$$
g_{0}=0\quad\mbox{ and }\quad\,\| g_{n}\|_{{2s_0+2\tau+4}}^{\textnormal{Lip},\nu}\leqslant C_{\ast}\varepsilon^{\mathtt{b}_0{-1}}\nu^{-1}\quad \hbox{for}\quad n\geqslant1.
$$
By setting 
$$
\quad {u}_{n} :=g_{n}-g_{n-1} \quad \hbox{for}\quad n\geqslant1,
$$
 we have 
\begin{align*}
    \forall s\in[s_0,S], \qquad\quad \| {u}_{1}\|_{s}^{{\textnormal{Lip},\nu}}&\leqslant\tfrac12 C_{\ast}\varepsilon^{\mathtt{b}_0-1}\nu^{-1},\\
     \forall\,\, 2\leqslant k\leqslant n, \quad \|{u}_{k}\|_{{{2s_0+2\tau+4}}}^{{\textnormal{Lip},\nu}}&\leqslant C_{\ast}\varepsilon^{\mathtt{b}_0{-1}}\nu^{-1}N_{k-1}^{-a_{2}}  .
\end{align*}

\item $(\mathcal{P}2)_{n}$ Set 
$$
\mathcal{A}_{0}=\mathcal{O}\quad \mbox{ and }\quad \mathcal{A}_{n+1}=\mathcal{A}_{n}\cap\mathcal{O}_{n}^1(g_n)\cap\mathcal{O}_{n}^2(g_n) \quad\forall n\in\mathbb{N}.
$$
Then we have the following estimate 
$$\|\mathcal{F}(g_{n})\|_{s_{0},\mathcal{A}_{n}}^{{\textnormal{Lip},\nu}}\leqslant C_{\ast}\varepsilon^{\mathtt{b}_0} N_{n-1}^{-a_{1}}.
$$
\item $(\mathcal{P}3)_{n}$ High regularity estimate: $\| g_{n}\|_{b_1}^{{\textnormal{Lip},\nu}}\leqslant C_{\ast}\varepsilon^{\mathtt{b}_0{-1}}\nu^{-1} N_{n-1}^{\mu_1}.$
\end{enumerate}
\end{pro}

\begin{proof}
The proof will be implemented  using an induction principle.

\medskip\noindent $\diamond$ \textit{Initialization.}  According to Theorem \ref{theo-approx1}, \eqref{main-eq1} and \eqref{choice-f1} one has 
 \begin{equation}\label{F-zero}
\|\mathcal{F}(0)\|_{s,\mathcal{O}}^{\textnormal{Lip},\nu}\lesssim \varepsilon^{\mathtt{b}_0+1}|\ln\varepsilon|.
 \end{equation}
The properties $(\mathcal{P}1)_{0},$ $(\mathcal{P}2)_{0}$ and $(\mathcal{P}3)_{0}$ then follow immediately.

\medskip\noindent $\diamond$ {\it{Induction step:}} Given $n\in\mathbb{N}$ and  assume that  we have constructed $g_n$ satisfying the assumptions 
$(\mathcal{P}1)_{k},$ $(\mathcal{P}2)_{k}$ and $(\mathcal{P}3)_{k}$   for all $k\in\llbracket 0,n\rrbracket$ and let us check them at the next order $n+1$. As we shall explain now,    the next  approximation $g_{n+1}$ will be performed through  a  Nash-Moser scheme. First, we verify each assumption in \eqref{hao11MM} required to apply Lemma \ref{lem-matrix-opF} for the invertibility of the linearized operator.
From \eqref{lambda-choice}, we have
\begin{align*}
\varepsilon^{2}|\ln\varepsilon|^{\frac12} \nu^{-1}+\nu \varepsilon^{{-}2}|\ln\varepsilon|^{-1}=|\ln\varepsilon|^{\frac12-\delta}+|\ln\varepsilon|^{\delta-1}\leqslant {\tfrac12{\epsilon}_0},
\end{align*}
provided that $\varepsilon\in(0,\varepsilon_0)$ with $\varepsilon_0$ small enough.  Similarly, using \eqref{lambda-choice}, we obtain
\begin{align}\label{smallness-5nu}
N_{0}^{\mu_{2}}\varepsilon^{{5}} {\nu^{-1}}=\varepsilon^{3-\upsilon\mu_{2}}|\ln\varepsilon|^{-\delta}\leqslant\tfrac12{{\epsilon}_0},
\end{align}
under the condition 
$$
0<\upsilon<\tfrac{3}{\mu_2},
$$
which follows from \eqref{Assump-DRP1}.

\medskip

For the first condition in \eqref{hao11MM}, we apply Sobolev embeddings and the interpolation inequality from Lemma \ref{Law-prodX1}:
\begin{align*}
\|u_k\|_{2s_{0}+2\tau+4+\frac32\mu_2}^{\textnormal{Lip},\nu}&\leqslant  \left( \|u_k\|_{2s_{0}+2\tau+4}^{\textnormal{Lip},\nu}\right)^{1-\theta}\left( \|u_k\|_{2s_{0}+2\tau+4+\beta_1}^{\textnormal{Lip},\nu}\right)^{\theta},
\end{align*}
provided that
\begin{align*}
\tfrac32\mu_2\leqslant \beta_1\quad\hbox{and}\quad \theta:=\frac{3\mu_2}{2\beta_1},
\end{align*}
which again follows from \eqref{Assump-DRP1}. Here $ b_1$ is defined such that
\begin{align*}
b_1=2s_{0}+2\tau+4+\beta_1 .
\end{align*}
Applying  $(\mathcal{P}1)_{k}$ for $1\leqslant k\leqslant n$, we have  
\begin{align*}
2\leqslant k\leqslant n,\quad \|u_k\|_{2s_{0}+2\tau+4+\frac32\mu_2}^{\textnormal{Lip},\nu}&\leqslant C_{\ast}\varepsilon^{\mathtt{b}_0{-1}}\nu^{-1}N_{k-1}^{\theta\mu_1-(1-\theta)a_2 },
\end{align*}
and for $k=1$,
\begin{align*}
\|u_1\|_{2s_{0}+2\tau+4+\frac32\mu_2}^{\textnormal{Lip},\nu}&\leqslant \tfrac12 C_{\ast}\varepsilon^{\mathtt{b}_0-1}\nu^{-1}.
\end{align*}
Consequently, by the triangle inequality together with \eqref{definition of Nm}, \eqref{lambda-choice}, and \eqref{choice-f1}, we infer 
\begin{align*}
 \|g_n\|_{2s_{0}+2\tau+4+\frac32\mu_2}^{\textnormal{Lip},\nu}&\leqslant \sum_{k=1}^n\|u_k\|_{2s_{0}+2\tau+4+\frac32\mu_2}^{\textnormal{Lip},\nu}\leqslant  C_{\ast}\varepsilon^{\mathtt{b}_0{-1}}\nu^{-1}\sum_{k=1}^nN_{k-1}^{-1}\\&\lesssim C_{\ast}\varepsilon^{N-\mu-2}|\ln\varepsilon|^{-\delta},
\end{align*}
provided  that 
$$
{a_2\geqslant \tfrac{\theta}{1-\theta}\mu_1+\tfrac{1}{1-\theta}},
$$
which is equivalent to 
$$
{0\leqslant \theta\leqslant \frac{a_2-1}{a_2+\mu_1}\Longleftrightarrow\beta_1\geqslant \frac32\frac{\mu_2(a_2+\mu_1)}{a_2-1}}.
$$
Indeed, this is also equivalent to
$$
{b_1\geqslant 2s_0+2\tau+4+\frac32\frac{\mu_2(a_2+\mu_1)}{a_2-1}}.
$$
According to \eqref{Assump-DRP1}, we have   $\mu<N-2$ and therefore we infer from  \eqref{choice-f1} and \eqref{lambda-choice} that for sufficiently small $\varepsilon$, 
where the last inequality follows from the smallness condition we get
\begin{align}\label{small-PP1}
\varepsilon^{\mathtt{b}_0-1}\nu^{-1}
= \varepsilon^{N-2-\mu}|\ln\varepsilon|^{-\delta}
\leqslant \epsilon_0,
\end{align}
which implies
$$
 \|g_n\|_{2s_{0}+2\tau+4+\frac32\mu_2}^{\textnormal{Lip},\nu}\leqslant 1.
$$
Hence, Lemma \ref{lem-matrix-opF}-(2) applied to the operator
$$\mathcal{L}_{n}:=\,\partial_g\mathcal{F}(g_n),$$
ensures the existence of an operator $\mathcal{T}_n$ well-defined on the entire parameter set $\mathcal{O}$, such that \begin{equation}\label{estimate Tm}	
				\forall \, s\in\,[ s_0, S],\quad\|{\mathcal{T}}_{n}h\|_{s}^{\textnormal{Lip},\nu}\leqslant \frac{C}{\nu\varepsilon|\ln\varepsilon|}\Big(\|h\|_{s+2\tau}^{\textnormal{Lip},\nu}+\|g_n\|_{s+2\tau+6}^{\textnormal{Lip},\nu}\|h\|_{s_0+2\tau}^{\textnormal{Lip},\nu}\Big).
			\end{equation}
			Moreover, on the Cantor set $\mathcal{A}_{n+1}$, we have
			$$
				{\mathcal{L}_n}\,{\mathcal{T}}_{n}=\textnormal{Id}+{\mathcal{E}}_{n},
			$$
			with suitable  estimates on the remainder  ${\mathcal{E}}_{n}$.
				According  to  $(\mathcal{P}1)_{n}$,  \eqref{estimate Tm} and \eqref{small-PP1}, we deduce
				\begin{align*}
					\nonumber \|{\mathcal{T}}_{n}h\|_{s_0}^{\textnormal{Lip},\nu}&\leqslant   \frac{C}{\nu\varepsilon|\ln\varepsilon|}\|h\|_{s_0+2\tau}^{\textnormal{Lip},\nu}\Big(1+C_{\ast}\varepsilon^{\mathtt{b}_0{-1}}\nu^{-1}\Big)\leqslant   \frac{C}{\nu\varepsilon|\ln\varepsilon|}\|h\|_{s_0+2\tau}^{\textnormal{Lip},\nu}.
	\end{align*}
    We define
				\begin{align}\label{def-un}
				{g}_{n+1}:= g_{n}+{u}_{n+1}\quad\mbox{ with }\quad {u}_{n+1}:= -{\Pi}_{n}\hbox{Ext}\mathcal{T}_{n}\Pi_{n}\mathcal{F}(g_{n})\in E_{n},
				\end{align}
				where   $\textnormal{Ext} f$ denotes a Lipschitz extension of $f$ from the set $\mathcal{A}_n$ to the full interval $(\lambda_*,\lambda^*)=\mathcal{O}$, as stated in Lemma \ref{thm-extend}. In particular, for every $\lambda\in \mathcal{A}_n$ we have
                \begin{align*}
u_{n+1} = -\Pi_{n}\mathcal{T}_{n}\Pi_{n}\mathcal{F}(g_{n}), \qquad
\|u_{n+1}\|_{s,\mathcal{O}}^{{\textnormal{Lip},\nu}} \lesssim \|{u}_{n+1}\|_{s,\mathcal{A}_n}^{{\textnormal{Lip},\nu}}.
\end{align*}
				It is worth noting that the following estimate from $(\mathcal{P}1)_{n}$, 
				$$
				\forall s\in[s_0,S], \,\| {u}_{1}\|_{s}^{{\textnormal{Lip},\nu}}\leqslant\tfrac12 C_{\ast}\varepsilon^{\mathtt{b}_0-1}\nu^{-1},
				$$
				can be derived from Lemma \ref{lem-matrix-opF}-(2) and \eqref{F-zero} as follows:
				\begin{align*}
				\| {u}_{1}\|_{s}^{{\textnormal{Lip},\nu}}&\leqslant C\| \mathcal{T}_{0}\Pi_{0}\mathcal{F}(0)\|_{s}^{{\textnormal{Lip},\nu}}\leqslant \frac{C}{\nu\varepsilon|\ln\varepsilon|}\| \mathcal{F}(0)\|_{s+2\tau}^{{\textnormal{Lip},\nu}}\leqslant \tfrac{C}{|\ln\varepsilon|}\nu^{-1}\varepsilon^{\mathtt{b}_0-1}\leqslant \tfrac12 C_{\ast}\varepsilon^{\mathtt{b}_0-1}\nu^{-1},
				\end{align*}
				provided that $\varepsilon$ is small enough.
                
				\medskip
                
				\noindent $\blacktriangleright$ 
                \textbf{Verification of $(\mathcal{P}2)_{n+1}$: estimates of $\mathcal{F}({g}_{n+1})$.}

               Let us introduce the quadratic remainder
				\begin{align}\label{Def-Qm}
					Q_{n}:= \mathcal{F}(g_{n}+{u}_{n+1})-\mathcal{F}(g_{n})-\mathcal{L}_{n}{u}_{n+1}.
				\end{align}
On the set $\mathcal{A}_n$, using straightforward transformations, we obtain
				\begin{align}\label{Decom-RTT1}
					\nonumber \mathcal{F}({g}_{n+1})& =  \mathcal{F}(g_{n})-\mathcal{L}_{n}{\Pi}_{n}\mathcal{T}_{n}\Pi_{n}\mathcal{F}(g_{n})+Q_{n}\\
					\nonumber& =  \mathcal{F}(g_{n})-\mathcal{L}_{n}\mathcal{T}_{n}\Pi_{n}\mathcal{F}(g_{n})+\mathcal{L}_{n}{\Pi}_{n}^{\perp}\mathcal{T}_{n}\Pi_{n}\mathcal{F}(g_{n})+Q_{n}\\
					\nonumber& =  \mathcal{F}(g_{n})-\Pi_{n}\mathcal{L}_{n}\mathcal{T}_{n}\Pi_{n}\mathcal{F}(g_{n})+(\mathcal{L}_{n}{\Pi}_{n}^{\perp}-\Pi_{n}^{\perp}\mathcal{L}_{n})\mathcal{T}_{n}\Pi_{n}\mathcal{F}(g_{n})+Q_{n}\\
					& =  \Pi_{n}^{\perp}\mathcal{F}(g_{n})-\Pi_{n}(\mathcal{L}_{n}\mathcal{T}_{n}-\textnormal{Id})\Pi_{n}\mathcal{F}(g_{n})+(\mathcal{L}_{n}{\Pi}_{n}^{\perp}-\Pi_{n}^{\perp}\mathcal{L}_{n})\mathcal{T}_{n}\Pi_{n}\mathcal{F}(g_{n})+Q_{n}.
				\end{align}	
				The first goal is  to prove the following estimate,
				$$\|\mathcal{F}({g}_{n+1})\|_{s_{0},\mathcal{A}_{n+1}}^{\textnormal{Lip},\nu}\leqslant C_{\ast}\varepsilon^{\mathtt{b}_0} N_{n}^{-a_{1}},$$
			and thus it suffices to verify it for each term on the right-hand side of \eqref{Decom-RTT1}.

				 \medskip\noindent $\diamond$ \textit{Estimate of $\Pi_{n}^{\perp}\mathcal{F}(g_{n}).$} 
				Applying Taylor’s formula together with \eqref{F-zero}, Corollary \ref{prop:asymp-lin-2}, and $(\mathcal{P}1)_{n},$ we obtain  
				\begin{align}\label{link mathcalF(Um) and Wm}
					\forall s\geqslant s_{0},\quad\|\mathcal{F}(g_{n})\|_{s}^{\textnormal{Lip},\nu}&\leqslant\|\mathcal{F}(0)\|_{s}^{\textnormal{Lip},\nu}+\|\mathcal{F}(g_{n})-\mathcal{F}(0)\|_{s}^{\textnormal{Lip},\nu}\lesssim\varepsilon^{\mathtt{b}_0}+\| g_{n}\|_{s+2}^{\textnormal{Lip},\nu},
				\end{align}
				where we have used  the estimate (from $(\mathcal{P}1)_{n}$ and \eqref{small-PP1})
				\begin{align*}
				\nonumber\| g_{n}\|_{s_{0}+2}^{\textnormal{Lip},\nu}&\leqslant C_{\ast}\varepsilon^{\mathtt{b}_0{-1}}\nu^{-1}\leqslant \epsilon_0.
				\end{align*}
				By the standard decay properties of the projectors $\Pi_n$ and the previous estimate, it follows that
				\begin{align*}
					\|\Pi_{n}^{\perp}\mathcal{F}(g_{n})\|_{s_0}^{\textnormal{Lip},\nu}&\leqslant N_{n}^{s_{0}-b_{1}}\|\mathcal{F}(g_{n})\|_{b_1}^{\textnormal{Lip},\nu}\lesssim  N_{n}^{s_0-b_{1}}\left(\varepsilon^{\mathtt{b}_0}+\| g_{n}\|_{b_1+2}^{\textnormal{Lip},\nu}\right).
				\end{align*}
				Using $(\mathcal{P}3)_{n}$, \eqref{definition of Nm} and $g_n\in E_{n-1}$, we deduce 
				\begin{align}\label{Wm in high norm}
			\nonumber\varepsilon^{\mathtt{b}_0}+\| g_{n}\|_{b_1+2}^{\textnormal{Lip},\nu}&\leqslant\varepsilon^{\mathtt{b}_0}+C_{\ast}\varepsilon^{\mathtt{b}_0{-1}} \nu^{-1}N_{n-1}^{\mu_{1}+2}\\
\nonumber&\leqslant\varepsilon^{\mathtt{b}_0}+C_{\ast}\varepsilon^{\mathtt{b}_0{-1}} \nu^{-1} N_{n}^{\frac{2}{3}\mu_{1}+\frac{4}{3}}
                        \\
					&\leqslant 2C_{\ast}\varepsilon^{\mathtt{b}_0{-1}}\nu^{-1} N_{n}^{\frac{2}{3}\mu_{1}+2}.
				\end{align}
				From \eqref{lambda-choice}, we have
				\begin{equation}\label{fra-1}
				\nu^{-1}=N_0^{\frac{2}{\upsilon}}|\ln\varepsilon|^{-\delta}, \qquad \varepsilon^{-1}=N_0^{\frac{1}{\upsilon}}.
				\end{equation}
				Combining \eqref{link mathcalF(Um) and Wm}, \eqref{Wm in high norm}, and \eqref{fra-1}, we obtain
				\begin{align}\label{final estimate PiperpF(Um)}
					\|\Pi_{n}^{\perp}\mathcal{F}(g_{n})\|_{s_0}^{\textnormal{Lip},\nu}&\lesssim  C_{\ast}\varepsilon^{\mathtt{b}_0} N_{n}^{s_{0}+\frac{2}{3}\mu_{1}+2+\frac{{3}}{\upsilon}-b_{1}}.
				\end{align}
				Similarly, combining \eqref{link mathcalF(Um) and Wm}, \eqref{Wm in high norm}, and \eqref{fra-1}, we find
				\begin{align}\label{HDP10}
				\|\mathcal{F}(g_{n})\|_{b_{1}}^{\textnormal{Lip},\nu} &\lesssim\varepsilon^{\mathtt{b}_0}+\| g_{n}\|_{b_{1}+2}^{\textnormal{Lip},\nu} \lesssim \varepsilon^{\mathtt{b}_0{-1}}\nu^{-1}N_{n}^{\frac{2}{3}\mu_{1}+2}
               \lesssim  		C_{\ast}\varepsilon^{\mathtt{b}_0} N_{n}^{\frac{2}{3}\mu_{1}+2+\frac{{3}}{\upsilon}}.
				\end{align}

				\medskip\noindent$\diamond$ \textit{Estimate of $\Pi_{n}(\mathcal{L}_{n}\mathcal{T}_{n}-\textnormal{Id})\Pi_{n}\mathcal{F}(g_{n})$.} 	By Lemma \ref{lem-matrix-opF}-(2), in the space  $\mathcal{A}_{n+1}$ we can decompose
				$$\Pi_{n}(\mathcal{L}_{n}\mathcal{T}_{n}-\textnormal{Id})\Pi_{n}\mathcal{F}(g_{n})=\Pi_{n}\mathcal{E}_{n,1}\Pi_{n}\mathcal{F}(g_{n})+\Pi_{n}\mathcal{E}_{n,2}\Pi_{n}\mathcal{F}(g_{n}):= \mathscr{E}_{n,1}+\mathscr{E}_{n,2},
				$$	 
				where for all $s\in [s_0,S],$
				\begin{align}\label{dimanche1}
                \|\mathscr{E}_{n,1}\|_{s_0,\mathcal{A}_{n}}^{\textnormal{Lip},\nu}
			&\leqslant C N_n^{s_0-s}\nu^{-1}\Big( \|\Pi_{n}\mathcal{F}(g_{n})\|_{s+2\tau+1}^{\textnormal{Lip},\nu}+\| g_n\|_{s+4\tau+6}^{\textnormal{Lip},\nu}\|\Pi_{n}\mathcal{F}(g_{n})\|_{s_{0}+2\tau}^{\textnormal{Lip},\nu} \Big),\\
			 \|\mathscr{E}_{n,2}\|_{s_0,\mathcal{A}_{n}}^{\textnormal{Lip},\nu} & \leqslant C\varepsilon^5\nu^{-1} N_{0}^{\mu_{2}}N_{n+1}^{-\mu_{2}}\|\Pi_{n}\mathcal{F}(g_{n})\|_{s_{0}+2\tau+2}^{\textnormal{Lip},\nu}.\nonumber
			\end{align}
				Consequently,
				\begin{equation}\label{e-ai-NM}
					\|\Pi_{n}(\mathcal{L}_{n}\mathcal{T}_{n}-\textnormal{Id})\Pi_{n}\mathcal{F}(g_{n})\|_{s_0,\mathcal{A}_{n+1}}^{\textnormal{Lip},\nu}\leqslant
					\|\mathscr{E}_{n,1}\|_{s_0,\mathcal{A}_{n}}^{\textnormal{Lip},\nu}+\|\mathscr{E}_{n,2}\|_{s_0,\mathcal{A}_{n}}^{\textnormal{Lip},\nu}.
				\end{equation}
		{\bf Case $n\geqslant 1$} :
    Applying \eqref{dimanche1} with $s=b_{1}$ and using \eqref{HDP10}, $(\mathcal{P}_2)_n$, $(\mathcal{P}_3)_n$, 
        \eqref{small-PP1} and \eqref{fra-1}, we obtain
				\begin{align}\label{Esc2n}
					\|\mathscr{E}_{n,1}\|_{s_0,\mathcal{A}_{n}}^{\textnormal{Lip},\nu}&\lesssim\nu^{-1}N_n^{s_0-b_1}\left(\|\Pi_{n}\mathcal{F}(g_n)\|_{b_1+2\tau+1,\mathcal{A}_{n}}^{\textnormal{Lip},\nu}+\|g_n\|_{b_1+4\tau+6}^{\textnormal{Lip},\nu}\|\Pi_{n}\mathcal{F}(g_n)\|_{s_0+2\tau,\mathcal{A}_{n}}^{\textnormal{Lip},\nu}\right)\nonumber\\
					&\lesssim \nu^{-1}N_n^{s_0-b_1}\left(N_n^{2\tau+1}\|\mathcal{F}(g_n)\|_{b_1,\mathcal{A}_{n}}^{\textnormal{Lip},\nu}+N_{n-1}^{4\tau+6}N_n^{2\tau}\|g_n\|_{b_1}^{\textnormal{Lip},\nu}\|\mathcal{F}(g_n)\|_{s_0,\mathcal{A}_{n}}^{\textnormal{Lip},\nu}\right)\nonumber
                    \\
					&\lesssim C_{\ast}\varepsilon^{\mathtt{b}_0} \nu^{-1}N_n^{s_0-b_1}\left(N_{n}^{2\tau+\frac{2}{3}\mu_{1}+{3}+\frac{{3}}{\upsilon}}+C_{\ast}\varepsilon^{\mathtt{b}_0{-1}}\nu^{-1} N_n^{2\tau} N_{n-1}^{4\tau+6+\mu_1-a_1} \right)\nonumber
                    \\
					&\lesssim C_{\ast}\varepsilon^{\mathtt{b}_0}N_n^{s_0-b_1+\frac{{2}}{\upsilon}}\left( N_{n}^{2\tau+\frac{2}{3}\mu_{1}+3+\frac{{3}}{\upsilon}}+ N_n^{2\tau+\frac23(4\tau+6+\mu_1-a_1)} \right)\nonumber
                    \\
					&\lesssim C_{\ast}\varepsilon^{\mathtt{b}_0} N_n^{s_0+\frac{2}{3}\mu_{1}+2\tau+{3}+\frac{{5}}{\upsilon}-b_1}\left(1+ N_{n}^{\frac23(4\tau-a_1+\frac32)}\right)\nonumber\\
					&\lesssim C_{\ast}\varepsilon^{\mathtt{b}_0} N_n^{s_0+2\tau+\frac{2}{3}\mu_{1}+3+\frac{{5}}{\upsilon}-b_1},
				\end{align}
                 where we have used  $a_1\geqslant 4\tau+\tfrac32$ in the last step.
        Similarly, from \eqref{dimanche1}, $(\mathcal{P}_2)_n$ and \eqref{smallness-5nu}, we find 
				\begin{align}\label{Esc02n}
			\nonumber \|\mathscr{E}_{n,2}\|_{s_0,\mathcal{A}_{n}}^{\textnormal{Lip},\nu}&\leqslant C\varepsilon^5\nu^{-1} N_{0}^{\mu_{2}}N_{n+1}^{-\mu_{2}}\|\Pi_{n}\mathcal{F}(g_{n})\|_{s_{0}+2\tau+2,\mathcal{A}_{n}}^{\textnormal{Lip},\nu}\\
		\nonumber &\leqslant CC_\ast\varepsilon^5\nu^{-1} N_{0}^{\mu_{2}}N_{n+1}^{-\mu_{2}}N_n^{2\tau+2}\varepsilon^{\mathtt{b}_0}N_{n-1}^{-a_1}\\
		\nonumber &\leqslant C C_\ast\underbrace{\varepsilon^5\nu^{-1} N_{0}^{\mu_{2}}}_{\leqslant 1}\,\varepsilon^{\mathtt{b}_0}N_{n}^{2\tau+2-\mu_{2}-\frac23 a_1}\\
		&\leqslant C C_\ast \varepsilon^{\mathtt{b}_0} N_{n}^{2\tau+2-\mu_{2}-\frac23 a_1}.
			\end{align}
				Combining \eqref{e-ai-NM}, \eqref{Esc2n}, and \eqref{Esc02n}, we conclude that for $n\geqslant 1$,
				\begin{equation}\label{est-aait}
					\|\Pi_{n}(\mathcal{L}_{n}\mathcal{T}_{n}-\textnormal{Id})\Pi_{n}\mathcal{F}(g_{n})\|_{s_0,\mathcal{A}_{n+1}}^{\textnormal{Lip},\nu}\leqslant CC_{\ast}\varepsilon^{\mathtt{b}_0}\left(N_n^{s_0+2\tau+\frac{2}{3}\mu_{1}+3+\frac{{5}}{\upsilon}-b_1}+N_{n}^{-\mu_{2}-\frac23 a_1+2+2\tau}\right).
				\end{equation}
				{\bf Case $n=0$}: From the first line of \eqref{Esc2n} and \eqref{F-zero}, we have
				\begin{align}\label{Esc2nL}
					\|\mathscr{E}_{0,1}\|_{s_0}^{\textnormal{Lip},\nu}&\lesssim\nu^{-1}N_0^{s_0-b_1}\|\Pi_{0}\mathcal{F}(0)\|_{b_1+2\tau+1}^{\textnormal{Lip},\nu}\lesssim\nu^{-1}N_0^{s_0-b_1}\varepsilon^{\mathtt{b}_0}.				
				\end{align}
				Moreover, in view of \eqref{F-zero} and \eqref{smallness-5nu}, inequality \eqref{Esc02n} becomes
				\begin{align}\label{Esc02N}
			\nonumber \|\mathscr{E}_{0,2}\|_{s_0}^{\textnormal{Lip},\nu}&\leqslant C\varepsilon^5\nu^{-1} N_{0}^{\mu_{2}}N_{1}^{-\mu_{2}}\|\Pi_{0}\mathcal{F}(0)\|_{s_{0}+2\tau+2}^{\textnormal{Lip},\nu}\\
		\nonumber &\leqslant C C_\ast\underbrace{\varepsilon^5\nu^{-1} N_{0}^{\mu_{2}}}_{\leqslant 1}N_{1}^{-\mu_{2}}\varepsilon^{\mathtt{b}_0}\\
		&\leqslant C C_\ast \varepsilon^{\mathtt{b}_0} N_{0}^{-\frac32\mu_{2}}.
			\end{align}
		Thus, combining \eqref{Esc2nL}, \eqref{Esc02N}, and \eqref{fra-1}, we deduce
							\begin{align}\label{e-ai-0}
					\|\Pi_{0}(\mathcal{L}_{0}\mathcal{T}_{0}-\textnormal{Id})\Pi_{0}\mathcal{F}(0)\|_{s_0}^{\textnormal{Lip},\nu}&\leqslant\|\mathscr{E}_{0,1}\|_{s_0}^{\textnormal{Lip},\nu}+\|\mathscr{E}_{0,2}\|_{s_0}^{\textnormal{Lip},\nu} \lesssim C_\ast\varepsilon^{\mathtt{b}_0} N_{0}^{s_{0}+\frac2\upsilon-b_{1}}+C_\ast \varepsilon^{\mathtt{b}_0} N_{0}^{-\frac32\mu_{2}}\cdot
				\end{align}

				\medskip\noindent $\diamond$  \textit{Estimate of $\big(\mathcal{L}_{n}{\Pi}_{n}^{\perp}-\Pi_{n}^{\perp}\mathcal{L}_{n}\big)\mathcal{T}_{n}\Pi_{n}\mathcal{F}(g_{n}).$} 		
				From the structure of the operator $\mathcal{L}_{n}$ described by Corollary  \ref{prop:asymp-lin-2}, straightforward computations yield 		
				$$
				\|(\mathcal{L}_{n}{\Pi}_{n}^{\perp}-\Pi_{n}^{\perp}\mathcal{L}_{n})h\|_{s_0}^{\textnormal{Lip},\nu}\lesssim\varepsilon N_{n}^{s_{0}-b_{1}}\left(\|h\|_{b_1+1}^{\textnormal{Lip},\nu}+\|g_{n}\|_{b_1+1}^{\textnormal{Lip},\nu}\|h\|_{s_0+1}^{\textnormal{Lip},\nu}\right).
				$$
				Consequently,
				\begin{align}\label{commu-t1}
			\nonumber	\|(\mathcal{L}_{n}{\Pi}_{n}^{\perp}-\Pi_{n}^{\perp}\mathcal{L}_{n})\mathcal{T}_{n}\Pi_{n}\mathcal{F}(g_{n})\|_{s_0,\mathcal{A}_{n+1}}^{\textnormal{Lip},\nu}\lesssim&\varepsilon N_{n}^{s_{0}-b_{1}}\|\mathcal{T}_{n}\Pi_{n}\mathcal{F}(g_{n})\|_{b_1+1,\mathcal{A}_{n+1}}^{\textnormal{Lip},\nu}\\
				&+\varepsilon N_{n}^{s_{0}-b_{1}}\|g_{n}\|_{b_1+1}^{\textnormal{Lip},\nu}\|\mathcal{T}_{n}\Pi_{n}\mathcal{F}(g_{n})\|_{s_0+1,\mathcal{A}_{n+1}}^{\textnormal{Lip},\nu}.
				\end{align}
				{\bf Case $n\geqslant 1$}: Using \eqref{estimate Tm}, for all  $s\in\,[ s_0, S]$			
				\begin{equation*}
                \|\mathcal{T}_{n}\Pi_{n}\mathcal{F}(g_{n})\|_{s,\mathcal{A}_{n}}^{\textnormal{Lip},\nu}\leqslant C\nu^{-1}\varepsilon^{-1}\Big(\|\Pi_{n}\mathcal{F}(g_{n})\|_{s+2\tau,\mathcal{A}_{n}}^{\textnormal{Lip},\nu}+\|g_n\|_{s+2\tau+6}^{\textnormal{Lip},\nu}\|\Pi_{n}\mathcal{F}(g_{n})\|_{s_0+2\tau,\mathcal{A}_{n}}^{\textnormal{Lip},\nu}\Big).
			\end{equation*}
			Thus,  for $s=s_0+1$ one has
			\begin{align*}
				\varepsilon\|\mathcal{T}_{n}\Pi_{n}\mathcal{F}(g_{n})\|_{s_0+1,\mathcal{A}_{n}}^{\textnormal{Lip},\nu}&\leqslant C\nu^{-1}N_n^{2\tau}\Big(N_n+\|g_n\|_{s_0+2\tau+7}^{\textnormal{Lip},\nu}\Big)\|\mathcal{F}(g_{n})\|_{s_0,\mathcal{A}_{n}}^{\textnormal{Lip},\nu}
                \\
	&\leqslant CC_*\nu^{-1}\varepsilon^{\mathtt{b}_0}N_n^{2\tau+1-\frac23 a_1}\Big(1 +\underbrace{C_* \varepsilon^{\mathtt{b}_0{-1}}\nu^{-1}}_{\leqslant   1}\Big)\\
	&\leqslant CC_*\nu^{-1}\varepsilon^{\mathtt{b}_0}N_n^{2\tau+1-\frac23 a_1}.			
			\end{align*}
			For $s=b_1+1$ we find, by virtue of \eqref{HDP10},  $(\mathcal{P}2)_{n}$,  $(\mathcal{P}3)_{n}$ and \eqref{small-PP1},
			\begin{align*}
				\varepsilon\|\mathcal{T}_{n}\Pi_{n}\mathcal{F}(g_{n})\|_{b_1+1,\mathcal{A}_{n}}^{\textnormal{Lip},\nu}&\leqslant C\nu^{-1} N_n^{2\tau}\Big(N_n\|\Pi_{n}\mathcal{F}(g_{n})\|_{b_1,\mathcal{A}_{n}}^{\textnormal{Lip},\nu}+\|g_n\|_{b_1+2\tau+7}^{\textnormal{Lip},\nu}\|\Pi_{n}\mathcal{F}(g_{n})\|_{s_0,\mathcal{A}_{n}}^{\textnormal{Lip},\nu}\Big)\\
			&\leqslant CC_*\nu^{-1}\varepsilon^{\mathtt{b}_0}N_n^{2\tau}\Big(N_{n}^{\frac{2}{3}\mu_{1}+{3}+\tfrac{{3}}{\upsilon}}+{C_*\nu^{-1}\varepsilon^{\mathtt{b}_0{-1}}} N_{n-1}^{2\tau+7+\mu_1-a_1} \Big)\\
            &\leqslant CC_*\nu^{-1}\varepsilon^{\mathtt{b}_0}N_n^{2\tau}\Big(N_{n}^{\frac{2}{3}\mu_{1}+3+\tfrac{{3}}{\upsilon}}+\underbrace{C_*\nu^{-1}\varepsilon^{\mathtt{b}_0{-1}}}_{\leqslant 1} N_{n}^{2+\frac23\mu_1} N_{n-1}^{2\tau+4-a_1}  \Big)\\
				&\leqslant CC_*\nu^{-1}\varepsilon^{\mathtt{b}_0}N_{n}^{2\tau+3+\frac{2}{3}\mu_{1}+\frac{{3}}{\upsilon}},
			\end{align*}
            where we used $a_1\geqslant 2\tau+4$.
			Plugging these bounds into \eqref{commu-t1} and using \eqref{small-PP1}, we obtain
			\begin{align}\label{final estimate commutator}
			\nonumber 	\|(\mathcal{L}_{n}{\Pi}_{n}^{\perp}-\Pi_{n}^{\perp}\mathcal{L}_{n})\mathcal{T}_{n}\Pi_{n}\mathcal{F}(g_{n})\|_{s_0,\mathcal{A}_{n+1}}^{\textnormal{Lip},\nu}&\lesssim \varepsilon^{\mathtt{b}_0}\nu^{-1} N_{n}^{s_{0}-b_{1}}\Big(N_{n}^{2\tau+3+\frac{2}{3}\mu_{1}+\frac{{3}}{\upsilon}}\\
  &\qquad\qquad\qquad \qquad  \nonumber+\varepsilon^{\mathtt{b}_0-1}\nu^{-1} N_n^{2+2\tau+\frac23\mu_1-\frac23 a_1}\Big)\\
				\nonumber&\lesssim \varepsilon^{\mathtt{b}_0}\nu^{-1}(1+\varepsilon^{\mathtt{b}_0-1} \nu^{-1})N_{n}^{s_{0}-b_{1}+2\tau+3+\frac{2}{3}\mu_{1}+\frac{{3}}{\upsilon}}\\
				&\lesssim \varepsilon^{\mathtt{b}_0} N_{n}^{s_{0}-b_{1}+2\tau+3+\frac{2}{3}\mu_{1}+\frac{{5}}{\upsilon}}.
				\end{align}
				{\bf Case $n=0$}: From \eqref{commu-t1}, \eqref{estimate Tm}, \eqref{F-zero}, and \eqref{fra-1}, we deduce
								\begin{align}\label{Thm-1}
	\nonumber			\|(\mathcal{L}_{0}{\Pi}_{0}^{\perp}-\Pi_{0}^{\perp}\mathcal{L}_{0})\mathcal{T}_{0}\Pi_{0}\mathcal{F}(g_{0})\|_{s_0}^{\textnormal{Lip},\nu}&\lesssim\varepsilon N_{0}^{s_{0}-b_{1}}\|\mathcal{T}_{0}\Pi_{0}\mathcal{F}(0)\|_{b_1+1}^{\textnormal{Lip},\nu}\\
\nonumber&\lesssim \nu^{-1}  N_{0}^{s_{0}-b_{1}}	\|\mathcal{F}(0)\|_{b_1+2\tau+1}^{\textnormal{Lip},\nu}\\
&\lesssim C_*\varepsilon^{\mathtt{b}_0}N_{0}^{s_{0}-b_{1}+\frac2\upsilon}			.
				\end{align}

				\medskip\noindent $\diamond$  \textit{Estimate of $Q_{n}$.} Applying Taylor’s formula together with \eqref{Def-Qm}, we have
				$$Q_{n}=\int_{0}^{1}(1-t)d_{g}^{2}\mathcal{F}(g_{n}+t{u}_{n+1})[{u}_{n+1},{u}_{n+1}]dt.$$
				From Corollary \ref{prop:asymp-lin-2}, it follows that
				\begin{align}\label{mahma-YDa1}
					\| Q_{n}\|_{s_0,\mathcal{A}_{n}}^{\textnormal{Lip},\nu}\lesssim \varepsilon^{\mu+1}\left(1+\varepsilon\|g_{n}\|_{s_0+2}^{\textnormal{Lip},\nu}+\varepsilon\| {u}_{n+1}\|_{s_0+2,\mathcal{A}_{n}}^{\textnormal{Lip},\nu}\right)\left(\| {u}_{n+1}\|_{s_0+2,\mathcal{A}_{n}}^{\textnormal{Lip},\nu}\right)^{2}.
				\end{align}
Combining \eqref{def-un} and \eqref{estimate Tm}, for all
				 $s\in[s_{0},S]$ 
				\begin{align}\label{Tu-L1}
					\| {u}_{n+1}\|_{s,\mathcal{A}_{n}}^{\textnormal{Lip},\nu} & =  \|{\Pi}_{n}\mathcal{T}_{n}\Pi_{n}\mathcal{F}(g_{n})\|_{s,\mathcal{A}_{n}}^{\textnormal{Lip},\nu}\nonumber\\
					& \lesssim  \nu^{-1}\varepsilon^{-1}|\ln\varepsilon|^{-1}\left(\|\Pi_{n}\mathcal{F}(g_{n})\|_{s+2\tau,\mathcal{A}_{n}}^{\textnormal{Lip},\nu}+\|g_{n}\|_{s+2\tau+6}^{\textnormal{Lip},\nu}\|\Pi_{n}\mathcal{F}(g_{n})\|_{s_0+2\tau,\mathcal{A}_{n}}^{\textnormal{Lip},\nu}\right)\nonumber\\
					& \lesssim  \nu^{-1}\varepsilon^{-1}N_{n}^{2\tau}\left(\|\Pi_{n}\mathcal{F}(g_{n})\|_{s,\mathcal{A}_{n}}^{\textnormal{Lip},\nu}+\| g_{n}\|_{s+2\tau+6}^{\textnormal{Lip},\nu}\|\Pi_{n}\mathcal{F}(g_{n})\|_{s_0,\mathcal{A}_{n}}^{\textnormal{Lip},\nu}\right).
				\end{align}
                {\bf Case $n\geqslant 1$}: Taking $s={s_0+2}$ in\eqref{Tu-L1} and using \eqref{small-PP1},   $(\mathcal{P}1)_{n}$ and $(\mathcal{P}2)_{n}$, 
                \begin{align}\label{Tu-L1ss02}
					\varepsilon\| {u}_{n+1}\|_{s_0+2,\mathcal{A}_{n}}^{\textnormal{Lip},\nu} 
				\nonumber & \lesssim  \nu^{-1}N_{n}^{2\tau}\left(N_{n}^{2}+{\| g_{n}\|_{s_0+2\tau+8}^{\textnormal{Lip},\nu}}\right)\|\Pi_{n}\mathcal{F}(g_{n})\|_{s_0,\mathcal{A}_{n}}^{\textnormal{Lip},\nu}
                \nonumber \\& \lesssim  \nu^{-1}N_{n}^{2\tau}\left(N_{n}^{2}+N_{n}^{2}\| g_{n}\|_{2s_0+2\tau+4}^{\textnormal{Lip},\nu}\right)\|\Pi_{n}\mathcal{F}(g_{n})\|_{s_0,\mathcal{A}_{n}}^{\textnormal{Lip},\nu}\\ \nonumber	& \lesssim  \nu^{-1}N_{n}^{2\tau+2}\left(1+C_{\ast}\varepsilon^{\mathtt{b}_0{-1}}\nu^{-1}\right)C_{\ast}\varepsilon^{\mathtt{b}_0}N_{n-1}^{-a_1}
                \nonumber\\
					& \nonumber\lesssim C_{\ast}\varepsilon^{\mathtt{b}_0}\nu^{-1}  N_{n}^{2\tau+2-\frac23 a_1}.
				\end{align}
                Choosing $\varepsilon$ sufficiently small and using  \eqref{small-PP1} with $
				{a_1\geqslant 3\tau+3}
				$, we ensure
								\begin{align*}
					\varepsilon\| {u}_{n+1}\|_{s_0+2,\mathcal{A}_{n}}^{\textnormal{Lip},\nu} 
					& \leqslant 2,
				\end{align*}
 Plugging \eqref{Tu-L1ss02} into \eqref{mahma-YDa1}, for all  $n\geqslant 1$
                \begin{align}\label{final estimate for Qm}
				 \nonumber\| Q_{n}\|_{s_0,\mathcal{A}_{n}}^{\textnormal{Lip},\nu}&\lesssim \varepsilon^{\mu-1}\left(\varepsilon\| {u}_{n+1}\|_{s_0+2,\mathcal{A}_{n}}^{\textnormal{Lip},\nu}\right)^{2}\\ &\nonumber\lesssim C_{\ast}^2\varepsilon^{2\mathtt{b}_0+\mu-1}\nu^{-2}  N_{n}^{4\tau+4-\frac43 a_1}\\ &\lesssim C_{\ast}^2\varepsilon^{\mathtt{b}_0}   N_{n}^{4\tau+4-\frac43 a_1},
				\end{align}
                provided that
                \begin{align*}
	 \varepsilon^{\mathtt{b}_0+\mu-1}\nu^{-2}&=\varepsilon^{N-4}|\ln\varepsilon|^{-2\delta}\leqslant \epsilon_0,
\end{align*}
 which follows from \eqref{choice-f1}, \eqref{lambda-choice}, and \eqref{Assump-DRP1}. \\     
 {\bf Case $n=0$}:
				From \eqref{Tu-L1} and \eqref{F-zero},  for all $s\in[s_{0},S]$,
				\begin{align}\label{H1}
					\|{u}_{1}\|_{s}^{\textnormal{Lip},\nu}&\lesssim\nu^{-1}\varepsilon^{-1}|\ln\varepsilon|^{-1} \|\Pi_{0}\mathcal{F}(0)\|_{s+2\tau}^{\textnormal{Lip},\nu}\leqslant C_{\ast}\nu^{-1}\varepsilon^{\mathtt{b}_0-1}.
				\end{align}
				Thus, by \eqref{mahma-YDa1} and \eqref{H1},
								\begin{align}\label{e-Q0}
				 \|Q_{0}\|_{s_{0}}^{\textnormal{Lip},\nu}&\lesssim C_{\ast} \varepsilon^{\mu-1} (\nu^{-1}\varepsilon^{\mathtt{b}_0})^2\lesssim C_{\ast}\varepsilon^{2\mathtt{b}_0+\mu-1}\nu^{-2}.
				\end{align}
				 \textit{Conclusion.}
				Inserting  \eqref{final estimate PiperpF(Um)}, \eqref{est-aait}, \eqref{final estimate commutator} and \eqref{final estimate for Qm}, into \eqref{Decom-RTT1} yields, for $n\in\mathbb{N}^{*}$,
				\begin{align*}
					\|\mathcal{F}({g}_{n+1})\|_{s_{0},\mathcal{A}_{n+1}}^{\textnormal{Lip},\nu}&\leqslant CC_{\ast}\varepsilon^{\mathtt{b}_0}\left(N_{n}^{s_{0}+\frac{2}{3}\mu_{1}+2+\frac{{3}}{\upsilon}-b_{1}}+N_n^{s_0+2\tau+\frac{2}{3}\mu_{1}+3+\frac{{5}}{\upsilon}-b_1}+N_{n}^{-\mu_{2}-\frac23 a_1+2+2\tau}\right)\\
					&\qquad +C C_\ast\varepsilon^{\mathtt{b}_0} N_{n}^{s_{0}-b_{1}+2\tau+\frac{2}{3}\mu_{1}+3+\frac{{5}}{\upsilon}}+C C_{\ast}\varepsilon^{\mathtt{b}_0} N_{n}^{4\tau+4-\frac43a_1}.
				\end{align*} 
				The parameter conditions in \eqref{Assump-DRP1} ensure
				\begin{equation}\label{Assump-DR1}\left\lbrace\begin{array}{rcl}
						s_{0}+2\tau+3+\frac{2}{3}\mu_{1}+\frac{{5}}{\upsilon}+a_1&< & b_{1},\\
						\frac{1}{3}a_{1}+2\tau+2 &< & \mu_{2},\\
						4\tau+4&< & \frac{1}{3}a_{1}.
					\end{array}\right.
				\end{equation}
				Thus, by choosing $N_{0}$ sufficiently large (equivalently, $\varepsilon$ sufficiently small), we obtain for all $n\in\mathbb{N}^*,$
				\begin{align*}
					\|\mathcal{F}({g}_{n+1})\|_{s_{0},\mathcal{A}_{n+1}}^{\textnormal{Lip},\nu}\leqslant C_{\ast}\varepsilon^{\mathtt{b}_0} N_{n}^{-a_{1}}.
				\end{align*}
				For  $n=0$, we plug \eqref{final estimate PiperpF(Um)}, \eqref{e-ai-0},  \eqref{Thm-1} and \eqref{e-Q0} into \eqref{Decom-RTT1} and use \eqref{lambda-choice} and \eqref{small-PP1}   to get
				\begin{align*}\|\mathcal{F}({g}_{1})\|_{s_{0}}^{\textnormal{Lip},\nu}&\leqslant CC_{\ast}\varepsilon^{\mathtt{b}_0}\left(N_{0}^{s_{0}+\frac{2}{3}\mu_{1}+2+\frac3\upsilon-b_{1}} +N_{0}^{s_{0}+\frac2\upsilon-b_{1}}+ N_{0}^{-\frac32\mu_{2}}+ N_{0}^{s_{0}+\frac2\upsilon-b_{1}}	+\varepsilon^{\mathtt{b}_0+\mu-1}\nu^{-2}\right)\\
				&\leqslant CC_{\ast}\varepsilon^{\mathtt{b}_0}\left(N_{0}^{s_{0}+\frac{2}{3}\mu_{1}+2+\frac3\upsilon-b_{1}}+  N_{0}^{-\frac32\mu_{2}}+N_0^{-\frac{{\mathtt{b}_0+\mu}-5}{\upsilon}}\right).
				\end{align*}
To guarantee this decay, we impose
\begin{equation}\label{Assump-DRL1}\left\lbrace\begin{array}{rcl}
						s_{0}+\frac23\mu_1+2+\frac3\upsilon+a_1&< & b_{1},\\
						a_1&<&\frac32\mu_2,\\
						a_1&< & \frac{{\mathtt{b}_0+\mu}-5}{\upsilon},
					\end{array}\right.
				\end{equation}
                which, for $\varepsilon$ small enough, gives
				\begin{align*}\|\mathcal{F}({g}_{1})\|_{s_{0}}^{\textnormal{Lip},\nu}&\leqslant C_{\ast}\varepsilon^{\mathtt{b}_0} N_0^{-a_1}.
				\end{align*}			
				Note that the first condition in \eqref{Assump-DRL1} follows from the first in \eqref{Assump-DR1}, while the third is equivalent to 
				\begin{align*}0<\upsilon< \tfrac{\mathtt{b}_0+\mu-5}{a_1}=\tfrac{N-4}{a_1}\cdot
				\end{align*}
				The second condition is immediate from the second in \eqref{Assump-DRP1}.
				This completes the proof of the estimates in $(\mathcal{P}2)_{n+1}.$ 

                \medskip
				
				\noindent $\blacktriangleright$ \textbf{Verification of $(\mathcal{P}1)_{n+1}$ and $(\mathcal{P}3)_{n+1}.$} 
Using  \eqref{Tu-L1}, $(\mathcal{P}3)_{n}$ and \eqref{small-PP1}, \eqref{HDP10} we deduce that
\begin{align}\label{Tu-LT1}
					\nonumber\| {u}_{n+1}\|_{b_1}^{\textnormal{Lip},\nu} 
					& \leqslant C\nu^{-1}\varepsilon^{-1}N_{n}^{2\tau}\left(\|\Pi_{n}\mathcal{F}(g_{n})\|_{b_1,\mathcal{A}_{n}}^{\textnormal{Lip},\nu}+C_{\ast}\varepsilon^{\mathtt{b}_0}N_{n-1}^{-a_1}\| g_{n}\|_{b_1+2\tau+6}^{\textnormal{Lip},\nu}\right)\\ & \nonumber \leqslant   C C_{\ast}\varepsilon^{\mathtt{b}_0-1} \nu^{-1}N_{n}^{2\tau}\left(N_{n}^{\frac23\mu_1+2+\frac3\upsilon}+C_{\ast}\varepsilon^{\mathtt{b}_0{-1}}\nu^{-1} N_{n}^{\frac43\tau+{4}+\frac23\mu_1-\frac23a_1}\right)\\
					&\leqslant C C_*\varepsilon^{\mathtt{b}_0-1} \nu^{-1}N_{n}^{{2\tau}+2+\frac23\mu_1+\frac{{3}}{\upsilon}}.
				\end{align}
Now gathering \eqref{Tu-LT1} and $(\mathcal{P}3)_{n}$ allows to write
				\begin{align*}
					\|g_{n+1}\|_{b_1}^{\textnormal{Lip},\nu}  & \leqslant   \| g_{n}\|_{b_1}^{\textnormal{Lip},\nu}+ \|u_{n+1}\|_{b_1}^{\textnormal{Lip},\nu} \\
					& \leqslant  C_{\ast} \varepsilon^{\mathtt{b}_0{-1}}\nu^{-1}N_{n}^{\frac23\mu_{1}}+CC_\ast\varepsilon^{\mathtt{b}_0{-1}}\nu^{-1}N_{n}^{2\tau}+2+\frac23\mu_1+\frac3\upsilon\\
					& \leqslant  C_{\ast} \varepsilon^{\mathtt{b}_0{-1}}\nu^{-1}N_{n}^{\mu_{1}},
				\end{align*}
				provided that
\begin{align*}
{2\tau}+2+\tfrac{3}{\upsilon}<\tfrac13\mu_1,
\end{align*}
which follows from \eqref{Assump-DRP1}.		
This achieves $(\mathcal{P}3)_{n+1}$. \\
In view of \eqref{Tu-L1ss02}, we get
	\begin{align*}
				\nonumber\| {u}_{n+1}\|_{{2s_0+2\tau+4},\mathcal{A}_{n}}^{\textnormal{Lip},\nu}& \lesssim N_{n}^{s_0+2\tau+2}\|{u}_{n+1}\|_{s_0+2,\mathcal{A}_{n}}^{\textnormal{Lip},\nu}\lesssim C_{\ast}\varepsilon^{\mathtt{b}_0{-1}}\nu^{-1} N_{n}^{{s_0+4\tau+4}-\frac23 a_1}\\
                &
                  \lesssim C_{\ast}\varepsilon^{\mathtt{b}_0{-1}}\nu^{-1} N_{n}^{-a_2},
				\end{align*}
				provided that
				\begin{align*}
				{s_0+4\tau+4+a_2}<\tfrac23 a_1,
				\end{align*}
				which is a consequence of  \eqref{Assump-DRP1}. 
								By $(\mathcal{P}1)_{n}$, we infer
								\begin{align*}
					\|g_{n+1}\|_{2s_{0}+2\tau+4}^{\textnormal{Lip},\nu}&\leqslant\|u_{1}\|_{2s_{0}+2\tau+4}^{\textnormal{Lip},\nu}+\sum_{k=2}^{n+1}\|u_{k}\|_{2s_{0}+2\tau+4}^{\textnormal{Lip},\nu}\\
					&\leqslant \tfrac{1}{2}C_{\ast}\varepsilon^{\mathtt{b}_0-1}\nu^{-1}+C_{\ast}\varepsilon^{\mathtt{b}_0{-1}}\nu^{-1}\sum_{k=0}^{\infty}N_{k}^{-a_2}\\
					&\leqslant \tfrac{1}{2}C_{\ast}\varepsilon^{\mathtt{b}_0-1}\nu^{-1}+CN_{0}^{-a_2}C_{\ast}\varepsilon^{\mathtt{b}_0{-1}}\nu^{-1}\\
					&\leqslant C_{\ast}\varepsilon^{\mathtt{b}_0{-1}}\nu^{-1}.
				\end{align*}
				This completes the proof of $(\mathcal{P}1)_{n+1}.$  The proof of Proposition \ref{Nash-Moser} is now complete.
\end{proof}
The next target is to study the convergence of Nash-Moser  scheme stated in Proposition \ref{Nash-Moser} and show that the limit is a solution to the problem \eqref{main-eq1}. For this aim, we need to introduce the final Cantor set,
\begin{equation}\label{Cantor-set01}\mathtt{C}_{\infty}:=\bigcap_{m\in\mathbb{N}}\mathcal{A}_{m}.
\end{equation}
\begin{cor}\label{prop-construction}
There exists $\lambda\in\mathcal{O}\mapsto g_\infty$ satisfying
$$\|g_\infty\|_{2s_{0}+2\tau+4}^{\textnormal{Lip},\nu}\leqslant C_{\ast}\varepsilon^{\mathtt{b}_0-1}\nu^{-1}\quad \hbox{and}\quad \|g_\infty-g_m\|_{2s_{0}+2\tau+4}^{\textnormal{Lip},\nu}\leqslant  C_{\ast}\varepsilon^{\mathtt{b}_0-1}\nu^{-1}N_{m}^{-a_{2}},
$$
such that  
$$\forall \lambda\in \mathtt{C}_{\infty}, \quad \mathcal{F}\big(g_{\infty}(\lambda)\big)=0.
$$
\end{cor}
\begin{proof}
According to Proposition \ref{Nash-Moser} one may write for each $m\geqslant1,$
$$
g_m=\sum_{n=1}^m u_n.
$$
Define the formal infinite sum
$$
g_\infty:=\,\sum_{n=1}^\infty u_n.
$$
By using the estimates of $(\mathcal{P}1)_{m}$, \eqref{definition of Nm} and \eqref{lambda-choice} we get
\begin{align*}
\|g_\infty\|_{2s_{0}+2\tau+4}^{\textnormal{Lip},\nu}&\leqslant\sum_{n=1}^\infty \|u_n\|_{2s_0+2\tau+4}^{\textnormal{Lip},\nu}\\
&\leqslant \tfrac12 C_{\ast}\varepsilon^{\mathtt{b}_0-1}\nu^{-1}+C_{\ast}\varepsilon^{\mathtt{b}_0{-1}}\nu^{-1}\sum_{n=1}^\infty N_{n}^{-a_{2}}\\
&\leqslant \tfrac12 C_{\ast}\varepsilon^{\mathtt{b}_0-1}\nu^{-1}+C_{\ast}\varepsilon^{\mathtt{b}_0{-1}}\nu^{-1}N_{0}^{-a_{2}}.
\end{align*}
Thus for $ \varepsilon$ small enough,  we deduce that
\begin{align*}
\|g_\infty\|_{2s_{0}+2\tau+4}^{\textnormal{Lip},\nu}
&\leqslant  C_{\ast}\varepsilon^{\mathtt{b}_0{-1}}\nu^{-1}.
\end{align*}
On the other hand, we get in a similar way
\begin{align*}
\|g_\infty-g_m\|_{2s_{0}+2\tau+4}^{\textnormal{Lip},\nu}&\leqslant\sum_{n=m+1}^\infty \|u_n\|_{2s_{0}+2\tau+4}^{\textnormal{Lip},\nu}\leqslant C_{\ast}\varepsilon^{\mathtt{b}_0{-1}}\nu^{-1}\sum_{n=m}^\infty N_{n}^{-a_{2}}\leqslant C_{\ast}\varepsilon^{\mathtt{b}_0{-1}}\nu^{-1}N_{m}^{-a_{2}}.
\end{align*}
It follows, from Sobolev embeddings,  that the sequence $(g_m)_{m\geqslant1}$ converges pointwise to $g_\infty$. Applying $(\mathcal{P}2)_{m}$ yields 
$$\|\mathcal{F}(g_{m})\|_{s_{0},\mathtt{C}_\infty}^{{\textnormal{Lip},\nu}}\leqslant C_{\ast}\varepsilon^{\mathtt{b}_0} N_{m-1}^{-a_{1}}.
$$
Passing to the limit, we conclude that
$$\forall \lambda\in \mathtt{C}_{\infty}, \quad \mathcal{F}\big(g_{\infty}(\lambda)\big)=0.
$$
This completes  the proof of the desired result.
\end{proof}

\subsection{Cantor set measure}\label{sec:cantorset}

This section is devoted to estimating the measure of the final Cantor set
$\mathtt{C}_{\infty}$ defined in \eqref{Cantor-set01}. Since $\mathtt{C}_{\infty}$
arises from the spectral conditions imposed on the linearized operators along
the approximate sequence $(g_m)_{m\in\mathbb{N}}$, it is convenient to isolate a
subset that depends solely on the limiting solution $g_\infty$ constructed in
Corollary~\ref{prop-construction}. The measure of this reduced set can then be
computed directly.\\
Consider the following sets, for $k=1,2$
$$
\mathcal{C}_{\infty}^{k}=\Big\{\lambda\in\mathcal{O}, \; \forall \ell\in\mathbb{Z},\;  \forall    |j|\geqslant k, \; \, \left|\varepsilon^2|\ln\varepsilon| \omega(\lambda)\ell+\mu_{j,k}(\lambda,g_\infty)\right|\geqslant  {2\nu}{| j |^{-\tau}}\Big\},
$$
where we recall from \eqref{def:muj2} and the Cantor set in Proposition \ref{Second-change}:
\begin{align}\label{mu-j1}
					 \nonumber\mu_{j,k}(\lambda,g)&= j\mathtt{c}_1(\lambda,g)+(k-1)\left(\mathtt{c}_2(\lambda,g)\tfrac{j}{|j|}	+\varepsilon^2\mathtt{c}_3(\lambda,g)\mathtt{d}_j\right)\\
                     &= j\left( -\tfrac12+\varepsilon^2|\ln\varepsilon|^\frac12\mathtt{c}_0(\lambda,g)\right)+(k-1)\left(\left( \tfrac12+\varepsilon^2\mathtt{c}_4(\lambda,g)\right)\tfrac{j}{|j|}	+\varepsilon^2\mathtt{c}_3(\lambda,g)\mathtt{d}_j\right),
								\end{align}
                                and $\mathtt{d}_j$ is defined in Proposition \ref{prop-diagonal-impo}.
The following result holds.
\begin{lem}\label{Lem-Cantor-measu}
Under \eqref{Assump-DRP1},  we have the inclusion
$$
\mathcal{C}_{\infty}^{1}\cap \mathcal{C}_{\infty}^{2}\subset \mathtt{C}_{\infty}.
$$
In addition, there exists $\varepsilon_0>0$ and $C>0$ such that for any $\varepsilon\in(0,\varepsilon_0)$ 
$$|\mathcal{O}\backslash  \mathtt{C}_{\infty}|\leqslant C|\ln\varepsilon|^{\delta-1},
$$
where $\delta$ is defined through \eqref{lambda-choice}-\eqref{Assump-DRP1}.
\end{lem}
\begin{proof}
It follows directly from the definitions of $\mathcal{O}_{n}^1(g)$ and $\mathcal{O}_{n}^2(g)$ given  in Propositions \ref{Second-change} and \ref{prop-perp-orth-1} that
\begin{align*}\mathtt{C}_{\infty}=\bigcap_{m\in\mathbb{N}}\mathcal{A}_{m}=\mathcal{O}_{\infty,1}\cap\mathcal{O}_{\infty,2},
\end{align*}
where 
$$
\mathcal{O}_{\infty,k}:=\Big\{\lambda\in \mathcal{O},\; \forall \ell\in\mathbb{Z},\; n\in\N,\;  k\leqslant |j|\leqslant N_n, \;\,\left|\varepsilon^2|\ln\varepsilon| \omega(\lambda)\ell+\mu_{j,k}(\lambda,g_n)\right|\geqslant  {\nu}{| j |^{-\tau}}\Big\}.
$$
Our goal is to prove that, for each  $k=1,2$, 
\begin{equation}\label{Mahma-l}
\mathcal{C}_{\infty}^{k}\subset \mathcal{O}_{\infty,k},
\end{equation}
which in turn implies
$
\mathcal{C}_{\infty}^{1}\cap \mathcal{C}_{\infty}^{2}\subset \mathtt{C}_{\infty}.
$\\
To this end, fix $\lambda\in \mathcal{C}_{\infty}^{k}$ and $n\in\mathbb{N}$. Then, for any integers $k\leqslant |j|\leqslant N_n$, we can apply  triangle inequality to obtain
\begin{align*}
\left|\varepsilon^2|\ln\varepsilon| \omega(\lambda)\ell+\mu_{j,k}(\lambda,g_n)\right|&\geqslant \left|\varepsilon^2|\ln\varepsilon| \omega(\lambda)\ell+\mu_{j,k}(\lambda,g_\infty)\right|-|\mu_{j,k}(\lambda,g_\infty)-\mu_{j,k}(\lambda,g_n)|\\
&\geqslant {2\nu}{| j |^{-\tau}}-|\mu_{j,k}(\lambda,g_\infty)-\mu_{j,k}(\lambda,g_n)|.
\end{align*}
By applying Propositions \ref{Second-change}-(4) and Corollary \ref{prop-construction}, we obtain, for $\varepsilon$ sufficiently small, that
\begin{align*}
|\mu_{j,k}(\lambda,g_\infty)-\mu_{j,k}(\lambda,g_n)|&\leqslant \varepsilon^2|\ln\varepsilon|^\frac12|j|\left|\mathtt{c}_0(\lambda,g_\infty)-\mathtt{c}_0(\lambda,g_n)\right|+\varepsilon^2(k-1)\left|\mathtt{c}_4(\lambda,g_\infty)-\mathtt{c}_4(\lambda,g_n)\right|\\ &\quad +\varepsilon^2(k-1)|\mathtt{d}_j|\left|\mathtt{c}_3(\lambda,g_\infty)-\mathtt{c}_3(\lambda,g_n)\right|
\\ &
\leqslant C\varepsilon^{2}|\ln\varepsilon|^{\frac12} |j|\|g_{\infty}-g_n\|_{2s_{0}+2\tau+3}^{\textnormal{Lip},\nu}\\
&  \leqslant C_{\ast} \varepsilon^{1+\mathtt{b}_0}|\ln\varepsilon|^{\frac12}\nu^{-1}|j|N_n^{-a_2}\\
&\leqslant  C_{\ast}\nu |j|^{-\tau}\varepsilon^{\mathtt{b}_0-3}|\ln\varepsilon|^{\frac12-2\delta}N_n^{1+\tau-a_2}.
\end{align*}
From \eqref{choice-f1}, together with the assumptions $1+\tau<a_2$ and  $\delta>\frac12$  (see \eqref{Assump-DRP1}) we deduce that,   for  $\varepsilon$ sufficiently small,
\begin{align*}
|\mu_{j,k}(\lambda_0,g_\infty)-\mu_{j,k}(\lambda,g_n)|&
\leqslant  C_{\ast}\nu|j|^{-\tau}\varepsilon^{N-\mu-2}|\ln\varepsilon|^{\frac12-2\delta}N_n^{1+\tau-a_2}\leqslant  \nu |j|^{-\tau}.
\end{align*} 
Consequently, for any $n\in\N$,  any $\ell\in\Z$ and $k\leqslant |j|\leqslant N_n$, we obtain
\begin{align*}
|\varepsilon^2 |\ln\varepsilon| \omega(\lambda)\ell+\mu_{j,k}(\lambda,g_n)|&\geqslant \nu |j|^{-\tau}\cdot
\end{align*}
This shows that $\lambda\in \mathcal{O}_{\infty,k}$, which completes the proof of \eqref{Mahma-l}.\\ We now turn to the measure estimate. First, observe that
$$
\big|\mathcal{O}\backslash  \mathtt{C}_{\infty}\big|\leqslant \big|\mathcal{O}\backslash  \mathcal{C}_{\infty}^{1}\big|+\big|\mathcal{O}\backslash \mathcal{C}_{\infty}^{2}\big|.
$$ 
From \eqref{lambda-choice}, we have
$$\mathcal{C}_{\infty}^{k}=\bigcap_{n\in\N}\bigcap_{\substack{(\ell,j)\in\mathbb{Z}^{2}\\ k\leqslant |j|\leqslant N_{n}}}A_{\ell,j}^k,
$$
where
$$A_{\ell,j}^k:=\left\lbrace \lambda\in \mathcal{O}=[\lambda_*,\lambda^*];\;\, \big|\varepsilon^2|\ln\varepsilon|\omega(\lambda)  \ell+\mu_{j,k}(\lambda,g_\infty)\big|\geqslant 2{\varepsilon^2|\ln\varepsilon|^{\delta}}{| j|^{-\tau}}\right\rbrace.
$$
Notice that
$$
\bigcap_{\substack{(\ell,j)\in\mathbb{Z}^{2}\\ k\leqslant |j|}}A_{\ell,j}^k\subset \mathcal{C}_{\infty}^{k}.
$$
By Proposition \ref{prop:period},  the function   $\lambda\in\mathcal{O}\mapsto \omega(\lambda)$ is positive and strictly decreasing. Therefore,
$$
\forall \lambda\in[\lambda_*,\lambda^*],\quad 0<\omega^*=\omega(\lambda^*)\leqslant \omega(\lambda)\leqslant \omega(\lambda_*)=\omega_*.
$$
Moreover, from \eqref{mu-j1}, Proposition \ref{prop-diagonal-impo}-(a) and Proposition \ref{prop:Phi4}-(b) we have the bounds: for all 
$0<\varepsilon\leqslant \varepsilon_0$ 
\begin{equation}\label{est-c_0}
-\tfrac23\leqslant \mathtt{c}_1(\lambda,g_\infty)\leqslant -\tfrac13, \quad
\tfrac13\leqslant \mathtt{c}_2(\lambda,g_\infty)\leqslant {\tfrac{13}{24}},  \quad\textnormal{and}\quad |\mathtt{d}_j|\leqslant \tfrac{1}{12}.
\end{equation}
In what follows, we distinguish several cases.

\medskip

\medskip\noindent $\diamond$ Case $\varepsilon^2 |\ln\varepsilon| \omega^*|\ell|\geqslant |j|\geqslant 1, k=1.$ One has by the triangle inequality 
\begin{align*}
 \big|\varepsilon^2|\ln\varepsilon|\omega(\lambda)  \ell+j \mathtt{c}_1(\lambda,g_\infty)\big|&\geqslant \varepsilon^2|\ln\varepsilon|\omega^* |\ell|-\tfrac23|j|\\
 &\geqslant\tfrac13|j|\\
 &\geqslant 
 2{\varepsilon^{2}|\ln\varepsilon|^\delta}{| j|^{-\tau}}\cdot
 \end{align*}
Hence
$$
 A_{\ell,j}^1=[\lambda_*,\lambda^*].
 $$

\medskip\noindent $\diamond$ Case $\varepsilon^2|\ln\varepsilon| \omega^*|\ell|\geqslant |j|\geqslant 2, k=2.$ By the triangle inequality, we infer for small $\varepsilon$ 
\begin{align*}
 \left|\varepsilon^2|\ln\varepsilon|\omega(\lambda)  \ell+j \mathtt{c}_1(\lambda,g_\infty)+\mathtt{c}_2(\lambda,g_\infty)\tfrac{j}{|j|}	+\varepsilon^2\mathtt{c}_3(\lambda,g_\infty)\mathtt{d}_j\right|&\geqslant \varepsilon^2|\ln\varepsilon|\omega^* |\ell|-\tfrac23|j|-{\tfrac{13}{24}}-\varepsilon^2|\mathtt{d}_j|\\
 &\geqslant\tfrac13|j|-{\tfrac{13}{24}}-\tfrac{\varepsilon^2}{12}\geqslant \tfrac{1}{12}\\
 &\geqslant {\varepsilon^{2}|\ln\varepsilon|^\delta}{| j|^{-\tau}}\cdot
 \end{align*}
 Therefore
 $$
 A_{\ell,j}^2=[\lambda_*,\lambda^*].
 $$

\medskip\noindent $\diamond$ Case $|j|\geqslant 24\varepsilon^2 |\ln\varepsilon| \omega_*|\ell|, k=1.$ One gets by the triangle inequality 
 \begin{align*}
 \forall \,|j|\geqslant 1,\quad \big|\varepsilon^2|\ln\varepsilon|\omega(\lambda)  \ell+j \mathtt{c}_1(\lambda,g_\infty)\big|&\geqslant \tfrac13|j|-\varepsilon^2|\ln\varepsilon|\omega_*|\ell|\\
 &\geqslant \big(\tfrac{1}{24}|j|-\varepsilon^2|\ln\varepsilon|\omega_*|\ell|\big) +\tfrac{7}{24}|j|\\
 &\geqslant\tfrac{7}{24}\cdot
 \end{align*}
  It follows that for small $\varepsilon$
  \begin{align*}
 \big|\varepsilon^2|\ln\varepsilon|\omega(\lambda)  \ell+j \mathtt{c}_1(\lambda,g_\infty)\big|
 &\geqslant  {\varepsilon^{2}|\ln\varepsilon|^{\delta}}{| j|^{-\tau}}\cdot
 \end{align*}
 Hence
 $$
  A_{\ell,j}^1=[\lambda_*,\lambda^*].
  $$

 \medskip\noindent $\diamond$ Case $|j|\geqslant  24\varepsilon^2 |\ln\varepsilon|\omega_*|\ell|, k=2.$ By the triangle inequality and \eqref{est-c_0}, we obtain for $|j|\geqslant 2$
 \begin{align*}
 \left|\varepsilon^2|\ln\varepsilon|\omega(\lambda)  \ell+j \mathtt{c}_1(\lambda,g_\infty)+\mathtt{c}_2(\lambda,g)\tfrac{j}{|j|}	+\varepsilon^2\mathtt{c}_3(\lambda,g)\mathtt{d}_j\right|&\geqslant \tfrac13|j|-\varepsilon^2|\ln\varepsilon|\omega_*|\ell|-{\tfrac{13}{24}}-\tfrac{1}{12}\varepsilon^2\\
 &\geqslant \big(\tfrac{1}{24}|j|-\varepsilon^2|\ln\varepsilon|\omega_*|\ell|\big) +\tfrac{7}{24}|j|-\tfrac{27}{48}\\
 &\geqslant\tfrac{1}{48}\cdot
 \end{align*}
 Thus, for small $\varepsilon$ we infer
 \begin{align*}
 \left|\varepsilon^2|\ln\varepsilon|\omega(\lambda)  \ell+j \mathtt{c}_1(\lambda,g_\infty)+\mathtt{c}_2(\lambda,g)\tfrac{j}{|j|}	+\varepsilon^2\mathtt{c}_3(\lambda,g)\mathtt{d}_j\right|
 &\geqslant  {\varepsilon^{2}|\ln\varepsilon|^{\delta}}{| j|^{-\tau}},
 \end{align*}
 implying that
 $$
 A_{\ell,j}^2=[\lambda_*,\lambda^*].
 $$
 It follows that
 \begin{align}\label{Cinfty}\bigcap_{\substack{|j|\leqslant 24\varepsilon^2 |\ln\varepsilon|\omega_*|\ell|\\ \varepsilon^2|\ln\varepsilon| \omega^*|\ell|\leqslant |j|}}\ A_{\ell,j}^k\subset \mathcal{C}_{\infty}^{k}.
 \end{align}
 Define the function
 \begin{align*}
     f_{\ell,j}(\lambda)&:=\varepsilon^2|\ln\varepsilon|\omega(\lambda)  \ell+j\left( -\tfrac12+\varepsilon^2|\ln\varepsilon|^\frac12\mathtt{c}_0(\lambda,g_\infty)\right)\\ &\quad +(k-1)\left(\left( \tfrac12+\varepsilon^2\mathtt{c}_4(\lambda,g_\infty)\right)\tfrac{j}{|j|}	+\varepsilon^2\mathtt{c}_3(\lambda,g_\infty)\mathtt{d}_j\right)\cdot
 \end{align*}
 Differentiating in $\lambda$ and using \eqref{mu-j1} yields
 $$
  f_{\ell,j}^\prime(\lambda)=\varepsilon^2|\ln\varepsilon|\omega^\prime(\lambda)  \ell+j\varepsilon^2|\ln\varepsilon|^\frac12\mathtt{c}_0^\prime(\lambda) +\varepsilon^{2}(k-1){\mathtt{c}}_4^\prime(\lambda)\tfrac{j}{|j|}+\varepsilon^2(k-1){\mathtt{c}}_3^\prime(\lambda)\mathtt{d}_j.
   $$
 Applying Proposition \ref{prop:period},  together with Propositions \ref{Second-change}-(4) and  Proposition \ref{prop:Phi4}-(b) we obtain two constants $\underline{c}>0,\overline{c}>0$ such that 
 $$\forall \lambda\in[\lambda_*,\lambda^*],\quad \underline{c}\leqslant |\omega_0^\prime(\lambda)|\quad\hbox{and}\quad |\mathtt{c}_k^\prime(\lambda)|\leqslant \overline{c}, \quad k\in\{0,3,4\},
 $$ and thus
 \begin{align*}
 \forall \lambda\in [\lambda_*,\lambda^*],\quad  |f_{j,\ell}^\prime(\lambda)|\geqslant \underline{c}\,\varepsilon^2 |\ln\varepsilon|  |\ell|-\overline{c}\,|j|\varepsilon^2|\ln\varepsilon|^\frac12-2\overline{c}\,\varepsilon^{2} .
  \end{align*}
  Hence, as $|j|\leqslant 24\varepsilon^2 |\ln\varepsilon|\omega_*|\ell|,$ we get
  \begin{align*}
 | f_{j,\ell}^\prime(\lambda)|\geqslant |j|\big(\tfrac{\underline{c}}{24\omega_*}-\overline{c}\,\varepsilon^2|\ln\varepsilon|^\frac12\big)-2\overline{c}\,\varepsilon^{2},
  \end{align*}
  which implies for  $\varepsilon$ small enough 
   \begin{align*}
  \forall \,\lambda\in\mathcal{O},\quad |f_{j,\ell}^\prime(\lambda)|\geqslant \tfrac{\overline{c}}{48\omega_*} |j|.
   \end{align*}
  Applying Lemma \ref{Piralt}, we deduce that
  \begin{align*}
  \left|\mathcal{O}\backslash A_{\ell,j}^k\right|\leqslant C\frac{\varepsilon^{2}|\ln\varepsilon|^{\delta}}{| j|^{1+\tau}}\cdot
  \end{align*}
  Consequently, by  \eqref{Cinfty} and  the fact $\tau>1$, it follows that
  \begin{align*}
  \big|[\lambda_*,\lambda^*]\backslash \mathcal{C}_{\infty}^{k}\big|&\leqslant C \sum_{ \varepsilon^2|\ln\varepsilon| \omega^*|\ell|\leqslant |j|}\frac{\varepsilon^{2}|\ln\varepsilon|^{\delta}}{| j|^{1+\tau}}\leqslant C |\ln\varepsilon|^{\delta-1} \sum_{   |j|\geqslant 1}\tfrac{1}{| j|^{\tau}}\leqslant C |\ln\varepsilon|^{\delta-1}.
  \end{align*}
  Finally, we get
  \begin{align*}
 \big|[\lambda_*,\lambda^*]\backslash \mathtt{C}_{\infty}\big|\leqslant  \big|[\lambda_*,\lambda^*]\backslash \mathcal{C}_{\infty}^{1}\big|+\big|[\lambda_*,\lambda^*]\backslash \mathcal{C}_{\infty}^{2}\big|
  &\leqslant C |\ln\varepsilon|^{\delta-1}.
  \end{align*}
  This achieves the proof of the desired result.
\end{proof}

\appendix

\section{Transport equation with constant coefficients}\label{apendix:transport}

In order to analyze a degenerate transport operator with constant coefficients, it is necessary to understand the solvability of the simple  equation
$$
\epsilon_1 \omega(\lambda)\,\partial_\varphi \varrho
+ \mathtt{c}(\lambda)\,\partial_\theta \varrho = h,
$$
posed on the torus $\T^2$.  
At the formal level, the inversion is straightforward in Fourier series. However, the denominators 
$$
\epsilon_1 \omega(\lambda) l + j \mathtt{c}(\lambda), \qquad (l,j)\in\Z^2\setminus\{0\},
$$
may become arbitrarily small, giving rise to the classical \emph{small-divisor problem}.  

To overcome this obstruction, we introduce Cantor-type conditions that exclude resonant parameters $\lambda$ for which the divisors vanish or are too small. On such nonresonant sets, the operator is invertible with a well-defined Fourier representation.  
Nevertheless, in the forthcoming analysis we require an operator that is defined smoothly on the whole parameter set $\mathcal{O}$, not just on the Cantor subset. For this reason, we construct a regularized extension of the inverse, obtained by introducing a smooth cutoff function. This modified inverse coincides with the true one on the Cantor set, but remains bounded and smooth for all $\lambda\in \mathcal{O}$.  

The following result formalizes this construction: it provides tame estimates for the extended inverse operator, and shows that on nonresonant truncated Cantor sets it recovers the exact inversion of the transport operator. This lays the foundation for its use in iterative schemes where uniform control in $\lambda$ is required.

Let $\nu, \epsilon_1, \epsilon_2 \in (0,1]$ and $\tau > 0$.  
Consider two smooth functions $\lambda \in \mathcal{O} \mapsto \mathtt{c}(\lambda), \omega(\lambda)\in \RR$.  
We introduce the Cantor-type set
$$
\mathtt{C}^{\tau}
=\Big\{ \lambda\in \mathcal{O}:\; 
|\epsilon_1 \omega(\lambda) \ell + j \mathtt{c}(\lambda)| > \tfrac{\epsilon_2}{\langle j\rangle^{\tau}}
\quad \forall (\ell,j)\in \Z^2\setminus\{(0,0)\}\Big\}.
$$
For each $N\in \N^\ast$, we also define the truncated Cantor set
$$
\mathtt{C}^{\tau,N}
=\Big\{ \lambda\in \mathcal{O}:\;
|\ell|\leqslant N,\quad 
|\epsilon_1 \omega(\lambda) \ell + j \mathtt{c}(\lambda)| > \tfrac{\epsilon_2}{\langle j\rangle^{\tau}}
\quad \forall (\ell,j)\in \Z^2\setminus\{(0,0)\}\Big\}.
$$
Let $h:\mathcal{O}\times \T^2 \to \RR$ be a smooth zero-average function, with Fourier expansion
$$
h(\lambda,\varphi,\theta) 
= \sum_{(\ell,j)\in\Z^{2}\setminus\{0\}} h_{l,j}(\lambda)\, {\bf e}_{\ell,j},
\qquad 
{\bf e}_{\ell,j}(\varphi,\theta) := e^{ i(\ell\varphi+j\theta)}.
$$
We study the transport equation
\begin{equation}\label{eqq-kk}
 \epsilon_1 \omega(\lambda)\, \partial_\varphi \varrho 
 + \mathtt{c}(\lambda)\, \partial_\theta \varrho 
 = h,
\end{equation}
posed on the torus $\T^2$, where $\varrho:\T^2 \to \RR$.  
If $\lambda\in \mathtt{C}^{\tau}$, equation \eqref{eqq-kk} can be solved by Fourier series inversion:
$$
\varrho(\lambda,\varphi,\theta)
= - i \sum_{(\ell,j)\in \Z^2\setminus\{0\}} 
\frac{h_{l,j}(\lambda)}{\epsilon_1 \omega(\lambda) \ell + j\mathtt{c}(\lambda)}\,
{\bf e}_{\ell,j}.
$$
To extend this inverse to all $\lambda\in \mathcal{O}$ in a smooth way, we introduce the following regularized inverse operator:
\begin{align}\label{ExtendMM00}
(\epsilon_1\omega(\lambda)\partial_\varphi + \mathtt{c}(\lambda)\partial_\theta)_{\mathrm{ext}}^{-1}h
&:= - i \sum_{(\ell,j)\in\Z^{2}\setminus\{0\}}
\frac{\chi\big( [\epsilon_1\omega(\lambda) \ell + j\mathtt{c}(\lambda)] \epsilon_2^{-1}\langle j\rangle^\tau\big)}
{\epsilon_1\omega(\lambda) \ell + j\mathtt{c}(\lambda)}\,
h_{\ell,j}(\lambda)\,{\bf e}_{\ell,j}\\
&= - i \epsilon_2^{-1} \sum_{(\ell,j)\in \Z^2\setminus\{0\}}
\chi_1\Big([\epsilon_1\omega(\lambda) \ell + j\mathtt{c}(\lambda)] \epsilon_2^{-1}\langle j\rangle^\tau\Big)
\,\langle j\rangle^{\tau} h_{\ell,j}(\lambda)\,{\bf e}_{\ell,j},\nonumber
\end{align}
where $\chi\in \mathscr{C}^\infty(\RR,\RR)$ is an even, nonnegative cutoff function such that
$$
\chi(\xi)=
\begin{cases}
0, & |\xi|\leqslant \tfrac13,\\
1, & |\xi|\geqslant \tfrac12,
\end{cases}
$$
and $\chi_1(x) := \tfrac{\chi(x)}{x}\in \mathscr{C}_b^\infty(\RR)$.  
The following lemma is a straightforward adaptation of classical results (see, e.g., \cite{Baldi-berti,HHM24}).  
\begin{lem}\label{L-InvertMM2}
Let $\nu, \epsilon_1, \epsilon_2 \in (0,1]$, $\tau>0$, and assume $\mathtt{c},\omega:\mathcal{O}\to\RR$ are smooth functions satisfying
$$
\inf_{\lambda\in \mathcal{O}} |\omega'(\lambda)| > 0,
\qquad 
\nu \epsilon_2^{-1}\leqslant 1,
\qquad 
\|(\omega,\mathtt{c})\|^{\mathrm{Lip},\nu}\leqslant 1.
$$
Then for all $s\geqslant 1$, one has the tame estimate
$$
\big\| (\epsilon_1\omega(\lambda)\partial_\varphi + \mathtt{c}(\lambda)\partial_\theta)_{\mathrm{ext}}^{-1} h \big\|_{s}^{\mathrm{Lip},\nu}
\;\leqslant\; C \epsilon_2^{-1} \,\||\partial_\theta|^{2\tau} h\|_{s}^{\mathrm{Lip},\nu}.
$$
Moreover, for any $N\in\N$ and any $\lambda\in \mathtt{C}^{\tau,N}$,
$$
(\epsilon_1\omega(\lambda)\partial_\varphi+\mathtt{c}(\lambda)\partial_\theta)
\,(\epsilon_1\omega(\lambda)\partial_\varphi+\mathtt{c}(\lambda)\partial_\theta)_{\mathrm{ext}}^{-1}\,\Pi_N
=\Pi_N,
$$
where $\Pi_N$ denotes the Fourier projector
$$
\Pi_N\!\Big(\sum_{(\ell,j)\in \Z^2} h_{\ell,j}{\bf e}_{\ell,j}\Big)
=\sum_{\substack{(\ell,j)\in \Z^2\\ |j|\leqslant N}} h_{l,j}{\bf e}_{\ell,j}.
$$
\end{lem}

\begin{proof}
Set $$L(\lambda):=\epsilon_{1}\,\omega(\lambda)\,\partial_\varphi+\mathtt{c}(\lambda)\,\partial_\theta,$$ and recall that $|\partial_\theta|^\sigma$ denotes the multiplier with symbol $\langle j\rangle^\sigma$ and that
$$\|u\|_{s}^{\mathrm{Lip},\nu}=\sup_{\lambda\in\mathcal{O}}\|u(\lambda)\|_{H^s}+\nu\,\sup_{\lambda_1\neq\lambda_2\in\mathcal{O}}\tfrac{\|u(\lambda_1)-u(\lambda_2)\|_{H^{s-1}}}{|\lambda_1-\lambda_2|}.$$
By the explicit extension formula \eqref{ExtendMM00},
$$
L_{\mathrm{ext}}^{-1}h
=- i\,\epsilon_2^{-1}\sum_{(\ell,j)\in\Z^2\setminus\{(0,0)\}}
\chi_1\!\Big([\epsilon_1\omega(\lambda)\ell+j\mathtt{c}(\lambda)]\,\epsilon_2^{-1}\langle j\rangle^{\tau}\Big)\,
\langle j\rangle^{\tau}\,h_{\ell,j}(\lambda)\,{\bf e}_{\ell,j}.
$$
Taking the $H^s$-norm and using that $\chi_1$ is bounded immediately yields the  bound
$$
\|L_{\mathrm{ext}}^{-1}h\|_{H^s}\;\lesssim\;\epsilon_2^{-1}\,\||\partial_\theta|^{\tau}h\|_{H^s}.
$$
For the Lipschitz seminorm, differentiate in $\lambda$ (equivalently argue by increments), we obtain
$$
\partial_\lambda(L_{\mathrm{ext}}^{-1}h)=I_1+I_2,
$$
with
\begin{align*}
I_1&=- i\,\epsilon_2^{-1}\!\!\sum_{(\ell,j)\neq(0,0)}
\chi_1([\epsilon_1\omega(\lambda)\ell+j\mathtt{c}(\lambda)]\,\epsilon_2^{-1}\langle j\rangle^{\tau})\,\langle j\rangle^{\tau}\,(\partial_\lambda h_{\ell,j})(\lambda)\,{\bf e}_{\ell,j},
\\
&
I_2=- i\,\epsilon_2^{-2}\!\!\sum_{(\ell,j)\neq(0,0)}
\chi_1'([\epsilon_1\omega(\lambda)\ell+j\mathtt{c}(\lambda)]\,\epsilon_2^{-1}\langle j\rangle^{\tau})\,\partial_\lambda\!\Big([\epsilon_1\omega \ell+j\mathtt{c}]\,\Big)\,
\langle j\rangle^{2\tau} h_{\ell,j}\,{\bf e}_{\ell,j}.
\end{align*}
Note that $\partial_\lambda([\epsilon_1\omega \ell+j\mathtt{c}])=\epsilon_1\omega'(\lambda)\ell+j\mathtt{c}'(\lambda)$.
The term $I_1$ is straightforward:
$$
\nu\|I_1\|_{H^{s-1}}\;\lesssim\;\epsilon_2^{-1}\,\|\,|\partial_\theta|^{\tau}h\|_{s}^{\mathrm{Lip},\nu}.
$$
To estimate $I_2$, use the standard symbol bound
$
|\chi_1'(x)|\lesssim x^{-2}\mathbf 1_{\mathrm{supp}\,\chi_1}(x)
$,
leading to
$$
\big|\chi_1'([\epsilon_1\omega(\lambda)\ell+j\mathtt{c}(\lambda)]\,\epsilon_2^{-1}\langle j\rangle^{\tau})\big|
\;\lesssim\;[\epsilon_1\omega \ell+j\mathtt{c}]^{-2}\,\epsilon_2^{2}\,\langle j\rangle^{-2\tau},
$$
and hence, after multiplying by the prefactors in $I_2$, using the triangle inequality and the fact that $\inf_{\lambda\in\mathcal{O}}|\omega^\prime(\lambda)|>0$
\begin{align*}
\epsilon_2^{-2}\,|\chi_1'([\epsilon_1\omega(\lambda)\ell+j\mathtt{c}(\lambda)]\,\epsilon_2^{-1}\langle j\rangle^{\tau})|\,\big|\epsilon_1\omega' \ell+j\mathtt{c}'\big|\,\langle j\rangle^{2\tau}
&\;\lesssim\;
\frac{\big|\epsilon_1\omega' \ell+j\mathtt{c}'\big|}{[\epsilon_1\omega \ell+j\mathtt{c}]^{2}}
\;\mathbf 1_{\mathrm{supp}\,\chi_1}\\
&
\;\lesssim\;
\Big(\frac{1}{|\epsilon_1\omega \ell+j\mathtt{c}|}+\frac{|j|}{|\epsilon_1\omega \ell+j\mathtt{c}|^{2}}\Big)\mathbf 1_{\mathrm{supp}\,\chi_1}.
\end{align*}
On the support of $\chi$, we get
$$
\frac{1}{|\epsilon_1\omega \ell+j\mathtt{c}|}\lesssim \epsilon_2^{-1}\langle j\rangle^{\tau},
\qquad
\frac{|j|}{|\epsilon_1\omega \ell+j\mathtt{c}|^{2}}\lesssim \epsilon_2^{-2}\langle j\rangle^{1+2\tau}.
$$
Consequently,
$$
\nu\|I_2\|_{H^{s-1}}\;\lesssim\;\nu\,\big(\epsilon_2^{-1}\||\partial_\theta|^{\tau}h\|_{H^{s-1}}
+\epsilon_2^{-2}\||\partial_\theta|^{2\tau+1}h\|_{H^{s-1}}\big)
\;\lesssim\;\nu\,\epsilon_2^{-2}\,\||\partial_\theta|^{2\tau}h\|_{H^{s}},
$$
and, under the  smallness condition $\nu\,\epsilon_2^{-1}\leqslant 1$,
$$
\nu\|I_2\|_{H^{s-1}}\;\lesssim\;\epsilon_2^{-1}\,\||\partial_\theta|^{2\tau}h\|_{s}^{\mathrm{Lip},\nu}.
$$
Combining the two estimates for $I_1$ and $I_2$ gives
$$
\|L_{\mathrm{ext}}^{-1}h\|_{s}^{\mathrm{Lip},\nu}\;\lesssim\;\epsilon_2^{-1}\,\||\partial_\theta|^{2\tau}h\|_{s}^{\mathrm{Lip},\nu}.
$$
For the second point, let $\Pi_N$ be the projector onto modes $|j|\leqslant N$. By \eqref{ExtendMM00},
$$
L_{\mathrm{ext}}^{-1}\Pi_N h
=- i\sum_{\substack{(\ell,j)\neq(0,0)\\ |j|\leqslant N}}
\frac{\chi\!\big([\epsilon_1\omega \ell+j\mathtt{c}]\,\epsilon_2^{-1}\langle j\rangle^\tau\big)}{\epsilon_1\omega \ell+j\mathtt{c}}\,
h_{\ell,j}\,{\bf e}_{\ell,j}.
$$
By construction of the Cantor-like set $\mathtt{C}_N^{\tau}$, for $\lambda\in\mathtt{C}_N^{\tau}$ and $|j|\leqslant N$ the cutoff is identically $1$.
Hence
$$
L_{\mathrm{ext}}^{-1}\Pi_N h
=- i\sum_{\substack{(\ell,j)\neq(0,0)\\ |j|\leqslant N}}\frac{h_{\ell,j}}{\epsilon_1\omega \ell+j\mathtt{c}}\,{\bf e}_{\ell,j},
$$
and therefore
$$
L\,\big(L_{\mathrm{ext}}^{-1}\Pi_N h\big)
= i\sum_{\substack{(\ell,j)\neq(0,0)\\ |j|\leqslant N}}h_{\ell,j}\,{\bf e}_{\ell,j}
=\Pi_N h.
$$
This establishes the claimed identity and completes the proof.
\end{proof}
In what follows, we invoke Lemma~\ref{L-InvertMM2} to analyze the model transport equation
\begin{equation}\label{Trans-model-ideal}
  \varepsilon^2|\ln\varepsilon|\,\omega(\lambda)\,\partial_\varphi \varrho
  \;+\;\mathtt{c}(\lambda)\,\partial_\theta \varrho \;=\; f,
\end{equation}
where the forcing term $f$ has zero $\theta$–average and is spectrally localized. That is a toy model of the linear operator analyzed in this paper.

  \begin{cor}\label{Lemma-transp-inver}
      Let $\epsilon_2=\nu\in (0,1],\tau>0, $ and assume $\mathtt{c},\omega:\mathcal{O}\to\RR$ are two smooth functions satisfying
$$
\inf_{\lambda\in \mathcal{O}} |\omega'(\lambda)| > 0,\quad\inf_{\lambda\in \mathcal{O}}|\mathtt{c}(\lambda)|>0,
\quad 
\|(\omega,\mathtt{c})\|^{\mathrm{Lip},\nu}\leqslant 1.
$$
Let  $f$ be smooth  and $\langle f\rangle_\theta=0$.
Set $\epsilon_1:=\varepsilon^2|\ln\varepsilon|$.
Then there exists a smooth map $\lambda\mapsto \varrho(\lambda)$ with $\langle \varrho\rangle_\theta=0,$ and admitting for any $M\in\mathbb{N}^\ast$ the finite expansion
      \begin{align*}
          \varrho=\mathtt{c}^{-1}\sum_{k=0}^{M-1}(-1)^k\big(\tfrac{\omega(\lambda)}{\mathtt{c}(\lambda)}\epsilon_1\big)^{k}\partial_\theta^{-k-1}f+\epsilon_1^M\varrho_M\quad\hbox{with}\quad  \|\varrho_M\|_{s}^{\textnormal{Lip},\nu}\lesssim \nu^{-1}\|f\|_{s+2\tau+M+1}^{\textnormal{Lip},\nu}.
      \end{align*}
      Moreover, for any $\lambda$ in the Cantor-like set
$$
\mathtt{C}_N^\tau
= \bigcap_{\substack{(\ell,j)\in\mathbb{Z}^2\backslash\{(0,0)\}\\ |j|\leqslant N}}
\Big\{ \lambda\in\mathcal{O}\;:\;
\big| \epsilon_1\,\omega(\lambda)\,\ell + j\,\mathtt{c}(\lambda)\big|
\ge \nu\,|j|^{-\tau}\Big\},
$$
the function $\varrho(\lambda)$ solves \eqref{Trans-model-ideal} provided that $\Pi_N f=f$.
  \end{cor}
  \begin{proof}
      We seek $\varrho$ as the truncated series
      $$
      \varrho=\sum_{k=0}^M\big(\omega(\lambda)\epsilon_1\big)^k\varrho_k.
      $$
      Inserting into \eqref{Trans-model-ideal} and matching powers of $\epsilon_1\,\omega$ yields
      $$
      \mathtt{c}_1\partial_\theta \varrho_0=f,\quad \mathtt{c}_1\partial_\theta \varrho_k=-\,\partial_\varphi \varrho_{k-1},\quad 1\leqslant k\leqslant M-1,
      $$
      and the remainder equation
    $$\varepsilon^2|\ln\varepsilon|\omega(\lambda)\,\partial_\varphi \varrho_M+\mathtt{c}\partial_\theta \varrho_M=-\,\partial_\varphi\varrho_{M-1}.
      $$
      Since $\langle f\rangle_\theta=0$, $\partial_\theta^{-1}$ is well-defined on $f$, and the first two equations integrate to
      \begin{align}\label{rhoo_k}
      \varrho_0=\frac{1}{\mathtt{c}}\partial_\theta^{-1}f,\quad \varrho_k=-\big(-\tfrac{1}{\mathtt{c}}\big)^{k+1}\partial_\theta^{-k-1}\big( f \big).
      \end{align}
      As $\Pi_Nf=f$ then $\Pi_N\varrho_k=\varrho_k$ for all $0\leqslant k\leqslant M-1.$
      The last function $\varrho_n$ will be solved using Cantor like set
      $$
\mathtt{C}_N^\tau
= \bigcap_{\substack{(\ell,j)\in\mathbb{Z}^2\backslash\{(0,0)\}\\ |j|\leqslant N}}
\Big\{ \lambda\in\mathcal{O}\;:\;
\big| \epsilon_1\,\omega(\lambda)\,\ell + j\,\mathtt{c}(\lambda)\big|
\ge \nu\,|j|^{-\tau}\Big\}.
$$
Using Lemma \ref{L-InvertMM2} and by setting
$$
\varrho_M=-L_{\mathrm{ext}}^{-1}[\partial_\varphi \varrho_{M-1}],
$$we find in view of \eqref{rhoo_k}
\begin{align*}
\|\varrho_M\|_{s}^{\mathrm{Lip},\nu} &\lesssim \nu^{-1}\||\partial_\theta|^{2\tau} \partial_\varphi \varrho_{M-1}\|_{s}^{\mathrm{Lip},\nu}\lesssim \nu^{-1}\|f\|_{s+2\tau+M+1}^{\mathrm{Lip},\nu}.
\end{align*}
Moreover, for $\lambda\in \mathtt{C}_N^\tau,$ we find that $\varrho_n(\lambda)$ solves the equation. Therefore, under this restriction on $\lambda$, the function  $g$ satisfies \eqref{Trans-model-ideal}. This completes the proof. 
  \end{proof}

We shall state a particular statement of Kirszbraun Theorem \cite{Kirsz}.
		\begin{lem}[Kirszbraun Theorem]\label{thm-extend}
	Given $a<b,$ $ U$ a subset of $[a,b]$ and $H$ a Hilbert space. Let $ f:U\to H$ be  a Lipschitz function, then $f$ admits a Lipschitz  extension $F=\textnormal{Ext}f: [a,b]\to H$ with the same Lipschitz constant.
	\end{lem}
		We recall the following classical result on bi-Lipschitz functions and measure theory, used in the proof to measure the Cantor-like set that appears before.
	\begin{lem}\label{Piralt}
Let  $(\alpha,\beta)\in(\mathbb{R}_{+}^{*})^{2}$ and $f:[a,b]\to \RR$ be  a bi-Lipschitz  function   such that 
$$\forall x,y \in[a,b],\quad |f(x)-f(y)|\geqslant\beta |x-y|.$$
Then  there exists $C>0$ independent of $ \|f\|_{\textnormal{Lip}}$ such that 
$$\Big|\big\lbrace x\in [a,b];\;\, |f(x)|\leqslant\alpha\big\rbrace\Big|\leqslant C\tfrac{\alpha}{\beta}\cdot $$
\end{lem}

\section{Symplectic change of coordinates system}\label{appendix:KAM}
The main goal of this section is to discuss useful  results related to some change of coordinates system. For the proofs we refer the reader to the papers \cite{Baldi-Montalto21,BFM,FGMP19}.
Let $\beta: \mathcal{O}\times \T^{2}\to \R$ be a smooth function such that $\displaystyle\sup_{\lambda\in \mathcal{O}}\|\beta(\lambda,\cdot,\centerdot)\|_{\textnormal{Lip}}<1$ 
then there exists $\widehat\beta: \mathcal{O}\times \T^{2}\to \R$ smooth 
such that
\begin{equation}\label{def betahat}
	y=\theta+\beta(\lambda,\varphi,\theta)\Longleftrightarrow \theta=y+\widehat\beta(\lambda,\varphi,y).
\end{equation}
Define the operators
\begin{equation}\label{definition symplectic change of variables}
	\mathscr{B}=(1+\partial_{\theta}\beta)\mathcal{B}, \qquad \mathcal{B}h(\lambda,\varphi,\theta)=h\big(\lambda,\varphi,\theta+\beta(\lambda,\varphi,\theta)\big).
\end{equation}
By straightforward  computations we obtain, see for instance \cite{HR21},
\begin{equation*}
	\mathscr{B}^{-1}=(1+\partial_{\theta}\widehat\beta)\mathcal{B}^{-1} ,\qquad \mathcal{B}^{-1} h(\lambda,\varphi,y)=h\big(\lambda,\varphi,y+\widehat{\beta}(\lambda,\varphi,y)\big).
\end{equation*}
	We shall now give some elementary algebraic properties for $\mathcal{B}^{\pm 1}$ and $\mathscr{B}^{\pm 1}$ which can be checked by straightforward computations, for more details we refer to \cite{Baldi-Montalto21,BFM,FGMP19}.
			\begin{lem}\label{algeb1}
				The following assertions hold true.
				\begin{enumerate}
					\item Let $\mathscr{B}_1,\mathscr{B}_2$  be two periodic  change of variables as  in \eqref{definition symplectic change of variables}, then
					$$
					\mathscr{B}_{{1}}\mathscr{B}_2=(1+\partial_{\theta}\beta)\mathcal{B},
					$$
					with
					$$
\beta(\varphi,\theta):= \beta_1(\varphi,\theta)+\beta_2\big(\varphi,\theta+\beta_1(\varphi,\theta)\big).
$$					\item
					 The conjugation of the transport operator by $\mathscr{B}$  keeps the same structure
					$$
					\mathscr{B}^{-1}\Big(\omega\cdot\partial_\varphi+\partial_\theta\big(V(\varphi,\theta)\cdot\big)\Big)\mathscr{B}=\omega\cdot\partial_\varphi+\partial_y\big(\mathscr{V}(\varphi,y)\cdot\big),
					$$
					with
					$$
					\mathscr{V}(\varphi,y):=\,\mathcal{B}^{-1}\Big(\omega\cdot\partial_{\varphi} \beta(\varphi,\theta)+V(\varphi,\theta)\big(1+\partial_\theta \beta(\varphi,\theta)\big)\Big).
					$$
				\end{enumerate}
			\end{lem}
			In what follows, and in the rest of this appendix, we assume that $(\lambda,s,s_0)$ satisfy \eqref{cond1} and we consider  $\beta\in \textnormal{Lip}_\nu(\mathcal{O},H^{s}(\T^{2})) $ satisfying the smallness condition 
				\begin{equation*}
					\|\beta \|_{2s_0}^{\textnormal{Lip},\nu}\leqslant \epsilon_0,
				\end{equation*}
			with $\epsilon_{0}$ small enough.\\
			The following result is proved in \cite{FGMP19}. We also refer to \cite[(A.2)]{BFM}.
			\begin{lem}\label{Compos1-lemm}
The following assertions hold true.
				\begin{enumerate}
					\item The linear operators $\mathcal{B},\mathscr{B}:\textnormal{Lip}_\nu(\mathcal{O},H^{s}(\T^{2}))\to \textnormal{Lip}_\nu(\mathcal{O},H^{s}(\T^{2}))$ are continuous and invertible, with 
					\begin{align*}
						 \|\mathcal{B}^{\pm1}h\|_{s}^{\textnormal{Lip},\nu}\leqslant \|h\|_{s}^{\textnormal{Lip},\nu}\left(1+C\|\beta\|_{s_{0}}^{\textnormal{Lip},\nu}\right)+C\|\beta\|_{s}^{\textnormal{Lip},\nu}\|h\|_{s_{0}}^{\textnormal{Lip},\nu}, 
\\
						\|\mathscr{B}^{\pm1}h\|_{s}^{\textnormal{Lip},\nu}\leqslant \|h\|_{s}^{\textnormal{Lip},\nu}\left(1+C\|\beta\|_{s_{0}}^{\textnormal{Lip},\nu}\right)+C\|\beta\|_{s+1}^{\textnormal{Lip},\nu}\|h\|_{s_{0}}^{\textnormal{Lip},\nu}. 
					\end{align*}
					\item The functions $\beta$ and $\widehat{\beta}$ defined through \eqref{def betahat} satisfy the estimates
					\begin{equation*}
					\|\widehat{\beta}\|_{s}^{\textnormal{Lip},\nu}\leqslant C\|\beta\|_{s}^{\textnormal{Lip},\nu}.
					\end{equation*}
				\end{enumerate}
			\end{lem}
			The next result was useful. 
			\begin{lem} \label{beta-inv-asym}
			Let $\varepsilon\in[0,1),$ $f:\mathcal{O}\times\R\to\R$ and $g:\mathcal{O}\times \T^2\to\R$ be two smooth  functions such that
			$$\|g\|_{2s_0}^{\textnormal{Lip},\nu}\leqslant 1,
			$$
			  and $f(\lambda,\cdot)$ or $f(\lambda,\cdot)-\textnormal{Id}$ is $2\pi-$periodic.
			Then the transformation $\Phi(\lambda):\R^2\to\R^2$ defined by
			$$
			\Phi(\lambda,\varphi,\theta)=\big(\varphi,\theta+f(\lambda,\varphi)+ \varepsilon g(\lambda,\varphi,\theta)\big),
			$$
			is a diffeomorphism, with
			$$
			\Phi^{-1}(\lambda,\varphi,\theta)=\big(\varphi,\theta-f(\lambda,\varphi)- \varepsilon g\big(\lambda,\varphi,\theta-f(\lambda,\varphi)\big)+\varepsilon^2 \mathtt{r}(\lambda,\varphi,\theta)\big),
			$$
			and $ \mathtt{r}$ satisfies  the estimate 
			$$
			\forall s\geqslant s_0,\quad \|\mathtt{r}\|_{s}^{\textnormal{Lip},\nu}\leqslant C(s,\|f\|_{s}^{\textnormal{Lip},\nu}) \|g\|_{s+1}^{\textnormal{Lip},\nu}.
			$$
			\end{lem}

\section{Asymptotics for the kernel}\label{sec:appendix}
In this section we derive the asymptotic expansions of the kernel $G$ in the scaling regime
relevant for vortex rings. Our computations follow the approach used for Green function
expansions in \cite{Gallay-Sverak-2}. We first introduce the auxiliary function
\begin{equation*}
J(s):=\int_{0}^{\pi}\frac{\cos\vartheta}{\big(2(1-\cos\vartheta)+s\big)^{1/2}}\,d\vartheta,
\qquad s>0,
\end{equation*}
and we study its behavior as $s\to0$.

The asymptotic expansion of $J$ for small $s$ is established in \cite[Lemma~3.1]{Gallay-Sverak-2}.
For convenience, we restate it below in the notation used throughout this paper.
\begin{lem}\label{prop-split-1}
For all $0<s<4$ we have the power series representation
\begin{align*}
J(s)=&|\ln(s)|\sum_{n\geqslant 0}A_n s^n+\sum_{n\geqslant 0}B_n s^n,
\end{align*}
where $A_m$, $B_m$ are real numbers. Moreover,
\begin{align*}
    A_0=&\frac12,\quad A_1=\frac{3}{32},\quad A_2=-\frac{15}{2048},\quad A_3=\frac{105}{98304},\\
    B_0=&-2+\ln8\quad B_1=-\frac{1}{16}+\frac{3}{16}\ln8,\quad B_2=\frac{31}{2048}-\frac{15}{1024}\ln8.
\end{align*}
\end{lem}

The next lemma provides the corresponding local expansion of $G$ under the vortex-ring scaling. 
It is the analogue of \cite[Lemma~3.3]{Gallay-Sverak-2} in our setting.

\begin{lem}[Local expansion of the kernel]\label{lem:kernel-expansion}
Let $Z=(z_1,z_2)\in\mathbb R^2$ with $z_1>0$, and let $X=(x_1,x_2)$, $Y=(y_1,y_2)\in\mathbb R^2$.
Define
$$
D := (x_1-y_1)^2+2z_1(x_2-y_2)^2,\qquad
E := (x_1-y_1)^2+z_1(x_2-y_2)^2.
$$
Then, for $0<\epsilon\ll1$, the following expansion holds as
$\epsilon\to0$:
$$
G(Z+\epsilon X,\, Z+\epsilon Y)
=
\ln\!\Big(\frac1\epsilon\Big)\sum_{n=0}^3 \mathcal A_n(X,Y,Z)\,\epsilon^n
+\sum_{n=0}^2 \mathcal B_n(X,Y,Z)\,\epsilon^n
+O(\epsilon^3),
$$
where the remainder is uniform for $(X,Y,Z)$ ranging in compact subsets of
$\{(X,Y,Z):\; z_1>0,\; D>0\}$. The coefficients are given by
\begin{align*}
\mathcal A_0(X,Y,Z) &:= \sqrt{z_1},\quad
\mathcal A_1(X,Y,Z) := \frac{x_1+y_1}{4\sqrt{z_1}},\\
\mathcal A_2(X,Y,Z) &:= \frac{1}{64z_1^{3/2}}
\Big(6z_1(x_2-y_2)^2-2x_1y_1-3x_1^2-3y_1^2\Big),\\
\mathcal A_3(X,Y,Z) &:= \frac{x_1+y_1}{256z_1^{5/2}}
\Big(5x_1^2-2x_1y_1+5y_1^2-6z_1(x_2-y_2)^2\Big),
\end{align*}
and
\begin{align*}
\mathcal B_0(X,Y,Z)
&:= \sqrt{z_1}\,(\ln 8-2)-\frac{\sqrt{z_1}}{2}\,\ln\left(\frac{D}{(2z_1)^2}\right),\\[0.4em]
\mathcal B_1(X,Y,Z)
&:= \frac{x_1+y_1}{4\sqrt{z_1}}
\left(\ln 8-2-\tfrac12\ln\left(\frac{D}{(2z_1)^2}\right)\right)
+\frac{x_1+y_1}{2\sqrt{z_1}}\,\frac{E}{D},\\[0.4em]
\mathcal B_2(X,Y,Z)
&:= \frac{1}{16z_1^{3/2}}
\Big(x_1y_1-\tfrac32x_1^2-\tfrac32y_1^2\Big)
\left(\ln 8-2-\tfrac12\ln\left(\frac{D}{(2z_1)^2}\right)\right)
+\frac{(x_1+y_1)^2}{8z_1^{3/2}}\,\frac{E}{D}\\
&\quad+\frac{1}{16}(3\ln 8-1)\sqrt{z_1}\,\frac{D}{(2z_1)^2}
+\frac{(x_1+y_1)^2}{4z_1^{3/2}}\,\left(\frac{E}{D}\right)^2
-\frac{3}{32}\sqrt{z_1}\,\frac{D}{(2z_1)^2}\ln\left(\frac{D}{(2z_1)^2}\right)\\
&\quad-\frac{1}{32z_1^{3/2}D}
\Big(15x_1^4-12x_1^3y_1-6x_1^2y_1^2-12x_1y_1^3+15y_1^4
\\ &\qquad\qquad\qquad\qquad\quad +4z_1(x_2-y_2)^2(3x_1^2+2x_1y_1+3y_1^2)\Big).
\end{align*}
\end{lem}
\begin{proof}
Set
$$
\mathcal P(\epsilon):=(z_1+\epsilon x_1)^{\frac14}(z_1+\epsilon y_1)^{\frac14},
\qquad
s(\epsilon):=\frac{2\big(\sqrt{z_1+\epsilon x_1}-\sqrt{z_1+\epsilon y_1}\big)^2+\epsilon^2(x_2-y_2)^2}
{2\sqrt{z_1+\epsilon x_1}\sqrt{z_1+\epsilon y_1}}.
$$
Then
$$
G(Z+\epsilon X,\,Z+\epsilon Y)=\mathcal P(\epsilon)\,J\big(s(\epsilon)\big).
$$
Fix compact subsets of $\{(X,Y,Z): z_1>0,\ D>0\}$, where
$$
D=(x_1-y_1)^2+2z_1(x_2-y_2)^2,\qquad E=(x_1-y_1)^2+z_1(x_2-y_2)^2.
$$
For $\epsilon>0$ sufficiently small we have $z_1+\epsilon x_1>0$, $z_1+\epsilon y_1>0$ and $s(\epsilon)\sim \epsilon^2 D/(4z_1^2)$, hence $s(\epsilon)\to0^+$.

\medskip

Using the binomial series $(1+t)^{1/4}=1+\frac14t-\frac{3}{32}t^2+\frac{7}{128}t^3+O(t^4)$, we obtain
$$
\mathcal P(\epsilon)
=\sqrt{z_1}+\epsilon\,\alpha_1+\epsilon^2\,\alpha_2+\epsilon^3\,\alpha_3+O(\epsilon^4),
$$
where
\begin{align*}
\alpha_1&=\frac{x_1+y_1}{4\sqrt{z_1}},\quad
\alpha_2=\frac{1}{z_1^{3/2}}\Big(\frac{1}{16}x_1y_1-\frac{3}{32}x_1^2-\frac{3}{32}y_1^2\Big),\\
\alpha_3&=\frac{1}{128z_1^{5/2}}\Big(7x_1^3 -3x_1^2y_1-3x_1y_1^2+7y_1^3\Big).
\end{align*}
On the other hand, a Taylor expansion of the square-root terms yields
$$
s(\epsilon)=\epsilon^2P_0\Big(1+\epsilon P_1+\epsilon^2P_2+O(\epsilon^3)\Big),
$$
with
$$
P_0=\frac{D}{(2z_1)^2},\qquad
P_1=-\frac{x_1+y_1}{z_1}\frac{E}{D},
$$
and
$$
P_2=\frac{1}{16z_1^2}\frac{1}{D}\Big(15x_1^4-12x_1^3y_1-6x_1^2y_1^2-12x_1y_1^3+15y_1^4
+4z_1(x_2-y_2)^2(3x_1^2+2x_1y_1+3y_1^2)\Big).
$$
Consequently, we have
\begin{align*}
\ln s(\epsilon)
&=\ln\!\big(\epsilon^2P_0\big)+\ln\!\Big(1+\epsilon P_1+\epsilon^2P_2+O(\epsilon^3)\Big)\\
&=-2\ln(1/\epsilon)+\ln P_0+\epsilon P_1+\epsilon^2\Big(P_2-\frac12P_1^2\Big)+O(\epsilon^3),
\end{align*}
and since $s(\epsilon)\to0^+$,
$$
|\ln s(\epsilon)|
=2\ln(1/\epsilon)-\ln P_0-\epsilon P_1-\epsilon^2\Big(P_2-\frac12P_1^2\Big)+O(\epsilon^3).
$$
By Lemma~\ref{prop-split-1}, and using $s(\epsilon)=O(\epsilon^2)$, we may truncate
$$
J(s)=A_0|\ln s|+B_0+s\big(A_1|\ln s|+B_1\big)+O(s^2|\ln s|),
$$
as $s\to0^+$, where $A_0=\frac12$, $A_1=\frac{3}{32}$, $B_0=\ln 8-2$ and
$B_1=-\frac{1}{16}+\frac{3}{16}\ln 8$.
Since $s(\epsilon)^2|\ln s(\epsilon)|=O(\epsilon^4L)=o(\epsilon^3)$ as $\epsilon\to0^+$, this remainder is admissible for our purposes.
Substituting the expansions of $|\ln s(\epsilon)|$ and $s(\epsilon)$ gives
$$
J\big(s(\epsilon)\big)
=
\ln(1/\epsilon)+C_0+\epsilon C_1+\epsilon^2\Big(\frac{3}{16}P_0\,\ln(1/\epsilon)+C_2\Big)
+\epsilon^3\Big(\frac{3}{16}P_0P_1\,\ln(1/\epsilon)\Big)+O(\epsilon^3),
$$
where
\begin{align*}
C_0&:=\ln 8-2-\frac12\ln P_0,\quad
C_1:=-\frac12P_1,\\
C_2&:=\frac{1}{16}(3\ln 8-1)P_0-\frac12\Big(P_2-\frac12P_1^2\Big)-\frac{3}{32}P_0\ln P_0.
\end{align*}
Multiplying the expansions of $\mathcal P(\epsilon)$ and $J(s(\epsilon))$ and collecting powers of
$\epsilon$ and $\ln(1/\epsilon)$ yields
$$
G(Z+\epsilon X,\,Z+\epsilon Y)
=
\ln(1/\epsilon)\sum_{n=0}^3\mathcal A_n\,\epsilon^n+\sum_{n=0}^2\mathcal B_n\,\epsilon^n+O(\epsilon^3),
$$
with
\begin{align*}
\mathcal A_0&=\sqrt{z_1},\quad
\mathcal A_1=\alpha_1,\quad
\mathcal A_2=\alpha_2+\sqrt{z_1}\,\frac{3}{16}P_0,\\
\mathcal A_3&=\alpha_3+\alpha_1\frac{3}{16}P_0+\sqrt{z_1}\,\frac{3}{16}P_0P_1,
\end{align*}
and
\begin{align*}
\mathcal B_0&=\sqrt{z_1}\,C_0,\quad
\mathcal B_1=\alpha_1C_0+\sqrt{z_1}\,C_1,\\
\mathcal B_2&=\alpha_2C_0+\alpha_1C_1+\sqrt{z_1}\,C_2.
\end{align*}
Finally, substituting $P_0=\frac{D}{(2z_1)^2}$, $P_1=-\frac{x_1+y_1}{z_1}\frac{E}{D}$,
and the explicit expressions for $\alpha_1,\alpha_2,\alpha_3$, $C_0,C_1,C_2$, a straightforward algebraic simplification
 gives the formulas stated in Lemma~\ref{lem:kernel-expansion}.
\end{proof}
\begin{lem}\label{prop:asympt-induced}
Let $Z=(z_1,z_2)$ with $z_1>0$. For $X=(x_1,x_2)$, $Y=(y_1,y_2)$ and for $\epsilon>0$
sufficiently small, one has
\begin{align*}
G\Big(&z_1+\epsilon (2z_1)^{\frac14}x_1,\ z_2+\epsilon (2z_1)^{-\frac14}x_2,\ 
      z_1+\epsilon (2z_1)^{\frac14}y_1,\ z_2+\epsilon (2z_1)^{-\frac14}y_2\Big) \\
&=|\ln(\epsilon)|\sum_{n=0}^3\mathscr{A}_n(X,Y)\,\epsilon^n
+\sum_{n=0}^2\mathscr{B}_n(X,Y)\,\epsilon^n
+O(\epsilon^3),
\end{align*}
where  the remainder is uniform for $(X,Y)$ in
compact subsets of $\{(X,Y):X\neq Y\}$. The coefficients are
\begin{align*}
\mathscr{A}_0(X,Y)&:=\sqrt{z_1},\quad
\mathscr{A}_1(X,Y):=\frac{2^{1/4}}{4 z_1^{1/4}}(x_1+y_1),\\ 
\mathscr{A}_2(X,Y)&:=\frac{2^{1/2}}{64z_1}\left(3(x_2-y_2)^2 -2x_1y_1-3x_1^2-3y_1^2\right),\\
\mathscr{A}_3(X,Y)&:=\frac{2^{3/4}}{256 z_1^{7/4}}(x_1 + y_1)\left(5 x_1^2 -2 x_1 y_1  + 5y_1^2 -3( x_2 -y_2)^2\right),
\end{align*}
and
\begin{align*}
\mathscr{B}_0(X,Y)&:=\sqrt{z_1}\big(\ln8-2\big)-\frac{\sqrt{z_1}}{2}\ln\left({|X-Y|^2}\right)+\frac{3\sqrt{z_1}}{4}\ln(2z_1),\\
\mathscr{B}_1(X,Y)&:=\frac{2^{1/4}}{4 z_1^{1/4}}(x_1+y_1)\left(\frac54\ln8+\frac34\ln(z_1)-1+\frac{(x_1-y_1)^2}{|X-Y|^2}-\frac12\ln\big(|X-Y|^2\big)\right),\\
\mathscr{B}_2(X,Y)&:=\frac{\sqrt{2}}{16z_1} \left(x_1y_1-\frac{3}{2}x_1^2-\frac{3}{2}y_1^2\right)
\left(\ln8-2-\frac12\ln\left(|X-Y|^2\right)+\frac34\ln(2z_1)\right)\\
&\quad +\frac{\sqrt{2}}{8z_1}(x_1+y_1)^2\frac{(x_1-y_1)^2+2^{-1}(x_2-y_2)^2}{|X-Y|^2}+ \frac{\sqrt{2}}{64 z_1}(3\ln8-1)|X-Y|^2 \\
&\quad
+\frac{\sqrt{2}}{4z_1}(x_1+y_1)^2\Big(\frac{(x_1-y_1)^2+2^{-1}(x_2-y_2)^2}{|X-Y|^2}\Big)^2
 - \frac{ 3\sqrt{2}}{128 z_1}|X-Y|^2 \ln\left(\frac{|X-Y|^2}{(2z_1)^{3/2}}\right)\\
&\quad
-\frac{\sqrt{2}}{32 z_1}\frac{1}{|X-Y|^2}\Big[15 x_1^4 +15 y_1^4- 12 x_1^3 y_1-12 x_1 y_1^3- 6 x_1^2 y_1^2
\\ &\qquad\qquad \qquad\qquad\qquad\qquad\qquad+\big(6 x_1^2  +4 x_1 y_1  +  6 y_1^2\big)(x_2 - y_2)^2\Big].
\end{align*}
\end{lem}
\begin{proof}
Define the rescaled vectors
$$
\widetilde X:=\big((2z_1)^{1/4}x_1,\ (2z_1)^{-1/4}x_2\big),
\qquad
\widetilde Y:=\big((2z_1)^{1/4}y_1,\ (2z_1)^{-1/4}y_2\big).
$$
Then the left-hand side of the lemma is exactly $G(Z+\epsilon \widetilde X,\ Z+\epsilon \widetilde Y)$.
Applying Lemma~\ref{lem:kernel-expansion}, yields
$$
G(Z+\epsilon \widetilde X,\ Z+\epsilon \widetilde Y)
=
|\ln\epsilon|\sum_{n=0}^3 \mathcal A_n(\widetilde X,\widetilde Y,Z)\,\epsilon^n
+\sum_{n=0}^2 \mathcal B_n(\widetilde X,\widetilde Y,Z)\,\epsilon^n
+O(\epsilon^3).
$$
We set
$$
\mathscr A_n(X,Y):=\mathcal A_n(\widetilde X,\widetilde Y,Z),
\qquad
\mathscr B_n(X,Y):=\mathcal B_n(\widetilde X,\widetilde Y,Z),
$$
and compute these coefficients explicitly.

First note that
$$
(\widetilde x_1-\widetilde y_1)^2=\sqrt{2z_1}\,(x_1-y_1)^2,\qquad
(\widetilde x_2-\widetilde y_2)^2=(2z_1)^{-1/2}(x_2-y_2)^2,
$$
hence
$$
(\widetilde x_1-\widetilde y_1)^2+2z_1(\widetilde x_2-\widetilde y_2)^2
=\sqrt{2z_1}\big((x_1-y_1)^2+(x_2-y_2)^2\big)=\sqrt{2z_1}\,|X-Y|^2.
$$
Therefore the logarithmic quantity in Lemma~\ref{lem:kernel-expansion} satisfies
$$
\frac{D}{(2z_1)^2}=\frac{(\widetilde x_1-\widetilde y_1)^2+2z_1(\widetilde x_2-\widetilde y_2)^2}{(2z_1)^2}
=\frac{|X-Y|^2}{(2z_1)^{3/2}},
$$
and the ratio
$$
\frac{E}{D}=\frac{(\widetilde x_1-\widetilde y_1)^2+z_1(\widetilde x_2-\widetilde y_2)^2}{(\widetilde x_1-\widetilde y_1)^2+2z_1(\widetilde x_2-\widetilde y_2)^2}
=\frac{(x_1-y_1)^2+\tfrac12(x_2-y_2)^2}{|X-Y|^2}.
$$
Substituting $\widetilde X,\widetilde Y$ into the explicit formulas for $\mathcal A_n$ and $\mathcal B_n$ and
using the identities above gives exactly the expressions stated for $\mathscr A_n$ and $\mathscr B_n$.
The $O(\epsilon^3)$ remainder is inherited from Lemma~\ref{lem:kernel-expansion} and is uniform on compact
sets away from the diagonal $X=Y$.
\end{proof}
As a consequence of the previous lemma, we obtain the following asymptotic expansion for the gradient of the kernel.
\begin{lem}\label{pro-decomp2-nabla}
Let $Z=(z_1,z_2)$ with $z_1>0$ and let $X=(x_1,x_2)$, $Y=(y_1,y_2)$.
For $\epsilon>0$ sufficiently small, one has
\begin{align*}
(\nabla_1 G)\Big(&z_1+\epsilon (2z_1)^{\frac14}x_1,\ z_2+\epsilon (2z_1)^{-\frac14}x_2,\ 
z_1+\epsilon (2z_1)^{\frac14}y_1,\ z_2+\epsilon (2z_1)^{-\frac14}y_2\Big)\\
&=|\ln(\epsilon)|\sum_{n=1}^3 \epsilon^{n-1}\mathscr{C}_n(X,Y,Z)
+\sum_{n=0}^2\epsilon^{n-1} \mathscr{D}_n(X,Y,Z)
+O(\epsilon^2),
\end{align*}
and, by symmetry,
\begin{align*}
(\nabla_2 G)\Big(&z_1+\epsilon (2z_1)^{\frac14}x_1,\ z_2+\epsilon (2z_1)^{-\frac14}x_2,\ 
z_1+\epsilon (2z_1)^{\frac14}y_1,\ z_2+\epsilon (2z_1)^{-\frac14}y_2\Big)\\
&=|\ln(\epsilon)|\sum_{n=1}^3 \epsilon^{n-1}\mathscr{C}_n(Y,X,Z)
+\sum_{n=0}^2\epsilon^{n-1} \mathscr{D}_n(Y,X,Z)
+O(\epsilon^2).
\end{align*}
The coefficients $\mathscr C_n$ are
\begin{align*}
\mathscr{C}_1(X,Y,Z)&:=\frac{1}{4 \sqrt{z_1}}(1,0),\quad 
\mathscr{C}_2(X,Y,Z):=\frac{1}{32z_1^{3/2}}\left( -(2z_1)^{1/4}(y_1+3x_1),\ 3 (2z_1)^{3/4} (x_2-y_2)\right),\\
\mathscr{C}_3(X,Y,Z)&:=\frac{\sqrt{2}}{256 z_1^{2}}\left(15 x_1^2+6 x_1y_1+3 y_1^2-3(x_2-y_2)^2,\ -6(2z_1)^{1/2}(x_1+y_1)(x_2-y_2)\right),
\end{align*}
and $\mathscr D_0,\mathscr D_1$ are
\begin{align*}
\mathscr{D}_0(X,Y,Z)&:=-\sqrt{z_1}\Big(\frac{(2z_1)^{-1/4}(x_1-y_1)}{|X-Y|^2},\ \frac{(2z_1)^{1/4}(x_2-y_2)}{|X-Y|^2}\Big),\\
\mathscr{D}_1(X,Y,Z)&:=\frac{1}{4 \sqrt{z_1}}\bigg(\tfrac54\ln8+\tfrac34\ln(z_1)-1-\tfrac12\ln\big(|X-Y|^2\big)\\
&\qquad\qquad\qquad\quad+ \frac{(x_1-y_1)^2+(x_1+y_1)(x_1-y_1)}{|X-Y|^2}
- \frac{2( x_1+y_1)( x_1-y_1)^3}{|X-Y|^4},\\
&\qquad\qquad\qquad\quad -\frac{\sqrt{2z_1}(x_1+y_1)(x_2-y_2)}{|X-Y|^2}
- \frac{2\sqrt{2z_1}( x_1+y_1)( x_1-y_1)^2(x_2-y_2)}{|X-Y|^4}\bigg).
\end{align*}
Finally, $\mathscr D_2$ is given by
$$
\mathscr D_2(X,Y,Z):=\mathrm{diag}\big((2z_1)^{-1/4},(2z_1)^{1/4}\big)\,\nabla_X \mathscr B_2(X,Y),
$$
where $\mathscr B_2$ is the coefficient from Lemma~\ref{prop:asympt-induced}.
\end{lem}
\begin{proof}
Define the scaled points
$$
p_\epsilon(X):=\big(z_1+\epsilon (2z_1)^{1/4}x_1,\ z_2+\epsilon (2z_1)^{-1/4}x_2\big),
\qquad
q_\epsilon(Y):=\big(z_1+\epsilon (2z_1)^{1/4}y_1,\ z_2+\epsilon (2z_1)^{-1/4}y_2\big),
$$
and set
$$
F_\epsilon(X,Y):=G\big(p_\epsilon(X),\,q_\epsilon(Y)\big).
$$
By Lemma~\ref{prop:asympt-induced},
$$
F_\epsilon(X,Y)=|\ln\epsilon|\sum_{n=0}^3 \epsilon^n \mathscr A_n(X,Y)
+\sum_{n=0}^2 \epsilon^n \mathscr B_n(X,Y)
+O(\epsilon^3).
$$
Differentiating $F_\epsilon$ with respect to $X$ gives
$$
\partial_{x_1}F_\epsilon(X,Y)=\epsilon(2z_1)^{1/4}\,(\partial_\rho G)\big(p_\epsilon(X),q_\epsilon(Y)\big),\quad
\partial_{x_2}F_\epsilon(X,Y)=\epsilon(2z_1)^{-1/4}\,(\partial_z G)\big(p_\epsilon(X),q_\epsilon(Y)\big).
$$
Equivalently,
$$
(\nabla_1 G)\big(p_\epsilon(X),q_\epsilon(Y)\big)
=\epsilon^{-1}M\,\nabla_XF_\epsilon(X,Y),
\qquad
M:=\mathrm{diag}\big((2z_1)^{-1/4},(2z_1)^{1/4}\big).
$$
Since the coefficients $\mathscr A_n,\mathscr B_n$ are smooth away from $X=Y$ and the remainder is uniform on
compact sets with $X\neq Y$, we may differentiate the expansion termwise:
$$
\nabla_XF_\epsilon(X,Y)
=|\ln\epsilon|\sum_{n=0}^3 \epsilon^n \nabla_X\mathscr A_n(X,Y)
+\sum_{n=0}^2 \epsilon^n \nabla_X\mathscr B_n(X,Y)
+O(\epsilon^3).
$$
Multiplying by $\epsilon^{-1}M$ yields
$$
(\nabla_1 G)\big(p_\epsilon(X),q_\epsilon(Y)\big)
=|\ln\epsilon|\sum_{n=1}^3 \epsilon^{n-1} \underbrace{M\nabla_X\mathscr A_n}_{=:~\mathscr C_n}
+\sum_{n=0}^2 \epsilon^{n-1} \underbrace{M\nabla_X\mathscr B_n}_{=:~\mathscr D_n}
+O(\epsilon^2).
$$
The explicit formulas for $\mathscr C_1,\mathscr C_2,\mathscr C_3$ follow by direct differentiation of
$\mathscr A_1,\mathscr A_2,\mathscr A_3$ (note that $\nabla_X\mathscr A_0\equiv 0$), and similarly
$\mathscr D_0,\mathscr D_1$ follow by differentiating $\mathscr B_0,\mathscr B_1$.
The coefficient $\mathscr D_2$ is $M\nabla_X\mathscr B_2$.

\medskip

Using symmetry of the kernel, $G(p,q)=G(q,p)$, we have
$$
(\nabla_2G)(p,q)=(\nabla_1G)(q,p),
$$
which gives the second expansion by swapping $X$ and $Y$.
\end{proof}

\begin{lem}\label{prop-asym-derivatives-G}
Given the scaling
$$
p_{1,1}=\kappa+\frac12 r_\varepsilon (2\kappa)^\frac14 \mathtt{x}_1,\quad p_{2,1}=\kappa-\frac12 r_\varepsilon (2\kappa)^\frac14\mathtt{x}_1,\quad p_{2,2}-p_{1,2}=r_\varepsilon (2\kappa)^{-\frac14}\mathtt{x}_2, \quad r_\varepsilon=\tfrac{(2\kappa)^\frac14}{|\ln\varepsilon|^\frac12},
$$
we have the asymptotics 
\begin{align*}
|\ln\varepsilon|^{-1}(\partial_{p_{1,1}}^2G)(P_1,P_2)&=\tfrac{1}{2\sqrt{\kappa}}\tfrac{\mathtt{x}_1^2-\mathtt{x}_2^2}{(\mathtt{x}_1^2+\mathtt{x}_2^2)^2}+O(|\ln\varepsilon|^{-\frac12}),\\
|\ln\varepsilon|^{-1}(\partial_{p_{1,2}}^2G)(P_1,P_2)&=-{\sqrt{\kappa}}\tfrac{\mathtt{x}_1^2-\mathtt{x}_2^2}{(\mathtt{x}_1^2+\mathtt{x}_2^2)^2}+O(|\ln\varepsilon|^{-\frac12}),\\
|\ln\varepsilon|^{-1}(\partial_{p_{1,1}}\partial_{p_{1,2}}G)(P_1,P_2)&={\sqrt{2}}\tfrac{\mathtt{x}_1\mathtt{x}_2}{(\mathtt{x}_1^2+\mathtt{x}_2^2)^2}+O(|\ln\varepsilon|^{-\frac12}).
\end{align*}
\end{lem}

\begin{proof}
    For the proof, we use Lemma \ref{prop:asympt-induced} with parameters
    \begin{equation}\label{parameters-0} 
    z_1=\kappa,\quad z_2=0, \quad \epsilon=\tfrac12 r_\varepsilon, \quad x_1=\mathtt{x}_1, \quad y_1=-\mathtt{x}_1, \quad x_2=\mathtt{x}_2,\quad y_2=-\mathtt{x}_2,
    \end{equation}
    and the chain rule. For example, note that
 \begin{align*}
    \partial_{x_1} \Big({G}\Big(z_1&+\epsilon (2z_1)^{\frac14}x_1, z_2+\epsilon (2z_1)^{-\frac14}x_2,z_1+\epsilon (2z_1)^{\frac14}y_1, z_2+\epsilon (2z_1)^{-\frac14}y_2\Big)\Big)\\
    =\epsilon (2z_1)^\frac14(\partial_{p_{1,1}} G)&\Big(z_1+\epsilon (2z_1)^{\frac14}x_1, z_2+\epsilon (2z_1)^{-\frac14}x_2,z_1+\epsilon (2z_1)^{\frac14}y_1, z_2+\epsilon (2z_1)^{-\frac14}y_2\Big),
    \end{align*}
    and then
     \begin{align*}
    \partial^2_{x_1} \Big({G}&\Big(z_1+\epsilon (2z_1)^{\frac14}x_1, z_2+\epsilon (2z_1)^{-\frac14}x_2,z_1+\epsilon (2z_1)^{\frac14}y_1, z_2+\epsilon (2z_1)^{-\frac14}y_2\Big)\Big)\\
    =\epsilon^2 (2z_1)^\frac12(\partial^2_{p_{1,1}} G)&\Big(z_1+\epsilon (2z_1)^{\frac14}x_1, z_2+\epsilon (2z_1)^{-\frac14}x_2,z_1+\epsilon (2z_1)^{\frac14}y_1, z_2+\epsilon (2z_1)^{-\frac14}y_2\Big)\Big).
    \end{align*}
    Now we use Lemma \ref{prop:asympt-induced} to compute $\partial^2_{x_1} G$ asymptotically in $\epsilon$ small using
  \begin{align*}
\partial^2_{x_1}\Big({G}\Big(z_1+&\epsilon (2z_1)^{\frac14}x_1, z_2+\epsilon (2z_1)^{-\frac14}x_2,z_1+\epsilon (2z_1)^{\frac14}y_1, z_2+\epsilon (2z_1)^{-\frac14}y_2\Big)\Big)\\ =&|\ln(\epsilon)|\sum_{n=0}^3\partial^2_{x_1}\mathscr{A}_n(X,Y) \epsilon^n+\sum_{n=0}^2\partial^2_{x_1} \mathscr{B}_n(X,Y) \epsilon^n+O(\epsilon^3)\\
=&\frac{\sqrt{z_1}}{2}\frac{4z_1(x_2-y_2)^2-2(x_1-y_1)^2}{[(x_1-y_1)^2+2z_1(x_2-y_2)^2]^2}+O(\epsilon).
\end{align*}  
Now, adding the values of the parameters \eqref{parameters-0} and dividing by $\epsilon^2(2z_1)^\frac12$ we get the first asymptotic expansion. The following expansions follow the same idea.
\end{proof}
\section{Fourier series expansion of the tangential part}\label{appendix:Fourier}

In the following lemma, we give more details about the operator $\mathcal{S}$ and $\mathcal{Q}$ defined in \eqref{def-S} and \eqref{operator-mode1}, respectively.

\begin{lem}\label{prop:SandQ}
The operator $\mathcal{S}$ acts on the first Fourier mode as follows:
\begin{align*}
\partial_\theta\mathcal{S}[\cos(\theta)]&= 
\mathtt{g}_{3}(\varphi)\sin(2\theta),\quad 
\partial_\theta\mathcal{S}[\sin(\theta)]= 
\mathtt{g}_{3}(\varphi)\cos(2\theta).
\end{align*}
Moreover, its projection onto the first Fourier mode satisfies
\begin{align*}
\Pi_{1,\mathbf{s}}\partial_\theta\mathcal{S}[h(\theta)]&=\partial_\theta\mathcal{S}[\cos(2\theta)]= 
\tfrac{1}{2}  \mathtt{g}_{3}(\varphi)\sin(\theta),\\
\Pi_{1,\mathbf{c}}\partial_\theta\mathcal{S}[h(\theta)]&=\partial_\theta\mathcal{S}[\sin(2\theta)]= 
\tfrac{1}{2} \mathtt{g}_{3}(\varphi)\cos(\theta).
\end{align*}
The operator $\mathcal{H}_{\mathtt{u},0}$ acts on the first a Fourier mode according to
\begin{align*}
    \mathcal{H}_{\mathtt{u},0}[\cos(\theta)]&=-3\mathtt{g}_{3}(\varphi)\sin(2\theta),
    \quad
   \mathcal{H}_{\mathtt{u},0}[\sin(\theta)]=3\mathtt{g}_{3}(\varphi)\cos(2\theta), 
\end{align*}
and its projection onto the first mode is given by
\begin{align*}
  \Pi_{1,\mathbf{s}}\mathcal{H}_{\mathtt{u},0}[h(\theta)]&=    \Pi_{1,\mathbf{s}} \mathcal{H}_{\mathtt{u},0}[\cos(2\theta)]=-\tfrac{3}{2} \mathtt{g}_{3}(\varphi) \sin(\theta), \\
   \Pi_{1,\mathbf{c}} \mathcal{H}_{\mathtt{u},0}[h(\theta)]&=\Pi_{1,\mathbf{c}} \mathcal{H}_{\mathtt{u},0}[\sin(2\theta)]=\tfrac{3}{2} \mathtt{g}_{3}(\varphi) \cos(\theta).
\end{align*}
Moreover, the operator $\mathcal{Q}[h]$ localizes on first Fourier mode and writes
\begin{align*}
       \partial_\theta\mathcal{Q}[h](\varphi,\theta)=&\big(\tfrac{1}{16 }(2\mathtt{p}_{1,1})^{-\frac32} h_{1,c}+\mathtt{q}_1(\varphi)h_{1,c,\star} +\mathtt{q}_2(\varphi)h_{1,s,\star} \big)\sin(\theta)\\ &+\big(-\tfrac{3}{16 }(2\mathtt{p}_{1,1})^{-\frac32} h_{1,s}+\mathtt{q}_3(\varphi)h_{1,c,\star} +\mathtt{q}_4(\varphi)h_{1,s,\star}\big)\cos(\theta).
  \end{align*}
  The auxilliary functions and their expansions are given in Lemma \ref{lem-functions}.
\end{lem}

\begin{proof}
 
By \eqref{operator-mode1},  one has
\begin{align*}
       \partial_\theta\mathcal{Q}[h]&(\varphi,\theta)=-\tfrac{1}{16\pi }(2\mathtt{p}_{1,1})^{-\frac32}\partial_\theta\left( \cos\theta\int_{0}^{2\pi} h(\eta) \cos\eta d\eta+3\sin\theta\int_{0}^{2\pi} h(\eta) \sin\eta d\eta\right)\\ &+\tfrac{1}{2\pi\sqrt{\mathtt{p}_{1,1}}}(2\mathtt{p}_{1,1})^{\frac14}(2\mathtt{p}_{2,1})^{\frac14}|\ln\varepsilon|^{-1}(\partial_{\mathtt{p}_{1,1}}\partial_{\mathtt{p}_{2,1}}G)(\mathtt{P}_1,\mathtt{P}_2)\partial_\theta\left(\cos(\theta)\int_{0}^{2\pi} h_\star(\eta) \cos(\eta)d\eta\right)
      \\ &+\tfrac{1}{2\pi\sqrt{\mathtt{p}_{1,1}}}(2\mathtt{p}_{1,1})^{\frac14}(2\mathtt{p}_{2,1})^{-\frac14}|\ln\varepsilon|^{-1}(\partial_{\mathtt{p}_{1,1}}\partial_{\mathtt{p}_{2,2}}G)(\mathtt{P}_1,\mathtt{P}_2)\partial_\theta\left(\cos(\theta)\int_{0}^{2\pi} h_\star(\eta) \sin(\eta)d\eta\right)
      \\ &+\tfrac{1}{2\pi\sqrt{\mathtt{p}_{1,1}}}(2\mathtt{p}_{1,1})^{-\frac14}(2\mathtt{p}_{2,1})^{\frac14}|\ln\varepsilon|^{-1}(\partial_{\mathtt{p}_{1,2}}\partial_{\mathtt{p}_{2,1}}G)(\mathtt{P}_1,\mathtt{P}_2)\partial_\theta\left(\sin(\theta)\int_{0}^{2\pi} h_\star(\eta)\cos(\eta)d\eta\right)
      \\ &+\tfrac{1}{2\pi\sqrt{\mathtt{p}_{1,1}}}(2\mathtt{p}_{1,1})^{-\frac14}(2\mathtt{p}_{2,1})^{-\frac14}|\ln\varepsilon|^{-1}(\partial_{\mathtt{p}_{1,2}}\partial_{\mathtt{p}_{2,2}}G)(\mathtt{P}_1,\mathtt{P}_2)\partial_\theta\left(\sin(\theta)\int_{0}^{2\pi} h_\star(\eta) \sin(\eta)d\eta\right).
  \end{align*}
Using the decomposition
$$
h=\sum_{k\geqslant 1} h_{k,c}\cos(k\theta)+h_{k,s}\sin(k\theta),
$$
we get
\begin{align*}
       \partial_\theta\mathcal{Q}[h](\varphi,\theta)=&\tfrac{1}{16 }(2\mathtt{p}_{1,1})^{-\frac32}\left(h_{1,c} \sin(\theta)-3 h_{1,s}\cos(\theta)\right)
       \\ &-\tfrac{1}{\sqrt{2}}(2\mathtt{p}_{1,1})^{-\frac14}(2\mathtt{p}_{2,1})^{\frac14}|\ln\varepsilon|^{-1}(\partial_{\mathtt{p}_{1,1}}\partial_{\mathtt{p}_{2,1}}G)(\mathtt{P}_1,\mathtt{P}_2)h_{1,c,\star}\sin(\theta) 
      \\ &-\tfrac{1}{\sqrt{2}}(2\mathtt{p}_{1,1})^{-\frac14}(2\mathtt{p}_{2,1})^{-\frac14}|\ln\varepsilon|^{-1}(\partial_{\mathtt{p}_{1,1}}\partial_{\mathtt{p}_{2,2}}G)(\mathtt{P}_1,\mathtt{P}_2)h_{1,s,\star}\sin(\theta) 
      \\ &+\tfrac{1}{\sqrt{2}}(2\mathtt{p}_{1,1})^{-\frac34}(2\mathtt{p}_{2,1})^{\frac14}|\ln\varepsilon|^{-1}(\partial_{\mathtt{p}_{1,2}}\partial_{\mathtt{p}_{2,1}}G)(\mathtt{P}_1,\mathtt{P}_2)h_{1,c,\star}\cos(\theta)
      \\ &+\tfrac{1}{\sqrt{2}}(2\mathtt{p}_{1,1})^{-\frac34}(2\mathtt{p}_{2,1})^{-\frac14}|\ln\varepsilon|^{-1}(\partial_{\mathtt{p}_{1,2}}\partial_{\mathtt{p}_{2,2}}G)(\mathtt{P}_1,\mathtt{P}_2)h_{1,s,\star}\cos(\theta).
  \end{align*}

Concerning the operator $\mathcal{S}$  in \eqref{shift-operator1}, its explicit action on the different modes can be obtained directly from its definition. Finally, to evaluate the action of the operator $\mathcal{H}_{\mathtt{u},0}$ in \eqref{def-Hu0} on the first and second Fourier modes, we use the identity
\begin{align*}
    \fint_{\T} \ln\!\big|\sin ({\tfrac{\eta}{2}})\big|\, e^{in\eta}\, d\eta=-\tfrac{1}{2|n|},
\end{align*}
which, after a change of variables, yields
\begin{align*}
   \fint_{\T} \ln\!\Big|\sin\!\Big(\tfrac{\theta-\eta}{2}\Big)\Big|\,\cos(n\eta)\, d\eta&=-\tfrac{1}{2|n|}\cos(n\theta),\\
    \fint_{\T} \ln\!\Big|\sin\!\Big(\tfrac{\theta-\eta}{2}\Big)\Big|\,\sin(n\eta)\, d\eta&=-\tfrac{1}{2|n|}\sin(n\theta).
\end{align*}
Using the previous identities, together with the trigonometric relation $2\cos^2(\eta)=1+\cos(2\eta)$ we obtain
\begin{align*}
   \mathcal{H}_{\mathtt{u},0}[\cos(\theta)]&=-\tfrac12(2\mathtt{p}_{1,1})^{-\frac34}\partial_\theta\left( \fint_{\T} \ln\!\Big|\sin\!\Big(\tfrac{\theta-\eta}{2}\Big)\Big|\, \Big(\cos(\theta) \cos(\eta)+\tfrac12+\tfrac12\cos(2\eta)\Big)\, d\eta\right)\\
   &=-\tfrac12(2\mathtt{p}_{1,1})^{-\frac34}\partial_\theta\left( -\tfrac12\cos^2(\theta) -\tfrac18\cos(2\theta)\right)\\
   &=-\tfrac38(2\mathtt{p}_{1,1})^{-\frac34} \sin(2\theta).
\end{align*}
Similarly, using the trigonometric identity $2\cos(\eta)\cos(2\eta)=\cos(\eta)+\cos(3\eta)$, we find
\begin{align*}
   \mathcal{H}_{\mathtt{u},0}[\cos(2\theta)]&=-\tfrac12(2\mathtt{p}_{1,1})^{-\frac34}\partial_\theta \fint_{\T} \ln\!\Big|\sin\!\Big(\tfrac{\theta-\eta}{2}\Big)\Big|\,\left(\cos(\theta) \cos(2\eta)+\tfrac12\cos(\eta)+\tfrac12\cos(3\eta)\right)\, d\eta\\
   &=-\tfrac12(2\mathtt{p}_{1,1})^{-\frac34}\partial_\theta\big(-\tfrac14 \cos(\theta) \cos(2\theta)-\tfrac14\cos(\theta)-\tfrac{1}{12}\cos(3\theta)\big)
   \\
   &=-\tfrac{1}{16}(2\mathtt{p}_{1,1})^{-\frac34}\big(3\sin(\theta)+5\sin(3\theta)\big).
\end{align*}
Proceeding analogously, we obtain
\begin{align*}
   \mathcal{H}_{\mathtt{u},0}[\sin(\theta)]
   &=\tfrac38(2\mathtt{p}_{1,1})^{-\frac34}\cos(2\theta),\\
    \mathcal{H}_{\mathtt{u},0}[\sin(2\theta)]&=\tfrac{1}{16}(2\mathtt{p}_{1,1})^{-\frac34}\big(3\cos(\theta)+5\cos(3\theta)\big).
\end{align*}
This concludes the proof of the lemma.
\end{proof}

In the following lemma we compute explicitly the operators \eqref{defL1Q1L2Q2}, arising from the speed modulation, expressed in Fourier variables.

\begin{lem}\label{lem-FourLQ}
Let $$
h(\theta)=\sum_{k\geqslant 1} \big(h_{k,c}\cos(k\theta)+h_{k,s}\sin(k\theta)\big).
$$ 
Then,
\begin{align*}
   \mathtt{L}_1[h]&=\tfrac{1}{32\mathtt{p}_{1,1}}(h_{4,c}-3h_{2,c}),\\
   \mathtt{L}_2[h]&=\tfrac12{(2\mathtt{p}_{1,1})^{-\frac14}}\left(\big(\tfrac{1}{8}(2\mathtt{p}_{1,1})^{-\frac32}+\mathtt{g}_2\big) h_{1,c}-2\mathtt{q}_1h_{1,c,\star} -2\mathtt{q}_2h_{1,s,\star} -\mathtt{f}_2h_{3,s}-\mathtt{f}_2 h_{1,s}+\mathtt{g}_2h_{3,c} \right),\\
   \mathtt{Q}_1[f,h]&= -\tfrac{1}{4}(2\mathtt{p}_{1,1})^{-\frac14}\sum_{k\geqslant 1}\big( f_{k,c}h_{k+1,c}+f_{k+1,c}h_{k,c}+f_{k,s}h_{k+1,s}+f_{k+1,s}h_{k,s}\big).
\end{align*}
    
\end{lem}
\begin{proof}
Straightforward Fourier computations yield
\begin{align*}
    \fint_{\T}\big(\cos(3\theta)-2\cos(\theta)\big) h(\theta)\cos(\theta)d\theta &=-\tfrac14(h_{4,c}-h_{2,c}), \\
    \fint_{\T}\big(\mathtt{f}_2\sin(2\theta)-\mathtt{g}_2\cos(2\theta)\big)h(\theta)\cos(\theta)d\theta &=\tfrac14\mathtt{f}_2(h_{3,s}+h_{1,s})-\tfrac14\mathtt{g}_2(h_{3,c}+h_{1,c}).
\end{align*}
By Lemma \ref{prop:SandQ}, 
\begin{align*}
   \fint_{\T}\big(\mathcal{H}_{\mathtt{u},0}[h](\theta)+{\partial_\theta}\mathcal{S}[h](\theta)\big)\sin(\theta)d\theta=-\tfrac12\mathtt{g}_3(\varphi)h_{2,c},  
\end{align*}
\begin{align*}
       \fint_{\T}\partial_\theta\mathcal{Q}[h](\theta)\sin(\theta)d\theta &=
       \tfrac{1}{32 }(2\mathtt{p}_{1,1})^{-\frac32} h_{1,c}-\tfrac{1}{2}\mathtt{q}_1h_{1,c,\star} -\tfrac{1}{2}\mathtt{q}_2 h_{1,s,\star}.
  \end{align*}
From \eqref{defL1Q1L2Q2}, we obtain
  \begin{equation*}
\begin{aligned}
  \mathtt{L}_1[h]&=-\tfrac{1}{8\mathtt{p}_{1,1}}\fint_{\T}\big(\cos(3\theta)-2\cos(\theta)\big) h(\theta)\cos(\theta)d\theta\\ &\quad +2{(2\mathtt{p}_{1,1})^{-\frac14}}\fint_{\T}\big(\mathcal{H}_{\mathtt{u},0}[h](\theta)+{\partial_\theta}\mathcal{S}[h](\theta)\big)\sin(\theta)d\theta\\ 
  &=\tfrac{1}{32\mathtt{p}_{1,1}}(h_{4,c}-h_{2,c})-\mathtt{g}_3(2\mathtt{p}_{1,1})^{-\frac14}h_{2,c}.
\end{aligned}
\end{equation*}
Using the identity relating $\mathtt{g}_3$ and $\mathtt{p}_{1,1}$ in \eqref{list-functions},
this simplifies to
 \begin{equation*}
  \mathtt{L}_1[h]=\frac{1}{32\mathtt{p}_{1,1}}(h_{4,c}-3h_{2,c}).
\end{equation*}
For $\mathtt{L}_2[h]$, we similarly obtain
\begin{equation*}
\begin{aligned}
&\mathtt{L}_2[h]=2(2\mathtt{p}_{1,1})^{-\frac14}\bigg(\fint_{\T}\partial_\theta\mathcal{Q}[h](\theta)\sin(\theta)d\theta  -\fint_{\T}\big(\mathtt{f}_2\sin(2\theta)-\mathtt{g}_2\cos(2\theta)\big)h(\theta)\cos(\theta)d\theta\bigg)\\
    &=(2\mathtt{p}_{1,1})^{-\frac14}\bigg(\big(\tfrac{1}{16}(2\mathtt{p}_{1,1})^{-\frac32}+\tfrac12\mathtt{g}_2\big) h_{1,c}-\mathtt{q}_1\, h_{1,c,\star} -\mathtt{q}_2\, h_{1,s,\star} -\tfrac12\mathtt{f}_2h_{3,s}-\tfrac12\mathtt{f}_2 h_{1,s}+\tfrac12\mathtt{g}_2h_{3,c} \bigg).
    \end{aligned}
\end{equation*}
Finally, 
\begin{equation*}
\begin{aligned}
\mathtt{Q}_1[f,h]&= -(2\mathtt{p}_{1,1})^{-\frac14}\fint_{\T} f(\theta)h(\theta) \cos(\theta)d\theta\\
&= -\tfrac{1}{4}(2\mathtt{p}_{1,1})^{-\frac14}\sum_{k\geqslant 1}\big( f_{k,c}h_{k+1,c}+f_{k+1,c}h_{k,c}+f_{k,s}h_{k+1,s}+f_{k+1,s}h_{k,s}\big).
\end{aligned}
\end{equation*}
This ends the proof of the Lemma.
\end{proof}

The next lemma records the Fourier-series expansions of the derivatives of $\mathbf F_1$ evaluated at the reference profile $f_0$. These expansions is used in Theorem~\ref{theo-approx1}.

\begin{lem}\label{prop:dF-computations3}
 Assume that $f_0$ admits the expansion
$$
f_0
 = g_0\cos\theta+\mathtt g_3\cos(3\theta)
 +\varepsilon|\ln\varepsilon|
 \bigl(-2\mathtt g_2\cos(2\theta)+2\mathtt f_2\sin(2\theta)\bigr)
 +O\!\bigl(\varepsilon|\ln\varepsilon|^{\frac12}\bigr).
$$
Then the following asymptotic expansions hold.
\begin{enumerate}
\item[{\rm (1)}] We have
\begin{align*}
\qquad \Pi_{1,{\bf{c}}}(\partial_f {\bf F}_1)(\varepsilon,f_0)[h]&= \Big\{\varepsilon^3 |\ln\varepsilon|\Big(\omega h_{1,c}'+\mathtt{q}_3 h_{1,c,\star}+ \mathtt{f}_2 h_{1,c}-\tfrac12 \mathtt{g}_3 g_0 h_{3,s}- \mathtt{f}_2 h_{3,c}-\mathtt{g}_2 h_{3,s}\Big)\\
 &\qquad +\tfrac14\varepsilon^2  \mathtt{g}_3 ( h_{2,s}+3h_{4,s})+\tfrac14\varepsilon^2 g_0 h_{2,s}
 +\varepsilon^2 \mathcal{M}_1(f_0)[h]\\
 &\qquad +\varepsilon^3|\ln\varepsilon| \mathcal{M}_2(f_0)[h]+O(\varepsilon^3|\ln\varepsilon|^{\frac12}h)\Big\}\cos(\theta).
\end{align*}
Moreover, for the quantities entering $\mathcal M_j$ (see \eqref{defMj}), one has
\begin{align*}
   \mathtt{L}_1[h]+   \mathtt{Q}_1[f_0,h]&= -\tfrac14 (2\mathtt{p}_{1,1})^{-\frac14}g_0 h_{2,c}+\tfrac{1}{64\mathtt{p}_{1,1}}(h_{4,c}-7 h_{2,c})\\
   &\quad -\tfrac{1}{2}(2\mathtt{p}_{1,1})^{-\frac14}\varepsilon|\ln\varepsilon|(-\mathtt{g}_2 h_{3,c}-\mathtt{g}_2h_{1,c}+\mathtt{f}_2 h_{3,s})+O(\varepsilon|\ln\varepsilon|^{\frac12} h),
\end{align*}
and
\begin{align*}
   \mathtt{L}_2[h]&=(2\mathtt{p}_{1,1})^{-\frac14}\bigg(\big(\tfrac{1}{16}(2\mathtt{p}_{1,1})^{-\frac32}+\tfrac12\mathtt{g}_2\big) h_{1,c}-\mathtt{q}_1 h_{1,c,\star}  -\tfrac12\mathtt{f}_2h_{3,s}+\tfrac12\mathtt{g}_2h_{3,c} \bigg)\\ &\quad+O(\varepsilon|\ln\varepsilon|^{\frac12}h).
\end{align*}
\item[{\rm (2)}] Restricting to the first cosine mode of $h$, we obtain \begin{align*}
\qquad \Pi_{1,{\bf{c}}}(\partial_f {\bf F}_1)(\varepsilon,f_0)&[\Pi_{1,{\bf{c}}} h]=\varepsilon^3 |\ln\varepsilon|\, \omega \Big\{  h_{1,c}'+\tfrac{T_0(\lambda,\kappa)}{2\pi}\check{\alpha}(h_{1,c}- h_{1,c,\star})\\
 &+ \tfrac{T_0^2(\lambda)\check{\mathtt{h}}_2 e^{\check{\mathtt{f}}_3(\varphi)}} {64\pi^2 \kappa}\, \int_0^\varphi e^{-\check{\mathtt{f}}_3(s)} (h_{1,c}-h_{1,c,\star}) (s)ds  +O(|\ln\varepsilon|^{-\frac12}h)\Big\}\cos(\theta),
\end{align*}
where $\check{\mathtt{f}}_3$ is given by \begin{equation}\label{check-f3}
  \check{\mathtt{f}}_3(\varphi):= \tfrac{T_0(\lambda,\kappa)}{\pi}\int_0^\varphi \tfrac{\mathtt{y}_1\mathtt{y}_2}{(\mathtt{y}_1^2+\mathtt{y}_2^2)^2}(s)ds,
\end{equation}.
\item[{\rm (3)}] Finally,
\begin{align*}
\qquad \Pi_{1,{\bf{c}}}\Big\{(\partial_f {\bf F}_1)(\varepsilon,0)&[f_0]+\tfrac12 (\partial_f^2 {\bf F}_1)(\varepsilon,0)[f_0,f_0]\Big\}= \varepsilon^3 |\ln\varepsilon| \omega\Big\{ g_0'+\tfrac{T_0}{2\pi}\,\check{\alpha} \,(g_0- g_{0,\star}) \\ &+\,\tfrac{ T_0^2\check{\mathtt{h}}_2(\varphi) e^{\check{\mathtt{f}}_3(\varphi)}} {64\pi^2\kappa}\, \int_0^\varphi e^{-\check{\mathtt{f}}_3(s)} (g_{0}-g_{0,\star})(s)ds+O(|\ln\varepsilon|^{-\frac12})\Big\}\cos(\theta).
\end{align*}
\end{enumerate}
\end{lem}

\begin{proof}
Let $f$ and $h$  be smooth functions with Fourier expansions
\begin{align*}
h(\theta)&=h_{1,c}\cos(\theta)+\sum_{k\geqslant 2} h_{k,c}\cos(k\theta)+h_{k,s}\sin(k\theta),
\\
f(\theta)&=f_{1,c}\cos(\theta)+\varepsilon|\ln\varepsilon|\big(f_{2,c}\cos(2\theta)+f_{2,s}\sin(2\theta)\big)+f_{3,c}\cos(3\theta)\\ &\quad+\varepsilon|\ln\varepsilon|^{\frac12}\Big(f_{3,s}\sin(3\theta)  +\sum_{k\geqslant 4} f_{k,c}\cos(k\theta)+f_{k,s}\sin(k\theta)\Big).
\end{align*}
By Corollary~\ref{cor:asymp-linF1},
\begin{align*}
 & \Pi_{1,{\bf{c}}}(\partial_f {\bf F}_1)(\varepsilon,f)[h]= \varepsilon^3 |\ln\varepsilon|\omega h_{1,c}'\cos(\theta)+\varepsilon^2\Pi_{1,{\bf{c}}} \partial_\theta\Big\{ h\, V_0(\varepsilon,f)\Big\}
      +\varepsilon^2\Pi_{1,{\bf{c}}} \partial_\theta \mathcal{S}[h]+\varepsilon^2\Pi_{1,{\bf{c}}} \mathcal{H}_{\mathtt{u},0}[h]\\ & \quad+\varepsilon^2 \mathcal{M}_1(f)[h]\cos(\theta)+\varepsilon^3|\ln\varepsilon| \mathcal{M}_2(f)[h]\cos(\theta)
  +{\varepsilon^3}|\ln\varepsilon| \Pi_{1,{\bf{c}}}\partial_\theta \mathcal{Q}[h]+\varepsilon^3 \Pi_{1,{\bf{c}}}{R}(\varepsilon,f)[h].
\end{align*}
Moreover, Lemma~\ref{prop:SandQ} yields
\begin{align*}
\nonumber \Pi_{1,{\bf{c}}}\partial_\theta\mathcal{Q}[h](\varphi,\theta)&=\mathtt{q}_3 h_{1,c,\star} \cos(\theta),
\\  \Pi_{1,{\bf{c}}}(\partial_\theta \mathcal{S}+\mathcal{H}_{\mathtt{u},0})[h]&=2\mathtt{g}_3 h_{2,s}\cos(\theta).
  \end{align*}
For the transport term we compute
\begin{align*}
 &\Pi_{1,{\bf{c}}}\partial_\theta\big\{ h(\theta) V_0(\varepsilon,\varepsilon)\big\}=\Big\{-\tfrac32 \mathtt{g}_3 h_{2,s}+\tfrac12 \mathtt{g}_3 h_{4,s}+\tfrac14f_{1,c} h_{2,s}+\tfrac14f_{3,c}(h_{4,s}-h_{2,s})\\ &+\tfrac12\varepsilon|\ln\varepsilon|\Big( h_{1,c}\mathtt{f}_2 -\mathtt{f}_2 h_{3,c}-\mathtt{g}_2 h_{3,s}+\tfrac12 h_{1,c} f_{2,s} +\tfrac12(f_{2,c}h_{3,s}-h_{3,c}f_{2,s})\Big) +O\big(\varepsilon|\ln\varepsilon|^{\tfrac12}\big) \Big\}\cos(\theta).
\end{align*}
Collecting these contributions, we find
\begin{align}\label{prop:dF-computations}
\nonumber \Pi_{1,{\bf{c}}}&(\partial_f {\bf F}_1)(\varepsilon,f)[h]= \Big[\varepsilon^3 |\ln\varepsilon|\Big(\omega h_{1,c}'+\mathtt{q}_3 h_{1,c,\star}+\tfrac12 \mathtt{f}_2 h_{1,c}-\tfrac12 \mathtt{g}_3 f_{1,c} h_{3,s}-\tfrac12 \mathtt{f}_2 h_{3,c}-\tfrac12\mathtt{g}_2 h_{3,s}\\ \nonumber
 &+\tfrac14(f_{2,c}h_{3,s}-h_{3,c}f_{2,s})+\tfrac14  h_{1,c} f_{2,s}\Big)+\varepsilon^2\Big( \tfrac12 \mathtt{g}_3 (h_{2,s}+h_{4,s})+\tfrac14 f_{1,c} h_{2,s}+\tfrac14f_{3,c}(h_{4,s}-h_{2,s})\Big)\nonumber \\
 &+\varepsilon^2 \mathcal{M}_1(f)[h]+\varepsilon^3|\ln\varepsilon| \mathcal{M}_2(f)[h]+O(\varepsilon^3|\ln\varepsilon|^{\frac12})\Big]\cos(\theta).
\end{align}
Here,
\begin{equation*}
 \begin{aligned}
     \mathcal{M}_j(f)[h](\varphi)&:=\tfrac{\mathtt{h}_2(\varphi) e^{\mathtt{f}_3(\varphi)}} {2\omega (2\mathtt{p}_{1,1})^{\frac14}}\, \int_0^\varphi e^{-\mathtt{f}_3(s)} \big(\mathtt{L}_j[h]-(\mathtt{L}_{j}[h])_{\star}\big)(s)ds +(2-j)\mathcal{M}_{3}[f,h],\\
     \mathcal{M}_{3}[f,h](\varphi)(\varphi)&:=\tfrac{\mathtt{h}_2(\varphi) e^{\mathtt{f}_3(\varphi)}} {2\omega (2\mathtt{p}_{1,1})^{\frac14}}\, \int_0^\varphi e^{-\mathtt{f}_3(s)} \big(\mathtt{Q}_{1}[f,h]-(\mathtt{Q}_{1}[f,h])_{\star}\big)(s)ds,
\end{aligned}   
\end{equation*}
 and the expressions of $\mathtt{L}_1[h]$, $\mathtt{L}_2[h]$ and $\mathtt{Q}_{1}[f,h]$ are given by Lemma \ref{lem-FourLQ}.
 
 \smallskip
 
\noindent\textbf{Proof of (1).}
Taking $f=f_0$ in \eqref{prop:dF-computations} and using
$$
(f_0)_{2,s}=2\mathtt f_2+O\!\bigl(\varepsilon|\ln\varepsilon|^{\frac12}\bigr),\qquad
(f_0)_{2,c}=-2\mathtt g_2+O\!\bigl(\varepsilon|\ln\varepsilon|^{\frac12}\bigr),\qquad
(f_0)_{3,c}=\mathtt g_3+O\!\bigl(\varepsilon|\ln\varepsilon|^{\frac12}\bigr),
$$
we obtain
\begin{align*}
 \Pi_{1,{\bf{c}}}&(\partial_f {\bf F}_1)(\varepsilon,f_0)[h]= \Big\{\varepsilon^3 |\ln\varepsilon|\Big(\omega h_{1,c}'+\mathtt{q}_3 h_{1,c,\star}+\tfrac12 \mathtt{f}_2 h_{1,c}-\tfrac12 \mathtt{g}_3 g_0 h_{3,s}-\tfrac12 \mathtt{f}_2 h_{3,c}-\tfrac12\mathtt{g}_2 h_{3,s}\\
 &-\tfrac12 (h_{3,c}\mathtt{f}_2+h_{3,s}\mathtt{g}_2)+\tfrac12  h_{1,c} \mathtt{f}_2\Big)+\varepsilon^2 \tfrac12 \mathtt{g}_3 (h_{2,s}+h_{4,s})+\tfrac14\varepsilon^2 g_0 h_{2,s}+\varepsilon^2\tfrac14 \mathtt{g}_3(h_{4,s}- h_{2,s})\\
 &
 +\varepsilon^2 \mathcal{M}_1(f_0)[h]+\varepsilon^3|\ln\varepsilon| \mathcal{M}_2(f_0)[h]+O(\varepsilon^3|\ln\varepsilon|^{\frac12}h)\Big\}\cos(\theta).
\end{align*}
It remains to identify the quantities entering $\mathcal M_j$. By Lemma~\ref{lem-FourLQ},
\begin{align*}
   \mathtt{L}_1[h]+   \mathtt{Q}_1[f_0,h]=& -\tfrac14 (2\mathtt{p}_{1,1})^{-\frac14}g_0 h_{2,c}+\frac{1}{64\mathtt{p}_{1,1}}(h_{4,c}-7 h_{2,c})\\
   &-\tfrac{1}{2}(2\mathtt{p}_{1,1})^{-\frac14}\varepsilon|\ln\varepsilon|(-\mathtt{g}_2 h_{3,c}-\mathtt{g}_2h_{1,c}+\mathtt{f}_2 h_{3,s})+O(\varepsilon|\ln\varepsilon|^{\frac12}h).
\end{align*}
\begin{align*}
   \mathtt{L}_2[h]=&(2\mathtt{p}_{1,1})^{-\frac14}\bigg(\big(\tfrac{1}{16}(2\mathtt{p}_{1,1})^{-\frac32}+\tfrac12\mathtt{g}_2\big) h_{1,c}-\mathtt{q}_1h_{1,c,\star}  -\tfrac12\mathtt{f}_2h_{3,s}+\tfrac12\mathtt{g}_2h_{3,c} \bigg)+O\!\bigl(|\ln\varepsilon|^{-\frac12}\bigr),\\
   \mathtt{Q}_1[f_0,h]=& -\tfrac{1}{4}(2\mathtt{p}_{1,1})^{-\frac14}(g_0 h_{2,c}+\mathtt{g}_3h_{4,c}+\mathtt{g}_3 h_{2,c})\\
   &-\tfrac{1}{2}(2\mathtt{p}_{1,1})^{-\frac14}\varepsilon|\ln\varepsilon|(-\mathtt{g}_2 h_{3,c}-\mathtt{g}_2h_{1,c}+\mathtt{f}_2 h_{3,s})+O(\varepsilon^3|\ln\varepsilon|^{\frac12}h).
\end{align*}
which is precisely the claim in (1).

\smallskip
\noindent\textbf{Proof of (2).} Applying (1) to $\Pi_{1,\mathbf c}h$, we get
\begin{align*}
 \Pi_{1,{\bf{c}}}(\partial_f {\bf F}_1)(\varepsilon,f_0)[\Pi_{1,{\bf{c}}} h]&= \omega\Big\{\varepsilon^3 |\ln\varepsilon|\Big( h_{1,c}'+\tfrac{\mathtt{q}_3}{\omega} h_{1,c,\star}+\tfrac{\mathtt{f}_2}{\omega} h_{1,c}+\tfrac{1}{\omega}\mathcal{M}_2(f_0)[\Pi_{1,{\bf{c}}} h]\Big) \\
 &\quad+\varepsilon^2 \tfrac{1}{\omega}\mathcal{M}_1(f_0)\Pi_{1,{\bf{c}}}[h]+O(\varepsilon^3|\ln\varepsilon|^{\frac12})\Big\}\cos(\theta).
\end{align*}
From Lemma \ref{lem-lam}, one has
$$
\omega(\varepsilon,\lambda)
=
\frac{2\pi}{\sqrt{2\kappa}\,T_0(\lambda,\kappa)}
+O\!\bigl(|\ln\varepsilon|^{-\frac12}\bigr)
,
$$
and by Lemma \ref{lem-functions},
$$
\tfrac{\mathtt{q}_3}{\omega}=\tfrac{T_0(\lambda,\kappa)}{2\pi} \check{\alpha}+O(|\ln\varepsilon|^{-\frac12}),\quad \tfrac{\mathtt{f}_2}{\omega}=\tfrac{T_0(\lambda,\kappa)}{2\pi} \check{\alpha}+O(|\ln\varepsilon|^{-\frac12}). 
$$
For the terms inside $\mathcal{M}_1+\varepsilon|\ln\varepsilon| \mathcal{M}_2$, using Lemma \ref{lem-functions},   one finds
\begin{align*}
   \mathtt{Q}[f_0,\Pi_{1,{\bf{c}}} h]&:=\big(\mathtt{L}_1\Pi_{1,{\bf{c}}}+\varepsilon|\ln\varepsilon| \mathtt{L}_2\big)[\Pi_{1,{\bf{c}}} h]+   \mathtt{Q}_1[f_0,\Pi_{1,{\bf{c}}} h]\\ &=\varepsilon|\ln\varepsilon|(2\mathtt{p}_{1,1})^{-\frac14}\Big( \mathtt{g}_2h_{1,c}+{\tfrac{1}{16}}(2\mathtt{p}_{1,1})^{-\frac32} h_{1,c}-\mathtt{q}_1 h_{1,c,\star}\Big)\\
   &=\varepsilon|\ln\varepsilon|(2\kappa)^{-\frac14}\Big(  \tfrac{1}{4}(2\kappa)^{-\frac32}h_{1,c}-\tfrac12(2\kappa)^{-\frac12}\check{\mathtt{h}}_2h_{1,c} -\tfrac12(2\kappa)^{-\frac12}\check{\mathtt{h}}_2h_{1,c,\star}\Big)\\ &\quad  +O(|\ln\varepsilon|^{-\frac12}).
\end{align*}
Subtracting the $\star$-shift and using the symmetry properties of the second part of Corollary \ref{cor-F0} gives 
\begin{align*}
  \mathtt{Q}[f_0,\Pi_{1,{\bf{c}}} h]-\big(\mathtt{Q}[f_0,\Pi_{1,{\bf{c}}} h]\big)_\star 
   &={\tfrac{1}{4}}\varepsilon|\ln\varepsilon|(2\kappa)^{-\frac74}(h_{1,c}-h_{0,\star}) +O(|\ln\varepsilon|^{-\frac12}).
\end{align*}
Moreover,
$$
\mathtt f_3(\varphi)
=
\frac{T_0(\lambda,\kappa)}{\pi}
\int_0^\varphi
\frac{\mathtt y_1\mathtt y_2}{(\mathtt y_1^2+\mathtt y_2^2)^2}\,d\tau
+O\!\bigl(|\ln\varepsilon|^{-\frac12}\bigr)
=
\check{\mathtt f}_3(\varphi)+O\!\bigl(|\ln\varepsilon|^{-\frac12}\bigr),
$$
and Lemma~\ref{lem-lam} gives
$$
\tfrac{1} {2\omega^2 }(2\mathtt{p}_{1,1})^{-\frac14}=\tfrac{T^2} {2(2\pi)^2 }(2\kappa)^{\frac34}
+O\!\bigl(|\ln\varepsilon|^{-\frac12}\bigr).
$$
Therefore,
\begin{equation*}
 \begin{aligned}
     \tfrac{1}{\omega}\big(\mathcal{M}_1(f_0)+\varepsilon|\ln\varepsilon|\mathcal{M}_2(f_0)\big)[h]&=\varepsilon|\ln\varepsilon|\tfrac{T_0^2(\lambda)\check{\mathtt{h}}_2 e^{\check{\mathtt{f}}_3(\varphi)}} {16(2\pi)^2 \kappa}\, \int_0^\varphi e^{-\check{\mathtt{f}}_3(s)} (h_{1,c}-h_{0,\star}) (s)ds \\ &+O(|\ln\varepsilon|^{-\frac12}).
\end{aligned}   
\end{equation*}
This concludes the proof of (2).

\smallskip

\noindent\textbf{Proof of (3).}
Now take $f=0$ and $h=f_0$ in \eqref{prop:dF-computations}. Since
$$
(f_0)_{2,s}=2\mathtt f_2+O\!\bigl(\varepsilon|\ln\varepsilon|^{\frac12}\bigr),\qquad
(f_0)_{2,c}=-2\mathtt g_2+O\!\bigl(\varepsilon|\ln\varepsilon|^{\frac12}\bigr),\qquad
(f_0)_{3,c}=\mathtt g_3+O\!\bigl(\varepsilon|\ln\varepsilon|^{\frac12}\bigr),
$$
we have
\begin{align*}
\Pi_{1,\mathbf c}(\partial_f \mathbf F_1)(\varepsilon,0)[f_0]
&=
\Big[
\varepsilon^3|\ln\varepsilon|
\Big(
\omega g_0'
+\mathtt q_3 g_{0,\star}
+\tfrac12\mathtt f_2 g_0
+\tfrac12\mathtt f_2\mathtt g_3
\Big)
\\
&\qquad
+\varepsilon^2\mathcal M_1(0)[f_0]
+\varepsilon^3|\ln\varepsilon|\,\mathcal M_2(0)[f_0]
+O\!\bigl(\varepsilon^3|\ln\varepsilon|^{\frac12}\bigr)
\Big]\cos\theta.
\end{align*}
Here
\begin{align*}
\mathcal M_1(0)[f_0]
&=
\tfrac{\mathtt h_2 e^{\mathtt f_3(\varphi)}}{2\omega(2\mathtt p_{1,1})^{\frac14}}
\int_0^\varphi e^{-\mathtt f_3(s)}
\bigl(\mathtt L_1[f_0]-(\mathtt L_1[f_0])_\star\bigr)(s)\,ds,
\\
\mathcal M_2(0)[f_0]
&=
\tfrac{\mathtt h_2 e^{\mathtt f_3(\varphi)}}{2\omega(2\mathtt p_{1,1})^{\frac14}}
\int_0^\varphi e^{-\mathtt f_3(s)}
\bigl(\mathtt L_2[f_0]-(\mathtt L_2[f_0])_\star\bigr)(s)\,ds.
\end{align*}
Using Lemmas~\ref{lem-functions} and \ref{lem-FourLQ}, we compute
\begin{align*}
\mathtt L_1[f_0]
&=
\varepsilon|\ln\varepsilon|\frac{3}{16\mathtt p_{1,1}}\mathtt g_2
=
\varepsilon|\ln\varepsilon|(2\kappa)^{-\frac32}\frac{3}{16}
\Big(\frac{3}{16\kappa}-\check{\mathtt h}_2\Big)
+O\!\bigl(|\ln\varepsilon|^{-\frac12}\bigr),
\end{align*}
and
\begin{align*}
\mathtt L_2[f_0]
&=
\tfrac12(2\mathtt p_{1,1})^{-\frac14}
\bigg[
\Big(\tfrac{1}{8}(2\mathtt p_{1,1})^{-\frac32}+\mathtt g_2\Big)g_0
-2\mathtt q_1 g_{0,\star}
+\mathtt g_2\mathtt g_3
\bigg]
+O\!\bigl(|\ln\varepsilon|^{-\frac12}\bigr)
\\
&=
\tfrac14(2\kappa)^{-\frac34}
\bigg[
\Big(\tfrac{5}{16\kappa}-\check{\mathtt h}_2\Big)g_0
-2\check{\mathtt h}_2 g_{0,\star}
+\tfrac18(2\kappa)^{-\frac34}
\Big(\tfrac{3}{16\kappa}-\check{\mathtt h}_2\Big)
\bigg]
+O\!\bigl(|\ln\varepsilon|^{-\frac12}\bigr).
\end{align*}
In particular, since $\check{\mathtt h}_2=\check{\mathtt h}_{2,\star}$, we have 
\begin{align*}
\mathtt L_1[f_0]-(\mathtt L_1[f_0])_\star
&=
+O\!\bigl(|\ln\varepsilon|^{-\frac12}\bigr),
\\
\mathtt L_2[f_0]-(\mathtt L_2[f_0])_\star
&=
\tfrac14(2\kappa)^{-\frac34}
\Big(\tfrac{5}{16\kappa}+\check{\mathtt h}_2\Big)
(g_0-g_{0,\star})
+O\!\bigl(\varepsilon|\ln\varepsilon|^{\frac12}\bigr).
\end{align*}
Moreover,
$$
\mathtt f_3(\varphi)
=
\frac{T_0(\lambda,\kappa)}{\pi}
\int_0^\varphi
\frac{\mathtt y_1\mathtt y_2}{(\mathtt y_1^2+\mathtt y_2^2)^2}\,d\tau
+O\!\bigl(|\ln\varepsilon|^{-\frac12}\bigr)
=
\check{\mathtt f}_3(\varphi)+O\!\bigl(|\ln\varepsilon|^{-\frac12}\bigr),
$$
and Lemma~\ref{lem-lam} gives
$$
\omega(\varepsilon,\lambda)
=
\frac{2\pi}{\sqrt{2\kappa}\,T_0(\lambda,\kappa)}
+O\!\bigl(|\ln\varepsilon|^{-\frac12}\bigr).
$$
Therefore,
$$
\frac{T_0(\lambda,\kappa)}{\omega\sqrt{2\kappa}}
=
\frac{T_0^2(\lambda)}{2\pi}
+O\!\bigl(|\ln\varepsilon|^{-\frac12}\bigr),
\qquad
\frac{1}{\omega\sqrt{2\kappa}}
=
\frac{T_0(\lambda,\kappa)}{2\pi}
+O\!\bigl(|\ln\varepsilon|^{-\frac12}\bigr).
$$
Hence
$$
\mathcal M_1(0)[f_0]
=
O\!\bigl(\varepsilon|\ln\varepsilon|^{\frac12}\bigr),
$$
and
$$
\mathcal M_2(0)[f_0]
=
\tfrac{T_0(\lambda,\kappa)\check{\mathtt h}_2(\varphi)e^{\check{\mathtt f}_3(\varphi)}}{16\pi\sqrt{2\kappa}}
\int_0^\varphi e^{-\check{\mathtt f}_3(s)}
\Big(\tfrac{5}{16\kappa}+\check{\mathtt h}_2\Big)
(g_0-g_{0,\star})(s)\,ds
+O\!\bigl(\varepsilon|\ln\varepsilon|^{\frac12}\bigr).
$$
Consequently,
\begin{align*}
\Pi_{1,\mathbf c}(\partial_f \mathbf F_1)(\varepsilon,0)[f_0]
&=
\varepsilon^3|\ln\varepsilon|\,\omega
\Big(
g_0'
-\tfrac{T_0}{2\pi}\check\alpha\,g_{0,\star}
+\tfrac{T_0}{4\pi}\check\alpha\,g_0
+\tfrac{T_0}{32\pi}(2\kappa)^{-\frac34}\check\alpha
\Big)\cos\theta
\\
&\qquad
+\varepsilon^3|\ln\varepsilon|
\tfrac{T_0^2\check{\mathtt h}_2(\varphi)e^{\check{\mathtt f}_3(\varphi)}}{32\pi^2}
\int_0^\varphi e^{-\check{\mathtt f}_3(s)}
\Big(\tfrac{5}{16\kappa}+\check{\mathtt h}_2\Big)
(g_0-g_{0,\star})(s)\,ds\;\cos\theta
\\
&\qquad
+O\!\bigl(\varepsilon^3|\ln\varepsilon|^{\frac12}\bigr)\cos\theta.
\end{align*}
On the other hand,
\begin{align*}
\Pi_{1,\mathbf c}(\partial_f^2 \mathbf F_1)(\varepsilon,0)[f_0,f_0]
&=
\varepsilon^3|\ln\varepsilon|
\bigl(\mathtt f_2 g_0-\mathtt f_2\mathtt g_3\bigr)\cos\theta
+\varepsilon^2\mathcal M_3[f_0,f_0]
+O\!\bigl(\varepsilon^3|\ln\varepsilon|^{\frac12}\bigr)\cos\theta
\\
&=
\varepsilon^3|\ln\varepsilon|
\Big(
\tfrac{1}{\sqrt{2\kappa}}\check\alpha\,g_0
-\tfrac18\check\alpha(2\kappa)^{-\frac54}
\Big)\cos\theta
+\varepsilon^2\mathcal M_3[f_0,f_0]
\\
&\qquad
+O\!\bigl(\varepsilon^3|\ln\varepsilon|^{\frac12}\bigr)\cos\theta,
\end{align*}
where
$$
\mathcal M_3[f_0,f_0]
=
\tfrac{\mathtt h_2 e^{\mathtt f_3(\varphi)}}{2\omega(2\mathtt p_{1,1})^{\frac14}}
\int_0^\varphi e^{-\mathtt f_3(s)}
\bigl(\mathtt Q_1[f_0,f_0]-(\mathtt Q_1[f_0,f_0])_\star\bigr)(s)\,ds.
$$
Furthermore,
\begin{align*}
\mathtt Q_1[f_0,f_0]
&=
\varepsilon|\ln\varepsilon|(2\mathtt p_{1,1})^{-\frac14}
\bigl(\mathtt g_2 g_0+\mathtt g_2\mathtt g_3\bigr)
+O\!\bigl(\varepsilon|\ln\varepsilon|^{\frac12}\bigr)
\\
&=
\tfrac12\varepsilon|\ln\varepsilon|(2\kappa)^{-\frac34}
\Big(\tfrac{3}{16\kappa}-\check{\mathtt h}_2\Big)
\Big(g_0+\tfrac18(2\kappa)^{-\frac34}\Big)
+O\!\bigl(\varepsilon|\ln\varepsilon|^{\frac12}\bigr),
\end{align*}
and therefore
\begin{align*}
\mathtt Q_1[f_0,f_0]-(\mathtt Q_1[f_0,f_0])_\star
&=
\tfrac12\varepsilon|\ln\varepsilon|(2\kappa)^{-\frac34}
\Big(\tfrac{3}{16\kappa}-\check{\mathtt h}_2\Big)
(g_0-g_{0,\star})
+O\!\bigl(\varepsilon|\ln\varepsilon|^{\frac12}\bigr).
\end{align*}
Thus
$$
\frac{1}{2\omega}\mathcal M_3[f_0,f_0]
=
\varepsilon|\ln\varepsilon|
\tfrac{T_0^2\check{\mathtt h}_2(\varphi)e^{\check{\mathtt f}_3(\varphi)}}{32\pi^2}
\int_0^\varphi e^{-\check{\mathtt f}_3(s)}
\Big(\tfrac{3}{16\kappa}-\check{\mathtt h}_2\Big)
(g_0-g_{0,\star})(s)\,ds
+O\!\bigl(\varepsilon|\ln\varepsilon|^{\frac12}\bigr).
$$
Combining the previous identities, we arrive at
\begin{align*}
&\Pi_{1,\mathbf c}(\partial_f \mathbf F_1)(\varepsilon,0)[f_0]
+\tfrac12\Pi_{1,\mathbf c}(\partial_f^2 \mathbf F_1)(\varepsilon,0)[f_0,f_0]
\\
&\qquad
=
\varepsilon^3|\ln\varepsilon|\,\omega
\Bigg(
g_0'
+\tfrac{T_0}{2\pi}\check\alpha\,(g_0-g_{0,\star})
+\frac{T_0^2\check{\mathtt h}_2(\varphi)e^{\check{\mathtt f}_3(\varphi)}}{64\pi^2\kappa}
\int_0^\varphi e^{-\check{\mathtt f}_3(s)}(g_0-g_{0,\star})(s)\,ds
\\
&\hspace{5.3cm}
+O\!\bigl(|\ln\varepsilon|^{-\frac12}\bigr)
\Bigg)\cos\theta,
\end{align*}
which proves (3).
\end{proof}

\section{Integral formulas}\label{appendix:integrals}
In this section, we provide useful integrals needed in  the computation of  linearized operator, see for instance \cite[Appendix A]{CCGS-active}.

\begin{lem}\label{lem:integrals}
The following integrals reads as
\begin{align}
\nonumber \int_{0}^{2\pi}\int_{0}^{1}\ln(|1-\rho e^{i\eta}|)\rho d\rho d\eta&=0,\\
\int_{0}^{2\pi}\int_{0}^{1}\ln(|1-\rho e^{i\eta}|^2)e^{i\eta}\rho^2 d\rho d\eta&=-\frac{\pi}{2},\label{id1}\\
\int_{0}^{2\pi}\int_{0}^{1}\ln(|1-\rho e^{i\eta}|^2)\cos(2\eta)\rho^3 d\rho d\eta&=-\frac{\pi}{6},\nonumber
\\ 
\label{id:integralV11}
   \int_0^1\int_0^{2\pi}\frac{1-\rho\cos(\eta)}{1+\rho^2-2\rho\cos(\eta)}\rho d\rho d\eta&=\pi,
\end{align}
\begin{align}
  \nonumber  
  \fint_{\T}\int_0^{1}\frac{(\cos(\theta)-\rho\cos(\eta))^2}{|e^{i\theta}-\rho e^{i\eta}|^2}\rho d\rho d\eta&=\frac{1}{4}+\frac{1}{8}\cos(2\theta),
  \\
   \nonumber
\fint_{\T}\int_0^{1}\frac{(\cos(\theta)+\rho\cos(\eta))(\sin(\theta)-\rho\sin(\eta))}{|e^{i\theta}-\rho e^{i\eta}|^2}\rho d\rho d\eta&=\frac{3}{8}\sin(2\theta),\\  
 \nonumber \fint_{\T}\int_{0}^{2\pi}\int_0^{1}\frac{(\cos(\theta)+\rho\cos(\eta))(\cos(\theta)-\rho\cos(\eta))^3}{|e^{i\theta}-\rho e^{i\eta}|^4}\rho d\rho d\eta&=\frac{3}{16}+\frac{1}{4}\cos(2\theta), 
 \\
 \nonumber \fint_{\T}\int_0^{1}\frac{(\cos(\theta)+\rho\cos(\eta))(\cos(\theta)-\rho\cos(\eta))(\sin(\theta)-\rho\sin(\eta))^2}{|e^{i\theta}-\rho e^{i\eta}|^4}\rho d\rho d\eta&=\frac{1}{16}+\frac18\cos(2\theta),
  \\
\nonumber  \fint_{\T}\int_0^{1}\frac{(\cos(\theta)+\rho\cos(\eta))(\cos(\theta)-\rho\cos(\eta))^2(\sin(\theta)-\rho\sin(\eta))}{|e^{i\theta}-\rho e^{i\eta}|^4}\rho d\rho d\eta&=\frac{1}{16}\sin(2\theta), 
\\
  \fint_{\T}\int_{0}^{1}\bigg(\frac{ (\cos(\theta)+\rho\cos(\eta))(\cos(\theta)-\rho\cos(\eta))^2}{|e^{i\theta}-\rho e^{i\eta}|^2}\bigg)\rho d\rho d\eta =\frac{1}{4}\cos(\theta)&+\frac{1}{12}\cos(3\theta).
  \label{id14}
\end{align}

\end{lem}

\end{document}

%% file: dibujo-puntos.pdf_tex
\begingroup%
  \makeatletter%
  \providecommand\color[2][]{%
    \errmessage{(Inkscape) Color is used for the text in Inkscape, but the package 'color.sty' is not loaded}%
    \renewcommand\color[2][]{}%
  }%
  \providecommand\transparent[1]{%
    \errmessage{(Inkscape) Transparency is used (non-zero) for the text in Inkscape, but the package 'transparent.sty' is not loaded}%
    \renewcommand\transparent[1]{}%
  }%
  \providecommand\rotatebox[2]{#2}%
  \newcommand*\fsize{\dimexpr\f@size pt\relax}%
  \newcommand*\lineheight[1]{\fontsize{\fsize}{#1\fsize}\selectfont}%
  \ifx\svgwidth\undefined%
    \setlength{\unitlength}{351.71747643bp}%
    \ifx\svgscale\undefined%
      \relax%
    \else%
      \setlength{\unitlength}{\unitlength * \real{\svgscale}}%
    \fi%
  \else%
    \setlength{\unitlength}{\svgwidth}%
  \fi%
  \global\let\svgwidth\undefined%
  \global\let\svgscale\undefined%
  \makeatother%
  \begin{picture}(1,0.65237199)%
    \lineheight{1}%
    \setlength\tabcolsep{0pt}%
    \put(0,0){\includegraphics[width=\unitlength,page=1]{dibujo-puntos.pdf}}%
    \put(0.85556276,0.28727227){\color[rgb]{0,0,0}\makebox(0,0)[lt]{\lineheight{1.25}\smash{\begin{tabular}[t]{l}$\varrho=\frac12 r^2$\end{tabular}}}}%
    \put(0.50041755,0.29190338){\color[rgb]{0,0,0}\makebox(0,0)[lt]{\lineheight{1.25}\smash{\begin{tabular}[t]{l}$\kappa$\end{tabular}}}}%
    \put(0.71071091,0.10840791){\color[rgb]{0,0,0}\makebox(0,0)[lt]{\lineheight{1.25}\smash{\begin{tabular}[t]{l}$P_1$\end{tabular}}}}%
    \put(0.30350582,0.54505603){\color[rgb]{0,0,0}\makebox(0,0)[lt]{\lineheight{1.25}\smash{\begin{tabular}[t]{l}$P_2$\end{tabular}}}}%
    \put(0.40611162,0.38936076){\color[rgb]{0,0,0}\makebox(0,0)[lt]{\lineheight{1.25}\smash{\begin{tabular}[t]{l}$(2\kappa)^\frac12\mathtt{x}_1(\tau)|\ln\varepsilon|^{-\frac12}$\end{tabular}}}}%
    \put(0.00797118,0.27044602){\color[rgb]{0,0,0}\rotatebox{89.96494205}{\makebox(0,0)[lt]{\lineheight{1.25}\smash{\begin{tabular}[t]{l}$\mathtt{x}_2(\tau)|\ln\varepsilon|^{-\frac12}$\end{tabular}}}}}%
    \put(0.0252591,0.62951179){\color[rgb]{0,0,0}\makebox(0,0)[lt]{\lineheight{1.25}\smash{\begin{tabular}[t]{l}$z$\end{tabular}}}}%
    \put(0,0){\includegraphics[width=\unitlength,page=2]{dibujo-puntos.pdf}}%
    \put(0.2570384,0.50776142){\color[rgb]{0,0,0}\makebox(0,0)[lt]{\lineheight{1.25}\smash{\begin{tabular}[t]{l}$\varepsilon$\end{tabular}}}}%
  \end{picture}%
\endgroup%